\documentclass[b5paper]{book}
\usepackage{geometry}
\usepackage[noinfo,width=176truemm,height=246truemm,center]{crop}

\usepackage{ccicons}
\usepackage{cite}
\usepackage{caption}
\usepackage{xfrac}
\usepackage{xparse}
\usepackage{relsize}
\usepackage{colortbl}
\usepackage{bussproofs}
\usepackage{xcolor}
\usepackage{graphicx}
\usepackage{makeidx}
\usepackage{xstring}
\usepackage{braket}
\usepackage{mathtools}
\usepackage{multicol}
\usepackage{marginnote}
\usepackage{ifthen}
\usepackage{xifthen}
\usepackage{paralist}
\usepackage{calc}
\usepackage{tikz}
\usepackage{tikz-cd}
\usetikzlibrary{calc,arrows,matrix,decorations.pathreplacing}
\usepackage{picins}
\usepackage{amssymb}
\usepackage{amsmath}
\usepackage{mathrsfs}
\usepackage{nicefrac}
\usepackage{stmaryrd}
\usepackage{sectsty}
\usepackage{xr}
\usepackage{tikz}
\usepackage[all]{xypic}
\usepackage{hyperref} 

\newcommand{\doi}[1]{\href{https://doi.org/#1}{doi:#1}}

\newcommand{\purl}[1]{\href{#1}{#1}}

\newcommand{\spacingfix}[0]{\vspace{-1\baselineskip}\ignorespaces}

\newcommand{\IGNORE}[1]{}

\newcommand{\colofon}[4]{%
\vspace*{8em}
\noindent
\begin{center}
    \textcolor{darkblue}{\textsf{\textbf{ Identifiers}}}
    \vspace{.5em}

\href{https://hdl.handle.net/#4}{hdl:~#4}
    \vspace{.5em}

\href{https://arxiv.org/abs/#3}{arXiv:~#3}
    \vspace{.5em}

\textsc{isbn:}~#2
\end{center}
\vspace{3em}
\begin{center}
    \textcolor{darkblue}{\textsf{\textbf{ Persistent links}}}

\vspace{.5em}
\purl{https://arxiv.org/abs/#3}

\vspace{.5em}
\purl{https://doi.org/#4}

\vspace{.5em}
\purl{https://hdl.handle.net/#4}
\end{center}
\vspace{3em}

\begin{center}
    \textcolor{darkblue}{\textsf{\textbf{Source code}}}
\vspace{.5em}

\begin{tabular}{rl}
\LaTeX{} & \purl{https://github.com/westerbaan/theses} \\
cover & \purl{https://github.com/westerbaan/ndpt}
\end{tabular}

\end{center}
\vfill{}

\noindent
Printed by GVO drukkers \& vormgevers B.V., Ede, 
\purl{https://proefschriften.nl}.

\vspace{1em}
\noindent
Where applicable,
\ccCopy{}~2019 #1,
\ccLogo{}\,\ccAttribution{}  available under \textsc{cc by}, \cite{ccby40}.
\newpage
}

\definecolor{pitchblack}{cmyk}{1 1 1 1}
\definecolor{pitchgray}{cmyk}{.5 .5 .5 .5}

\colorlet{darkgreen}{green!50!pitchblack}
\colorlet{darkblue}{blue!75!pitchblack}
\colorlet{lightblue}{blue!50!white}
\colorlet{lightgray}{pitchgray!50!white}
\colorlet{darkgray}{pitchgray!75!pitchblack}
\colorlet{darkred}{red!75!pitchblack}

\allsectionsfont{\sffamily\color{darkblue}}

\newcommand{\textPointHeaderI}[1]{\textcolor{darkblue}{\textbf{\textsf{#1}}}}
\newcommand{\textPointHeaderII}[1]{\textcolor{darkblue}{\textsf{#1}}}
\newcommand{\textPointHeaderIII}[1]{\textcolor{lightgray}{\textsf{#1}}}
\newcommand{\textParsecNumber}[1]{\textcolor{darkblue}{\textbf{\textsf{#1}}}}
\newcommand{\textPointNumberI}[1]{\textcolor{darkblue}{\textsf{#1}}}
\newcommand{\textPointNumberII}[1]{\textcolor{lightblue}{\textsf{#1}}}
\newcommand{\textPointNumberIII}[1]{\textcolor{lightgray}{\textsf{#1}}}
\newcommand{\textDefine}[1]{\textcolor{darkblue}{#1}}
\newcommand{\textSref}[1]{\textsf{#1}}

\newcommand{\Define}[1]{\textDefine{#1}}

\newcommand{\grayed}[1]{\textcolor{darkgray}{#1}}

\renewcommand{\leq}{\leqslant}
\renewcommand{\geq}{\geqslant}
\renewcommand{\nleq}{\not\leqslant}

\newcommand{\Cat}[1]{\mathbf{#1}}

\newcommand{\haW}[1]{\Cat{haW}^*_{\text{\textsc{#1}}}}
\newcommand{\W}[1]{\Cat{W}^*_{\text{\textsc{#1}}}}
\newcommand{\dW}[1]{\Cat{dW}^*_{\text{\textsc{#1}}}}
\newcommand{\cCstar}[1]{\Cat{cC}^*_{\text{\textsc{#1}}}}
\newcommand{\Cstar}[1]{\Cat{C}^*_{\text{\textsc{#1}}}}
\newcommand{\CH}{\Cat{CH}}
\newcommand{\op}[1]{{#1}{^{\mathsf{op}}}}

\newcommand{\ketbra}[2]{\left|#1\right>\!\left<#2\right|}
\newcommand{\uleq}{\mathbin{\rotatebox[origin=c]{90}{$\leq$}}}

\newcommand{\wn}{\mathop{\mathrm{wn}}}

\newcommand{\bsp}{\scrB}

\NewDocumentCommand{\vmleq}{o}{\mathrel{%
\IfNoValueTF{#1}{\lesssim}{\lesssim_{#1}}}}

\newcommand{\C}{\mathbb{C}}
\newcommand{\N}{\mathbb{N}}
\newcommand{\Z}{\mathbb{Z}}
\newcommand{\R}{\mathbb{R}}
\newcommand{\I}{\mathbb{I}}
\newcommand{\spec}{\mathrm{sp}}
\newcommand{\scrA}{\mathscr{A}}
\newcommand{\scrB}{\mathscr{B}}
\newcommand{\scrC}{\mathscr{C}}
\newcommand{\scrD}{\mathscr{D}}
\newcommand{\scrE}{\mathscr{E}}
\newcommand{\scrG}{\mathscr{G}}
\newcommand{\scrH}{\mathscr{H}}
\newcommand{\scrF}{\mathscr{F}}
\newcommand{\scrP}{\mathscr{P}}
\newcommand{\scrK}{\mathscr{K}}
\newcommand{\scrL}{\mathscr{L}}
\newcommand{\scrR}{\mathscr{R}}
\newcommand{\scrS}{\mathscr{S}}
\newcommand{\scrT}{\mathscr{T}}
\newcommand{\scrX}{\mathscr{X}}
\newcommand{\scrY}{\mathscr{Y}}
\newcommand{\scrZ}{\mathscr{Z}}

\newcommand{\ceil}[1]{\left\lceil#1\right\rceil}
\newcommand{\cceil}[1]{\left\lceil\!\!\left\lceil#1\right\rceil\!\!\right\rceil}

\newcommand{\floor}[1]{\left\lfloor#1\right\rfloor}
\newcommand{\ceill}[1]{\left\lceil#1\right)}
\newcommand{\ceilr}[1]{\left(#1\right\rceil}

\newcommand{\linf}{\ell^\infty}

\newcommand{\limp}{\mathbin{\multimap}}
\newcommand{\Mon}{\mathrm{Mon}}
\newcommand{\CMon}{\mathrm{cMon}}
\newcommand{\nsp}{\mathrm{nsp}}
\newcommand{\qbit}{\mathsf{qubit}}
\newcommand{\bit}{\mathsf{bit}}
\DeclarePairedDelimiter{\sem}{\llbracket}{\rrbracket}

\newcommand{\Real}[1]{#1_{\R}}
\newcommand{\Imag}[1]{#1_{\I}}
\newcommand{\sa}[1]{#1_{\R}}
\newcommand{\pos}[1]{#1_{+}}
\newcommand{\bang}{\mathord{!}}

\DeclareMathOperator{\dom}{dom}
\DeclareMathOperator{\length}{length}

\newcommand\BOX{%
    {\rlap{\kern .5pt\tikz[baseline]{\draw[line width=0.08ex] (0,0) rectangle (0.6ex,0.6ex)}}%
    \phantom{\diamond}}}

\newcommand{\cmpr}[2]{\ensuremath{\{#1|{\kern.2ex}#2\}}}
\newcommand{\id}{\mathrm{id}}

\DeclareMathOperator{\Proj}{Proj}
\DeclareMathOperator{\supp}{supp}
\DeclareMathOperator{\Ker}{Ker}
\DeclareMathOperator{\Ran}{Ran}

\newcommand{\uwlim}{\qopname\relax m{uwlim}}

\makeatletter
\newcommand{\bigperp}{%
  \mathop{\mathpalette\bigp@rp\relax}%
  \displaylimits
}

\newcommand{\bigp@rp}[2]{%
  \vcenter{
    \m@th\hbox{\scalebox{\ifx#1\displaystyle2.1\else1.5\fi}{$#1\perp$}}
  }%
}
\makeatother

\makeatletter
\def\ourrawref#1{%
    \expandafter\expandafter\expandafter
    \@car\csname r@#1\endcsname\@nil
}
\makeatother

\newcommand{\qed}{\hfill\textcolor{darkblue}{\ensuremath{\square}}}

\newcounter{tmptmp}

\newcounter{parsec} 
\newcounter{parsecMajor}
\newcounter{parsecMinor} 

\newcommand\refsforparsec{%
    \setcounter{parsecMajor}{\value{parsec}/10}%
    \setcounter{parsecMinor}{\value{parsec}-10*\value{parsecMajor}}%
    \setcounter{parsec}{\value{parsec}-1}%
    \refstepcounter{parsec}%
}

\NewDocumentEnvironment{parsec}{g o}{%
	\leavevmode\unskip%
	\par\vskip1em\noindent%
    \renewcommand{\theparsec}{\the\value{parsecMajor}\alph{parsecMinor}}%
    \setcounter{parsec}{#1}%
    \refsforparsec{}%
	\setcounter{tmptmp}{2*\value{parsec}}%
	\markboth{\the\value{tmptmp}}{\the\value{tmptmp}}%
	\IfValueT{#2}{\label{#2}}%
    \renewcommand{\thepoint}{\the\value{parsecMajor}\alph{parsecMinor}}%
    \setcounter{point}{\value{point}-1}%
    \refstepcounter{point}%
    \label{parsec-\the\value{parsec}}%
        \marginnote{\makebox[3em][c]{\textParsecNumber{%
            \the\value{parsecMajor}%
            \alph{parsecMinor}%
            }}}%
	\ignorespaces%
}{%
	\leavevmode\unskip%
	\setcounter{tmptmp}{2*\value{parsec}+1}%
	\markboth{\the\value{tmptmp}}{\the\value{tmptmp}}%
	\ignorespaces%
}

\newcounter{point} 
\newcounter{pointMajor}
\newcounter{pointMinor}
\numberwithin{point}{parsec}
\newcounter{pointdepth}  

\newcommand\refsforpoint{%
    \setcounter{pointMajor}{\value{point}/10}%
    \setcounter{pointMinor}{\value{point}-10*\value{pointMajor}}%
    \setcounter{point}{\value{point}-1}%
    \refstepcounter{point}%
}

\NewDocumentEnvironment{point}{g o g}{%
	\leavevmode\unskip%
	\setcounter{pointdepth}{\value{pointdepth}+1}%
	\refstepcounter{point}
    \setcounter{point}{#1}%
    \IfValueT{#2}{%
        \renewcommand{\thepoint}{\the\value{parsec}}%
        \refsforpoint{}%
        \label{#2::parsec}%
        \renewcommand{\thepoint}{\Roman{pointMajor}\alph{pointMinor}}%
        \refsforpoint{}%
        \label{#2::point}%
        \renewcommand{\thepoint}{\theparsec\,\Roman{pointMajor}\alph{pointMinor}}%
        \refsforpoint{}%
        \label{#2}%
    }%
    \renewcommand{\thepoint}{\theparsec\,\Roman{pointMajor}\alph{pointMinor}}%
    \refsforpoint{}%
    \label{parsec-\the\value{parsec}.\the\value{point}}%
	\ifthenelse{\equal{\value{point}}{10}}{}{%
		\ifthenelse{\equal{\value{pointdepth}}{1}}{%
			\par\penalty-100\vskip.6em\noindent%
		}{%
			\ifthenelse{\equal{\value{pointdepth}}{2}}{%
				\par\penalty-50\vskip.2em\noindent%
			}{%
				\par\penalty-25\vskip.1em\noindent%
			}%
		}%
		\marginnote{\makebox[2em][c]{\small%
                    \ifthenelse{\equal{\value{pointdepth}}{1}}{%
                        \textPointNumberI{\Roman{pointMajor}\alph{pointMinor}}%
                    }{%
                        \ifthenelse{\equal{\value{pointdepth}}{2}}{%
                            \textPointNumberII{\Roman{pointMajor}\alph{pointMinor}}%
                        }{%
                            \textPointNumberIII{\Roman{pointMajor}\alph{pointMinor}}%
                        }%
		}}}%
	}%
	\IfValueT{#3}{%
		\ifthenelse{\equal{\value{pointdepth}}{1}}{%
			\textPointHeaderI{#3}%
		}{%
			\ifthenelse{\equal{\value{pointdepth}}{2}}{%
				\textPointHeaderII{#3}%
			}{%
				\textPointHeaderIII{(#3)}%
			}}%
	\ \ }%
\ignorespaces%
}{%
\leavevmode\unskip%
\setcounter{pointdepth}{\value{pointdepth}-1}%
\ignorespaces%
}

\NewDocumentCommand{\sref}{m}{\textSref{%
	\ifthenelse{%
		\equal{\value{parsec}}{\ourrawref{#1::parsec}}%
	}{%
		\ref{#1::point}%
	}{%
        \ref{#1}%
	}%
}}

\usepackage{fancyhdr}
\renewcommand{\chaptermark}[1]{}  
\renewcommand{\sectionmark}[1]{}

\newcounter{firstParsec}
\newcounter{lastParsec}
\newcounter{firstParsecF}
\newcounter{lastParsecF}
\newcounter{firstParsecMajor}
\newcounter{firstParsecMinor}
\newcounter{lastParsecMajor}
\newcounter{lastParsecMinor}

\newcounter{parsecToBeContinued}  
\newcommand\ourfancyfooters{%
    \fancyfoot[CE]{%
        \IfInteger{\rightmark}{%
            \setcounter{tmptmp}{\leftmark+0}%
            \ifthenelse{\equal{\value{tmptmp}}{\value{lastParsecF}}}{%
                \setcounter{firstParsecF}{0}%
            }{%
                \setcounter{firstParsecF}{\rightmark+0}%
            }
        }{%
        }%
    }%
    \fancyfoot[RO]{%
        \textsf{\footnotesize\textcolor{lightgray}{\thepage}}}%
    \fancyfoot[CO]{%
        \ifthenelse{\equal{\value{firstParsecF}}{0}}{%
            \setcounter{firstParsecF}{\rightmark+0}%
        }{%
        }
        \setcounter{firstParsec}{\value{firstParsecF}/2}%
        \setcounter{lastParsecF}{\leftmark+0}%
        \ifthenelse{\equal{\value{lastParsecF}}{0}}{%
        }{%
            \setcounter{lastParsec}{\value{lastParsecF}/2}%
            \textPointNumberI{%
                \ifthenelse{\equal{\value{parsecToBeContinued}}{1}}{..}{}%
                \setcounter{firstParsecMajor}{\value{firstParsec}/10}%
                \setcounter{firstParsecMinor}{\value{firstParsec}-\value{firstParsecMajor}*10}%
                \the\value{firstParsecMajor}\alph{firstParsecMinor}%
                \ifthenelse{\equal{\value{firstParsec}}{\value{lastParsec}}}{%
                }{
                    \setcounter{tmptmp}{\value{firstParsec}+1}%
                    \ifthenelse{\equal{\value{tmptmp}}{\value{lastParsec}}}{, }{--}
                    \setcounter{lastParsecMajor}{\value{lastParsec}/10}%
                    \setcounter{lastParsecMinor}{\value{lastParsec}-\value{lastParsecMajor}*10}%
                    \the\value{lastParsecMajor}\alph{lastParsecMinor}%
                }%
                \setcounter{tmptmp}{\value{lastParsec}*2}%
                \ifthenelse{\equal{\value{tmptmp}}{\value{lastParsecF}}}{%
                    ..\setcounter{parsecToBeContinued}{1}%
                }{%
                    \setcounter{parsecToBeContinued}{0}%
                }%
            }%
        }%
    }%
}

\ourfancyfooters
\fancypagestyle{plain}{\ourfancyfooters}

\makeatletter
\def\@wrindex#1{%
   \protected@write\@indexfile{}%
   {\string\indexentry{#1|parsechyperlink{\the\value{parsec}.\the\value{point}}}{\the\value{parsec}.\the\value{point}}}%
 \endgroup
 \@esphack}%
\makeatother

\let\oldchapter\chapter%
\renewcommand\chapter[1]{%
    \oldchapter{#1}%
    \setcounter{tmptmp}{(\value{parsec}/10)*10+10}
    \addtocontents{parsectoc}{\protect\contentsline{chapter}{\numberline{\thechapter} #1}{\sref{parsec-\the\value{tmptmp}}}{chapter.\thechapter}}%
}
\let\oldsection\section%
\renewcommand\section[1]{%
    \oldsection{#1}%
    \setcounter{tmptmp}{(\value{parsec}/10)*10+10}
    \addtocontents{parsectoc}{\protect\contentsline{section}{\numberline{\thesection} #1}{\sref{parsec-\the\value{tmptmp}}}{section.\thesection}}%
}
\let\oldsubsection\subsection%
\renewcommand\subsection[1]{%
    \oldsubsection{#1}%
    \setcounter{tmptmp}{(\value{parsec}/10)*10+10}
    \addtocontents{parsectoc}{\protect\contentsline{subsection}{\numberline{\thesubsection} #1}{\sref{parsec-\the\value{tmptmp}}}{subsection.\thesubsection}}%
}
\let\oldsubsubsection\subsubsection%
\renewcommand\subsubsection[1]{%
    \oldsubsubsection{#1}%
    \setcounter{tmptmp}{(\value{parsec}/10)*10+10}
    \addtocontents{parsectoc}{\protect\contentsline{subsubsection}{\numberline{\thesubsubsection} #1}{\sref{parsec-\the\value{tmptmp}}}{subsubsection.\thesubsubsection}}%
}

\newcommand\backmattertitle[1]{{%
    \sffamily\color{darkblue}\Huge\bfseries #1\vspace{1em}}}



\usepackage{subfiles}
\makeindex

\begin{document}

\newcommand{\titelpagina}[8]{%
\vspace*{4em}
\begin{center}
{\huge\sffamily\color{darkblue}The Category of Von Neumann Algebras}
\end{center}
\vspace{2em}
\begin{center}
{\par\noindent\sffamily\color{darkblue}\large\textbf{#1}}
\end{center}
\begin{center}
#2
\end{center}
\vspace{1em}
\begin{center}
#3
\end{center}
\begin{center}
#4
\end{center}
\begin{center}
#5
\end{center}
\begin{center}
#6
\end{center}
\vspace{10em}
\begin{center}%
#7
\end{center}%
\begin{center}%
\large\color{darkblue}Abraham Anton \textsc{Westerbaan}
\end{center}%
\begin{center}%
#8
\end{center}
\newpage
}

\newcommand{\achterkanttitelpagina}[3]{
\vspace*{8em}
\begin{center}
\textbf{\large\sffamily\color{darkblue}#1:}
\end{center}
\begin{center}
Prof.~dr.~B.P.F.~\textsc{Jacobs}
\end{center}
\vspace{5em}
\begin{center}
\textbf{\large\sffamily\color{darkblue}#2:}
\end{center}
\begin{center}
Prof.~dr.~J.D.M.~\textsc{Maassen}\\
\vspace{1em}
    Prof.~dr.~P.~\textsc{Panangaden} \\ {\footnotesize(McGill University, Canada)} \\
\vspace{1em}
    Prof.~dr.~P.~\textsc{Selinger} \\ {\footnotesize(Dalhousie University, Canada)} \\
\vspace{1em}
    Dr.~C.J.M.~\textsc{Heunen} \\ {\footnotesize(University of Edinburgh, #3)} \\
\vspace{1em}
Dr.~A.R.~\textsc{Kissinger}
\end{center}
\newpage
}

\newcommand{\whiteout}[1]{{\color{white}#1}}

\titelpagina{\whiteout{P}$\,$}{$\,$\whiteout{ter}\\
\whiteout{aan}$\,$\\
\whiteout{op}$\,$\\
\whiteout{volgens}$\,$\\
\whiteout{in}$\,$}{\whiteout{op}$\,$}{\whiteout{dinsdag}$\,$}{\whiteout{om}$\,$}{\whiteout{10.30}$\,$}{\whiteout{door}}{$\,$\whiteout{geboren}$\,$\\
\whiteout{te}$\,$}

\colofon{A.A.~Westerbaan}{978-94-6332-484-7}{1804.02203}{2066/201611}

\titelpagina{Proefschrift}{ter verkrijging van de graad van doctor\\
aan de Radboud Universiteit Nijmegen\\
op gezag van de rector magnificus prof. dr. J.H.J.M. \textsc{van Krieken},\\
volgens besluit van het college van decanen \\
in het openbaar te verdedigen}{op}{\textbf{dinsdag 14 mei 2019}}{om}{\textbf{10.30 uur precies}}{door}{geboren op 30 augustus 1988\\
te Nijmegen}

\achterkanttitelpagina{Promotor}{Manuscriptcommissie}{Verenigd Koninkrijk}

\titelpagina{Doctoral Thesis}{to obtain the degree of doctor \\
from  Radboud University Nijmegen \\
on the authority of the Rector Magnificus prof. dr. J.H.J.M. \textsc{van Krieken},\\
according to the decision of the Council of Deans \\
to be defended in public}{on}{\textbf{Tuesday, May 14, 2019}}{at}{\textbf{10.30 hours}}{by}{born on August 30, 1988\\
in Nijmegen (the Netherlands)}

\achterkanttitelpagina{Supervisor}{Doctoral Thesis Committee}{United Kingdom}

\vspace*{.1em}

\newpage
\makeatletter\@starttoc{parsectoc}\makeatother
\chapter{Introduction}


\begin{parsec}{10}
\begin{point}{10}
What does this Ph.D.~thesis offer?
Proof, perhaps,
to the doctoral thesis committee
of passable academic work;
an advertisement, as it may be,
of my school's perspective
to colleagues;
a display, even,
of intellectual achievement
to friends and family.
But I believe such narrow and selfish goals \emph{alone}
barely serve to keep a writer's spirits 
energised---and are definitely detrimental to that of the readers.
That is why I have foolhardily
challenged
myself
not just 
to drily list contributions,
but to write this thesis 
as the introduction,
that I would have liked to read
when I started
research for this thesis
back in May~2014.

The topic is von Neumann algebras,
the category they form,
and how they may be used
to model aspects of quantum computation.
Let us just say for now that a von Neumann algebra
is a special type of complex vector
space endowed with
a multiplication operation among some other additional structure.
An important example is the complex vector space~$M_2$
of~$2\times 2$ complex matrices,
because it models (the predicates on) a qubit;
but all~$N\times N$-complex matrices form a von Neumann algebra~$M_N$ as well.
Using von Neumann algebras
(and their little cousins, $C^*$-algebras) 
to describe quantum data types 
seems to be quite a recent idea
(see e.g.~\cite{jacobs2013block,rennela2015operator,furber2013kleisli}, 
	and~\cite{cho2016semantics} for an overview)
and has two distinct features.
Firstly, classical data types
are neatly incorporated:
$\C^2 \equiv \C\oplus \C$
models a bit,
and the direct sum $M_2\oplus M_3$
models the union type of a qubit and a qutrit.
Secondly,
von Neumann algebras
allow for infinite data types as well
	such as~$\scrB(\ell^2(\Z))$,
which represents a ``quantum integer.''\footnote{Though
	other methods of modelling infinite dimensional
	quantum computing have been proposed as well
	e.g.~using non-standard analysis\cite{Gogioso2017},
	pre-sheaves\cite{malherbe2013categorical},
	the geometry of interaction\cite{hasuo2017semantics},
	and quantitative semantics\cite{pagani2014applying}.}
It should be said that this last feature
is both a boon and a bane:
it brings with it all the inherent
intricacies of dealing with infinite dimensions;
and it is no wonder that
most authors choose 
to restrict themselves
to finite dimensions,
especially since
this seems to be enough to describe quantum algorithms,
see e.g.~\cite{nielsen2002quantum}.
\end{point}
\begin{point}{20}%
In this thesis, however,
we do face infinite dimensions,
because the two main results demand it:
\begin{enumerate}
\item
For the first result,
that von Neumann algebras
form  a model of Selinger and Valiron's quantum lambda calculus,
		as Cho and I explained in~\cite{model}
and for which I'll provide the foundation here,
we need to interpret function types,
some of which are essentially infinite dimensional.
\item
The second result,
an axiomatisation
of the map $a\mapsto \sqrt{p}a\sqrt{p}\colon \scrA\to\scrA$
representing measurement
of an element $p\in[0,1]_\scrA$
of a von Neumann algebra~$\scrA$
was tailored by B.E.~Westerbaan (my twin brother) and myself to work for
both finite and infinite dimensional~$\scrA$.
\end{enumerate}
These results
are part of a line of research that
tries to find patterns
in the category of von Neumann algebras,
that may also be cut from other categories
modelling computation---ideally in order to arrive at categorical axioms
for (probabilistic) computation in general.
When I joined the fray
the notion of \emph{effectus}\cite{newdirections} had already
been established by Jacobs,
and the two results
above offer potential additional axioms.
The work in this area has largely been a collaborative effort,
primarily between Jacobs, Cho, my twin brother, and myself,
and many of their insights have ended up in this thesis.

Of this I'd say no more
than that my work appears conversely, and proportionally,
in their writings too, except that the close cooperation
with my brother begs further explanation.
Our efforts on certain topics have been like interleaving 
of the pages of two phone books:
separating them  would be nigh impossible,
especially the work on the axiomatisation
of~$a\mapsto \sqrt{p}a\sqrt{p}$ and Paschke dilations.
So that's why we decided to write our theses
as two volumes of the same work;
preliminaries on von Neumann algebras,
and the axiomatisation of $a\mapsto \sqrt{p}a\sqrt{p}$
appear in this thesis,
while the work on dilations,
and effectus theory appear in my brother's thesis,
\cite{bas}.
\end{point}
\begin{point}{30}
The two results mentioned above only
make up 
about a third of this thesis;
the rest of it is devoted to 
the introduction to the theory of von Neumann algebras
needed to understand these results.
My aim is that anyone 
with, say, a bachelor's degree in mathematics
(more specifically, basic knowledge of linear algebra,
analysis\cite{rudin1964principles}, 
topology\cite{willard}
and set theory\cite{devlin2012joy})
should at least be able to follow the lines of reasoning
with only minimal recourse to external sources.
But I hope that they will gain some deeper understanding
of the material as well.
To this end, and because I wanted to gain some of this insight
for myself too,
I've not just mixed
and matched
results from the literature,
but I tailored a thorough treatise
of everything that's needed,
including proofs (except~\sref{intersection-tensor-proof}).
Whenever possible,
I've taken shortcuts
(e.g.~avoiding for example
the theory of Banach algebras
and locally convex spaces entirely)
to  prevent the mental tax
the added concepts
(and pages) would have brought.
For the same reasons
I've refrained from putting
everything in its proper abstract (and categorical\cite{maclane}) context
trusting that it'll shine through of its own accord.
I've however not been able to restrain
myself in making perhaps frivolous variations on the existing
theory whenever not strictly necessary,
taking for example Kadison's characterisation\cite{kadison1956}
of von Neumann algebras
as my definition,
and developing the elementary theory for it;
in my defence I'll just say this adds to
the original element that is expected of a thesis.
\end{point}
\begin{point}{40}{Advertisements}%
Due to space--time constraints
this thesis is based only on a selection
\cite{model,cho2015quotient,cho2016duplicable,qpakm,westerbaan2016universal}
of the works
I produced under supervision of Jacobs,
and while~\cite{wwpaschke,effintro,statesofconvexsets}
are incorporated in my brother's thesis,
this means~\cite{jacobs2015effect,jacobs2017distances} 
are unfortunately left out.
If you like this thesis,
then
you might also want to take a look
at these~\cite{rennela2017infinite,
rennela2015complete,
furber2013kleisli,
kornell2012,
heunen2015domains,
Maassen2010} recent works on von Neumann algebras,
and $C^*$-algebras.
If you're curious
about effectus theory
and related matters,
please have a look at~\cite{jacobs2017quantum,
cho2017disintegration,
jacobs2016hyper,
jacobs2017channel,
jacobs2017formal,
cho2017efprob,
jacobs2017probability,
jacobs2017recipe,
jacobs2016effectuses,
jacobs2016affine,
jacobs2016relating,
effintro,
statesofconvexsets,
cho2015quotient,
jacobs2017distances,
jacobs2015effect,
jacobs2016expectation,
jacobs2016predicate,
newdirections}.
But if you'd like more pictures instead,
I'd suggest~\cite{coecke2017picturing}.
\end{point}
\begin{point}{50}[on-writing-style]{Writing style}
I've replaced page numbers by
paragraph numbers
such as~\sref{on-writing-style}
for this paragraph.
The numbers after~\sref{final-bram} refer to paragraphs
	in my twin brother's thesis\cite{bas}.
\Define{Definitions}\index{definitions} are set like that
(i.e.~in blue),
and can be found in the index.
Proofs of certain facts
that are easily obtained on the back of an envelope,
and would clutter this manuscript,
have been left out.
Instead these facts have been phrased as exercises
as a challenge to the reader.
\end{point}
\begin{point}{51}
    The symbol ``\Define{$:=$}''\index{(((defequal@$:=$, is defined to be}
should be interpreted as ``is defined to be'',
    while ``\Define{$\equiv$}''\index{(((equiv@$\equiv$, being of the form}
should be read as ``being of the form''.
Sometimes ``$\equiv$'' is
used to define something on its right-hand side,
as in ``let $A\equiv\left(
\begin{smallmatrix} a&b\\b^* &c\end{smallmatrix}
\right)$ be a self-adjoint matrix.''
Other times~``$\equiv$''
indicates a simple rewrite step, as in
``since~$a=2$, we have $a+2=2+2\equiv 4$,''
where it's not suggested $a=2$ implies $2+2=4$.
\end{point}
\begin{point}{60}{Acknowledgements}
The work in this thesis specifically
has benefited greatly from
discussions
with John van de Wetering,
Robert Furber,
Kenta Cho,
and Bas Westerbaan,
but I've also had the pleasure
of discussing a variety
of other topics 
with 
Aleks Kissinger,
Andrew Polonsky,
Bert Lindenhovius,
Frank Roumen, 
Hans Maassen,
Henk Barendregt,
Joshua Moerman,
Martti Karvonen,
Robbert Krebbers,
Robin Adams, 
Robin Kaarsgaard,
Sam Staton, 
Sander Uijlen,
Sebastiaan Joosten,
and many others.
I'm especially honoured to have been received
in Edinburgh by Chris Heunen 
and in Oberwolfach
by Jianchao Wu.
I'm very grateful to
Arnoud van Rooij,
Bas Westerbaan
and
John van de Wetering
for proofreading large parts of
this manuscript,
without whose efforts
even more shameful errors would have remained.
I'm very grateful too for the manuscript committee's members' 
various suggestions and comments, and hope the improvements I made to this text
do them justice.
I should of course not forget to mention
the contribution
of friends (both close and distant),
family,
and colleagues---too numerous to name---of keeping me sane
these past years.

This is the second dissertation topic
I've worked on;
my first attempt
under different supervision
was unfortunately cut short after $1\sfrac{1}{2}$ years.
When Bart Jacobs graciously offered
me a second chance,
I initially had my reservations,
but accepted on account of the challenging topic.
Little did I know 
that behind the ambition and suit
	one finds a man
of singular moral fibre,
embodying
what was said
	about von Neumann himself:
	``[he] had to understand and accept much that most 
of us do not want to accept and do not even wish to understand.''%
\footnote{An excerpt from Eugene P.~Wigner's writings,
see page~130 of~\cite{wigner2013collected}.}
\end{point}
\begin{point}{70}{Funding} was received from the 
European Research Council under grant agreement \textnumero~320571.
\end{point}
\end{parsec}

\chapter{$\text{C}^*$-algebras}
\begin{parsec}{20}
\begin{point}{10}
We redevelop the essentials of the theory of (unital) $C^*$-algebras
in this chapter.
Since we are ultimately interested
in von Neumann algebras
(a special type of $C^*$-algebras)
we will evade
delicate
topics such as tensor products (of $C^*$-algebras), 
quotients, approximate identities,
and $C^*$-algebras without a unit.
The zenith of this chapter
is \emph{Gelfand's representation theorem} (see~\sref{gelfand}),
the fact that every commutative (unital) $C^*$-algebra
is isomorphic
to the $C^*$-algebra
$C(X)$ of continuous functions on some compact Hausdorff space~$X$
--- it yields a duality
	between the category~$\CH$
	of compact Hausdorff spaces (and continuous maps)
and the category~$\cCstar{miu}$ of commutative $C^*$-algebras (and
unital $*$-homomorphisms,
the appropriate structure preserving maps), see~\sref{gelfand-equivalence}.

As the road to Gelfand's representation theorem 
is a bit winding ---
involving intricate relations between technical concepts --- 
we have put emphasis on the invertible and  positive elements
so that the important
theorems about them may serve as landmarks along the way:
\begin{enumerate}
\item
first we show that the norm
on a $C^*$-algebra
is determined by the invertible elements
(via the \emph{spectral radius}), see~\sref{norm-spectrum};

\item
then we construct a \emph{square root} of a positive element in~\sref{sqrt};

\item
and finally we
show that an element of a commutative $C^*$-algebra
is not invertible iff it is mapped to~$0$
by some multiplicative state, see~\sref{inv-mult-state}.
\end{enumerate}
At every step along the way
the positive and invertible elements 
(and the norm, multiplicative states, multiplication
and other structure on a $C^*$-algebra)
are bound more tightly together
until Gelfand's representation theorem emerges.

To make this chapter
more accessible
we have removed
much material
from the ordinary development
of $C^*$-algebras
such as the more general theory of Banach algebras
(and its pathology).
This forces us
to take a slightly different path than is usual in the literature 
(see e.g.~\sref{gelfand-mazur-predicament}).

After Gelfand's representation theorem
we deal with two smaller topics:
that a $C^*$-algebra
may be represented as a concrete $C^*$-algebra
of bounded operators on a Hilbert space (see~\sref{gns}),
and that the $N\times N$-matrices with entries drawn from a 
$C^*$-algebra~$\scrA$
form a $C^*$-algebra~$M_N(\scrA)$
(see~\sref{cstar-matrices}).
We end
with an overture to von Neumann algebras---the
topic of the next chapter.
\end{point}%
\end{parsec}
\section{Definition and Examples}
\begin{parsec}{30}
\begin{point}{10}{Definition}
A \Define{$C^*$-algebra}\index{Cstar-algebra@$C^*$-algebra}
is a complex vector space~$\scrA$
endowed with
\begin{enumerate}
\item
a binary operation,
called \Define{multiplication}
(and denoted as such),
which is associative, and linear in both coordinates;
\item
an element~$1$, called 
\Define{unit}\index{unit!of a {$C^*$-algebra}},
such that $1\cdot a = a = a\cdot 1$
for all~$a\in \scrA$;
\item
a unary operation $\Define{(\,\cdot\,)^*}$,
called \Define{involution} %
\index{$(\,\cdot\,)^*$ !involution on a $C^*$-algebra}%
\index{involution!on a $C^*$-algebra}
such that $(a^*)^*=a$,
$(ab)^*=b^*a^*$,
$(\lambda a)^* = \bar\lambda a^*$,
and $(a+b)^* = a^*+b^*$
for all~$a,b\in\scrA$ and~$\lambda\in \C$;
\item
a complete \Define{norm}%
\index{$"\"|\,\cdot\,"\"|$, norm!on a $C^*$-algebra}
$\Define{\|\,\cdot\,\|}$
such that
$\|ab\|\leq\|a\|\|b\|$
for all~$a,b\in\scrA$,
and 
\begin{equation*}
\label{eq:Cstar-identity}
\|a^*a\|\ =\ \|a\|^2
\end{equation*}
holds; this equality is called the \Define{$C^*$-identity}.%
\index{Cstar-identity@$C^*$-identity}
\end{enumerate}
The $C^*$-algebra $\scrA$ is called \Define{commutative}%
\index{Cstar-algebra@$C^*$-algebra!commutative}
if $ab=ba$ for all~$a,b\in\scrA$.
\begin{point}{20}{Warning}%
In the literature it is usually not
required that a $C^*$-algebra
possess a unit; but when it does it is called
a \Define{unital $C^*$-algebra}.%
\index{unital!$C^*$-algebra}
\end{point}
\end{point}
\begin{point}{30}{Example}%
The vector space~$\C$ of \Define{complex numbers}%
\index{C@$\C$, the complex numbers!as a $C^*$-algebra}
forms a commutative  $C^*$-algebra
in which
multiplication and~$1$
have their usual meaning.
Involution is given by conjugation ($z^*=\bar{z}$),
and norm by modulus ($\|z\|=|z|$).
\end{point}
\begin{point}{40}{Example}%
A \Define{$C^*$-subalgebra}%
\index{Cstar-subalgebra@$C^*$-subalgebra}
of a $C^*$-algebra~$\scrA$
is a subset~$\scrB$ of~$\scrA$,
which is a linear subspace of~$\scrA$,
contains the unit, $1$, is closed under multiplication
and involution, 
and is closed with respect to the norm of~$\scrA$;
such a $C^*$-subalgebra of~$\scrA$
is itself a $C^*$-algebra
when endowed with the operations and norm
of~$\scrA$.
\end{point}
\begin{point}{50}[cstar-product]{Example}%
One can form products (in the categorical sense,
see~\sref{cstar-product-2}) of $C^*$-algebras as follows.
Let~$\scrA_i$ be a $C^*$-algebra
for every element~$i$ of some index set~$I$.
The \Define{direct sum}%
\index{direct sum!of $C^*$-algebras}
\index{$\bigoplus$, direct sum!$\bigoplus_i \scrA_i$, of $C^*$-algebras}
of the family $(\scrA_i)_i$
is the $C^*$-algebra
denoted by \Define{$\bigoplus_{i\in I}\scrA_i$} on the set
\begin{equation*}
\textstyle
\bigl\{\ 
a\in \prod_{i\in I}\scrA_i\colon\  \sup_{i \in I} \|a(i)\|< \infty \ 
\bigr\}
\end{equation*}
whose operations are defined coordinatewise,
and whose norm is a \Define{supremum norm}%
\index{supremum norm}%
\index{$"\"|\,\cdot\,"\"|$, norm!supremum $\sim$}
given by $\|a\|=\sup_{i}\|a(i)\|$.
If each~$\scrA_i$ is commutative,
then~$\bigoplus_{i\in I}\scrA_i$
is commutative.

In particular,
taking~$\scrA_i\equiv \C$,
we see that
the vector space~\Define{$\ell^\infty(X)$}%
\index{linfty@$\ell^\infty(X)$}%
\index{linfty@$\ell^\infty(X)$!as a $C^*$-algebra}
of bounded complex-valued functions
on a set~$X$ forms a commutative $C^*$-algebra
with pointwise operations and supremum norm.
\end{point}
\begin{point}{60}{Example}%
The \Define{bounded continuous functions
on a topological space}~$X$
form a commutative $C^*$-subalgebra~\Define{$BC(X)$}%
\index{$BC(X)$! as a $C^*$-algebra}
of~$\ell^\infty(X)$ (see above).
In particular,
since a continuous function on a compact Hausdorff space is 
automatically bounded,
we see that the \Define{continuous functions
on a compact Hausdorff space} $X$
form a commutative $C^*$-algebra~\Define{$C(X)$}%
\index{$C(X)$}%
\index{$C(X)$!as a $C^*$-algebra}
with pointwise operations and sup-norm.
We'll see that every commutative $C^*$-algebra
is isomorphic to a~$C(X)$
in~\sref{gelfand}.
\end{point}
\begin{point}{70}[cstar-matrices-example]{Example}%
An example of a non-commutative
$C^*$-algebra
is
the vector space~\Define{$M_n$}%
\index{$M_n$, the $n\times n$-matrices!as a $C^*$-algebra}
of \Define{$n\times n$-matrices} ($n>1$) over~$\C$
with the usual (matrix) multiplication,
the identity matrix as unit,
and conjugate transpose
as involution
(so~$(A^*)_{ij} = \overline{A_{ji}}$).
The norm~$\|A\|$ of a matrix~$A$ in~$M_n$
is less obvious,
being
the \emph{operator norm}
(cf.~\sref{bounded-linear-maps})
of the associated linear map~$v\mapsto Av,\ \C^n\to\C^n$,
that is,
$\|A\|$ is
the least number~$r\geq 0$
with $\|Av\|_2\leq r\|v\|_2$
for all~$v\in \C^n$
(where $\|w\|_2=(\sum_i \left|w_i\right|^2)^{\nicefrac{1}{2}}$
denotes the $2$-norm
of~$w\in \C^n$).

It is not entirely obvious that~$\|A^*A\|=\|A\|^2$
holds
and that $M_n$ is complete.
We will prove these facts in the more general setting
of bounded operators between Hilbert spaces, 
see~\sref{adjoinables-cstar-algebra}.
Suffice it to say, $\C^n$ is a Hilbert space
with~$\left<v,w\right>=\sum_i \overline{v}_iw_i$
as inner product,
each matrix gives a (bounded) linear map $v\mapsto Av,\C^n\to \C^n$,
and the conjugate transpose $A^*$ is \emph{adjoint} to~$A$
in the sense that $\left<v,Aw\right> = \left<A^*v,w\right>$
for all~$v,w\in\C^n$.
\end{point}
\begin{point}{80}{Remark}%
\index{Cstar-algebra@$C^*$-algebra!finite dimensional}
Combining~\sref{cstar-product}
and~\sref{cstar-matrices-example}
we see that
$\bigoplus_k M_{n_k}$
is a finite-dimensional
$C^*$-algebra
for any tuple $n_1,\dotsc,n_K$
of natural numbers.
In fact,
any finite-dimensional $C^*$-algebra
is of this form
as we'll see in~\sref{fdcstar}.\footnote{Although
clearly related to the 
Wedderburn--Artin theorem,
see e.g.~\cite{nicholson1993},
this description of finite-dimensional
$C^*$-algebras
does not seem to be an immediate consequence of it.}
\end{point}
\end{parsec}
\subsection{Operators}
\begin{parsec}{40}[hilb]%
\begin{point}{10}[example-hilb]{Example}%
Let us now turn to perhaps the most important
and difficult example:
we'll show that the vector space~\Define{$\scrB(\scrH)$}%
\index{BH@$\scrB(\scrH)$!as a $C^*$-algebra}
of \Define{bounded operators
on a Hilbert space}~$\scrH$ forms a $C^*$-algebra
when endowed with the operator 
norm.
Multiplication is given by composition,
involution by taking the \emph{adjoint} (see~\sref{hilb-def}),
and unit by the identity operator.
A \Define{concrete $C^*$-algebra}%
\index{Cstar-algebra@$C^*$-algebra!concrete}
or
a \Define{$C^*$-algebra of bounded operators}%
\index{Cstar-algebra@$C^*$-algebra!of bounded operators}
refers to a $C^*$-subalgebra of~$\scrB(\scrH)$.
We will eventually see that every $C^*$-algebra is isomorphic to a $C^*$-algebra
of bounded operators in~\sref{gelfand-naimark}.
\end{point}
\begin{point}{20}[bounded-linear-maps]{Definition}%
Let~$\scrX$ and~$\scrY$ be normed
vector spaces.
We say that~$r\in [0,\infty)$
is a \Define{bound}%
\index{bound!for a linear map}
for a linear map (=\Define{operator}%
\index{operator}) 
$T\colon \scrX\to\scrY$
when  $\|Tx\|\leq r\|x\|$ for all~$x\in \scrX$,
and we say that~$T$ is \Define{bounded}%
\index{operator!bounded}
when there is such a bound.
In that case~$T$ has a least bound,
which is called the \Define{operator norm}%
\index{operator norm}%
\index{$"\"|\,\cdot\,"\"|$, norm!of an operator} of~$T$,
and is denoted by~$\Define{\|T\|}$.
The vector space of bounded operators
from~$\scrX$ to~$\scrY$
is denoted by~$\Define{\scrB(\scrX,\scrY)}$,%
\index{BXY@$\scrB(\scrX,\scrY)$}
and the vector space of bounded operators
from~$\scrX$ to itself is denoted by~$\scrB(\scrX)$.%
\index{BX@$\scrB(\scrX)$}
\end{point}
\begin{point}{30}[bounded-operators-basic]{Exercise}%
Let~$\scrX$, $\scrY$ and~$\scrZ$ be normed complex vector spaces.
\begin{enumerate}
\item
Show that the operator norm on~$\scrB(\scrX,\scrY)$
is, indeed, a norm.
\item
Let~$T\colon \scrX\to \scrY$ and~$S\colon \scrY\to\scrZ$
be bounded operators.
Show that $ST$ is bounded by~$\|S\|\|T\|$,
so that~$\|ST\|\leq\|S\|\|T\|$.
\item
Show that the identity operator $\id\colon \scrX\to \scrX$
is bounded by~$1$.%
\end{enumerate}%
\spacingfix{}
\end{point}%
\begin{point}{40}[operator-norm-ball]{Exercise}%
Let $T\colon \scrX\to\scrY$
be a bounded operator between normed vector spaces,
	and let~$r\in[0,\infty)$.
Show that 
	\begin{equation*}
		\textstyle
		r\|T\|\ =\ \sup_{x\in (\scrX)_r} \|Tx\|,
	\end{equation*}
where $\Define{(\scrX)_r}=\{x\in \scrX\colon \|x\|\leq r\}$.%
\index{*ballr@$(\scrX)_r$, $r$-ball}%
\index{*ball@$(\scrX)_1$, unit ball}
	(The set~$(\scrX)_1$
	is called the \Define{unit ball} of~$\scrX$.)%
\index{unit ball}
\end{point}
\begin{point}{50}[operator-norm-complete]{Lemma}%
The operator norm on~$\scrB(\scrX,\scrY)$ is complete
when~$\scrY$ is a complete normed vector space.
\begin{point}{60}{Proof}%
Let~$(T_n)_n$ be a Cauchy sequence in~$\scrB(\scrX,\scrY)$.
We must show that~$(T_n)_n$ converges to some
bounded operator $T\colon \scrX\to\scrY$.
Let~$x\in \scrX$ be given.
Since 
\begin{equation*}
\|\,T_nx - T_mx\,\|\ =\ \|\,(T_n-T_m)\,x\,\|\ \leq\  \|T_n-T_m\|\,\|x\|
\end{equation*}
and~$\|T_n-T_m\|\to 0$ as~$n,m\to \infty$ 
(because~$(T_k)_k$ is Cauchy),
we see that $\|\,T_nx-T_mx\,\|\to 0$ as $n,m\to \infty$,
and so $(T_nx)_n$ is a Cauchy sequence in~$\scrY$.
Since~$\scrY$ is complete,
 $(T_nx)_n$ converges,
and  we may define $Tx:=\lim_n T_nx$,
giving a map $T\colon \scrX\to \scrY$,
which is easily seen to be linear
(by continuity of addition and scalar multiplication).

It remains to be shown that~$T$ is bounded,
and that~$(T_n)_n$ converges to~$T$ with respect to the operator norm.
Let~$\varepsilon>0$ be given, and pick~$N$ such that
$\|T_n-T_m\|\leq \frac{1}{2}\varepsilon$ for all~$n,m\geq N$.
Then for every~$x\in \scrX$
we can find~$M\geq N$ with 
$\|T x - T_m x\|\leq \frac{1}{2}\varepsilon\|x\|$ for all $m\geq M$,
and so,
for $n\geq N$, $m\geq M$,
\begin{equation*}
\|(T - T_n) x\| \ \leq\ \|T x - T_mx\|\,+\,\|T_m x - T_n x\|
\ \leq\  \varepsilon\|x\|
\end{equation*}
giving that~$T-T_n$ is bounded
and $\|T-T_n\|\leq \varepsilon$ for all~$n\geq N$.
Whence~$T$ is bounded too,
and $(T_n)_n$ converges to~$T$.\qed
\end{point}
\end{point}
\begin{point}{70}[bounded-operators-banach-algebra]%
From~\sref{bounded-operators-basic}
and~\sref{operator-norm-complete}
it is clear that the complex vector space
of bounded operators~$\scrB(\scrX)$
on a complete normed vector space~$\scrX$
with composition as multiplication
and the identity operator as unit
satisfies all the requirements
to be a $C^*$-algebra that do not involve the involution, $(\,\cdot\,)^*$
(that is, $\scrB(\scrX)$ is a \Define{Banach algebra}).
To get an involution,
we need the additional structure
provided by a Hilbert space as follows.
\end{point}
\begin{point}{80}[hilb-def]{Definition}%
An \Define{inner product}%
\index{inner product!*C-valued@$\C$-valued}
on a complex vector space~$V$ 
is a map $\left<\,\cdot\,,\,\cdot\,\right>\colon V\times V\to \C$%
\index{$\left<\,\cdot\,,\,\cdot\,\right>$, inner product!$\C$-valued}
such that,
for all~$x,y\in V$,
$\left<x,\,\cdot\,\right>\colon V\to V$ is linear;
$\left<x,x\right>\geq 0$;
and
$\left<x,y\right>=\overline{\left<y,x\right>}$.
We say that the inner product is \Define{definite}%
\index{inner product!*C-valued@$\C$-valued!definite}
when~$\left<x,x\right>=0\implies x=0$ for~$x\in V$.
A \Define{pre-Hilbert space}~$\scrH$%
\index{pre-Hilbert space}
is a complex vector space endowed with a definite inner product.
We'll shortly see that every such~$\scrH$
carries a norm
given by
 $\|x\|:= \left<x,x\right>^{\nicefrac{1}{2}}$;
if~$\scrH$ is complete with respect to this norm,
we say that~$\scrH$ is a \Define{Hilbert space}.%
\index{Hilbert space}

Let~$\scrH$ and~$\scrK$ be pre-Hilbert spaces.
We say that an operator~$T\colon \scrH\to \scrK$
is \Define{adjoint}%
\index{adjoint!of an operator}
to an operator
$S\colon \scrK\to \scrH$ 
when
\begin{equation*}
\left<Tx,y\right> \ = \ \left<x,Sy\right>
\qquad\text{for all $x\in \scrH$ and $y\in \scrK$.}
\end{equation*}
In that case, we call~$T$ \Define{adjointable}.%
\index{adjointable!operator}%
\index{operator!adjointable}
We'll see (in~\sref{uniqueness-adjoint})
that such adjointable~$T$ is adjoint to exactly one~$S$,
which we denote by~\Define{$T^*$}.%
\index{$(\,\cdot\,)^*$!adjoint of an operator}
\end{point}
\begin{point}{90}[hilb-basic-examples]{Example}%
We endow $\C^N$%
\index{C@$\C$, the complex numbers!as a Hilbert space}
(where~$N$ is a natural number)
with the inner product
given by
$\left<x,y\right>=\sum_i \overline{x}_iy_i$,
making it a Hilbert space.

The space~$\Define{c_{00}}$%
\index{c00@$c_{00}$!as a pre-Hilbert space}
of sequences $x_1,x_2,\dotsc$
for which~$x_n$ is non-zero
for finitely many~$n$'s
is an example of a 
pre-Hilbert
which is not complete
when  endowed with $\left<x,y\right>=\sum_{n=0}^\infty \overline{x}_ny_n$
as inner product.

For an example
of an infinite-dimensional Hilbert
space,
we'll have to wait until~\sref{hilb-sum}
where 
we'll show
that the sequences $x_1,x_2,\dotsc$
with $\sum_n \left|x_n\right|^2<\infty$
form a Hilbert space~$\Define{\ell^2}$%
\index{l2@$\ell^2$!as a Hilbert space}
with $\left<x,y\right>=\sum_{n=0}^\infty \overline{x}_ny_n$
as its inner product,
because at this point it is not
even clear that this sum converges.
\end{point}
\begin{point}{100}[uniqueness-adjoint]{Exercise}%
Let~$x$ and~$x'$ be elements of a pre-Hilbert space~$\scrH$
with $\left<y,x\right>=\left<y,x'\right>$
for all~$y\in\scrH$.
Show that~$x=x'$ (by taking $y=x-x'$).
Conclude that every operator between pre-Hilbert spaces
has at most one adjoint.
\begin{point}{110}{Remark}%
Note that we did not require that
an adjointable operator $T\colon \scrH\to\scrK$
between pre-Hilbert spaces be bounded,
and in fact, it might not be.
Take for example
the operator
$T\colon c_{00}\to c_{00}$
given by~$(T x)_n = nx_n$,
which is adjoint to itself,
and not bounded.
On the other hand,
if either~$\scrH$ or~$\scrK$ is complete,
then both~$T$ and~$T^*$ are automatically bounded
as we'll see in~\sref{hellinger-toeplitz}.
\end{point}
\end{point}
\begin{point}{120}{Exercise}%
Let~$S$ and~$T$ be adjointable operators on a pre-Hilbert space.
\begin{enumerate}
\item
Show that~$T^*$ is adjoint to~$T$ (and so $T^{**}=T$).
\item
Show that~$(T+S)^*=T^*+S^*$
and $(\lambda S)^*=\overline{\lambda}S^*$
for every~$\lambda\in \C$.
\item
Show that~$ST$ is adjoint to $T^*S^*$ (and so $(ST)^*=T^*S^*$).
\end{enumerate}
We will, of course, show
that every bounded operator on a Hilbert space is adjointable,
see~\sref{bounded-operator-adjoinable}.
But let us first show that~$\|\,\cdot\,\|$
defined in~\sref{hilb-def} is a norm,
which boils down to the following fact
about  $2\times 2$-matrices.
\end{point}
\begin{point}{130}[positive-2x2matrix]{Lemma}%
For a positive matrix $A\equiv 
\left(\begin{smallmatrix}p & \overline{c} \\ c & q\end{smallmatrix}\right)$
(i.e.~$\left(
\begin{smallmatrix}\overline{u}&\overline{v}\end{smallmatrix}\right)
A
\left(\begin{smallmatrix}u \\ v \end{smallmatrix}\right) \,\geq \, 0$
for all~$u,v\in \C$),
we have
$p,q\geq 0$, and $\left|c\right|^2 \leq pq$.
\begin{point}{140}{Proof}%
Let~$u,v\in\C$ be given.
We have
\begin{equation*}
0\ \leq\ 
\left(\begin{smallmatrix}\overline{u}&\overline{v}\end{smallmatrix}\right)
A
\left(\begin{smallmatrix}u \\ v \end{smallmatrix}\right)
\ = \ 
\left|u\right|^2 p\,+\, 
\overline{u}v\,\overline{c} \,+\,
u\overline{v}\,c \,+\,
\left|v\right|^2 q.
\end{equation*}
By taking~$u=1$ and $v=0$, we see that~$p\geq 0$,
and similarly $q\geq 0$.

The trick to see that~$\left|c\right|^2\leq pq$
is to
take~$v=1$ and $u=t\overline{c}$ with~$t\in \R$:
\begin{equation*}
0 \ \leq\ p\left|c\right|^2t^2
\,+\,2\left|c\right|^2t 
\,+\, q.
\end{equation*}
If~$p=0$, then~$-2\left|c\right|^2t \leq q $
for all~$t\in \R$,
which implies that~$\left|c\right|^2=0=pq$.

Suppose that~$p>0$.
Then taking~$t=-p^{-1}$ we see that
\begin{equation*}
0 \ \leq\ \left|c\right|^2p^{-1}
\,-\,2\left|c\right|^2p^{-1} 
\,+\, q \ = \ -\left|c\right|^2p^{-1}\,+\,q.
\end{equation*}
Rewriting gives us
 $\left|c\right|^2\leq pq$.\qed
\end{point}
\end{point}
\begin{point}{150}[inner-product-basic]{Exercise}%
Let~$\left<\,\cdot\,,\,\cdot\,\right>$
be an inner product on a vector space~$V$.
Show that
the formula~$\Define{\|x\|}=\smash{\sqrt{\left<x,x\right>}}$%
\index{$"\"|\,\cdot\,"\"|$, norm!on a pre-Hilbert space}
defines a seminorm on~$V$,
that is,
$\|x\|\geq 0$,
$\|\lambda x\|=\left|\lambda\right|\|x\|$,
and---the \Define{triangle inequality}---$\|x+y\|\leq \|x\|+\|y\|$
for all~$\lambda\in \C$ and~$x,y\in V$.

Moreover, prove that~$\|\,\cdot\,\|$
is a norm when~$\left<\,\cdot\,,\,\cdot\,\right>$
is definite;
and for~$x,y\in V$:
\begin{enumerate}
\item
The \Define{Cauchy--Schwarz inequality}:%
\index{Cauchy--Schwarz inequality!for $\C$-valued inner products}
$\left|\left<x,y\right>\right|^2\,\leq\, \left<x,x\right>
\,\left<y,y\right>$;
\item
\Define{Pythagoras' theorem}:%
\index{Pythagoras' theorem}
$\|x\|^2+\|y\|^2\,=\,\|x+y\|^2$ when~$\left<x,y\right>=0$;
\item
The \Define{parallelogram law}:%
\index{parallelogram law}
$\|x\|^2\,+\,
\|y\|^2
\,= \,
\frac{1}{2}(\,\|x+y\|^2\,+\,\|x-y\|^2\,)$;
\item
\label{polarization-identity}%
The \Define{polarisation identity}:%
\index{polarisation identity!for an inner product}
$\left<x,y\right> \,=\, \frac{1}{4}\sum_{n=0}^3i^n\|i^nx+y\|^2$.
\end{enumerate}

(Hint: prove the Cauchy--Schwarz inequality
before the triangle inequality
by applying~\sref{positive-2x2matrix} to the matrix
$\smash{\bigl(\begin{smallmatrix}
\smash{\left<x,x\right>} & \smash{\left<x,y\right>} \\
\smash{\left<y,x\right>} & \smash{\left<y,y\right>}
\end{smallmatrix}\bigr)}$.
Then prove $\|x+y\|^2\leq (\|x\|+\|y\|)^2$
using the inequalities~$\left<x,y\right>+\left<y,x\right>
\leq 2\left|\left<x,y\right>\right| \leq 2\|x\|\|y\|$.)
\end{point}
\begin{point}{160}[operators-cstar-identity]{Lemma}%
For an adjointable operator~$T$ on a pre-Hilbert space~$\scrH$
\begin{equation*}
\|T^*T\|\ =\ \|T\|^2\qquad\text{and}\qquad\|T^*\|\ =\ \|T\|.
\end{equation*}%
\spacingfix{}%
\begin{point}{170}{Proof}%
If~$T=0$, then~$T^*=0$, and the statements are surely true.

Suppose~$T\neq 0$ (and so~$T^*\neq 0$).
Since $\|Tx\|^2=\left<Tx,Tx\right>=\left<x,T^*Tx\right>
\leq \|x\|\,\|T^*Tx\|\leq \|x\|^2\|T^*T\|$
for every~$x\in \scrH$
by Cauchy--Schwarz,
we have $\|T\|^2\leq \|T^*T\|$.
Since~$\|T^*T\|\leq \|T^*\|\|T\|$
and $\|T\|\neq 0$,
it follows that~$\|T\|\leq \|T^*\|$.
Since by a similar reasoning $\|T^*\|\leq \|T\|$,
we get~$\|T\|=\|T^*\|$.
But then $\|T\|^2\leq \|T^*T\|\leq \|T^*\|\|T\|=\|T\|^2$,
and so $\|T\|^2=\|T^*T\|$.\qed
\end{point}
\end{point}
\begin{point}{180}{Exercise}%
Given a Hilbert space~$\scrH$
show that the adjointable operators
form a closed subspace of~$\scrB(\scrH)$.
\end{point}
\begin{point}{190}[ketbra]{Exercise}%
Let~$x$ and~$y$ be vectors from a Hilbert space~$\scrH$.
\begin{enumerate}
\item
Show that $\Define{\ketbra{x}{y}}\colon\, z\mapsto\left<y,z\right>x$
\index{*ketbra@$\ketbra{x}{y}$, with $x,y\in\scrH$}
defines a bounded operator
$\scrH\to\scrH$,
and, moreover, that~$\|\,\ketbra{x}{y}\,\|=\|x\|\|y\|$.
\item
Show that~$\ketbra{x}{y}$
is adjointable,
and~$(\ketbra{x}{y})^*=\ketbra{y}{x}$.
\end{enumerate}
\spacingfix{}
\end{point}%
\end{parsec}%
\begin{parsec}{50}[hilb-adjoint]%
\begin{point}{10}[adjoinables-cstar-algebra]%
At this point
it is clear that the vector space of adjointable operators
on a Hilbert space forms a $C^*$-algebra.
So to prove that $\scrB(\scrH)$
is a $C^*$-algebra,
it remains to be shown that every bounded operator
is adjointable (which we'll do in~\sref{bounded-operator-adjoinable}).
We first show that each bounded functional $f\colon \scrH\to \C$
has an adjoint, see~\sref{riesz-representation-theorem},
for which we need the (existence and) properties of ``projections''
on (closed) linear subspaces:
\end{point}
\begin{point}{20}[projection-on-closed-linear-subspace]{Definition}
Let~$x$ be an element of a pre-Hilbert space~$\scrH$.
We say that an element~$y$ of a linear subspace~$C$
of~$\scrH$ is a \Define{projection of~$x$ on~$C$}%
\index{projection!of~$x$ on~$C$}
if
\begin{equation*}
\|x-y\|\,=\,\min\{\,\|x-y'\|\colon \,y'\in C\,\}.
\end{equation*}
(In other words,~$y$ is one of the elements of~$C$ closest to~$x$.)
\end{point}
\begin{point}{30}{Exercise}%
We'll see in~\sref{projection-theorem}
that on a \emph{closed}
linear subspace
every vector has a projection.
For arbitrary linear subspaces this
isn't so:
show that the only vectors in~$\ell_2$
having 
a projection on the linear subspace~$c_{00}$
(from \sref{hilb-basic-examples})
are the vectors in~$c_{00}$ themselves.
\end{point}
\begin{point}{40}{Lemma}%
Let~$\scrH$ be a pre-Hilbert space,
and let $x,e\in\scrH$ with
$\|e\|=1$.

Then~$y=\left<e,x\right>e$ is the unique projection of~$x$ on~$e\C$.
\begin{point}{50}{Proof}%
Let~$y'\in e\C$
with~$y'\neq y$
be given.
To prove that~$y$
is the unique projection of~$x$ on $e\C$
it suffices to show that $\|x-y\|<\|x-y'\|$.
Since~$y'\neq y\equiv \left<e,x\right>e$,
there is~$\lambda\in \C$, $\lambda\neq 0$ 
with $y'=(\lambda+\left<e,x\right>)e$.

Note that $\left<e,y\right>=\left<e,\left<e,x\right>e\right>=
\left<e,x\right>\left<e,e\right>
= \left<e,x\right>$,
and so~$\left<e,x-y\right>=0$.
Then~$y'-y\equiv \lambda e$ and~$x-y$ are orthogonal too,
and thus, by Pythagoras'~theorem (see~\sref{inner-product-basic}),
we have $\|y'-x\|^2
=\|y'-y\|^2+\|y-x\|^2\equiv \left|\lambda\right|^2+\|x-y\|^2
>\|x-y\|^2$, because~$\lambda\neq 0$.
Hence~$\|y'-x\|>\|y-x\|$.\qed
\end{point}
\end{point}
\begin{point}{60}[hilb-projection-basic]{Exercise}%
Let~$y$ be a projection of an element~$x$ of a pre-Hilbert space~$\scrH$
on a linear subspace~$C$.
Show that~$y$ is a projection of~$x$ on $y\C$.
Conclude that~$y$ is the unique projection of~$x$ on~$C$,
and that~$\left<y,x-y\right>=0$.
Show that~$y+c$ is the projection of~$x+c$ on~$C$
for every~$c\in C$.
Conclude that~$\left<y',x-y\right>\equiv\left<y',(x+y'-y)-y'\right>=0$ 
for every~$y'\in C$.
\end{point}
\begin{point}{70}[projection-theorem]{Projection Theorem}%
\index{Projection Theorem}%
Let~$C$ be a closed linear subspace
of a Hilbert space~$\scrH$.
Each~$x\in \scrH$
has a unique projection~$y$ on~$C$,
and $\left<y',y\right>=\left<y',x\right>$ for~$y'\in C$.
\begin{point}{80}{Proof}%
We only need to show that there is a projection~$y$
of~$x$ on~$C$,
because~\sref{hilb-projection-basic}
gives us that such~$y$ is unique and satisfies
$\left<y',y\right> = \left<y',x\right>$ for all~$y'\in C$.

Write~$r:=\inf\{\,\|x-y'\|\colon\, y'\in C\,\}$,
and pick a sequence $y_1,y_2,\dotsc \in C$
such that $\|x-y_n\|\rightarrow r$.
We will show that~$y_1,y_2,\dotsc$ is Cauchy.
Let~$\varepsilon >0$
be given,
and pick~$N$ such that $\|y_n-x\|^2\leq r^2+\frac{1}{4}\varepsilon$
for all~$n\geq N$.
Let~$n,m\geq N$ be given.
Then since $\frac{1}{2}(y_n+y_m)$
is in~$C$, we have
$\|y_n+y_m-2x\|\equiv 
2\|\frac{1}{2}(y_n+y_m)-x\|\geq 2r$,
and so by the parallelogram law (see \sref{inner-product-basic}),
\begin{alignat*}{3}
\|y_n-y_m\|^2
\ &\equiv\ 
\|(y_n-x)-(y_m-x)\|^2\\
\ &=\ 
2\|y_n-x\|^2 + 2\|y_m-x\|^2
- \|y_n+y_m-2x\|^2\\
\ &\leq\ 
4r^2 + \varepsilon - 4r^2 \ \leq \ \varepsilon.
\end{alignat*}
Hence~$y_1,y_2,\dotsc$ is Cauchy,
and converges to some~$y\in C$,
because~$\scrH$ is complete and~$C$ is closed.
It follows easily that~$\|x-y\|=r$,
and thus~$y$ is the projection of~$x$ on~$C$.\qed
\end{point}
\end{point}
\begin{point}{90}[riesz-representation-theorem]{Riesz'~Representation Theorem}%
\index{Riesz' Representation Theorem}%
Let~$\scrH$ be a Hilbert space.
For every bounded linear map~$f\colon \scrH\to\C$
there is a unique vector~$x\in \scrH$
with $\left<x,\,\cdot\,\right>=f$.
\begin{point}{100}{Proof}%
If~$f=0$, then $x=0$ does the job.
Suppose that~$f\neq 0$.
There is an~$x'\in\scrH$ with~$f(x')=1$.
Note that~$\ker(f)$ is closed, because~$f$
is bounded.
So by~\sref{projection-theorem},
we know that~$x'$
has a projection~$y$ on~$\ker(f)$,
and $\left<x',z\right>=\left<y,z\right>$
for all~$z\in \ker(f)$.
Then for~$x'':=x'-y$,
we have $f(x'')=1$ and~$\left<x'',y'\right>=0$
for all~$y'\in \ker(f)$.
Given $z\in \scrH$,
we have $f(\,z-f(z)x''\,)=0$,
so~$z-f(z)x''\in \ker(f)$,
and thus~$0=\left<x'',z-f(z)x''\right>\equiv \left<x'',z\right>-f(z)\|x''\|^2$.
Hence writing $x:=x''\|x''\|^{-2}$
we have~$f(z)=\left<x,z\right>$
for all~$z\in \scrH$.

Finally, uniqueness of~$x$ follows from~\sref{uniqueness-adjoint}.\qed
\end{point}
\end{point}
\begin{point}{110}[bounded-operator-adjoinable]{Exercise}%
Prove that every bounded operator~$T$ on a Hilbert space~$\scrH$
is adjointable, as follows.
Let~$x\in \scrH$ be given.
Prove that~$\left<x,T(\,\cdot\,)\right>\colon \scrH\to \C$
is a bounded linear map.
Let~$Sx$ be the unique vector with $\left<Sx,\,\cdot\,\right>
=\left<x,T(\,\cdot\,)\right>$,
which exists by~\sref{riesz-representation-theorem}.
Show that~$x\mapsto Sx$
gives a bounded linear map $S$, which is adjoint to~$T$.
\end{point}
\begin{point}{120}%
Thus the bounded operators
on a Hilbert space~$\scrH$
form a $C^*$-algebra~$\scrB(\scrH)$%
\index{BH@$\scrB(\scrH)$!as a $C^*$-algebra}
as described in~\sref{example-hilb}.
We will return to Hilbert spaces
in~\sref{gelfand-naimark},
where we show that every $C^*$-algebra
is isomorphic to a $C^*$-subalgebra of
a $\scrB(\scrH)$.
\end{point}
\end{parsec}
\begin{parsec}{60}%
\begin{point}{10}%
Here is a non-trivial
example of a Hilbert space
that will be used later on.
\end{point}
\begin{point}{20}[hilb-sum]{Proposition}%
Given a family  $(\scrH_i)_{i\in I}$ 
of Hilbert spaces,
the vector space
\begin{equation*}
	\textstyle
	\Define{\bigoplus_i \scrH_i} \ =\ \{\ 
		x\in \prod_i \scrH_i\colon\ 
	\sum_i \|x_i\|^2 <\infty \ \}.
\end{equation*}%
\index{direct sum!of Hilbert spaces}%
\index{$\bigoplus$, direct sum!$\bigoplus_i\scrH_i$, of Hilbert spaces}
is a Hilbert space
when endowed with the inner product
$\left<x,y\right>=\sum_i \left<x_i,y_i\right>$.
\begin{point}{30}{Proof}%
To begin with
we must show that~$\sum_i \left<x_i,y_i\right>$
converges for~$x,y\in\bigoplus_i \scrH_i$.
Given~$\varepsilon>0$
we must
find a finite subset~$G$ of~$I$
such that~$ \left|\sum_{i \in F} \left<x_i,y_i\right>\right|
\leq \varepsilon$ for all finite $F\subseteq I\backslash G$.
Since an obvious application
of the Cauchy--Schwarz inequality
gives
us that for every finite subset~$F$ of~$I$
\begin{equation*}
	\Bigl|\sum_{i\in F}
	\left<x_i,y_i\right>\Bigr|^2
	\ \leq\ 
	\sum_{i\in F}\|x_i\|^2
	\, \sum_{i\in F}\|y_i\|^2,
\end{equation*}
any~$G\subseteq I$
with $\sum_{i\in I\backslash G} \|x_i\|^2 \leq \sqrt{\varepsilon}$
	and~$\sum_{i\in I\backslash G} \|y_i\|^2 \leq \sqrt{\varepsilon}$
will do.

It is easily seen that
$\left<x,y\right>:=\sum_i \left<x_i,y_i\right>$
gives a definite inner product on~$\bigoplus_i \scrH_i$; the
remaining difficulty
lies in showing that the resulting norm is complete.
To this end, let $x_1,x_2,\dotsc$ be a Cauchy sequence 
in~$\bigoplus_{i\in I}\scrH_i$;
we must show that it converges to some~$x_\infty\in \bigoplus_i \scrH_i$.
We do the obvious thing:
since for every~$i\in I$
the sequence $(x_1)_i, (x_2)_i,\dotsc$
is Cauchy in~$\scrH_i$
we may define $(x_\infty)_i:=\lim_n (x_n)_i$,
and thereby get an element $x_\infty$ of~$\prod_i \scrH_i$.
Since for each finite subset~$F$ of~$I$
we have $\sum_{i\in F} \|(x_\infty)_i\|^2
=\lim_n \sum_{i\in F} \|(x_n)_i\|^2
\leq \lim_n \|x_n\|^2$,
we have $\sum_{i\in I} \|(x_\infty)_i\|^2 
\leq \lim_n \|x_n\|^2 <\infty$,
and so~$x_\infty\in\bigoplus_i \scrH_i$.

It remains to be shown that~$x_1,x_2,\dotsc$
converges to~$x_\infty$
(not only coordinatewise but also)
with respect to the inner product on~$\bigoplus_i \scrH_i$.
Given~$\varepsilon >0$
pick $N$ such that $\|x_n - x_m\|\leq \frac{1}{2\sqrt{2}}\varepsilon$
for all~$n,m\geq N$.
We claim that for such~$n$
we have $\|x_\infty -x\|\leq \varepsilon$.
Indeed, first note that since the sum
\begin{equation*}
	\sum_{i\in I} \|(x_\infty)_i - (x_n)_i \|^2
\ \equiv\ 
\sum_{i\in F} \|(x_\infty)_i - (x_n)_i \|^2
\ +\ 
\sum_{i\in I\backslash F} \|(x_\infty)_i - (x_n)_i \|^2
\end{equation*}
converges (to $\|x_\infty - x_n\|^2$),
we can find
a finite subset~$F$ (depending on~$n$)
such that second term in the right-hand side above
is smaller than~$\frac{1}{2}\varepsilon^2$.
To see that the first term is also below~$\frac{1}{2}\varepsilon^2$,
begin by noting that for every~$m$,
\begin{equation*}
\Bigl(\,\sum_{i\in F} \|(x_\infty)_i - (x_n)_i \|^2\,\Bigr)^{\nicefrac{1}{2}}
\ \leq \ 
\Bigl(\,\sum_{i\in F} \|(x_\infty)_i - (x_m)_i \|^2\,\Bigr)^{\nicefrac{1}{2}}
\ +\ 
\Bigl(\,\sum_{i\in F} \|(x_m)_i - (x_n)_i \|^2\,\Bigr)^{\nicefrac{1}{2}}.
\end{equation*}
Since~$F$ is finite,
and~$(x_m)$ converges
to~$x_\infty$ coordinatewise
we can find an~$m$ large enough
that the first term on the right-hand side above
is below~$\smash{\frac{1}{2\sqrt{2}}\varepsilon}$.
If we choose~$m\geq N$
we see that the second term is below $\smash{\frac{1}{2\sqrt{2}}\varepsilon}$
as well,
and we conclude that~$\|x_\infty-x_n\|\leq \varepsilon$.\qed
\end{point}
\end{point}

\end{parsec}
\section{The Basics}

\begin{parsec}{70}%
\begin{point}{10}%
Now that we have seen the most important examples
of $C^*$-algebras,
we can begin developing the theory.
We'll start easy with the self-adjoint elements:
\end{point}
\begin{point}{20}{Definition}%
Given an element $a$ of a $C^*$-algebra $\scrA$, 
\begin{enumerate}
\item we say that $a$ is \Define{self-adjoint}%
\index{self adjoint} if $a^* =a$, and
\item we write $\Define{\Real{a}}:= \frac{1}{2}(a+a^*)$
and $\Define{\Imag{a}}:=\frac{1}{2i}(a-a^*)$
for the \Define{real} and \Define{imaginary part}%
\index{real part!of an element of a $C^*$-algebra}%
        \index{$\Real{(\,\cdot\,)}$, real part!$\Real{a}$, of an element of a $C^*$-algebra}%
\index{imaginary part!of an element of a $C^*$-algebra}%
        \index{$\Imag{(\,\cdot\,)}$, imaginary part!$\Imag{a}$, of an element of a $C^*$-algebra}
of~$a$, respectively.
\end{enumerate}
The set of self-adjoint elements of~$\scrA$
is denoted by~\Define{$\sa{\scrA}$}.%
\index{$\Real{(\,\cdot\,)}$, real part!$\Real{\scrA}$, of a $C^*$-algebra}
\end{point}
\begin{point}{30}[cstar-involution-basic]{Exercise}%
Let~$a$ be an element of a $C^*$-algebra~$\scrA$.
\begin{enumerate}
\item 
Show that $\Real{a}$ and $\Imag{a}$ are self-adjoint,
and  $a= \Real{a}+i\Imag{a}$.
\item
Show that if $a\equiv b+ic$ for self-adjoint elements $b$, $c$ of~$\scrA$,
then $b=\Real{a}$ and~$c=\Imag{a}$.
\item
Show that $\Real{(a^*)}=\Real{a}$ and $\Imag{(a^*)}=-\Imag{a}$.
\item 
Show that~$a$ is self-adjoint iff $\Real{a}=a$ iff $\Imag{a}=0$.
\item
Show that $a\mapsto \Real{a}$ and $a\mapsto \Imag{a}$
give $\R$-linear maps $\scrA\to\scrA$.
\item
Show that $\Imag{a} = -\Real{(ia)}$ and $\Real{a}=\Imag{(ia)}$.
\item
Show that $a^*a$ is self-adjoint,
and  $a^*a=\Real{a}^2+\Imag{a}^2+i(\Real{a}\Imag{a}-\Imag{a}\Real{a})$.
\item
Give an example of~$\scrA$ and~$a$ 
with  $\Real{a}\Imag{a} \neq \Imag{a}\Real{a}$.

(This inequality is a source of many technical difficulties.)
\item
Show that $a^*a+aa^* = 2(\Real{a}^2+\Imag{a}^2)$.
\item
The product of self-adjoint elements $b$, $c$ need not be self-adjoint;
show that, in fact, $bc$ is self-adjoint iff $bc=cb$.
\item
Show that $\|a^*\| = \|a\|$. (Hint:  $\|a\|^2=\|a^*a\|\leq \|a^*\|\|a\|$.)

\item
Show that $\|\Real{a}\|\leq \|a\|$ and $\|\Imag{a}\|\leq \|a\|$.
\item
Show that $\|a^2\|=\|a\|^2$ when~$a$ is self-adjoint.

However,
show that $\|a^2\|\neq \|a\|^2$ might occur
when~$a$ is not self-adjoint.
(Hint: $\bigl(
\begin{smallmatrix}
	0&1\\
	0&0
\end{smallmatrix}
\bigr)$.)
\end{enumerate}%
\spacingfix{}
\end{point}%
\end{parsec}%
\begin{parsec}{80}%
\begin{point}{10}{Notation}%
Recall that (in this text) every $C^*$-algebra~$\scrA$ has a unit, $1$.
Thus, for every scalar $\lambda\in \C$,
we have an element $\lambda\cdot 1$ of~$\scrA$,
which we will simply denote by~$\lambda$.
This should hardly cause any confusion,
for while an expression of an element of~$\scrA$
such as $i+2+5a$ (where $a\in \scrA$) 
may be interpreted in several ways,
the result is always the same.
\end{point}
\begin{point}{20}{Exercise}%
There is a subtle point regarding
the norm~$\|\lambda\|$ of a
scalar~$\lambda\in \C$ inside a $C^*$-algebra~$\scrA$:
	we do not always have~$\|\lambda\|=\left|\lambda\right|$
	on the nose.
\begin{enumerate}
\item
Indeed, show 
that $\|1\|=0\neq 1$ when~$\scrA=\{0\}$ is the trivial $C^*$-algebra.
\item 
Show that $\|\lambda\|\leq \left| \lambda\right|$ (in~$\C$).
\item
Show that~$\|\lambda\|=\left|\lambda\right|$
when~$\|\lambda\|$ and~$\left|\lambda\right|$
are interpreted as elements of~$\scrA$.
\end{enumerate}%
\spacingfix{}%
\end{point}%
\end{parsec}%
\begin{parsec}{90}%
\begin{point}{10}%
Let us now generalise the notion of a positive function
in~$C(X)$
to a positive element of a $C^*$-algebra.
There are several descriptions of
positive functions in~$C(X)$ in terms of the $C^*$-algebra structure
(see~\sref{cx-positive}) on which we can base such a  generalisation,
and while we will eventually see that these all yield the same notion
of positive element of a $C^*$-algebra (see~\sref{cstar-positive-final})
we base our definition of positive element (\sref{cstar-positive-def})
on the description that is perhaps
not most familiar,
but does give us the richest structure at this stage.
\end{point}
\begin{point}{20}[cx-positive]{Exercise}%
Let~$X$ be a compact Hausdorff space.
Show that for self-adjoint $f\in C(X)$, the following are equivalent.
\begin{enumerate}
\item \label{cx-positive-1}
$f(X)\subseteq [0,\infty)$;
\item
$f\equiv g^2$ for some $g\in \sa{C(X)}$;
\item
$f\equiv g^* g$ for some~$g\in C(X)$;
\item
$\|f-t\|\leq t$ for some $t\in \R$;
\item
$\|f-t\|\leq t$ for all $t\geq \frac{1}{2}\|f\|$.
\end{enumerate}
(Hint: $\|f-t\|\leq t$ iff $-t\leq f-t\leq t$
iff $0\leq f\leq 2t$.)
\begin{point}{30}{Exercise}%
To see how condition~\ref{cx-positive-1}
can be expressed in terms of the $C^*$-algebra structure of~$C(X)$,
prove that  $\lambda\in f(X)$ iff $f-\lambda$
is not invertible.
\end{point}
\end{point}
\begin{point}{40}[cstar-positive-def]{Definition}%
A self-adjoint element~$a$ of a $C^*$-algebra~$\scrA$ is called
\Define{positive}%
\index{positive!element of a $C^*$-algebra}
if $\|a-t\|\leq t$
for some~$t\in \R$.
We write $\Define{a\leq b}$%
\index{((leq@$\leq$, order!on a $C^*$-algebra}
for $a,b\in\scrA$ when $b-a$ is positive,
and we denote the set of positive elements of~$\scrA$
by~$\Define{\pos{\scrA}}$.%
\index{$(\,\cdot\,)_+$, positive part!$\scrA_+$,
of a $C^*$-algebra}
\begin{point}{41}
Given elements $a$ and~$b$ of a $C^*$-algebra~$\scrA$
we denote by $\Define{[a,b]_\scrA}$,
or sometimes simply $\Define{[a,b]}$,%
\index{$[a,b]$, interval!$[a,b]_\scrA$, in a $C^*$-algebra}
the set of elements~$c$ of~$\scrA$
with $a\leq c\leq b$.
\end{point}
\begin{point}{50}{Remark}%
One advantage of this
definition
over, say, taking the elements of the form~$a^*a$
to be positive,
is that it is immediately clear
that an element~$b$ of a $C^*$-subalgebra~$\scrB$
of a $C^*$-algebra~$\scrA$ is positive in~$\scrB$
iff~$b$ is positive in~$\scrA$---that is,
`positive permanence' comes for free (cf.~\sref{spectral-permanence}).
Another advantage is
that it's also pretty easy to see 
that the sum of such positive elements
is again positive, see~\sref{cstar-positive-sum}.
\end{point}
\begin{point}{51}{Remark}%
Note that when an element~$a$
of a $C^*$-algebra
is positive on the grounds that $\|a-t\|\leq t$
for some~$t\in \R$,
then this number~$t$ must be positive,
and we even have $t\geq \frac{1}{2}\|a\|$,
since $\|a\|-\|t\|\leq \|a-t\|\leq t$.
There's nothing special about this~$t$:
we'll see in~\sref{cstar-positive-1}
that $\|a-s\|\leq s$
    for all~$s\geq \frac{1}{2} \|a\|$ and positive~$a$.
\end{point}
\end{point}
\begin{point}{60}{Example}%
We'll see in~\sref{hilb-positive-operators},
that a bounded operator~$T$ on a Hilbert space~$\scrH$
is positive iff~$\left<x,Tx\right>\geq 0$ for all~$x\in\scrH$.
\end{point}
\begin{point}{70}[cstar-positive-sum]{Lemma}%
Let~$a,b$ be positive elements of a $C^*$-algebra.
Then $a+b$ is positive.
\begin{point}{80}{Proof}
Since~$a\geq 0$,
there is~$t\in \R$
    with $\|a-t\|\leq t$.
Similarly, there is~$s\in \R$
with $\|b-s\|\leq s$.
Then $\|a+b-(t+s)\|\leq \|a-t\|+\|b-s\|\leq t+s$.\qed
\end{point}
\end{point}
\begin{point}{90}{Exercise}%
Given an element~$a$
of a $C^*$-algebra~$\scrA$
with~$0\leq a\leq 1$
(which is called an \Define{effect})%
\index{effect!in a $C^*$-algebra}
show that 
the \Define{orthosupplement} $\Define{a^\perp} :=1-a$%
\index{orthosupplement!operation in a $C^*$-algebra}%
\index{$(\,\cdot\,)^\perp$!$a^\perp$, orthosupplement of an effect}
is an effect too.
\end{point}
\begin{point}{100}[cstar-positive]{Exercise}%
Let~$\scrA$ be a $C^*$-algebra.
\begin{enumerate}
\item
Show that~$\pos{\scrA}$ is a \emph{cone}:
$0\in \pos{\scrA}$,
$a+b\in \pos{\scrA}$ for all $a,b\in\pos{\scrA}$,
and
$\lambda a\in \pos{\scrA}$  
for all $a\in \pos{\scrA}$ and $\lambda\in [0,\infty)$.
Conclude that~$\leq$ is a preorder.
\item
Show that~$1$ is positive, and  $-\|a\|\leq a \leq \|a\|$
for every self-adjoint element~$a$ of~$\scrA$.
(Thus $1$ is an \emph{order unit} of~$\sa{\scrA}$.)
\item
The behaviour of positive elements may be surprising:
give an example of positive elements $a$ and~$b$
from a $C^*$-algebra
such that $ab$ is not positive.
\item
Given a self-adjoint element~$a$ of~$\scrA$ define
\begin{equation*}
\|a\|_o \ = \ \inf\{\ \lambda\in[0,\infty)\colon \ 
-\lambda\leq a\leq \lambda\ \}.
\end{equation*}
Show that $\|-\|_o$ is a seminorm on~$\sa{\scrA}$,
and that~$\|a\|_o\leq \|a\|$
for all~$a\in\sa{\scrA}$.

Prove that $0\leq a\leq b$ implies that~$\|a\|_o\leq\|b\|_o$
for $a,b\in\sa{\scrA}$.

\item
There is not much more that can easily be
proven about positive elements, at this point,
but don't take my word for it:
try to prove the following facts
about a self-adjoint element~$a$ of~$\scrA$ directly.
\begin{enumerate}
\item $a^2$ is positive;
\item if $a$ is the limit of positive $a_n\in\scrA$,
then $a$ is positive;
\item if $a\geq -\frac{1}{n}$ for all~$n\in \N$, then $a\geq 0$;
\item  $\|a\|=\|a\|_o$;
\item $a=0$ when~$0\leq a\leq 0$.
\end{enumerate}
We will prove these facts
when we return to the positive elements in~\sref{cstar-positive-2}.
\end{enumerate}%
\spacingfix{}%
\end{point}%
\end{parsec}%
\begin{parsec}{100}%
\begin{point}{10}%
Let us spend some words
on the morphisms between $C^*$-algebras.
\end{point}
\begin{point}{20}[maps]{Definition}
A linear map $f\colon \scrA \to \scrB$
between $C^*$-algebras
is called
\begin{enumerate}
\item
\Define{\textbf{m}ultiplicative}%
\index{multiplicative!map between $C^*$-algebras}
if $f(ab)=f(a)f(b)$ for all $a,b\in\scrA$;
\item
\Define{\textbf{i}nvolution preserving}%
\index{involution preserving!map between $C^*$-algebras}
if $f(a^*)=f(a)^*$ for all~$a\in\scrA$;
\item
\Define{\textbf{u}nital}%
\index{unital!map between $C^*$-algebras}
if $f(1)=1$;
\item
\Define{\textbf{s}ub\textbf{u}nital}%
\index{subunital map between $C^*$-algebras}
if $f(1)\leq 1$;
\item
\Define{\textbf{p}ositive}%
\index{positive!map between $C^*$-algebras}
if $f(a)$ is positive
for every positive $a\in\scrA$, and
\item
\Define{\textbf{c}ompletely \textbf{p}ositive}%
\index{completely positive!map between $C^*$-algebras}%
\index{positive!completely~$\sim$ map between $C^*$-algebras}
if $\sum_{i,j} b_i^*\,f(\,a_i^*a_j\,)\,b_j$ is positive
for all~$a_1,\dotsc,a_n\in \scrA$, and $b_1,\dotsc,b_n\in\scrB$
(see Remark~5.1 of~\cite{paschke}).
\item
(For \emph{\textbf{n}ormal} maps,
we refer to~\sref{bh-normal} and~\sref{p-uwcont}.)
\end{enumerate}%
\spacingfix{}%
\begin{point}{30}%
We use the bold letters as abbreviations,
so for instance,
$f$ is \Define{pu}%
\index{pu-map} if it is positive and unital,
and a \Define{miu-map}%
\index{miu-map}
is a multiplicative, involution preserving,
unital linear map between $C^*$-algebras
(which is usually called a \Define{unital $*$-homomorphism}%
\index{homomorphism@$*$-homomorphism}).

We'll denote the category of $C^*$-algebras
and miu-maps by~$\Define{\Cstar{miu}}$,%
\index{Cstar@$\Cstar{}$: $\Cstar{miu}$, $\Cstar{cpu}$, \dots}
and
the subcategory of commutative $C^*$-algebras
by~$\Define{\cCstar{miu}}$.%
\index{cCstar@$\cCstar{}$: $\cCstar{miu}$, $\cCstar{pu}$,\dots}
We'll use similar notation
for the other classes of maps,
but
will,
naturally, only mention $\Cstar{cpu}$
after having established that cp-maps are closed under composition.

The advantages of completely positive maps
become apparent
only later on 
when we start dealing with matrices (see~\sref{n-pos})
and the tensor product (see~\sref{tensor-functorial}).
\end{point}
\end{point}
\begin{point}{40}[cstar-p-implies-i]{Lemma (``p$\Rightarrow$i'')}
A positive map $f\colon \scrA\to\scrB$ between
$C^*$-algebras is involution preserving.
\begin{point}{50}{Proof}%
Let~$a\in \scrA$ be given. We must show that~$f(a^*)=f(a)^*$.

But first we'll show that if~$a$ is self adjoint,
then so is~$f(a)$.
Indeed, since $\|a\|$ and $\|a\|-a$ are positive (see~\sref{cstar-positive}),
we see that $f(\|a\|)$ and $f(\|a\|-a)$ are positive,
and so~$f(a)=f(\|a\|)-f(\|a\|-a)$ being positive is self adjoint.

It follows that $\Real{f(a)}=f(\Real{a})$
and $\Imag{f(a)}=f(\Imag{a})$ (for~$a\in\scrA$),
because $f(a)\equiv f(\Real{a})+if(\Imag{a})$,
and~$f(\Real{a})$ and~$f(\Imag{a})$
are self adjoint
(see~\sref{cstar-involution-basic}).

Hence $f(a^*)\equiv f(\Real{a}-i\Imag{a})
=\Real{f(a)}-i\Imag{f(a)}\equiv f(a)^*$.\qed
\end{point}
\end{point}
\begin{point}{60}{Remark}%
Other important relations between these types of morphisms
can only be established later on
once we have a firmer grasp on the positive elements.
We will then see  
that every mi-map 
is completely positive (in~\sref{cp}),
and that every completely positive map is positive 
(in~\sref{astara-pos-basic-consequences}).
\begin{point}{61}%
Note that we didn't bother to include
an abbreviation for bounded linear maps in our list, \sref{maps}.
That's because we'll see in~\sref{weak-russo-dye} that any positive
map between $C^*$-algebras is automatically bounded.
\end{point}
\end{point}
\begin{point}{70}{[Moved to~\sref{cstar-product-2}.]}%
\end{point}
\begin{point}{80}{[Moved to~\sref{cstar-equaliser-1}.]}%
\begin{point}{90}{[Moved to~\sref{cstar-no-pu-equalisers}.]}%
\end{point}%
\end{point}%
\end{parsec}%
%
%
\begin{parsec}{110}%
\begin{point}{10}%
After having visited the positive elements,
let us explore our second landmark,
the  invertible elements
of a $C^*$-algebra,
whose role 
is as important as it is technical.
This paragraph culminates in what is essentially
 \emph{spectral permanence} (\sref{spectral-permanence}):
the fact that if an element $a$ of a $C^*$-subalgebra $\scrB$
of a $C^*$-algebra~$\scrA$
is invertible in~$\scrA$,
then~$a$ is already invertible in~$\scrB$,
see~\sref{inverse-permanence}.
\end{point}
\begin{point}{20}[geometric]{Lemma}%
\index{geometric series}%
Let~$a$ be an element of a $C^*$-algebra~$\scrA$ with~$\|a\|<1$.
Then~$a^\perp\equiv 1-a$ has an inverse,
namely~$(a^\perp)^{-1}= \sum_{n=0}^\infty\, a^n$.
Moreover, this series converges absolutely,
that is,
$\sum_{n=0}^\infty \|a^n\|<\infty$.
\begin{point}{30}{Proof}%
Note that
$(1-\|a\|)\,(1+\|a\|+\|a\|^2+\dotsb+\|a\|^N) \,=\, 1-\|a\|^{N+1}$,
and so 
\begin{equation*}
\sum_{n=0}^N \|a\|^n \ =\  \frac{1-\|a\|^{N+1}}{1-\|a\|}
\end{equation*}
for every~$N$.
Thus,
since $\|a\|^N$ converges to~$0$
(because\footnote{
In case you've never seen the argument:
the limit $b:=\lim_N \|a\|^N$
exists, because $\|a\|\geq \|a\|^2\geq \dotsb\geq 0$,
and is zero
because $\|a\|b=\lim_N \|a\|^{N+1}=b$
and~$\|a\| < 1$.}
 $\|a\|<1$),
we  get $\sum_{n=0}^\infty \|a\|^n = (1-\|a\|)^{-1}$.
Note that since~$\|a^n\|\leq \|a\|^n$ for every~$n$,
    this entails that $\sum_{n=0}^\infty \|a^n\|\leq (1-\|a\|)^{-1}
    <\infty$.
\begin{point}{40}%
Note that $a^N$ norm converges to~$0$,
because $\|a\|^N$ converges to~$0$.
Also (but slightly less obvious),
$\sum_n a^n$ norm converges,
because~$\sum_n \|a\|^n$ converges.
\end{point}
\begin{point}{50}%
Thus, taking the norm limit
on both sides of $(1-a)(1+a+a^2+\dotsb a^N) = 1-a^{N+1}$,
gives us $(1-a)(\sum_n a^n) = 1$.
Since we can derive $(\sum_n a^n)(1-a) = 1$
in a similar manner, 
we see that $\sum_n a^n$ is the inverse of~$1-a$.\qed
\end{point}
\end{point}
\end{point}
\begin{point}{60}[spectrum-bounded]{Exercise}%
\index{invertible!element of a $C^*$-algebra}%
Let~$a$ be an element of a $C^*$-algebra~$\scrA$.
\begin{enumerate}
\item
Show that $a-\lambda$ is invertible
for every~$\lambda\in\C$ with~$\|a\|< \left|\lambda\right|$.
\item
Show that $a-b$ is invertible
when~$b\in\scrA$ is invertible and $\|a\| < \|b\|$.
\item
Show that $U:=\{\ b\in\scrA\colon\ \text{$b$ is invertible}\ \}$
is an open subset of~$\scrA$.
\end{enumerate}%
\spacingfix{}
\end{point}%
\begin{point}{70}[geometric-convergence]{Lemma}%
\index{geometric series}%
For a self-adjoint element~$a$ of~$\scrA$
the series $\sum_n a^n$ 
converges iff~$\|a\|<1$;
and in that case converges absolutely.
\begin{point}{80}{Proof}%
We have already seen in~\sref{geometric}
that~$\sum_n a^n$
converges absolutely when~$\|a\|<1$.
Now, if $\sum_n a^n$ converges,
then~$\|a^n\|$
(being the norm of the difference
between consecutive partial sums of~$\sum_n a^n$)
converges to~$0$.
    In particular, $\smash{\|a\|^{2^n}}$
    (being equal to $\|\smash{a^{2^n}}\|$
    by the $C^*$-identity)
converges to~$0$ too,
which only happens when~$\|a\|<1$.\qed
\end{point}
\begin{point}{90}[geometric-non-self-adjoint]{Remark}%
For non-self-adjoint elements~$a$ of~$\scrA$,
the convergence of~$\sum_n a^n$
is a more delicate matter.
Take for example
the matrix $A:=\bigl(\begin{smallmatrix}0&2\\0&0\end{smallmatrix}\bigr)$
for which the series $\sum_n A^n$ converges (to~$1+A$),
while~$\|A\|=2$ --- the problem being that
$\|A^2\|^{\nicefrac{1}{2}}$ differs from $\|A\|$.
In fact,
we'll see from~\sref{hadamard}
(although we won't need it)
that~$\sum_n a^n$
converges absolutely when $1>\limsup_n \|a^n\|^{\nicefrac{1}{n}}$,
and diverges when $1<\limsup_n \|a^n\|^{\nicefrac{1}{n}}$.
This begs the question
what happens when $1=\limsup_n \|a^n\|^{\nicefrac{1}{n}}$
--- which I do not know.
\end{point}
\end{point}
\begin{point}{100}[cstar-inv-continuous]{Lemma}%
Let~$\scrA$ be a $C^*$-algebra.
The assignment $a\mapsto a^{-1}$
gives a  continuous map
(from the set $\{\,b\in \scrA\colon\, \text{$b$ is invertible}\,\}$
to~$\scrA$.)
\begin{point}{110}[cstar-inv-continuous-1]{Proof}
(Based on Proposition 3.1.6 of~\cite{kr}.)

First we establish continuity at~$1$:
let~$a\in\scrA$ with $\|1-a\|\leq \frac{1}{2}$ be given;
we claim that~$a$ is invertible,
and~$\|1-a^{-1}\| \leq 2\|1-a\|$.

Indeed, since~$\|1-a\|\leq \frac{1}{2}<1$,
$a$ is invertible by~\sref{geometric},
and $a^{-1}=\sum_{n=0}^\infty (1-a)^n$.
Then~$\|1-a^{-1}\|=\|\sum_{n=1}^\infty (1-a)^n\|\leq \sum_{n=1}^\infty \|1-a\|^n
= \|1-a\|\, (1-\|1-a\|)^{-1}$.
Thus, as $\|1-a\|\leq\frac{1}{2}$,
we get $(1-\|1-a\|)^{-1}\leq 2$,
and so $\|1-a^{-1}\|\leq 2\|1-a\|$.
\begin{point}{120}%
Let~$a$ be an invertible element of~$\scrA$,
and let~$b\in\scrA$ with~$\|a-b\|\leq\frac{1}{2}\|a^{-1}\|$.
We claim that~$b$ is invertible,
and~$\|a^{-1}-b^{-1}\|\leq 2\|a-b\|\,\|a^{-1}\|^2$.
Since $\|a-b\|\leq \frac{1}{2}\|a^{-1}\|$
we have
$\|1-a^{-1}b\|\leq \|a^{-1}\|\,\|a-b\|\leq \frac{1}{2}$.
By~\sref{cstar-inv-continuous-1}, $a^{-1}b$ is invertible,
and $\|1-(a^{-1}b)^{-1}\|\leq 2\|1-a^{-1}b\|\leq 2\|a-b\|\,\|a^{-1}\|$.
Hence $\|a^{-1}-b^{-1}\| = \|(1-(a^{-1}b)^{-1})a^{-1}\|
\leq \|1-(a^{-1}b)^{-1}\|\,\|a^{-1}\|\leq 2 \|a-b\|\,\|a^{-1}\|^2$.\qed
\end{point}
\end{point}
\end{point}
%
%
\begin{point}{130}{Lemma}%
For a self-adjoint element~$a$ from a $C^*$-algebra,
$a-i$ is invertible.
\begin{point}{140}{Proof}%
(Based on Proposition 4.1.1(ii) of~\cite{kr}.)

The trick
is to 
write~$a-i\equiv (a+ni)\,-\,(n+1)i$
for sufficiently large~$n$,
because  
then
by~\sref{spectrum-bounded}
$a-i$
is invertible provided that~$n+1 > \|a+ni\|$.
Indeed, for~$n$ such that~$\|a\|<2n+1$,
we have $\|a+ni\|^2 = \|(a+ni)^*(a+ni)\|
= \|a^2+n^2\|
\leq \|a\|^2+n^2 < 2n+1+n^2 = (n+1)^2$,
and so $\|a+ni\| < n+1$.\qed
\end{point}
\end{point}
\begin{point}{150}[spectrum-self-adjoint-real]{Exercise}%
Let~$a$ be a self-adjoint element of a $C^*$-algebra.
\begin{enumerate}
\item
Show that~$a-\lambda$ is invertible for all $\lambda\in \C\backslash \R$.
\item
Show that $a^2-\lambda$ is invertible for all 
$\lambda\in \C\backslash[0,\infty)$.\\
(Hint: first prove that
 $a^2+1 \equiv (a+i)(a-i)$ is invertible.)

Conclude that $a^n-\lambda$ is invertible for all 
$\lambda\in\C\backslash[0,\infty)$ and \emph{even} $n\in\N$.
\item
Let~$n\in \N$ be \emph{odd}.
Show that $a^n-\lambda$ is invertible
for all~$\lambda\in \C\backslash[0,\infty)$
if and only if $a-\lambda$ is invertible
for all~$\lambda\in \C\backslash[0,\infty)$.\\
(Hint: show that
$a^n+1= \prod_{k=1}^n a+\zeta^{2k+1}$
where $\zeta=e^{\frac{\pi i}{n}}$.)
\end{enumerate}%
\spacingfix{}
\end{point}%
\begin{point}{160}[inverse-permanence]{Proposition}%
Let~$\scrA$ be a $C^*$-subalgebra
of a $C^*$-algebra $\scrB$.
Let~$a$ be a self-adjoint element of~$\scrA$,
which has an inverse, $a^{-1}$, in~$\scrB$.
Then~$a^{-1}\in\scrA$.
\begin{point}{170}{Proof}%
While we do not know yet that~$a$ is invertible in~$\scrA$,
we do know that~$a+\nicefrac{i}{n}$ 
has an inverse $(a+\nicefrac{i}{n})^{-1}$ in~$\scrA$
by~\sref{spectrum-self-adjoint-real}
for each~$n$
(using that $a$ is self-adjoint.)
Since~$a+\nicefrac{i}{n}$ converges to~$a$ in~$\scrB$ as~$n$ increases,
we see that $(a+\nicefrac{i}{n})^{-1}$ converges to~$a^{-1}$
in~$\scrB$ by~\sref{cstar-inv-continuous}.
Thus, as all~$(a+\nicefrac{i}{n})^{-1}$ are in~$\scrA$,
and~$\scrA$ is closed in~$\scrB$,
we see that~$a^{-1}$ is in~$\scrA$.\qed
\end{point}
\end{point}
\begin{point}{180}[improved-inverse-permanence]{Exercise}%
	Show that the assumption in~\sref{inverse-permanence} 
	that~$a$ is self-adjoint
may be dropped. 

(Hint: consider $a^*a$, see Proposition VIII.1.14 of~\cite{conway2013}.)
\end{point}
\begin{point}{190}[spectrum-of-element]{Definition}%
The \Define{spectrum},
\Define{$\spec(a)$},%
\index{sp@$\spec$, spectrum!$\spec(a)$, of an element of a $C^*$-algebra}
of an element $a$
of a $C^*$-algebra
is the set 
of complex numbers~$\lambda$
for which~$a-\lambda$ is not invertible.
\end{point}
\begin{point}{200}{Exercise}%
Verify the following examples.
\begin{enumerate}
\item
The spectrum of a continuous function~$f\colon X\to \R$
on a compact Hausdorff space~$X$
being an element of the $C^*$-algebra $C(X)$
is the image of~$f$, that is,
$\spec(f) = \{f(x)\colon x\in X\}$.
\item
The spectrum of a square matrix~$A$
from the $C^*$-algebra $M_n$
is the set of eigenvalues of~$A$.
\end{enumerate}%
\spacingfix%
\end{point}%
\begin{point}{210}[spectrum-basic]{Exercise}%
Let~$a$ be an element of a $C^*$-algebra $\scrA$.
\begin{enumerate}
\item
Prove that $\spec(a)\subseteq \R$ when $a$ is self-adjoint
(see~\sref{spectrum-self-adjoint-real}).

The reverse implication does not hold:
show that~$\spec(\bigl(
\begin{smallmatrix}0&2\\0&0\end{smallmatrix}\bigr))=\{0\}$.

\item
Show that $\spec(a^2)\subseteq [0,\infty)$ when $a$ is self-adjoint
(see~\sref{spectrum-self-adjoint-real}).

\item
Show that $|\lambda|\leq \|a\|$ for all~$\lambda\in\spec(a)$
using~\sref{spectrum-bounded}.

In fact, we will see in~\sref{norm-spectrum},
that $\|a\|=\sup\{\left|\lambda\right|\colon \lambda\in \spec(a)\}$.
\item
Show that $\spec(a)$ is closed (using~\sref{spectrum-bounded}).\\
Conclude that~$\spec(a)$ is compact.
\item
Show that $\spec(a+z)=\{\lambda+z\colon \lambda\in\spec(a)\}$
for all~$z\in \C$.
\item
Prove that~$\spec(a^{-1})=\{\lambda^{-1}\colon \lambda\in \spec(a)\}$
if~$a$ is invertible (and~$0\notin \spec(a)$).
\end{enumerate}%
\spacingfix{}%
\end{point}%
\begin{point}{220}%
On first sight,
the spectrum $\spec(a)$
of an element~$a$ of a $C^*$-algebra~$\scrA$ 
depends not only on~$a$,
but also on the surrounding $C^*$-algebra~$\scrA$ for it determines
for which~$\lambda\in\C$ the operator $a-\lambda$ is invertible.
Thus we should perhaps write $\spec_\scrA(a)$ instead
of~$\spec(a)$.
However, such careful bookkeeping turns out 
be unnecessary
by the following result.
\end{point}
\begin{point}{230}[spectral-permanence]{Theorem (Spectral Permanence)}%
\index{Spectral Permanence}%
Let~$\scrB$ be a $C^*$-subalgebra of a $C^*$-algebra $\scrA$.
Then~$\spec_{\scrA}(a)=\spec_\scrB(a)$
for every element~$a$ of~$\scrB$.
\begin{point}{240}{Proof}%
Let~$a$ be an element of~$\scrB$,
and let~$\lambda\in \C$.
We must show that $a-\lambda$ is invertible in~$\scrA$
iff $a-\lambda$ is invertible in~$\scrB$.
Surely,
if $a-\lambda$ has an inverse $(a-\lambda)^{-1}$ in~$\scrB$,
then~$(a-\lambda)^{-1}$ is also an inverse of~$a-\lambda$ in~$\scrA$,
since~$\scrB\subseteq \scrA$.
The other, non-trivial, direction follows
    directly from~\sref{inverse-permanence} 
    (and~\sref{improved-inverse-permanence}.)\qed%
\end{point}
\end{point}
\end{parsec}
\section{Positive Elements}
\subsection{Holomorphic Functions}
\begin{parsec}{120}%
\begin{point}{10}%
The next order of business
is to show that the spectrum~$\spec(a)$ of an element~$a$
of a $C^*$-algebra contains enough points, so to speak.
One incarnation of this idea 
is that~$\spec(a)$ is non-empty
(see~\sref{spectrum-non-empty}), but
we will need more,
and prove that  $\|a\|=\left|\lambda\right|$
for some~$\lambda\in\spec(a)$
(provided that~$a$ is self-adjoint).
Somewhat bafflingly,
the canonical and apparently
easiest way to derive this fact is by considering the power series
expansion of a cleverly chosen $\scrA$-valued function
(see~\sref{norm-spectrum}).
To this end,
we'll first quickly redevelop some complex analysis
for~$\scrA$-valued functions
(instead of $\C$-valued functions),
which will only be needed to prove this fact.
\end{point}
\begin{point}{20}{Setting}%
Fix a $C^*$-algebra~$\scrA$ for the remainder of this paragraph.
For brevity,
we'll say that a \Define{function}%
\index{function!$\scrA$-valued \& partial}
is a partially defined map $f\colon \C\to \mathscr{A}$
whose domain of definition $\dom(f)$%
\index{dom@$\dom(f)$, domain of an $\scrA$-valued partial function}
is an open subset of~$\C$.
Such a function is called \Define{holomorphic} at a point~$x\in \C$%
\index{function!holomorphic (at~$x$)}%
\index{holomorphic function}
if $f$ is defined on~$x$,
and 
\begin{equation*}
\frac{f(x)-f(y)}{x-y}
\end{equation*}
converges (with respect to the norm on~$\scrA$)
to some element~$f'(x)$ of~$\scrA$
as $y\in \dom(f)\backslash\{x\}$
converges to~$x$.

We say that~$f$ is \Define{holomorphic}
if~$f$ is holomorphic at~$x$ for all~$x\in \dom(f)$,
and the function $z\mapsto \Define{f'}(z)$
with $\dom(f')=\dom(f)$
is called its \Define{derivative}.%
\index{derivative of a holomorphic function}%
\index{$(\,\cdot\,)'$!$f'$, derivative of a holomorphic function}
\end{point}
\begin{point}{30}{Exercise}%
Verify the following examples of holomorphic functions.
\begin{enumerate}
\item
If~$f$ and $g$ are holomorphic functions with $\dom(f)=\dom(g)$,
then $f+g$ and $f\cdot g$ are holomorphic,
and $(f+g)'=f'+g'$ and $(f\cdot g)' = f'g+g'f$.

\item
The function~$f$ given by $f(z)=z$ and~$\dom(f)=\C$
is holomorphic, and $f'(z)=1$ for all $z\in\C$.

\item
Let~$a\in \scrA$. The constant function $f$ given by $f(z)=a$
for all~$z\in \C$ is holomorphic, and $f'(z)=0$ for all~$z\in \scrA$.

\item
Any polynomial,
that is, function~$f$ of the form $f(z)\equiv a_n z^n+\dotsb+a_1 z+a_0$
with~$a_i\in \scrA$
is holomorphic with $f'(z)=na_nz^{n-1}+\dotsb+2a_2z+a_1$.
\end{enumerate}%
\spacingfix{}%
\end{point}%
\end{parsec}%
\begin{parsec}{130}%
\begin{point}{10}%
We now turn
to perhaps the most important example
of a holomorphic $\scrA$-valued function ---
or at the very least the very source from
which (as we'll see) all holomorphic functions
draw their interesting and pleasant
properties:
the holomorphic $\scrA$-valued function
given by a power series  $\sum_n a_n z^n$.
\end{point}
\begin{point}{20}[hadamard]{Theorem}%
\index{power series}%
Let~$a_0,a_1,a_2,\dotsc\in\scrA$
be given,
and write~$R:=(\limsup_n \|a_n\|^{\nicefrac{1}{n}})^{-1}$.
Then for every~$z\in\C$,
\begin{enumerate}
\item
$\sum_n a_n z^n$
converges absolutely
when~$\left|z\right| < R$, and 
\item
if~$\sum_n a_n z^n$ converges,
then~$\left|z\right|\leq R$.
\end{enumerate}
(The number~$R\in[0,\infty]$ is called the \Define{radius of convergence}%
\index{radius of convergence}
of the series $\sum_n a_n z^n$.)
\begin{point}{30}{Proof}%
Suppose that $\left|z\right|<R$.
To show that the series 
$\sum_n a_nz^n$ converges absolutely,
we must show that $\sum_n \left\|a_n\right\|\left|z\right|^n
\equiv \sum_n (\,\left\|a_n\right\|^{\nicefrac{1}{n}}\left|z\right|\,)^n
<\infty$.
If~$z=0$, this is obvious,
so we'll assume that~$\left|z\right| > 0$.
Then, since~$\left|z\right|<R$,
we have~$R^{-1}\left|z\right|<1$
(and $R^{-1}<\infty$).
Note that
there is $\varepsilon>0$ with $(R^{-1}+\varepsilon)\left|z\right|<1$.
The point of this~$\varepsilon$
is that~$\limsup_n \|a_n\|^{\nicefrac{1}{n}} 
< R^{-1}+\varepsilon$,
so that we can find~$N$ 
with $\|a_n\|^{\nicefrac{1}{n}} \leq R^{-1}+\varepsilon$
for all~$n\geq N$.
Then $\|a_n\|^{\nicefrac{1}{n}}\left|z\right|
\leq (R^{-1}+\varepsilon)\left|z\right|<1$
for all~$n\geq N$,
and so
$\sum_n \|a_n\|\left|z\right|^n 
\leq\sum_{n=0}^{N-1} \|a_n\|\left|z\right|^n+ \sum_{n=N}^\infty 
(\,(R^{-1}+\varepsilon)\left|z\right|\,)^n < \infty$
by  convergence
of the geometric series (c.f.~\sref{geometric}).

Suppose now instead that $\sum_n a_n z^n$ converges.
Then~$\|a_n\|\left|z\right|^n$
converges to~$0$.
In particular,
there is~$N$ with $\|a_n\|\left|z\right|^n \leq  1$
for all~$n\geq N$.
Then~$\|a_n\|^{\nicefrac{1}{n}} \left|z\right| \leq 1$,
and $\|a_n\|^{\nicefrac{1}{n}} \leq \left|z\right|^{-1}$
for all~$n\geq N$,
so that $R^{-1}\equiv \limsup_n \|a_n\|^{\nicefrac{1}{n}}
\leq \left|z\right|^{-1}$,
giving $\left|z\right|\leq R$.\qed
\end{point}
\end{point}
\begin{point}{40}{Proposition}%
The $\scrA$-valued function~$f$
given by a series $\sum_n a_n z^n$
with radius of convergence~$\smash{R:=(\,\limsup_n \|a_n\|^{\nicefrac{1}{n}}\,)^{-1}}$
is holomorphic
when defined
on the disk $\dom(f)=\{z\in\C\colon \left|z\right|<R\}$,
and $f'(z)=\sum_{n=1}^\infty n a_n z^{n-1}$
for all~$z\in \dom(f)$.
\begin{point}{50}{Proof}%
If~$R=0$,
the statement is rather dull, but clearly true,
so we assume that~$R\neq 0$,
that is, $\smash{\limsup_n \|a_n\|^{\nicefrac{1}{n}}}<\infty$.

Note
that the radius of convergence
of~$\sum_{n=1}^\infty na_n z^{n-1}
\equiv \sum_{n=0}^\infty (n+1)a_{n+1}z^n$
is also~$R$,
because
\begin{equation*}
	\smash{\bigl\|\,(n+1)\, a_{n+1}\,\bigr\|^{\nicefrac{1}{n}}}
\ =\  \smash{(n+1)^{\nicefrac{1}{n}}}
\ \smash{\|a_{n+1}\|^{\frac{1}{n+1}}}
\ \smash{\bigl(\|a_{n+1}\|^{\frac{1}{n+1}}\bigr)^{\nicefrac{1}{n}}},
\end{equation*}
and
$R^{-1}=\limsup_n \|a_{n+1}\|^{\frac{1}{n+1}}$,
and both 
$(n+1)^{\nicefrac{1}{n}}$
and  $\bigl(\|a_{n+1}\|^{\frac{1}{n+1}}\bigr)^{\nicefrac{1}{n}}$
converge to~$1$ as~$n\to\infty$
(using here that $\smash{\limsup_n\|a_n\|^{\nicefrac{1}{n}}}<\infty$).

Hence~$\sum_{n=1}^\infty n a_n z^{n-1}$
converges absolutely for every~$z\in\C$ with $\left|z\right|<R$.
Let~$z\in \C$ with $\left|z\right|<R$
be given. We must show that~$f$
is holomorphic at~$z$ with~$f'(z)=\sum_n na_nz^{n-1}$.
For this it suffices to show that
\begin{equation}
\label{power-series-derivative-0}
\sum_{n=0}^\infty
\|a_n\|\left|\,\frac{(z+h)^n-z^n}{h}-nz^{n-1}\,\right|
\end{equation}
converges to~$0$
as $h\in\C$ (with $h\neq 0$ and $\left|z+h\right|<R$)
tends to~$0$.

Pick~$r>0$ with $\left|z\right| < r < R$.
With the appropriate algebraic gymnastics
(involving the identity
$a^n-b^n=(a-b)\sum_{k=1}^n a^{n-k}b^{k-1}$
and the inequalities $\left|z+h\right| \leq r$
and $\left|z\right|\leq r$)
we get, for every~$n$
and~$h\in \C$ with
$h\neq 0$ and~$\left|z+h\right|<r$,
\begin{align}
\left|\  \frac{(z+h)^n-z^n}{h}-nz^{n-1}\ \right|
\ &= \ 
\biggl|\ \sum_{k=1}^n \bigl(\,(z+h)^{n-k}-z^{n-k}\,\bigr)z^{k-1} \ \biggr|
\label{power-series-derivative-1}
\\ 
\label{power-series-derivative-2}
\ &\leq\ 
2nr^{n-1}.
\end{align}
On the one hand,
we see from~\eqref{power-series-derivative-1}
that any term
--- and thus any partial sum ---
of the series from~\eqref{power-series-derivative-0}
converges to~$0$ as~$h$ tends to~$0$.
On the other hand,
we see from~\eqref{power-series-derivative-2}
that the series 
from~\eqref{power-series-derivative-0}
is dominated by~$2\sum_n \|a_n\|nr^{n-1}$
(which converges
because the radius of convergence of $\sum_n a_n nz^{n-1}$
is $R>r$),
so that the tails of the series in~\eqref{power-series-derivative-0}
vanish uniformly in~$h$.
All in all, the sum of the infinite series
from~\eqref{power-series-derivative-0}
converges to~$0$ as~$h$ tends to~$0$.\qed
\end{point}
\end{point}
\begin{point}{60}[powerseries-uniqueness-coeffients]{Exercise}%
Let~$\sum_n a_n z^n$
be a power series over~$\scrA$
with radius of convergence~$R>0$
such that~$\sum_n a_n z^n=0$
for all~$z$ from
some disk around~$0$ with radius~$r<R$.
Show that~$0=a_0=a_1=a_2=\dotsb$.

(Hint: clearly~$a_0=0$.  Show that the derivative
of the power series also vanishes on the disk around~$0$
with radius~$r$.)
\end{point}
\end{parsec}%
\begin{parsec}{140}%
\begin{point}{10}%
All holomorphic functions
are power series
in the sense
that any $\scrA$-valued holomorphic
function~$f$ defined on~$0$
is given by some power series $\sum_n a_n z^n$
on the largest disk around~$0$
that fits in~$\dom(f)$.
This fact,
which follows from~\sref{taylor}
and~\sref{rigid-expansion} below,
is all the more remarkable,
because here the pointwise (``local'') property
of being holomorphic
entails
the uniform (``global'') property
of being equal to a power
series (on some disk).
The device
that bridges
this
gap
is integration of $\scrA$-valued
holomorphic functions along
line segments.
\end{point}
\begin{point}{20}{Exercise}%
We're going to define as quickly as possible
an integral~$\int f $
for every continuous map $f\colon [0,1]\to\scrA$.%
\index{S@$\int$, integral!$\int f$, of continuous $f\colon [0,1]\to \scrA$}
Any interval~$I$
in~$[0,1]$
is of one of the following forms
\begin{equation*}
	[s,t]\qquad[s,t)\qquad(s,t]\qquad(s,t)
\end{equation*}
where~$0\leq s\leq t\leq 1$;
we'll denote the length of an interval~$I$ --- being~$t-s$ 
in the four cases above --- by $\Define{\left|I\right|}$.
An \Define{$\scrA$-valued step function}
is a function $f\colon [0,1]\to\scrA$
of the form
$f\equiv \sum_n a_n \mathbf{1}_{I_n}$
for some~$a_1,\dotsc,a_N\in\scrA$
and intervals $I_1,\dotsc,I_N$
(where~$\mathbf{1}_{I_n}$ is~$1$ is the 
\emph{indicator function}%
\index{$\mathbf{1}_A$, indicator function!on $[0,1]$} of~$I_n$
which is~$1$
on~$I_n$
and~$0$ elsewhere);
and the set of $\scrA$-valued step functions
is denoted by~$\Define{S_\scrA}$,%
\index{SA@$S_\scrA$}
which is a subset
of the space of all bounded functions
$f\colon [0,1]\to\scrA$
which we'll denote by~$\Define{B_\scrA}$.%
    \index{BA@$B_\scrA$}
\begin{enumerate}
\item
Show that there is a unique
linear map $\int\colon S_\scrA\to \scrA$
with~$\int a \mathbf{1}_{I}=\left| I \right|a$
for every interval~$I$ in~$[0,1]$
and~$a\in\scrA$.

(Hint:  the difficulty
here is to show that no contradiction
arises in the sense that 
$\sum_n a_n\left|I_n\right| = \sum_m a_m' \left|I_m'\right|$
when 
$\sum_n a_n \mathbf{1}_{I_n}=\sum_m a_m' \mathbf{1}_{I_n'}$
for intervals $I_1,\dotsc,I_N,I_1',\dotsc,I_M'$ in~$[0,1]$
and $a_1,\dotsc,a_N,a_1',\dotsc,a_M'\in\scrA$.)

\item
We endow~$B_\scrA$
with the supremum norm,
viz.~$\|f\|=\sup_{t\in[0,1]} \|f(t)\|$
for all~$f\in B_\scrA$.

Show that every $\scrA$-valued step function~$f$
may be written as
$f\equiv \sum_n a_n \mathbf{1}_{I_n}$
where $I_1,\dotsc,I_N$
are \emph{disjoint and non-empty}
intervals in~$[0,1]$.

Show that 
for such a representation
$\|f\|=\sup_n \|a_n\|$, and
$\sum_n \left|I_n\right|\leq 1$.
Deduce that
$\|\int f\| \leq \sum_n \|a_n\|\left|I_n\right|
\leq \|f\|$.

Conclude that~$\int\colon S_\scrA\to\scrA$
is a bounded linear map
and can therefore
be uniquely extended to a bounded linear map
$\int\colon \overline{S}_\scrA\to\scrA$
on the closure~$\overline{S}_\scrA$
of~$S_\scrA$.

\item
Show that every continuous function $f\colon [0,1]\to\scrA$
is the supremum norm limit
of a sequence $g_1,g_2,\dotsc$
of $\scrA$-valued step functions
 (i.e.~$f\in\overline{S}_\scrA$).

\item
Show that~$\int af = a\int f$
when~$f\colon [0,1]\to\C$
is continuous and~$a\in \scrA$.
\end{enumerate}%
\spacingfix{}%
\end{point}%
\begin{point}{30}{Definition}%
The integral of
a holomorphic $\scrA$-valued function~$f$
along a line segment $[w,w']\subseteq\dom(f)$
	(where~$w$ and~$w'$ are thus \emph{complex numbers})
is now defined as
\begin{equation*}
\Define{\int_{w}^{w'}f}
\ = \ 
(w'-w)\int_0^1f(\,w+t(w'-w)\,)\,dt.
\end{equation*}
\index{S@$\int$, integral!of continuous $f\colon \C\to\scrA$!$\int_w^{w'} f$, over an interval}%
We'll also need integration along a triangle~$T$,%
\index{triangle, for our purposes}
which is for this  purpose a triple of 
complex numbers~$w_0,w_1,w_2$
(of which the order \emph{does} matter)
called the
\Define{vertices} of~$T$.
The \Define{boundary} of
such a triangle~$T$
is $\Define{\partial T}
:=[w_0,w_1]\cup [w_1,w_2]\cup [w_2,w_0]$,
and given any $\scrA$-valued
holomorphic function~$f$
with $\partial T\subseteq \dom(f)$
we define
\begin{equation*}
	\Define{\int_T f}\ = \ \int_{w_0}^{w_1} f
\,+\, \int_{w_1}^{w_2} f
\,+\, \int_{w_2}^{w_0} f.
\end{equation*}
\index{S@$\int$, integral!of continuous $f\colon \C\to\scrA$!$\int_T f$, over a triangle}%
We'll need some more terminology
relating to our triangle~$T$.
Its \Define{closure},
written $\Define{\mathrm{cl}(T)}$,
is the convex hull of~$w_0,w_1,w_2$,
and its \Define{interior}
is simply
$\Define{\mathrm{in}(T)}=\mathrm{cl}(T)\backslash \partial T$.
The length
of~$T$ is given by
$\Define{\mathrm{length}(T)}:=
\left|w_1-w_0\right|\,+\,
\left|w_2-w_1\right|\,+\,
\left|w_0-w_2\right|$.

The number of times the triangle~$T$ winds
around a point $z\in \C\backslash \partial T$
in the counterclockwise direction
is
called the 
\Define{winding number}, and
is
written $\Define{\mathrm{wn}_T(z)}$,%
\index{wn@$\mathrm{wn}_T$, winding number}
is either $1$ or~$-1$ when~$z\in\mathrm{in}(T)$
(depending on the order of the vertices),
is~$0$ when~$z\notin\mathrm{cl}(T)$,
and undefined on~$\partial T$.
It is defined formally for $z\in \C\backslash \partial T$ by
\begin{equation*}
	\textstyle
2\pi \wn_T(z)\ = \ 
\measuredangle(w_0,z,w_1)
\,+\,\measuredangle(w_1,z,w_2)
\,+\, \measuredangle(w_2,z,w_0),
\end{equation*}
where~$\Define{\measuredangle(w_0,z,w_1)}$%
\index{$\measuredangle(w_0,z,w_1)$, angle between complex numbers}
denotes the number of radians in~$(-\pi,\pi)$
needed to rotate the line through~$z$ and~$w_0$ counterclockwise
around~$z$ to hit~$w_1$,
that is, the angle of the corner on the right when travelling
from~$w_0$ to~$w_1$ via~$z$.

The winding number~$\mathrm{wn}_T$
pops
up in the value of the integral $\int_T (z-z_0)^{-1}dz$
later on,
see~\sref{invint}.
\end{point}
\begin{point}{40}[goursat]{Goursat's Theorem}%
\index{Goursat's Theorem}%
Let~$f$ be a holomorphic function,
and let~$T$ be a triangle whose closure
is entirely contained in~$\dom(f)$.
Then~$\int_T f = 0$.
\begin{point}{50}[goursat-1]{Proof}%
(Based on~\cite{moore1900}.)
If two vertices of~$T$ coincide
the result is obviously true,
so we may assume that they're all distinct,
that is, $\mathrm{in}(T)\neq \varnothing$.

Note that if~$f$ has an antiderivative,
that is, $f\equiv g'$ for some holomorphic function~$g$,
then one can show that~$\int_T f=0$
(after deriving the fundamental theorem of calculus).
Although it is true that every holomorphic function
with simply connected domain has a antiderivative,
this result is not yet available 
(and in fact usually depends on this very theorem).
Instead we will approximate~$f$
by an affine function
(which does have an antiderivative)
using the derivative of~$f$.
But since such an approximation only
concerns a single point,
we first need to zoom in.
\begin{point}{60}[goursat-2]%
If we split~$T$ into four similar triangles
$T^\text{i}$, $T^\text{ii}$,
$T^\text{iii}$, $T^\text{iv}$
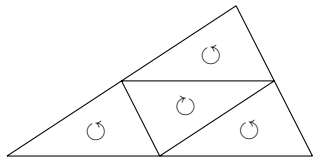
\begin{equation*}
\begin{tikzpicture}
\coordinate (A) at (0,0){};
\coordinate (B) at (4,0){};
\coordinate (C) at (3,2){};
\coordinate (AB) at ($(A)!0.5!(B)$) {};
\coordinate (BC) at ($(B)!0.5!(C)$) {};
\coordinate (CA) at ($(C)!0.5!(A)$) {};
\node  at ($(AB)!0.5!(BC)!0.33!(CA)$) {$\circlearrowright$};
\node  at ($(A)!0.5!(AB)!0.33!(CA)$) {$\circlearrowleft$};
\node  at ($(B)!0.5!(BC)!0.33!(AB)$) {$\circlearrowleft$};
\node  at ($(C)!0.5!(BC)!0.33!(CA)$) {$\circlearrowleft$};
\draw
	(A) -- (AB) 
	(AB) -- (B)
	(B)  -- (BC)
	(BC) -- (C) 
	(C) -- (CA)
	(CA) -- (A)
	(CA) -- (BC)
	(BC) -- (AB)
	(AB) -- (CA);
\end{tikzpicture}
\end{equation*}
we have $\smash{\int_Tf = \sum_{n={\text{i}}}^{\text{iv}} \int_{T^n}f}$.
There is $T'$ among
$T^\text{i}$, $T^\text{ii}$,
$T^\text{iii}$, $T^\text{iv}$
with 
 $\|\int_Tf\|\leq 4 \|\int_{T'} f\|$.
Clearly, $\length(T)=2\length(T')$.
Write~$T_0 := T$ and $T_1 := T'$. 

From this it is clear how to
 get a sequence of similar triangles $T_0, T_1, T_2, \dotsc$
with $\|\int_Tf\|\leq 4^n \|\int_{T_n} f\|$,
and $\length(T)=2^n\length(T_n)$.
\end{point}
\begin{point}{70}%
If we pick a point on the closure $\mathrm{cl}(T_n)$
of each triangle~$T_n$ 
we get a Cauchy sequence
that converges to some point~$z_0\in\C$
which lies in~$\bigcap_n \mathrm{cl}(T_n)$.
We can approximate $f$ by an affine
function at~$z_0$ as follows.
For $z\in \dom(f)$,
\begin{equation*}
f(z)\ = \ f(z_0)\,+\,f'(z_0)\,(z-z_0)\,-\,r(z)\,(z-z_0),
\end{equation*}
where~$r\colon \dom(f)\to \C$
is given by $r(z)=f'(z_0)-(f(z)-f(z_0))(z-z_0)^{-1}$ for $z\neq z_0$
and $r(z_0)=0$.
We see that~$r(z)$ converges to~$0$ as~$z\to z_0$.

Let~$\varepsilon >0$ be given.
There is~$\delta>0$
such that $z\in\dom(f)$
and $\|r(z)\|\leq \varepsilon$
for all~$z\in \C$ with $\|z-z_0\|<\delta$.
There is~$n$ such that the triangle~$T_n$ is contained
in the ball around~$z_0$ of radius~$\delta$.
Note that $\int_{T_n} f(z_0)+f'(z_0)(z-z_0)\,dz=0$
by the discussion in~\sref{goursat-1}, because
the integrated function is affine.
Thus
\begin{equation*}
\textstyle
\int_{T_n} f \  = \ -\int_{T_n}r(z)\,(z-z_0)\,dz.
\end{equation*}
Note that for $z\in T_n$,
we have  $\|z-z_0\|\leq \length(T_n)$,
and $\|r(z)\|\leq \varepsilon$ (because $\|z-z_0\|< \delta$),
and so $\|r(z)(z-z_0)\|\leq \varepsilon\,\length(T_n)$.
Thus:
\begin{equation*}
\textstyle
\|\int_{T_n} f\| \  = \ \|\int_{T_n}r(z)\,(z-z_0)\,dz\|
\ \leq\ \varepsilon\length(T_n)^2.
\end{equation*}
Using the inequalities from~\sref{goursat-2},
we get
\begin{equation*}
\textstyle
\|\int_T f\|\ \leq\ 4^n\, \|\int_{T_n} f\|
\ \leq\ \varepsilon \,4^n\,\length(T_n)^2 
\ \equiv\ \varepsilon \length(T)^2.
\end{equation*}
Since~$\varepsilon>0$ was arbitrary,
we see that~$\int_T f=0$.\qed
\end{point}
\end{point}
\end{point}%
\begin{point}{80}[invint]{Exercise}%
	The assumption in Goursat's Theorem (\sref{goursat})
that the holomorphic function~$f$
is defined not only on the boundary~$\partial T$
of the triangle~$T$
but also on the interior $\mathrm{in}(T)$
is essential,
for if only a single hole in~$\dom(f)$ is allowed within $\mathrm{in}(T)$
the integral~$\int_T f$ can become non-zero---which we will demonstrate
here by computing
$\int_T (z-z_0)^{-1}dz$.
\begin{enumerate}%
\item
Show that for a non-zero
complex number~$z$ we have
\begin{equation*}
	z^{-1}\ =\  \frac{\Real{z}-i\Imag{z}}{\Real{z}^2+\Imag{z}^2}.
\end{equation*}
\item
Given real numbers~$a\neq 0$ and~$b$,
show that
\begin{alignat*}{3}
	\int_{a}^{a+ib}\ 
	z^{-1}\,dz
	\ &=\ 
	i\int_{0}^t \frac{a-it}{a^2+t^2}dt
	\\
	\ &=\ 
	i\int_0^b \frac{a}{a^2+t^2}\,dt 
	\ +\ 
	\int_0^b \frac{t}{a^2+t^2}\,dt\\
	\ &=\ 
	\textstyle
	i\,\arctan(\,b/a\,)
	\,+\, \log\left|a+ib\right| - \log\left|ia\right|,
\end{alignat*}
and similarly, show that for real numbers~$a$ and $b\neq 0$
\begin{equation*}
	\int_{a+ib}^{ib} z^{-1}\,dz
	\ = \ 
	i\arctan(\,a/b\,) \ +\ 
	\log\left|ib\right| \,-\,
	\log\left|a+ib\right|.
\end{equation*}
	\item
Show that for complex numbers~$w$, $w'$ and~$z_0$
with~$z_0\notin [w,w']$
\begin{equation*}
	\int_{w}^{w'}\,(z-z_0)^{-1}\,dz
\ = \ 
i\, \measuredangle(w,z_0,w')\ +\ 
\log\,\frac{\left|w'-z_0\right|}{\left|w-z_0\right|},
\end{equation*}
where~$\measuredangle(w,z_0,w')$
denotes
the number of radians
in~$(-\pi,\pi)$
needed
to rotate the line through~$z_0$ and~$w$
counterclockwise around~$z_0$ to hit~$w'$.

(Hint: 
using Goursat's Theorem, \sref{goursat},
one may reduce the problem
to integration along horizontal and vertical line segments.)
\item
Given a triangle~$T$ and~$z_0\in\C\backslash \partial T$,
show that
\begin{equation*}
	\frac{1}{2\pi i}\int_T (z-z_0)^{-1}\,dz
	\ =  \mathrm{wn}_T(z_0).
\end{equation*}%
\end{enumerate}%
\spacingfix{}%
\end{point}%
\begin{point}{90}%
Thus integration of 
$z\mapsto (z-z_0)^{-1}$
along a triangle~$T$ detects
the number of times~$T$ winds
around~$z_0$.
There is nothing special about a triangle:
a similar result---not needed here---holds
for a broad class of curves
(c.f.~Thm~2.9 of~\cite{conway2013}).

Integration along a curve can also be used
to probe the value of a holomorphic function at a point~$z_0$.
On this occasion
we restrict ourselves
to regular $N$-gons.
\end{point}
\end{parsec}%
\begin{parsec}{150}%
\begin{point}{10}[cauchy-formula]{Theorem (Cauchy's Integral Formula)}%
\index{Cauchy's Integral Formula}%
Let~$f$ be a holomorphic $\scrA$-valued function
which is defined on the interior and boundary
of some regular $N$-gon 
with centre~$c\in\C$,
circumradius~$r$
and vertices $w_n := c+r\cos(2\pi/n)+ir\sin(2\pi/n)$.
Then for any complex number~$z_0$ in
the interior of the~$N$-gon
we have
\begin{equation*}
	f(z_0)\ = \ \frac{1}{2\pi i}\,\sum_{n=0}^{N-1}\int_{w_n}^{w_{n+1}}
	\frac{f(z)}{z-z_0}\,dz
\end{equation*}
\begin{point}{20}{Proof}%
	Since~$\sum_{n=0}^{N-1} \int_{w_n}^{w_{n+1}} \frac{f(z_0)}{z-z_0}\,dz
	= 2\pi i f(z_0)$ by~\sref{invint}
it suffices to show that
\begin{equation}
\label{eq:cauchy-formula-1}
\sum_{n=0}^{N-1}\int_{w_n}^{w_{n+1}} \frac{f(z)-f(z_0)}{z-z_0}\,dz \ = \ 0.
\end{equation}
\begin{point}{30}[cauchy-formula-1]%
Let~$\varepsilon>0$ be given.
Since~$f$ is holomorphic at~$z_0$
we can find $\delta>0$ with
\begin{equation*}
\left\|\frac{f(z)-f(z_0)}{z-z_0}\right\|
\ \leq \ \,\|f'(z_0)\|\,+\,37
\end{equation*}
for all~$z\in\dom(f)$ with $\left|z-z_0\right|\leq \delta$. 
\end{point}
\begin{point}{40}%
To use~\sref{cauchy-formula-1},
we must restrict our attention to a smaller polygon.
Let~$T$ be a triangle 
that is entirely inside the~$N$-gon
such that  $\mathrm{wn}_T(z_0)=-1$,
$\length(T)\leq \varepsilon$,
and $\|z_0-z\|\leq \delta$ for all~$z\in \partial T$.
By partitioning the area
between~$T$ and the~$N$-gon
in the obvious manner
into triangles~$T_1,\dotsc,T_M$
(for which~$\int_{T_m}f=0$ for all~$m$ by~\sref{goursat})
we see that
\begin{equation}
\label{eq:cauchy-formula-2}
\sum_{n=0}^{N-1}\int_{w_n}^{w_{n+1}} \frac{f(z)-f(z_0)}{z-z_0}\,dz
\ = \ 
\int_T \frac{f(z)-f(z_0)}{z-z_0}\,dz.
\end{equation}
Hence by~\sref{cauchy-formula-1}
we have
\begin{alignat*}{3}
	\left\|\,\sum_{n=0}^{N-1} 
	\int_{w_n}^{w_{n+1}} \frac{f(z)-f(z_0)}{z-z_0}\,dz\,\right\|
	\ &\leq \ \length(T)\,\cdot\,
\sup_{z\in\partial T} \,\left\|\,\frac{f(z)-f(z_0)}{z-z_0}\,\right\|
\\
	\ &\leq \ \|f'(z_0)\|\varepsilon\,+\,37\varepsilon.
\end{alignat*}
Since~$\varepsilon>0$ was arbitrary,
 \eqref{eq:cauchy-formula-1}
follows from \eqref{eq:cauchy-formula-2}.\qed
\end{point}
\end{point}
\end{point}
\begin{point}{50}[taylor]{Proposition}%
Let~$f$ be a holomorphic $\scrA$-valued function
defined on the boundary and interior
of a regular $K$-gon
with vertices $w_0,\dotsc,w_{K-1},w_K=w_0$
as in~\sref{cauchy-formula}.
Then for every element~$z$ of an open disk in the interior of the $K$-gon
with centre~$w$,
\begin{equation*}
f(z)\ = \ 
	\sum_{n=0}^\infty \ 
	\left(\frac{1}{2\pi i}\sum_{k=0}^{K-1}\int_{w_k}^{w_{k+1}} 
	\frac{f(u)}{(u-w)^{n+1}}\,du\right)
\ (z-w)^n.
\end{equation*} 
\begin{point}{60}{Proof}%
By~\sref{cauchy-formula} and some easy algebra we have
\begin{alignat*}{3}
2\pi if(z)\ &=\  
	\sum_{k=0}^{K-1}\int_{w_k}^{w_{k+1}}
	\frac{f(u)}{u-z}\,du
\ =\ 
\sum_{k=0}^{K-1}\int_{w_k}^{w_{k+1}}
   \frac{f(u)}{u-w}\,\frac{1}{1-\frac{z-w}{u-w}}\,du
\end{alignat*}
Note that~$\left|z-w\right|<\left|u-w\right|$
for all~$u\in [w_k,w_{k+1}]$ and~$k$,
because the open disk with centre~$w$
from which~$z$ came lies entirely in the~$K$-gon.
Hence,
by~\sref{geometric},
\begin{alignat*}{3}
	2\pi if(z) \ &= \ 
\sum_{k=0}^{K-1}\int_{w_k}^{w_{k+1}}
   \frac{f(u)}{u-w}\, \sum_{n=0}^\infty 
\frac{(z-w)^n}{(u-w)^n}
 \,du\\
 \ &= \ 
  \sum_{n=0}^\infty \ 
\sum_{k=0}^{K-1}\int_{w_k}^{w_{k+1}}
  \frac{f(u)}{(u-w)^{n+1}}du \ (z-w)^n,
\end{alignat*}
where the interchange of ``$\sum$'' and ``$\int$''
is allowed
because the partial sum 
$\sum_{n=0}^Nf(u)\frac{(z-w)^n}{(u-w)^{n+1}}$
converges uniformly in~$u$ as~$N\to\infty$.\qed
\end{point}
\end{point}
\begin{point}{70}[rigid-expansion]{Proposition}%
Let~$f$ be an $\scrA$-valued holomorphic
function
that 
can be written as a power
series $f(z)=\sum_n a_n (z-w)^n$
where~$a_0,a_1,\dotsc\in\scrA$
for all~$z$ from some disk in~$\dom(f)$ around~$w$
with radius~$r>0$.

Then the formula $f(z)=\sum_n a_n (z-w)^n$
holds also for any $z$ from a larger disk 
with radius~$R>r$
around~$w$ that still fits in~$\dom(f)$.
\begin{point}{80}{Proof}%
Let~$z$ with $\left|z-w\right|<R$ be given.
By choosing~$K$ large enough
we can fit the boundary of a regular $K$-gon
centred around~$w$
with vertices $w_0,\dotsc,w_{K-1},w_{K}\equiv w_0$
inside the difference between the two disks,
and we can moreover, by~\sref{taylor},
choose the polygon
in such a way that
$f(z')=\sum_n b_n (z'-w)^n$
for all~$z'\in\C$
with $\left|z'-w\right|\leq \left|z-w\right|$
where $b_n = \sum_{k=0}^{K-1}
\int_{w_k}^{w_{k+1}}
\frac{f(u)}{(u-w)^{n+1}}\,du$.

Thus to show that $f(z)=\sum_n a_n (z-w)^n$
it suffices to show that~$a_n=b_n$ for all~$n$.
This in turn
follows by~\sref{powerseries-uniqueness-coeffients} from the fact
that $\sum_n a_n (z'-w)^n
= \sum_n b_n (z'-w)^n$
for all~$z'\in \C$ with $\left|z'-w\right|<r$.\qed
\end{point}
\end{point}
\end{parsec}
\subsection{Spectral Radius}
\begin{parsec}{160}%
\begin{point}{10}%
Our analysis of $\scrA$-valued
holomorphic functions
allows us to expose
the following connection
between the norm
and the invertible elements
in a $C^*$-algebra.
\end{point}
\begin{point}{20}[norm-spectrum]{Proposition}%
For a self-adjoint element~$a$ of a $C^*$-algebra~$\scrA$,
we have
\begin{equation*}
\|a\|\,=\,\sup\{\,\left|\lambda\right|\colon 
\,\lambda\in \spec(a)\,\}.
\end{equation*}
(The quantity on the right hand-side above
is called the \Define{spectral radius} of~$a$.)%
\index{spectral radius}
\begin{point}{30}{Proof}%
Write~$r=
\sup\{\left|\,\lambda\right|\colon\, \lambda\in \spec(a)\backslash\{0\}\,\}$
where the supremum is computed
in~$[0,\infty]$ so that~$\sup\varnothing=0$.
Since~$\left|\lambda\right| \leq \|a\|$
for all~$\lambda\in\spec(a)$
(\sref{spectrum-bounded})
we see that~$r\leq \|a\|$,
and so we only need to show that~$\|a\|\leq r$. 
Note that this is clearly true if~$\|a\|=0$,
so we may assume that~$\|a\|\neq 0$.

The trick is to consider
the power series expansion
around~$0$ of the holomorphic function~$f$ defined
on~$G:=\{\,z\in \C\colon 1-az\text{ is invertible}\,\}$ 
by  $f(z)=z(1-az)^{-1}$.
More specifically,
we are interested in the distance~$R$
of~$0$ to the complement of~$G$,
viz.~$R= \inf\{\left|\lambda\right|\colon \lambda\in \C\backslash G\}$
(where the infimum is computed in~$[0,\infty]$
so that~$\inf\varnothing=\infty$)
because since $0\in G$
and $z\notin G\iff z^{-1}\in \spec(a)$,
we have~$R=r^{-1}$
(using the convention $0^{-1}=\infty$).

Note that $f$ has the power series expansion
$f(z) = \sum_n a^nz^{n+1}$
for all~$z\in \C$ with $\|z\|<\|a\|^{-1}$,
because for such~$z$
we have $\sum_n (az)^n=(1-az)^{-1}$
by~\sref{geometric},
and thus~$f(z)=z(1-az)^{-1}=z\sum_n (az)^n = \sum_n a^nz^{n+1}$.

By~\sref{rigid-expansion}
we know that~$f(z)=\sum_n a^nz^{n+1}$
is valid not only for~$z\in \C$ with $\left|z\right|< \|a\|^{-1}$,
but for all~$z$ with $\left|z\right|< R$.
However, $R$ cannot be strictly larger than~$\|a\|^{-1}$,
because for every $z\in\C$ with $\left|z\right|>\|a\|^{-1}$
the series $\sum_n(az)^n$ 
and thus $\sum_n a^n{z}^{n+1}$ diverges (see~\sref{geometric-convergence})
--- using here that~$a$ is self-adjoint.
Hence~$R=\|a\|^{-1}$, and so~$r=\|a\|$.\qed
\end{point}
\end{point}
\begin{point}{40}{Remark}%
For an arbitrary (possibly non-self-adjoint)
element~$a$ of a $C^*$-algebra~$\scrA$
the formula in~\sref{norm-spectrum}
might be incorrect, e.g.~$\bigl\|\,\bigl(\begin{smallmatrix}
0& 1\\
0 & 0
\end{smallmatrix}\bigr)\,\|=1$
while 
$\spec(\,\bigl(\begin{smallmatrix}
0& 1\\
0 & 0
\end{smallmatrix}\bigr)\,)=\{0\}$
cf.~\sref{geometric-non-self-adjoint}.
For such~$a$
the formula
$\sup\{\,\left|\lambda\right|\colon\, \lambda\in \spec(a)\,\}
\,=\, \limsup_n \|a^n\|^{\nicefrac{1}{n}}$
	can be derived (see e.g.~Theorem~3.3.3 of~\cite{kr}) --- 
	which we won't need here.
\end{point}
\begin{point}{50}[spectrum-non-empty]{Exercise}%
Given a self-adjoint element~$a$ of a $C^*$-algebra show that
$\spec(a)\neq \varnothing$.
\end{point}
\begin{point}{60}{Exercise}%
Given a self-adjoint element~$a$ of a $C^*$-algebra
and~$\lambda\in\R$
show that $\spec(a) =\{\lambda\}$ iff $a=\lambda$.
\end{point}
\begin{point}{61}{Exercise}%
Use the previous exercise to prove the following theorem.
\end{point}
\begin{point}{70}{Theorem (Gelfand--Mazur for $C^*$-algebras)}%
\index{Gelfand--Mazur's Theorem}%
If every non-zero element of a $C^*$-algebra~$\scrA$
is invertible, then~$\scrA=\C$ or~$\scrA=\{0\}$.
\end{point}
\begin{point}{80}[gelfand-mazur-predicament]{Remark}%
A logical next step
towards Gelfand's representation theorem
is to show that if~$\lambda\in\spec(a)$
for some element~$a$ of a \emph{commutative} $C^*$-algebra~$\scrA$,
then there is an miu-map $f\colon \scrA\to \C$
with~$f(a)=\lambda$.
Here we have moved ourselves into a tight spot
by evading Banach algebras,
because the mentioned result is usually obtained
by finding a maximal ideal~$I$ of~$\scrA$
(by Zorn's Lemma) that contains~$\lambda-a$,
and then forming the \emph{Banach algebra} quotient~$\scrA/I$.
One then applies Gelfand--Mazur's Theorem for \emph{Banach algebras}, 
to see that
$\scrA/I= \C$,
and thereby obtain an miu-map~$f\colon \scrA\to \C$ with~$f(a-\lambda)=0$.
The problem here is that while $\scrA/I$
will turn out to be a $C^*$-algebra (indeed, be $\C$)
the formation of the $C^*$-algebra quotient
is non-trivial and depends on Gelfand's representation theorem
(see e.g. \S{}VIII.4 of~\cite{conway2013}) 
which is the very theorem we are working towards.
The way out of this predicament
is to avoid ideals and quotients of $C^*$- and Banach algebras
altogether,
and instead work 
with order ideals (and what are essentially
 quotients of Riesz and order unit spaces).
To this end,
we develop the theory
of the positive elements of a $C^*$-algebra
farther than is usually done
for Gelfand's representation theorem.
\end{point}
\end{parsec}
\begin{parsec}{170}[cstar-positive-2]%
\begin{point}{10}%
We return to the positive elements 
in a $C^*$-algebra (see~\sref{cstar-positive-def}).
We'll see that the connection we have established
between the norm and invertible elements
of a $C^*$-algebra
via the spectral radius (\sref{norm-spectrum})
affects the positive elements as well, see~\sref{cstar-positive-1}.
\end{point}
\begin{point}{20}[real-pos-ineq]{Exercise}%
Show that 
$\left|\,\lambda-t\,\right| \,\leq\, t$ iff  $\lambda \in[0,2t]$,
where $\lambda,t\in\R$.
\end{point}
\begin{point}{30}[pos-spectrum]{Proposition}%
For a self-adjoint element $a$ from a $C^*$-algebra,
and $t\in [0,\infty]$, 
\begin{equation*}
\|a-t\|\,\leq\, t\qquad\iff\qquad \spec(a)\subseteq [0,2t].
\end{equation*}%
\spacingfix{}%
\begin{point}{40}{Proof}%
To begin, note that~$\spec(a-t)=\spec(a)-t\subseteq \R$ 
by~\sref{spectrum-basic},
because~$a$ is self-adjoint.
Thus $\|a-t\|=\sup\{\,\left|\lambda-t\right|\colon \lambda\in \spec(a)\,\}$
by~\sref{norm-spectrum}.
Hence $\|a-t\|\leq t$
iff $\left|\lambda-t\right|\leq t$ for all~$\lambda\in\spec(a)$
iff $\spec(a)\subseteq [0,2t]$ (by \sref{real-pos-ineq}).\qed
\end{point}
\end{point}
\begin{point}{50}[cstar-positive-1]{Exercise}%
\index{positive!element of a $C^*$-algebra}%
Show
(using~\sref{pos-spectrum} and~\sref{spectrum-basic})
that
for any self-adjoint element $a$ of a $C^*$-algebra~$\scrA$,
the following are equivalent.
\begin{enumerate}
\item 
\label{cstar-pos-1}
$\|a-t\|\leq t$
for some $t\geq \frac{1}{2}\|a\|$;
\item 
\label{cstar-pos-2}
$\|a-t\|\leq t$
for all $t\geq \frac{1}{2}\|a\|$;
\item 
\label{cstar-pos-3}
$\spec(a)\subseteq[0,\infty)$;
\item
$a$ is positive.
\end{enumerate}
We will complete this list in~\sref{cstar-positive-final}.
\end{point}
\begin{point}{60}[positive-basic-2]{Exercise}%
Let~$\scrA$ be a $C^*$-algebra.
\begin{enumerate}
\item
Show that $0\leq a\leq 0$ entails that~$a=0$
for all~$a\in\scrA$.
\item
Show that~$\pos{\scrA}$ is closed.
\item
Let~$a$ be a self-adjoint element of~$\scrA$.
Show that
 $-\lambda \leq a\leq \lambda$
iff $\|a\|\leq \lambda$,
for $\lambda\in [0,\infty)$.
Conclude that $\|a\| = \inf\{ \lambda \in \R\colon 
-\lambda \leq a\leq \lambda\}$.

(In other words
$\sa{\scrA}$ is a \emph{complete Archimedean order unit space},
see Definition~1.12 of~\cite{alfsen2012}---a 
type of structure first studied in~\cite{kadison1951}.)

Show that $0\leq a \leq b$ entails $\|a\|\leq \|b\|$
for $a,b\in\sa{\scrA}$.

\item 
Recall that $ab$ need not be positive if~$a,b\geq 0$. However:

Show that $a^2$ is positive for every self-adjoint element~$a$ of~$\scrA$.

Show that $a^n$ is positive for \emph{even} $n\in \N$ and~$a\in\sa{\scrA}$.

Show that $a^n$ is positive iff $a$ is positive for \emph{odd} $n\in \N$
and $a\in\sa{\scrA}$.

Show that $a^n$ is positive
for every positive $a$ from~$\scrA$ and~$n\in \N$.
\item
Let~$a$ be an invertible element of~$\scrA$.
Show that $a\geq 0$ iff $a^{-1}\geq 0$.

\item
Show that a positive element~$a$ of~$\scrA$ is invertible
iff $a\geq \frac{1}{n}$ for some~$n>0$.
(Hint: show that $\spec(a)\subseteq [\frac{1}{n},\infty)$
		when~$a\geq \frac{1}{n}$.)
\end{enumerate}%
\spacingfix{}%
\end{point}%
\end{parsec}%
\begin{parsec}{180}%
\begin{point}{10}{[Moved to \sref{cstar-product-2}.]}
\end{point}
\end{parsec}
\begin{parsec}{190}
\begin{point}{10}%
Although we can't quite yet see that~$a^*a$
is positive---for this we need the existence of the square root, 
    \sref{sqrt},---we
    can already prove that~$a^*a$ can't be negative,
    see~\sref{astara-non-negative}.
\end{point}
\begin{point}{11}[prod-spec]{Lemma}%
For elements $a$ and $b$ from a $C^*$-algebra,
we have
\begin{equation*}
\spec(ab)\backslash\{0\}\ =\ \spec(ba)\backslash\{0\}.
\end{equation*}%
\spacingfix{}%
\begin{point}{20}{Proof}%
Let~$\lambda\in \C$ with $\lambda\neq 0$ be given.
We must show that $\lambda - ab$ is invertible
iff $\lambda - ba$ is invertible.
Suppose that $\lambda-ab$ is invertible.
Then using the equality $a(\lambda-ba)=(\lambda-ab)a$
one sees that $(1+b(\lambda-ab)^{-1}a)(\lambda-ba)=\lambda$.
Since similarly $(\lambda-ba)(1+b(\lambda-ab)^{-1}a)=\lambda$,
we see that $\lambda^{-1}(1+b(\lambda-ab)a)$
is the inverse of~$\lambda-ba$.\qed
\end{point}
\end{point}
\begin{point}{30}[astara-non-negative]{Lemma}%
We have $a^*a  \leq 0\implies a=0$
for every element~$a$ of a $C^*$-algebra.
\begin{point}{40}{Proof}%
Suppose that $a^*a\leq 0$.
Then~$\spec(a^*a)\subseteq (-\infty,0]$, almost by definition,
and so $\spec(aa^*)\subseteq (-\infty,0]$ by~\sref{prod-spec},
giving $aa^*\leq 0$.
Thus $a^*a+aa^*\leq 0$.

But on the other hand, 
$a^*a+aa^* = 2(\Real{a}^2 + \Imag{a}^2) \geq 0$,
and so~$a^*a+aa^*=0$.
Then $0\geq a^*a=-aa^*\geq 0$ gives $a^*a=0$,
and $a=0$.\qed
\end{point}
\end{point}
\end{parsec}
\begin{parsec}{200}%
\begin{point}{10}%
Observe that the norm and order on 
(the self-adjoint elements of a) $C^*$-algebra~$\scrA$
completely determine one another (using the unit):
on the one hand
$\|a\|=\inf\{\lambda\geq 0\colon -\lambda\leq a\leq \lambda\}$
by~\sref{positive-basic-2},
and on the other hand
$a\geq 0$ iff $\|a-s\|\leq s$ for some $s\in\R$
by definition (\sref{cstar-positive-def}).
This has some useful consequences.
\end{point}
\begin{point}{20}[weak-russo-dye]{Lemma}%
A positive map~$f\colon \scrA\to\scrB$
between $C^*$-algebras
is bounded.
More specifically,
we have
$\|f(a)\|\leq \|f(1)\|\,\|a\|$
for all self-adjoint~$a\in\sa{\scrA}$,
and we have $\|f(a)\|\leq 2\|f(1)\|\,\|a\|$
for arbitrary $a\in \scrA$.
\begin{point}{30}{Proof}%
Given~$a\in\sa{\scrA}$
we have~$-\|a\|\leq a \leq \|a\|$,
and $-\|a\|\,f(1)\leq f(a)\leq \|a\|\,f(1)$
(because $f$ is positive),
and thus~$\|f(a)\|\leq f(1)\,\|a\|\leq \|f(1)\|\,\|a\|$ 
by~\sref{positive-basic-2}.

For an arbitrary element $a\equiv \Real{a}+i\Imag{a}$
of~$\scrA$
we have 
$\|f(a)\|\leq \|f(\Real{a})\|+\|f(\Imag{a})\|\leq 2\|f(1)\|\,\|a\|$.\qed
\end{point}
\begin{point}{40}[russo-dye-remark]{Remark}%
It is a non-trivial theorem (see~\sref{russo-dye})
that the factor ``2'' in the statement above
can be dropped, i.e.~$\|f\|=\|f(1)\|$
(c.f.~Corollary~1 of~\cite{russodye}).
We'll be using this improved bound mostly
for completely positive maps,
for which it's much easier to obtain (see~\sref{cp-russo-dye}).

For miu-maps we can already obtain the improved bound here:
\end{point}
\end{point}
\begin{point}{50}[norm-mi-map]{Lemma}%
Any miu-map $\varrho\colon\scrA\to\scrB$ 
    between $C^*$-algebras $\scrA$ and~$\scrB$
is positive, bounded,
    and, in fact, $\|\varrho\|\leq 1$.
\begin{point}{51}{Proof}%
Let~$a$ be a positive element of~$\scrA$,
    so~$\spec(a)\subseteq[0,\infty)$ by~\sref{cstar-positive-1},
To show that~$\varrho$ is positive,
we must prove that~$\varrho(a)\geq 0$,
that is, $\spec(\varrho(a))\subseteq[0,\infty)$.
This follows immediately
from the observation that~$\spec(\varrho(a))\subseteq \spec(a)$:
 when $a-\lambda$ is invertible,
so is~$\varrho(\,a-\lambda\,)\equiv\varrho(a)-\lambda$,
for any~$\lambda\in \C$.
Hence~$\varrho$ is positive.

It follows by~\sref{weak-russo-dye}
that~$\varrho$ is bounded, and~$\|\varrho(b)\|\leq \|b\|$
for \emph{self-adjoint} $b\in \scrA$.
It remains to be shown that~$\|\varrho(a)\|\leq \|a\|$
for arbitrary~$a\in\scrA$.
But since~$a^*a$ is self-adjoint for such~$a$,
we have $\|\varrho(a)\|^2\equiv\|\varrho(a^*a)\|\leq \|a^*a\|=\|a\|^2$
by the $C^*$-identity and using that~$\varrho$
    is an miu-map. Whence $\|\varrho(a)\|\leq\|a\|$
    for all~$a\in\scrA$, and so~$\|\varrho\|\leq 1$.\qed
\end{point}
\end{point}
\begin{point}{60}[cstar-isometry]{Lemma}%
For a pu-map $f\colon \scrA\to\scrB$
the following are equivalent.
\begin{enumerate}
\item\label{cstar-isometry-1}
$f$ is \Define{bipositive}%
\index{bipositive!map between $C^*$-algebras}%
, that is, $f(a)\geq 0$ iff $a\geq 0$
for all~$a\in\scrA$;
\item\label{cstar-isometry-2}%
$f$ is an isometry on~$\sa{\scrA}$, 
that is, $\|f(a)\|=\|a\|$ for all~$\in \sa{\scrA}$;
\item\label{cstar-isometry-3}
	$f$ is an isometry on~$\pos{\scrA}$.
\end{enumerate}%
\spacingfix{}%
\begin{point}{70}{Proof}%
It is clear that \ref{cstar-isometry-2} implies~\ref{cstar-isometry-3}.
\begin{point}{80}{\ref{cstar-isometry-1}$\Longrightarrow$\ref{cstar-isometry-2}}%
Let~$a\in \sa{\scrA}$ be given.
Note that $-\lambda \leq a\leq \lambda$
iff $-\lambda \leq f(a) \leq \lambda$
for all~$\lambda \geq 0$,
because~$f$ is bipositive and unital.
In particular,
since~$-\|a\|\leq a\leq \|a\|$,
we have $-\|a\|\leq f(a)\leq \|a\|$,
and so~$\|f(a)\|\leq \|a\|$.
On the other hand,
$-\|f(a)\|\leq f(a)\leq \|f(a)\|$
implies $-\|f(a)\|\leq a\leq \|f(a)\|$,
and so $\|a\|\leq \|f(a)\|$.
Thus $\|a\|=\|f(a)\|$,
and $f$ is an isometry on~$\sa{\scrA}$.
\end{point}
\begin{point}{90}{\ref{cstar-isometry-3}$\Longrightarrow$\ref{cstar-isometry-1}}%
Let~$a\in \scrA$ be given.
We must show that~$f(a)\geq 0$ iff $a\geq 0$.
Since~$f$ is involution preserving (\sref{cstar-p-implies-i})
$a$ is self-adjoint iff $f(a)$ is self-adjoint,
and so we might as well assume that~$a$ is self-adjoint
to start with.
Since~$f$ is an isometry on~$\scrA_+$,
$\|a\|-a$ is positive,
and $f$ is unital,
we have $\|\,\|a\|-a\,\|=\|f(\|a\|-a)\|=\|\,\|a\|-f(a)\,\|$.
Now,
observe that
$0\leq a$
iff
$ \|\,\|a\|-a\,\|\leq \|a\|$,
and that
$\|\,\|a\|-f(a)\,\|\leq \|a\|$ 
iff $0\leq f(a)$,
by~\sref{positive-basic-2},
because $\frac{1}{2}\|a\|\leq \|a\|$
and $\frac{1}{2}\|f(a)\|\leq  \|a\|$
(by~\sref{weak-russo-dye}).\qed
\end{point}
\end{point}
\begin{point}{100}{Warning}%
Such a map~$f$ need not preserve the norm of arbitrary elements:
the map $A\mapsto \frac{1}{2}A+\frac{1}{2}A^T\colon M_2\to M_2$
is bipositive and unital,
but
\begin{equation*}
\left\|\left(\begin{matrix}0&1\\0&0\end{matrix}
\right)\right\|
\ = \ 1 \ \neq \ \frac{1}{2}\ = \ 
\left\|\,\left(\begin{matrix}0 & \nicefrac{1}{2} \\ 0 & 0
\end{matrix}\right)
\,+\,\left(\begin{matrix}0 & 0\\ \nicefrac{1}{2} & 0
\end{matrix}\right)\,\right\|.
\end{equation*}
(Even if~$f$ is completely positive, \sref{cp},
it might still only preserve the norm of self-adjoint elements
cf.~\sref{warning-norm-states}.)
\end{point}
\end{point}
\end{parsec}
\begin{parsec}{201}
\begin{point}{10}[cstar-product-2]{Exercise}%
\index{product!in $\Cstar{miu}$ and $\cCstar{miu}$}
\index{product!in $\Cstar{pu}$}%
Show that
the product~$\bigoplus_{i\in I}\scrA_i$
of
a family $(\scrA_i)_{j\in I}$
of $C^*$-algebras
defined in~\sref{cstar-product}
is also the categorical product 
of these $C^*$-algebras
    in~$\Cstar{miu}$ and~$\cCstar{miu}$
with as projections
the maps~$\Define{\pi_j}\colon \bigoplus_{i \in I}\scrA_i\to\scrA_j$%
\index{pi@$\pi_j$, projection!in $\Cstar{miu}$}
given by~$\pi_j(a)=a(j)$.

(Hint: use here that
the projections $\pi_j$ are bounded by~\sref{norm-mi-map}.)

Show that the same description applies to~$\cCstar{pu}$
and $\cCstar{pu}$.
(Hint: first show that an element~$a$
of $\bigoplus_{i\in I}\scrA_i$ is positive
    iff $a(i)$ is positive for every~$i\in I$.)

We'll return to the product of $C^*$-algebras
a final time in~\sref{cstar-product-4}.
\end{point}
\begin{point}{20}[cstar-equaliser-1]{Exercise}%
\index{equaliser!in~$\Cstar{miu}$ and $\Cstar{pu}$}
Show that given miu-maps $f,g\colon \scrA\to\scrB$
between $C^*$-algebras
the collection $\scrE:=\{a\in\scrA\colon f(a)=g(a)\}$
is a $C^*$-subalgebra
    of~$\scrA$ (using the fact that~$f$ and~$g$
    are bounded by~\sref{norm-mi-map} to show that~$\scrE$ is closed.)
Show that the inclusion $e\colon \scrE\to\scrA$
is a (positive) miu-map
that is in fact the equaliser of~$f$ and~$g$
in~$\Cstar{miu}$ and~$\Cstar{pu}$.
Show that the same description applies
to~$\cCstar{miu}$ and~$\cCstar{pu}$.
\begin{point}{30}[cstar-no-pu-equalisers]{Remark}%
The assumption here that~$f$ and~$g$ are miu-maps
is essential:
the pair of pu-maps $f,g\colon \C^4\to\C$
given by
\begin{equation*}
    \textstyle
f(a,b,c,d)\,=\, \frac{1}{2}(a+b),
    \quad \text{and}\quad
g(a,b,c,d)\,=\, \frac{1}{2}(c+d),
\end{equation*}
for example,
has no equaliser
in~$\Cstar{pu}$,
as we'll show in~\sref{cstar-no-pu-equalisers-example}.
\end{point}%
\end{point}%
\end{parsec}%
\begin{parsec}{210}%
\begin{point}{10}%
We just saw in~\sref{cstar-isometry} that a
map on a $C^*$-algebra~$\scrA$
that preserves and reflects the order
determines the norm of the self-adjoint
--- but not all --- elements of~$\scrA$.
This theme, to what extent a linear map
(or a collection of linear maps)
on a $C^*$-algebra
determines its structure,
while tangential at the moment,
will
grow ever more important 
until it is essential for the theory 
of von Neumann algebras.
That's why we introduce
the four levels
of discernment
that a collection of maps on a $C^*$-algebra
might have already here.
\end{point}
\begin{point}{20}[separating]{Definition}%
A collection~$\Omega$ of linear maps on a $C^*$-algebra~$\scrA$
will be called
\begin{enumerate}
\item
\label{separating-1}
\Define{order separating}%
\index{order separating collection!of maps on a $C^*$-algebras}
if an element~$a$ of~$\scrA$
is positive iff $0\leq \omega(a)$
for all~$\omega\in \Omega$;
\item
    \label{separating-2}
\Define{separating}%
\index{separating collection!of maps on a $C^*$-algebra} 
if an element~$a$ of~$\scrA$
is zero iff $\omega(a)=0$ for all~$\omega\in\Omega$;
\item
    \label{separating-3}
\Define{faithful} if an element~$a$ of~$\scrA_+$%
\index{faithful collection!of maps on a $C^*$-algebra}
is zero iff~$\omega(a)=0$ for all~$\omega\in\Omega$; and
\item
    \label{separating-4}
\Define{centre separating}%
\index{centre separating collection!of maps on a $C^*$-algebra}
if $a\in\scrA_+$
is zero iff $\omega(b^*ab)=0$ for all~$\omega\in\Omega$
and~$b\in \scrA$.

(The ``centre'' in ``centre separating''
will be explained in~\sref{vn-center-separating}.)
\end{enumerate}
(Note that 
    $\text{\eqref{separating-1}}
\implies
    \text{\eqref{separating-2}}
\implies
    \text{\eqref{separating-3}}
\implies
    \eqref{separating-4}$.)
\end{point}
\begin{point}{30}{Examples}%
We'll see later on that the following
collections
are order separating.
\begin{enumerate}
\item
The set of all pu-maps $\omega\colon \scrA\to\C$
    (called \Define{states}\index{state of a $C^*$-algebra})
on a $C^*$-algebra
(see~\sref{states-order-separating}).
\item
The set of all miu-maps $\omega\colon \scrA\to\C$
on a commutative $C^*$-algebra
(see~\sref{gelfand-representation-isometry}).
\item
The set of functionals on~$\scrB(\scrH)$,
where~$\scrH$ is a Hilbert space,
of the form 
$\left<x,(\,\cdot\,)x\right>
\colon \scrB(\scrH)\to\C$ where $x\in \scrH$
(see~\sref{hilb-vector-states-order-separating}).

We'll call these functionals
\Define{vector functionals}%
\index{vector functional!for a Hilbert space}%
\index{functional!vector}.
(They are clearly bounded
and involution preserving linear maps,
and once we know that each positive element of a $C^*$-algebra
is a square, in~\sref{sqrt},
it'll be obvious that vector functionals
are positive too.)

The unital vector functionals (called \Define{vector states})%
\index{vector state = unital vector functional}
are order separating too.
\end{enumerate}%
\spacingfix%
\begin{point}{40}
None of the four levels of separation
coincide.
This follows from the following examples,
that we'll just mention here,
but can't verify yet.
\begin{enumerate}
\item
A single non-zero vector~$x$ from a Hilbert space~$\scrH$
gives a vector functional $\left<x,(\,\cdot\,)x\right>$
on~$\scrB(\scrH)$ that is centre separating
on its own, but is not faithful when~$\scrH$ has dimension~$\geq 2$.
\item
Given an orthonormal basis~$\scrE$
of a Hilbert space~$\scrH$
		the collection 
\begin{equation*}
	\{\,\left<e,(\,\cdot\,)e\right>\colon\,e\in\scrE\,\}
\end{equation*}
		of vector functionals on~$\scrB(\scrH)$
		is faithful, but not separating when~$\scrE$
		has more than one element.
\item
Given Hilbert spaces~$\scrH$
and~$\scrK$
the set of vector functionals
\begin{equation*}
	\{\,\left<\,x\otimes y,\,(\,\cdot\,)\,x\otimes y\,\right>\colon\,
x\in\scrH,\,y\in\scrK\,\}
\end{equation*}
on~$\scrB(\scrH\otimes\scrK)$
is separating,
but not order separating
when both~$\scrH$ and~$\scrK$
are at least two dimensional.
\end{enumerate}%
\spacingfix{}%
\end{point}%
\end{point}%
\begin{point}{50}[separating-self-adjoint]{Exercise}%
One use for 
a separating collection~$\Omega$
of involution preserving maps
on a $C^*$-algebra~$\scrA$
is checking whether an element~$a\in\scrA$
is self-adjoint:
show that $a\in\scrA$ is self-adjoint
iff $\omega(a)$ is self-adjoint for all~$\omega\in\Omega$.
\end{point}
\begin{point}{60}%
An order separating
collection senses the norm 
of a self-adjoint element:
\end{point}
\begin{point}{70}[order-separating-norm]{Proposition}%
For a collection~$\Omega$ of pu-maps on a $C^*$-algebra~$\scrA$
the following are equivalent.
\index{order separating collection!of pu-maps on a $C^*$-algebra}
\begin{enumerate}
\item 
	$\Omega$ is order separating;
\item
	$\|a\|= \sup_{\omega\in\Omega} \left\|\omega(a)\right\|$
	for all $a\in \sa{\scrA}$;
\item
	$\|a\| = \sup_{\omega\in\Omega} \left\|\omega(a)\right\|$
	for all~$a\in \pos{\scrA}$.
\end{enumerate}%
\spacingfix{}%
\begin{point}{80}{Proof}%
Denoting the codomain of~$\omega\in\Omega$
by~$\scrB_\omega$
(so that $\omega\colon \scrA\to\scrB_\omega$),
apply~\sref{cstar-isometry}
to the pu-map $\left<\omega\right>_{\omega\in\Omega}\colon 
\scrA\to\bigoplus_{\omega\in\Omega}
\scrB_\omega$ (see~\sref{cstar-product-2}).\qed
\end{point}
\begin{point}{90}[warning-norm-states]{Warning}%
The formula
$\|a\|=\sup_{\omega\in\Omega} \|\omega(a)\|$
need not be correct
for an arbitrary (not necessarily self-adjoint)
element~$a$.
Indeed,
consider the matrix $A:=\smash{%
\bigl(\begin{smallmatrix}0&1\\0&0\end{smallmatrix}\bigr)}$,
and the collection $\Omega=\{\,\left<x,(\,\cdot\,)x\right>\colon 
x\in \C^2,\,\|x\|=1\,\}$,
which will turn out to be order separating.
We have $\|A\|=1$,
while $\left| \left<x,\omega(A)x\right>\right|
 =\left|x_1\right|\left|x_2\right|$
 never exceeds~$\nicefrac{1}{2}$
for $x\equiv (x_1,x_2)\in \scrH$ with $1=\|x\|$.
\end{point}
\end{point}
\begin{point}{100}[order-separating-dense-subset]{Exercise}%
Show that any operator norm dense subset~$\Omega'$
of an order separating collection~$\Omega$
of positive functionals
on a $C^*$-algebra~$\scrA$
is order separating too.
\end{point}
\end{parsec}
\begin{parsec}{220}%
\begin{point}{10}%
We'll use~\sref{order-separating-norm}
to show 
that the pu-maps $\omega\colon \scrA\to\C$
on a $C^*$-algebra~$\scrA$
(called states
of~$\scrA$ for short)
are order separating
by showing that
for every self-adjoint element~$a\in \scrA$
there is a state~$\omega$ of~$\scrA$ with $\omega(a)=\|a\|$ or 
$\omega(a)=-\|a\|$.
To obtain such a state
we first find its kernel,
which leads us to the following definitions.
\end{point}
\begin{point}{20}{Definition}%
An \Define{order ideal}%
\index{order ideal of a $C^*$-algebra}
of a $C^*$-algebra~$\scrA$
is a linear subspace~$I$ of~$\scrA$
with $b\in I\implies b^*\in I$
and $b\in I\cap\pos{\scrA}\implies [-b,b]
    \equiv \{\,a\in\scrA\colon \,-b\leq a\leq b\, \} \ \subseteq\, I$.

The order ideal~$I$ is called \Define{proper}%
\index{order ideal of a $C^*$-algebra!proper}
if~$1\notin I$,
and \Define{maximal} 
\index{order ideal of a $C^*$-algebra!maximal}
if it is maximal among all proper order ideals.
\begin{point}{21}{Warning}%
``Order ideals'' like ``subspaces'' appear in relation
to other structures as well,
with appropriately varying meanings.
Our definition for $C^*$-algebras is based on
to the order ideals for order unit spaces
from Definition~2.2 of~\cite{kadison1951}.
\end{point}
\end{point}
\begin{point}{30}[order-ideal-basic]{Exercise}%
Let~$\scrA$ be a $C^*$-algebra.
\begin{enumerate}
\item
Show that the kernel of a state is a maximal order ideal.

(Hint: the kernel of a state is already maximal as linear subspace.)
\item
Let~$I$ be a proper order ideal of~$\scrA$.
Show that there is a maximal 
order ideal~$J$ of~$\scrA$ with $I\subseteq J$.
(Hint: Zorn's Lemma may be useful.)
\item
Let~$a\in \sa{\scrA}$.
Show that there is a least order ideal~\Define{$(a)$}%
\index{$(a)$, order ideal generated by $a$}
that contains~$a$,
and that given~$b\in\Real{\scrA}$
we have $b\in (a)$
iff there are $\lambda,\mu\in \R$
with~$\lambda a\leq b\leq \mu a$.

Show that~$(a)=\C a$
when~$0\nleq a\nleq 0$.

Show that~$1\in (a)$ if and only if $a$ is invertible
and either $0\leq a$ or $a\leq 0$.

\item
Let~$a$ be a self-adjoint element of~$\scrA$ which
is not invertible.
Show that there is a maximal order ideal~$J$
of~$\scrA$
with $a\in J$.

\item
Let~$a$ be a self-adjoint element of~$\scrA$.
Show that  $\|a\|-a$
or $\|a\|+a$ is not invertible
(perhaps by considering the spectrum of~$a$.)
\end{enumerate}%
\spacingfix{}%
\end{point}%
\begin{point}{40}[maximal-ideal-state]{Lemma}%
For every maximal order ideal~$I$ of a $C^*$-algebra~$\scrA$,
 there is a state $\omega \colon \scrA\to \C$
with $\ker(\omega)=I$.
\begin{point}{50}{Proof}%
Form the quotient vector space $\scrA/I$
with quotient map $q\colon \scrA\to \scrA/I$.
Note that since~$1\notin I$
we have $q(1)\neq 0$
and so we may regard~$\C$ 
to be a linear subspace of~$\scrA/I$
via $\lambda\mapsto q(\lambda)$.
We will, in fact, show that~$\C=\scrA/I$.

But let us first put an order on~$\scrA/I$:
we say that $\mathfrak{a}\in \scrA/I$ is positive
if $\mathfrak{a}\equiv q(a)$ for some~$a\in\pos{\scrA}$,
and write $\mathfrak{a}\leq \mathfrak{b}$ 
if $\mathfrak{b}-\mathfrak{a}$ is positive
for $\mathfrak{a},\mathfrak{b}\in\scrA/I$.
Note that the definition of ``order ideal'' is such
that if both~$\mathfrak{a}$ and $-\mathfrak{a}$ are positive,
then~$\mathfrak{a}=0$.
We leave it to the reader to verify 
that~$\scrA/I$ becomes a partially ordered vector space
with the order defined above.
There is, however,
one detail we'd like to draw attention to,
namely that a scalar $\lambda$ is positive in~$\scrA/I$
iff $\lambda$ is positive in~$\C$.
Indeed, if~$\lambda\geq 0$ in~$\C$,
then $\lambda\geq 0$ in~$\scrA$, and so~$\lambda \geq 0$ in~$\scrA/I$.
On the other hand,
if~$\lambda\geq 0$ in~$\scrA/I$, but~$\lambda\leq 0$ in~$\C$,
then $\lambda\leq 0$ in~$\scrA/I$,
and so $\lambda=0$.
This detail
has the pleasant consequence
that once we have shown that~$\scrA/I=\C$,
we automatically get that~$q\colon \scrA\to\C$ is positive.

\begin{point}{60}[pos-hahn-banach-1]%
Let~$a\in \sa{\scrA}$ be given.
Define~$\alpha := \inf\{\,\lambda\in\R\colon\, q(a)\leq \lambda\,\}$.
Note that $-\|a\| \leq \alpha\leq \|a\|$.
We will prove that~$q(a)=\alpha$
by considering the order ideal
\begin{alignat*}{3}
	J\ := \ \{\,b\in \scrA\colon\, 
		& \exists\lambda,\mu\in\R\,[\ 
	\lambda (\alpha-q(a))\,\leq\, \Real{b}\,\leq\, \mu 
(\alpha-q(a))\ ]\,\wedge\,
		\\
		&\exists\lambda,\mu\in \R\,[\ 
\lambda (\alpha-q(a))\,\leq\, \Imag{b}\,\leq\, \mu (\alpha-q(a)) \ ] \,\}.
\end{alignat*}
We claim that $1\notin J$.
Indeed, suppose not---towards a contradiction.
Then there is~$\mu\in \R$
with $1\leq \mu (\alpha-q(a))$.
What can we say about~$\mu$?
If~$\mu <0$,
then $0\geq \nicefrac{1}{\mu}\geq \alpha-q(a)$,
so~$\alpha-\nicefrac{1}{\mu} \leq q(a)$,
but $q(a)\leq \alpha+\varepsilon$
for every~$\varepsilon>0$,
and so~$\alpha-\nicefrac{1}{\mu}\leq q(a)\leq \alpha-\nicefrac{1}{2\mu}$,
which is absurd.
If $\mu=0$,
then we get $1\leq \mu (\alpha-q(a))\equiv 0$, which is absurd.
If $\mu> 0$,
then $\nicefrac{1}{\mu}\leq \alpha-q(a)$,
or in other words,
 $q(a) \leq \alpha - \nicefrac{1}{\mu}$,
 giving $\alpha \leq \alpha-\nicefrac{1}{\mu}$
by definition of~$\alpha$,
which is absurd.
Hence~$1\notin J$.

But then since~$I\subseteq J$,
we get~$I=J$, by maximality of~$I$.
Thus, as $\alpha-a\in J$, we have $\alpha-a\in I$,
and so $q(a)=\alpha$, as desired.
\end{point}
\begin{point}{70}%
Let~$a\in \scrA$ be given.
Then~$a=\Real{a}+i\Imag{a}$.
By~\sref{pos-hahn-banach-1},
there are $\alpha,\beta\in \R$ with $q(\Real{a})=\alpha$,
and $q(\Imag{a})=\beta$.
Thus~$q(a)=\alpha+i\beta$.
Hence~$\scrA/I=\C$.
Since the quotient map $q\colon \scrA\to \scrA/I\equiv \C$
is pu, and $\ker(q)=I$, we are done.\qed
\end{point}
\end{point}
\end{point}
\begin{point}{80}[states-order-separating]{Exercise}%
\index{state of a $C^*$-algebra!order separating}%
Show using~\sref{maximal-ideal-state} that given a self-adjoint element~$a$
of a $C^*$-algebra~$\scrA$
there is a state~$\omega$ with $\left|\omega (a)\right| = \|a\|$.
Conclude that the set of states of a $C^*$-algebra
is order separating (see~\sref{separating}).
\end{point}
\end{parsec}
\subsection{The Square Root}
\begin{parsec}{230}%
\begin{point}{10}%
The key that unlocks the remaining basic facts 
about the (positive) elements of a  $C^*$-algebra
is the existence of the square root~$\sqrt{a}$ of a positive element~$a$,
and its properties.
For technical reasons,
we will assume $\|a\|\leq 1$,
and construct
 $1-\sqrt{1-a}$ instead of~$\sqrt{a}$.
\end{point}
\begin{point}{20}{Lemma}%
Let $a$ be an element of a $C^*$-algebra $\scrA$
with $0\leq a\leq 1$.
Then there is a unique element~$b\in\scrA$ 
with, $0\leq b\leq 1$,
$ab=ba$,
and~$(1-b)^2 = 1-a$.
To be more specific,
$b$ is the norm limit of
the sequence $b_0\leq b_1\leq \dotsb$
given by $b_0=0$ and $b_{n+1} = \frac{1}{2}(a+b_n^2)$.
Moreover,
if~$c\in\scrA$ commutes with~$a$, then~$c$ commutes with~$b$,
and if in addition $c^2\leq 1-a$ and $c^*=c$,
we have $c\leq 1-b$.
\begin{point}{30}{Proof}%
When discussing $b_n$ it 
is convenient to write~$b_n \equiv q_n(a)$
where~$q_0,q_1,\dotsc$ are the polynomials over~$\R$ given by
$q_0=0$ and $q_{n+1}=\frac{1}{2}(x + q_n^2)$.
For example,
we have~$b_n\geq 0$, 
because all coefficients of~$q_n$ are all positive,
and $a,a^2,a^3,\dotsc$ are positive by~\sref{positive-basic-2}.
With a similar argument we can see that
 $b_0 \leq b_1\leq b_2\leq \dotsb$.
Indeed, 
the coefficients of~$q_{n+1}-q_n$
are positive,
by induction,
because
\begin{alignat*}{3}
q_{n+2}-q_{n+1} \ &=\ \textstyle \frac{1}{2}(x+ q_{n+1}^2)
\,-\, \textstyle\frac{1}{2}(x+q_n^2) \\
&=\ \textstyle\frac{1}{2}(q_{n+1}^2- q_n^2) \\
&=\ \textstyle\frac{1}{2}(q_{n+1}+q_n)(q_{n+1}-q_n) \\
&=\ (q_n+\textstyle\frac{1}{2}(q_{n+1}-q_n))(q_{n+1}-q_n),
\end{alignat*}
has positive coefficients
if~$q_{n+1}-q_n$ has positive coefficients,
and $q_1-q_0\equiv \frac{1}{2}x$ clearly has positive coefficients.
Hence~$b_{n+1}-b_{n} = q_{n+1}(a)- q_n(a)$ is positive.
(Note that we have carefully avoided
using the fact here that the product of positive 
commuting elements is positive,
which is not available to us until~\sref{ineq-square-root}.)

Let us now show that~$b_0\leq b_1\leq \dotsb$ converges.
Let~$n\geq N$ from~$\N$ be given.
Since the coefficients of $q_n-q_N$ are positive,
and $\|a\|\leq 1$,
the triangle inequality gives us
$\|b_n-b_N\|\equiv \|(q_n-q_N)(a)\|\leq q_n(1)-q_N(1)$,
and
so it suffices to 
show that the ascending sequence
 $q_0(1)\leq q_1(1)\leq \dotsb$
of real numbers
converges,
i.e.~is bounded.
Indeed,
we have $q_n(1)\leq 1$,
by induction,
because $q_{n+1}(1)\equiv \frac{1}{2}(1+q_n(1)^2)
\leq 1$ if $q_n(1)\leq 1$,
and clearly $0\equiv q_0(1)\leq 1$.

Let~$b$ be the limit of $b_0\leq b_1\leq\dotsb$.
Then~$b$ being the limit of positive elements
is positive
(see~\sref{positive-basic-2}),
and if $c\in \scrA$ commutes with~$a$,
then $c$ commutes with all powers of~$a$,
and therefore with all~$b_n$,
and thus with~$b$.
Further, 
from the recurrence relation $q_{n+1} = \frac{1}{2}(a+q_n^2)$
we get $b=\frac{1}{2}(a+b^2)$,
and so $-a = -2b+b^2$, 
giving us  $(1-b)^2 = 1-2b+b^2 = 1-a$.

Let us prove that~$b\leq 1$.
To begin, note that~$\|b_n\|\leq 1$ for all~$n$, by induction,
because $0\equiv \|b_0\|\leq 1$,
and if $\|b_n\|\leq 1$, then $\|b_{n+1}\|\leq \frac{1}{2}(\|a\|+\|b_n\|^2)
\leq 1$, since $\|a\|\leq 1$.
Since~$b_n\geq 0$, we get $-1\leq b_n\leq 1$ for all~$n$,
and so $b\leq 1$.
\begin{point}{40}[square-commuting-monotone]%
Let us take a step back for the moment.
From what we have proven so far
we see that each positive $c\in\scrA$
is of the form $c\equiv d^2$ for some positive~$d\in\scrA$
which commutes with all~$e\in \scrA$ that commute with~$c$.

From this we can see that $c_1c_2\geq 0$
for  
 $c_1,c_2 \in\pos{\scrA}$
with $c_1c_2 = c_2c_1$.
Indeed, writing $c_i\equiv d_i^2$ with $d_i$ as above,
we have $d_1c_2=c_2d_1$ (because $c_1c_2=c_2c_1$), and thus 
$d_1d_2=d_2d_1$. It follows that $d_1d_2$ is self-adjoint,
and $c_1c_2 = (d_1d_2)^2$. Hence $c_1c_2\geq 0$.

We will also need the following corollary.
For~$c,d\in\pos{\scrA}$ with $c\leq d$ and $cd=dc$,
we have $c^2\leq d^2$.
Indeed, $d^2-c^2 \equiv d(d-c)+c(d-c)$
is positive by the previous paragraph.
\end{point}
\begin{point}{50}[ineq-square-root]%
Let~$c\in\sa{\scrA}$ be such that~$ca=ac$ and $c^2\leq 1-a$, 
    that is, $a\leq 1-c^2$.
We must show that $c\leq 1-b$,
that is, $b\leq 1-c$.
Of course,
since~$b$ is the limit of $b_1,b_2,\dotsc$,
it suffices to show that~$b_n\leq 1-c$,
and we'll do this by induction.
Since $0\leq c^2 \leq 1-a$,
 we have $\|c\|^2\leq \|1-a\|\leq 1$,
and so $-1\leq c\leq 1$.
Thus $b_0\equiv 0\leq 1-c$.
Now, suppose that~$b_n\leq 1-c$ for some~$n$.
Then $b_{n+1} = \frac{1}{2}(a+b_n^2)
\leq \frac{1}{2}( (1-c^2)+(1-c)^2) = 1-c$,
where we have used that $b_n^2 \leq (1-c)^2$,
because $b_n\leq 1-c$
by~\sref{square-commuting-monotone}.
\begin{point}{60}%
We'll now show that~$b$ is unique
in the sense that $b=b'$
for any~$b'\in \scrA$ with $0\leq b'\leq 1$,
 $b'a=ab'$ and $(1-b')^2=1-a$.
Note that $b'\leq 1$,
because $\|1-b'\|^2=\|1-a\|\leq 1$,
From $a=1-(1-b')^2$,
we immediately get $b \leq 1-(1-b')=b'$ by~\sref{ineq-square-root}.
For the other direction,
note that
$(1-b')^2= (1-b)^2 \equiv (1-b'+(b'-b))^2 = (1-b')^2+2(1-b')(b'-b)+(b'-b)^2$,
which gives $0=2(1-b')(b'-b)+(b'-b)^2$.
Now, since~$1-b'$ and $b'-b$ are positive,
and commute, we see that $(1-b')(b'-b)$ is positive 
by~\sref{ineq-square-root}, and so 
 $0=2(1-b')(b'-b)+(b'-b)^2\geq (b'-b)^2 \geq 0$,
which entails $(b'-b)^2=0$, and so $\|(b'-b)^2\|=\|b'-b\|^2=0$,
yielding $b=b'$.\qed
\end{point}
\end{point}
\end{point}
\end{point}
\begin{point}{70}[sqrt]{Exercise}%
\index{*sqrt@$\sqrt{a}$, square root!in a $C^*$-algebra}%
Let~$a$ be a positive element of a $C^*$-algebra~$\scrA$.
Show that there is a unique 
positive element of~$\scrA$
denoted by $\Define{\sqrt{a}}$ 
(and by~$\Define{a^{\nicefrac{1}{2}}}$)
with $\smash{\sqrt{a}^2}=a$
and $a\sqrt{a}=\sqrt{a}a$.
Show that if~$c\in\scrA$ commutes with~$a$,
then $c\sqrt{a}=\sqrt{a}c$,
and if in addition $c^*=c$ and $c^2\leq a$,
then $c\leq \sqrt{a}$.
Using this, verify:
\begin{enumerate}
\item
If~$a,b\in \scrA$ are positive,
and~$ab=ba$,
then $ab\geq 0$.

\item
Let~$a\in\pos{\scrA}$.
If $b,c\in \sa{\scrA}$ commute with~$a$,
then $b\leq c$ implies $ab\leq ac$.

\item
If~$a,b\in\Real{\scrA}$ commute, and~$a\leq b$, then~$a^2\leq b^2$.

\item
The requirement in the previous item  that~$a$ and~$b$ commute is essential:
there are positive elements $a$, $b$ of a $C^*$-algebra~$\scrA$
with $a\leq b$, but $a^2 \nleq b^2$.

In other words, the square $a\mapsto a^2$
on the positive elements of a $C^*$-algebra
need not be monotone,
(but $a\mapsto \sqrt{a}$ \emph{is} monotone, see~\sref{sqrt-monotone}).

(Hint: take $a=(\begin{smallmatrix}1&0\\0&0\end{smallmatrix})$
and $b=a+\frac{1}{2}(\begin{smallmatrix}1&1\\1&1\end{smallmatrix})$
from~$M_2$.)
\end{enumerate}%
\spacingfix{}%
\end{point}%
\end{parsec}%
\begin{parsec}{240}
\begin{point}{10}{Definition}
Given a self-adjoint element~$a$ of a $C^*$-algebra $\scrA$,
we write
\begin{equation*}
\textstyle
\Define{\left|a\right|}\ :=\ \sqrt{a^2}
\qquad
\Define{\pos{a}}\ :=\ \frac{1}{2}(\left|a\right| + a)
\qquad
\Define{a_{-}}\ :=\ \frac{1}{2}(\left|a\right| - a).
\end{equation*}%
\index{$(\,\cdot\,)_+$, positive part!$a_+$, of a self-adjoint element of a $C^*$-algebra}
We call $a_+$ the \Define{positive part} of~$a$,
and $a_-$ the \Define{negative part}.
\index{$(\,\cdot\,)_-$, negative part!$a_-$, of a self-adjoint element of a $C^*$-algebra}
\end{point}
\begin{point}{20}[cstar-pos-neg-part]{Exercise}%
Let~$a$ be a self-adjoint element of a
 $C^*$-algebra $\scrA$.
\begin{enumerate}
\item
Show that $-\left|a\right| \leq a \leq \left| a \right|$,
and $\|\,\left|a\right|\,\|= \|a\|$.
\item
Prove that $a_+$ and $a_-$ are positive,  $a=a_+-a_-$
and $a_+a_-=a_-a_+=0$.
\item
One should not read too much into the notation
$\left|\,\cdot\,\right|$
in the non-commutative case:
give an example of
self-adjoint elements~$a$ and~$b$ of a $C^*$-algebra with
 $\left|a+b\right|\nleq \left|a\right|+ \left|b\right|$.

(Hint: one may take  
$a=\frac{1}{2}\left(\begin{smallmatrix}1 & 1 \\ 1 & 1\end{smallmatrix}\right)$
and $b=-\left(\begin{smallmatrix}1 & 0 \\ 0 & 0 \end{smallmatrix}\right)$.)
\end{enumerate}
\spacingfix%
\end{point}
\begin{point}{30}%
The existence of positive and negative parts
in a $C^*$-algebra
has many pleasant and subtle consequences
of which we'll now show one.
\end{point}
\begin{point}{40}[astara-positive]{Lemma}%
Given an element $a$ of a $C^*$-algebra $\scrA$,
we have $a^*a\geq 0$.
\begin{point}{50}{Proof}%
Writing $b:=a((a^*a)_-)^{\nicefrac{1}{2}}$,
we have $b^*b=
((a^*a)_-)^{\nicefrac{1}{2}} a^*a
((a^*a)_-)^{\nicefrac{1}{2}}
= (a^*a)_- \,a^*a
=
-((a^*a)_-)^2\leq 0$,
and so
$b=0$
by~\sref{astara-non-negative}.
    Hence~$((a^*a)_-)^2=0$,
    and thus~$(a^*a)_-=0$ (by, say, the $C^*$-identity,)
    giving us $a^*a=(a^*a)_+\geq 0$.\qed
\end{point}
\end{point}
\end{parsec}
\begin{parsec}{250}
\begin{point}{10}[cstar-positive-final]{Exercise}%
\index{positive!element of a $C^*$-algebra}%
Round up our results regarding positive elements
to 
prove that
the following are equivalent
for a self-adjoint element $a$ of a $C^*$-algebra~$\scrA$.
\begin{enumerate}
\item 
$a$ is positive, that is,  $\|a-t\|\leq t$
for some $t\in \R$;
\item
$\|a-t\|\leq t$ for all~$t\geq \frac{1}{2}\|a\|$;
\item
$a\equiv b^2$ for some self-adjoint $b\in\scrA$;
\item
$a\equiv c^* c$ for some $c\in\scrA$;
\item
$\spec(a)\subseteq [0,\infty)$.
\end{enumerate}%
\spacingfix%
\end{point}%
\begin{point}{20}[astara-pos-basic-consequences]{Exercise}%
The fact that $a^*a$ is positive
for an element~$a$ of a $C^*$-algebra~$\scrA$
has some nice consequences
of its own needed later on.
\begin{enumerate}
\item
Show that $b\leq c\implies a^*ba \leq a^*ca$
for all~$b,c\in\sa{\scrA}$ and~$a\in\scrA$.
\item
Show that every mi-map and cp-map is positive.
\item
Show that~$a\leq b^{-1}$ 
iff $\sqrt{b}a\sqrt{b}\leq 1$
iff $\|\sqrt{a}\sqrt{b}\|\leq 1$
iff $b\leq a^{-1}$
for positive invertible elements $a$, $b$ of~$\scrA$
(and so $a\leq b$ entails $b^{-1}\leq a^{-1}$).
\item
Prove that $(1+a)^{-1}a\leq (1+b)^{-1}b$
for $0\leq a\leq b$ from~$\scrA$.\\
(Hint: add $(1+a)^{-1} + (1+b)^{-1}$
to both sides of the inequality.)
\end{enumerate}%
\spacingfix%
\end{point}%
\begin{point}{30}[hilb-vector-states-order-separating]{Proposition}%
The vector states
of~$\scrB(\scrH)$
are order separating (see~\sref{separating})
for every Hilbert space~$\scrH$.
\begin{point}{40}{Proof}%
By~\sref{order-separating-norm}
it suffices to show
that~$\|T\|=  \sup_{x\in (\scrH)_1} 
\left|\left< x,Tx\right>\right|$
for given~$T\in\scrB(\scrH)_+$.
Since $\left|\left<x,Tx\right>\right|
=\left<T^{\nicefrac{1}{2}}x,T^{\nicefrac{1}{2}}x\right>
=\|T^{\nicefrac{1}{2}}x\|^2$
for all~$x\in \scrH$, 
we have $
\|T\| = \|T^{\nicefrac{1}{2}}\|^2
=(\,\sup_{x\in (\scrH)_1}\left\|T^{\nicefrac{1}{2}}x\right\|\,)^2
=\sup_{x\in (\scrH)_1} \left|\left<x,Tx\right>\right|$.\qed
\end{point}
\end{point}
\begin{point}{50}[hilb-positive-operators]{Corollary}%
For a bounded operator~$T$
on a Hilbert space~$\scrH$, we have
\begin{enumerate}
\item
$T$ is self-adjoint iff $\left<x,Tx\right>$
is real for all~$x\in (\scrH)_1$;
\item
$0\leq T$ iff $0\leq\left<x,Tx\right>$
for all~$x\in (\scrH)_1$;
\item
	$\|T\|=\sup_{x\in (\scrH)_1}\left|\left<x,Tx\right>\right|$
when~$T$ is self-adjoint.
\end{enumerate}%
\spacingfix%
\begin{point}{60}{Proof}%
This follows from~\sref{separating-self-adjoint} 
and~\sref{order-separating-norm}
because the vector states on~$\scrB(\scrH)$
are order separating 
by~\sref{hilb-vector-states-order-separating}.\qed
\end{point}
\end{point}
\end{parsec}
\begin{parsec}{260}%
\begin{point}{10}%
The interaction between the multiplication and order
on a $C^*$-algebra can be subtle,
but
when the $C^*$-algebra is commutative
almost all peculiarities disappear.
This is to be expected
as any commutative $C^*$-algebra
is isomorphic to a $C^*$-algebra
of continuous functions on a compact Hausdorff space
(as we'll see in~\sref{gelfand}).
\end{point}
\begin{point}{20}[commutative-cstar-basic]{Exercise}%
Let~$\scrA$ be a \emph{commutative} $C^*$-algebra.
Let~$a,b,c\in\sa{\scrA}$.
\begin{enumerate}
\item
Show that $\left| a\right|$ is the supremum of~$a$ and~$-a$
in~$\sa{\scrA}$.
\item
Show that if~$a$ and $b$ have a supremum, $a\vee b$, in $\sa{\scrA}$,
then~$c\,+\,a\vee b$ is the supremum of~$a+c$ and $b+c$.
\item
Show that~$\sa{\scrA}$ is a \Define{Riesz space},
that is,  a lattice ordered vector space.\\
(Hint: prove that $\frac{1}{2}(a+b+\left|a-b\right|)$
is the supremum of~$a$ and~$b$ in~$\sa{\scrA}$.)
\item
Show that an miu-map $f\colon \scrA\to\scrB$
between commutative $C^*$-algebras
preserves finite suprema and infima.
\end{enumerate}%
\spacingfix%
\end{point}%
\begin{point}{30}[riesz-decomposition-lemma]{Exercise}%
Prove the \Define{Riesz decomposition lemma}:%
\index{Riesz decomposition lemma}
For positive elements~$a,b,c$ of a commutative $C^*$-algebra~$\scrA$
with~$c\leq a+b$
we have $c\equiv a'+b'$
where  $0\leq a'\leq a$ and $0\leq b'\leq b$.
\end{point}
\end{parsec}
\section{Representation}
\subsection{\dots by Continuous Functions}
\begin{parsec}{270}%
\begin{point}{10}%
Now that we have have a firm grip
on the positive elements of a $C^*$-algebra
we turn to what is perhaps the most important
fact about commutative $C^*$-algebras:
that they are isomorphic to $C^*$-algebras
of continuous functions on a compact Hausdorff space,
via the \emph{Gelfand representation}.
\end{point}
\begin{point}{20}{Setting}%
$\scrA$ is a commutative $C^*$-algebra.
\end{point}
\begin{point}{30}[gelfand-representation]{Definition}%
The \Define{spectrum} of~$\scrA$,%
\index{sp@$\spec$, spectrum!$\spec(\scrA)$, of a $C^*$-algebra}
denoted by \Define{$\spec(\scrA)$},
is the set of all miu-maps $f\colon \scrA\to \C$.
We endow~$\spec(\scrA)$
with the topology of pointwise convergence.

The \Define{Gelfand representation}
of~$\scrA$
is the miu-map~$\gamma\colon \scrA\to C(\spec(\scrA))$%
\index{Gelfand representation@$\gamma$, Gelfand representation}
given by $\gamma(a)(f)=f(a)$.
\end{point}
\begin{point}{40}[gelfand-representation-basic]{Exercise}%
Verify that 
 the map $\spec(\scrA)\to \C,\ f\mapsto f(a)$ is indeed
continuous for every~$a\in\scrA$,
and that~$\gamma$ is miu.
\end{point}
\begin{point}{50}{Remark}
One might wonder if there is any connection between
the spectrum~$\spec(\scrA)$
of a commutative $C^*$-algebra,
and the spectrum~$\spec(a)$
of one of~$\scrA$'s elements (from~\sref{spectrum-of-element});
and indeed there is
as we'll see in~\sref{spectrum-miu}
(and~\sref{functional-calculus}).
\end{point}
\begin{point}{60}%
Our program for this paragraph is to show that
the Gelfand representation~$\gamma$ is 
an miu-isomorphism.
In fact,
we will show that it gives the unit
of an equivalence between the category of commutative $C^*$-algebras
(with miu-maps)
and the opposite of the category of compact Hausdorff spaces
(with continuous maps).
The first hurdle we take is the injectivity of~$\gamma$
--- that there are sufficiently many
points in the spectrum of a commutative $C^*$-algebra,
so to speak ---,
and involves
the following special type of order ideal.
\end{point}
\begin{point}{70}{Definition}%
A \Define{Riesz ideal}%
\index{Riesz ideal}
of~$\scrA$
is an order ideal~$I$
such that $a\in I\cap\sa{\scrA}\implies \left|a\right|\in I$.
A \Define{maximal Riesz ideal}%
	\index{Riesz ideal!maximal}
is a proper Riesz ideal which is maximal among
proper Riesz ideals.
\end{point}
\begin{point}{80}[riesz-ideal-ring-ideal]{Lemma}%
Let~$I$ be a Riesz ideal of~$\scrA$.
For all~$a\in \scrA$ and $x\in I$ we have $ax\in I$.
\begin{point}{90}{Proof}%
Since~$x=\Real{x}+i\Imag{x}$,
it suffices to show that~$a\Real{x}\in I$ and $a\Imag{x}\in I$.
Note that~$\Real{x},\Imag{x}\in I$,
so we might as well assume that~$x$ is self-adjoint to begin with.
Similarly, using that
 $\pos{x}\in I$ (because $\pos{x}=\frac{1}{2}(\left|x\right|+x)$
and~$\left|x\right|\in I$) and $x_-\in I$,
we can reduce the problem to the case that~$x$ is positive.
We may also assume that~$a$ is self-adjoint.
Now, since~$x\geq 0$ and $-\|a\|\leq a\leq \|a\|$,
we have $-\|a\|x \leq ax\leq \|a\|x$
by~\sref{sqrt},
and so~$ax\in I$,
because $\|a\|x\in I$.\qed
\end{point}
\end{point}
\begin{point}{100}[riesz-ideal-basic]{Exercise}%
Verify the following facts about Riesz ideals.
\begin{enumerate}
\item
The least Riesz ideal that contains a self-adjoint element~$a$
of~$\scrA$ is
\begin{equation*}
(a)_m\ :=\ \{\,b\in \scrA\colon\, 
\exists n\in \N\,[\ \left|\Real{b}\right|,\,\left|\Imag{b}\right|
\,\leq\, n\left|a\right| \ ]\,\}.
\end{equation*}
Moreover,  $(a)_m=\scrA$ iff $a$ is invertible,
and we have~$(a)=(a)_m$ when~$a\geq 0$
(where $(a)$ is the least order ideal that contains~$a$,
see~\sref{order-ideal-basic}).
For non-positive~$a$, however, we may have~$(a)\neq (a)_m$.
\item
$I+J$ is a Riesz ideal of~$\scrA$
when $I$ and~$J$ are Riesz ideals. (Hint: use~\sref{riesz-decomposition-lemma}.)
But~$I+J$ might not be an order ideal
when~$I$ and~$J$ are order ideals.

\item
Each proper Riesz ideal is contained in a maximal Riesz ideal.
\end{enumerate}%
\spacingfix%
\end{point}%
\begin{point}{110}[maximal-riesz-ideal-maximal-order-ideal]{Lemma}%
A maximal Riesz ideal~$I$ of~$\scrA$
is a maximal order ideal.
\begin{point}{120}{Proof}%
Let~$J$ be a proper order ideal with $I\subseteq J$.
We must show that $J=I$.
Let~$a\in J$ be given;
we must show that $a\in I$.
Since~$\Real{a},\Imag{a}\in J$,
it suffices to show that~$\Real{a},\Imag{a}\in I$,
and so we might as well assume that~$a$ is self-adjoint
to begin with.
Similarly,
since~$\left|a\right|\in J$,
and it suffices to show that~$\left|a\right|\in I$,
because then $-\left|a\right|\leq a\leq \left|a\right|$
entails $a\in I$,
we might as well assume that~$a$ is positive.

Note that the least ideal~$(a)$ that contains~$a$
is also a Riesz ideal by~\sref{riesz-ideal-basic}.
Hence  $I+(a)$ is a Riesz ideal by~\sref{riesz-ideal-basic}
Since~$a\in J$, we have $(a)\subseteq J$,
and so~$I+(a)\subseteq J$ is proper.
It follows that $a\in I+(a)=I$ by maximality of~$I$.\qed
\end{point}
\end{point}
\begin{point}{130}[riesz-ideal-miu-map]{Lemma}%
Let~$I$ be a maximal Riesz ideal of~$\scrA$.
Then there is an miu-map $f\colon \scrA\to \C$
with $\ker(f)=I$.
\begin{point}{140}{Proof}%
Since~$I$ is a maximal order ideal 
by~\sref{maximal-riesz-ideal-maximal-order-ideal},
there is a pu-map $f\colon \scrA\to \C$
with~$\ker(f)=I$ by~\sref{maximal-ideal-state}.
It remains to be shown that~$f$ is multiplicative.
Let~$a,b\in \scrA$ be given;
we must show that $f(ab)=f(a)f(b)$.
Surely, since~$f$ is unital,
we have $f(b-f(b))=f(b)-f(b)=0$,
an so $b-f(b)\in \ker(f)\equiv I$.
Now, since~$I$ is a Riesz ideal,
we have $a(b-f(b))\in I\equiv \ker(f)$ by~\sref{riesz-ideal-ring-ideal},
and so~$0=f(\,a(b-f(b))\,)=f(ab)-f(a)f(b)$.
Hence~$f$ is multiplicative.\qed
\end{point}
\end{point}
\begin{point}{150}[inv-mult-state]{Proposition}%
Let~$a$ be a self-adjoint element of a $C^*$-algebra.
Then~$a$ is not invertible
iff there is $f\in\spec(\scrA)$ 
with~$f(a)=0$.
\begin{point}{160}{Proof}%
Note that if~$a$ is invertible,
then~$f(a^{-1})$ is the inverse of~$f(a)$---and so~$f(a)\neq 0$---for 
every~$f\in\spec(\scrA)$.
For the other, non-trivial, direction,
assume that~$a$ is not invertible.
Then 
by~\sref{riesz-ideal-basic}
the least Riesz ideal $(a)_m$
that contains~$a$ is proper,
and can be extended to a maximal Riesz ideal~$I$.
By~\sref{riesz-ideal-miu-map}
there is an miu-map $f\colon \scrA\to\C$
with~$\ker(f)=I$.
Then~$f\in\spec(\scrA)$
and~$f(a)=0$.\qed
\end{point}
\end{point}
\begin{point}{170}[spectrum-miu]{Exercise}%
Show that $\spec(a)=\{f(a)\colon f\in\spec(\scrA)\}$
for each self-adjoint $a\in\scrA$.
\end{point}
\begin{point}{180}[gelfand-representation-isometry]{Exercise}%
Prove that $\|\gamma(a)\|=\|a\|$
for each~$a\in\scrA$
where~$\gamma$ is from~\sref{gelfand}.

(Hint: first assume that~$a$ is self-adjoint,
and use \sref{spectrum-miu} and~\sref{norm-spectrum}.
For the general case,
use the $C^*$-identity.)

Conclude that the Gelfand representation $\gamma\colon \scrA\to C(\spec(\scrA))$
is injective,
and that its range $\{\gamma(a)\colon a\in\scrA\}$
is a $C^*$-subalgebra of~$C(\spec(\scrA))$.
\end{point}
\begin{point}{190}%
To show that~$\gamma$ is surjective,
we use the following special case of
the Stone--Weierstra\ss{} theorem.%
\index{Stone--Weierstra\ss{}' Theorem}
\end{point}
\begin{point}{200}[stone-weierstrass]{Theorem}%
Let~$X$ be a compact Hausdorff space,
and let~$\scrS$ be a $C^*$-subalgebra of~$C(X)$
which `separates the points of~$X$',
that is, for all~$x,y\in X$ with~$x\neq y$
there is~$f\in \scrS$ with $f(x)\neq f(y)$.
Then~$\scrS=C(X)$.
\begin{point}{210}{Proof}%
Let~$g\in \pos{C(X)}$ and $\varepsilon >0$.
To prove that~$\scrS=C(X)$,
it suffices to show that~$g\in \scrS$,
and for this,
it suffices to find~$f\in \scrS$ with $\|f-g\|\leq \varepsilon$,
because~$\scrS$ is closed.
It is convenient to assume that~$g(x)> 0$ for all~$x\in X$,
which we may, without loss of generality,
by replacing~$g$ by~$1+g$.
\begin{point}{220}[stone-weierstrass-1]%
Let~$x,y\in X$ with~$x\neq y$
be given.
We know there is~$f\in \scrS$ with $f(x)\neq f(y)$.
Note that we can assume that~$f(x)=0$ (by replacing~$f$ by~$f-f(x)$),
and that~$f$ is self-adjoint (by replacing~$f$
by either~$\Real{f}$ or~$\Imag{f}$),
and that~$f$ is positive
(by replacing~$f$ by~$f_+$ or~$f_-$),
and that~$f(y)=g(y)>0$
(by replacing $f$ by $\frac{g(y)}{f(y)} f$),
and that~$f\leq g(y)$
(by replacing $f$ by $f\wedge g(y)$).
\end{point}
\begin{point}{230}%
Let~$y\in X$ be given.
We will show that there is~$f\in\scrS$
with $0\leq f\leq g+\varepsilon$
and~$f(y)=g(y)$.
Indeed,
since~$g$ is continuous
there is an open neighbourhood~$V$ of~$y$
with~$g(y) \leq  g(x)+\varepsilon$
for all~$x\in V$.
For each~$x\in X\backslash V$ there is $f_x \in [0,f(y)]_{\scrS}$
with $f_x(x)=0$ and~$f_x(y)=g(y)$ by~\sref{stone-weierstrass-1}.
Since the open subsets
$U_x := \{\,z\in X\colon f_x(z)\leq \varepsilon\,\}$
with~$x\in X\backslash V$
form an open cover of the closed (and thus compact) subset $X\backslash V$,
there are $x_1,\dotsc,x_N\in X\backslash U$
with $U_{x_1}\cup\dotsb\cup U_{x_N}\supseteq X\backslash V$.
Define $f:=f_{x_1}\wedge \dotsb \wedge f_{x_N}$.
Then~$f\in \scrS$, $0\leq f\leq g(y)$, $f(y)=g(y)$,
and $f(x)\leq \varepsilon$
for every~$x\in X\backslash V$.

We claim that $f\leq g+\varepsilon$.
Indeed,
if~$x\in X\backslash V$,
then $f(x)\leq \varepsilon\leq g(x)+\varepsilon$.
If~$x\in V$,
then $f(x)\leq g(y)\leq g(x)+\varepsilon$
(by definition of~$V$).
Hence $f\leq g+\varepsilon$.
\end{point}
\begin{point}{240}%
Thus for each~$y\in X$
there is $f_y\in \scrS$ with $0\leq f_y \leq g+\varepsilon$
and~$f_y(y)=g(y)$.
Since~$f_y$ is continuous at~$y$,
and~$f_y(y)=g(y)$,
there is an open neighbourhood~$U_y$ of~$y$
with $g(y)-\varepsilon\leq f_y(x)$
for all~$x\in U_y$.
Since these open neighbourhoods cover~$X$,
and~$X$ is compact,
there are $y_1,\dotsc,y_N\in X$
with $U_{y_1}\cup\dotsb\cup U_{y_N} = X$.
Define $f:=f_{y_1}\vee \dotsb\vee f_{y_N}$.
Then~$f\in\scrS$,
and $g-\varepsilon \leq f\leq g+\varepsilon$,
giving $\|f-g\|\leq \varepsilon$.\qed
\end{point}
\end{point}
\end{point}
\begin{point}{250}[spectrum-calg-compact-hausdorff]{Lemma}%
The spectrum $\spec(\scrA)$ of~$\scrA$ is a compact Hausdorff space.
\begin{point}{260}{Proof}%
Since for each~$a\in \scrA$
and~$f\in \spec(\scrA)$
we have  $\|f(a)\|\leq \|a\|$ 
by~\sref{norm-mi-map}
we see that~$f(a)$ is an element of the compact set
$\{\,z\in \C\colon\, \left|z\right|\leq\|a\|\,\}$,
and so~$\spec(\scrA)$ is a subset of
\begin{equation*}
\textstyle
\prod_{a\in \scrA}\, \{\,z\in \C\colon\, \left|z\right|\leq \|a\|\,\},
\end{equation*}
which is a compact Hausdorff space
(by Tychonoff's theorem, under the product topology
it inherits
from the space of all functions $\scrA\to \C$).
So to prove that~$\spec(\scrA)$
is a compact Hausdorff space,
it suffices to show that~$\spec(\scrA)$
is closed.
In other words,
we must show that if~$f\colon \scrA\to \C$
is the pointwise limit of a net of miu-maps $(f_i)_i$,
then~$f$ is an miu-map as well.
But this is easily achieved
using the continuity of addition, involution and multiplication on~$\C$,
because, for instance, 
for~$a,b\in\scrA$, we have $f(ab)
= \lim_i f_i(ab)=\lim_i f_i(a)f_i(b)
 = (\lim_i f_i(a))\,(\lim_i f_i(b))
= f(a) \,f(b)$.\qed
\end{point}
\end{point}
\begin{point}{270}[gelfand]{Gelfand's Representation Theorem}%
\index{Gelfand's Representation Theorem}%
\index{Cstar-algebra@$C^*$-algebra!commutative}
For a commutative $C^*$-algebra~$\scrA$,
the Gelfand representation, 
 $\gamma\colon \scrA\to C(\spec(\scrA))$
defined in~\sref{gelfand-representation}
is an miu-isomorphism.
\begin{point}{280}{Proof}%
We already know that~$\gamma$ is an injective miu-map
(see~\sref{gelfand-representation-basic} 
and~\sref{gelfand-representation-isometry}).
So to prove that~$\gamma$ is an miu-isomorphism,
it remains to be shown that~$\gamma$ is surjective.
Since~$\spec(\scrA)$ is a compact Hausdorff space 
(by~\sref{spectrum-calg-compact-hausdorff}),
and~$\gamma(\scrA)\equiv \{\gamma(a)\colon a\in \scrA\}$
is a $C^*$-subalgebra of~$C(\spec(\scrA))$
(by~\sref{gelfand-representation-isometry}),
it suffices to show that~$\gamma(\scrA)$
separates the points of~$\spec(X)$
by~\sref{stone-weierstrass}.
This is obvious,
because
for~$f,g\in \spec(\scrA)$ with~$f\neq g$
there is~$a\in \scrA$ with~$f(a)\equiv \gamma(a)(f)
\neq \gamma(a)(g)\equiv g(a)$.\qed
\end{point}
\end{point}
\end{parsec}
\begin{parsec}{280}%
\begin{point}{10}%
While Gelfand's representation theorem
is a result about commutative $C^*$-algebras,
it tells us a lot about non-commutative $C^*$-algebras too,
via their commutative $C^*$-subalgebras.
\end{point}
\begin{point}{20}[functional-calculus]{Exercise}%
Let~$a$ be an element of a (not necessarily commutative)
$C^*$-algebra~$\scrA$.
We are going to use Gelfand's representation
theorem to define~$f(a)$
for every continuous map $f\colon \spec(a)\to\C$
whenever~$a$ is contained in some commutative $C^*$-algebra.
This idea is referred to as the \Define{continuous functional calculus}.%
\index{functional calculus@$f(a)$, continuous functional calculus}
\begin{enumerate}
\item
Show that there is a least $C^*$-subalgebra
$\Define{C^*(a)}$%
\index{Cstara@$C^*(a)$, $C^*$-subalgebra generated by~$a$}
of~$\scrA$
that contains~$a$.

Given~$b\in C^*(a)$
show that~$bc=cb$ for all~$c\in\scrA$
with~$ac=ca$.
\item
We call~$a\in\scrA$ \Define{normal}%
\index{normal!element of a $C^*$-algebra}
when~$C^*(a)$ is commutative.

Show that~$a$ is normal iff
$aa^*=a^*a$
iff~$\Real{a}\Imag{a}=\Imag{a}\Real{a}$.
\item
From now on assume~$a$ is normal
so that~$C^*(a)$ is commutative.

Show that~$j\colon\, \varrho\mapsto \varrho(a),\, \spec(C^*(a))\to \spec(a)$
is a continuous map.

Denoting the composition of
the miu-maps
\begin{equation*}
\xymatrix@C=4em{
	C(\spec(a))
	\ar[r]^-{f\mapsto f\circ j}
	&
	C(\spec(C^*(a)))
	\ar[r]^-{\cong,\, \text{\sref{gelfand}}}
	&
	C^*(a)
	\ar[r]^-{\text{inclusion}}
	&
	\scrA.
}
\end{equation*}
by~$\Phi$,
we write $\Define{f(a)}:=\Phi(f)$
for all~$f\in C(\spec(a))$.

We have hereby defined, for example, $a^\alpha$ when $a\geq 0$
and~$\alpha\in (0,\infty)$.

From the fact that~$\Phi$ is miu some properties
of~$f(a)$ can be derived.
Show, for example,
that~$a^\alpha a^\beta = a^{\alpha+\beta}$
for all~$\alpha,\beta\in(0,\infty)$
when~$a\geq 0$.
\item
Given $f\in C(\spec(a))$,
show that~$f(a)$ is the unique element of~$C^*(a)$
with 
\begin{equation*}
	\varphi(f(a))
	 \ =\ f(\varphi(a))
\end{equation*}
for all~$\varphi\in\spec(C^*(a))$.
\item
(\Define{Spectral mapping thm.})%
\index{Spectral Mapping Theorem}
Show that~$\spec(f(a))=f(\spec(a))$
for~$f\in C(\spec(a))$.
\item
Show that~$\spec(\varrho(a))\subseteq
\spec(a)$
and $\varrho(f(a))=f(\varrho(a))$
for every $f\in C(\spec(a))$ and  miu-map
$\varrho\colon \scrA\to\scrB$
into a $C^*$-algebra~$\scrB$.
\item
Given~$f\in C(\spec(a))$ and~$g\in C(f(\spec(a)))$
show that $g(f(a))=(g\circ f)(a)$.

Show that $(a^\alpha)^\beta = a^{\alpha\beta}$
for~$\alpha,\beta\in (0,\infty)$
and~$a\in\scrA_+$.
\end{enumerate}%
\spacingfix{}%
\end{point}%
\begin{point}{30}[sqrt-monotone]{Theorem}%
We have $0\leq a\leq b \implies a^\alpha \leq b^\alpha$
for all positive elements $a$ and~$b$
of a $C^*$-algebra~$\scrA$,
and $\alpha\in (0,1]$.
\begin{point}{40}{Proof}%
(Based on~\cite{pedersen1972}.)
Note that the result is trivial if~$a$ and~$b$ commute.

It suffices to show that
$(a+\frac{1}{n})^\alpha \leq (b+\frac{1}{n})^\alpha$
for all~$n$,
because~$(a+\frac{1}{n})^\alpha$
norm converges to~$a^\alpha$
as~$n\to\infty$.
In other words,
it suffices to prove $a^\alpha \leq b^\alpha$
under the additional assumption that~$a$ and~$b$ are invertible.
Note that~$a^0$ and~$b^0$ are defined for such
invertible~$a$ and~$b$,
because the
function~$(\,\cdot\,)^0\colon [0,1]\to\C$
is only discontinuous at~$0$.
Writing~$E$ for the set of all~$\alpha\in[0,1]$
for which $b\mapsto b^\alpha$ is monotone
on positive,
\emph{invertible} elements of~$\scrA$
we must prove that~$E\supseteq(0,1]$,
and we will in fact show that~$E=[0,1]$.
Since clearly~$0,1\in E$
it suffices to show that~$E$ is convex.
We'll do this by showing that~$E$
is closed,
and $\alpha,\beta\in E\implies \frac{1}{2}\alpha 
+ \frac{1}{2}\beta \in E$.
\begin{point}{50}{$E$ is closed}
Let~$b$ be a positive and invertible element of~$\scrA$.
A moment's thought reveals it suffices to 
prove that $\alpha\mapsto b^\alpha, \,[0,1]\to\scrA$
is continuous.
And indeed it is
being
the composition of the map $\alpha \mapsto b^\alpha\colon\,[0,1]
\to C(\spec(b))$, 
	which is norm continuous,
and the functional
	calculus $f\mapsto f(b)\colon\,
	C(\spec(b))\to \scrA$, which being an miu-map is norm continuous
	as well.
\end{point}
\begin{point}{60}{$\alpha,\beta\in E\implies \frac{1}{2}\alpha+\frac{1}{2}\beta
	\in E$}
Let~$\alpha,\beta\in E$. Let~$a,b\in\scrA$ be positive
and invertible with $a\leq b$.
We must show that $a^{\frac{\alpha+\beta}{2}}\leq 
b^{\frac{\alpha+\beta}{2}}$.
Since the map $b^{\frac{\alpha+\beta}{4}}(\,\cdot\,)
b^{\frac{\alpha+\beta}{4}}$
is positive (by~\sref{astara-pos-basic-consequences}),
it suffices to show that
$b^{-\frac{\alpha+\beta}{4}}\,a^{\frac{\alpha+\beta}{2}}\,
b^{-\frac{\alpha+\beta}{4}} \leq 1$,
that is, 
$\|b^{-\frac{\alpha+\beta}{4}}\,a^{\frac{\alpha+\beta}{2}}\,
b^{-\frac{\alpha+\beta}{4}} \| \leq 1$.

For this, it seems, we must take a look under the hood
of the theory of $C^*$-algebras:
writing $\varrho(c):=\sup_{\lambda\in \spec(c)}\left|\lambda \right|$
for~$c\in \scrA$,
we know that $\varrho(c)\leq \|c\|$ for any~$c$,
and $\varrho(c)=\|c\|$ for self-adjoint~$c$ by \sref{norm-spectrum}.
Moreover, recall from~\sref{prod-spec}
that $\spec(cd)\backslash\{0\}
=\spec(dc)\backslash\{0\}$,
and so~$\varrho(cd)=\varrho(dc)$
for all $c,d\in \scrA$.
Hence
\begin{alignat*}{3}
\|\,b^{-\frac{\alpha+\beta}{4}}\,a^{\frac{\alpha+\beta}{2}}\,
b^{-\frac{\alpha+\beta}{4}}\,\| 
\ &=\ 
\varrho(\,b^{-\frac{\alpha+\beta}{4}}\,a^{\frac{\alpha+\beta}{2}}\,
b^{-\frac{\alpha+\beta}{4}} \,) \\
&=\ 
\varrho(\,b^{-\frac{\alpha+\beta}{4}}\,a^{\frac{\alpha+\beta}{2}}\,
b^{-\frac{\alpha+\beta}{4}} \,b^{-\frac{\alpha-\beta}{4}}
\,b^{\frac{\alpha-\beta}{4}}\,) \\
&=\ 
\varrho(
\,b^{\frac{\alpha-\beta}{4}}\,
\,b^{-\frac{\alpha+\beta}{4}}\,a^{\frac{\alpha+\beta}{2}}\,
b^{-\frac{\alpha+\beta}{4}} \,b^{-\frac{\alpha-\beta}{4}}\,)\\
&=\ 
\varrho(
\,b^{-\nicefrac{\beta}{2}}\,a^{\nicefrac{\beta}{2}}\,
a^{\nicefrac{\alpha}{2}}\,
b^{-\nicefrac{\alpha}{2}} \,) \\
&\leq\ 
\|\,b^{-\nicefrac{\beta}{2}}\,a^{\nicefrac{\beta}{2}}\,\|
\,\|\,
a^{\nicefrac{\alpha}{2}}\,
b^{-\nicefrac{\alpha}{2}} \,\|\\
&=\ 
\|\,b^{-\nicefrac{\beta}{2}}\,a^{\beta}\,
b^{-\nicefrac{\beta}{2}}\,\|^{\nicefrac{1}{2}}
\,\|\,
b^{-\nicefrac{\alpha}{2}} \,a^\alpha\,
b^{-\nicefrac{\alpha}{2}} \,\|^{\nicefrac{1}{2}}\\
\ &\leq\ 
\|\,b^{-\nicefrac{\beta}{2}}\,b^{\beta}\,
b^{-\nicefrac{\beta}{2}}\,\|^{\nicefrac{1}{2}}
\,\|\,
b^{-\nicefrac{\alpha}{2}} \,b^\alpha\,
b^{-\nicefrac{\alpha}{2}} \,\|^{\nicefrac{1}{2}} \ = \ 1,
\end{alignat*}
and so we're done.\qed
\end{point}%
\end{point}%
\end{point}%
\end{parsec}%
\begin{parsec}{290}[gelfand-equivalence]%
\begin{point}{10}%
As a cherry on the cake,
we use Gelfand's representation theorem~\sref{gelfand}
to get an equivalence between the categories $\op{(\cCstar{miu})}$
and~$\Define{\CH}$%
\index{CH@$\CH$}
of continuous maps between compact Hausdorff spaces.

To set the stage,
we extend $X\mapsto C(X)$ to a functor
$\CH\to \op{(\cCstar{miu})}$
by sending a continuous function~$f\colon X\to Y$
to the miu-map $C(f)\colon C(Y)\to C(X)$
given by~$C(f)(g)=g\circ f$ for $g\in C(Y)$,
and we extend $\scrA\mapsto \spec(\scrA)$
to a functor $\spec\colon \op{(\cCstar{miu})}\to \CH$
by sending an miu-map $\varphi \colon \scrA\to\scrB$
to the continuous map~$\spec(\varphi)\colon \spec(\scrB)\to\spec(\scrA)$
given by~$\spec(\varphi)(f)=f\circ \varphi$.

The Gelfand representations $\gamma_\scrA\colon \scrA\to C(\spec(\scrA))$
form a natural isomorphism
from $C\circ \spec$ to the identity functor on~$\op{(\cCstar{miu})}$.
So to get an equivalence,
it suffices to find a natural isomorphism
from the identity on~$\CH$ to~$\spec\circ C$,
which is provided by the following lemma.
\end{point}
\begin{point}{20}{Lemma}%
Let~$X$ be a compact Hausdorff space,
and let~$\tau \colon C(X)\to \C$ be an miu-map.
Then there is~$x\in X$ with $\tau(f)=f (x)$
for all~$f\in C(X)$.
\begin{point}{30}{Proof}%
Define
$Z\,= \, \{\,x\in X\colon \ h(x)\neq 0\text{ for some~$h\in \pos{C(X)}$
with $\tau(h)=0$}\,\}$.
We'll prove~$X\backslash Z$ contains
exactly one point, $x_0$, and $\tau(f)=f(x_0)$ for all~$f$.
\begin{point}{40}
To see that~$X\backslash Z$ contains no more than one point,
let~$x,y\in X$ with $x\neq y$ be given;
we will show that either~$x\in Z$ or~$y\in Z$.
By the usual topological trickery,
we can find~$f,g\in \pos{C(X)}$
with $fg=0$, $f(x)=1$ and~$g(y)=1$.
Then~$0=\tau(fg)=\tau(f)\,\tau(g)$,
so either~$\tau(f)=0$ (and~$x\in Z$), or~$\tau(g)=0$
(and~$y\in Z$).

That~$X\backslash Z$ is non-empty
follows from the following result (by taking~$f=1$).
\end{point}
\begin{point}{50}[multiplicative-state-on-cx-1]%
For~$f\in \pos{C(X)}$
with~$f(x)> 0 \implies x\in Z$ for all~$x\in X$
we have~$\tau(f)=0$.
Indeed, for each~$x\in X$ with~$f(x)>0$
(and so~$x\in Z$)
we can find~$h\in \pos{C(X)}$
with $\tau(h)=0$ and~$h(x)\neq 0$.
Then~$f(x)< g(x)$
and~$\tau(g)=0$
for $g:=(\frac{f(x)}{h(x)}+1)h$.
By compactness,
we can find $g_1,\dotsc,g_N\in \pos{C(X)}$
with~$\tau(g_n)=0$,
such that for every~$x\in X$
there is~$n$ with $g(x)<f_n(x)$.
Writing $g:=g_1\vee \dotsb \vee g_N$,
we have $0\leq f\leq g$ and~$\tau(g)=0$
(because by~\sref{commutative-cstar-basic}
$\tau$ preserves finite infima).
It follows that~$\tau(f)=0$.
\end{point}
\begin{point}{60}%
We now know that~$X\backslash Z$ contains exactly
one point, say~$x_0$.
To see that~$\tau(f)=f(x_0)$
for~$f\in C(X)$,
write $g:=(f-f(x_0))^*(f-f(x_0))$
and note that $g(x)>0\implies x\neq x_0\implies  x\in Z$.
Thus by~\sref{multiplicative-state-on-cx-1},
we get $0=\tau(g)=\left|\tau(f)-f(x_0)\right|^2$,
and so $\tau(f)=f(x_0)$.\qed
\end{point}
\end{point}
\end{point}
\begin{point}{70}{Exercise}%
Let~$X$ be a compact Hausdorff space.
Show that for every~$x\in X$
the map $\delta_x\colon C(X)\to \C,\ f\mapsto f(x)$
is miu,
and that the map $X\to \spec(C(X)),\ x\mapsto \delta_x$
is a continuous bijection
onto a compact Hausdorff space,
and thus a homeomorphism.
\end{point}
\begin{point}{80}[injective-miu-isometry]{Exercise}%
\index{miu-map!injective!is isometry}
As an application of the equivalence
between $\op{(\cCstar{MIU})}$
and~$\CH$,
we will show that every injective miu-map
between $C^*$-algebras
is an isometry.

Show that an arrow $f\colon X\to Y$
in~$\CH$ is mono iff injective, and epi iff surjective
(using complete regularity of~$Y$).
Conclude that~$f$ is both epi and mono in~$\CH$
only if~$f$ is an isomorphism (i.e.~homeomorphism).

Let~$\varrho\colon \scrA\to\scrB$
be an injective miu-map between $C^*$-algebras.
Let~$a$ be a self-adjoint element of~$\scrA$.
Show that~$\varrho$ can be restricted
to an miu-map $\sigma\colon C^*(a)\to C^*(\varrho(a))$,
which is both epi and mono in~$\cCstar{MIU}$.
Conclude that~$\sigma$ is an isomorphism,
and thus~$\|\varrho(a)\|=\|a\|$.
Use the $C^*$-identity
to extend the equality $\|\varrho(a)\|=\|a\|$ 
to (not necessarily self-adjoint) $a\in \scrA$.
\end{point}
\begin{point}{90}[injective-miu-iso-on-image]{Exercise}%
Let~$\varrho\colon \scrA\to\scrB$ 
be an injective miu-map.
Show that~$\varrho(\scrA)$
is closed (using~\sref{injective-miu-isometry}).
Conclude that~$\varrho(\scrA)$
is a
 $C^*$-subalgebra of~$\scrB$
isomorphic to~$\scrA$.
\end{point}
\end{parsec}
\subsection{Representation by Bounded Operators}
\begin{parsec}{300}%
\begin{point}{10}[completion-inner-product-space]%
Let us prove that every $C^*$-algebra~$\scrA$
is isomorphic
to a $C^*$-algebra
of bounded operators on some Hilbert space.
We proceed as follows.
To each p-map $\omega\colon \scrA\to\C$
(see~\sref{maps})
we assign a inner product $[\,\cdot\,,\,\cdot\,]_\omega$ on~$\scrA$,
which can be ``completed'' to a Hilbert space $\scrH_\omega$.
Every element~$a\in \scrA$ gives a bounded operator on~$\scrH_\omega$
via the action $b\mapsto ab$, which in turn gives a 
miu-map $\varrho_\omega\colon \scrA\to \scrB (\scrH_\omega)$.
In general $\varrho_\omega$ is not injective,
but if~$\Omega$ is a set of p-maps which separates the
points of~$\scrA$,
then the composition
\begin{equation*}
	\xymatrix@C=6em{
		\scrA\ar[r]^-{\left<\varrho_\omega\right>_{\omega\in \Omega}}
		&
		\bigoplus_{\omega\in\Omega} \scrB(\scrH_\omega)
		\ar[r]
		&
		\scrB(\,\bigoplus_{\omega\in\Omega}\scrH_\omega\,)
	}
\end{equation*}
does give an injective miu-map~$\varrho$,
which restricts to an isomorphism 
(\sref{injective-miu-iso-on-image})
from~$\scrA$
to the $C^*$-algebra~$\varrho(\scrA)$
of bounded operators
on $\bigoplus_{\omega\in \Omega} \scrH_\omega$,
see~\sref{hilb-sum}.

The creation of~$\varrho_\omega$ from~$\omega$
is known as the \emph{Gelfand--Naimark--Segal (GNS) construction}
and will make a reappearance in the theory of von Neumann algebras
(in~\sref{normal-functionals-lemma}).

We take a somewhat utilitarian stance towards the GNS construction here,
but there is much more that can be said about it:
in the first chapter of my twin brother's thesis, \cite{bas},
you'll see that the GNS construction has a certain universal property,
and that it can be generalised to apply
not only to maps of the form $\omega\colon \scrA\to\C$,
but also to maps of the form $\varphi\colon \scrA\to\scrB$.
\end{point}
\begin{point}{20}[state-inner-product]{Lemma}%
For every p-map~$\omega\colon \scrA\to \C$ on a
$C^*$-algebra~$\scrA$,
$\Define{[a,b]_\omega} = \omega(a^*b)$
defines an inner product~$\Define{[\,\cdot\,,\,\cdot\,]_\omega}$%
\index{*innerprodomega@$[\,\cdot\,,\,\cdot\,]_\omega$, given np-functional $\omega$}
on~$\scrA$
(see~\sref{hilb-def}).
\begin{point}{30}{Proof}%
Note that $[a,a]_\omega\equiv \omega(a^*a)\geq 0$ for each~$a\in\scrA$,
because $a^*a\geq 0$ (by~\sref{astara-positive});
and  $\smash{\overline{[a,b]}_\omega}=[b,a]_\omega$
for $a,b\in\scrA$,
because $\omega$ is involution preserving (by~\sref{cstar-p-implies-i}).
Finally, it is clear that $[a,\,\cdot\,]_\omega\equiv\omega(a^*\,\cdot\,)$
is linear for each~$a\in\scrA$.\qed
\end{point}
\end{point}
\begin{point}{40}[omega-norm-basic]{Exercise}%
Let~$\omega\colon \scrA\to\C$
be a p-map on a $C^*$-algebra.
Let us for a moment study
the 
semi-norm
$\Define{\|\,\cdot\,\|_\omega}$%
\index{*seminormomega@$\|\,\cdot\,\|_\omega$, given np-functional $\omega$}
	on~$\scrA$
induced by the inner product $[\,\cdot\,,\,\cdot\,]_\omega$
(so~$\smash{\|a\|_\omega = \omega(a^*a)^{\nicefrac{1}{2}}}$),
because it plays an important role
here,
and all throughout the next chapter.
\begin{enumerate}
\item
Use Cauchy--Schwarz
(\sref{inner-product-basic})
to prove \Define{Kadison's inequality}%
\index{Kadison's inequality}: 
for all~$a,b\in\scrA$,
\begin{equation*}
\left|\omega(a^*b)\right|^2\ \leq\ \omega(a^*a)\ \omega(b^*b).
\end{equation*}
\item
Show that $\|ab\|_\omega \leq \|\omega\| \,\|a\|\,\|b\|_\omega$
for all $a,b\in\scrA$
(using $a^*a\leq \|a\|^2$).

Show that we do \emph{not} always have 
$\|ab\|_\omega\leq \|\omega\|\|a\|_\omega \|b\|$.

(Hint:
take $a=(\begin{smallmatrix}0&0\\0&1\end{smallmatrix})$
and $b=\frac{1}{2}(\begin{smallmatrix}1&1\\1&1\end{smallmatrix})$
from~$\scrA=M_2$,
and $\omega(\,(\begin{smallmatrix}c & d\\e&f\end{smallmatrix})\,)=c$.)

Show that neither always
$\|ab\|_\omega \leq \|a\|_\omega \|b\|_\omega$,
or
$\|a^*a\|_\omega = \|a\|^2_\omega$.

(Hint: 
take~$a=b=\frac{1}{2}(\begin{smallmatrix}1 & 1 \\ 1 & 1\end{smallmatrix})$
from~$\scrA= M_2$,
and 
$\omega((\,(\begin{smallmatrix}c&d\\e&f\end{smallmatrix})\,)=c$.)

Give a counterexample to $\|a^*\|_\omega = \|a\|_\omega$.
\end{enumerate}%
\spacingfix%
\end{point}%
\begin{point}{50}[inner-product-completion]{Exercise}%
\index{inner product!*C-valued@$\C$-valued!completion}
Let us begin by showing how a complex vector space~$V$
with inner product
$[\,\cdot\,,\,\cdot\,]$ can be ``completed'' to a Hilbert space~$\scrH$.

We will take for~$\scrH$ the set of Cauchy sequences on~$V$
modulo the following equivalence relation.
Two Cauchy sequences $(a_n)_n$ and~$(b_n)_n$ in~$V$
are considered equivalent
iff $\lim_n \|a_n-b_n\|=0$.
We ``embed'' $V$ into~$\scrH$ via the map $\eta\colon V\to \scrH$
which sends~$a$ to
the constant sequence $a,a,a,\dotsc$.
Note, however, that $\eta$ need not be injective:
show that $\eta(a)=\eta(b)$ iff $\|a-b\|=0$ for all $a,b\in V$.

Show that $d(\,(a_n)_n,\,(b_n)_n\,) = \lim_n \|a_n-b_n\|$
defines a metric on~$\scrH$,
that~$\scrH$ is complete with respect to this metric,
and that if $(a_n)_n$ is a Cauchy sequence in~$V$,
then $(\eta(a_n))_n$ converges to the \emph{element}~$(a_n)_n$ of~$\scrH$
(so $V$ is dense in~$\scrH$).

Show that every uniformly continuous 
map $f\colon V\to X$ to a complete metric space~$X$
can be uniquely extended to a uniformly continuous map $g\colon \scrH\to X$.
(We say that~$g$ extends~$f$ when $f=g\circ \eta$.)

Show that addition, scalar multiplication, and inner product on~$V$
(being uniformly continuous)
can be uniquely extended to uniformly continuous operations on~$\scrH$,
and turn~$\scrH$ into a Hilbert space.
(Also verify that the extended inner product agrees with the complete
metric we've already put on~$\scrH$.)

Show that every bounded linear map $f\colon V\to\scrK$
to a Hilbert space~$\scrK$
can be uniquely extended to a bounded linear map $g\colon \scrH\to\scrK$.

(Categorically speaking,
Hilbert spaces
form a reflexive subcategory of
the category of bounded linear maps between
complex vector spaces with an inner product.)
\end{point}
\begin{point}{60}[gns]{Definition (Gelfand--Naimark--Segal construction)}%
\index{Gelfand--Naimark--Segal (GNS)}%
	\\
Let $\omega\colon \scrA\to\C$ be a p-map on a $C^*$-algebra~$\scrA$.

Let~$\Define{\scrH_\omega}$%
\index{Homega@$\scrH_\omega$} 
	denote the completion
of~$\scrA$ endowed with the inner product $[\,\cdot\,,\,\cdot\,]_\omega$
(see~\sref{state-inner-product})
to a Hilbert space as discussed in~\sref{inner-product-completion}.
Recall that we have an ``embedding''
$\Define{\eta_\omega}\colon \scrA\to\scrH_\omega$%
\index{etaomega@$\eta_\omega$}
with $\left<\eta_\omega(a),\eta_\omega(b)\right>
= [a,b]_\omega$ for all~$a,b\in \scrA$.

Since given~$a\in \scrA$
the map $b\mapsto ab,\ \scrA\to\scrA$ is
bounded with respect to~$\|\,\cdot\,\|_\omega$
(because $\|ab\|_\omega\leq \|\omega\|\|a\|\|b\|_\omega$
by~\sref{omega-norm-basic}),
it can be uniquely extended to a bounded linear map
$\scrH_\omega\to\scrH_\omega$
(by the universal property of~$\scrH_\omega$, 
see~\sref{inner-product-completion}),
which we'll denote by~$\Define{\varrho_\omega}(a)$.%
\index{rhoomega@$\varrho_\omega$}
So~$\varrho_\omega(a)$ is the unique
bounded linear map $\scrH_\omega\to\scrH_\omega$
with $\varrho_\omega(a)(\eta_\omega(b)) = \eta_\omega(ab)$
for all~$b\in\scrA$.
\end{point}
\begin{point}{70}{Proposition}%
The map $\varrho_\omega\colon \scrA\to\scrB(\scrH_\omega)$
given by~\sref{gns} is an miu-map.
\begin{point}{80}{Proof}%
Let~$a_1,a_2\in\scrA$ be given.
Since $\varrho_\omega(a_1+a_2)\,\eta_\omega(b)
= \eta_\omega((a_1+a_2)b)
= \eta_\omega(a_1b)+\eta_\omega(a_2b)
= (\varrho_\omega(a_1) + \varrho_\omega(a_2))\,\eta_\omega(b)$
for all $b\in\scrA$,
and~$\{\eta_\omega(b)\colon b\in\scrA\}$
is dense in~$\scrH_\omega$,
we see that $\varrho_\omega(a_1+a_2)
=\varrho_\omega(a_1)+\varrho_\omega(a_2)$.
Since similarly $\varrho_\omega(\lambda a)
= \lambda\varrho_\omega(a)$
for $\lambda\in\C$ and~$a\in\scrA$,
we see that~$\varrho_\omega$ is linear.

Since $\varrho_\omega(1)\,\eta_\omega(b)
= \eta_\omega(b)$ for all~$b\in\scrA$,
we have $\varrho_\omega(1)\,x=x$
for all~$x\in\scrH_\omega$,
and so~$\varrho_\omega$ is unital,
$\varrho_\omega(1)=1$.

To see that~$\varrho_\omega$ is multiplicative,
note that
$(\varrho_\omega(a_1)\,\varrho_\omega(a_2))\,\eta_\omega(b)
= \eta_\omega(a_1a_2b)=\varrho_\omega(a_1a_2)\,\eta_\omega(b)$
for all~$a_1,a_2,b\in\scrA$.

Let~$a\in\scrA$ be given.
To show that~$\varrho_\omega$ is involution preserving
it suffices to prove that~$\varrho_\omega(a^*)$
is the adjoint of~$\varrho_\omega(a)$.
Since~$\left<\varrho_\omega(a^*)\,\eta_\omega(b),\eta_\omega(c)\right>
\equiv [a^*b,c]_\omega = \omega(b^*ac)=[b,ac]_\omega
\equiv \left<\eta_\omega(b),\varrho_\omega(a)\,\eta_\omega(c)\right>$
for all~$b,c\in\scrA$,
and~$\{\eta_\omega(b)\colon b\in\scrA\}$
is dense in~$\scrH_\omega$,
we get~$\left<\varrho_\omega(a^*)x,y\right>=\left<x,\varrho_\omega(a)y\right>$
for all~$x,y\in\scrH_\omega$,
and so~$\varrho_\omega(a^*)=\varrho_\omega(a)^*$.\qed
\end{point}
\end{point}
\begin{point}{90}[gelfand-naimark-representation]{Definition}%
Given a collection~$\Omega$ of p-maps $\omega\colon \scrA\to\C$
on a $C^*$-algebra~$\scrA$,
let $\Define{\varrho_\Omega}\colon \scrA\to \scrB(\scrH_\Omega)$%
\index{rhoOmega@$\varrho_\Omega$}
be the miu-map given by~$\varrho_\Omega(a)x 
= \sum_{\omega\in\Omega} \varrho_\omega(a)x(\omega)$,
where~$\Define{\scrH_\Omega}=\bigoplus_{\omega\in\Omega}\scrH_\omega$%
\index{HOmega@$\scrH_\Omega$}
(and $\varrho_\omega$ is as in~\sref{gns}).
\end{point}
\begin{point}{100}[proto-gelfand-naimark]{Proposition}%
For a collection~$\Omega$ of positive maps $\scrA\to \C$
on a $C^*$-algebra~$\scrA$,
the following are equivalent.
\begin{enumerate}
\item
\label{proto-gelfand-naimark-1}
$\varrho_\Omega\colon \scrA\to\scrB(\scrH_\Omega)$
is injective;
\item
\label{proto-gelfand-naimark-2}
$\Omega$ is centre separating on~$\scrA$
(see~\sref{separating});
\item
\label{proto-gelfand-naimark-3}
$\Omega'=\{\,\omega(b^*(\,\cdot\,)b)\colon \, b\in\scrA,\,\omega\in\Omega\,\}$
is order separating on~$\scrA$.
\end{enumerate}
In that case, $\varrho_\Omega(\scrA)$ is a $C^*$-subalgebra
of~$\scrB(\scrH_\Omega)$,
and~$\varrho_\Omega$
restricts to an miu-isomorphism from~$\scrA$ to~$\varrho_\Omega(\scrA)$.
\begin{point}{110}{Proof}%
It is clear that~\ref{proto-gelfand-naimark-3}
entails~\ref{proto-gelfand-naimark-2}.
\begin{point}{120}{\ref{proto-gelfand-naimark-2}$\Longrightarrow$%
\ref{proto-gelfand-naimark-1}}%
Let~$a\in \scrA$ with $\varrho_\Omega(a)=0$ be given.
We must show that~$a=0$ (in order to show that~$\varrho_\Omega$
is injective),
and for this it is enough to prove that~$a^*a=0$.
Let~$b\in\scrA$ and~$\omega\in\Omega$ be given.
Since~$\Omega$ is centre separating,
it suffices to show that $0=\omega(b^*a^*ab) \equiv \|ab\|_\omega^2$.
Since~$\varrho_\Omega(a)=0$,
we have $\varrho_\omega(a)=0$,
thus $0=\varrho_\omega(a)\,\eta_\omega(b)
=\eta_\omega(ab)$,
and so $\|ab\|_\omega=0$.
Hence~$\varrho_\Omega$ is injective.
\end{point}
\begin{point}{130}{\ref{proto-gelfand-naimark-1}$\Longrightarrow$%
\ref{proto-gelfand-naimark-3}}%
Let~$a\in\scrA$ with $\omega(b^*a b)\geq 0$
for all~$\omega\in\Omega$ and~$b\in\scrA$
be given.
We must show that~$a\geq 0$.
Since~$\varrho_\Omega$ is injective,
we know by~\sref{injective-miu-iso-on-image}
that~$\varrho_\Omega(\scrA)$ is a $C^*$-subalgebra
of~$\scrB(\scrH_\Omega)$,
and~$\varrho_\Omega$ restricts to an miu-isomorphism
from~$\scrA$ to~$\varrho_\Omega(\scrA)$.
So in order to prove that~$a\geq 0$,
it suffices to show that $\varrho_\Omega(a)\geq 0$,
and for this we must prove that $\varrho_\omega(a)\geq 0$
for given $\omega\in \Omega$.
Since the vector states on~$\scrH_\omega$ are order separating
by~\sref{hilb-vector-states-order-separating}, it suffices to show that 
$\left<x,\varrho_\omega(a)x\right>\geq 0$
for given~$x\in \scrH_\omega$.
Since~$\{\eta_\omega(b)\colon b\in\scrA\}$
is dense in~$\scrH_\omega$,
we only need to prove 
that~$0\leq \left<\eta_\omega(b),\varrho_\omega(a)\eta_\omega(b)\right>
\equiv \omega(b^*ab)$ for given~$b\in \scrA$,
but this is true
by assumption.\qed
\end{point}
\end{point}
\end{point}
\begin{point}{140}[gelfand-naimark]{Theorem (Gelfand--Naimark)}%
\index{Gelfand--Naimark's Theorem}
Every $C^*$-algebra~$\scrA$ is miu-isomorphic
to a $C^*$-algebra of operators on a Hilbert space.
\begin{point}{150}{Proof}%
Since the states on~$\scrA$
are separating
(\sref{states-order-separating}),
and therefore centre separating,
the miu-map $\varrho_\Omega\colon \scrA\to\scrB(\scrH_\Omega)$
(defined in~\sref{gelfand-naimark-representation})
restricts to an miu-isomorphism
from~$\scrA$ onto the $C^*$-subalgebra
$\varrho(\scrA)$ of~$\scrB(\scrH_\Omega)$
by~\sref{proto-gelfand-naimark}.\qed
\end{point}
\end{point}
\end{parsec}
\section{Matrices over $C^*$-algebras}
\begin{parsec}{310}%
\begin{point}{10}%
We have seen (in~\sref{hilb}) that the $N\times N$-matrices
($N$ being a natural number) over the complex numbers~$\C$
form a $C^*$-algebra (denoted by~$M_N$) by interpreting
them as bounded operators on the Hilbert space $\C^N$,
and proving
that the bounded operators~$\scrB(\scrH)$
on any Hilbert space~$\scrH$ form a $C^*$-algebra.

In this paragraph, we'll prove the analogous
and more general
result that the 
$N\times N$-matrices \emph{over a $C^*$-algebra~$\scrA$}
form a $C^*$-algebra by interpreting them
as \emph{adjointable module maps} on
the \emph{Hilbert $\scrA$-module} $\scrA^N$,
see~\sref{chilb-basic} and~\sref{bax-cstar}.
\end{point}
\end{parsec}
\begin{parsec}{320}%
\begin{point}{10}[chilb-basic]{Definition}%
	An ($\scrA$-valued)
	\Define{inner product}%
\index{inner product!$\scrA$-valued}
	on a right $\scrA$-module~$X$
($\scrA$ being a $C^*$-algebra) is a map
$\left<\,\cdot\,,\,\cdot\,\right>\colon X\times X\to\scrA$%
\index{$\left<\,\cdot\,,\,\cdot\,\right>$, inner product!*avalued@$\scrA$-valued}
such that, for all $x,y\in X$,
$\left<x,\,\cdot\,\right>\colon X\to \scrA$
is a module map, $\left<x,x\right>\geq 0$,
and $\left<x,y\right>=\left<y,x\right>^*$.
We say that such an inner product is \Define{definite}%
	\index{$\left<\,\cdot\,,\,\cdot\,\right>$, inner product!*avalued@$\scrA$-valued!definite}
if~$\left<x,x\right>=0\implies x=0$ for all~$x\in X$.

A \Define{pre-Hilbert $\scrA$-module}%
\index{pre-Hilbert $\scrA$-module}
$X$ (where~$\scrA$ is always assumed to be a $C^*$-algebra)
is a right $\scrA$-module endowed with a definite inner product.
Such~$X$ is called
a \Define{Hilbert $\scrA$-module}%
\index{Hilbert $\scrA$-module}
when it is complete
with respect to
the norm we'll define in~\sref{chilb-norm-basic}.

Let~$X$ and~$Y$ be pre-Hilbert $\scrA$-modules.
We say that a map $T\colon X\to Y$
is adjoint to a map $S\colon Y\to X$
when
\begin{equation*}
\left<Tx,y\right>\ =\ \left<x,Sy\right>
\qquad \text{for all $x\in X$ and $y\in Y$}.
\end{equation*}
In that case, we call~$T$ \Define{adjointable}.%
\index{adjointable!map between pre-Hilbert $\scrA$-modules}
It is not difficult to see that~$T$
must be linear, and a module map, and 
adjoint to exactly one~$S$, which we denote by~$\Define{T^*}$.%
\index{adjoint!of a adjointable map between pre-Hilbert $\scrA$-modules}

(Note that we did not require that~$T$
is bounded, and in fact, it need not be, 
see~\sref{hellinger-toeplitz-needs-complete}.
However, if~$T$ is bounded, then so is~$T^*$, 
see~\sref{chilb-form-bounded},
and if either~$X$ or~$Y$ is complete,
then~$T$ is automatically bounded, see~\sref{hellinger-toeplitz}.)

The vector space of adjointable bounded module maps~$T\colon X\to Y$ 
is denoted
by~$\Define{\scrB^a(X,Y)}$,%
\index{BaXY@$\scrB^a(X,Y)$}
and we write $\Define{\scrB^a(X)}=\scrB^a(X,X)$.%
\index{BaX@$\scrB^a(X)$}
\end{point}
\begin{point}{20}{Example}%
We endow $\scrA^N$
(where~$\scrA$ is a $C^*$-algebra and~$N$ is a natural number)
with the inner product $\left<x,y\right>=\sum_n x_n^*y_n$,
making it a Hilbert $\scrA$-module.
\end{point}
\begin{point}{30}{Exercise}%
Let~$S$ and~$T$ be adjointable operators on a 
pre-Hilbert $\scrA$-module.
\begin{enumerate}
\item
	Show that~$T^*$ is adjoint to~$T$ (and so~$T^{**}=T$).
\item
Show that $(T+S)^*=T^*+S^*$ 
and $(\lambda S)^*=\overline{\lambda}S^*$ for $\lambda\in \C$.
\item
Show that $ST$ is adjoint to $T^*S^*$
(and so $(ST)^*=T^*S^*$).
\end{enumerate}%
\spacingfix%
\end{point}%
\begin{point}{40}{Exercise}%
Although a bounded linear map between Hilbert spaces
is always adjointable (see~\sref{hilb-adjoint}),
a bounded module map
between Hilbert $\scrA$-modules
might have no adjoint
as is revealed by the following example
(based on~\cite{paschke}, p.~447).

Prove that~$J:=\{\,f\in C[0,1]\colon\, f(0)=0\,\}$
is a closed right ideal of~$C[0,1]$, and thus a
Hilbert $C[0,1]$-module.

Show that the inclusion $T\colon J\to C[0,1]$
is a bounded module map,
which has no adjoint
by proving that there is no~$b\in J$
with $\left<b,a\right>=Ta\equiv a$ for all~$a\in J$
(for if~$T$ had an adjoint~$T^*$,
then $\left<T^*1,a\right>=\left<1,Ta\right>=a$
for all~$a\in J$).
\begin{point}{50}{Remark}%
Note that part of the problem here is the lack 
of the obvious analogue to
Riesz'~representation theorem (\sref{riesz-representation-theorem})  
for Hilbert $\scrA$-modules.
One solution (taken in the literature) is to simply 
add Riesz'~representation theorem as axiom
giving us the \emph{self-dual} Hilbert $\scrA$-modules.
For those
who like to keep Riesz'~representation theorem a theorem,
I'd like to mention that 
it is also possible to assume instead that the Hilbert $\scrA$-module
is complete with respect to a suitable uniformity,
as in done in my twin brother's thesis, \cite{bas}, see~\sref{dils-selfdual}.
\end{point}
\end{point}
\begin{point}{60}[chilb-cs]{Proposition (Cauchy--Schwarz)}%
\index{Cauchy--Schwarz inequality!for A-valued@for $\scrA$-valued inner products}
	We have
$\left<x,y\right>\,\left<y,x\right>
\,\leq\,\left\|\left<y,y\right>\right\|\,\left<x,x\right>$
for every inner product $\left<\,\cdot\,,\,\cdot\,\right>$
on a right $\scrA$-module~$X$,
and $x,y\in X$.
\begin{point}{70}{Remark}%
The symmetry-breaking norm symbols ``$\|$'' cannot simply 
be removed from this version of Cauchy--Schwarz,
because 
$0\leq \left<x,y\right>\,\left<y,x\right>
\leq \left<y,y\right>\left<x,x\right>$
would
imply
that $\left<y,y\right>\left<x,x\right>$
is positive, and self-adjoint,
and thus that $\left<y,y\right>$
and~$\left<x,x\right>$ commute,
which is not always the case.
\end{point}
\begin{point}{80}{Proof}%
Let~$\omega\colon \scrA\to \C$ be a state of~$\scrA$.
Since the states on~$\scrA$
are order separating (\sref{states-order-separating}), 
it suffices to show that
$\omega(\,\left<x,y\right>\,\left<y,x\right>\,)
\,\leq\,\left\|\left<y,y\right>\right\|\,\omega(\left<x,x\right>)$.
Noting that $(u,v)\mapsto \omega(\left<u,v\right>)$
is a complex-valued inner product on~$X$,
we compute
\begin{alignat*}{3}
	\omega&(\,\left<x,y\right>\,\left<y,x\right>\,)^2\\
\ &= \ 
\omega(\,\left<x,\,y\left<y,x\right>\right>\,)^2
\\
&\leq\ 
\omega(\left<x,x\right>)\ 
\omega(\,\left<\,y\left<y,x\right>,\, y\left<y,x\right>\,\right>\,)
\qquad
&&\text{by Cauchy--Schwarz, \sref{inner-product-basic}}
\\
&=\ 
\omega(\left<x,x\right>)\ 
\omega(\,\left<x,y\right> \,\left<y,y\right>\, \left<y,x\right>\,)
\\
&\leq\ 
\omega(\left<x,x\right>)\ 
\omega(\,\left<x,y\right>\left<y,x\right>\,)
\  \left\|\left<y,y\right>\right\|
\qquad
&&\text{since $\left<y,y\right>\leq \left\|\left<y,y\right>\right\|$.}
\end{alignat*}
It follows
(also when~$\omega(\,\left<x,y\right>\,\left<y,x\right>\,)=0$),
that 
\begin{equation*}
\omega(\,\left<x,y\right>\,\left<y,x\right>\,)\ \leq\ 
\left\|\left<y,y\right>\right\|\,
\omega(\left<x,x\right>),
\end{equation*}
and so we're done.\qed
\end{point}
\end{point}
\begin{point}{90}[chilb-norm-basic]{Exercise}%
Let~$X$ be a pre-Hilbert $\scrA$-module.
Verify that
\begin{enumerate}
	\item
$\Define{\|x\|} = \left\|\left<x,x\right>\right\|^{\nicefrac{1}{2}}$
defines a norm~$\left\|\,\cdot\,\right\|$%
\index{$"\"|\,\cdot\,"\"|$, norm!on a pre-Hilbert $\scrA$-module}
on~$X$, and
\item
$\left\|xb\right\|\leq \left\|x\right\|\left\|b\right\|$
and $\left\|\left<x,y\right>\right\|\leq \left\|x\right\|
\left\|y\right\|$
for all~$x,y\in X$ and $b\in \scrA$.
\end{enumerate}%
\spacingfix%
\end{point}%
\begin{point}{100}[chilb-form-bounded]{Lemma}%
For a linear map~$T\colon X\to Y$
between pre-Hilbert $\scrA$-modules,
and $B>0$,
the following are equivalent.
\begin{enumerate}
\item 
\label{chilb-form-bounded-1}
$\|Tx\|\leq B\,\|x\|$ for all~$x\in X$
(that is, $T$ is bounded by~$B$);
\item
\label{chilb-form-bounded-2}
$\left\|\left<y,Tx\right>\right\|\leq B\,\|y\|\|x\|$
for all~$x\in X$, $y\in Y$.
\end{enumerate}
Moreover,
if~$T$ is adjointable,
and bounded, then $\|T^*\|=\|T\|$.
\begin{point}{110}{Proof}%
If~$\|Tx\|\leq B\|x\|$ for all~$x\in X$,
then~$T$ is bounded, $\|T\|\leq B$, and therefore
$\left\|\left<y,Tx\right>\right\|
\leq \|y\|\,\|Tx\|\leq B \|y\|\|x\|$
for all~$x\in X$ and~$y\in Y$ using~\sref{chilb-cs}.

On the other hand,
if~\ref{chilb-form-bounded-2} holds,
and~$x\in X$ is given,
then we have $\|Tx\|^2=\left\|\left<Tx,Tx\right>\right\|
\leq B \,\|Tx\|\|x\|$,
entailing $\|Tx\|\leq B\|x\|$
(also when~$\|Tx\|=0$).

If~$T$ is adjointable, and bounded,
then~$\left\|\left<x,T^*y\right>\right\|=\left\|\left<y,Tx\right>\right\|
\leq \|T\|\|y\|\|x\|$ for all~$x\in X$, $y\in Y$,
so~$\|T^*\|\leq \|T\|$,
giving us that~$T^*$ is bounded.
Since by a similar reasoning $\|T\|\leq \|T^*\|$,
we get $\|T\|=\|T^*\|$.\qed
\end{point}
\end{point}
\begin{point}{120}[module-maps-cstar-identity]{Exercise}%
Show that $\|T^*T\|=\|T\|^2$
for every adjointable bounded map~$T$ on a pre-Hilbert $\scrA$-module.
(Hint: adapt the proof of~\sref{operators-cstar-identity}.)
\end{point}
\begin{point}{130}[bax-cstar]{Proposition}%
The adjointable bounded module maps
on a Hilbert $\scrA$-module
form a $C^*$-algebra%
\index{BaX@$\scrB^a(X)!as $C^*$-algebra}
$\scrB^a(X)$
with composition as multiplication,
adjoint as involution,
and the operator norm as norm.
\begin{point}{140}{Proof}%
Considering~\sref{bounded-operators-banach-algebra}
and~\sref{module-maps-cstar-identity},
the only thing that remains to be shown is that~$\scrB^a(X)$
is closed (with respect to the operator norm)
in the set of all bounded \emph{linear} maps $\scrB(X)$.
So let~$T\colon X\to X$ be a bounded linear map
which is the limit of a sequence $T_1,T_2,\dotsc$
of adjointable bounded module maps.

To see that~$T$ has an adjoint,
note that~$\left\|T_n^*-T_m^*\right\|
=\left\|(T_n-T_m)^*\right\|
=\left\|T_n-T_m\right\|$
for all~$n,m$, and so $T_1^*,\,T_2^*,\,\dotsc$
is a Cauchy sequence,
and converges to some bounded operator~$S$ on~$X$.
Since for~$x,y\in X$ and~$n$,
\begin{alignat*}{3}
\left\|\left<Sx,y\right>-\left<x,Ty\right>\right\|
\ &\leq\ 
\left\|\left<(S-T^*_n)x,y\right>\right\|
\,+\,
\left\|\left<x,(T_n-T)y\right>\right\|
\\
\ &\leq\ 
\|S-T^*_n\|\|x\|\|y\|\,+\,\|T_n-T\|\|x\|\|y\|,
\end{alignat*}
we see that $\left<Sx,y\right>=\left<x,Ty\right>$,
so~$S$ is the adjoint of~$T$,
and~$T$ is adjointable.
\qed
\end{point}
\end{point}
\begin{point}{150}[chilb-vector-states-order-separating]{Exercise}%
Let~$X$ be a Hilbert~$\scrA$-module.
Show that 
the vector states%
\index{vector functional!for a Hilbert $\scrA$-module}
of~$\scrB^a(X)$
are order separating (see~\sref{separating}).
Conclude that 
for an adjointable operator~$T$ on~$X$
\begin{enumerate}
\item
$T$ is self-adjoint iff $\left<x,Tx\right>$
is self-adjoint for all~$x\in (X)_1$;
\item
$0\leq T$ iff $0\leq\left<x,Tx\right>$
for all~$x\in (X)_1$;
\item
	$\|T\|=\sup_{x\in (X)_1}\|\left<x,Tx\right>\|$
when~$T$ is self-adjoint.
\end{enumerate}
(Hint:~adapt the proofs
of~\sref{hilb-vector-states-order-separating}
and~\sref{hilb-positive-operators}.)
\end{point}
\begin{point}{160}{Corollary}%
The operator
$T^*T$ is positive
in~$\scrB^a(X)$
for every adjointable operator~$T\colon X\to Y$
between Hilbert $\scrA$-modules.
\begin{point}{170}{Proof}%
$\left<x,T^*Tx\right>=
\left<Tx,Tx\right> \geq 0$
for all~$x\in X$,
and so~$T^*T\geq 0$ by~\sref{hilb-positive-operators}.\qed
\end{point}
\end{point}
\end{parsec}
\begin{parsec}{330}
\begin{point}{10}[cstar-matrices]{Exercise}%
Let us consider matrices over a $C^*$-algebra $\scrA$.
\begin{enumerate}
\item
	Show that every $N\times M$-matrix~$A$ (over~$\scrA$)
gives a bounded module map~$\underline{A}\colon \scrA^N\to\scrA^M$ 
via $\underline{A}(a_1,\dotsc,a_N)= A(a_1,\dotsc,a_N)$,
which is adjoint to~$\underline{A^*}$
(where $\Define{A^*}= (A_{ji}^*)_{ij}$ is conjugate transpose).

\item
Show that $A\mapsto \underline{A}$
gives a linear bijection between the vector 
space of $N\times M$-matrices 
over~$\scrA$ and the vector space of adjointable bounded
module maps~$\scrB^a(\scrA^N,\scrA^M)$.

\item
Show that~$\underline{A}\circ \underline{B} = \underline{AB}$
when $A$ is an $N\times M$ and~$B$ an $M\times K$ matrix.

\item
Conclude that the vector space $\Define{M_N\scrA}$%
\index{$M_n\scrA$, the $n\times n$-matrices over~$\scrA$!as a $C^*$-algebra}
of $N\times N$-matrices over~$\scrA$
is a $C^*$-algebra
with matrix multiplication (as multiplication),
conjugate transpose as involution,
and the operator norm (as norm, so~$\|A\|=\|\underline{A}\|$).
\end{enumerate}
\spacingfix%
\end{point}%
\begin{point}{20}[when-a-matrix-over-a-cstar-algebra-is-positive]{Exercise}%
Let us describe the positive  $N\times N$ matrices
over a $C^*$-algebra~$\scrA$.
\begin{enumerate}
\item
Show that an $N\times N$ matrix~$A$ over~$\scrA$
is positive iff $0\leq \sum_{i,j} a_i^* A_{ij} a_j$
for all~$a_1,\dotsc,a_N\in\scrA$.
(Hint: use~\sref{hilb-vector-states-order-separating}.)
\item
Show that the matrix $(\,\left<x_i,x_j\right>\,)_{ij}$
is positive for all vectors $x_1,\dotsc,x_N$
from a pre-Hilbert $\scrA$-module~$X$.
\item
Show that the matrix $(a^*_ia_j)_{ij}$
is positive for all $a_1,\dotsc,a_N\in\scrA$.
\end{enumerate}
\spacingfix
\end{point}%
\begin{point}{30}[mnf]{Exercise}%
Let~$f\colon \scrA\to\scrB$ be a linear map between $C^*$-algebras.
\begin{enumerate}
\item
Show that applying~$f$ entry-wise to an $N\times N$ matrix~$A$
over~$\scrA$ (yielding the matrix $(f(A_{ij}))_{ij}$ over~$\scrB$)
gives a linear map,
which we'll denote by~$\Define{M_Nf}\colon M_N\scrA\to M_N\scrB$.%
\index{$M_nf$}
\item
The map~$M_Nf$ inherits some traits of~$f$:
show that if~$f$ is unital, then~$M_Nf$ unital;
if~$f$ is multiplicative, then $M_Nf$ is multiplicative; and
if~$f$ is involution preserving, then so is~$M_Nf$.
\item
However,
show that $M_nf$ need not be positive when~$f$ is positive,
and that~$M_nf$ need not be bounded, when~$f$ is.
\end{enumerate}%
\spacingfix%
\end{point}%
\end{parsec}%
\begin{parsec}{340}%
\begin{point}{10}%
Let us briefly return
to the completely positive maps (defined in~\sref{maps}),
to show that a map $f$ between $C^*$-algebras
is completely positive precisely
when~$M_Nf$ is positive for all~$N$,
and to give some examples of completely positive maps.

We also prove two lemmas
stating special properties of completely positive maps (setting
them apart from plain positive maps),
that'll come in very handy later on.
The first one is a variation on Cauchy--Schwarz
(\sref{cp-cs}),
and the second one concerns
the points at which a cpu-map is multiplicative (\sref{choi}).

Completely positive maps are often touted as 
good models for quantum processes
(over plain positive maps)
with an argument involving the tensor product,
and while we agree,
we submit that the absence of analogues of \sref{cp-cs} and~\sref{choi}
for positive maps
is already enough to make complete positivity indispensable.
\end{point}
\begin{point}{20}[n-pos]{Lemma}%
\index{completely positive!map between $C^*$-algebras}
For a linear map $f\colon \scrA\to\scrB$
between $C^*$-algebras,
and natural number~$N$,
the following are equivalent.
\begin{enumerate}
\item
\label{n-pos-1}
$M_Nf\colon M_N\scrA\to M_N\scrB$
is positive;
\item
\label{n-pos-2}
	$\sum_{ij} b^*_if(a^*_ia_j)b_j \geq 0$
	for all~$a\equiv(a_1,\dotsc,a_N)\in \scrA^N$
	and $b\in \scrB^N$;
\item
\label{n-pos-3}
the matrix $(\,f(a_i^*a_j)\,)_{ij}$
is positive in $M_N\scrB$ for all $a\in\scrA^N$.
\end{enumerate}
\spacingfix%
\begin{point}{30}{Proof}%
Recall that~$M_Nf$ is positive
iff $(M_Nf)(C)$ is positive for all $C\in \pos{(M_N\scrA)}$.
The trick is to note that such~$C$ can be written as $C\equiv A^*A$
for some~$A\in M_N\scrA$,
and thus as $C \equiv (a_1^T)^* a_1^T+\dotsb+(a_N^T)^*a_N^T$,
where $a_n\equiv(A_{n1},\dotsc,A_{nN})$ is the $n$-th row of~$A$.
Hence~$M_Nf$ is positive
iff $(M_Nf)(\,(a^T)^*a^T\, )\equiv(\,f(a_i^*a_j)\,)_{i,j}$ is positive
for all tuples~$a\in\scrA^N$.
Since~$B\in M_N\scrB$ is positive iff $\left<b,Bb\right>\geq 0$
for all~$b\in \scrB^N$,
we conclude:
$M_Nf$ is positive iff 
$0\leq\left<b,(M_Nf)(\, (a^T)^*a^T\,) b\right>
= \sum_{ij} b_i^*f(a_i^*a_j)b_j$
for all~$a\in\scrA^N$ and~$b\in\scrB^N$.\qed
\end{point}
\end{point}
\begin{point}{40}[cp]{Exercise}%
Conclude from~\sref{n-pos}
that a linear map~$f$ between $C^*$-algebras
is completely positive iff~$M_Nf$ is positive for all~$N$
iff 
for all~$N$ and~$a\in \scrA^N$
the matrix $(\,f(a_i^*a_j)\,)_{i,j}$ 
is positive in~$M_N\scrB$.

Deduce that the composition of cp-maps is
completely positive.

Show that a mi-map~$f$ is completely positive.
(Hint: $M_Nf$ is a mi-map too.)
\end{point}
\begin{point}{50}[ad-cp]{Exercise}%
Show that
given a $C^*$-algebra~$\scrA$,
the following maps are completely positive:
\begin{enumerate}
\item
$b\mapsto a^*ba\colon \scrA\to\scrA$
for every~$a\in\scrA$;%
\index{$a^*(\,\cdot\,)a\colon \scrA\to\scrA$!is completely positive}
\item
$T\mapsto S^* T S\colon \scrB^a(X)
\to\scrB^a(Y)$
\index{$A^*(\,\cdot\,) A\colon \scrB^a(X)
\to\scrB^a(Y)$!is completely positive}
for every adjointable operator $S\colon Y\to X$
between Hilbert $\scrA$-modules;
\item
$T\mapsto \left<x,Tx\right>,\scrB^a(X)\to \scrA$%
		\index{vector functional!is completely positive}
for every element~$x$ of a Hilbert $\scrA$-module~$X$.
\end{enumerate}
\spacingfix%
\end{point}%
\begin{point}{60}[cstar-product-4]{Exercise}%
\index{product!in $\Cstar{cpsu}$}%
\index{equaliser!in $\Cstar{cpsu}$}%
Show that the product 
of a family of $C^*$-algebras $(\scrA_i)_i$
in the category~$\Cstar{cpsu}$
(see~\sref{maps}) 
is given by~$\bigoplus_i \scrA_i$
with the same projections as in~\sref{cstar-product-2}.

Show that the equaliser
of miu-maps $f,g\colon\scrA\to\scrB$
in~$\Cstar{cpsu}$
is the inclusion of
the $C^*$-subalgebra
$\{\,a\in\scrA\colon\, f(a)=g(a)\,\}$
of~$\scrA$ into~$\scrA$.
\end{point}
\begin{point}{70}[ccstar-pos-mat]{Lemma}%
Let~$\scrA$ be a commutative $C^*$-algebra,
and let~$N$ be a natural number.
The set of  matrices of the form $\sum_k a_k B_k$,
where $a_1,\dotsc,a_K\in \scrA_+$
and $B_1,\dotsc,B_K\in M_N(\C)_+$,
is norm dense in~$(M_N\scrA)_+$.
\begin{point}{80}{Proof}%
Since~$\scrA$ is isomorphic to~$C(X)$ for some compact
Hausdorff space~$X$ (by~$\sref{gelfand})$),
we may as well assume that~$\scrA\equiv C(X)$.

Let~$A\in M_N(C(X))_+$ and~$\varepsilon>0$ be given.
We're looking for $g_1,\dotsc,g_K\in C(X)_+$
and $B_1,\dotsc,B_K\in (M_N)_+$
with $\|A-\sum_k g_k B_k\|\leq \varepsilon$.
Since $A(x):=(A_{ij}(x))_{ij}$
gives a continuous map $X\to M_N$,
the sets
$U_x = \{\,y\in X\colon \, \|A(x)-  A(y)\| < \varepsilon\,\}$
form an open cover of~$X$.
By compactness of~$X$
this cover has a finite subcover;
there are $x_1,\dotsc,x_K\in X$ with
$U_{x_1}\cup\dotsb\cup U_{x_K}=X$.

Let~$y\in X$ be given. Since $y\in U_{x_k}$
for some~$k$, there is, by complete regularity of~$X$,
a function $f_y\in (C(X))_+$
with $f_y(y)>0$
and $\supp(f_y)\subseteq U_{x_k}$.
Since the open subsets~$\supp(f_y)$
cover~$X$
there are (by compactness of~$X$) finitely many $y_1,\dotsc y_L$
with $X = \supp(f_{y_1})\cup \dotsb \cup \supp(f_{y_L})$,
and so~$\sum_\ell f_{y_\ell} > 0$.
Let us group together the $f_{y_\ell}$s:
pick for each~$\ell$ an $k_\ell$ with $\supp(f_{y_\ell})\subseteq 
U_{x_{k_\ell}}$,
and let $g_k:= \sum\{f_\ell\colon k_\ell = k\}$.
Then $g_k\in (C(X))_+$,
$\supp(g_k)\subseteq U_k$,
and $\sum_k g_k >0$.
	Upon replacing $g_k$ with $(\sum_\ell g_\ell)^{-1} g_k$ if necessary,
we see that $\sum_k g_k=1$.

Since~$\supp(g_k)\subseteq U_{x_k}$,
we have $-\varepsilon \leq A(x)-A(x_k)\leq \varepsilon$
for all~$x\in \supp(g_k)$,
and so  $-g_k(x) \varepsilon
\,\leq\, g_k(x) A(x) - g_k(x) A(x_k)\,\leq\, g_k(x) \varepsilon$
for all~$x\in X$,
that is,  $-g_k \varepsilon
\,\leq\, g_k A - g_k A(x_k)\,\leq\, g_k \varepsilon$.
Summing yields
$-\varepsilon \,\leq\, A- \sum_k g_k A(x_k)\,\leq\, \varepsilon$,
and so $\|A-\sum_k g_k A(x_k)\|\leq \varepsilon$.\qed
\end{point}
\end{point}
\begin{point}{90}[cp-commutative]{Proposition}%
Let~$f\colon \scrA\to\scrB$ be a 
positive map between $C^*$-algebras.
If either~$\scrA$ or~$\scrB$ is commutative,
then~$f$ is completely positive.
\begin{point}{100}{Proof}%
Suppose that~$\scrB$ is commutative,
and let~$a_1,\dotsc,a_N\in \scrA$,
$b_1,\dotsc,b_N\in\scrB$ be given.
We must
show that $\sum_{i,j} b_i^*f(a_i^*a_j)b_j$ is positive.
This follows from the observation that
$\omega(\,\sum_{i,j} b_i^*f(a_i^*a_j)b_j\,)
= \omega(f(\,\sum_{i,j}(a_i\omega(b_i))^*\,a_j \omega(b_j)\,))\,\geq \,0$
for every~$\omega\in\spec(\scrB)$.
\begin{point}{110}%
Suppose instead that~$\scrA$ is commutative,
and let $A\in (M_N\scrA)_+$
be given for some natural number~$N$.
We must show that~$(M_Nf)(A)$ is positive in~$M_N\scrB$.
By~\sref{ccstar-pos-mat},
the problem reduces to the case that~$A\equiv a B$
where~$a\in \scrA_+$ and~$B\in (M_N)_+$.
Since $(M_Nf)(aB)\equiv f(a)B$ is clearly positive
in~$M_N\scrB$,
we are done.\qed
\end{point}
\end{point}
\end{point}
\begin{point}{120}[cstar-positive-2x2matrix]{Lemma}%
For a positive
matrix $A\equiv \bigl(\begin{smallmatrix} p & a \\ a^* & q
\end{smallmatrix}\bigr)$
over a $C^*$-algebra~$\scrA$
we have
\begin{equation*}
	a^*a\ \leq\  \|p\|q
	\quad\text{ and }\quad
	aa^*\leq \|q\|p.
\end{equation*}
In particular,
if $p=0$ or~$q=0$, then~$a=a^*=0$.
\begin{point}{130}{Proof}%
Since $(x,y)\mapsto \left<x,Ay\right>$
gives an $\scrA$-valued inner product
on~$\scrA^2$,
{
\newcommand\twovect[2]{%
\left(\begin{smallmatrix}#1\\#2\end{smallmatrix}\right)}
\newcommand\onezero{\twovect{1}{0}}
\newcommand\zeroone{\twovect{0}{1}}
\begin{alignat*}{3}
	aa^*
	\ &=\  \left<\,\onezero,\,A\zeroone\,\right> \ 
\left<\,\zeroone,\,A\onezero\,\right> \\
\ &\leq\  
\left\|\left<\,\zeroone,\,A\zeroone\right>\right\| \ 
\left<\,\onezero,\,A\onezero\,\right>
\ =\  \|q\|\  p
\end{alignat*}
}
by Cauchy--Schwarz (see \sref{chilb-cs}).

By a similar reasoning, we get $a^*a\leq \|p\|q$.\qed
\end{point}
\end{point}
\begin{point}{140}[cp-cs]{Lemma}%
We have $f(a^*b) f(b^*a)\leq \|f(b^*b)\|\,f(a^*a)$
for every p-map $f\colon \scrA\to\scrB$
between $C^*$-algebras
and $a,b\in\scrA$,
provided that $M_2f$ is positive.
\begin{point}{150}{Proof}%
Since writing $x\equiv (a,b)\in \scrA^2$,
the $2\times 2$ matrix $(x^T)^* x^T\equiv 
	\bigl(
\begin{smallmatrix}
a^*a & a^*b \\
b^*a & b^* b
\end{smallmatrix} \bigr)$
in $M_2\scrA$
is positive,
the $2\times 2$ matrix $T:=\bigl(
\begin{smallmatrix}
	f(a^*a) & f(a^*b) \\
	f(b^*a) & f(b^* b)
\end{smallmatrix}\bigr)$
in~$M_2\scrB$ is positive.
Thus we get $f(a^*b) f(b^*a)\leq \|f(b^*b)\|\,f(a^*a)$
by~\sref{cstar-positive-2x2matrix}.\qed
\end{point}
\end{point}
\begin{point}{160}[cp-russo-dye]{Corollary}%
$\|f\|= \|f(1)\|$
for every cp-map $f\colon \scrA\to\scrB$ between $C^*$-algebras.
\begin{point}{170}{Proof}%
Let~$a\in\scrA$ be given.
It suffices to show that $\|f(a)\|\leq \|f(1)\|\,\|a\|$
so that~$\|f\|\leq\|f(1)\|$,
because we already know that~$\|f(1)\|\leq \|f\|\,\|1\| = \|f\|$.
Since $\|f(a^*a)\|\leq \|f(1)\|\,\|a^*a\|$
by~\sref{weak-russo-dye},
we have
 $\|f(a)\|^2=\|f(a)^*f(a)\|=\|f(a^*1)f(1^*a)\|
\leq \|f(1^*1)\|\,\|f(a^*a)\|
\leq \|f(1)\|\, \|f(1)\|\|a^*a\|
= \|f(1)\|^2 \|a\|^2$
by~\sref{cp-cs},
and so~$\|f(a)\|\leq \|f(1)\|\,\|a\|$.\qed
\end{point}
\end{point}
\begin{point}{180}[choi]{Lemma (Choi\cite{choi})}%
\index{Choi's Theorem}%
We have
$f(a)^*f(a) \leq f(a^* a)$ for
every
cpu-map~$f\colon \scrA\to\scrB$ between $C^*$-algebras,
and~$a\in\scrA$.
Moreover, if $f(a^*a)=f(a)^*f(a)$
for some~$a\in\scrA$,
then~$f(ba)=f(b)f(a)$
for all~$b\in \scrA$.
\begin{point}{190}{Proof}%
By~\sref{cp-cs}
we have $f(a)^*f(a)=f(a^* 1)f(1^* a) \leq
\|f(1^*1)\| f(a^*a)=f(a^*a)$,
where we used that~$f$ is unital, viz.~$f(1)=1$.

Let~$a,b\in \scrA$ be given,
and assume that $f(a^*a)=f(a)^*f(a)$.
Instead of~$f(ba)=f(b)f(a)$
we'll prove that $f(a^*b)=f(a)^*f(b)$
(but this is nothing more than  a reformulation).
Since~$M_2f$ is cp,
we have, writing 
$A\equiv\bigl(\begin{smallmatrix}a&b\\0&0\end{smallmatrix}\bigr)$,
\begin{alignat*}{3}
\left(\,\begin{matrix}f(a)^*f(a)&f(a)^*f(b)\\
f(b)^*f(a)&f(b)^*f(b)\end{matrix}\,\right) 
	\ &=\ (M_2f)(A)^*\,(M_2f)(A)\\
\ &\leq\ (M_2f)(A^*A) \ =\ 
\left(\,\begin{matrix}f(a^*a)&f(a^*b)\\
f(b^*a)&f(b^*b)\end{matrix}\,\right).
\end{alignat*}
Hence
(using that $f(a^*a)=f(a)^*f(a)$)
the following matrix is positive.
\begin{equation*}
\left(\,\begin{matrix}
0 & f(a^*b) - f(a)^*f(b) \\
f(b^*a)-f(b)^*f(a) & f(b^*b)-f(b)^*f(b)
\end{matrix}\,\right)
\end{equation*}
But then
by~\sref{cstar-positive-2x2matrix}
we have
$f(a^*b)-f(a)^*f(b)=0$.\qed
\end{point}
\end{point}
\end{parsec}
\begin{parsec}{341}%
\begin{point}{10}%
We've just seen in~\sref{cp-russo-dye}
that the norm of a \emph{completely} positive map 
    $f\colon \scrA\to\scrB$ between $C^*$-algebras
is given by $\left\|f\right\|=\left\|f(1)\right\|$.
We'll show here that the same result holds
when~$f$ is just positive.
This result will be useful at the end of this
thesis in~\sref{lem:sef-instrument}, where we'll
try to consider the broadest possible class
of duplicators
$\delta\colon \scrA \otimes \scrA\to\scrA$
    (see~\sref{def:duplicator})
being a priori just positive,
not completely positive.
The proof consists of two ingredients:
    the fact, \sref{normal-russo-dye},
    that $\left\|f(a)\right\|\leq \|f(1)\| \|a\|$
    for all \emph{normal}~$a\in\scrA$ (see~\sref{functional-calculus}),
and the result,
known as Russo--Dye's theorem,
\sref{russo-dye},
that the convex combinations
of unitaries (being normal)
are norm dense in the unit ball~$(\scrA)_1$
of~$\scrA$.
\end{point}
\begin{point}{20}[normal-russo-dye]{Lemma}%
We have $\left\|f(a)\right\| \leq \left\|f(1)\right\|\,\left\|a\right\|$
for every p-map $f\colon \scrA\to\scrB$ between $C^*$-algebras,
and \emph{normal} $a\in\scrA$.
    \begin{point}{30}{Proof}%
Since~$a$ is normal,
the $C^*$-subalgebra $C^*(a)$
of~$\scrA$ generated by~$a$
is commutative (see~\sref{functional-calculus}),
and so the restriction of~$f$
to a map $f\colon C^*(a)\to\scrB$
is completely positive by~\sref{cp-commutative}.
Thus~$\|f(a)\|\leq \|f(1)\|\,\|a\|$ by~\sref{cp-russo-dye}.\qed
\end{point}
\end{point}
\begin{point}{40}[cstar-unitary]{Definition}%
An element~$u$ of a $C^*$-algebra
is \Define{unitary}\index{unitary!in a $C^*$-algebra}
when $u^*u=1$ and~$uu^*=1$.

    In that case we also say that~$u$ is \Define{\emph{a} unitary}.
\end{point}
\begin{point}{50}{Exercise}%
Let~$\scrA$ be a $C^*$-algebra.
\begin{enumerate}
\item
Show that any~$\lambda\in\C$ with $\left|\lambda\right|=1$
is unitary in~$\scrA$.

In particular, the unit, $1$, of~$\scrA$ is unitary.
\item
Show that a unitary~$u\in \scrA$ is invertible with inverse~$u^{-1}=u^*$,
and that~$u^*$ is a unitary as well.
\item
Show that the product $uv$ of unitaries $u,v\in\scrA$
is unitary.
\item
Show that every unitary~$u$ of~$\scrA$ is normal,
that is, $uu^*=u^*u$
(see~\sref{functional-calculus}).

Show that a normal element~$a$ of~$\scrA$ is unitary
        iff $\Real{a}^2 + \Imag{a}^2 = 1$.
\item
Show that every self-adjoint element~$a$ of~$\scrA$
with~$\|a\|\leq 1$
is the real part of some unitary~$u$, 
        so~$a=\Real{u}\equiv \frac{1}{2}(u+u^*)$.
(Hint:
        take~$u := a + i\sqrt{1-a^2}$.)
\item
Show that every invertible element~$a$ of~$\scrA$ can be written
as $a=u\sqrt{a^*a}$,
where~$u$ is a unitary.
(Hint: take $u=\sqrt{a^{-1}(a^{-1})^*}$.)

This is a variation on the polar decomposition
we'll see in~\sref{polar-decomposition}.
\end{enumerate}
\spacingfix%
\end{point}%
\begin{point}{60}{Exercise}
(Based on II.3.2.14--17 of~\cite{blackadar2006operator}.)
Let~$\scrA$ be a $C^*$-algebra.
\begin{enumerate}
\item
Show that every invertible element~$a$ of~$\scrA$
with~$\|a\|\leq 2$ can be written
as the sum of two unitaries.
(Hint: write $a=u\sqrt{a^*a}$ with~$u$ as above.)
\item
Let~$u\in \scrA$ be a unitary, and~$a\in\scrA$ with $\|a\|< 1$.

Show that~$u+a$ is the sum of two unitaries.

(Hint: write~$u+a=u(1+u^*a)$, and
        note that $1+u^*a$ is invertible by~\sref{geometric}.)
\item
Let~$a\in \scrA$ be given,
and let~$N$ be a natural number with $\|a\| < N$.

Show that~$a$ is the sum of~$N+2$ unitaries. 

(Hint: write $a=1+(N+1)b$ where $b:=\frac{a-1}{N+1}$,
and show that $\|b\|< 1$.)

\item
Prove the following theorem.
\end{enumerate}
\spacingfix%
\end{point}%
\begin{point}{70}[russo-dye]{Theorem (Russo--Dye)}%
\index{Russo--Dye's Theorem}
An element~$a$ of a $C^*$-algebra $\scrA$
with $\|a\|< 1-\frac{2}{N}$
for some natural number~$N$
can be written as~$a=\frac{1}{N}(u_1+\dotsb + u_N)$
for some unitaries $u_1,\dotsc,u_N\in \scrA$.
\end{point}
\begin{point}{80}[russo-dye-cor]{Corollary}%
The operator norm of a positive
    linear map~$f\colon\scrA\to\scrB$ between
    $C^*$-algebras is given by  $\|f\|=\|f(1)\|$.
\begin{point}{90}{Proof}%
We must show that $\|f(a)\|\leq \|f(1)\|$
for every~$a\in\scrA$ with~$\|a\|\leq 1$.
Since by Russo--Dye's theorem 
every~$a\in\scrA$ with~$\|a\|\leq 1$ 
may be approximated with respect to the norm by a sequence
of elements of the form  $b:=\frac{1}{N}(u_1+\dotsb+u_N)$,
where~$u_1,\dotsc,u_N$ are unitaries,
it suffices to show that $\|f(b)\|\leq \|f(1)\|$
for such~$b$.
Since~$u_n$ is normal,
    and thus $\|f(u_n)\|\leq \|f(1)\|\,\|u_n\|\leq \|f(1)\|$
    by~\sref{normal-russo-dye},
    we get $\|f(b)\|\leq \frac{1}{N}(\|f(u_1)\|\,+\,\dotsb\,+\,\|f(u_N)\|)
    \leq \|f(1)\|$, and so $\|f\|=\|f(1)\|$.\qed
\end{point}
\end{point}
\end{parsec}

\section{Towards von Neumann Algebras}
\begin{parsec}{350}%
\begin{point}{10}%
Let us work towards
the subject of the next chapter, von Neumann algebras,
by pointing out two special properties
of~$\scrB(\scrH)$
on which the definition of a von Neumann algebra is based,
namely that
\begin{enumerate}
\item
any norm-bounded directed subset of
self-adjoint operators on~$\scrH$
has a supremum (in~$\Real{\scrB(\scrH)}$), and
\item
all vector functionals
$\left<x,(\,\cdot\,)x\right>\colon \scrB(\scrH)\to\C$ 
preserve these suprema.
\end{enumerate}
We'll end
the chapter
by showing
in~\sref{bh-np}
 that every  functional on~$\scrB(\scrH)$
that  preserves the aforementioned
suprema
is a (possibly infinite) sum of vector functionals.
\end{point}

\subsection{Directed Suprema}
\begin{point}{20}[pub]{Theorem (Uniform Boundedness)}%
\index{Principle of Uniform Boundedness}%
\index{Uniform Boundedness Theorem}%
A set~$\scrF$ of bounded linear maps 
from a complete normed vector space~$\scrX$
to a normed vector space~$\scrY$
is bounded
in the sense that $\sup_{T\in \scrF} \|T\|<\infty$
provided that 
 $\sup_{T\in \scrF} \|Tx\|<\infty$
 for all~$x\in \scrX$.
\begin{point}{30}{Proof}%
Based on~\cite{sokal}.
\begin{point}{40}[sokal-lemma]%
Let $r>0$ and~$T\in\scrF$ be given.
Writing~$B_r(x)=\{\,y\in\scrX\colon \|x-y\|\leq r\,\}$
for the ball around~$x\in\scrX$ with radius~$r$,
note that $r\|T\|=\sup_{\xi\in B_r(0)} \|T \xi\|$
almost by definition of the operator norm.
We will need the less obvious fact
that $r\|T\|\leq \sup_{\xi \in B_r(x)}\|T \xi\|$
for every~$x\in \scrX$.

To see why this is true,
note that for~$\xi\in B_r(0)$
either $\|T\xi\|\leq \|T(x+\xi)\|$
or $\|T\xi\|\leq \|T(x-\xi)\|$,
because we would otherwise have
$2\|T\xi\| = \|T(x+\xi)-T(x-\xi)\|
\leq \|T(x+\xi)\|+\|T(x-\xi)\|<2\|T\xi\|$.
Hence
$r\|T\|=\sup_{\xi\in B_r(0)} \|T\xi\|\leq  
\sup_{\xi \in B_r(x)} \|T\xi \|$.
\end{point}
\begin{point}{50}%
Suppose towards a contradiction
that $\sup_{T\in\scrF}\|T\|=\infty$,
and pick~$T_1,T_2,\dotsc$ with $\|T_n\|\geq n3^{n}$.
Using~\sref{sokal-lemma},
choose $x_1,x_2,\dotsc$ in~$\scrX$
with $\|x_{n}-x_{n-1}\|\leq 3^{-n}$
and~$\|T_{n} x_{n}\|\geq \frac{2}{3}3^{-n}\|T_{n}\|$,
so that~$(x_n)_n$ is a Cauchy sequence, 
and therefore converges to some
$x\in\scrX$.
Note that~$\|x-x_n\|\leq \frac{1}{2}3^{-n}$
(because $\sum_{k=0}^\infty 3^{-k}=\frac{3}{2}$),
and so $\|T_n x\|\geq  \|T_nx_n\| - \|T_n(x_n-x)\|
\geq \frac{2}{3}3^{-n}\|T_n\|-\frac{1}{2}3^{-n}\|T_n\|
\geq \frac{1}{6}n$,
which contradicts
the assumption that $\sup_{T\in \scrF} \|Tx\| <\infty$.\qed
\end{point}
\end{point}
\end{point}
\begin{point}{60}[hellinger-toeplitz]{Theorem}%
Let~$T\colon X\to Y$ be an adjointable map
between pre-Hilbert $\scrA$-modules.
If either~$X$ or~$Y$ is complete,
then~$T$ and~$T^*$ are bounded.
\begin{point}{70}{Proof}%
We may assume without loss of generality
that~$X$ is complete (by swapping~$T$ for~$T^*$
and~$X$ with~$Y$ if necessary).

Note that for every~$y\in Y$,
the linear map $\left<y,T\,\cdot\,\right>\equiv
\left<T^*y,\,\cdot\,\right>\colon Y\to \scrA$
is bounded,
because $\left\|\left<T^*y,x\right>\right\| \leq \|T^*y\|\|x\|$
for all~$x\in X$ (see~\sref{chilb-cs}).

Since 
on the other hand,
$\left\|\left<y,Tx\right>\right\|
\leq \|y\|\,\|Tx\|\leq \|Tx\|$
for all~$x\in X$ and~$y\in Y$ with $\|y\|\leq 1$,
we have $\sup_{\|y\|\leq 1} \|\left<y,Tx\right>\| \leq \|Tx\|<\infty$
for all~$x\in X$,
and thus $B:=\sup_{\|y\|\leq 1} \|\left<y,T\,\cdot\,\right>\|<\infty$
by~\sref{pub}.

It follows that~$\|\left<y,Tx\right>\|\leq B\|y\|\|x\|$
for all~$y\in Y$ and~$x\in X$,
and thus~$T$ and~$T^*$ are bounded, by~\sref{chilb-form-bounded}.\qed
\end{point}
\begin{point}{80}{Remark}%
As a special case of the preceding theorem
we get the fact,
known as the \Define{Hellinger--Toeplitz theorem},%
\index{Hellinger--Toeplitz's Theorem}
that every symmetric
operator on a Hilbert space is bounded.
\end{point}
\begin{point}{90}[hellinger-toeplitz-needs-complete]{Example}%
The condition that either~$X$ or~$Y$ be complete may not be dropped:
the linear map $T\colon c_{00}\to c_{00}$
given by $T\alpha = (n\alpha_n)_n$ for $\alpha\in c_{00}$
is self-adjoint,
but not bounded,
because~$T$ maps $(1,\frac{1}{2},\dotsc,\frac{1}{n},0,0,\dotsc)$
having 2-norm below~$\frac{\pi}{\sqrt{6}}$
to $(1,1,\dotsc,1,0,0,\dotsc)$,
which has $2$-norm equal to~$\sqrt{n}$.
\end{point}
\end{point}
\end{parsec}
\begin{parsec}{360}%
\begin{point}{10}[self-dual]{Definition}%
A Hilbert $\scrA$-module~$X$ is \Define{self-dual}%
\index{Hilbert $\scrA$-module!self dual}
when every bounded module map $r\colon X\to \scrA$
is of the form $r\equiv \left<y,(\,\cdot\,)\right>$
for some~$y\in X$.
\end{point}
\begin{point}{20}{Example}%
By Riesz' representation theorem (\sref{riesz-representation-theorem})
every Hilbert space is self-dual.
\end{point}
\begin{point}{30}{Exercise}%
Show that given a $C^*$-algebra~$\scrA$
the Hilbert $\scrA$-module $\scrA^{N}$
of $N$-tuples is self dual.
\end{point}
\begin{point}{40}[chilb-form]{Definition}%
Let us say that a \Define{(bounded) form}%
\index{form, between Hilbert $\scrA$-modules}%
\index{form, between Hilbert $\scrA$-modules!bounded}
on Hilbert $\scrA$-modules
$X$ and~$Y$
is a map $[\,\cdot\,,\,\cdot\,]\colon X\times Y\to \scrA$
such that $[x,\,\cdot\,]\colon Y\to \scrA$
and $[\,\cdot\,,y]^*\colon X\to \scrA$
are (bounded) module maps for all~$x\in X$ and~$y\in Y$.
\end{point}
\begin{point}{50}[chilb-form-representation]{Proposition}%
For every bounded form  $[\,\cdot\,,\,\cdot\,]\colon X\times Y\rightarrow \scrA$
on self-dual Hilbert $\scrA$-modules
$X$ and~$Y$
there is a unique adjointable bounded module map
$T\colon X\to Y$.
with
$[x,y]\equiv \left<Tx,y\right>$
for all $x\in X$ and~$y\in Y$.
\begin{point}{60}{Proof}%
Let $x\in X$ be given.
Since~$[x,\,\cdot\,]\colon Y\to \scrA$ is a
a bounded module map,
and~$Y$ is self-dual,
there is a unique $Tx\in Y$ with
$[x,y]=\left<Tx,y\right>$
for all~$y\in Y$,
giving a map $T\colon X\to Y$.
For a similar reason
we get a map $S\colon Y\to X$
with $\left<Sy,x\right>=[x,y]^*$ 
for all~$x\in X$ and~$y\in Y$.
Since $S$ and~$T$ are clearly adjoint,
they are bounded module maps by~\sref{hellinger-toeplitz}.\qed
\end{point}
\end{point}
\end{parsec}
\begin{parsec}{370}%
\begin{point}{10}%
	Another consequence of~\sref{pub}
	is this:
\end{point}
\begin{point}{20}[hilb-weakly-bounded-complete]{Proposition}%
Given a net~$(y_\alpha)_\alpha$
in a Hilbert space~$\scrH$
for which $\left<y_\alpha,x\right>$
is Cauchy \emph{and bounded}\footnote{Recall that while every Cauchy
\emph{sequence} is bounded,
a Cauchy net need only be eventually bounded.}
for every~$x\in \scrH$,
there is a unique~$y\in\scrH$
with $\left<y,x\right>=\lim_\alpha \left<y_\alpha,x\right>$
for all~$x\in\scrH$.
\begin{point}{30}{Proof}%
To obtain~$y$,
we want to apply  Riesz' representation theorem
(\sref{riesz-representation-theorem})
to the linear map $f\colon \scrH\to\C$
defined by~$f(x)=\lim_\alpha\left<y_\alpha,x\right>$,
but must first show that~$f$ is bounded.
For this it suffices to show
that~$\sup_\alpha \left \|\left<y_\alpha,(\,\cdot\,)\right>\right\|<\infty$,
and this follows by~\sref{pub}
from the assumption 
that $\sup_{\alpha} \left|\left<y_\alpha,x\right>\right| <\infty$
for every~$x\in \scrH$.

By Riesz' representation theorem (\sref{riesz-representation-theorem}),
there is a unique~$y\in\scrH$ with 
$\left<y,x\right>=f(x)\equiv \lim_\alpha \left<y_\alpha,x\right>$
for all~$x\in \scrH$,
and so we're done.\qed
\end{point}
\begin{point}{40}{Remark}%
The condition in~\sref{hilb-weakly-bounded-complete} 
that the net~$(\,\left<y_\alpha,x\right>\,)_\alpha$
be bounded for every~$x$ may not be omitted
(even though $(\,\left<y_\alpha,x\right>\,)_\alpha$
being Cauchy is eventually bounded).

To see this,
consider a linear map $f\colon \scrH\to\C$ on a Hilbert space~$\scrH$
which is not bounded.
We claim that there is a net~$(y_\alpha)_\alpha$ in~$\scrH$
with $f(x)=\lim_\alpha \left<y_\alpha,x\right>$ for all~$x\in\scrH$,
and so there can be no~$y\in \scrH$ 
with $\left<y,x\right> = \lim_\alpha \left<y_\alpha,x\right>$
for all~$x\in \scrH$, because 
that would imply that~$f$ is bounded.

To create this net,
note that~$f$ is bounded
on the span $\left<F\right>$ of every 
finite subset $F\equiv \{x_1,\dotsc,x_n\}$
of vectors from~$\scrH$,
and so by Riesz' representation theorem~\sref{riesz-representation-theorem}
applied to~$f$ restricted to closed subspace~$\left<F\right>$
of~$\scrH$ there is a unique $y_F\in \left<F\right>$
such that~$f(x)=\left<y_F,x\right>$
for all~$x\in\left<F\right>$.

These $y_F$'s form a net in~$\scrH$
(when we order the finite subsets~$F$ of~$\scrH$ by inclusion),
which approximates~$f$ in the
sense that~$f(x)=\lim_F \left<y_F,x\right>$
for every~$x\in \scrH$,
(because 
$f(x)=\left<y_F,x\right>$
for every~$F$
with $\{x\}\subseteq F$).
\end{point}
\end{point}
\begin{point}{50}[swot]{Definition}%
Let~$\scrH$ be a Hilbert space.
\begin{enumerate}
\item
	The \Define{weak operator topology (WOT)}%
\index{WOT, weak operator topology}
on~$\scrB(\scrH)$ is the least topology
with respect to which $T\mapsto \left<x,Tx\right>,\,\scrB(\scrH)\to\C$
is continuous for every~$x\in\scrH$.

So a net $(T_\alpha)_\alpha$
		converges to~$T$ in $\scrB(\scrH)$
		with respect to the weak operator topology
		iff $\left<x,T_\alpha x\right>\to \left<x,Tx\right>$
		as~$\alpha\to\infty$ for all~$x\in \scrH$.
\item
	The \Define{strong operator topology (SOT)}%
\index{SOT, strong operator topology}
on~$\scrB(\scrH)$ is the least topology
with respect to which $T\mapsto \|Tx\|\equiv \smash{\left<x,T^*Tx\right>^{%
\nicefrac{1}{2}}}$
is continuous for every~$x\in\scrH$.

		So a net $(T_\alpha)_\alpha$
		converges to~$T$
		in $\scrB(\scrH)$
		with respect to the strong operator topology iff
		$\|T_\alpha x -Tx \| \to 0$ as $\alpha\to\infty$
		for all~$x\in\scrH$.
\end{enumerate}
\spacingfix%
\begin{point}{60}{Remark}%
Although we'll only make use of the weak operator
topology we have nonetheless included
the definition of the strong operator topology here
for comparison
with the \emph{ultrastrong
topology} that appears in the next chapter.
\end{point}
\end{point}
\begin{point}{70}[bh-wot-bounded-complete]{Lemma}%
Let~$(T_\alpha)_\alpha$ be a net of bounded operators
on a Hilbert space~$\scrH$
such that $(\,\left<x,T_\alpha x \right>\,)$ is
Cauchy and bounded for every~$x\in \scrH$.

Then~$(T_\alpha)_\alpha$
WOT-converges to some bounded operator~$T$ in $\scrB(\scrH)$.
\begin{point}{80}{Proof}%
Let~$x,y\in \scrH$ be given.
Since by a simple computation
(c.f.~\sref{inner-product-basic}\eqref{polarization-identity})
\begin{equation*}
	\textstyle
	\left<y,T_\alpha x\right>
	\ = \ \frac{1}{4}\sum_{k=0}^3
	i^k\left<\,i^ky+x,\,T_\alpha (i^ky+x)\,\right>,
\end{equation*}
 $(\,\left<y,T_\alpha x\right>\,)_\alpha$
is bounded for every~$y\in \scrH$,
and so by~\sref{hilb-weakly-bounded-complete} there is~$Tx\in \scrH$ 
with $\left<y,Tx\right>=\lim_\alpha \left<y,T_\alpha x\right>$
for all~$y\in\scrH$,
giving us a linear map $T\colon \scrH\to \scrH$.
It is clear that~$(T_\alpha)_\alpha$
WOT-converges to~$T$,
provided that~$T$ is bounded.

So to complete the proof,
we must show that~$T$ is bounded,
and we'll do this by showing that~$T$ has an adjoint
(see~\sref{hellinger-toeplitz}).
Note that $\left<x,T_\alpha^* x\right>=\overline{\left<x,T_\alpha x\right>}$
is Cauchy and bounded (with~$\alpha$ running),
so by a similar reasoning as before (but with~$T^*_\alpha$
instead of~$T_\alpha$)
we get a map
$S\colon \scrH\to\scrH$
with $\left<x,Sy\right>=\lim_\alpha \left<x,T^*_\alpha y\right>$
for all~$x,y\in\scrH$, which will be adjoint to~$T$,
which is therefore bounded.\qed
\end{point}
\end{point}
\begin{point}{90}[hilb-suprema]{Proposition}%
Let~$\scrH$ be a Hilbert space,
and~$\scrD$ an upwards directed subset of~$\Real{\scrB(\scrH)}$
with $\sup_{T\in \scrD} \left<x,Tx\right> <\infty$
for all~$x\in \scrH$. Then
\begin{enumerate}
\item
$(T)_{T\in\scrD}$
converges 
in the weak operator topology
to some~$T'$ in~$\Real{(\scrB(\scrH))}$,
\item
$T'$ is the supremum of~$\scrD$
in $\Real{(\scrB(\scrH))}$,
and 
\item
$\left<x,T'x\right> = 
\sup_{T\in\scrD}\left<x,Tx\right> $
for all~$x\in \scrH$.
\end{enumerate}
\spacingfix
\begin{point}{100}{Proof}%
Let~$x\in \scrH$.
Since $\left<x,(\,\cdot\,)x\right>\colon \scrB(\scrH)\to \C$
is positive
we see that
$(\left<x,Tx\right>)_{T\in\scrD}$
is an increasing net in~$\R$, 
bounded from above (by assumption),
and therefore converges to~$\sup_{T\in\scrD}\left<x,Tx\right>$.
In particular,~$(T)_{T\in\scrD}$
is WOT-Cauchy,
and ``WOT-bounded'',
and thus
(by~\sref{bh-wot-bounded-complete})
WOT-converges to some self-adjoint~$T'$ from~$\scrB(\scrH)$.

Since $(\,\left<x,Tx\right>\,)_{T\in\scrD}$
converges both to~$\left<x,T'x\right>$,
and to~$\sup_{T\in\scrD} \left<x,Tx\right>$,
we conclude that $\left<x,T'x\right>=\sup_{T\in\scrD}\left<x,Tx\right>$
for every~$x\in\scrH$.
In particular,  $\left<x,Tx\right>\leq \left<x,T'x\right>$
for all~$x\in\scrH$ and $T\in\scrD$, and thus $T\leq T'$
for all~$T\in\scrD$.

Let~$S$ be a self-adjoint bounded operator on~$\scrH$ with $T\leq S$
for all~$T\in\scrD$.
To prove that~$T'$ is the supremum of~$\scrD$,
we must show that~$T'\leq S$.
Let~$x\in \scrH$ be given.
Since $\left<x,Tx\right>\leq \left<x,Sx\right>$
for each~$T\in \scrD$ (because $T\leq S$),
we have $\left<x,T'x\right>\equiv \sup_{T\in\scrD} \left<x,Tx\right>
\leq \left<x,Sx\right>$,
and therefore $T'\leq S$ by~\sref{hilb-positive-operators}.\qed
\end{point}
\end{point}
\begin{point}{110}{Definition}%
Let~$\scrH$ be a Hilbert space.
The supremum of a (norm) bounded directed subset~$\scrD$ 
in~$\Real{(\scrB(\scrH))}$
(which exists by~\sref{hilb-suprema})
is denoted by~$\Define{\bigvee\scrD}$.%
\index{*infsup@$\bigvee D$, supremum of~$D$!in $\scrB(\scrH)$}
\end{point}
\end{parsec}
\subsection{Normal Functionals}
\begin{parsec}{380}%
\begin{point}{10}[bh-normal]{Definition}%
Given a Hilbert space~$\scrH$
a p-map $\omega\colon \scrB(\scrH)\to\C$
is called \Define{\textbf{n}ormal}%
\index{normal!positive functional!on $\scrB(\scrH)$}
when
$\omega(\bigvee \scrD)=\bigvee_{T\in\scrD} \omega(T)$
for every bounded directed subset $\scrD$ of~$\Real{\scrB(\scrH)}$.
\begin{point}{11}[bh-normal-abbreviation]{Notation}%
We use the letter ``n'' to abbreviate ``normal'' in line with~\sref{maps}.
So an npu-map $\omega\colon \scrB(\scrH)\to\C$
is a normal positive unital linear functional
on~$\scrB(\scrH)$.
\end{point}
\end{point}
\begin{point}{20}{Example}%
All vector functionals%
\index{vector functional!for a Hilbert space!is normal}
$\left<x,(\,\cdot\,)x\right>$ are normal by~\sref{hilb-suprema}.
\end{point}
\begin{point}{30}[bh-normal-effects]{Exercise}%
To show that a positive linear functional is normal, it suffices to show 
that it preserves directed suprema of \emph{effects}: 
show that given a Hilbert space~$\scrH$
a  positive map $\omega\colon \scrB(\scrH)\to\C$
is normal 
provided that $\omega(\bigvee \scrD) = \bigvee_{T\in\scrD} \omega(T)$
for every directed subset $\scrD$ of $[0,1]_{\scrB(\scrH)}$.
\end{point}
\begin{point}{40}[bh-functional-lemma]{Lemma}%
Every sequence $x_1,x_2,\dotsc $ in a Hilbert space~$\scrH$
with $\sum_n \|x_n\|^2 < \infty$
gives an np-map $\omega\colon\scrB(\scrH)\to\C$
defined by~$\omega(T)=\sum_n \left<x_n,Tx_n\right>$.
\begin{point}{50}{Proof}%
Given $T\in\scrB(\scrH)$ 
we have $\left|\left<x_n,Tx_n\right>\right|\leq \|x_n\|^2\|T\|$ 
by Cauchy--Schwarz (\sref{inner-product-basic}),
so $\sum_n \left|\left<x_n,Tx_n\right>\right|
\leq \|T\| \sum_n \|x_n\|^2$,
which means that~$\sum_n \left<x_n,Tx_n\right>$
converges, 
and so we may define~$\omega$ as above.

It is easy to see that~$\omega$ is linear and positive,
so we'll only show that~$\omega$ is normal.
We must prove that $\omega(\bigvee \scrD)=\bigvee_{T\in\scrD} \omega(T)$
for every bounded directed subset of~$\Real{(\scrB(\scrH))}$.
By~\sref{bh-normal-effects}
we may assume without loss of generality that 
$\scrD\subseteq [0,1]_{\scrB(\scrH)}$.
This has the benefit that $\left<x_n,T x_n\right>$
is positive for all~$n$ and~$T\in\scrD$,
so that their sum (over~$n$) is given by
a supremum over partial sums, viz.~$\sum_n\left<x_n,Tx_n\right>
=\bigvee_N\sum_{n=1}^N\left<x_n,Tx_n\right>$.
Completing the proof is now simply a matter of
interchanging suprema,
\begin{alignat*}{3}
	\textstyle \bigvee_{T\in \scrD} \omega(T)
	\ &=\ 
	\textstyle\bigvee_{T\in \scrD} \bigvee_N \sum_{n=1}^N 
	\left<x_n,Tx_n\right>\\
	\ &=\ 
	\textstyle\bigvee_N \bigvee_{T\in \scrD}\sum_{n=1}^N 
	\left<x_n,Tx_n\right>\\
	\ &=\ 
	\textstyle\bigvee_N \sum_{n=1}^N \left<x_n,(\bigvee \scrD )\,x_n\right>
	\ =\ \textstyle\omega(\bigvee\scrD),
\end{alignat*}
where we used that~$\sum_{n=1}^N \left<x_n,(\,\cdot\,)x_n\right>$
is normal.\qed
\end{point}
\end{point}
\begin{point}{60}[vector-functional-convergence]{Exercise}%
The following observations
regarding
a net~$(x_\alpha)_\alpha$ in a Hilbert space~$\scrH$
will be useful later on.
\begin{enumerate}
\item
Show that~$\sum_\alpha \|x_\alpha\|^2<\infty$
if and only if~$\sum_\alpha \left<x_\alpha,(\,\cdot\,)x_\alpha\right>$
converges with respect to the operator norm
to some bounded functional on~$\scrB(\scrH)$.
\item
Given some~$x\in \scrH$,
show that~$x_\alpha$ converges to~$x$
if and only if  $\left<x_\alpha,(\,\cdot\,)x_\alpha\right>$
operator-norm converges to~$\left<x,(\,\cdot\,)x\right>$.

(For the ``if'' part it may be convenient
to first prove that $\left<x_\alpha,x\right>\to \left<x,x\right>$
by considering the bounded operator
$\ketbra{x}{x}$ on~$\scrB(\scrH)$.)
\end{enumerate}
\spacingfix
\end{point}%
\end{parsec}%
\begin{parsec}{390}%
\begin{point}{10}%
The final project of this chapter
is to show that each normal positive functional~$\omega$ 
on a~$\scrB(\scrH)$
	is of the form 
	$\omega\equiv \sum_{n=0}^\infty\left<x_n,(\,\cdot\,)x_n\right>$
	for some~$x_1,x_2,\dotsc$
	with~$\sum_n\|x_n\|^2<\infty$.
For this we'll need some more nuggets
from the theory of Hilbert spaces.
\end{point}
\begin{point}{20}{Definition}%
A subset~$\scrE$ of a Hilbert space
is called \Define{orthonormal}%
\index{orthonormal, subset of a Hilbert space}
if $\left<e,e'\right>=0$
for all~$e,e'\in\scrE$ with~$e\neq e'$,
and~$\left<e,e\right>=1$
for all~$e\in\scrE$.
We say that~$\scrE$ is \Define{maximal}
\index{orthonormal, subset of a Hilbert space!maximal}
when~$\scrE$ is maximal
among all orthonormal  subsets of~$\scrH$
ordered by inclusion,
and in that case
we call~$\scrE$ an \Define{orthonormal basis}%
\index{orthonormal basis, for a Hilbert space}
for~$\scrH$
for reasons that will be become clear in~\sref{orthonormal}
below.
\begin{point}{30}{Remark}%
Clearly, by Zorn's lemma,
each Hilbert space has an orthonormal basis.
\end{point}
\end{point}
\begin{point}{40}[orthonormal]{Proposition}%
Given an orthonormal subset~$\scrE$
of a Hilbert space~$\scrH$,
and~$x\in \scrH$,
\begin{enumerate}
\item
\label{orthonormal-1}
\Define{(Bessel's inequality)}%
\index{Bessel's inequality}
\ 
 $\sum_{e\in \scrE}\left|\left<e,x\right>\right|^2
\leq \|x\|^2$;
\item
\label{orthonormal-2}
$\sum_{e\in \scrE} \left<e,x\right>e$
converges in~$\scrH$,
\item
\label{orthonormal-3}
$\sum_{e\in \scrE} \left<e,x\right>e=x$
if~$\scrE$ is maximal, and
\item
\label{orthonormal-4}
\Define{(Parseval's identity)}%
\index{Parseval's identity}
$\sum_{e\in\scrE}\left|\left<e,x\right>\right|^2 = \|x\|^2$
if~$\scrE$ is maximal.
\end{enumerate}
\spacingfix%
\begin{point}{50}{Proof}%
\ref{orthonormal-1}\ 
Since for any finite subset $\scrF$ of $\scrE$
we have $0\leq \|x-\sum_{e\in \scrF} \left<e,x\right>e\|^2
= \|x\|^2-2\sum_{e\in \scrF} \left<e,x\right>\left<x,e\right>
+ \sum_{e,e'\in\scrF} \left<x,e'\right>\left<e',e\right>\left<e,x\right>
= \|x\|^2-\sum_{e\in\scrF}\left|\left<e,x\right>\right|^2$,
and so~$\sum_{e\in\scrF} \left|\left<e,x\right>\right|^2\leq \|x\|^2$,
we get~$\sum_{e\in\scrE}\left|\left<e,x\right>\right|^2\leq \|x\|^2$.

\ref{orthonormal-2}\ 
From the observation
that~$\|\sum_{e\in\scrF} \left<e,x\right>e\|^2
= \sum_{e\in \scrF} \left|\left<e,x\right>\right|^2$
for any finite~$\scrF\subseteq \scrE$,
and the fact that~$\sum_{e\in \scrE} \left|\left<e,x\right>\right|^2$
converges (by the previous point),
one deduces that~$(\sum_{e\in\scrF} \left<e,x\right>e)_\scrF$
is Cauchy,
and so~$\sum_{e\in\scrE} \left<e,x\right>e$
converges.

\ref{orthonormal-3}\
Writing~$y:=\sum_{e\in\scrE} \left<e,x\right>e$ we must show that~$x=y$.
If it were not so,
that is,~$x\neq y$,
then~$e':=\|x-y\|^{-1}(x-y)$
satisfies
$\left<e',e'\right>=1$
and
$\left<e',e\right>=0$
for all~$e\in\scrE$,
and so may be added to~$\scrE$
to yield an orthonormal basis~$\scrE\cup\{e'\}$
extending~$\scrE$
contradicting~$\scrE$s  maximality.

\ref{orthonormal-4}\ 
Finally,
$\|x\|^2=\left<x,x\right>
= \sum_{e,e'\in\scrE}
\left<x,e'\right>\left<e',e\right>\left<e,x\right>
= \sum_{e\in \scrE} \left|\left<e,x\right>\right|^2$.\qed
\end{point}
\end{point}
\begin{point}{60}[sum-ketbras]{Exercise}%
Let~$\scrE$ be an orthonormal basis of a Hilbert space~$\scrH$.
\begin{enumerate}
\item
Show that~$\sum_{e\in\scrE} \ketbra{e}{e}$
converges to~$1$
in the weak operator topology.
\item
Show that $\sum_{e\in\scrE}\ketbra{e}{e}=1$
also in the sense
that the directed
set of partial sums
$\sum_{e\in \scrF} \ketbra{e}{e}$,
where~$\scrF$ is a finite subset of~$\scrE$,
has~$1$ as its supremum.
\item
Conclude that~$\omega(1)=\sum_{e\in\scrE} \omega(\ketbra{e}{e})$
for every np-map $\omega\colon \scrB(\scrH)\to\C$.
\end{enumerate}
\spacingfix%
\end{point}%
\begin{point}{70}[bh-np-lemma]{Lemma}%
Given a Hilbert space~$\scrH$
with orthonormal basis~$\scrE$,
we have
\begin{equation*}
	\omega(A)\ = \ 
	\sum_{e,e'\in\scrE} \left<e,Ae'\right>\ \omega(\,\ketbra{e}{e'}\,).
\end{equation*}
for every normal p-map $\omega\colon \scrB(\scrH)\to\C$
and~$A\in\scrB(\scrH)$.
\begin{point}{80}{Proof}%
Let~$\scrF$ be a finite subset of~$\scrE$,
and write $P=\sum_{e\in \scrF} \ketbra{e}{e}$.
Since $PAP
= \sum_{e,e'\in\scrF}\left<e,Ae'\right>\,\ketbra{e}{e'}$
it suffices
to show that~$\omega(A-PAP)$
vanishes as~$\scrF$ increases.
Note that~$P^*P=P$ and
$(P^\perp)^*P^\perp=P^\perp$.
Further,
since
$\|P\|\leq 1$,
and~$A-PAP=P^\perp A + PAP^\perp$,
we have,
by Kadison's inequality,
\begin{alignat*}{3}
\left|\omega(A-PAP)\right|
\ &\leq\  \left|\omega(P^\perp A)\right| \,+\,\left|\omega(PAP^\perp)\right| \\
\ &\leq\  
\omega(P^\perp)^{\nicefrac{1}{2}}\,
\omega(A^*A)^{\nicefrac{1}{2}}
\ +\ 
\omega(PAA^*P)^{\nicefrac{1}{2}}\,
\omega(P^\perp)^{\nicefrac{1}{2}}\\
\ &\leq\  
2\|A\| 
\omega(1)^{\nicefrac{1}{2}}\ 
\omega(P^\perp)^{\nicefrac{1}{2}}.
\end{alignat*}
But since~$\sum_{e\in\scrE} \omega(\ketbra{e}{e})=\omega(1)$
by~\sref{sum-ketbras}
we see that~$\omega(P^\perp)\to0$
as~$\scrF\to\infty$.\qed
\end{point}
\end{point}
\begin{point}{90}[bh-np]{Theorem}%
\index{normal!positive functional!on~$\scrB(\scrH)$}%
Let~$\scrH$ be a Hilbert space.
Every normal p-map $\omega\colon \scrB(\scrH)\to \C$
is of the form $\omega = \sum_n\left<x_n,(\,\cdot\,)x_n\right>$
where $x_1,x_2,\dotsc\in \scrH$
with~$\sum_n \|x_n\|^2=\|\omega\|$.
\begin{point}{100}{Proof}%
By~\sref{chilb-form-representation}
there is a
unique
$\varrho\in\scrB(\scrH)$
with $\omega(\ketbra{y}{x})=\left<x,\varrho y\right>$
for all~$x,y\in\scrH$,
because
$(x,y)\mapsto \omega(\ketbra{y}{x}),\,
\scrH\times\scrH\to\C$
is a bounded form in the sense of~\sref{chilb-form}.
Note that~$\varrho$ 
is positive by~\sref{hilb-positive-operators}
because $\left<x,\varrho x\right>
=\omega(\ketbra{x}{x})\geq 0$
for all~$x\in\scrH$.
Now, let~$\scrE$ be an orthonormal basis for~$\scrH$.
Since~$\omega$ is normal,
\sref{sum-ketbras} gives us
$\omega(1)=\sum_{e\in\scrE} \omega(\ketbra{e}{e})
= \sum_{e\in\scrE} \left<e,\varrho e\right>
= \sum_{e\in \scrE} \|\sqrt{\varrho} e\|^2$,
so that $\omega':=\sum_{e\in\scrE} \left<\sqrt{\varrho}e,(\,\cdot\,)
\sqrt{\varrho}e\right>$
defines a normal positive functional on~$\scrB(\scrH)$
by~\sref{vector-functional-convergence}.
Thus,
we are done if
can show that~$\omega'=\omega$,
(because $\sqrt{\varrho}e$
is non-zero for at most countably many~$e\in\scrE$).
To this end,
note that
$\omega(\ketbra{x}{x})
= \left<\sqrt{\varrho}x,\sqrt{\varrho}x\right>
= \sum_{e\in \scrE}
\left<\sqrt{\varrho}x,e\right>
\left<e,\sqrt{\varrho}x\right>
= \sum_{e\in \scrE} \left< \sqrt{\varrho}e,
\ketbra{x}{x} \sqrt{\varrho}e \right>
=\omega'(\ketbra{x}{x})$
for each~$x\in\scrH$,
and so $\omega(\ketbra{x}{y})=\omega'(\ketbra{x}{y})$
for all~$x,y\in\scrH$
by polarisation,
and thus~$\omega=\omega'$
by~\sref{bh-np-lemma}.\qed
\end{point}
\end{point}
\end{parsec}
\begin{parsec}{400}
\begin{point}{10}
In this chapter
we've studied
the algebraic  structure of 
the space~$\scrB(\scrH)$
of bounded operators on a Hilbert space~$\scrH$
abstractly 
via the notion
of
a $C^*$-algebra.
We've seen not only that every $C^*$-algebra
is miu-isomorphic to a $C^*$-subalgebra
	of such a $\scrB(\scrH)$ (in~\sref{gelfand-naimark}),
but also that any commutative $C^*$-algebra
is miu-isomorphic
to the space~$C(X)$ of continuous functions
	on some compact Hausdorff space (in~\sref{gelfand}).
But there's more to~$\scrB(\scrH)$
than just being a $C^*$-algebra:
it has the two additional properties
of having suprema of bounded directed subsets (see~\sref{hilb-suprema}),
and having a faithful collection 
of normal functionals (viz.~the vector functionals,
\sref{hilb-vector-states-order-separating}).
This leads us to the study of von Neumann algebras---the topic
of the next chapter.
\end{point}
\end{parsec}

\chapter{Von Neumann Algebras}
\begin{parsec}{410}%
\begin{point}{10}%
We have arrived at the main subject of this thesis,
the special class of $C^*$-algebras
called von Neumann algebras (see definition~\sref{vna} below)
that are characterised by the existence
of certain directed suprema 
and an abundance  of functionals that preserve
these suprema.
While all $C^*$-algebras
and the cpsu-maps
between them
may perhaps serve as models for
quantum data types and processes, respectively,
we focus
for the purposes of this thesis
 our attention on
the subcategory~$\W{cpsu}$ of von Neumann algebras
and the cpsu-maps between them that preserve these suprema
(called \emph{normal} maps, see~\sref{p-uwcont}),
because 
\begin{enumerate}
\item
$\W{cpsu}$
is a model of the quantum lambda calculus
(in a way that~$\Cstar{cpsu}$ is not,
see~\sref{cstar-no-model}), and
\item
we were able to axiomatise
the sequential product ($b\mapsto \sqrt{a}b\sqrt{a}$)
in~$\W{cpsu}$ 
(but not in~$\Cstar{cpsu}$)
see~\sref{uniqueness-sequential-product}.
\end{enumerate}
Both these are reserved for the next chapter;
in this chapter we'll (re)develop the theory
we needed to prove them.

The archetypal von Neumann algebra
is the $C^*$-algebra~$\scrB(\scrH)$
of bounded operators on a Hilbert space~$\scrH$.
In fact,
the
original~\cite{vn1930,mvn1936}
and common~\cite{kr,conway2000} 
definition of a von Neumann algebra
is a $C^*$-subalgebra~$\scrA$
of a~$\scrB(\scrH)$
that is closed in a ``suitable topology''
such as the strong or weak operator topology
(see~\sref{swot}).
Most authors make the distinction
between such rings of operators
(called von Neumann algebras)
and the $C^*$-algebras
miu-isomorphic to them
(called \emph{$W^*$-algebras}),
but we won't bother and call them all von Neumann algebras.
Partly because it seems difficult
to explain 
to someone
picturing a quantum data type
the meaning of the weak operator topology
and the Hilbert space~$\scrH$,
we'll use Kadison's characterisation~\cite{kadison1956}
of von Neumann algebras
as $C^*$-algebras
with a certain dcpo-structure (c.f.~\sref{hilb-suprema})
and sufficiently many Scott-continuous functionals (c.f.~\sref{bh-normal})
as our definition instead, see~\sref{vna}.

But we also use Kadison's definition
just to see
to what extent the representation
of von Neumann algebras
as rings of operators (see~\sref{ngns}) can be avoided
when erecting the basic theory.
Instead we'll put the directed suprema and normal 
positive functionals
on centre stage.
All the while
our treatment doesn't stray too far
from the beaten path,
and borrows
many arguments
from
the standard texts~\cite{sakai,kr};
but  most of them had to be tweaked in places, and
some demanded a complete overhaul.

The material on von Neumann algebras
is less tightly knit as the theory of $C^*$-algebras,
and so after the basics
we deal with four topics
more or less in linear order
(instead of intertwined.)

The great abundance of projections
(elements~$p$ with $p^*p=p$)
in von Neumann algebras---a definite advantage
over $C^*$-algebras---is
the first topic.
We'll see for example that
the existence of norm bounded directed suprema
in a von Neumann algebra~$\scrA$
allows us to show
that
there is a least projection~$\ceil{a}$
above any effect~$a$ from~$\scrA$
given by~$\ceil{a}=\bigvee_n a^{\nicefrac{1}{2}^n}$
(see~\sref{vna-ceil});
and also that any element of a von Neumann algebra
can be written as a norm limit
of linear combinations of projections (in~\sref{projections-norm-dense}).
Many a result about von Neumann algebras
can be proven by an appeal to projections.

The second topic concerns
two topologies that are instrumental
for the more delicate results and constructions:
the \emph{ultraweak topology}
induced by the normal positive functionals~$\omega\colon\scrA\to\C$,
and the \emph{ultrastrong topology}
induced by the associated seminorms~$\|\,\cdot\,\|_\omega$
(see~\sref{vna}).
We'll show among other things that a von Neumann algebra
is complete with respect to the ultrastrong topology
and \emph{bounded} complete with respect to the ultraweak topology
(see~\sref{vn-complete}).

This completeness allows us to
define,
for example,
for any pair~$a$, $b$ of elements
from a von Neumann algebra~$\scrA$
with~$a^*a\leq b^*b$
an element $a/b$
with~$a=(a/b) \, b $
(see~\sref{division})---this is the third topic.
Taking~$b=\sqrt{a^*a}$ we
obtain
the famous
polar decomposition~$a = (a/\sqrt{a^*a}) \, \sqrt{a^*a}$
(see~\sref{polar-decomposition},
which is usually proven 
for a bounded operator on a Hilbert space first).

The fourth, and final topic,
is ultraweakly continuous functionals
on a von Neumann algebra:
we'll show
in~\sref{vn-center-separating-fundamental}
that any centre separating collection (\sref{separating})
of normal positive functionals~$\Omega$
on a von Neumann algebra completely
determines the normal positive functionals,
which will be important for the definition of the tensor
product of von Neumann algebras in the next chapter,
see~\sref{tensor}.
\end{point}
\end{parsec}
\section{The Basics}
\subsection{Definition and Counterexamples}
\begin{parsec}{420}[vna]%
\begin{point}{10}{Definition (Kadison~\cite{kadison1956})}%
A $C^*$-algebra~$\scrA$
is a \Define{von Neumann algebra}%
\index{von Neumann algebra}
when
\begin{enumerate}
\item
every bounded directed subset~$D$
of self-adjoint elements of~$\scrA$ (so $D\subseteq \sa{\scrA}$) 
    has a supremum \Define{$\bigvee D$}%
\index{*infsup@$\bigvee D$, supremum of~$D$!in a von Neumann algebra}
    in $\sa{\scrA}$, and
\item
if $a$ is a positive element of~$\scrA$
with $\omega(a)=0$ for every \emph{normal} (see below) positive 
linear map $\omega\colon \scrA\to \C$,
then~$a=0$.\footnote{In other words,
the collection of normal positive functionals should be faithful
(see~\sref{separating}).
Interestingly,
it's already enough for the normal positive 
functionals to be centre separating,
but since we have encountered no example
of a von Neumann algebra
where it wasn't already clear that the normal positive 
functionals are faithful
instead of just centre separating
we did not use this weaker albeit more complex condition.}
\end{enumerate}
\spacingfix%
\begin{point}{20}[def-np-functional]%
A positive linear map $\omega\colon \scrA\to \C$
is called \Define{\textbf{n}ormal}%
\index{normal!functional}
if $\omega(\bigvee D) = \bigvee_{d\in D} \omega(d)$
for every bounded directed subset of self-adjoint elements of~$D$
which has a supremum $\bigvee D$ in $\sa{\scrA}$.
\begin{point}{21}
Recall that we use the letter ``n'' as abbreviation
for ``normal'', see~\sref{bh-normal-abbreviation}.
\end{point}
\end{point}%
\begin{point}{30}%
The \Define{ultraweak topology}
on $\scrA$
is the least topology
that makes all normal positive linear maps $\omega\colon \scrA\to \C$
continuous; the ultraweakly open subsets of~$\scrA$
are exactly the unions
of finite intersections of
sets of the form $\omega^{-1}(U)$, where
$\omega\colon \scrA\to\C$ is an np-map, 
and $U$ is an open subset of~$\C$.
One can verify that a net $(b_\alpha)_\alpha$
in~$\scrA$ converges ultraweakly to some~$b$ in~$\scrA$
    iff $\omega(b_\alpha)\to b$ for all
    np-maps $\omega\colon \scrA\to\C$.
The \Define{ultrastrong topology}%
\index{ultraweak and ultrastrong}
	on~$\scrA$
is the topology
induced by the seminorms
$\|\,\cdot\,\|_\omega$
associated to the np-maps $\omega\colon \scrA\to\C$
    (given by
$\|a\|_\omega \equiv \omega(a^*a)^{\nicefrac{1}{2}}$,
    see~\sref{omega-norm-basic});
a subset of~$\scrA$ is ultrastrongly open
iff it is the union of a finite intersections
of sets of the form
$\{\,a\in\scrA\colon\, \|a-b\|_\omega \leq \varepsilon\,\}$,
where~$b\in \scrA$, $\omega\colon \scrA\to\C$ is an np-map,
and~$\varepsilon>0$.
One can prove that a net~$(b_\alpha)_\alpha$
in~$\scrA$
converges ultrastrongly to an element~$b$ of~$\scrA$
    iff~$\|b_\alpha-b\|_\omega\to 0$ for all np-maps 
    $\omega\colon \scrA\to\C$.
\end{point}
\end{point}
\begin{point}{40}{Remark}%
We work with the ultraweak and ultrastrong topology in tandem,
because neither is ideal, and they tend to be complementary:
for example, $a\mapsto a^*$ is ultraweakly continuous
but not ultrastrongly (see~\sref{vn-counterexamples},
point~\ref{vn-counterexamples-4}), 
while $a\mapsto \left|a\right|$
is ultrastrongly continuous (see~\sref{abs-us-cont})
but not ultraweakly (\sref{vn-counterexamples}, 
point~\ref{vn-counterexamples-6}).
This doesn't prevent 
the ultraweak topology 
from being weaker than the ultrastrong topology:
a net that converges ultrastrongly converges ultraweakly as well,
see~\sref{uwweaker}.
\end{point}
\begin{point}{50}[von-neumann-examples]{Examples}%
\begin{enumerate}
\item
\index{C@$\C$, the complex numbers!as a von Neumann algebra}%
$\C$ and~$\{0\}$ are clearly von Neumann algebras.
\item
\index{BH@$\scrB(\scrH)$!as a von Neumann algebra}
The $C^*$-algebra $\scrB(\scrH)$
of bounded operators on a Hilbert space~$\scrH$
is a von Neumann algebra:
$\scrB(\scrH)$ has bounded directed suprema
of self-adjoint elements
by~\sref{hilb-suprema},
and the vector states
(and thus all normal functionals)
are order separating
(and thus faithful)
by~\sref{hilb-vector-states-order-separating}.
\item%
\index{direct sum!of von Neumann algebras}%
\index{$\bigoplus$, direct sum!$\bigoplus_i \scrA_i$, of von Neumann algebras}
The direct sum $\bigoplus_i \scrA_i$
(see~\sref{cstar-product})
of a family $(\scrA_i)_i$
of von Neumann algebras
is itself a von Neumann algebra.

(While we're not quite ready to define morphisms
between von Neumann algebras,
we can already spoil that the direct sum
gives the categorical product of von Neumann algebras
once we do,
see~\sref{vn-products}.)
\item
A $C^*$-subalgebra~$\scrB$
of a von Neumann algebra~$\scrA$
is called a \Define{von Neumann subalgebra}%
\index{von Neumann subalgebra}
(and is itself a von Neumann algebra)
if for every bounded directed subset~$D$
of self-adjoint elements from~$\scrB$
we have $\bigvee D\in\scrB$
(where the supremum is taken in~$\sa{\scrA}$).

\item[4a.]
Let~$S$ be a subset of a von Neumann algebra~$\scrA$.
Since the intersection of an arbitrary collection of von Neumann subalgebras
of~$\scrA$ is a von Neumann subalgebra of~$\scrA$ as well,
there is
        a least von Neumann subalgebra, $\Define{W^*(S)}$,%
\index{$W^*(S)$, von Neumann subalgebra generated by~$S$}
that contains~$S$.

\item
We'll see in~\sref{commutant-basic}
that given a subset~$S$ of a von Neumann algebra~$\scrA$
the set~$S^\square = \{\,a\in\scrA\colon\, \forall s\in S\,[\ as=sa\ ]\,\}$
called the \emph{commutant} of~$S$
is a von Neumann subalgebra of~$\scrA$
when~$S$ is closed under involution.
\item
We'll see in~\sref{mn-vna}
that the $N\times N$-matrices over a von Neumann algebra~$\scrA$
form a von Neumann algebra.
\item
We'll see in~\sref{Linfty-vn}
that the bounded measurable functions
on a finite complete measure space~$X$ 
(modulo the negligible ones)
form
a commutative von Neumann algebra~$L^\infty(X)$.

(Recall that a measure space~$X$ is called finite
when $\mu(X)<\infty$.)
\end{enumerate}
\spacingfix%
\end{point}%
\end{parsec}%
\begin{parsec}{430}%
\begin{point}{10}[uwweaker]{Exercise}%
Let~$\scrA$ be a von Neumann algebra.
\begin{enumerate}
\item
Show that 
$\left|\omega(a)\right|\leq \|a\|_\omega \|\omega\|^{\nicefrac{1}{2}}$
for every np-map $\omega\colon \scrA\to\C$
and~$a\in\scrA$.
\item
Show that when a net $(a_\alpha)_\alpha$
in~$\scrA$ converges ultrastrongly to~$a\in \scrA$
it does so ultraweakly, too.
\item
Show that an ultraweakly closed subset~$C$ of~$\scrA$
is also ultrastrongly closed.
\end{enumerate}%
\spacingfix%
\end{point}%
\begin{point}{11}[infima-in-vna]{Exercise}%
Note that given a von Neumann algebra~$\scrA$
the map $a\mapsto -a\colon \scrA\to\scrA$ 
is an order reversing isomorphism.

Deduce from this that any bounded filtered\footnote{`Filtered' is the
order dual of `directed':
    $F$ is
filtered when
for all~$a,b\in F$ there is~$c\in F$ with $c\leq a$ and~$c\leq b$.}
subset~$F$
of self-adjoint elements of~$\scrA$
has as infimum $\Define{\bigwedge F} := -\bigvee\{\,-d\colon \,d\in F\,\}$.%
\index{*inf@$\bigwedge F$, infimum of~$F$!in a von Neumann algebra}
\end{point}
\begin{point}{20}[vn-counterexamples]{Exercise}%
We give some counterexamples in $\scrB(\ell^2)$
to plausible propositions
to sharpen your understanding of the ultrastrong and ultraweak topologies,
and so that you may better appreciate
the strange manoeuvres we'll need to pull off later on.
\begin{enumerate}
\item
First some notation: given~$n,m\in \N$,
we denote by $\Define{\ketbra{n}{m}}$%
\index{*ketbranm@$\ketbra{n}{m}$, with $n,m\in\N$}
the bounded operator on~$\ell^2$
given by $(\ketbra{n}{m})(f)(n)=f(m)$
and~$(\ketbra{n}{m})(f)(k)=0$ for $k\neq n$
and $f\in \ell^2$.

Verify the following computation rules,
where $k,\ell,m,n\in \N$.
\begin{equation*}
(\ketbra{n}{m})^*\ =\ \ketbra{m}{n},
\qquad
\ketbra{n}{m}\ketbra{\ell}{k}\ =\ 
\begin{cases}
\ \ \ketbra{n}{k} & \text{if $m=\ell$} \\
\ \ 0 & \text{otherwise}
\end{cases}
\end{equation*}
\item
Show that $\bigvee_N \sum_{n=0}^N \ketbra{n}{n}=1$.

Conclude that~$(\,\ketbra{n}{n}\,)_n$
converges ultrastrongly (and ultraweakly) to~$0$.

Thus ultrastrong (and ultraweak) convergence does not imply norm convergence,
which isn't unexpected.
But we also see that if a sequence~$(b_n)_n$ converges ultrastrongly
(or ultraweakly) to some~$b$,
then $(\|b_n\|)_n$ doesn't even have to converge to~$\|b\|$.

(Note that~$(\ketbra{n}{n})_n$ resembles a `moving bump'.)
\item
Note that when a net $(a_\alpha)_\alpha$
converges ultrastrongly to~$a$,
then $(\,a_\alpha^*a_\alpha\,)_\alpha$
is norm-bounded and
converges ultraweakly to~$a^*a$.

The converse does not hold:
show that (already in~$\C$)
$e^{in}$
does not converge ultraweakly 
(nor ultrastrongly) as $n\to \infty$,
while $1\equiv e^{-in} e^{in}$
is norm-bounded and
converges ultraweakly to~$1$ as~$n\to\infty$.
\item
\label{vn-counterexamples-4}
Show that~$(\,\ketbra{0}{n}\,)_n$ converges ultrastrongly 
(and ultraweakly) to~$0$.

Deduce that $(\,\ketbra{n}{0}\,)_n$ converges ultraweakly to~$0$,
but doesn't converge ultrastrongly at all.

Conclude that~$a\mapsto a^*$ is not ultrastrongly continuous 
on~$\scrB(\ell^2)$.

(This has the annoying side-effect
that it is not immediately clear that the ultrastrong
closure of a $C^*$-subalgebra of a von Neumann algebra 
is a von Neumann subalgebra; we'll deal with this
by showing that the ultrastrong closure coincides
with the ultraweak closure in~\sref{ultraclosed}.)
\item
Show that the unit ball~$(\,\scrB(\ell^2)\,)_1$
of~$\scrB(\ell^2)$ is not ultrastrongly compact
by proving that $(\,\ketbra{0}{n}\,)_n$
has no ultrastrongly convergent subnet.

(But we'll see in~\sref{vn-ball-compact} that
the unit ball of a von Neumann algebra
is ultraweakly compact.)

\item
\label{vn-counterexamples-6}
Show that $\ketbra{n}{0}+\ketbra{0}{n}$
converges ultraweakly to~$0$ as $n\to \infty$,
while $(\ketbra{n}{0}+\ketbra{0}{n})^2\equiv \ketbra{0}{0}+\ketbra{n}{n}$
converges ultraweakly to~$\ketbra{0}{0}$.

Conclude that~$a\mapsto a^2$ is not ultraweakly continuous on~$\scrB(\ell^2)$.

Conclude that $a,b\mapsto ab$ is not jointly ultraweakly continuous
on~$\scrB(\ell^2)$.

Prove that~$\left|\, \ketbra{n}{0}+\ketbra{0}{n}\,\right|
= \ketbra{0}{0}+\ketbra{n}{n}$.

Conclude that~$a\mapsto \left|a\right|$
is not ultraweakly continuous on~$\sa{(\scrB(\ell^2))}$.

(We'll see in~\sref{proto-kaplansky} that $a\mapsto \left|a\right|$
is ultrastrongly continuous on self-adjoint elements.)

\item
Let us consider the two extensions of~$\left|\,\cdot\,\right|$
to arbitrary elements, namely
$a\mapsto \sqrt{a^*a}=:\Define{\left|a\right|_s}$ and 
$a\mapsto \sqrt{aa^*}=:\Define{\left|a\right|_r}$
(for \textbf{s}upport and \textbf{r}ange,
c.f.~\sref{hilb-ceil}).

Prove that $\ketbra{0}{0} + \ketbra{0}{n}$
converges ultrastrongly to~$\ketbra{0}{0}$ as $n\to \infty$.

Show that $\left|\,\ketbra{0}{0}+\ketbra{0}{n}\,\right|_s
= \ketbra{0}{0}
+\ketbra{0}{n}
+\ketbra{n}{0}
+\ketbra{n}{n}$
converges ultraweakly to~$\left|\,\ketbra{0}{0}\,\right|_s
\equiv \ketbra{0}{0}$
as $n\to \infty$,
but not ultrastrongly.

Show that $\left|\,\ketbra{0}{0}+\ketbra{0}{n}\,\right|_r
= \sqrt{2} \ketbra{0}{0}$.

Conclude that $\left|\,\cdot\,\right|_s$
and $\left|\,\cdot\,\right|_r$
are not ultrastrongly continuous on~$\scrB(\ell^2)$.

\item
Show that $1+\ketbra{n}{0}+\ketbra{0}{n}$
is positive,
and 
converges ultraweakly to~$1$ as~$n\to\infty$,
while the squares
$1+\ketbra{n}{n}+\ketbra{0}{0}+2\ketbra{n}{0}+2\ketbra{0}{n}$
converge ultraweakly to $1+\ketbra{0}{0}$
(as $n\to\infty$).

Hence~$a\mapsto a^2$
and $a\mapsto \sqrt{a}$
are not ultraweakly continuous on $\pos{\scrB(\ell^2)}$.

\item
\label{vn-counterexamples-9}
For the next counterexample,
we need a growing moving bump,
which still converges ultraweakly.
Sequences won't work here:

Show that $n\ketbra{n}{n}$ does not converge ultraweakly as~$n\to\infty$.

Show that $n \ketbra{f(n)}{f(n)}$ does not converge ultraweakly
as $n\to\infty$
for every strictly monotone (increasing) map~$f\colon \N\to\N$.

So we'll resort to a net.
Let~$D$ be the directed set which consists of pairs $(n,f)$,
where $n\in \N\backslash\{0\}$ and $f\colon \N\to\N$
is monotone, ordered by $(n,f)\leq (m,g)$ iff $n\leq m$ and $f\leq g$.

Show that the net $(\, n\ketbra{f(n)}{f(n)}\,)_{n,f\in D}$
converges ultrastrongly to~$0$.

So a net which converges ultrastrongly need not be bounded!
(The cure for this pathology is Kaplansky's density theorem, 
see~\sref{kaplansky}.)

Show that $\frac{1}{n} \ketbra{f(n)}{0}$
converges ultrastrongly to~$0$ as $D\ni(n,f)\to \infty$.

Show that the product
$\ketbra{f(n)}{0} = (\,n\ketbra{f(n)}{f(n)}\,)\,(\,\frac{1}{n}
\ketbra{f(n)}{0}\,)$
does not converge ultrastrongly 
as $D\ni(n,f)\to\infty$.

Conclude that multiplication $a,b\mapsto ab$
is not jointly ultrastrongly continuous on~$\scrB(\ell^2)$,
even when~$b$ is restricted to a bounded set.

(Nevertheless we'll see that multiplication is ultrastrongly continuous
when~$a$ is restricted to a bounded set in~\sref{mult-jus-cont}.)

\item
Show that
$a_{n,f} = \frac{1}{n}(\ketbra{f(n)}{0}+\ketbra{0}{f(n)})
\,+\, n\ketbra{f(n)}{f(n)}$
converges ultrastrongly to~$0$
as $D\ni(n,f)\to\infty$,
while $a_{n,f}^2$ does not.

Hence~$a\mapsto a^2$ is not ultrastrongly continuous on~$\sa{\scrB(\ell^2)}$.

\item
Let us show that~$\scrB(\ell^2)$
is not ultraweakly complete.

Show that there is an unbounded linear map~$f\colon \ell^2\to\C$
(perhaps using the fact that every vector space
has a basis by the axiom of choice),
and that for each finite dimensional linear subspace~$S$ of~$\ell^2$
there is a unique vector~$x_S\in S$ 
with
$f(x)=\left<x_S,y\right>$ for all~$y\in S$
(using~\sref{riesz-representation-theorem}).

Consider the net~$(\,\ketbra{e}{x_S}\,)_S$
where~$S$ ranges over the finite dimensional subspaces of~$\ell^2$
ordered by inclusion,
and~$e$ is some fixed vector in~$\ell^2$ with~$\|e\|=1$.

Let~$\omega\colon \scrB(\ell^2)\to\C$
be an np-map,
so $\omega\equiv \sum_n \left<y_n,(\,\cdot\,)y_n\right>$
for $y_1,y_2,\dotsc \in \ell^2$ with $\sum_n \|y_n\|^2 <\infty$,
see~\sref{bh-np}.

Show that $\omega(\,\ketbra{e}{x_S}-\ketbra{e}{x_T}\,)
= \left<\,x_S-x_T,\,\sum_n y_n\left<y_n,e\right>\,\right> = 0$
when $S$ and~$T$ are finite dimensional linear subspaces of~$\ell_2$
which contain the vector $\sum_n y_n\left<y_n,e\right>$.

Conclude that~$(\,\ketbra{e}{x_S}\,)_S$
is ultraweakly Cauchy.

Show that if~$(\,\ketbra{e}{x_S}\,)_S$
converges ultraweakly to some~$A$ in~$\scrB(\ell^2)$,
then we have~$\left<e,Ay\right>=f(y)$
for all~$y\in\ell^2$.

Conclude that~$(\,\ketbra{e}{x_S}\,)_S$
does not converge ultraweakly,
and that~$\scrB(\ell^2)$ is not ultraweakly complete.

(Nevertheless, we'll see that every von Neumann algebra
is ultrastrongly complete, and that
every norm-bounded ultraweakly Cauchy net
in a von Neumann converges, in~\sref{vn-complete}.)
\end{enumerate}
\spacingfix%
\end{point}%
\end{parsec}%
\subsection{Elementary Theory}
%
%
\begin{parsec}{440}%
\begin{point}{10}%
The basic facts concerning von Neumann algebras
we'll deal with first mostly involve the
relationship
between
multiplication
and the order structure.
For example,
while it is clear that translation and scaling
on a von Neumann algebra
are ultraweakly (and ultrastrongly) continuous,
the fact
that multiplication is ultraweakly (and ultrastrongly)
continuous in each coordinate is
less obvious (see~\sref{mult-uws-cont}).
Quite surprisingly,
this problem reduces to the ultraweak continuity
of $b\mapsto a^*ba$ by the following identity.
\end{point}
\begin{point}{20}[mult-polarization]{Exercise}%
\index{polarisation identity!in a von Neumann algebra}
Show that for elements~$a,b,c$ of a $C^*$-algebra,
\begin{equation*}
\textstyle
a^*\,c\,b\ =\ \frac{1}{4}\,\sum_{k=0}^3\ i^k\  (i^ka+b)^*\,c\,(i^ka+b).
\end{equation*}
(Note that this identity is a variation on the polarisation
identity for inner products,
see~\sref{inner-product-basic}.)
\end{point}
\begin{point}{30}[vanishing-effects]{Lemma}%
Let~$(x_\alpha)_{\alpha\in D}$ be 
a net of effects of a von Neumann algebra~$\scrA$,
which converges ultraweakly to~$0$.
Let~$(b_\alpha)_{\alpha\in D}$ be a 
net of elements with~$\|b_\alpha\| \leq 1$ for all~$\alpha$.
Then $(x_\alpha b_\alpha)_\alpha$ converges ultraweakly
to~$0$.
\begin{point}{40}{Proof}%
Let~$\omega\colon \scrA\to \C$ be an np-map.
We have, for each~$\alpha$,
\begin{alignat*}{3}
\left|\,\omega(x_\alpha b_\alpha)\,\right|^2
\ &=\ 
\left|\, \omega(\,\sqrt{x_\alpha}\,\sqrt{x_\alpha}\,b_\alpha\,)\, \right|^2
\qquad&&\text{since $x_\alpha\geq 0$}\\
\ &\leq\ 
\omega(x_\alpha)\  \omega(\,b_\alpha^* x_\alpha b_\alpha\,) 
\qquad&&\text{by Kadison's inequality, \sref{omega-norm-basic}}\\
\ &\leq\ 
\omega(x_\alpha)\ \omega(b_\alpha^* b_\alpha)
\qquad&&\text{since $x_\alpha\leq 1$}\\
\ &\leq\ 
\omega(x_\alpha)\ \omega(1)
\qquad&&\text{since $b_\alpha^*b_\alpha\leq 1$}.
\end{alignat*}
Thus,
since $(\omega(x_\alpha))_\alpha$
converges to~$0$,
we see that $(\omega(x_\alpha b_\alpha))_\alpha$
converges to~$0$,
and so $(x_\alpha b_\alpha)_\alpha$ converges ultraweakly to~$0$.\qed
\end{point}
\end{point}
\begin{point}{50}{Exercise}%
Let~$D$ be a bounded directed set of self-adjoint
elements of a von Neumann algebra~$\scrA$,
and let~$a\in \scrA$.
\begin{point}{60}[vna-supremum-uwlimit]%
Show that the net~$(d)_{d\in D}$ converges ultraweakly to~$\bigvee D$.
\end{point}
\begin{point}{70}[vna-supremum-mult]%
Use~\sref{vanishing-effects}
to show that $(da)_d$ converges ultraweakly to~$(\bigvee D)a$,
and that~$(a^*d)_d$ converges ultraweakly to~$a^* (\bigvee D)$.
\end{point}
\end{point}
%
%
\begin{point}{80}[ad-normal]{Proposition}%
Let~$a$ be an element of a von Neumann algebra~$\scrA$.
Then
\begin{equation*}
	\textstyle
	\bigvee_{d\in D} a^*\,d\,a \ =\  a^*\,(\bigvee D)\, a
\end{equation*}
for every bounded directed subset~$D$ of self-adjoint
elements of~$\scrA$.
\begin{point}{90}[ad-normal-1]{Proof}%
If~$a$ is invertible,
then the (by~\sref{astara-pos-basic-consequences}) order preserving map $b\mapsto a^*ba$
has an order preserving inverse (namely $b\mapsto (a^{-1})^* b a^{-1}$),
and therefore preserves all suprema.
\begin{point}{100}%
The general case reduces to the case that~$a$ 
is invertible
in the following way.
There is (by~\sref{spectrum-bounded})
 $\lambda>0$ such that $\lambda+a$ is invertible.
Then as $d$ increases 
\begin{equation*}
a^*\,d\,a \ \equiv\  (\lambda+a)^*\,d\,(\lambda+a) \,-\,
 \lambda^2d \,-\, \lambda a^*d \,-\, \lambda da
\end{equation*}
converges ultraweakly
to~$a^* \,(\bigvee D)\,a$,
because $(\ (\lambda+a)^*\,d\,(\lambda+a)\ )_d$
converges ultraweakly to $(\lambda+a)^*\,(\bigvee D)\,(\lambda+a)$
by~\sref{ad-normal-1} and~\sref{vna-supremum-uwlimit},
and $(a^*d+da)_d$ converges ultraweakly to $a^*(\bigvee D)+(\bigvee D)a$
by~\sref{vna-supremum-mult}.
Since~$(a^*da)_d$ converges to~$\bigvee_{d\in D} a^*d a$ too,
we could conclude that
$\bigvee_{d\in D} a^* \,d\, a = a^*\,(\bigvee D)\,a$
if we would already know that the ultraweak topology is Hausdorff.
At the moment, however,
we must content ourselves with
the conclusion that
	$\omega(\,a^*(\bigvee D) a\,-\, \bigvee_{d\in D} a^* d a\,)=0$
for every np-functional~$\omega$ on~$\scrA$.
But since
	$a^*(\bigvee D) a - 
\bigvee_{d\in D} a^* da$
happens to be positive,
we conclude that
	$a^*(\bigvee D) a  
- \bigvee_{d\in D} a^* d a =0$
nonetheless.\qed
\end{point}
\end{point}
\end{point}
\begin{point}{110}[vn-positive-basic]{Exercise}%
Show that the set of np-functionals
on a von Neumann algebra~$\scrA$
is not only faithful 
but also order separating
using~\sref{proto-gelfand-naimark}.
Deduce
\begin{enumerate}%
\item%
\index{ultraweak and ultrastrong!topologies are Hausdorff}%
that the ultraweak and ultrastrong topologies
are Hausdorff,
\item
that~$\scrA_+$, 
$\sa{\scrA}$ and~$[0,1]_\scrA$  are ultraweakly 
(and ultrastrongly) closed, 
\item
and that the unit ball
$(\scrA)_1$
is ultrastrongly closed.

(We'll see only later on, in~\sref{ultraclosed},
that~$(\scrA)_1$
is ultraweakly closed as well.)
\end{enumerate}
\spacingfix%
\end{point}%
\begin{point}{120}{Exercise}%
Let~$D$ be a directed subset of self-adjoint elements
of a von Neumann algebra~$\scrA$,
and let~$a\in\scrA$.
\begin{point}{130}[vna-supremum-commutes]%
Show that if~$ad=da$ for all~$d\in D$,
then $a(\bigvee D) = (\bigvee D)a $.
\end{point}
\begin{point}{140}[vna-supremum-uslimit]%
Use~\sref{vanishing-effects}
to show that $(\bigvee D-d)^2$ converges ultraweakly to~$0$
as $D\ni d\to\infty$.

Conclude that~$(d)_{d\in D}$ converges ultrastrongly to~$\bigvee D$.
\end{point}
\end{point}
\begin{point}{150}[p-uwcont]{Exercise}%
Show that for a positive linear map $f\colon \scrA\to\scrB$
between von Neumann algebras,
the following are equivalent.
\begin{enumerate}
\item
$f$ is ultraweakly continuous;
\item
$f$ is ultraweakly continuous on~$[0,1]_\scrA$;
\item
$f(\bigvee D)=\bigvee_{d\in D}f(d)$ for each bounded 
directed~$D\subseteq\sa{\scrA}$;
\item 
$\omega\circ f\colon \scrA\to\C$ is normal 
for each np-map $\omega\colon \scrB\to\C$.
\end{enumerate}
In that case we say that~$f$ is \Define{\textbf{n}ormal}.%
\index{normal!positive map between von Neumann algebras}
(Note that this definition of ``normal'' extends the one
for positive functionals
from~\sref{def-np-functional}.)

Conclude that $b\mapsto a^*ba,\,\scrA\to\scrA$%
\index{$a^*(\,\cdot\,)a\colon \scrA\to\scrA$!is normal}
is ultraweakly
continuous for every element~$a$ of a von Neumann 
algebra~$\scrA$.
\end{point}
\end{parsec}
\begin{parsec}{450}%
\begin{point}{10}{Exercise}%
Show that if a positive linear map $f\colon \scrA\to\scrB$
between von Neumann algebras 
is ultrastrongly continuous (on~$[0,1]_\scrA$), 
then~$f$ is normal.
(Hint: use that a bounded directed set 
$D\subseteq \sa{\scrA}$ converges ultrastrongly to~$\bigvee D$.)

The converse does not hold: give an example of a map~$f$ 
which is normal, but 
not ultrastrongly continuous. (Hint: transpose.)
\end{point}
\begin{point}{20}[cp-uscont]{Proposition}%
An ncp-map $f\colon \scrA\to\scrB$
between von Neumann algebras is 
ultrastrongly continuous.
\begin{point}{30}{Proof}%
Note that $f$ is ultrastrongly continuous at~$a\in\scrA$
iff $f((\,\cdot\,)+a)\equiv f + f(a)$ is ultrastrongly continuous at~$0$.
Thus to show that~$f$
is ultrastrongly continuous
it suffices to show that~$f$ is ultrastrongly continuous at~$0$.
So let~$(b_\alpha)_\alpha$ be a net in~$\scrA$
which converges ultrastrongly to~$0$;
we must show that $(f(b_\alpha))_\alpha$
converges ultrastrongly to~$0$, viz.~that
$(\,f(b_\alpha)^*f(b_\alpha)\,)_\alpha$ converges ultraweakly to~$0$.
Since
$f(b_\alpha)^*f(b_\alpha) \leq f(b_\alpha^*b_\alpha) \|f(1)\|$
by~\sref{cp-cs}, 
it suffices to show that~$(\,f(b_\alpha^*b_\alpha)\,)_\alpha$
converges ultraweakly to~$0$,
but this follows from the
facts that~$f$ is ultraweakly continuous
and~$(b_\alpha^*b_\alpha)_\alpha$
converges ultraweakly to~$0$
(since~$(b_\alpha)_\alpha$ converges ultrastrongly to~$0$).\qed
\end{point}
\end{point}
\begin{point}{40}[mult-uws-cont]{Exercise}%
Let~$\scrA$ be a von Neumann algebra.
Conclude (using~\sref{cp-uscont} and~\sref{ad-cp})
that the map $a\mapsto b^*ab,\,\scrA\to\scrA$
is ultrastrongly continuous for every
 element~$b\in\scrA$.

Use this,
and~\sref{mult-polarization},
to
show that $b\mapsto ab,\,ba\colon\ \scrA\to\scrA$
are ultraweakly and ultrastrongly continuous
for every element~$a$ of a von Neumann algebra~$\scrA$.
\end{point}
\begin{point}{50}%
We saw in~\sref{vn-counterexamples}
that the multiplication on a von Neumann algebra 
is not jointly ultraweakly continuous,
even on a bounded set.
Neither is $a,b\mapsto ab$ jointly ultrastrongly continuous,
even when~$b$ is restricted to a bounded set;
but it \emph{is} jointly 
ultrastrongly continuous when~$a$ is restricted to a bounded set:
\end{point}
\begin{point}{60}[mult-jus-cont]{Proposition}%
Let $(a_\alpha)_\alpha$
and~$(b_\alpha)_\alpha$
be nets
in a von Neumann algebra~$\scrA$
with the same index set
that converge ultrastrongly to~$a,b\in\scrA$, respectively.
Then the net~$(a_\alpha b_\alpha)_\alpha$
converges ultrastrongly to~$ab$
provided that~$(a_\alpha)_\alpha$
is bounded.
\begin{point}{70}{Proof}%
Let $\omega\colon \scrA\to\C$
be an np-functional.
Since
\begin{alignat*}{3}
\|ab-a_\alpha b_\alpha\|_\omega
\ &\leq\ 
	\|(a -a_\alpha)b\|_\omega
	\,+\, 
	\|a_\alpha(b-b_\alpha)\|_\omega
	\\
\ &\leq\ 
	\|a -a_\alpha\|_{\omega(b^*(\,\cdot\,)b)}
	\,+\, 
	\|a_\alpha\|\|b-b_\alpha\|_\omega
\end{alignat*}
vanishes as~$\alpha\to\infty$,
	we see that~$(a_\alpha b_\alpha)_\alpha$
converges ultrastrongly to~$ab$.\qed
\end{point}
\end{point}
\end{parsec}
\begin{parsec}{460}%
\begin{point}{10}%
We can now prove a bit more 
about the ultrastrong and ultraweak topologies.
\end{point}
\begin{point}{20}[usconv]{Exercise}%
Show that a net $(b_\alpha)_\alpha$ 
in a von Neumann algebra~$\scrA$
converges ultrastrongly to an element~$b$
of~$\scrA$
if and only if
both $b_\alpha^*b_\alpha\longrightarrow b^*b$
and~$b_\alpha\longrightarrow b$
ultraweakly as~$\alpha\to\infty$.
\end{point}
\begin{point}{30}[npuws]{Exercise}%
\index{normal!positive functional}%
Show that for a positive linear map $\omega \colon \scrA\to\C$
on a von Neumann algebra~$\scrA$
the following are equivalent
\begin{enumerate}
\item
	$\omega$ is normal;
\item
	$\omega$ is ultraweakly continuous;
\item
	$\omega$ is ultrastrongly continuous.
\end{enumerate}
(Hint: combine~\sref{p-uwcont} and \sref{cp-uscont}.)
\end{point}
\end{parsec}
\begin{parsec}{470}%
\begin{point}{10}%
Enter the eponymous hero(s) of this thesis.
\end{point}
\begin{point}{20}{Definition}%
We denote 
the category of
\emph{normal} cpsu-maps
by~$\Define{\W{cpsu}}$,
and its subcategory of nmiu-maps
by~$\Define{\W{miu}}$.
\index{Wmiu@$\W{miu}$, $\W{cpsu}$, \dots}%
(We omit the ``N'' for the sake of brevity.)
\begin{point}{30}%
Though arguably~$\W{miu}$
is a good candidate
for being called \emph{the} category of von Neumann algebra,
the title of this thesis refers to~$\W{cpsu}$.%
\index{category of von Neumann algebras}%
\index{von Neumann algebra!category of}
Indeed, it's the ncpsu-maps between von Neumann algebras
that stand to model the arbitrary quantum processes,
and it's the category of these quantum processes
we want to mine for abstract structure.
This is mostly a task for the next chapter,
though.
For now we'll just establish that~$\W{cpsu}$
has all products, \sref{vn-products}, certain equalisers,
\sref{vn-equalisers},
and that $\op{(\W{cpsu})}$ is an \emph{effectus}, see~\sref{vn-effectus}.
\end{point}
\end{point}
\begin{point}{40}[vn-products]{Exercise}%
\index{product!in $\W{miu}$ and $\W{cpsu}$}
Show that
given a family $(\scrA_i)_i$
of von Neumann algebras
the direct sum
$\bigoplus_i \scrA_i$
from~\sref{cstar-product}
is a von Neumann algebra
and  the projections
$\pi_j \colon \bigoplus_i \scrA_i\to\scrA_j$
are normal.
Moreover, show
that this makes~$\bigoplus_i \scrA_i$
into the  product of the~$\scrA_i$
in the categories~$\W{miu}$ and~$\W{cpsu}$
(see~\sref{cstar-product-2} and~\sref{cstar-product-4}).
\end{point}
\begin{point}{50}[vn-equalisers]{Exercise}%
\index{equaliser!in $\W{miu}$ and $\W{cpsu}$}
Show that given nmiu-maps $f,g\colon \scrA\to\scrB$
between von Neumann algebras
the set~$\scrE:=\{\,a\in\scrA\colon\, f(a)=g(a)\,\}$
is a von Neumann subalgebra of~$\scrA$,
and the inclusion~$e\colon \scrE\to\scrA$
is the equaliser of~$f$ and~$g$ in the 
categories~$\W{miu}$ and~$\W{cpsu}$
(see~\sref{cstar-equaliser-1} and~\sref{cstar-product-4}).
\end{point}
\begin{point}{60}[vn-effectus]
Let us briefly indicate what makes~$\op{(\W{cpsu})}$ an effectus;
for a precise formulation and proof of this fact
we refer to~\cite{effintro,kentapartial}
(or~\sref{effectus-vn}, 
\sref{effectus-in-partial-form}, and \sref{cho-thm} ahead).
Note that the sum $f+g$ of two ncpsu-maps
$f,g\colon \scrA\to\scrB$ between von Neumann algebras
is again an ncpsu-map
iff~$f(1)+g(1)\leq 1$.
The partial addition on ncpsu-maps
thereby defined
has, aside from some fairly obvious properties
(summarised by the
fact that the category~$\W{cpsu}$
    is $\Cat{PCM}$-enriched, see~\cite{kentapartial}),
the following special trait:
given ncpsu-maps $f\colon \scrA\to\scrD$
and $g\colon \scrB\to\scrD$
with $f(1)+g(1)\leq 1$
we may form an ncpsu-map $[f,g]\colon \scrA\times \scrB\to \scrD$
by $[f,g](a,b)=f(a)+g(b)$,
and, moreover,
every ncpsu-map $\scrA\times\scrB\to\scrD$ is of this form.
This observation, which gives the product of~$\W{cpsu}$ a coproduct-like quality
without forcing it to be a biproduct (which it's not),
	makes $\op{(\W{cpsu})}$ a \emph{FinPAC}\index{FinPAC} 
	(see~\sref{effectus-in-partial-form}).

	For~$\op{(\W{cpsu})}$ to be an effectus,\index{effectus}
we need a second ingredient:
 the complex numbers,~$\C$.
Since the ncpsu-maps $p\colon \C\to \scrA$
are all of the form $\lambda \mapsto \lambda a$
for some effect $a\in [0,1]_\scrA$,
the ncpsu-maps $p\colon \C\to\scrA$
(called \emph{predicates} in this context)
are not only endowed with a partial addition,
but even form an \emph{effect algebra}.
This,
combined with the observation that
an ncpsu-map $f\colon \scrA\to\scrB$
is constant zero iff~$f(1)=0$,
makes~$\op{(\W{cpsu})}$
an \emph{effectus in partial form} (see~\sref{effectus-in-partial-form}).

As you can see, there's nothing deep underlying~$\op{(\W{cpsu})}$
being an effectus.
In that respect effectus theory resembles topology:
just as a topology 
provides a basis for  
notions such as compactness, connectedness, meagerness, and homotopy,
so does an
effectus provide a framework to study
 aspects of  computation
such as side effects (\sref{sefp}) and purity
(\sref{pure-effectus}).
\end{point}
\end{parsec}
\begin{parsec}{480}%
\begin{point}{10}%
Let us quickly prove that every von Neumann algebra
is isomorphic to a von Neumann algebra of operators on a 
Hilbert space (see~\sref{ngns}).
\end{point}
\begin{point}{20}[normal-faithful]{Exercise}%
\index{normal!positive map between von Neumann algebras}
Let~$\Omega$ be a collection of np-functionals
on a von Neumann algebra~$\scrB$
that is faithful (see~\sref{separating}).
Show that a positive linear map $f\colon \scrA\to\scrB$
is normal iff~$\omega\circ f$
is normal for all~$\omega\in\Omega$.
\end{point}
\begin{point}{30}{Proposition}%
\index{rhoomega@$\varrho_\omega$!is normal}
Given an np-map~$\omega\colon \scrA \to \C$
on a von Neumann algebra~$\scrA$,
the map $\varrho_\omega\colon \scrA\to\scrB(\scrH_\omega)$
from~\sref{gns} is normal.
\begin{point}{40}{Proof}%
Since by definition of~$\scrH_\omega$
the vectors of the form~$\eta_\omega(a)$
where~$a\in\scrA$
are dense in~$\scrH_\omega$,
the vector functionals
$\left<\eta_\omega(a),(\,\cdot\,)\eta_\omega(a)\right>$
form a faithful collection
of np-functionals on~$\scrB(\scrH_\omega)$.
Thus by~\sref{normal-faithful}
it suffices to show given~$a\in\scrA$
that $\left<\eta_\omega(a),\varrho_\omega(\,\cdot\,)
\eta_\omega(a)\right>\equiv \omega(a^*(\,\cdot\,)a)$ is normal,
which it is, by~\sref{ad-normal}.\qed
\end{point}
\end{point}
\begin{point}{50}[varrho-Omega-normal]{Exercise}%
\index{rhoOmega@$\varrho_\Omega$!is normal}
Show that the map $\varrho_\Omega$ 
from~\sref{gelfand-naimark-representation}
is normal for every collection~$\Omega$ of np-maps $\scrA\to\C$
on a von Neumann algebra $\scrA$.
\end{point}
\begin{point}{60}[injective-nmiu-iso-on-image]{Lemma}%
Let~$f\colon \scrA\to\scrB$ be an injective nmiu-map
between von Neumann algebras.
Then the image~$f(\scrA)$ is a von Neumann subalgebra of~$\scrB$,
and~$f$ restricts to an nmiu-isomorphism from~$\scrA$
to~$f(\scrA)$.
\begin{point}{70}{Proof}%
We already know by~\sref{injective-miu-iso-on-image}
that~$f(\scrA)$ is a $C^*$-subalgebra of~$\scrA$,
and that~$f$ restricts to an miu-isomorphism~$f'\colon \scrA\to f(\scrA)$.
The only thing left to show is that~$f(\scrA)$
is a von Neumann subalgebra of~$\scrB$,
because an miu-isomorphism between von Neumann algebras
(being an order isomorphism)
will automatically be an nmiu-isomorphism.
Let~$D$ be a bounded directed subset of~$f(\scrA)$.
Note that~$S:=(f')^{-1}(D)$ is a bounded
directed subset of~$\scrA$,
and so~$\bigvee D\equiv  \bigvee f(\,S\,)
= f(\bigvee S)$, because~$f$ is normal.
Thus~$\bigvee f(D)\in f(\scrA)$,
and so $f(\scrA)$ is a von Neumann subalgebra of~$\scrB$.\qed
\end{point}
\end{point}
\begin{point}{80}[ngns]{Theorem (normal Gelfand--Naimark)}%
\index{Gelfand--Naimark's Theorem!for von Neumann algebras}%
Every von Neumann algebra~$\scrA$ is nmiu-isomorphic
to von Neumann algebra of operators on a Hilbert space.
\begin{point}{90}[ngns-proof]{Proof}%
Recall that an element $a\in \scrA$ is zero iff $\omega(a)=0$
for all np-maps $\omega\colon \scrA\to\C$.
It follows that the collection~$\Omega$
of all np-maps $\scrA\to\C$
obeys the condition of~\sref{proto-gelfand-naimark},
and so the miu-map $\varrho_\Omega\colon \scrA\to\scrB(\scrH_\Omega)$
(from~\sref{gelfand-naimark-representation})
is injective.
Since~$\varrho_\Omega$
is also normal by~\sref{varrho-Omega-normal},
we see by~\sref{injective-nmiu-iso-on-image} that~$\varrho_\Omega$
restricts to an nmiu-isomorphism
from~$\scrA$ to the von Neumann subalgebra~$\varrho_\Omega(\scrA)$
of~$\scrB(\scrH_\Omega)$.\qed
\end{point}
\end{point}
\end{parsec}
\subsection{Examples}
\subsubsection{Matrices over von Neumann algebras} 
\begin{parsec}{490}%
\begin{point}{10}%
We'll show that the $C^*$-algebra
of $N\times N$-matrices~$M_N(\scrA)$
over a von Neumann algebra~$\scrA$
is itself a von Neumann algebra,
and to this end,
we prove something a bit more more general.
\end{point}
\begin{point}{20}[bah-vn]{Theorem}%
Given a von Neumann algebra~$\scrA$,
the~$C^*$-algebra $\scrB^a(X)$ (\sref{bax-cstar})%
\index{BaX@$\scrB^a(X)$!as a von Neumann algebra}
of bounded adjointable module maps on
a self-dual (\sref{self-dual}) Hilbert $\scrA$-module~$X$
is a  von Neumann algebra,
and
$\left<x,(\,\cdot\,)x\right>\colon \scrB^a(X)\to\scrA$
is normal for every~$x\in X$.%
\index{vector functional!is normal}
\begin{point}{30}{Proof}%
We'll first show that
a  bounded directed subset~$\scrD$ of~$\Real{\scrB^a(X)}$
has a supremum (in~$\Real{\scrB^a(X)}$).
To obtain a candidate for this supremum,
we first define
a bounded form~$[\,\cdot\,,\,\cdot\,]\colon X\times X\to\scrA$
in the sense of~\sref{chilb-form}
and apply~\sref{chilb-form-representation}.
To this end note that given~$x\in X$
the subset~$\{\,\left<x,Tx\right>\colon \, T\in\scrD\,\}$
of~$\Real{\scrA}$
is bounded and directed,
and so 
(since~$\scrA$ is a von Neumann algebra)
has a supremum.
Since the 
net $(\,\left<x,Tx\right>\,)_{T\in \scrD}$
converges ultraweakly to this supremum
by~\sref{vna-supremum-uwlimit},
we see that
$\left<y,Tz\right>\,= \,\frac{1}{4}
\sum_{k=0}^3
i^k\left<y+i^kz,\smash{T(y+i^kz)}\right>$
converges ultraweakly to 
some element~$[y,z]$ of~$\scrA$
as~$T\to\infty$ for all~$y,z\in X$,
giving us a form~$[\,\cdot\,,\,\cdot\,]$ on~$X$.
Since~$\left\|\left<y,T z \right>\right\|
\leq \sup_{T'\in \scrD} \|T'\| \|y\|\|z\|$
for all~$T\in\scrD$
by~\sref{chilb-form-bounded},
and thus~$\|[y,z]\|\leq \sup_{T'\in\scrD} \|T'\|\|y\|\|z\|$
for all~$y,z\in X$,
we see that the form~$[\,\cdot\,,\,\cdot\,]$
is bounded.
Since~$X$ is self dual,
there is
by~\sref{chilb-form-representation}
$S\in \scrB^a(X)$
with
$[y,z]=\left<y,Sz\right>$
for all~$y,z\in X$;
we'll show that~$S$ is the supremum of~$\scrD$.

To begin,
given~$T\in\scrD$
we have~$\left<x,Tx\right>\leq\bigvee_{T'\in\scrD} \left<x,T'x\right>
=[x,x]=\left<x,Sx\right>$
for all~$x\in X$,
and so~$T\leq S$
by~\sref{chilb-vector-states-order-separating},
that is, $S$ is an upper bound for~$\scrD$.
Given another upper bound~$S'\in\Real{\scrB^a(X)}$
of~$\scrD$
(so  $T\leq S'$ for all~$T\in\scrD$)
we have~$\left<x,Tx\right>\leq \left<x,S'x\right>$
and so~$\left<x,Sx\right>
=[x,x]=\bigvee_{T\in \scrD}\left<x,Tx\right>\leq \left<x,S'x\right>$
for all~$x\in X$
implying that~$S\leq S'$.
Hence~$S$ is the supremum of~$\scrD$
in $\Real{\scrB^a(X)}$.
Note that since~$\left<x,Sx\right>=\bigvee_{T\in\scrD} \left<x,Tx\right>$
we immediately see
that $\left<x,(\,\cdot\,)x\right>\colon \scrB^a(X)\to\scrA$
preserves bounded directed suprema
for every~$x\in X$.

It remains to be shown that
there are sufficiently many np-functionals on $\scrB^a(X)$
in the sense that~$T\in(\scrB^a(X))_+$ is zero
when~$\omega(T)=0$ for every np-functional $\omega\colon \scrB^a(X)\to\C$.
This is indeed
the case for such an operator~$T$,
because~$\xi(\left<x,(\,\cdot\,)x\right>)$
is an np-functional on~$\scrB^a(X)$
for every~$x\in X$
and an np-functional $\xi\colon \scrA\to\C$,
implying that~$\xi(\left<x,Tx\right>)=0$,
and~$\left<x,Tx\right>=0$,
and so~$T=0$.\qed
\end{point}
\end{point}
\begin{point}{40}[mn-vna]{Exercise}%
\index{$M_n\scrA$, the $n\times n$-matrices over~$\scrA$!as a von Neumann algebra}%
Let~$\scrA$ be a von Neumann algebra,
and let~$N$ be a natural number.
\begin{enumerate}
\item
Show 
that the $C^*$-algebra
$M_N(\scrA)$ of $N\times N$-matrices over~$\scrA$ (see~\sref{cstar-matrices})
is a von Neumann algebra.
\item
Show that
the map $A\mapsto \sum_{ij} a_i^* A_{ij} a_j\colon\, M_N\scrA\to \scrA$
is normal and completely positive,
and that
the map $A\mapsto \sum_{ij} a_i^* A_{ij} b_j\colon\, M_N\scrA\to\scrA$
is ultrastrongly and ultraweakly continuous
for all~$a_1,\dotsc,a_N,b_1,\dotsc,b_N\in\scrA$.

In particular, $A\mapsto A_{ij}\colon\,M_N\scrA\to\scrA$
is ultraweakly and ultrastrongly continuous
for all~$i,j$.

Show that a net~$(A_\alpha)_\alpha$
in~$M_N\scrA$
converges ultraweakly (ultrastrongly)
to $B\in M_N\scrA$
iff $(A_\alpha)_{ij}$ converges ultraweakly (ultrastrongly)
to~$B_{ij}$ as~$\alpha\to\infty$ for all~$i,j$.
\item
Given an ncp-map $f\colon \scrA\to\scrB$
between von Neumann algebras,
show that the cp-map $M_N f\colon M_N\scrA\to M_N\scrB$
from~\sref{mnf} is normal.%
\index{$M_nf$!is normal}
\end{enumerate}
\spacingfix%
\end{point}%
\end{parsec}%
\subsubsection{Commutative von Neumann algebras}
\begin{parsec}{500}%
\begin{point}{10}[linfty-example]%
Another important source of examples of von 
Neumann algebras is measure theory:
we'll show that the bounded measurable
functions on a 
finite complete
measure space~$X$
form a commutative von Neumann algebra~$L^\infty(X)$
when functions that are equal almost everywhere are identified
(see~\sref{Linfty-vn}).
In fact, 
we'll see in~\sref{cvn}
that
every commutative von Neumann algebra is nmiu-isomorphic
to a direct sum of~$L^\infty(X)$s.
This is not only interesting in its own right,
but will also
be used in the next chapter
to show that the only von Neumann algebras
that can be endowed with a `duplicator'
are of the form~$\ell^\infty(X)$
for some set~$X$ (see~\sref{duplicable}).

We should probably mention that~$L^\infty(X)$
can be defined for any measure space~$X$,
and is a von Neumann algebra
precisely when~$X$ is \emph{localisable},
see~\cite{segal1951}.
This has the advantage that any commutative von Neumann algebra
is nmiu-isomorphic
to a single~$L^\infty(X)$
for some localisable measure space~$X$,
but since it has no other advantages relevant
to this text
we restrict ourselves to complete finite measure spaces.

We'll assume the reader is reasonably familiar with 
the basics of measure theory,
and we'll only show a selection of results
that we deemed important.
For the other details,
we refer to volumes~1 and~2 of~\cite{fremlin}.
Nevertheless,
we'll  recall
some basic definitions
to fix terminology,
which is sometimes simpler than in~\cite{fremlin}
(because we're dealing with finite complete measure spaces),
and sometimes modified to the complex-valued case
(c.f.~133C of~\cite{fremlin}).
A motivated reader will have no problem adapting
the results from~\cite{fremlin} 
to our setting.
\end{point}
\end{parsec}%
\begin{parsec}{510}[measure-theory-recap]%
\begin{point}{10}%
Let~$X$ be a finite and complete measure space.
We'll denote the $\sigma$-algebra
of measurable subsets of~$X$ by~$\Define{\Sigma_X}$,%
\index{SigmaX@$\Sigma_X$, measurable subsets}
and the measure by~$\Define{\mu_X}\colon \Sigma_X\to[0,\infty)$ (or~$\mu$%
\index{muX@$\mu_X$, measure}
when no confusion is expected).
That~$X$ is \Define{finite}%
\index{measure!finite}
means that~$\mu(X)<\infty$
(which doesn't mean that the set~$X$ is finite),
and that~$X$ is \Define{complete}%
\index{measure!complete}
means that every subset~$A$
of a negligible subset~$B$ of~$X$
is itself negligible.
(Recall that~$N\subseteq X$
is \Define{negligible}%
\index{negligible, subset of a measure space}
when~$N\in \Sigma_X$ and~$\mu(N)=0$.)
A function $f\colon X\to\C$ is \Define{measurable}%
\index{measurable function}
when the inverse image~$f^{-1}(U)$
of any open subset~$U$ of~$\C$
is measurable
(which happens precisely
when both~$x\mapsto \Real{f(x)},\,x\mapsto\Imag{f(x)}\colon\,
X\to\R$
are measurable in the sense of~121C of~\cite{fremlin}).
An important example of a measurable function on~$X$
is the indicator function~$\mathbf{1}_A$
of a measurable subset~$A$ of~$X$
(which is equal to~$1$ on~$A$ and~$0$ elsewhere.)
\end{point}
\begin{point}{20}%
The bounded measurable functions $f\colon X\to\C$
form a $C^*$-subalgebra of~$\C^X$ 
that we'll denote by~$\Define{\mathcal{L}^\infty(X)}$.%
\index{LcalinftyX@$\mathcal{L}^\infty(X)$}
The space $\mathcal{L}^\infty(X)$
is not only closed with respect to the (supremum) norm
on~$\C^X$,
but also with respect to coordinatewise limits
of bounded \emph{sequences}
(c.f.~121F of~\cite{fremlin}).
As a result,
the coordinatewise
(countable)
supremum~$\bigvee_n f_n$
of a bounded ascending sequence $f_1\leq f_2\leq \dotsb$
in~$\Real{\mathcal{L}^\infty(X)}$
is again in~$\mathcal{L}^\infty(X)$,
and is fact the supremum of~$(f_n)_n$ in~$\mathcal{L}^\infty(X)$.
However $\mathcal{L}^\infty(X)$
might still not be a von Neumann algebra
because not every bounded directed subset of~$\Real{\mathcal{L}^\infty(X)}$
might have a supremum
as we'll show presently; 
this is why we'll move from~$\mathcal{L}^\infty(X)$
to~$L^\infty(X)$ in a moment.
\end{point}
\begin{point}{30}%
For a counterexample
to~$\mathcal{L}^\infty(X)$ being always
a von Neumann algebra
we take~$X$ to be the unit interval~$[0,1]$
with the Lebesgue measure.
Let~$A$ be a non-measurable
subset of~$[0,1]$ (see 134B of~\cite{fremlin}).
The indicator functions~$\mathbf{1}_F$
where~$F$ is a finite subset of~$A$
form a bounded directed subset~$D$ of~$\Real{\mathcal{L}^\infty([0,1])}$
that---so we claim---has no supremum.
Indeed, note that
since $f\in \Real{\mathcal{L}^\infty([0,1])}$
is an upper bound for~$D$ iff~$\mathbf{1}_A\leq f$,
the least upper bound~$h$ for~$D$
would be the least bounded measurable function above~$\mathbf{1}_A$.
Surely, $h\neq \mathbf{1}_A$
for such~$h$
(because otherwise~$A$ would be measurable),
so~$h(x)>\mathbf{1}_A(x)$ for some~$x\in [0,1]$.
But then $h - (h(x)-\mathbf{1}_A(x))\mathbf{1}_{\{x\}}<h$
is an upper bound for~$D$ too
contradicting the minimality of~$h$.
Whence~$\mathcal{L}^\infty([0,1])$
is not a von Neumann algebra.
\end{point}
\begin{point}{40}%
To deal with~$L^\infty(X)$
we need to know a bit more about~$\mathcal{L}^\infty(X)$,
namely that the measure on~$X$ can be extended to a
an integral $\int\colon \mathcal{L}^\infty(X)\to \C$
(see~122M of~\cite{fremlin})%
\footnote{Note that every element of~$\mathcal{L}^\infty(X)$
being bounded is integrable by~122P of~\cite{fremlin}.}
with the following properties.
\begin{enumerate}
\item
$\int(\mathbf{1}_A) = \mu(A)$
for every measurable subset~$A$ of~$X$.
\item
$\int\colon \mathcal{L}^\infty(X)\to \C$
is a positive linear map
(see 122O of~\cite{fremlin}).
\item
$\int \bigvee_n f_n = \bigvee_n \int f_n$
for every bounded sequence
$f_1\leq f_2\leq \dotsb$
in~$\Real{\mathcal{L}^\infty(X)}$.
(This is a special case
of Levi's theorem, see 123A of~\cite{fremlin}.)
\end{enumerate}
Unsurprisingly, the integral
interacts poorly
with the uncountable directed suprema
that do exist in~$\mathcal{L}^\infty(X)$:
for example, the set~$D:=\{\,f\in [0,1]_{\mathcal{L}^\infty(X)}\colon \,
\int f = 0\,\}$
is directed, bounded, and has supremum~$1$,
but~$\bigvee_{f\in D} \int f=0 < 1=\int \bigvee D$.
What \emph{is} surprising
is that the lifting of~$\int$
to~$L^\infty(X)$ will be normal.
\end{point}
\begin{point}{50}%
But let us first define~$L^\infty(X)$.
We say that $f,g\in \mathcal{L}^\infty(X)$
are \Define{equal almost everywhere}
and write $\Define{f\approx g}$%
\index{((approx@$\approx$! $f\approx g$, for measurable functions}
when~$f(x)=g(x)$ for almost all~$x\in X$
(that is, $\{\,x\in X\colon\, f(x)\neq g(x)\,\}$
is negligible).
It is easily seen that~$\approx$ is an equivalence relation;
we denote the equivalence class
of a function~$f\in \mathcal{L}^\infty(X)$
by $\Define{f^\circ}$,%
\index{*supcirc@$f^\circ$, equivalence class of~$f$}
and
the set of equivalence classes by~$\Define{L^\infty(X)}
:= \{\, f^\circ\colon\, f\in \mathcal{L}^\infty(X)\,\}$,
	\index{LinftyX@$L^\infty(X)$}
which
becomes a commutative $C^*$-algebra
when
endowed with the same operations
as~$\mathcal{L}^\infty(X)$, but with
a slightly modified norm given by,
for~$\mathfrak{f}\equiv f^\circ \in L^\infty(X)$,
\begin{alignat*}{3}
	\|\mathfrak{f}\| \ &= \ 
\min\{\, \|g\|\colon \, g\in\mathcal{L}^\infty(X)\text{ and }
g^\circ = \mathfrak{f}\,\} \\
&=\ \min\{\ \lambda\geq 0\colon\,
		\left|f(x)\right|\,\leq\,\lambda
	\text{ for almost all~$x\in X$}\ \}.
\end{alignat*}
This is called the \Define{essential supremum norm}.%
\index{essential supremum norm}
To see that~$L^\infty(X)$ is complete
one can use the fact that~$\mathcal{L}^\infty(X)$
is complete in a slightly more general sense than discussed before:
when a bounded sequence $f_1,f_2,\dotsc$
in~$\mathcal{L}^\infty(X)$
converges coordinatewise for almost all~$x\in X$
to some bounded function~$f\colon X\to\C$,
this function~$f$ is itself measurable (and so~$f\in\mathcal{L}^\infty(X)$,
c.f.~121F of~\cite{fremlin}).

Another consequence of this
is that a bounded ascending sequence
$f_1^\circ \leq f_2^\circ \leq \dotsb$
in~$L^\infty(X)$
(so~$f_1,f_2,\dotsc\in\mathcal{L}^\infty(X)$,
and $f_1(x)\leq f_2(x)\leq\dotsb$ for almost all~$x\in X$)
has a supremum~$\bigvee_n f_n^\circ$ in~$\mathcal{L}^\infty(X)$.
Indeed, we'll have
$\bigvee_n f_n^\circ = g^\circ$
for any bounded map~$g\colon X\to\C$
with $g(x)=\bigvee_n f_n(x)$
for almost all~$x\in X$.
\end{point}%
\begin{point}{60}%
Now, let us return to the integral.
Since~$\int f = \int g$
for all~$f,g\in \mathcal{L}^\infty(X)$
with $f\approx g$
we get a map
$\int \colon L^\infty(X)\to \C$
given
by~$\int  f^\circ = \int f$.
Clearly, $\int$ is positive and linear,
and by (a slightly less special case
of) Levi's theorem (123A of~\cite{fremlin})
we see that $\int \bigvee_n \mathfrak{f}_n
= \bigvee_n \int \mathfrak{f}_n$
for any bounded ascending sequence
$\mathfrak{f}_1\leq \mathfrak{f}_2\leq\dotsb$
in~$\Real{L^\infty(X)}$.
Note that~$\int\colon L^\infty(X)\to\C$
is also faithful,
because if~$\int f^\circ=\int f = 0$
for some $f\in \mathcal{L}^\infty(X)_+$,
then~$f(x)=0$ for almost all~$x\in X$,
and so~$f^\circ=0$.
Now,
the fact that~$L^\infty(X)$
is a von Neumann algebra follows from
the following general and rather surprising observation.
\end{point}
\begin{point}{70}{Proposition}%
\index{von Neumann algebra!with a faithful np-functional}
Let~$\scrA$ be a $C^*$-algebra,
and let~$\tau\colon \scrA\to\C$
be a faithful positive map.
If every bounded ascending sequence~$a_1\leq a_2\leq \dotsb$
of self-adjoint elements from~$\scrA$ has a supremum~$\bigvee_n a_n$
(in~$\Real{\scrA}$)
and
$\tau(\bigvee_n a_n)=\bigvee_n \tau (a_n)$,
then~$\scrA$ is a von Neumann algebra,
and~$\tau$ is normal.
\begin{point}{80}{Proof}%
Our first task is to show that a bounded
directed subset~$D$ of self-adjoint elements of~$\scrA$
has a supremum~$\bigvee D$ in~$\Real{\scrA}$.
Since~$\bigvee_{d\in D} \tau(d)$
is a supremum in~$\R$
we can find~$a_1\leq a_2\leq \dotsb$ in~$D$
with~$\bigvee_n \tau(a_n) = \bigvee_{d\in D} \tau(d)$.
We'll show that~$\bigvee_n a_n$
is the supremum of~$D$.
Surely, any upper bound of~$D$ being also an upper bound 
for~$a_1\leq a_2\leq\dotsb$
is above~$\bigvee_n a_n$,
so the only thing that we need to show is that
$\bigvee_n a_n$ is an upper bound of~$D$.
So let~$b\in D$ be given.
The trick is to pick a sequence
$b_1\leq b_2\leq \dotsb$ in~$D$
with $b\leq b_1$ and~$a_n\leq b_n$ for all~$n$
(which exists on account of~$D$'s directedness).
Then~$\bigvee_n a_n \leq \bigvee_n b_n$,
and~$\bigvee_{d\in D} \tau(d) =\bigvee_n \tau(a_n)
=\tau(\,\bigvee_n a_n\,)
\leq \tau(\,\bigvee_n b_n\,)
= \bigvee_n \tau(b_n) \leq \bigvee_{d\in D} \tau(d)$,
so~$\tau(\,\bigvee_n a_n\,) = \tau(\,\bigvee_n b_n\,)$,
which implies that~$\bigvee_n a_n = \bigvee_n b_n$
as~$\tau$ is faithful.
Since then $b\leq b_1\leq \bigvee_n b_n = \bigvee_n a_n$
we see that~$\bigvee_n a_n$ is an upper bound
(and thus the supremum) of~$D$.
Moreover,
since
$\bigvee_{d\in D} \tau(d) \leq
\tau(\bigvee D)
=\tau(\bigvee_n a_n)
=\bigvee_n\tau(a_n)
\leq \bigvee_{d\in D} \tau(d)$,
we see that~$\bigvee_{d\in D} \tau(d)
= \tau(\bigvee D)$,
and so~$\tau$ is normal.
Since~$\tau$ is faithful
and normal, $\scrA$ is a von Neumann algebra.\qed
\end{point}
\end{point}
\begin{point}{90}[Linfty-vn]{Corollary}%
Given a finite complete measure space~$X$ 
the $C^*$-algebra $L^\infty(X)$
is a commutative von Neumann algebra,
and the assignment $f\mapsto \int f$
gives a faithful normal positive
map~$\int\colon L^\infty(X)\to\C$.
\end{point}
\end{parsec}
\begin{parsec}{520}[classification-cvn]%
\begin{point}{10}%
We'll show that any commutative von Neumann algebra~$\scrA$
that admits a faithful np-functional $\omega\colon \scrA\to \C$
is nmiu-isomorphic  to~$L^\infty(X)$
for some finite complete measure space~$X$.
It makes sense to regard this result
as a von Neumann algebra analogue
of Gelfand's theorem for commutative $C^*$-algebras,
(see~\sref{gelfand}---that any commutative $C^*$-algebra
is miu-isomorphic to~$C(Y)$ for some compact Hausdorff space~$Y$.)
But one should not take the
comparison too far too lightly:
while Gelfand's theorem readily
yields a clean equivalence 
between commutative $C^*$-algebras
and compact Hausdorff spaces
(see~\sref{gelfand-equivalence}),
the fact that  $L^\infty(X_1)\cong L^\infty(X_2)$
for finite complete measure spaces~$X_1$
and~$X_2$ does not even imply that~$X_1$ and~$X_2$ have
	the same cardinality.\footnote{Indeed, one may take
	$X_1$ to be a measure space consisting of a
	single non-negligible point~$*$ (so $X_1=\{*\}$
	and $\mu(X_1)\neq 0$),
	while letting $X_2$ be a measure space
	on an uncountable set
	formed by taking for the measurable subsets
	of~$X_2$ the countable subsets
	and their complements,
	by making the countable subsets negligible,
	and by giving all cocountable subsets
	the same non-zero measure.
	Then all measurable functions
	on~$X_1$ and on~$X_2$
	are constant almost everywhere,
	(because in $X_1$ and $X_2$ there are no 
	two non-negligible disjoint measurable subsets,)
	so that~$L^\infty(X_1)\cong \C\cong L^\infty(X_2)$.}
Obtaining an equivalence between commutative von Neumann algebras
and measure spaces is nonetheless
possible after a suitable non-trivial modification
to the category of measure spaces
(as is shown by Robert Furber 
in as of yet unpublished work.)

We obtain our finite complete measure space~$X$
from the commutative von Neumann algebra~$\scrA$
by taking for~$X$ the compact Hausdorff space~$\spec(\scrA)$
of all miu-functionals on~$\scrA$,
and declaring that a subset~$A$ of~$X\equiv \spec(\scrA)$ 
is measurable when~$A$ is clopen up to
a \emph{meagre}
subset
(defined below,~\sref{meagre}).
It takes some effort to show
that this yields a $\sigma$-algebra
in~$\spec(\scrA)$,
and that the faithful np-functional $\omega\colon\scrA\to \C$
gives a finite complete measure on~$\spec(\scrA)$,
but once this is achieved
it's easily seen that
$\scrA\cong C(\spec(\scrA))\cong L^\infty(\spec(\scrA))$.
\end{point}
\begin{point}{20}[meagre]{Definition}%
Let~$X$ be a topological space.
\begin{enumerate}
\item
A subset~$A$ of~$X$
is called \Define{meagre}%
\index{meagre}
when~$A\subseteq \bigcup_n B_n$
for some closed subsets $B_1\subseteq B_2\subseteq \dotsb$
of~$X$ with empty interior (so  $B_n^\circ=\varnothing$
for all~$n$.)
\item
Given $A,B\subseteq X$ 
we write $\Define{A\approx B}$
\index{((approx@$\approx$!$A\approx B$, for subsets of a topological space} when
$A\cup B\,\backslash\,A\cap B$
is meagre.
\item
We say that~$A\subseteq B$
is \Define{almost clopen}%
\index{almost clopen}
when~$A\approx C$
for some clopen~$C\subseteq X$.
\end{enumerate}
\spacingfix%
\end{point}%
\begin{point}{30}[meagre-basic]{Exercise}%
Given a topological space~$X$,
verify the following facts.
\begin{enumerate}
\item
A countable union~$\bigcup_n A_n$
of meagre subsets $A_1,A_2,\dotsc\subseteq X $ is meagre.
\item
A subset of a meagre set is meagre.
\item
$\overline{U}\approx U$
for every open subset~$U$ of~$X$.

(Hint: show that $\overline{U}\backslash U$
is closed with empty interior.)
\item
$\bigcup_n A_n\approx \bigcup_n B_n$
for all~$A_1,A_2,\dotsc, B_1,B_2,\dotsc\subseteq X$
with $A_n\approx B_n$.
\item
$A\backslash B \approx A'\backslash B'$
for all $A,A',B,B'\subseteq X$
with~$A\approx A'$ and~$B\approx B'$.
\item
If~$A,B\subseteq X$ are almost clopen,
then~$A\cup B$ and~$A\backslash B$ 
are almost clopen.
\end{enumerate}
\spacingfix
\end{point}%
\begin{point}{31}%
Meagerness can be thought of as a topological analogue of negligibility.
In fact, with respect to the measure we'll put on $\spec(\scrA)$
in~\sref{cvn-faithful},
    meagerness and negligibility actually coincide.
In general, however, 
the notions are disparate,
as is demonstrated rather dramatically by the following example.
\end{point}
\begin{point}{32}{Example}
of a meagre subset~$A$ of~$[0,1]$
with Lebesgue measure~$1$.
Given an enumeration~$q_1,q_2,\dotsc$
    of the rational numbers in~$[0,1]$,
the set
\begin{equation*}
    \textstyle
    B_m \ :=\  \bigcup_n 
    \bigl(\,q_n-\frac{1}{2m}2^{-n},\,q_n+\frac{1}{2m}2^{-n}\,\bigr)
    \cap[0,1]
\end{equation*}
is open and dense in~$[0,1]$, with a Lebesgue measure 
of at most~$\nicefrac{1}{m}$.
So the intersection~$B:=\bigcap_m B_m$ is Lebesgue negligible.
On the other hand, $[0,1]\backslash B_m$
is closed and has empty interior,
so that $A:=[0,1]\backslash B\equiv \bigcup_m B_m$
is meagre with Lebesgue measure~$1$.
\end{point}
\end{parsec}
\begin{parsec}{530}%
\begin{point}{10}%
The fact that the almost clopen subsets
of the spectrum~$\spec(\scrA)$
of a commutative von Neumann algebra~$\scrA$
are closed under countable unions
(and thus form a $\sigma$-algebra)
relies on a special topological property of~$\spec(\scrA)$
that is described in~\sref{vn-spectrum-extremally-disconnected} below.
\end{point}
\begin{point}{20}[ngelfand]{Exercise}%
\index{Gelfand's Representation Theorem!for von Neumann algebras}
Let~$\scrA$ be a commutative von Neumann algebra.
Using the fact that
the Gelfand representation~$\gamma_\scrA\colon \scrA\to C(\spec(\scrA))$
from~\sref{gelfand-representation}
is an miu-isomorphism 
by~\sref{gelfand}
and thus an order isomorphism,
show that~$C(\spec(\scrA))$
is a commutative von Neumann algebra
that is nmiu-isomorphic
to~$\scrA$ via~$\gamma_\scrA$.
\end{point}
\begin{point}{30}[vn-spectrum-extremally-disconnected]{Proposition}%
	\index{sp@$\spec$, spectrum!$\spec(\scrA)$, of a $C^*$-algebra!is extremally disconnected for a von Neumann algebra}
The spectrum~$\spec(\scrA)$
of a commutative von Neumann algebra~$\scrA$
is \Define{extremally disconnected}:%
\index{extremal disconnectedness} 
the closure~$\overline{U}$
of an open subset~$U$ of~$\spec(\scrA)$
is open.
\begin{point}{40}{Proof}%
(Based on \S6.1 of~\cite{riesz}.)

Let~$U$ be an open subset of~$\spec(\scrA)$,
and let~$\mathbf{1}_U$ be the indicator function
of~$U$.
The set $D=\{\,f\in C(\spec(\scrA))\colon\, f\leq \mathbf{1}_U\,\}$
is directed and bounded 
and so has a supremum~$\bigvee D$
in~$C(\spec(\scrA))$
since~$C(\spec(\scrA))$ is a von Neumann algebra by~\sref{ngelfand}.
Note that~$0\leq \bigvee D\leq 1$.
We'll prove that~$\bigvee D = \mathbf{1}_{\overline{U}}$,
because this entails
that~$\mathbf{1}_{\overline{U}}$
is continuous,  so that~$\overline{U}$ is both open and closed.

Let~$x\in U$ be given.
By Urysohn's lemma (see~15.6 of~\cite{willard},
using here that~$\spec(\scrA)$ being a compact Hausdorff
space, \sref{spectrum-calg-compact-hausdorff}, 
is normal by~17.10 of~\cite{willard})
there is $f\in [0,1]_{C(\spec(\scrA))}$
with~$f(x)=1$ and~$f(y)=0$ for all~$y\in \spec(X)\backslash U$.
It follows that~$f\in D$,
and~$f\leq \bigvee D\leq 1$,
so that $1=f(x)\leq (\bigvee D)(x)\leq 1$,
and~$(\bigvee D)(x)=1$.
By continuity of~$\bigvee D$,
we get $(\bigvee D)(x)=1$ for all~$x\in \overline{U}$.

Let~$y\in \spec(\scrA)\backslash U$ be given.
Again by Urysohn's lemma
there is $f\in [0,1]_{C(\spec(\scrA))}$
with~$f(y)=0$ and~$f(x)=1$ for all~$x\in \overline{U}$.
Since $g\leq \mathbf{1}_U \leq f$
for every~$g\in D$,
we get $\bigvee D\leq f$,
and so~$0\leq (\bigvee D)(y)\leq f(y) =0$,
which implies that~$(\bigvee D)(y)=0$.
Hence~$(\bigvee D)(y)=0$ for all~$y\in \spec(\scrA)\backslash U$.

All in all we have $\bigvee D = \mathbf{1}_{\overline{U}}$,
and so~$\overline{U}$ is open.\qed
\end{point}
\end{point}
\begin{point}{50}{Corollary}%
The almost clopen subsets
of an extremally disconnected topological
space~$X$ form a $\sigma$-algebra.
\begin{point}{60}[almost-clopen-sigma-algebra-proof]{Proof}%
In light of~\sref{meagre-basic}
it remains only to be shown that
the union~$\bigcup_n A_n$
of almost clopen subsets~$A_1,A_2,\dotsc$
is almost clopen.
Let~$C_1,C_2,\dotsc\subseteq X$
be clopen with $A_n\approx C_n$ for each~$n$.
Then $\bigcup_n A_n \approx \bigcup_n C_n$,
and~$C:=\bigcup_n C_n$ is open
(but not necessarily closed).
Since~$C\approx \overline{C}$
(by~\sref{meagre-basic}),
and~$\overline{C}$ is clopen (as $X$ is extremally disconnected)
we get~$\bigcup_n A_n \approx \overline{C}$,
so~$\bigcup_n A_n$ is almost clopen.\qed
\end{point}
\end{point}
\end{parsec}
\begin{parsec}{540}%
\begin{point}{10}%
The final ingredient
we need to prove the main result, \sref{cvn-faithful},
of this section
is the observation
that an almost clopen subset of a compact Hausdorff space
is equivalent to precisely one clopen,
which follows from the following famous theorem.
\end{point}
\begin{point}{20}[baire-category-theorem]{Baire category theorem}%
\index{Baire's Category Theorem}
A meagre subset of a compact Hausdorff space
has empty interior.
\begin{point}{30}{Proof}%
Let~$A$ be a meagre subset of a compact Hausdorff space~$X$.
So there are closed $B_1\subseteq B_2\subseteq \dotsc$
with~$A\subseteq \bigcup_n B_n$
and~$B_n^\circ=\varnothing$ for all~$n$.
Then~$U_n:= X\backslash B_n$
is an open dense subset of~$X$
for each~$n$.
Since~$A^\circ \subseteq (\bigcup_n B_n)^\circ  =
X\backslash (\,\smash{\overline{\bigcap_n U_n}}\,)$
it suffices to show that
$\bigcap_n U_n$ is dense in~$X$.
That is, given a non-empty open subset~$V$
of~$X$ we must show that~$V\cap\bigcap_n U_n\neq \varnothing$.

Write~$V_1:=V$.
Since~$U_1$ is open and dense, and~$V_1$ is open and not empty,
we have $U_1\cap V_1\neq \varnothing$.
	Since~$X$ is regular (see e.g.~\cite{willard})
we can find an open and non-empty subset $V_2$ of~$X$
with~$\overline{V}_2\subseteq U_1\cap V_1$.
Continuing this process
we obtain non-empty open subsets $V\equiv V_1 \supseteq V_2\supseteq \dotsb$
	of~$X$ with $\overline{V}_{n+1} \subseteq U_n\cap V_n$
for all~$n$,
and so~$
\overline{V}_1 \supseteq V_1 
\supseteq \overline{V}_2 \supseteq V_2 
\supseteq \dotsb$.
Since~$X$ is compact,
$\bigcap_n \overline{V}_n$
can not be empty,
and neither will be~$V\cap \bigcap_n U_n
\supseteq \bigcap_n \overline{V}_n$.\qed
\end{point}
\end{point}
\begin{point}{40}[approx-closure]{Lemma}%
For open subsets~$U$ and~$V$ of a compact Hausdorff
space~$X$,
\begin{equation*}
U\approx V\quad\iff\quad
\overline{U}\approx \overline{V}
\quad\iff\quad
\overline{U}=\overline{V}.
\end{equation*}
\spacingfix%
\begin{point}{50}{Proof}%
As~$U\approx \overline{U}$ by~\sref{meagre-basic}
the only thing that is not obvious
is that~$\overline{U}\approx \overline{V} \implies \overline{U}=\overline{V}$.
So suppose that~$\overline{U}\approx \overline{V}$.
Then~$U\backslash \overline{V}$
is empty, because it is an open subset of the
meagre set~$\overline{U}\cup\overline{V}\,\backslash\,
\overline{U}\cap\overline{V}$
(which has empty interior by~\sref{baire-category-theorem}.) In other words,
we have~$U\subseteq \overline{V}$, 
and thus~$\overline{U}\subseteq \overline{V}$.
Similarly, $\overline{V}\subseteq\overline{U}$,
and so~$\overline{V}=\overline{U}$.\qed
\end{point}
\end{point}
\begin{point}{60}{Corollary}%
Given an almost clopen subset~$A$ of a compact Hausdorff
space~$X$ there is precisely one clopen~$C$ with~$A\approx C$.
\begin{point}{70}[almost-meagre-fundamental]{Proof}%
When~$C\approx A \approx C'$
for clopen subsets~$C,C'\subseteq X$,
we have~$C\approx C'$,
and so~$C=C'$ by~\sref{approx-closure}.\qed
\end{point}
\end{point}
\begin{point}{80}%
Interestingly,
a compact Hausdorff space is extremally disconnected
iff each of its open subsets is ``measurable''
in the sense of being almost clopen:
\end{point}
\begin{point}{90}[open-almost-clopen]{Proposition}%
A compact Hausdorff space~$X$
is extremally disconnected
iff every open subset of~$X$
is almost clopen.
\begin{point}{100}{Proof}%
If~$X$ is extremally disconnected,
and~$U$ is open subset of~$X$,
then~$\overline{U}$ is clopen,
and~$\overline{U}\approx U$ by~\sref{meagre-basic}
giving us that~$U$ is almost clopen.

Conversely,
suppose that each open subset of~$X$ is almost clopen.
To show that~$X$ is extremally disconnected
we must show that~$\overline{U}$ is open given
an open subset~$U$ of~$X$.
Pick a clopen~$C$ with~$U\approx C$.
Then~$\overline{U}\approx U\approx C$
(by~\sref{meagre-basic}),
and so~$\overline{U}=C$ by~\sref{approx-closure}.\qed
\end{point}
\end{point}
\begin{point}{110}[cvn-faithful]{Theorem}%
Let~$\scrA$ be a commutative von Neumann algebra~$\scrA$.
Recall that
the Gelfand representation $\gamma_\scrA\colon 
\scrA\to C(\spec(\scrA))$
is an nmiu-isomorphism (by~\sref{ngelfand}),
$C(\spec(\scrA))$ is a von Neumann algebra,
and that the almost clopen subsets (see~\sref{meagre}) of~$\spec(\scrA)$
form a $\sigma$-algebra.

Given a faithful np-functional $\omega\colon \scrA\to\C$
there is a (unique) measure
$\mu$ on the almost clopen subsets of~$\spec(\scrA)$
such that~$\mu(A)=0$ iff $A$ is meagre,
and~$\mu(C)=\omega(\gamma_\scrA^{-1} (\mathbf{1}_{C}))$
for every clopen subset~$C$ of~$\spec(\scrA)$;
and this turns~$\spec(\scrA)$
into a finite complete measure space.

With respect to this measure space a
bounded function $f\colon \spec(\scrA)\to\C$
is measurable iff~$f$ is continuous almost everywhere.
Moreover, 
$f\mapsto f^\circ\colon C(\spec(\scrA))\to L^\infty(\spec(\scrA))$
is an nmiu-isomorphism,
and~$\int f^\circ = \omega(\gamma_\scrA^{-1}(f))$ for all~$f\in C(X)$.
All in all,
we get the following commuting diagram.
\begin{equation*}
\xymatrix@C=3em@R=3em{
\scrA
\ar[r]^-{\gamma_\scrA}_-\cong
\ar[rd]_\omega
&
C(\spec(\scrA))
\ar[r]^-{f\mapsto f^\circ}_-\cong
&
L^\infty(\spec(\scrA))
\ar[ld]^\int
\\
&
\C
}
\end{equation*}
\begin{point}{120}{Proof}%
By~\sref{almost-meagre-fundamental}
we know that given an almost clopen subset~$A$
of~$\spec(\scrA)$
there is a unique clopen~$C_A$
with~$A\approx C_A$,
and 
so we may define
$\mu(A):=\omega(\gamma_\scrA^{-1}(\mathbf{1}_{C_A}))$.
It is easily seen that~$\mu$
is finitely additive.
Further $\mu(A)=0$ for every meagre~$A\subseteq X$,
and so $\mu(A)=\mu(B)$ when~$A\approx B$.
Conversely, 
an almost clopen subset~$A$
of~$\scrA$
with~$\mu(A)=0$
is meagre,
because for the unique clopen~$C$ with~$A\approx C$,
we have~$\omega(\gamma_\scrA^{-1}(\mathbf{1}_C))=\mu(A)=0$,
so that~$\mathbf{1}_C=0$ and thus~$C=\varnothing$---using here
that~$\omega$ is faithful.

To show that~$\mu$
is a measure,
it suffices to prove
that~$\bigwedge_n \mu(A_n)=0$
given~$A_1\supseteq A_2\supseteq \dotsb$
with~$\bigcap_n A_n=\varnothing$.
To do this, pick clopen subsets~$C_1,C_2,\dotsc$
of~$\spec(\scrA)$
with~$A_n\approx C_n$ for all~$n$.
Then~$\bigwedge_n \mu(A_n)
=\bigwedge_n \mu(C_n)
= \omega(\gamma_\scrA^{-1}(\bigwedge_n \mathbf{1}_{C_n}))$---using 
here that~$\omega$ is normal.
So to prove that~$\bigwedge_n \mu(A_n)=0$
it suffices to show that~$\bigwedge_n \mathbf{1}_{C_n}=0$,
that is,
given a lower bound~$f$ of the~$\mathbf{1}_{C_n}$
in~$\Real{C(\spec(\scrA))}$
we must show that~$f\leq 0$.
Note that for such~$f$
we have~$f(x)\leq 0$
for all~$x\in X\backslash\bigcap_n C_n$.
Then~$f(x)\leq 0$ for all~$x\in X$
if we can show that~$X\backslash \bigcap_n C_n$ is dense in~$X$.
But this indeed the case
since~$\bigcap_n C_n\approx \bigcap_n A_n=\varnothing$
is meagre, and therefore has empty interior
(by~\sref{baire-category-theorem}).
Whence~$\mu$ is a measure.
Note that~$\mu$ is finite,
because~$\mu(\spec(\scrA))
= \omega(1)<\infty$,
and complete,
because a subset of a meagre set is meagre.

Let~$h\colon \spec(\scrA)\to\C$ be a bounded function.
We'll show that~$h$ is continuous almost everywhere
iff~$h$ is measurable.
Surely,
if~$h$ is continuous (everywhere),
then~$h$ is measurable
(since every open subset~$U$ of~$\spec(\scrA)$
is almost clopen, \sref{open-almost-clopen}).
So if~$h$ is continuous almost everywhere,
then~$h$ is measurable too.
For the converse,
it suffices to show
that $\varrho \colon h\mapsto h^\circ\colon C(\spec(\scrA))\to
L^\infty(\spec(\scrA))$ is surjective.
To this end, note first that~$\varrho$ is injective,
because a continuous function on~$\spec(\scrA)$
that is zero almost everywhere,
is non-zero on a meagre set,
and by~\sref{baire-category-theorem} zero on a dense
subset, and so is zero everywhere.
Since the image of the injective miu-map $\varrho$ is norm closed
in order to show that~$\varrho$ is surjective
it suffices to show that image of~$\varrho$ is 
norm dense in~$L^\infty(X)$.
This is indeed the case
since the elements of~$L^\infty(\spec(\scrA))$
of the form~$\sum_{n} \lambda_n \mathbf{1}_{A_n}^\circ$
where~$\lambda_1,\dotsc,\lambda_N\in\C$
and~$A_1,\dotsc,A_N$ are measurable (i.e.~almost clopen)
subsets of~$\spec(\scrA)$
are easily seen to be norm dense in~$L^\infty(\spec(\scrA))$
(c.f.~243I of~\cite{fremlin}),
and are in the range of~$\varrho$,
because given an almost clopen~$A\subseteq \spec(\scrA)$
and a clopen~$C$ with $A\approx C$
we have~$\mathbf{1}_A^\circ = \mathbf{1}_C^\circ$
and~$\mathbf{1}_{C}\in C(\spec(\scrA))$.
Hence~$\varrho$ is surjective.

It remains to be show that
$\int f^\circ = \omega(\gamma_\scrA^{-1}(f))$
for all~$f\in C(\spec(\scrA))$,
that is, $\int = \omega\circ \gamma_\scrA^{-1}\circ \varrho^{-1}$.
By the previous discussion the linear span
of the elements of~$L^\infty(\spec(\scrA))$
of the form~$\mathbf{1}_C^\circ$,
where~$C$ is (not just measurable but) clopen,
is norm dense in $L^\infty(\spec(\scrA))$.
Since~$\int \mathbf{1}_C
= \mu(C)=\omega(\gamma_\scrA^{-1} (\varrho^{-1}(\mathbf{1}_C^\circ ))$
for all clopen~$C$,
and both~$\int$ and~$\omega\circ\gamma_\scrA^{-1}\circ \varrho^{-1}$
are linear and bounded,
we conclude that~$\int = \omega\circ \gamma_\scrA^{-1} \circ \varrho^{-1}$,
and so we are done.\qed
\end{point}
\begin{point}{130}%
To deduce from this that all commutative von Neumann algebras
(and not just the ones admitting a faithful np-functional)
are nmiu-isomorphic
to direct sums of the form  $\bigoplus_i L^\infty(X_i)$
where the~$X_i$ are finite complete measure spaces
we first need some basic facts concerning the
\emph{projections} of a commutative von Neumann algebra.
\end{point}
\end{point}
\end{parsec}%
\section{Projections}
\begin{parsec}{550}%
\begin{point}{10}%
One pertinent feature
of von Neumann algebras
is an abundance of projections:
above each effect~$a$ 
there is a least projection~$\ceil{a}$
we call the ceiling of~$a$ (\sref{vna-ceil});
for every np-map $\omega\colon \scrA\to \scrB$
between von Neumann algebras
there is a least projection~$p$ with~$\omega(p^\perp)=0$
called the carrier of~$\omega$ (see~\sref{carrier});
the directed supremum of projections is again a projection;
the partial order of projections is complete
(see~\sref{ceil-floor-basic});
and each element of a von Neumann algebra is the norm limit
of linear combinations of projections
(see~\sref{projections-norm-dense}).
We'll prove all this and more in this section.
\end{point}
\begin{point}{20}{Definition}%
An element~$p$ of a $C^*$-algebra
is a \Define{projection}%
\index{projection!in a $C^*$-algebra}
when~$p^*p=p$.
\end{point}
\begin{point}{30}{Examples}%
\begin{enumerate}
\item
The only projections in~$\C$ are~$0$ and~$1$.
\item
Given a measurable
subset~$A$ of a finite complete measure space~$X$
the indicator function~$\mathbf{1}_A$
is a projection in~$L^\infty(X)$,
and every projection in~$L^\infty(X)$
is of this form.
\item
Given a closed linear subspace~$C$ of a Hilbert space~$\scrH$
the inclusion $E\colon C\to \scrH$
is a bounded linear map,
and  $\Define{P_C}:=EE^*\colon \scrH\to\scrH$
is a projection in~$\scrB(\scrH)$,
and
every projection in~$\scrB(\scrH)$ is of this form.
\end{enumerate}
\spacingfix%
\end{point}
\begin{point}{40}[projection-basic]{Exercise}%
Show that in a $C^*$-algebra:
\begin{enumerate}
\item
$0$ and~$1$ are projections.
\item
A projection~$p$ is an effect,
that is, $p^*=p$
and $0\leq p\leq 1$.
\item
    The orthocomplement $\Define{p^\perp} \equiv 1-p$ of a projection~$p$
    \index{$(\,\cdot\,)^\perp$!$p^\perp$, orthocomplement of a projection}
is a projection.
\item
An effect~$a$ is a projection iff $aa^\perp=0$.
\end{enumerate}
\spacingfix%
\end{point}%
\begin{point}{50}[ad-contraposed]{Lemma}%
Let~$a$ be an element of a $C^*$-algebra~$\scrA$
with $\|a\|\leq 1$,
and let~$p$ and~$q$ be projections on~$\scrA$.
Then 
$a^* p a \leq q^\perp$
iff $paq=0$
iff  $aqa^*\leq p^\perp$.
\begin{point}{60}{Proof}%
Suppose that~$a^*pa\leq q^\perp$.
Then we have $q a^*pa q \leq qq^\perp q = 0$
(see \sref{astara-pos-basic-consequences})
and so $paq=0$,
because $\|paq\|^2=\|(paq)^*paq\|=0$
by the $C^*$-identity.
Applying $(\,\cdot\,)^*$,
we get $qa^*p=0$, and so both $qa^* = qa^*p^\perp$
and $aq = p^\perp aq$, giving
us $aqa^* = p^\perp a q a^* p^\perp 
\leq p^\perp$,
where we used that $aqa^*\leq aa^*\leq \|aa^*\|=\|a\|^2\leq 1$.
By a similar reasoning,
we get $aqa^*\leq p^\perp \implies paq=0\implies a^*pa\leq q^\perp$.\qed
\end{point}
\end{point}
\begin{point}{70}{Exercise}%
Let~$a$ be an effect of a $C^*$-algebra~$\scrA$,
and~$p$ be a projection from~$\scrA$.
\begin{point}{80}[projection-above-effect]%
Show that $a\leq p$
iff $p\sqrt{a} = \sqrt{a}$
iff $\sqrt{a}p = \sqrt{a}$
iff $p^\perp\sqrt{a} = 0$
iff $\sqrt{a}p^\perp = 0$
iff $a^2\leq p$
iff $p a  = a$
iff $ a p = a $
iff $p^\perp a  = 0$
iff $ap^\perp = 0$
iff $\sqrt{a}\leq p$.
\end{point}
\begin{point}{90}[projection-below-effect]%
Show that $p\leq a$
iff $p \sqrt{a} = p$
iff $\sqrt{a} p = p$
iff $ p\sqrt{a}^\perp = 0$
iff $\sqrt{a}^\perp p = 0$
iff $p\leq a^2$
iff $ap=p$
iff $pa = p$
iff $pa^\perp =0$
iff $a^\perp p =0$
iff $p\leq \sqrt{a}$.
\end{point}
\end{point}
\begin{point}{100}[projection-order-sharp]{Lemma}%
An effect~$a$ of a $C^*$-algebra~$\scrA$
is a projection iff the only effect
below~$a$ and~$a^\perp$ is~$0$.
\begin{point}{110}{Proof}%
On the one hand,
if~$a$ is a projection,
and~$b$ is an effect with~$b\leq a$
and~$b\leq a^\perp$,
then~$a^\perp b=0$ and~$ab=0$ by~\sref{projection-above-effect},
and so~$b=ab+a^\perp b = 0$.
On the other hand,
if~$0$ is the only effect below both~$a$ and~$a^\perp$,
then~$aa^\perp\equiv \sqrt{a}a^\perp \sqrt{a}$
being an effect below~$a$, and below~$a^\perp$,
is zero, and so~$a$ is projection, by~\sref{projection-basic}.\qed
\end{point}
\end{point}
\begin{point}{120}{Definition}%
We say that projections~$p$ and~$q$
from a $C^*$-algebra~$\scrA$ are \Define{orthogonal}%
\index{orthogonal projections}
when~$pq=0$,
and we say that a subset of projections
from~$\scrA$ is \Define{orthogonal}
(and its elements are \Define{pairwise orthogonal})
when all~$p$ and~$q$ from~$E$
are either equal or orthogonal.
\end{point}
\begin{point}{130}[orthogonal-tuple-of-projections]{Exercise}%
Let~$\scrA$ be a $C^*$-algebra.
\begin{enumerate}
\item
Show that projections~$p$ and~$q$ from~$\scrA$
are orthogonal iff $pq=0$ iff $qp=0$ iff $pqp=0$
iff $p+q\leq 1$ iff $p\leq q^\perp$
iff $p+q$ is a projection.
\item
Show that a finite set of  projections $p_1,\dotsc,p_n$
from~$\scrA$ is orthogonal
iff~$\sum_i p_i \leq 1$
iff $\sum_i p_i$ is a projection.

Show that, in that case, $\sum_i p_i$ is the least projection
above~$p_1,\dotsc,p_n$.
\end{enumerate}
\spacingfix%
\end{point}%
\begin{point}{140}[projection-below-projection]{Exercise}%
Let~$p$ and~$q$ be projections from a $C^*$-algebra
with~$p\leq q$.\\
Show that~$q-p$ is a projection
(either directly, or using~\sref{orthogonal-tuple-of-projections}).
\end{point}
\end{parsec}
\subsection{Ceiling and Floor}
\begin{parsec}{560}%
\begin{point}{10}[vna-ceil]{Proposition}%
Above every effect~$b$ of a von Neumann algebra~$\scrA$,
there is a smallest projection, \Define{$\ceil{b}$},%
\index{*ceil@$\ceil{\,\cdot\,}$!$\ceil{a}$, ceiling}
which we call the \Define{ceiling}%
\index{ceiling}
of~$b$,
 given by $\ceil{b}=\bigvee_{n=0}^\infty b^{\nicefrac{1}{2^n}}$.\\
Moreover, if $a\in \scrA$ commutes with $b$,
then~$a$ commutes with~$\ceil{b}$.
\begin{point}{20}{Proof}
Since~$0\leq b\leq b^{\nicefrac{1}{2}}\leq
b^{\nicefrac{1}{4}}\leq\dotsb\leq 1$,
    we may define~$p:=\bigvee_n b^{\nicefrac{1}{2^n}}$.
\begin{point}{30}[vna-ceil-point-1]%
To begin,
note that if~$a\in \scrA$
commutes with~$b$,
then~$a$ commutes with~$p$.
Indeed, for such~$a$ we have~$a\sqrt{b}=\sqrt{b}a$
by~\sref{sqrt},
and so $a b^{\nicefrac{1}{2^n}} = b^{\nicefrac{1}{2^n}} a$
for each~$n$
by induction.
Thus~$ap=pa$ by~\sref{vna-supremum-commutes}.
\end{point}
\begin{point}{40}%
Let us prove that~$p$ is a projection, i.e.~$p^2=p$. 
Since~$p\leq 1$, we already have $p^2\equiv \sqrt{p}p\sqrt{p}\leq p$
by~\sref{astara-pos-basic-consequences},
and so we only need to show that $p\leq p^2$. We have:
\begin{alignat*}{3}
 p^2 \ &=\  \textstyle \bigvee_m \sqrt{p} \,b^{\nicefrac{1}{2^m}} \,\sqrt{p}
\qquad&&\text{by \sref{ad-normal}} \\
&=\ \textstyle\bigvee_m b^{\nicefrac{1}{2^{m+1}}}\, p\,
b^{\nicefrac{1}{2^{m+1}}} 
\qquad&&\text{by \sref{vna-ceil-point-1} and \sref{sqrt}} \\
&=\ \textstyle \bigvee_m \bigvee_n \, 
b^{\nicefrac{1}{2^{m+1}}}\, b^{\nicefrac{1}{2^n}}\,
b^{\nicefrac{1}{2^{m+1}}} \qquad && \text{by \sref{ad-normal}}
\end{alignat*}
Thus $p^2 \geq b^{\nicefrac{1}{2^k}}$
for each~$k$ (taking $n=m=k+1$,)
and so~$p^2 \geq p$.
\end{point}
\begin{point}{50}%
It remains to be shown that~$p$ is the \emph{least} projection
above~$b$.
Let~$q$ be a projection in~$\scrA$ with $b\leq q$;
we must show that~$q\leq p$.
We have $b^{\nicefrac{1}{2}}\leq q$
by~\sref{projection-above-effect},
and so $b^{\nicefrac{1}{2^n}}\leq q$ for each~$n$ by induction.
Hence $p\leq q$.\qed
\end{point}
\end{point}
\end{point}
\begin{point}{60}[vna-floor]{Proposition}%
Below every effect~$b$ of a von Neumann algebra~$\scrA$,
there is greatest projection, \Define{$\floor{b}$},%
\index{*floor@$\floor{\,\cdot\,}$!$\floor{a}$, of an effect}
	we call the \Define{floor}%
\index{floor!of an effect}
of~$b$,
given by~$\floor{b} = \bigwedge_{n=0}^\infty b^{2^{n}}$.\\
Moreover, if~$a\in \scrA$ commutes with~$b$,
then~$b$ commutes with~$\floor{b}$.
\begin{point}{70}{Proof}%
Note that $1\geq b\geq b^2 \geq b^4 \geq  \dotsb \geq 0$,
    and define $p:=\bigwedge_n b^{2^n}$
    (see~\sref{infima-in-vna}.)
\begin{point}{80}[vna-floor-point-1]%
If~$a\in \scrA$ commutes with~$b$,
then~$a$ commutes with~$p$.
Indeed, such~$a$ commutes with~$b^2$ (because
$ab^2 = bab = b^2a$,)
and so~$a$ commutes with~$b^{2^n}$ for each~$n$ by induction.
Thus~$a$ commutes with~$p\equiv\bigwedge_n b^{2^n}$ 
(by a variation on~\sref{vna-supremum-commutes}.)
\end{point}
\begin{point}{90}%
To see that~$p$ is a projection, i.e.~$p^2=p$,
we only need to show that~$p\leq p^2$,
because we get $p^2\equiv \sqrt{p}\,p\,\sqrt{p}\leq p$
from $p\leq 1$ (using~\sref{astara-pos-basic-consequences}.)
Now, since
\begin{alignat*}{3}
p^2 \ &=\ \textstyle \bigwedge_m\  \sqrt{p}\, b^{2^m} \sqrt{p}\qquad
&&\text{by a variation on~\sref{ad-normal}}\\
&=\ \textstyle \bigwedge_m \ b^{2^{m-1}} p\, b^{2^{m-1}}\qquad
&&\text{by~\sref{vna-floor-point-1} and~\sref{sqrt}}\\
&=\ \textstyle \bigwedge_m \bigwedge_n \ 
b^{2^{m-1}}\, b^{2^n}\, b^{2^{m-1}}\qquad
&&\text{by~\sref{ad-normal},}
\end{alignat*}
and $p\leq b^{2^{m-1}}\, b^{2^n}\,b^{2^{m-1}}$
for all~$n,m$, we get~$p\leq p^2$.
\end{point}
\begin{point}{100}%
It remains to be shown that~$p$ is the greatest projection below~$b$.
Let~$q$ be a projection in~$\scrA$ with~$q\leq b$.
We must show that~$q\leq p$.
Since~$q\leq b$,
we have~$q\leq b^2$ (by~\sref{projection-below-effect}),
and so~$q\leq b^{2^n}$ for each~$n$ by induction.
Thus~$q\leq p\equiv\bigwedge_n b^{2^n}$.\qed
\end{point}
\end{point}
\end{point}
\begin{point}{110}[ceil-floor-second-property]{Exercise}%
Show that given an effect~$a$ and a projection~$p$
in a von Neumann algebra~$\scrA$ we have
\begin{enumerate}
\item
$pa=a$ iff $ap=a$ iff $\ceil{a}\leq p$, and
\item
$pa=p$ iff $ap=p$ iff $p\leq \floor{a}$.
\end{enumerate}
Conclude that~$\ceil{a}$
is the least projection~$p$ with $a=a p$
(or, equivalently, $a=pa$),
and that $\floor{a}$
is the greatest projection~$p$ with $p=a p$
(or, equivalently, $p=pa$.)

In particular,
$a=a\ceil{a}=\ceil{a}a$
and $\floor{a}=a\floor{a}=\floor{a}a$.
\end{point}
\begin{point}{120}{Example}%
Given a finite complete measure space~$X$
we have
\begin{equation*}
\ceil{f^\circ}
	\,=\,\mathbf{1}_{\{x\in X\colon f(x)>0\}}^\circ
	\qquad
	\text{and}\qquad
	\floor{f^\circ}
	\,=\, \mathbf{1}_{\{x\in X\colon  f(x)=1\}}^\circ
\end{equation*}
for every~$f\in\mathcal{L}^\infty(X)$
with~$0\leq f^\circ \leq 1$.
\end{point}
\begin{point}{130}[ceil-floor-basic]{Exercise}%
Let~$a,b$ be effects of a von Neumann algebra~$\scrA$,
and let~$\lambda\in [0,1]$.
\begin{enumerate}
\item
Show that $\ceil{a}^\perp = \floor{a^\perp}$
and $\floor{a}^\perp = \ceil{a^\perp}$.
\item
Show that~$\ceil{\lambda a} = \ceil{a}$
when~$\lambda\neq 0$.

Use this to prove that~$\ceil{\lambda a+\lambda^\perp b}$
is the supremum of~$\ceil{a}$ and~$\ceil{b}$
in the poset of projections of~$\scrA$
when~$\lambda\neq 0$ and~$\lambda\neq 1$.
\item
Show that $\floor{a}=\floor{a^2}$
and $\ceil{a}=\ceil{a^2}$.
\end{enumerate}%
\spacingfix%
\end{point}%
\begin{point}{140}[vna-directed-supremum-projections]{Lemma}%
The supremum of a directed set~$D$ of projections
from a von Neumann algebra~$\scrA$ is a projection.
\begin{point}{150}{Proof}%
Writing $p=\bigvee D$,
we must show that $p^2=p$.
Note that $dp=d$ for all~$d\in D$
(by~\sref{projection-below-effect} because~$d\leq p$.)
Now, on the one hand, $(d)_{d\in D}$
converges ultraweakly to~$p$.
On the other hand,
$(dp)_{d\in D}$
converges ultraweakly to~$p^2$ by~\sref{vna-supremum-mult}.
Hence~$p=p^2$ by uniqueness of ultraweak limits.
\end{point}
\end{point}
\begin{point}{160}{Exercise}%
Deduce from this result
 that every set~$A$ of projections from~$\scrA$
has a supremum $\Define{\bigcup A}$%
\index{*bigcupA@$\bigcup A$, supremum of projections}
and an infimum $\Define{\bigcap A}$%
\index{*bigcapA@$\bigcap A$, infimum of projections}
\emph{in the poset of projections from~$\scrA$}.\\
(Hint: use~\sref{ceil-floor-basic},
and the fact that $p\mapsto p^\perp$ 
is an order anti-isomorphism on the poset of projections on~$\scrA$.)
\end{point}
\begin{point}{170}[ceil-supremum]{Exercise}
Let~$\scrA$ be a von Neumann algebra.
\begin{enumerate}
\item
Show that $\ceil{\bigvee D}=\bigcup_{d\in D} \ceil{d}$
for every directed set~$D$ of effects from~$\scrA$.
\item
Show that $\floor{\bigwedge D} = \bigcap_{d\in D} \floor{d}$
for every filtered set~$D$ of effects from~$\scrA$.
\item
Show that $\ceil{\,\cdot\,}$
does not preserve filtered infima,
and~$\floor{\,\cdot\,}$
does not preserve directed suprema.
(Hint: $1,\frac{1}{2},\frac{1}{3},\dotsc$.)

Conclude that $\ceil{\,\cdot\,}$
and $\floor{\,\cdot\,}$
are neither ultraweakly, ultrastrongly nor norm  continuous
as maps from $[0,1]_\scrA$ to~$[0,1]_\scrA$.
\end{enumerate}
\spacingfix
\end{point}%
\begin{point}{180}[sum-of-orthogonal-projections]{Exercise}%
Show that for a family~$(p_i)_{i\in I}$ 
of pairwise orthogonal projections
(with~$I$ potentially uncountable)
the series $\sum_i p_i$
converges ultrastrongly to~$\bigcup_i p_i$.

(Hint: 
use the fact that $\sum_{i\in F} p_i = \bigcup_{i \in F} p_i$
for finite subsets~$F$ of~$I$ by~\sref{orthogonal-tuple-of-projections}.)
\end{point}
\end{parsec}
\begin{parsec}{570}[floor-sequential-product]%
\begin{point}{10}{Lemma}%
Let~$a,b$ be effects of a von Neumann algebra~$\scrA$.
Then~$\floor{\sqrt{a}b\sqrt{a}}$ is the greatest projection
below~$a$ and~$b$, that is, 
$\floor{\sqrt{a}b\sqrt{a}}=\floor{a}\cap \floor{b}$.
\begin{point}{20}{Proof}%
Surely, $\floor{\sqrt{a}b\sqrt{a}}\leq \sqrt{a}b\sqrt{a} \leq a$.
Let us prove that~$\floor{\sqrt{a}b\sqrt{a}}\leq b$.
To this end,
recall
that (by~\sref{projection-below-effect})
a projection~$e$ is below an effect~$c$
iff $ec=e$ iff $e\sqrt{c}=e$.
In particular,
since~$\floor{\sqrt{a}b\sqrt{a}}\leq \sqrt{a}b\sqrt{a}$ and 
$\floor{\sqrt{a}b\sqrt{a}}\leq a$,
we get
\begin{equation*}
\floor{\sqrt{a}b\sqrt{a}}
\ =\ \floor{\sqrt{a}b\sqrt{a}}\sqrt{a}b\sqrt{a}\floor{\sqrt{a}b\sqrt{a}} \ =\ 
\floor{\sqrt{a}b\sqrt{a}}b\floor{\sqrt{a}b\sqrt{a}},
\end{equation*}
and so $\floor{\sqrt{a}b\sqrt{a}}b^\perp\floor{\sqrt{a}b\sqrt{a}}=0$,
which implies that
$\floor{\sqrt{a}b\sqrt{a}}\leq b$ by~\sref{ad-contraposed}.
\begin{point}{30}%
Now,
let~$e$ be a projection below~$a$ and~$b$,
that is, $e\sqrt{a}=e$ and~$eb=e$.
We must show that~$e\leq \floor{\sqrt{a}b\sqrt{a}}$,
or equivalently, $e\leq \sqrt{a}b\sqrt{a}$,
or put yet differently, $e\sqrt{a}b\sqrt{a}=e$.
But this is obvious: $e=e\sqrt{a}=eb\sqrt{a}=e\sqrt{a}b\sqrt{a}$.\qed
\end{point}
\end{point}
\end{point}
\end{parsec}
\begin{parsec}{580}%
\begin{point}{10}%
Having seen that~$\floor{\sqrt{a}b\sqrt{a}} = \floor{a}\cap\floor{b}$
in~\sref{floor-sequential-product}
one might wonder whether
there is a similar expression for $\ceil{\sqrt{a}b\sqrt{a}}$,
but this doesn't seem to exist.
However,
for projections
$p$ and~$q$
we have
$\ceil{pqp}= p\cap (p^\perp \cup q)$
as we'll show below.
\end{point}
\begin{point}{20}[floor-difference]{Lemma}%
Let~$p$ be a projection,
and let~$a$ be an effect of a von Neumann algebra
with $a\leq p$.
We have $p-\ceil{a}=\floor{p-a}$.
\begin{point}{30}{Proof}%
We must show that $p-\ceil{a}$ is the greatest projection below $p-a$.
To begin, $p-\ceil{a}\leq p-a$,
because $a\leq \ceil{a}$.
Further, since~$a\leq p$, we have $\ceil{a}\leq p$,
and so~$p-\ceil{a}$ is a projection
(by~\sref{projection-below-projection}).

Let~$q$ be a projection below~$p-a$.
We must show that~$q\leq p-\ceil{a}$.
The trick is to note that~$a\leq p-q$.
Since~$p-q$ is a projection (by~\sref{projection-below-projection}
because $q\leq p-a\leq p$),
we have $\ceil{a}\leq p-q$,
and so $q\leq p-\ceil{a}$.\qed
\end{point}
\end{point}
\begin{point}{40}[ceil-sequential-product]{Proposition}%
We have $\ceil{pqp}=p\cap (p^\perp \cup q)$
for all projections~$p$ and~$q$ from a von Neumann algebra.
\begin{point}{50}[ceil-sequential-product-1]{Proof}%
Observe that $(\ p\cap (p^\perp \cup q)\ )^\perp 
= p^\perp \cup(p\cap q^\perp)$.
Since~$p^\perp$ and $p\cap q^\perp$ are disjoint,
we have $p^\perp \cup (p\cap q^\perp) = p^\perp + p\cap q^\perp$,
and so $p\cap (p^\perp \cup q) = p-p\cap q^\perp$.

By point~\sref{ceil-sequential-product-1}, 
it suffices to show that~$\ceil{pqp}=p- p\cap q^\perp$,
that is, $p-\ceil{pqp}=p\cap q^\perp$.
Since $p-\ceil{pqp} = \floor{p-pqp}$
by~\sref{floor-difference} and $\floor{pq^\perp p}=p\cap q^\perp$
by~\sref{floor-sequential-product} we are done.\qed
\end{point}
\end{point}
\end{parsec}
\subsection{Range and Support}
\begin{parsec}{590}%
\begin{point}{10}[ceill]{Notation}%
Let~$\scrA$ be a von Neumann algebra.
Because it will be very convenient
we extend the definition of~$\ceil{b}$
to all positive~$b$ from~$\scrA$
by 
$\Define{\ceil{b}}:=\ceil{\|b\|^{-1} b}$%
\index{ceiling}%
\index{*ceil@$\ceil{\,\cdot\,}$!$\ceil{a}$, ceiling}
when~$b\nleq 1$.
    Note that---contrary to what the notation suggests---we
    do \emph{not} have  $b\leq \ceil{b}$, 
    for those~$b\nleq 1$.

Now, given an arbitrary element~$b$ of~$\scrA$,
we'll call $\Define{\ceill{b}}:=\ceil{b^*b}$
the \Define{support (projection)} of~$b$,%
\index{*ceill@$\ceill{a}$, support}
	and~$\Define{\ceilr{b}}:=\ceil{bb^*}$%
\index{*ceilr@$\ceilr{a}$, range}
the \Define{range (projection)} of~$b$.
\end{point}
\begin{point}{20}{Remark}%
Some explanation is in order here.
We did not just introduce
the range
and support notation for its own sake,
but will use it extensively in~\S\ref{S:division}
thanks to calculation  rules
such as $ab=0\iff \ceill{a}\ceilr{b}=0$
(see~\sref{mult-cancellation}).
The notation was chosen such 
that $\ceilr{b} b=b= b\ceill{b}$
(see~\sref{ceill-basic}).
Good examples are 
\begin{equation*}
	\ceill{\,\ketbra{x}{y}\,}
\,= \,\ketbra{y}{y}
\qquad\text{and}\qquad
\ceilr{\,\ketbra{x}{y}\,}
\,=\, \ketbra{x}{x}
\end{equation*}
for unit vectors~$x$ and~$y$ from a Hilbert space~$\scrH$.
\end{point}
\begin{point}{30}[ceil-basic]{Exercise}%
Let~$a$ and~$b$ be positive elements of a von Neumann algebra~$\scrA$.
\begin{enumerate}
\item
Given a projection~$p$ in~$\scrA$
show that $pa=a$ iff $ap=a$ iff $\ceil{a}\leq p$.

(In particular, $\ceil{a}$ is the least projection~$p$ of~$\scrA$
with $a p=a$.)

\item
Show that $\ceil{a}a=a\ceil{a}$,
and if fact, if $b\in\scrA$ commutes with~$a$
then~$b$ commutes with~$\ceil{a}$.

\item
Show that~$a=0$ iff $\ceil{a}=0$.

\item
Show that $\ceil{a}=\ceil{\lambda a}$
for every~$\lambda>0$.

Show that $\ceil{a+b}=\ceil{a}\cup\ceil{b}$.
\item
Show that $\ceil{a^2}=\ceil{a}$.
\end{enumerate}
\spacingfix%
\end{point}%
\begin{point}{40}[ceil-pos-part]{Exercise}%
Let~$a$ be a self-adjoint element of a von Neumann algebra.
\begin{enumerate}
\item
Show that~$\ceil{a_+}\ceil{a_-}= 0$.
(Hint: recall from~\sref{cstar-pos-neg-part} that~$a_+a_-=0$.)
\item
	Show that~$\ceil{a_+}a = a\ceil{a_+} = a_+$
	and~$\ceil{a_-}a=a\ceil{a_-} = -a_-$.
\end{enumerate}
\spacingfix%
\end{point}
\begin{point}{50}[ceil-suprema]{Exercise}%
Show that $\ceil{\bigvee D} = \bigcup_{d\in D}\ceil{d}$
for every bounded directed set of \emph{positive}
elements of a von Neumann algebra~$\scrA$.
\end{point}
\begin{point}{60}[ceill-basic]{Exercise}%
Let~$a$ and~$b$  be elements of a von Neumann algebra.
\begin{enumerate}
\item
Show that $\ceill{a}\equiv \ceil{a^*a}$
is the least projection~$p$ of~$\scrA$
with $ap =a$.

\item
Show that $\ceilr{a}\equiv \ceil{aa^*}$
is the least projection~$p$ of~$\scrA$
with $pa=a$.

\item
Show that $\ceill{a^*}=\ceilr{a}$
and $\ceilr{a^*}=\ceill{a}$.

\item
Show that~$\ceill{ab}\leq\ceill{b}$
and~$\ceilr{ab}\leq\ceilr{a}$.
\end{enumerate}
\spacingfix%
\end{point}
\begin{point}{70}[hilb-ceil]{Exercise}%
Let~$T$ be a bounded operator on a Hilbert space~$\scrH$.
\begin{enumerate}
\item
Show that $\ceilr{T}$
is the projection onto the closure
$\overline{\Ran(T)}$ of the range of~$T$.
\item
Show that $\ceill{T}$
is the projection onto the \emph{support}
of~$T$, i.e.~the orthocomplement
$\Ker(T)^\perp$ of the kernel of~$T$.
\item
Show that $\floor{T}$ is the projection
on $\{\,x\in\scrH\colon\, Tx=x\,\}$
when~$T$ is an effect.
\end{enumerate}
\spacingfix%
\end{point}%
\end{parsec}%
\begin{parsec}{600}%
\begin{point}{10}[ceil-functionals-lemma]{Lemma}%
Given a positive element~$a$
of a von Neumann algebra~$\scrA$ and an
np-functional~$\omega\colon \scrA\to\C$
we have~$\omega(a)=0$ iff~$\omega(\ceil{a})=0$.
\begin{point}{20}{Proof}%
Note that if~$a=0$,
the stated result is clearly correct,
and the other case, when~$\|a\| \neq 0$,
the problem reduces to the case
that $0\leq a\leq 1$
by replacing $a$ by~$\frac{a}{\|a\|}$.
So let us just assume that~$a\in [0,1]_\scrA$
to begin with.
For similar reasons, we may assume that $\omega (1)\leq 1$.

Now, since~$0\leq a\leq \ceil{a}$
we have $0\leq \omega(a)\leq \omega(\ceil{a})$,
so~$\omega(\ceil{a})=0\implies \omega(a)=0$
is obvious.
It remains to be shown that $\omega(\ceil{a})=0$
given~$\omega(a)=0$.
Since  $\ceil{a}=\bigvee_n a^{\nicefrac{1}{2^n}}$ 
(by~\sref{vna-ceil})
and~$\omega$ is normal,
we have $\omega(\ceil{a})=\bigvee_n \omega(a^{\nicefrac{1}{2^n}})$,
and so it suffices to show that~$\omega(a^{\nicefrac{1}{2^n}})=0$
for each~$n$.
As a result of Kadison's inequality
(see~\sref{omega-norm-basic})
we have~$\omega(\sqrt{a})^2\leq \omega(a)=0$,
and so~$\omega(\sqrt{a})=0$.
Since then~$\omega(\smash{\sqrt{\sqrt{a}}})=0$
by the same token, and so on,
we get $\omega(a^{\nicefrac{1}{2^n}})=0$
for all~$n$ by induction.\qed
\end{point}
\end{point}
\begin{point}{30}[ceil-functionals]{Proposition}%
For positive elements~$a$ and~$b$ of a von Neumann algebra~$\scrA$,
\begin{equation*}
\ceil{a}\leq \ceil{b}
\qquad\iff\qquad
\forall \omega\,[\quad\omega(b)=0\  \implies  \ \omega(a)=0\quad ],
\end{equation*}
where~$\omega$ ranges over all np-functionals~$\scrA\to\C$.
\begin{point}{40}{Proof}%
When~$\ceil{a}\leq\ceil{b}$
and~$\omega$ is an np-functional on~$\scrA$
with~$\omega(b)=0$,
then~$0\leq \omega(\ceil{a})\leq \omega(\ceil{b})=0$
(by~\sref{ceil-functionals-lemma}),
and so~$\omega(\ceil{a})=0$,
so that~$\omega(a)=0$
(again
by~\sref{ceil-functionals-lemma}).

For the other direction,
assume that~$\omega(b)=0\implies \omega(a)=0$
for every np-functional~$\omega$ on~$\scrA$;
we must show that~$\ceil{a}\leq \ceil{b}$,
or in other words, $\smash{\ceil{b}}^\perp \ceil{a}\smash{\ceil{b}}^\perp =0$.
Let~$\omega\colon \scrA\to\C$ be an arbitrary np-functional;
it suffices to show that 
$\omega(\,\smash{\ceil{b}}^\perp \ceil{a}\smash{\ceil{b}}^\perp \,)=0$.
Since~$\smash{\ceil{b}}^\perp b \smash{\ceil{b}}^\perp =0$
we have~$\omega(\smash{\ceil{b}}^\perp b\smash{\ceil{b}}^\perp)=0$
and so~$\omega(\smash{\ceil{b}}^\perp a \smash{\ceil{b}}^\perp)=0$
(by assumption, because
$\omega(\smash{\ceil{b}}^\perp(\,\cdot\,)\smash{\ceil{b}}^\perp)$
is an np-functional on~$\scrA$ as well),
which implies that
$\omega(\smash{\ceil{b}}^\perp\ceil{a}\smash{\ceil{b}}^\perp)=0$
by~\sref{ceil-functionals-lemma}.\qed
\end{point}
\end{point}
\begin{point}{50}[ncp-ceil]{Proposition}%
Let $f\colon \scrA\to\scrB$ be an np-map
between von Neumann algebras.
Then $\ceil{f(a)}=\ceil{f(\ceil{a})}$
for every $a\in\scrA_+$.
\begin{point}{60}{Proof}%
By~\sref{ceil-functionals}
it suffices to show that
$\omega(f(a))=0$ iff $\omega(f(\ceil{a}))=0$
for every np-functional $\omega\colon \scrB\to\C$,
and this is indeed the case by~\sref{ceil-functionals-lemma}.\qed
\end{point}
\end{point}
\begin{point}{70}[ceil-fundamental]{Exercise}%
Let~$a$ and~$b$ be elements of a von Neumann algebra~$\scrA$.
\begin{enumerate}
\item
Deduce from~\sref{ncp-ceil}
that
$\ceil{a^*ba} = \ceil{a^*\ceil{b}a}$
when~$b\geq 0$.
\item
Conclude that~$\ceill{ab}=\ceill{\ceill{a}b}$
and~$\ceilr{ab}=\ceilr{a\ceilr{b}}$
(see~\sref{ceill}).
\end{enumerate}
\spacingfix%
\end{point}%
\begin{point}{80}[mult-cancellation]{Exercise}%
Let~$a$ and~$b$ be elements of a von Neumann algebra~$\scrA$.
\begin{enumerate}
\item
Show that $cb=0$ iff $\ceill{c}\ceilr{b}=0$
iff $\ceill{c}\leq \ceilr{b}^\perp$
for~$c\in\scrA$.

(Hint: if $cb=0$,
then $\ceil{b^*c^*cb}\equiv \ceil{b^*\ceil{c^*c}b}=0$
by~\sref{ceil-fundamental}.)
\item
Show that~$c_1b=c_2b \implies c_1=c_2$
for all~$c_1,c_2\in\scrA$
with~$\ceill{c_i}\leq \ceilr{b}$.
\item
Show that $b^*c_1b=b^*c_2b\implies c_1=c_2$
for all~$c_1,c_2\in\ceilr{b}\scrA\ceilr{b}$
\end{enumerate}
\spacingfix%
\end{point}%
\begin{point}{90}[ncp-union]{Exercise}%
Let~$f\colon \scrA\to\scrB$ be an np-map
between von Neumann algebras.
\begin{enumerate}
\item
Show that
$\ceil{f(p\cup q)}
= \ceil{f(p)}\cup\ceil{f(q)}$
for all projections $p$ and~$q$ in~$\scrA$.

(Hint: recall from~\sref{ceil-floor-basic} 
that $p\cup q = \ceil{\frac{1}{2}p+\frac{1}{2}q}$.)

\item
Deduce from this and~\sref{ncp-ceil} that $\ceil{f(\bigcup A)}=\bigcup_{a\in A}\ceil{f(a)}$
for every set of projections~$A$ from~$\scrA$.

\item
Show that there is a greatest projection~$e$
in~$\scrA$ with~$f(e)=0$.
\end{enumerate}
\spacingfix%
\end{point}%
\end{parsec}%
\begin{parsec}{610}%
\begin{point}{10}%
Given the rule $\ceil{f(\ceil{a})}=\ceil{f(a)}$
for an np-map~$f$ and self-adjoint~$a$
one might surmise that
the equation $\ceil{f(\ceill{a})}=\ceill{f(a)}$
holds 
for arbitrary~$a$;
but
one would be mistaken to do so.
We can, however,
recover an inequality
by assuming that~$f$ is completely positive, see~\sref{ncp-ceill}.
One of its corollaries is
that ncpsu-isomorphisms
are in fact nmiu-isomorphisms (see~\sref{iso}).
\end{point}
\begin{point}{20}[ncp-ceill]{Proposition}%
Given an ncp-map $f\colon \scrA\to\scrB$
between von Neumann algebras
we have,
for all
$a\in\scrA$,
\begin{equation*}
	\ceil{f(\,\ceill{a}\,)}\,\leq\, \ceill{f(a)\,}
	\qquad\text{and}\qquad
	\ceil{f(\,\ceilr{a}\,)}\,\leq\, \ceilr{\,f(a)}.
\end{equation*}
\spacingfix%
\begin{point}{30}{Proof}%
Since~$f(a)^*f(a)\leq \|f(1)\|^2\,f(a^*a)$
 by~\sref{cp-cs},
we get~$\ceill{f(a)\,}\equiv \ceil{f(a)^*f(a)}
	\leq \ceil{\,\smash{\|f(1)\|^2f(a^*a)}\,}
	\leq \ceil{f(a^*a)}
	=\ceil{f(\ceil{a^*a})}
	\equiv \ceil{f(\ceill{a})}$.

One obtains $\ceil{f(\,\ceilr{a}\,)}\leq \ceilr{\,f(a)}$
along similar lines.\qed
\end{point}
\end{point}
\end{parsec}%
\begin{parsec}{620}%
\begin{point}{10}{Proposition}%
Let~$f\colon \scrA\to\scrB$ be a ncpsu-map
between von Neumann algebras.
Then~$\floor{f(a)}=\floor{f(\floor{a})}$
for every effect~$a$ from~$\scrA$.
\begin{point}{20}{Proof}%
Since~$\floor{a}\leq a$,
we have~$\floor{f(\floor{a})}\leq \floor{f(a)}$.
Thus we only need to show that~$\floor{f(a)}\leq \floor{f(\floor{a})}$,
or equivalently, $\floor{f(a)}\leq f(\floor{a})$.
We have
\begin{equation*}
\floor{f(a)}
\ \overset{\sref{ceil-floor-basic}}{=\joinrel=\joinrel=}\ 
\floor{f(a)^2}
\ \stackrel{\sref{inner-product-basic}}{\leq}\  
\floor{f(a^2)} \ \leq\ \floor{f(a)},
\end{equation*}
and so~$\floor{f(a)}=\floor{f(a^2)}$.
By induction,
and similar reasoning,
we get~$\floor{f(a)}=\floor{f(a^{2^n})}\leq f(a^{2^n})$
for every~$n$,
and so
$\floor{f(a)}\leq \bigwedge_n f(a^{2^n})
= f(\bigwedge_n a^{2^n})=f(\floor{a})$,
where we used that~$f$ is normal,
and~$\floor{a}=\bigwedge_n a^{2^n}$ (see~\sref{vna-floor}).\qed
\end{point}
\end{point}
\end{parsec}
\subsection{Carrier and Commutant}
\begin{parsec}{630}%
\begin{point}{10}[carrier]{Definition}%
The \Define{carrier}%
\index{carrier}
of an np-map $f\colon \scrA\to\scrB$
between von Neumann algebras
(written $\Define{\ceil{f}}$)%
\index{*ceil@$\ceil{\,\cdot\,}$!$\ceil{f}$, carrier of an np-map}
is the least projection~$p$
with~$f(p^\perp)=0$
(which exists by~\sref{ncp-union}.)
\end{point}
\begin{point}{20}[carrier-basic]{Exercise}%
Let~$f,g\colon \scrA\to\scrB$
be np-maps between von Neumann algebras.
\begin{enumerate}
\item
Show that $\ceil{\lambda f}=\ceil{f}$
for all~$\lambda>0$.
\item
Show that~$\ceil{f+g}=\ceil{f}\cup\ceil{g}$.
\item
Show that~$\ceil{f}=1$ iff~$f$ is faithful.
\item
Assuming~$f$ is multiplicative
show that~$\ceil{f}=1$ iff~$f$ is injective.

(There is more to be said about
the carrier of an nmiu-map, see~\sref{carrier-miu}.)
\end{enumerate}
\spacingfix%
\end{point}
\begin{point}{30}{Exercise}%
\begin{enumerate}
\item
Given an element~$a$ of a von Neumann algebra~$\scrA$
show that 
\begin{equation*}
\ceil{a^*(\,\cdot\,)a}\ =\ \ceil{aa^*}\equiv \ceilr{a}
\end{equation*}
where~$a^*(\,\cdot\,)a$
is interpreted as an np-map~$\scrA\to\scrA$.
\item
Given a bounded operator~$T\colon \scrH\to\scrK$
between Hilbert spaces
show that~$\ceil{T^*(\,\cdot\,)T}$
is the projection onto~$\overline{\Ran(T)}$
when~$T^*(\,\cdot\,)T$
is interpreted 
as a map
$\scrB(\scrK)\to\scrB(\scrH)$.
\item
Show that~$\ceil{\left<x,(\,\cdot\,)x\right>}=\ketbra{x}{x}$
for any unit vector~$x$ from a Hilbert space~$\scrH$
when~$\left<x,(\,\cdot\,)x\right>$
is interpreted as a map~$\scrB(\scrH)\to\C$.

(But be warned: when~$\scrA$ is a von Neumann subalgebra of~$\scrB(\scrH)$
the carrier of the restriction 
$\left<x,(\,\cdot\,)x\right>\colon \scrA\to\C$
might differ from~$\ketbra{x}{x}$ because the former is in~$\scrA$,
while the latter may not be, see~\sref{carrier-vector-state}.)
\end{enumerate}
\spacingfix%
\end{point}%
\begin{point}{40}[cp-comprehension]{Lemma}%
Let~$f\colon \scrA\to\scrB$ be a p-map between
$C^*$-algebras,
and let~$p$ be an effect of~$\scrA$ with~$f(p^\perp)=0$.
Then $f(a)=f(pa)=f(ap)=f(pap)$
for all~$a\in\scrA$.
\begin{point}{50}{Proof}%
Assume~$\scrB=\C$ for now.
Since~$p^\perp \leq 1$,
we have $(p^\perp)^2=\sqrt{p^\perp}p^\perp \sqrt{p^\perp}
\leq p^\perp$,
and so~$0\leq f(\,(p^\perp)^2\,) \leq f(p^\perp) = 0$,
giving us $f(\,(p^\perp)^2\,)=0$.
Since
$\left|\smash{f(p^\perp a)}\right|^2
\leq f(\,(p^\perp)^2\,)\  f(a^*a)\,=\, 0$
by Kadison's inequality, \sref{omega-norm-basic}, 
we get~$f(p^\perp a)=0$, and so~$f(pa)=f(a)$ for all~$a\in\scrA$.
In particular,
$f(ap)=f(pa^*)^*=f(a^*)^*=f(a)$
for all~$a\in\scrA$,
and so~$f(pap)=f(pa)=f(a)$
for all~$a\in \scrA$.

Letting~$\scrB$ be again arbitrary,
and given~$a\in\scrA$,
note that
since the states on~$\scrB$ 
are separating (by~\sref{states-order-separating})
it suffices 
to show that $\omega(f(a))=\omega(f(ap))=\omega(f(pa))=
\omega(f(pap))$ for all
states $\omega \colon \scrB\to\C$.
But this follows from the previous paragraph
since~$\omega\circ f$
is a p-map into~$\C$.\qed
\end{point}
\end{point}
\begin{point}{60}[carrier-fundamental]{Corollary}%
Given an np-map $f\colon \scrA\to\scrB$
between von Neumann algebras
we have $f(a)=f(\ceil{f}a) = f(a\ceil{f})=f(\ceil{f}a\ceil{f})$
for all~$a\in\scrA$.
\end{point}
\end{parsec}
\begin{parsec}{640}%
\begin{point}{10}%
We turn
to the task of showing that every element
of a von Neumann
algebra is the norm limit of linear combinations
of projections in~\sref{projections-norm-dense}.
We'll deal with the  commutative case first
(see~\sref{abelian-projections-norm-dense}).
\end{point}
\begin{point}{20}[abelian-projections-norm-dense]{Proposition}%
Every element~$a$ of a commutative von Neumann algebra~$\scrA$
is the norm limit
of linear combinations of projections.
\begin{point}{30}{Proof}%
By~\sref{ngelfand}
it suffices to show that the linear span
of projections is norm dense
in~$C(\spec(\scrA))$.
For this, in turn, it suffices 
by Stone--Weierstra\ss{}' theorem
(see~\sref{stone-weierstrass})
to show that the projections in~$C(\spec(\scrA))$
separate the points of~$\spec(\scrA)$
in the sense
that given $x,y\in\spec(\scrA)$
with~$x\neq y$
there is a projection~$f$ in~$C(\spec(\scrA))$
with~$f(x)\neq f(y)$.
Since~$\spec(\scrA)$
is Hausdorff
there are for such~$x$ and~$y$
disjoint open subsets $U$ and~$V$
of~$\spec(\scrA)$
with~$x\in U$ and~$y\in V$.

Then
$f:=\mathbf{1}_{\overline{U}}$
is a projection in~$C(\spec(\scrA))$
(continuous
because~$\overline{U}$
is clopen by~\sref{vn-spectrum-extremally-disconnected})
with~$f(x)=0\neq 1=f(y)$
since $x\in \overline{U}\subseteq \spec(\scrA)\backslash V$,
and so~$y\notin \overline{U}$.\qed
\end{point}
\end{point}
\end{parsec}
\begin{parsec}{650}%
\begin{point}{10}%
To reduce the general case
to the commutative case
we need the following tool
(that will be useful 
later on too for different reasons).
\end{point}
\begin{point}{20}[commutant]{Definition}%
Given a subset~$S$ of a von Neumann algebra~$\scrA$
	the \Define{commutant}\index{commutant} of~$S$
is the set, denoted by~$\Define{\smash{S^\square}}$,%
\index{*square@$S^\square$, commutant of~$S$}
of all~$a\in\scrA$ with $as=sa$ for all~$s\in S$.

The commutant of~$\scrA$ itself
is denoted by~$\Define{Z(\scrA)}:=\scrA^\square$%
\index{ZA@$Z(\scrA)$, centre of~$\scrA$}
and is called the \Define{centre}%
\index{centre of a von Neumann algebra} of~$\scrA$.
(Its elements, called \emph{central}, are the subjects of the next section.)
\end{point}
\begin{point}{30}[commutant-basic]{Exercise}%
Let~$S$ and~$T$ be subsets of a von Neumann algebra~$\scrA$.
\begin{enumerate}
\item
Show that $S \subseteq T^\square$ iff $T \subseteq S^\square$.

Show that $S\subseteq T$ entails $T^\square \subseteq S^\square$.

Show that $S\subseteq S^{\square\square}$,
and  $S^{\square\square\square}=S^\square$.
\item
Show that $S^\square$ is closed under addition,
(scalar) multiplication,
contains the unit of~$\scrA$,
and is ultraweakly closed.
\item
Show that the commutant $S^\square$ need not be closed under involution.\\
(Hint: compute 
$\{\bigl(\begin{smallmatrix}0&1\\0&0\end{smallmatrix}\bigr)\}^\square$
in $M_2$.)

Suppose~$S$ is closed under involution.

Show~$S^\square$ is closed under involution as well,
and conclude that in that case~$S^\square$
is a von Neumann subalgebra of~$\scrA$.

Show that~$Z(\scrA)$ is a von Neumann subalgebra of~$\scrA$.

Show that~$S^{\square\square}$
is a von Neumann subalgebra of~$\scrA$
with~$S\subseteq S^{\square\square}$.

Show that if~$S$ is commutative (i.e.~$S\subseteq S^\square$), 
then so is~$S^{\square\square}$.
\item
In particular,
if~$\scrB$ is a von Neumann subalgebra of~$\scrA$,
then~$\scrB^{\square\square}$
is a von Neumann subalgebra of~$\scrA$
with $\scrB\subseteq \scrB^{\square\square}$.

Show that~$(\,\scrA\cap \C\,)^\square=\scrA$,
and so~$(\,\scrA\cap\C\,)^{\square\square} = Z(\scrA)$.

So in general~$\scrB^{\square\square}$ needn't equal~$\scrB$.
Nevertheless,
we'll see in \sref{proto-double-commutant} that $\scrB^{\square\square}=\scrB$
when~$\scrA$ is of the form $\scrA=\scrB(\scrH)$
for some Hilbert space~$\scrH$.
\item
Given a von Neumann subalgebra~$\scrB$
of~$\scrA$
verify that $Z(\scrB)=\scrB\cap \scrB^{\square}$.
\end{enumerate}
\spacingfix%
\end{point}%
\begin{point}{40}[projections-norm-dense]{Proposition}%
Every self-adjoint element~$a$ of a von Neumann algebra~$\scrA$
is the norm limit
of linear combinations
of projections from~$\smash{\{a\}}^{\square\square}$.
\begin{point}{50}{Proof}%
Since~$a$ is an element
of the by~\sref{commutant-basic}
commutative von Neumann subalgebra~$\smash{\{a\}}^{\square\square}$
of~$\scrA$,
$a$ is the norm limit of linear combinations
of projections from~$\smash{\{a\}}^{\square\square}$
by~\sref{abelian-projections-norm-dense}.\qed
\end{point}
\end{point}
\end{parsec}
\begin{parsec}{660}%
\begin{point}{10}%
The carriers of np-functionals
play such an important role in the theory
that we decided to give them a name.
\end{point}
\begin{point}{20}{Definition}%
We call a projection~$p$ of a von Neumann algebra~$\scrA$
\Define{ultracyclic}%
\index{ultracyclic projection}
if~$p=\ceil{\omega}$
for some np-map $\omega\colon \scrA\to\C$.
\begin{point}{30}{Remark}%
Some explanation of this terminology
is in order.
A projection~$E$
in a von Neumann subalgebra~$\scrR$
of~$\scrB(\scrH)$
is usually defined to be \Define{cyclic}%
\index{cyclic projection}
when~$E$ is the projection 
onto~$\smash{\overline{\scrR^\square x}}$
for some~$x\in \scrH$
(see Definition~5.5.8~\cite{kr}).
With~\sref{carrier-vector-state} and~\sref{double-commutant}
we'll be able to see that
this amounts to requiring that~$E$
be
the carrier of the vector functional
$\left<x,(\,\cdot\,)x\right>\colon \scrR\to\C$.
So, loosely speaking,
a cyclic projection
is the carrier of a vector functional
with respect to some fixed Hilbert space,
while an ultracyclic projection
is the carrier of a vector functional
with respect to some arbitrary Hilbert space.
\end{point}
\end{point}
\begin{point}{40}[ultracyclic-basic]{Exercise}%
Let~$\scrA$ be a von Neumann algebra.
Verify the following facts.
\begin{enumerate}
\item
If~$p$ and~$q$ are ultracyclic projections in~$\scrA$,
then~$p\cup q$ is ultracyclic.
\item
If~$p \leq q$ are projections in~$\scrA$,
and~$q$ is ultracyclic,
then~$p$ is ultracyclic.
\item
Every projection~$p$ in~$\scrA$ is a directed supremum
of ultracyclic projections.
In fact, $p=\bigvee_\omega \ceil{\omega}$
where~$\omega$ ranges over the np-functionals on~$\scrA$
with $\omega(p^\perp)=0$.
(Hint: first consider~$p=1$.)
\item
Every projection~$p$ in~$\scrA$ is the 
sum of ultracyclic projections:
there are np-functionals $(\omega_i)_i$
on~$\scrA$ with $p=\sum_i \ceil{\omega_i}$.
\end{enumerate}
\spacingfix%
\end{point}%
\end{parsec}%
\subsection{Central Support and Central Carrier}
\begin{parsec}{670}%
\begin{point}{10}{Definition}%
An element~$a$ of a von Neumann algebra~$\scrA$
is called~\Define{central}%
\index{central!element of a von Neumann algebra}
when~$ab=ba$ for all~$b\in\scrA$
(that is, when~$a\in Z(\scrA)$, see~\sref{commutant-basic}).
\end{point}
\begin{point}{20}[central-examples]{Examples}%
\begin{enumerate}
\item
In a commutative von Neumann algebra
every element is central.
\item
An element~$a$ of a direct sum~$\bigoplus_i \scrA_i$
of von Neumann algebras
is central iff~$a_i$ is central for each~$i$.
\item
In~$\scrB(\scrH)$,
where~$\scrH$ is a Hilbert space,
only the scalars are central.

Indeed,
given a positive central element~$A$ 
of~$\scrB(\scrH)$,
we have $\left<x, A\|y\|^2 x\right>
        = \left<x,\,(A \ketbra{x}{y})\, y\right>
        = \left<x,\,(\ketbra{x}{y}A)\, y\right>
= \left<x,\smash{\|\sqrt{A}y\|^2 }x\right>$
for all~$x,y\in\scrH$,
and so~$A\|y\|^2= \|\sqrt{A}y\|^2$
for all~$y\in\scrH$.
Hence~$A$ is (zero or) a scalar.
\end{enumerate}
\spacingfix%
\end{point}%
\begin{point}{30}{Remark}%
A von Neumann algebra
in which only the scalars are central
--- of which a~$\scrB(\scrH)$ is but the simplest example ---
	is called a \Define{factor}.\index{factor}
The classification of these factors
is an important part
of the theory of von Neumann algebras
that we did not need in this thesis.
\end{point}
\begin{point}{40}[central-projections-sums]{Exercise}%
Note that if a von Neumann algebra
$\scrA$ can be written as
a direct sum $\scrA\cong \scrB_1\oplus \scrB_2$,
then~$(1,0)\in\scrB_1\oplus\scrB_2$ gives
a central projection in~$\scrA$.
The converse also holds:
\begin{enumerate}
\item
Given a central projection~$c$ in~$\scrA$,
show that~$c\scrA \equiv \{\,ca\colon\, a\in\scrA\,\}$ 
is a von Neumann subalgebra of~$\scrA$
for all but the fact that $1$ need not be in~$c\scrA$.

Show~$c\scrA$
is a von Neumann algebra with~$c$ as unit,
and that $a\mapsto (ca,c^\perp a)$
gives an nmiu-isomorphism
$\scrA\to c\scrA\oplus c^\perp \scrA$.
\item
Given a family of central projections $(c_i)_i$ in~$\scrA$
with $\sum_i c_i=1$
show that $a\mapsto (c_ia)_i$
gives an nmiu-isomorphism $\scrA\to\bigoplus_i c_i \scrA$.
\end{enumerate}
\spacingfix%
\end{point}%
\end{parsec}%
\begin{parsec}{680}%
\begin{point}{10}[cceil-fundamental]{Proposition}%
Given a projection~$e$ of a von Neumann algebra~$\scrA$
\begin{equation*}
	\Define{\cceil{e}}\ :=\ 
	\bigcup_{a\in\scrA} \ceil{a^* e a}
\end{equation*}%
\index{*cceil@$\cceil{\,\cdot\,}$!$\cceil{a}$, central support}%
is the least central projection above~$e$.
\begin{point}{20}{Proof}%
Let us first show that~$\cceil{e}$ is central.
Given~$b\in\scrA$
we have $\ceill{\cceil{e}b}=\ceil{b^*\cceil{e}b}
=\bigcup_{a\in\scrA} \ceil{b^*\ceil{a^*ea}b}
= \bigcup_{a\in \scrA} \ceil{(ab)^* eab}
\leq \cceil{e}$
by~\sref{ncp-union},
which implies that~$\cceil{e}b\cceil{e}=\cceil{e}b$.
Since similarly (or consequently)
$\cceil{e}b\cceil{e} = b\cceil{e}$
we get~$b\cceil{e}=\cceil{e}b\cceil{e}
=\cceil{e}b$,
and so~$\cceil{e}$ is central.

Clearly~$e\leq \cceil{e}$.
It remains to be shown that~$\cceil{e}\leq c$
given a central projection~$c$ with~$e\leq c$.
For this it suffices to show that $\ceill{ea}\equiv \ceil{a^*ea}\leq c$
given~$a\in\scrA$.
Now, since~$e\leq c$
we have $ec=e$
and so~$eac=eca=ea$
which implies that~$\ceill{ea}\leq c$.
Thus~$\cceil{e}\leq c$.\qed
\end{point}
\end{point}
\begin{point}{30}[central support]{Definition}%
Let~$a$ be an element of a von Neumann algebra~$\scrA$.
Since given a central projection~$c$ of~$\scrA$
we have
$\cceil{\ceill{a}}\leq c$
iff $\ceill{a}\leq c$
iff $ac=a$ 
iff
$ca=a$ iff $\cceil{\ceilr{a}}\leq c$,
we see that~$\Define{\cceil{a}}:=
\cceil{\ceill{a}}=\cceil{\ceilr{a}}$
is the smallest central projection~$p$
with~$p a=a$,
which we'll call the
\Define{central support}%
\index{central support} 
	of~$a$.
\end{point}
\begin{point}{40}[cceil-basic]{Exercise}%
Let~$\scrA$ be a von Neumann algebra.
\begin{enumerate}
\item
Show that $\cceil{a}=\cceil{a^*}
=\cceil{a^*a}=\cceil{aa^*}$
for all~$a\in\scrA$.
\item
Show that~$\cceil{ \bigvee D } = \bigcup_{d\in D} \cceil{d}$
for any bounded directed subset of~$\scrA$.

Show that~$\cceil{\bigcup E} = \bigcup_{e\in E} \cceil{e}$
for any collection of projections from~$\scrA$.

Show that~$\cceil{a+b}=\cceil{\ceil{a}\cup \ceil{b}}
=\cceil{a}\cup\cceil{b}$
for all~$a,b\in \scrA$.
\item
Given~$a\in\scrA$
and a central projection~$c$ of~$\scrA$
show that~$\cceil{a}c=\cceil{ac}$.

Conclude that~$\cceil{a}\cceil{b}=\cceil{a\cceil{b}}
= \cceil{\cceil{a}b}
=\cceil{a}\cap\cceil{b}$
for all~$a,b\in\scrA$.
\end{enumerate}
\spacingfix%
\end{point}%
\end{parsec}%
\begin{parsec}{690}%
\begin{point}{10}[cceil-map-def]{Definition}%
Let $f\colon \scrA\to\scrB$
be an np-map
between von Neumann algebras.
Show that given a central effect~$c$
of~$\scrA$ we have
$f(c^\perp)=0$
iff $\ceil{f}\leq c$
iff $\cceil{\ceil{f}}\leq c$,
and so~$\Define{\cceil{f}}:=\cceil{\ceil{f}}$%
\index{*cceil@$\cceil{\,\cdot\,}$!$\cceil{f}$, central carrier}
is the least central effect
(and central projection) $p$
with~$f(p^\perp)=0$,
which we'll call the \Define{central carrier}%
\index{central carrier}
of~$f$.
\end{point}
\begin{point}{20}[prop:weakly-closed-ideal]{Proposition}%
Every two-sided ideal~$\scrD$ of a von Neumann
algebra~$\scrA$
that is closed under bounded directed suprema of self-adjoint 
elements --- for example when~$\scrA$ is ultrastrongly closed ---
is of the form~$c\scrA$
for some unique central projection~$c$ of~$\scrA$.
Moreover, $c$ is the greatest projection in~$\mathscr{D}$.
\begin{point}{30}{Proof}%
We'll obtain~$c$ as the supremum over all effects in~$\scrD$,
and to this end we'll show first that~$\scrD\cap [0,1]_\scrA$
is directed.
Since $\ceil{a}\cup \ceil{b} \equiv \ceil{\frac{1}{2}a+\frac{1}{2}b}$
(see~\sref{ceil-floor-basic})
is an upper bound for~$a,b\in\scrD\cap [0,1]_\scrA$
it suffices to show that~$\ceil{a}\in \scrD$
for all~$a\in \scrD\cap[0,1]_\scrA$,
which, in turn,
follows from
$\ceil{a}=\bigvee_n a^{\nicefrac{1}{2^n}}$,
see~\sref{vna-ceil}.

Hence~$\scrD\cap [0,1]_\scrA$ is directed,
and so we may define~$c:=\bigvee \scrD\cap[0,1]_\scrA$.
Since~$\scrD$ is a von Neumann subalgebra of~$\scrA$,
we'll have~$c\in\scrD\cap[0,1]_\scrA$,
and so~$c$ is the greatest element of~$\scrD\cap[0,1]_\scrA$.
In particular,
$c$ will be above~$\ceil{c}$ implying~$\ceil{c}=c$
and making~$c$ a projection---the greatest projection in~$\scrD$.

Given~$a\in\scrA$ we claim that~$a\in \scrD$
iff~$ca=a$.
Surely, if~$a=ca$,
then~$a=ca\in \scrD$,
because~$\scrD$ is a two-sided ideal of~$\scrA$.
Concerning the other direction,
note that given~$a\in \scrD$
the equality~$ac=a$ holds
when~$a$ is an effect by~\sref{projection-above-effect} (because~$a\leq c$),
and thus when~$a$ is self-adjoint too
(by scaling),
and hence for arbitrary~$a\in\scrD$
by  writing~$a\equiv \Real{a}+i\Imag{a}$
where~$\Real{a}$ and~$\Imag{a}$ are self-adjoint.

Note that this claim entails that~$\scrD \subseteq c\scrA$.
Since $\scrD$ is an ideal
we also have~$c\scrA\subseteq \scrD$,
and so~$\scrD= c\scrA$.
The claim also entails that~$c$ is central.
Indeed,
given~$a\in\scrA$
we have $ac\in \scrD$ 
(because~$\scrD$ is an ideal)
and so~$c(ac)=ac$ by the claim.
Since similarly~$(ca)c=ca$,
we get~$ac=ca$.

The only thing that remains to be shown is that~$c$
is unique.
To this end let~$c$ and~$c'$ be central projections 
with~$c\scrA = \mathscr{D}=c'\scrA$.
As~$c'\in\mathscr{D}=c \scrA$,
there is~$a\in \scrA$
with~$c' = ca$.
Then~$c' = c'(c')^* = caa^*c^*\leq 
cc^*\|aa^*\|=c\|a\|^2$,
and so~$c'\leq c$.
Since similarly $c\leq c'$, we get~$c=c'$.\qed
\end{point}
\end{point}
\begin{point}{40}[carrier-miu]{Corollary}%
The carrier $\ceil{f}$
of an nmiu-map $f\colon \scrA\to\scrB$
between von Neumann algebras
is central, so
$\ceil{f}=\cceil{f}$.
    Moreover, $\smash{\ker(f)=\cceil{f}^\perp\!\scrA}$.
\end{point}
\begin{point}{41}[nmiu-factors]{Exercise}%
Show using~\sref{carrier-miu} and~\sref{central-projections-sums}
 that an nmiu-map $f\colon \scrA\to\scrB$
factors as
\begin{equation*}
\xymatrix{
\scrA 
    \ar[rr]^f
    \ar[rd]_{g\colon a\mapsto \cceil{f}a}
	&
    &
    \scrB
    \\
&
    \cceil{f}\scrA \ar[ru]_{h\colon a\mapsto f(a)}
},
\end{equation*}
where the~$g$ an nmiu-surjection, and~$h$ is an nmiu-injection.
\begin{point}{42}[nmiu-image]
Use this, and~\sref{injective-nmiu-iso-on-image}, to show that~$f(\scrA)$
is a von Neumann subalgebra of~$\scrB$.
\end{point}
\end{point}
\begin{point}{50}[proto-gns-ceil]{Lemma}%
We have~$\cceil{\omega}=
\ceil{\varrho_\omega}$
for every np-functional $\omega\colon \scrA\to\C$
on a von Neumann algebra~$\scrA$,
where~$\varrho_\omega$
is as in~\sref{gns}.
\begin{point}{60}{Proof}%
Let~$e$ be a projection in~$\scrA$.
Note that $0=\|\varrho_\omega(e)(\eta_\omega(a))\|^2
\equiv \omega(a^*ea)$
iff~$\ceil{a^*ea}\leq \smash{\ceil{\omega}}^\perp$
iff~$\ceil{a\ceil{\omega}a^*}\leq e^\perp$
for all~$a\in\scrA$.
So since the~$\eta_\omega(a)$'s lie dense in~$\scrH_\omega$,
we have~$\varrho_\omega(e)=0$
iff $\varrho_\omega(e)(\eta_\omega(a))=0$ for all~$a\in\scrA$
iff $\bigcup_{a\in\scrA} \ceil{a\ceil{\omega}a^*}\leq e^\perp$.
Hence $\ceil{\varrho_\omega}=
\bigcup_{a\in\scrA}
\ceil{a\ceil{\omega}a^*} \equiv
\bigcup_{a\in\scrA}
\ceil{a^*\ceil{\omega}a}
=\cceil{\ceil{\omega}}
=\cceil{\omega}$
by~\sref{cceil-fundamental}.\qed
\end{point}
\end{point}
\begin{point}{70}[gns-ceil]{Proposition}%
Given a collection of np-functionals~$\Omega$
on a von Neumann algebra~$\scrA$
we have $\ceil{\varrho_\Omega}
=\bigcup_{\omega\in\Omega} \cceil{\omega}$
for~$\varrho_\Omega\colon \scrA\to\scrB(\scrH_\Omega)$
from~\sref{gns}.
\begin{point}{80}{Proof}%
Let~$e$ be a projection of~$\scrA$.
Since
$\varrho_\Omega(e)(x)
= \sum_{\omega\in\Omega} \varrho_\omega(x_\omega)$
by~\sref{gns}
for all~$x\in\scrH_\Omega\equiv\bigoplus_{\omega\in\Omega}\scrH_\omega$,
we have~$\varrho_\Omega(e)=0$
iff~$\varrho_\omega(e)=0$ for all~$\omega\in\Omega$
iff~$e\leq \smash{\ceil{\varrho_\omega}^\perp \equiv 
\smash{\cceil{\omega}}^\perp}$
iff~$e\leq \bigcap_{\omega\in\Omega}\smash{\cceil{\omega}}^\perp 
\equiv (\bigcup_{\omega\in\Omega} \cceil{\omega})^\perp$.
Hence~$\ceil{\varrho_\Omega}=\bigcup_{\omega\in\Omega}\cceil{\omega}$.\qed
\end{point}
\end{point}
\begin{point}{90}[vn-center-separating]{Corollary}%
For a collection~$\Omega$ of np-functionals
on a von Neumann algebra,
the following are equivalent.
\begin{enumerate}
\item
\label{vn-center-separating-1}
$\Omega$ is centre separating (see~\sref{separating}).
\item
\label{vn-center-separating-2}
A central projection~$z$ of~$\scrA$ is zero
when~$\omega(z)=0$ for all~$\omega\in\Omega$.
\item
\label{vn-center-separating-3}
The map $\varrho_\Omega\colon\scrA\to\scrB(\scrH_\Omega)$
from~\sref{gns} 
is injective.
\end{enumerate}
\spacingfix%
\begin{point}{100}{Proof}%
We've seen in~\sref{proto-gelfand-naimark}
that
\ref{vn-center-separating-1}$\iff$\ref{vn-center-separating-3},
and~\ref{vn-center-separating-1}$\Rightarrow$\ref{vn-center-separating-2}
is trivial,
which leaves us with 
\ref{vn-center-separating-2}$\Rightarrow$\ref{vn-center-separating-3}.
So assume that~$\forall \omega\in\Omega\,[\,\omega(z)=0\,]\implies z=0$ 
for every central projection~$z$
of~$\scrA$.
Then since~$\ceil{\varrho_\Omega}^\perp$
is a central projection
by~\sref{gns-ceil}
with
$\ceil{\varrho_\Omega}^\perp
=\smash{\bigl(\bigcup_{\omega\in\Omega} \smash{\cceil{\omega}\bigr)}}^\perp
=\bigcap_{\omega\in \Omega} \cceil{\omega}^\perp
\leq \smash{\cceil{\omega}}^\perp
\leq \smash{\ceil{\omega}}^\perp$
and thus
$\omega(\smash{\ceil{\varrho_\Omega}}^\perp)
\leq \omega(\ceil{\omega}^\perp)=0$
for all~$\omega\in\Omega$
we get~$\ceil{\varrho_\Omega}^\perp =0$,
and so~$\varrho_\Omega$ is injective
by~\sref{carrier-basic}.\qed
\end{point}
\end{point}
\end{parsec}
\begin{parsec}{700}%
\begin{point}{10}%
With our new-found knowledge on central elements
we can complete the classification
of commutative von Neumann algebras
we started in~\sref{classification-cvn}.
\end{point}
\begin{point}{20}[central-projection-central-carrier]{Exercise}%
Show that every central projection~$c$
of a von Neumann algebra is of
the form~$c\equiv \sum_i \cceil{\omega_i}$
for some family of np-functionals $(\omega_i)_i$ on~$\scrA$.
(Hint: take
    $(\omega_i)_i$ to be a maximal set
    of np-functionals for which the $\cceil{\omega_i}$'s are orthogonal.)
\end{point}
\begin{point}{30}[cvn]{Theorem}%
\index{von Neumann algebra!commutative}
Every commutative von Neumann algebra
is nmiu-isomorphic
to a direct sum of the form  $\bigoplus_i L^\infty(X_i)$
where~$X_i$ are finite complete measure spaces.
\begin{point}{40}{Proof}%
By~\sref{central-projection-central-carrier}
we have $1\equiv \sum_i \cceil{\omega_i} $
for some np-functionals $\omega_i\colon \scrA\to\C$,
and so~$\scrA\cong \bigoplus_i \cceil{\omega_i}\!\scrA$
by~\sref{central-projections-sums}.
Since~$\scrA$ is commutative,
and so~$\cceil{\omega_i}=\ceil{\omega_i}$,
we see that restricting~$\omega_i$
gives a faithful functional on~$\cceil{\omega_i}\!\scrA$,
which is therefore by~\sref{cvn-faithful} nmiu-isomorphic to~$L^\infty(X_i)$
for some finite complete measure space~$X_i$.
From this the stated result follows.\qed
\end{point}
\end{point}
\end{parsec}

\section{Completeness}
\begin{parsec}{710}%
\begin{point}{10}%
We set to work on the ultrastrong and bounded ultraweak completeness
of von Neumann algebras (see~\sref{vn-complete}) and their precursors:
\begin{enumerate}
\item
A linear (not necessarily positive)
functional on a von Neumann algebra
is ultraweakly continuous iff it is ultrastrongly continuous
(see~\sref{luws}).
\item
A convex subset of a von Neumann algebra
is ultraweakly closed iff it is ultrastrongly closed
(see~\sref{ultraclosed}).
\item
\emph{(Kaplansky's density theorem)}\ 
The unit ball $(\scrA)_1$
of a $C^*$-subalgebra~$\scrA$
of a von Neumann algebra~$\scrB$
is ultrastrongly dense in~$(\bar{\scrA})_1$
where~$\bar{\scrA}$ is the ultrastrong (=ultraweak,
\sref{ultraclosed}) closure of~$\scrA$
(see~\sref{kaplansky}).
\item
Any von Neumann subalgebra~$\scrA$
of~$\scrB$ is ultraweakly and ultrastrongly
closed in~$\scrB$
(see~\sref{vnsac}).
\item
The von Neumann algebra~$\scrB(\scrH)$ 
of bounded operators on a Hilbert space~$\scrH$
is ultrastrongly 
(\sref{bh-us-complete})
and bounded ultraweakly complete
(\sref{bh-bounded-uw-complete}).
\end{enumerate}
\spacingfix%
\end{point}%
\end{parsec}%
\subsection{Closure of a Convex Subset}
\begin{parsec}{720}%
\begin{point}{10}%
We saw in~\sref{npuws}
that a \emph{positive} linear functional~$f$
on a von Neumann algebra 
is ultrastrongly continuous iff it is ultraweakly continuous.
In this section, we'll show that the same result holds
for an arbitrary linear functional~$f$.
Note that if~$f$ is ultraweakly continuous,
then~$f$ is automatically ultrastrongly continuous
(because ultrastrong convergence implies ultraweak convergence).
For the other direction,
we'll show that if~$f$ is ultrastrongly continuous,
then~$f$ can be written as a linear combination
$f\equiv \sum_{k=0}^3 i^k f_k$
of np-maps $f_0,\dotsc,f_3$,
and must therefore be ultraweakly continuous.
We'll need the following tool.
\end{point}
\begin{point}{20}[bstaromega]{Definition}%
Let~$\scrA$ be a von Neumann algebra.
Given an np-map $\omega\colon \scrA\to\C$,
and~$b\in \scrA$,
define~$\Define{b*\omega}\colon \scrA\to \C$%
\index{*astaromega@$a*\omega$}
	by
$(b*\omega)(a)=\omega(b^*ab)$ for all~$a\in \scrA$.
\end{point}
\begin{point}{30}[bstaromega-basic]{Exercise}%
Let~$\omega\colon \scrA\to \C$ be an np-map on a von Neumann algebra~$\scrA$.
\begin{enumerate}
\item
Note that $b*\omega\colon \scrA\to\C$
is an np-map for all~$b\in \scrA$.

Show that $\left|\omega(a^*bc)\right| 
\,\leq\, \|\omega\|\,\|a\|_\omega\, \|b\|\, \|c\|_\omega$
for all~$a,b,c\in\scrA$.

Deduce that $\| b*\omega - b'*\omega\|
\,\leq\, \|\omega\| \,\|b-b'\|_\omega\, (\|b\|_\omega + \|b'\|_\omega)$
for all~$b,b'\in\scrA$.

\item
Let~$b_1,b_2,\dotsc$ be a sequence in~$\scrA$
which is Cauchy with respect to~$\|\,\cdot\,\|_\omega$.
Show that the sequence~$b_1*\omega,\,b_2*\omega,\,\dotsc$ 
is Cauchy (in the operator norm
on bounded linear functionals $\scrA\to\C$),
and converges to a bounded linear map~$f\colon \scrA\to\C$.
Show that~$f$ is an np-map.
\end{enumerate}
\spacingfix%
\end{point}%
\begin{point}{40}{Exercise}%
Let~$f\colon \scrA\to \C$ be an ultrastrongly continuous linear
functional on a von Neumann algebra~$\scrA$.
Show that there are an np-map
$\omega\colon \scrA\to \C$
and $\delta>0$
with $\left|f(a)\right|\leq 1$
for all~$a\in\scrA$ with $\|a\|_\omega \leq \delta$.

(Keep this in mind when reading the following lemma.)
\end{point}
\begin{point}{50}[normal-functionals-lemma]{Lemma}%
Let~$\omega\colon \scrA\to\C$ be an np-map,
and let~$f\colon \scrA\to \C$ be a linear map.
The following are equivalent.
\begin{enumerate}
\item\label{normal-functionals-lemma-0}
$\left|f(a)\right|\leq B$ for all~$a\in \scrA$
with $\|a\|_\omega\leq \delta$, for some $\delta,B>0$;
\item\label{normal-functionals-lemma-1}
$\left|f(a)\right| \leq B \|a\|_\omega$ for all~$a\in\scrA$,
for some~$B>0$;
\item\label{normal-functionals-lemma-2}
$f(a)=[b,a]_\omega$ for all~$a\in\scrA$, 
for some~$b\in\scrH_\omega$
(where $\scrH_\omega$ is the Hilbert space completion of~$\scrA$
with respect to the inner-product
$[\,\cdot\,,\,\cdot\,]_\omega$).
\item\label{normal-functionals-lemma-3}
$f\equiv f_0+if_1-f_2-if_3$
where $f_0,\dotsc,f_3\colon \scrA\to \C$
are np-maps for which there is~$B>0$
such that~$f_k(a)\leq B \omega(a)$ for all~$a\in\pos{\scrA}$ 
and~$k$. 
\end{enumerate}%
\spacingfix%
\begin{point}{60}{Proof}%
We make a circle.
\begin{point}{70}{%
\ref{normal-functionals-lemma-3}$\Longrightarrow$%
\ref{normal-functionals-lemma-0}}%
For $a\in \scrA$ and~$k$, we have
$\left|f_k(a)\right|^2 \leq
f_k(1)\,f_k(a^*a) \leq f_k(1)B \,\omega(a^*a)$,
giving~$\left|f_k(a)\right| \leq (f_k(1)B)^{\nicefrac{1}{2}} \|a\|_\omega$,
and so~$\left|f(a)\right|\leq \tilde{B} \|a\|_\omega$,
where 
\begin{equation*}
	\textstyle
	\tilde{B} \,=\, B^{\nicefrac{1}{2}}\sum_{k=0}^3f_k(1)^{\nicefrac{1}{2}}.
\end{equation*}
Hence~$\left|f(a)\right|\leq \tilde{B}$
for all~$a\in\scrA$ with $\|a\|_\omega\leq 1$.
\end{point}
\begin{point}{80}{\ref{normal-functionals-lemma-0}$\Longrightarrow$%
\ref{normal-functionals-lemma-1}}%
Let~$a\in\scrA$, and~$\varepsilon>0$ be given.
Then for~$\tilde{a}:=\delta(\varepsilon+\|a\|_\omega)^{-1}\,a$,
we have $\|\tilde{a}\|_\omega\leq \delta$,
and so~$\left|f(\tilde{a})\right|
\equiv \delta(\varepsilon +\|a\|_\omega)^{-1} \,\left|f(a)\right|
\leq B$,
which entails $\left|f(a)\right|\leq 
B\delta^{-1}(\varepsilon+\|a\|_\omega)$.
Since~$\varepsilon>0$ was arbitrary, we get~$\left|f(a)\right|\leq
B\delta^{-1}\|a\|_\omega$.
\end{point}
\begin{point}{90}{%
\ref{normal-functionals-lemma-1}$\Longrightarrow$%
\ref{normal-functionals-lemma-2}}%
Since~$\left|f(a)\right|\leq B\|a\|_\omega$ for all~$a\in\scrA$,
the map~$f$ can be extended to a bounded linear map 
$\tilde{f}\colon \scrH_\omega\to \C$.
Then by Riesz' representation theorem, \sref{riesz-representation-theorem},
there is~$b\in \scrH_\omega$ with $\tilde{f}(x)=[b,x]_\omega$
for all~$x\in \scrH_\omega$.  
In particular,
$f(a)=[b,a]_\omega$ for all~$a\in \scrA$.
\end{point}
\begin{point}{100}{\ref{normal-functionals-lemma-2}$\Longrightarrow$%
\ref{normal-functionals-lemma-3}}%
We know that~$f(a)\equiv [b,a]_\omega$ for all~$a\in \scrA$,
for some~$b\in\scrH_\omega$.
Then, by definition of~$\scrH_\omega$,
there is a sequence~$b_1,b_2,\dotsc$ in~$\scrA$
which converges to~$b$ in~$\scrH_\omega$.
Then the maps $[b_n,\,\cdot\,]_\omega\colon \scrA\to \C$
approximate~$f=[b,\,\cdot\,]_\omega$
in the sense that 
$\left|f(a)-[b_n,a]_\omega\right|=\left|[b-b_n,a]_\omega\right|
\leq \|b-b_n\|_\omega \|a\|_\omega
	\leq \|b-b_n\|_\omega \|\omega\|^{\nicefrac{1}{2}} \|a\|$
for all~$a\in \scrA$.
In particular, $[b_1,\,\cdot\,]_\omega,\,[b_2,\,\cdot\,]_\omega,\,\dotsc$
converges to~$f$ (in the operator norm).
By ``polarisation'' (c.f.~\sref{mult-polarization}),
we have $[b_n,a]_\omega = \frac{1}{4}\sum_{k=0}^3 i^kf_{k,n}(a)$,
where $f_{k,n} := (i^kb_n+1)*\omega$ is an np-map.
Since~$(i^kb_n+1)_n$ is Cauchy with respect to~$\|\,\cdot\,\|_\omega$,
we see by~\sref{bstaromega-basic} that 
$(f_{k,n})_n$ converges to an np-map $f_k\colon \scrA\to\C$
(with respect to the operator norm).
It follows that~$f=\frac{1}{4}\sum_{k=0}^3 i^k f_k$.

It remains to be shown that there is~$B>0$ with $f_k(a)\leq B\omega(a)$
for all~$k$ and~$a\in\pos{\scrA}$.
Note that since $f_{k,n}(a) \leq \|i^kb_n+1\|_\omega \,\omega(a)
\leq (\|b_n\|_\omega+1) \,\omega(a)$,
for all~$n$, $k$, and~$a\in\pos{\scrA}$,
the number $B:=\lim_n \|b_n\|_\omega +1 $ will do.\qed
\end{point}
\end{point}
\end{point}
\begin{point}{110}[luws]{Corollary}%
For a linear map~$f\colon \scrA\to \C$
on a von Neumann algebra~$\scrA$ the following are equivalent.
\begin{enumerate}
\item
$f$ is ultrastrongly continuous;
\item
$f$ is ultraweakly continuous;
\item
$f\equiv f_0+if_1-f_2-if_3$
for some~np-maps $f_0,\dotsc,f_3\colon \scrA\to\C$;
\item
``$f$ is bounded on some $\|\,\cdot\,\|_\omega$-ball,''
that is,
\begin{equation*}
	\sup\{\ \left|f(a)\right|\colon a\in \scrA\colon 
		\|a\|_\omega\leq \delta\ \}\ <\ \infty
\end{equation*}
for some $\delta>0$ and  np-map $\omega\colon \scrA\to\C$;
\item
$\left|f(a)\right|\leq \|a\|_\omega$
for all~$a\in \scrA$, for some np-map $\omega\colon \scrA\to\C$.
\end{enumerate}
\spacingfix%
\end{point}%
\end{parsec}%
%
%
\begin{parsec}{730}%
\begin{point}{10}%
We'll show that the ultrastrong and ultraweak closure
of a convex set agree. 
For this we need the following proto-Hahn--Banach separation theorem,
which concerns the following notion of openness.
\end{point}
\begin{point}{20}{Definition}%
A subset~$A$ of a real vector space~$V$ 
is called \Define{radially open}%
\index{radially open set}
if for all~$a\in A$ and~$v\in V$
there is $t\in (0,\infty)$
with $a+sv\in A$ for all~$s\in [0,t)$.
\end{point}
\begin{point}{30}{Exercise}%
Let~$V$ be a vector space.
\begin{enumerate}
\item
Show that the radially open subsets of~$V$ form a topology.
\item
Show that with respect to this topology,
scalar multiplication and translations $x\mapsto x+a$
by a fixed vector~$a\in V$ are continuous.
\item
Show that 
the subset of~$\R^2$ depicted below in blue
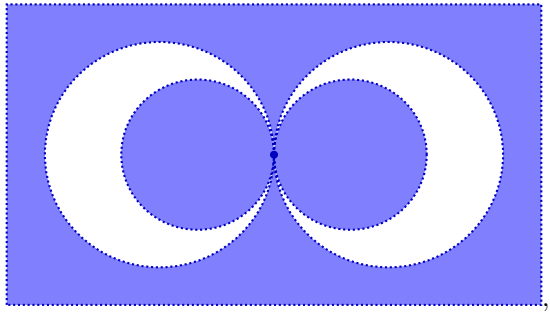
\begin{center}
\begin{tikzpicture}
    \fill[style={fill=lightblue,draw=darkblue,densely dotted,thick}, even odd rule]
    (-3.5cm,-2.0cm) rectangle (3.5cm,2.0cm)
    (-1.5cm,0) circle (1.5cm)
    (1.5cm,0) circle (1.5cm)
    (-1cm,0) circle (1cm)
    (1cm,0) circle (1cm)
    (0,0) node[circle,inner sep=0pt,minimum size=3pt,fill=darkblue]{};
\end{tikzpicture},
\end{center}
including the point in the middle
but not the dashed borders,
is radially open, but not open in the usual topology on~$\R^2$.
\item
Show that addition on~$\R^2$ is not jointly radially continuous.
\item
	Show that nevertheless $\{s\in \R \colon sx+(1-s) y\in A\}$
is open for every radially open~$A\subseteq V$, and $x,y\in V$.
\item
Show that $A+B$ is radially open when~$A,B\subseteq V$ are radially open.

Show that $\{\lambda a\colon a\in A,\lambda>0\}$ is radially open
when~$A$ is radially open.
\end{enumerate}%
\spacingfix%
\end{point}%
\begin{point}{40}[hahn-banach]{Theorem}%
\index{Hahn--Banach's Theorem}
For every radially open
convex subset~$K$ of a real vector space~$V$
with~$0\notin K$
there is a linear map $f\colon V\to\R$
with $f(x)>0$ for all~$x\in K$.
\begin{point}{50}{Proof}%
(Based on Theorem~1.1.2 of~\cite{kr}.)

By Zorn's Lemma we may assume without loss of generality that~$K$ is maximal
among radially open convex subsets of~$V$ that do not contain~$0$.

We also assume that~$K$ is non-empty,
because if~$K=\varnothing$, the result is trivial.

We will show in a moment that~$H:=\{x\in V\colon -x,x\notin K\}$
is a linear subspace and~$V/H$ is one-dimensional.
From this we see that there is a linear map $f\colon V\to\R$
with~$\ker(f)=H$.
Since~$f(K)$ is a convex subset which does not contain~$0$
(because $H\cap K=\varnothing$)
we either have $f(K)\subseteq (0,\infty)$
or $f(K)\subseteq(-\infty,0)$.
Thus, by replacing $f$ by $-f$ if necessary,
we see that there is a linear map $f\colon V\to \R$
with $f(x)>0$ for all~$x\in K$.
\begin{point}{60}{$H$ is a linear subspace}%
Note that~$x\in K,\,\lambda>0\implies \lambda x\in K$,
because the subset 
$\{\lambda x\colon x\in K,\lambda\in(0,\infty)\}\supseteq K$
is radially open, convex, doesn't contain~$0$,
and is thus~$K$ itself.
Furthermore,
$x,y\in K\implies x+y\in K$, because
$x+y=2(\frac{1}{2}x + \frac{1}{2}y)$, and~$K$ is convex.

Let~$\overline{K}$ be the set of all~$x\in V$
with $x+y\in K$ for all~$y\in K$.
Then it is not difficult to check that~$\overline{K}$ is a cone:
 $0\in\overline{K}$,
and
$x\in \overline{K},\lambda\geq 0\implies \lambda x\in \overline{K}$, and
$x,y\in\overline{K}\implies x+y\in \overline{K}$.

We claim that~$x\in \overline{K}$ iff $-x\notin K$.
Indeed, if~$x\in\overline{K}$, then $-x\notin K$, because otherwise
$-x\in K$ and so
$0=x+(-x)\in K$, which is absurd.
For the other direction, suppose that $-x\notin K$.
Then $x+y\in K$ for all~$y\in K$,
because
$\{\lambda x+y\colon y\in K,\lambda\geq0\}\supseteq K$
is radially open, convex, doesn't contain~$0$,
and is thus~$K$.

It follows that~$H=\overline{K}\cap -\overline{K}$.
Since~$\overline{K}$ is a cone, $-\overline{K}$ is a cone,
and thus~$H$ is a cone.  But then~$-H=H$ is a cone too,
and thus~$H$ is a linear subspace.
\end{point}
\begin{point}{70}{$V/H$ is one-dimensional}%
Note that~$H\neq V$, because~$K\cap H=\varnothing$
and~$K$ is (assumed to be) non-empty.
So to show that~$V/H$ is one-dimensional,
it suffices to show that
any~$x,y\in V$ 
 are linearly dependent in~$V/H$.
We may assume that~$x\in K$ and~$y\in -K$.
It suffices to find~$s\in [0,1]$ with $0=sx+s^\perp y$.
The trick is to consider the sets
 $S_0 = \{s\in [0,1]\colon sx+s^\perp y \in -K\}$
and~$S_1 = \{s\in [0,1]\colon sx+s^\perp y \in K\}$,
which are open (because~$K$ and~$-K$ are radially open),
non-empty (because $0\in S_0$ and~$1\in S_1$),
and therefore cannot cover~$[0,1]$
(because~$[0,1]$ is connected).
So there must be~$s\in (0,1)$ 
such that
$sx+s^\perp y $ 
is neither in~$K$
nor in~$-K$,
and thus $sx+s^\perp y \in H$ (by definition of~$H$).
Whence~$x$ and~$y$ are linearly dependent in~$V/H$
(since~$s\neq 0$).\qed
\end{point}
\end{point}
\end{point}
\begin{point}{80}[ultraclosed]{Exercise}%
\index{ultraweak and ultrastrong!convex $\sim$ly closed subset}
We will use~\sref{hahn-banach}
to prove that 
an ultrastrongly closed convex subset~$K$ of a 
von Neumann algebra~$\scrA$
is ultraweakly closed as well.

Let us first simplify the problem a bit.
If~$K$ is empty, the result is trivial,
so we may as well assume that $K\neq \varnothing$.
Note that we must show that no net in~$K$ converges ultraweakly
to any element~$a_0\in\scrA$ outside~$K$,
but by considering~$K-a_0$ instead of~$K$,
we see that it suffices to show that
no net in~$K$ converges ultraweakly to~$0$
under the assumption that~$0\notin K$.
To this end we'll find an ultraweakly continuous linear map
$g\colon \scrA\to \C$ and~$\delta>0$ 
with~$\Real{g(k)}\geq \delta$ for all~$k\in K$---if
a net $(k_\alpha)_\alpha$ in~$K$ were to converge ultraweakly to~$0$,
then~$\Real{g(k_\alpha)}$ would converge to~$0$ as well,
which is impossible.
\begin{enumerate}
\item
	Show that 
	there is an np-map~$\omega \colon\scrA\to\C$
	and~$\varepsilon>0$ 
	with $\|k\|_\omega \geq \varepsilon$ for all~$k\in K$.
	(Hint: use that~$K$ is ultrastrongly closed).
\item
	Show that~$B:= \{ b\in \scrA\colon \|b\|_\omega < \varepsilon\}$
	 is convex, radially open, $B\cap K=\varnothing$.

	Show that $B-K$ is convex, radially open, and $0\notin B-K$.
\item
	Use~\sref{hahn-banach} to show that
	 there is an $\R$-linear map $f\colon \scrA\to \R$
	with $f(b)<f(k)$ for all~$b\in B$ and~$k\in K$.
	Show that~$f$
	can be extended to a $\C$-linear map
	$g \colon \scrA\to \C$
	by $g(a)= f(a)-if(ia)$ for all~$a\in\scrA$.
\item
	Show that $\left| f(b)\right| \leq f(k)$ 
	and $\left|g(b)\right|\leq 2f(k)$
	for all~$b\in B$ and $k\in K$.\\
	(Hint: $b\in B\implies -b\in B$.)

	Conclude that~$g$ is ultraweakly continuous
	(using~\sref{luws} and $K\neq \varnothing$).
\item
	It remains to be shown that
	there is $\delta>0$ with $f(k)\equiv \Real{g(k)}\geq \delta$
	for all~$k\in K$.
	Show that in fact there is $b_0\in B$
	with $f(b_0) >0$,
	and  that $f(k)\geq f(b_0)>0 $ for all~$k\in K$.
\end{enumerate}
\spacingfix
\end{point}%
\end{parsec}%
\subsection{Kaplansky's Density Theorem}
\begin{parsec}{740}%
\begin{point}{10}[proto-kaplansky]{Proposition}%
Let~$\scrA$ be a von Neumann algebra,
and let~$f\colon \R\to\R$ be a continuous map 
with $f(t)=O(t)$,
that is,
there are~$n\in \N$ and~$b\in [0,\infty)$
such that $\left|f(t)\right|\leq b\left|t\right|$
for all~$t\in \R$ with~$\left|t\right| \geq n$.

Then the map~$a\mapsto f(a),\,\sa{\scrA}\to\sa{\scrA}$,
see~\sref{functional-calculus},
is ultrastrongly continuous.
\begin{point}{20}{Proof}%
(An adaptation of Lemma~44.2 from~\cite{conway2000}.)

Let~$S$ denote the set of all continuous $g\colon \R\to\R$
such that $a\mapsto g(a),\,\sa{\scrA}\to\sa{\scrA}$
is ultrastrongly continuous.
We must show that~$f\in S$.

Let us first make some general observations.
The identity map $t\mapsto t$ is in~$S$,
any constant function is in~$S$,
and~$S$ is closed under addition,
and scalar multiplication.
In particular, any affine transformation ($t\mapsto at+b$)
is in~$S$.
Moreover, we have~$g\circ h\in S$ when $g,h\in S$,
and also~$gh\in S$
provided that~$g$ is bounded.
Finally, $S$ is closed with respect to uniform convergence.

Now,
as $f(t)=f(t)\,\smash{\frac{1}{1+t^2}\,+\, f(t)\,\frac{t^2}{1+t^2}}$
 one can see from the remarks above
that it suffices
to show
that~$t\mapsto f(t)\,\smash{\frac{1}{1+t^2}}$ is in~$S$
--- here we use that $t\mapsto f(t) \,\smash{\frac{t}{1+t^2}}$ is bounded.
In other words,
we may assume without loss of generality,
that~$f$ vanishes at infinity, i.e.~$\lim_{t\to \infty}f(t)=0$.

Suppose for the moment
that there is $e\in S$, $e\neq 0$,
which vanishes at infinity.
Let~$a,b\in \R$.
Then $e_{a,b}\colon \R\to\R, t\mapsto e(at+b)$
--- an affine transformation followed by~$e$ ---
is also in~$S$,
vanishes at infinity,
and can be extended to a continuous real-valued
function on the one-point compactification $\R\cup \{ \infty\}$
of~$\R$
(by defining $e_{a,b}(\infty)=0$).
It is easy to see that the $C^*$-subalgebra
of~$C(\R\cup\{\infty\})$
generated by these extended~$e_{a,b}$'s 
separates the points of~$\R\cup\{\infty\}$,
and is thus~$C(\R\cup\{\infty\})$ itself
    by the Stone--Weierstra\ss{} theorem (see~\sref{stone-weierstrass}).
Since~$f$ vanishes at infinity,
$f$ can be extended to an element of $C(\R\cup\{\infty\})$,
 and can thus be obtained
(by taking real parts if necessary)
from the extended $e_{a,b}$'s and real constants 
via uniform limits, addition and (real scalar)
multiplication. 
Since~$S$ contains the $e_{a,b}$'s and constants
and is closed under these operations (acting on bounded functions),
we see that~$f\in S$.

To complete the proof,
we show that such~$e$ indeed exists.
Let $e,s\colon \R\to\R$ 
be given by $e(t)=ts(t)$ and $s(t)=\smash{\frac{1}{1+t^2}}$.
Clearly~$e$ and~$s$
are continuous
and
vanish at infinity.
To see that~$e$ is ultrastrongly continuous,
let $(b_\alpha)_\alpha$ be a net of self-adjoint elements of~$\scrA$
which converges ultrastrongly to $a\in\sa{\scrA}$,
and let~$\omega\colon \scrA\to \C$
be an npu-map. 
Unfolding the definitions
of~$e$ and~$s$ yields
the following equality.
\begin{equation*}
e(b_\alpha)-e(a) \ =\ s(b_\alpha)\,(b_\alpha-a)\,s(a)
\,-\, e(b_\alpha)\,(b_\alpha-a)\,e(a).
\end{equation*}
Since $\|s(b_\alpha)\|\leq 1$,
we have $\|s(b_\alpha)(b_\alpha-a)s(a)\|_\omega 
\leq  \,\|(b_\alpha-a)s(a)\|_\omega
\equiv \|b_\alpha-a\|_{s(a)*\omega}$.
Similarly, since~$\|e(b_\alpha)\|\leq 1$,
we get
\begin{alignat*}{3}
\|e(b_\alpha)-e(a)\|_\omega
\ \leq\ \|b_\alpha-a\|_{s(a)*\omega}\,+\,\|b_\alpha-a\|_{e(a)*\omega}.
\end{alignat*}
Thus~$e(b_\alpha)$ converges ultrastrongly to~$e(a)$,
and so~$e$ is ultrastrongly continuous.\qed
\end{point}
\end{point}

\begin{point}{30}[abs-us-cont]{Corollary}%
Given a von Neumann algebra~$\scrA$
the map $a\mapsto\left|a\right|\colon\,\Real{\scrA}\to\Real{\scrA}$
is ultrastrongly continuous.
\end{point}
\begin{point}{40}[kaplansky]{Kaplansky's Density Theorem}%
\index{Kaplansky's Density Theorem}
Let~$b$ be an element of a von Neumann algebra~$\scrB$
which is the ultrastrong limit of a net
of elements
from a $C^*$-subalgebra $\scrA$ of~$\scrB$.
Then
$b$ is the ultrastrong limit of a net~$(a_\alpha)_\alpha$
in~$\scrA$ with~$\|a_\alpha\|\leq\|b\|$ for all~$\alpha$.
Moreover,
\begin{enumerate}
\item
if~$b$ is self-adjoint,
then the~$a_\alpha$ can be chosen to be self-adjoint as well;
\item
if~$b$ is positive,
then the~$a_\alpha$ can be chosen to be positive as well, and
\item
if~$b$ is an effect,
then the~$a_\alpha$ can be chosen to be effects as well.
\end{enumerate}
\spacingfix%
\begin{point}{50}{Proof}%
Let~$(a_\alpha)_\alpha$
be a net in~$\scrA$ that converges ultrastrongly
to~$b$.

Assume for the moment that~$b$ is self-adjoint.
Then~$\Real{(a_\alpha)}$
converges ultraweakly (but perhaps not ultrastrongly)
to~$\Real{b}=b$
as~$\alpha\to\infty$,
and so~$b$ is in the ultraweak
closure of the convex set~$\Real{\scrA}$.
Since the ultraweak and ultrastrong closure
of convex subsets of~$\scrA$
coincide (by~\sref{ultraclosed}),
we see that~$b$ is also the ultrastrong limit
of some net $(a_\alpha')_\alpha$ in~$\Real{\scrA}$.
Since the map~$-\|b\|\vee(\,\cdot\,)\wedge \|b\|\colon
\Real{\scrB}\to\Real{\scrB}$
is ultrastrongly continuous by~\sref{proto-kaplansky}
we see that 
$a_\alpha'' := -\|b\|\vee a_\alpha'\wedge \|b\|$
gives a net $(a_\alpha'')_\alpha$
in~$[-\|b\|,\|b\|]_\scrA$
that
converges ultrastrongly
to~$b$.

If we assume in addition
that~$b$ is positive,
then~$a_\alpha''' := (a_\alpha'')_+$
gives a net~$(a_\alpha''')_\alpha$
in~$[0,\|b\|]_\scrA$
that converges ultrastrongly to $b_+=b$,
because the map $(\,\cdot\,)_+\colon \Real{\scrB}\to\Real{\scrB}$
is ultrastrongly continuous by~\sref{proto-kaplansky}.
Note that if~$b$ is an effect,
then so are the $a_\alpha'''$.

This takes care of all the special cases.
The general case
in which~$b$ is an arbitrary element of~$\scrB$ requires a trick:
since the element
$B:=\smash{\left(\begin{smallmatrix} 0 & b \\
b^* & 0 \end{smallmatrix}\right)}$
of the von Neumann algebra~$M_2(\scrB)$
is self-adjoint,
and the ultrastrong limit
of the net
$\smash{\bigl(\begin{smallmatrix} 0 & a_\alpha \\
a_\alpha^* & 0 \end{smallmatrix}\bigr)}$
from the  $C^*$-subalgebra $M_2(\scrA)$ of~$M_2(\scrB)$,
there is, as we've just seen, a net~$(A_\alpha)_\alpha$
in~$M_2(\scrA)$
that converges ultrastrongly to~$B$
with~$\|A_\alpha\|\leq \|B\|\equiv \|b\|$ for all~$\alpha$.
Since the upper-right entries
$ (A_\alpha)_{12}$ will then converge ultrastrongly 
to~$B_{12}\equiv b$ as~$\alpha\to\infty$,
and~$\|(A_\alpha)_{12}\|\leq \|A_\alpha\|\leq \|b\|$
for all~$\alpha$, we are done.\qed
\end{point}
\end{point}
\begin{point}{60}[dense-subalgebra]{Corollary}%
Given~$\varepsilon>0$  and an ultraweakly dense $*$-subalgebra~$\scrS$
of a von Neumann algebra~$\scrA$
each element~$a$ of~$\scrA$
is the ultrastrong limit of a net~$(s_\alpha)_\alpha$
from~$\scrS$
with $\|s_\alpha\|\leq \|a\|(1+\varepsilon)$
for all~$\alpha$.
\begin{point}{70}{Proof}%
As the norm closure~$\scrC$ of~$\scrS$ in~$\scrA$
is an ultraweakly
(and thus by~\sref{ultraclosed}
ultrastrongly) dense $C^*$-subalgebra of~$\scrA$,
the element $a$ of~$\scrA$ is by~\sref{kaplansky}
the ultrastrong limit
of net~$(c_\alpha)_{\alpha\in D}$
in~$\scrC$ 
with $\|c_\alpha\|\leq \|a\|$ for all~$\alpha$.
Each element~$c_\alpha$ 
is in its turn 
the norm (and thus ultrastrong)
limit of a sequence~$s_{\alpha1},\,s_{\alpha2},\,\dotsc$
in~$\scrS$,
and if we choose the~$s_{\alpha n}$
such that $\|c_\alpha - s_{\alpha n}\|\leq 2^{-n}$,
then~$s_{\alpha n}$
converge ultrastrongly to~$b$
as~$D\times \N\ni (\alpha,n)\to\infty$.
Finally, since~$\lim_n\|s_{\alpha n}\| = \|c_\alpha \|\leq \|c\| 
\leq (1+\varepsilon)\|c\|$
we have~$\|s_{\alpha n}\|\leq (1+\varepsilon)\|c\|$
for sufficiently large~$n$,
and thus for all~$n$
if we replace~$(s_{\alpha n})_n$ by the appropriate  subsequence.\qed
\end{point}
\end{point}
\end{parsec}
\subsection{Closedness of Subalgebras}
\begin{parsec}{750}%
\begin{point}{10}%
Recall that according to our definition (\sref{von-neumann-examples})
a von Neumann subalgebra~$\scrB$ 
of a von Neumann algebra~$\scrA$
is a $C^*$-subalgebra of~$\scrA$
which is closed under suprema
of bounded directed sets of self-adjoint elements.
We will show that such~$\scrB$ is ultrastrongly closed in~$\scrA$.
\end{point}
\begin{point}{20}[sequence-separation-lemma]{Lemma}%
Let~$\scrB$ be a von Neumann subalgebra
of a von Neumann algebra~$\scrA$.
Let~$\omega_0,\omega_1\colon \scrA\to\C$
be npu-maps,
which are separated
by a net $(b_\alpha)_\alpha$
of effects of~$\scrB$ 
in the sense that~$\lim_\alpha \omega_0(b_\alpha)=0$
and $\lim_\alpha \omega_1(b_\alpha^\perp)=0$.
Then~$\omega_0$ and~$\omega_1$ are separated by a 
projection~$q$ of~$\scrB$ 
in the sense that~$\omega_0(q)= 0 = \omega_1(q^\perp)$.
\begin{point}{30}{Proof}%
(Based on Lemma~45.3 and Theorem~45.6 of~\cite{conway2000}.)

Note that it suffices to find an effect~$a$ in~$\scrB$
with $\omega_0(a) = 0 = \omega_1(a^\perp)$,
because then~$\omega_0(\ceil{a}) = 0 = \omega_1(\ceil{a}^\perp)$
by~\sref{ceil-functionals-lemma}
and~$\ceil{a}\in\scrB$.

Note that we can find a subsequence~$(b_n)_n$ of~$(b_\alpha)_\alpha$
such that $\omega_0( b_n ) 
	\leq n^{-1}2^{-n}$
and $\omega_1(b_n^\perp)\leq n^{-1}$
for all~$n$.
For~$n < m$, define
\begin{equation*}
\textstyle
a_{nm}\ = \ (1+\sum_{k=n}^m kb_k)^{-1} \,\sum_{k=n}^m kb_k.
\end{equation*}
Since we have seen in~\sref{astara-pos-basic-consequences}
that the map~$d\mapsto (1+d)^{-1}d$ is order preserving
(on~$\pos{\scrB}$),
we have $0\leq a_{nm}\leq \frac{1}{2}$
and we get the formation
\begin{equation*}
\renewcommand{\labelstyle}{\textstyle}
\xymatrix{
a_{12}\ar@{}|-{\leq}[r] & 
a_{13}\ar@{}|-{\leq}[r] & 
a_{14}\ar@{}|-{\leq}[r]  & 
\dotsb & \ar@{}|-{\leq}[r]& a_1 \\
&
a_{23} \ar@{}|-{\leq}[r] \ar@{}[u]|-{\uleq} &
a_{24} \ar@{}|-{\leq}[r]  \ar@{}[u]|-{\uleq} &
\dotsb &\ar@{}|-{\leq}[r]& a_2  \ar@{}[u]|-{\uleq} \\
&
&
a_{34} \ar@{}|-{\leq}[r]  \ar@{}[u]|-{\uleq} & 
\dotsb & \ar@{}|-{\leq}[r] & a_3  \ar@{}[u]|-{\uleq} 
\\
&&&\ddots&&\vdots \ar@{}[u]|-{\uleq}& \\
&&&&&a\ar@{}[u]|-{\uleq}
},
\end{equation*}
where~$a_n:= \bigvee_{m\geq n} a_{nm}$
and~$a := \bigwedge_n a_n$.
We'll prove that~$\omega_0(a)=0=\omega_1(a^\perp)$.
\begin{point}{40}{$\omega_0(a)=0$}%
	Since~$\omega_0(b_n)\leq n^{-1}2^{-n}$ 
and $a_{nm}\leq \sum_{k=n}^m k b_k$,
we get~$\omega_0(a_{nm})\leq \sum_{k=n}^m k\omega_0(b_k) \leq 2^{1-n}$,
and so $\omega_0(a)=\bigwedge_n\bigvee_{m\geq n} \omega_0(a_{nm})
\leq \bigwedge_n 2^{1-n} = 0$.
\end{point}
\begin{point}{50}{$\omega_1(a^\perp )=0$}%
Let~$m> n$ be given.
Since~$\sum_{k=n}^m kb_k \geq mb_m$
and $d\mapsto (1+d)^{-1}d$
is monotone on~$\pos{\scrB}$
we get~$a_{nm} \geq (1+mb_m)^{-1} mb_m$,
and so~$a_{nm}^\perp \leq (1+mb_m)^{-1}$.

Observe that for a real number $t\in[0,1]$,
we have $tt^\perp \geq 0$,
and so $(1+mt)(1+mt^\perp) = 1+m+m^2tt^\perp \geq 1+m$.
This yields the inequality $(1+mt)^{-1}\leq (1+m)^{-1}(1+mt^\perp)$
for real numbers~$t\in[0,1]$.
The corresponding inequality for effects of a $C^*$-algebra
(obtained via Gelfand's representation theorem, \sref{gelfand})
gives us $\omega_1(a_{nm}^\perp)\leq \omega_1((1+mb_m)^{-1})
\leq (1+m)^{-1}(1+m\omega_1(b_m^\perp))\leq \frac{2}{1+m}$,
where we have used that~$\omega_1(b_m^\perp)\leq \frac{1}{m}$.
Hence~$\omega_1(a_n^\perp)=\bigwedge_{m\geq n} \omega_1(a_{nm}^\perp)
\leq \bigwedge_{m\geq n} \frac{2}{1+m}=0$ for all~$n$,
and so~$\omega_1(a^\perp)=\bigvee_n\omega_1(a_n^\perp)=0$.\qed
\end{point}
\end{point}
\end{point}
\begin{point}{60}[kadisons-lemma]{Lemma}%
Let~$\scrB$ be a von Neumann subalgebra
of a von Neumann algebra~$\scrA$.
Let~$p$ be a projection of~$\scrA$,
which is the ultrastrong limit of a net in~$\scrB$.

For all npu-maps $\omega_0,\omega_1\colon \scrA\to\C$
with~$\omega_0(p)=0= \omega_1(p^\perp)$
there is a projection~$q$ of~$\scrB$
with~$\omega_0(q)=0=\omega_1(q^\perp)$.
\begin{point}{70}{Proof}%
Let~$(b_\alpha)_\alpha$ be a net in~$\scrB$
which converges ultrastrongly to~$p$.
We may assume that all~$b_\alpha$
are effects
by Kaplansky's density theorem (\sref{kaplansky}).
Note that~$(\omega_0(b_\alpha))_\alpha$ converges to~$\omega_0(p)\equiv 0$,
and $(\omega_1(b^\perp_\alpha))_\alpha$ converges 
to~$\omega_1(p^\perp)\equiv 0$.
Now apply~\sref{sequence-separation-lemma}.\qed
\end{point}
\end{point}
\begin{point}{80}[vnsac]{Theorem}%
\index{von Neumann subalgebra!is ultraweakly closed}
A von Neumann subalgebra~$\scrB$ of a von Neumann algebra~$\scrA$
is ultrastrongly and ultraweakly closed.
\begin{point}{90}{Proof}%
It suffices to show that~$\scrB$ is ultrastrongly closed,
because then, by~\sref{ultraclosed}, $\scrB$ will be ultraweakly closed
as well.

Let~$p$ be a projection of~$\scrA$ which is the ultrastrong limit
of a net from~$\scrB$. It suffices to show that~$p\in\scrB$,
because the ultrastrong closure of~$\scrB$
being a von Neumann subalgebra of~$\scrA$
is generated by its projections, see~\sref{projections-norm-dense}.
Note that given an np-map $\omega\colon \scrA\to\C$,
the carrier~$\ceil{\omega}$ of~$\omega$
need not be equal to the carrier
of~$\omega$ restricted to~$\scrB$,
which we'll therefore denote by~$\ceil{\omega}_\scrB$;
but we do have $\ceil{\omega}\leq \ceil{\omega}_\scrB$.
Then by~\sref{ultracyclic-basic}
\begin{equation}
\label{ultracyclic-proof}
\textstyle
\bigvee_{\omega_1}
\ceil{\omega_1}_\scrB \ \geq\ 
\bigvee_{\omega_1}
\ceil{\omega_1}\ =\ p\ =\ 
\bigwedge_{\omega_0} \ceil{\omega_0}^\perp
\ \geq\ \bigwedge_{\omega_0} \ceil{\omega_0}_{\scrB}^\perp,
\end{equation}
where~$\omega_0$ ranges over np-maps $\omega_0\colon \scrA\to\C$
with~$\omega_0(p)=0$,
 and $\omega_1$ ranges over
np-maps $\omega_1\colon \scrA\to\C$ with~$\omega_1(p^\perp)=0$.
Since for such~$\omega_0$ and~$\omega_1$
there is 
by~\sref{kadisons-lemma}
a projection~$q$ in~$\scrB$
with $\omega_0(q)=0=\omega_1(q^\perp)$,
we get 
$\ceil{\omega_1}_\scrB \leq q \leq 
\ceil{\omega_0}_\scrB^\perp$,
and so~$\bigvee_{\omega_1} \ceil{\omega_1}_\scrB
\leq 
\bigwedge_{\omega_0} \ceil{\omega_0}_\scrB^\perp$.
It follows that the inequalities in~\eqref{ultracyclic-proof}
are in fact equalities,
and so~$p=\bigvee_{\omega_1}\ceil{\omega_1}_\scrB \in \scrB$.\qed
\end{point}
\end{point}
\end{parsec}
\subsection{Completeness}
\begin{parsec}{760}%
\begin{point}{10}[bh-us-complete]{Proposition}%
The von Neumann algebra~$\scrB(\scrH)$
of bounded operators on a Hilbert space~$\scrH$
is ultrastrongly complete.
\begin{point}{20}{Proof}
Let~$(T_\alpha)_\alpha$ be an ultrastrongly Cauchy net
in~$\scrB(\scrH)$
(which must be shown to converge ultrastrongly to
some operator~$T$ in~$\scrB(\scrH)$).

Note that given~$x\in \scrH$,
the net~$(T_\alpha x)_\alpha$ in~$\scrH$
is norm Cauchy,
because $\|(T_\alpha-T_\beta) x\|
= \| T_\alpha-T_\beta \|_{\left<x,(\,\cdot\,)x\right>}$
vanishes for sufficiently large~$\alpha,\beta$,
and so we may define~$Tx :=\lim_\alpha T_\alpha x$,
giving a map~$T\colon \scrH\to\scrH$.

It is clear that~$T$ will be linear,
but the question is whether~$T$ is bounded,
and whether in that 
case~$(T_\alpha)_\alpha$ converges ultrastrongly to~$T$.

Suppose towards a contradiction that~$T$ is not bounded.
Then we can find~$x_1,x_2,\dotsc\in\scrH$
with $\|x_n\|^2\leq 2^{-n}$
and $\|Tx_n\|^2\geq 1$ for all~$n$.
Since~$\omega:=\sum_n \left<x_n,(\,\cdot\,)x_n\right>\colon 
\scrB(\scrH)\to\C$
is an np-map by~\sref{bh-functional-lemma}, 
it follows that~$\|T_\alpha\|_\omega^2\equiv \sum_{n=1}^\infty 
\|T_\alpha x_n\|^2$
converges to some positive number~$R$.
Since any partial sum $\sum_{n=1}^N \|T_\alpha x_n\|^2
\leq \|T_\alpha\|_\omega^2$
converges to~$\sum_{n=1}^N \|T x_n\|^2\geq N$,
we must conclude that~$R\geq N$,
for all natural numbers~$N$,
which is absurd.
Hence~$T$ is bounded.

It remains to be shown that~$(T_\alpha)_\alpha$
converges ultrastrongly to~$T$.
So let $\omega\colon \scrB(\scrH)\to\C$ be an arbitrary
np-map,
being of the form~$\omega\equiv \sum_n\left<x_n,(\,\cdot\,)x_n\right>$
for some $x_1,x_2,\dotsc\in\scrH$ with $\sum_n \|x_n\|^2 <\infty$
by~\sref{bh-np}.
We must show that~$\|T-T_\alpha\|_\omega
\equiv (\sum_n \|(T-T_\alpha)x_n\|^2)^{\nicefrac{1}{2}}$ 
converges to~$0$ as $\alpha\to 0$.

Let~$\varepsilon>0$ be given,
and pick~$\alpha_0$
such that $\|T_\alpha-T_\beta\|_\omega \leq 
\smash{\frac{1}{2\sqrt{2}}}\,\varepsilon$
for all~$\alpha,\beta\geq \alpha_0$
--- this is possible  because~$(T_\alpha)_\alpha$ is ultrastrongly Cauchy.
We claim that $\|T-T_\alpha\|_\omega \leq \varepsilon$
for any~$\alpha\geq \alpha_0$.
Since for such~$\alpha$
the sum
\begin{equation*}
	\sum_{n=1}^\infty \|(T-T_\alpha)x_n\|^2 
	\ = \ 
	\sum_{n=1}^{N-1} \|(T-T_\alpha)x_n\|^2
	\ +\ 
	\sum_{n=N}^\infty
	\|(T-T_\alpha)x_n\|^2
\end{equation*}
converges (to~$\|T-T_\alpha\|_\omega^2$),
we can find~$N$ such that the second term in the bound above
is below~$\frac{1}{2}\varepsilon^2$.
The first term will also be below~$\frac{1}{2}\varepsilon^2$,
because
\begin{equation*}
	\bigl(\sum_{n=1}^{N-1} \|(T - T_\alpha)x_n\|^2\,\bigr)^{\nicefrac{1}{2}}
	\ \leq\ 
\bigl(\,
\sum_{n=1}^{N-1} \|(T-T_{\beta})x_n\|^2
\,\bigr)^{\nicefrac{1}{2}}
\ +\ 
\bigl(\,
\sum_{n=1}^{N-1} \|(T_\beta-T_\alpha)x_n\|^2
\,\bigr)^{\nicefrac{1}{2}}
\end{equation*}
for any~$\beta$,
and in particular for~$\beta$
large enough that the first term on the right-hand side above
is below~$\smash{\frac{1}{2\sqrt{2}}}\,\varepsilon$.
If we choose $\beta\geq \alpha_0$
the second term will be below~$\smash{\frac{1}{2\sqrt{2}}}\,\varepsilon$ too,
and we get $\|T-T_\alpha\|_\omega^2 \leq \frac{1}{2}\varepsilon^2
+ (\smash{\frac{1}{2\sqrt{2}}\,\varepsilon 
+ \frac{1}{2\sqrt{2}}\,\varepsilon})^2
\equiv \varepsilon^2$ 
all in all.
(This reasoning is very similar to that in~\sref{hilb-sum}.)

Hence~$\scrB(\scrH)$ is ultrastrongly complete.\qed
\end{point}
\end{point}
\begin{point}{30}[bh-bounded-uw-complete]{Proposition}%
The von Neumann algebra~$\scrB(\scrH)$
of bounded operators on a Hilbert space~$\scrH$
is bounded ultraweakly complete.
\begin{point}{40}{Proof}%
Let~$(T_\alpha)_\alpha$ be a norm-bounded ultraweakly Cauchy net
in~$\scrB(\scrH)$.
We must show that~$(T_\alpha)_\alpha$
converges ultraweakly
to some bounded operator~$T$ on~$\scrH$.

Note that given $x,y\in\scrH$
the net $(\,\left<x,T_\alpha y\right>\,)_\alpha$
is Cauchy
(because $\left<x,(\,\cdot\,)y\right>
\equiv \frac{1}{4}\sum_{k=0}^3 i^k\left<i^kx+y,(\,\cdot\,)(i^kx+y)\right>$
is ultraweakly continuous),
and so we may define $[x,y] = \lim_\alpha \left<x,T_\alpha y\right>$.
The resulting `form'
$[\,\cdot\,,\,\cdot\,]\colon \scrH\times\scrH\to \C$
(see~\sref{chilb-form})
is bounded, 
because~$\left\|[x,y]\right\| \leq (\sup_\alpha \|T_\alpha\|)\|x\|\|y\|$
for all~$x,y\in\scrH$
and $\sup_\alpha\|T_\alpha\|<\infty$
since $(T_\alpha)_\alpha$
is norm bounded.
By~\sref{chilb-form-representation},
there is a unique bounded operator~$T$
with $\left<x,Ty\right>=[x,y]$
for all~$x,y\in\scrH$.

By definition of~$T$ it is clear 
that~$\lim_\alpha \left<x,(T-T_\alpha)x\right>=0$
for any~$x\in\scrH$,
but it is not yet clear that~$(T_\alpha)_\alpha$ converges ultraweakly to~$T$.
For this we must show that $\lim_\alpha \omega(T-T_\alpha)=0$
for any np-map~$\omega\colon \scrB(\scrH)\to \C$.
By~\sref{bh-np},
we know that such~$\omega$ is of the form
$\omega=\sum_n \left<x_n,(\,\cdot\,)x_n\right>$
for some $x_1,x_2,\dotsc\in\scrH$ with $\sum_n\|x_n\|^2<\infty$.
Now, given~$N$ and~$\alpha$ we easily obtain the 
following bound.
\begin{equation*}
	|\,\omega(T-T_\alpha)| \ \leq\  
	\sum_{n=1}^{N-1} \left|\left<x_n(T-T_\alpha),x_n\right>\right|
\ +\ \bigl(\,\|T\|+\sup_\alpha\|T_\alpha\|\,\bigr)\,\sum_{n=N}^\infty \|x_n\|^2 
\end{equation*}
Since the first term of this bound converges to~$0$ as~$\alpha\to\infty$,
we get, for all~$N$,
\begin{equation*}
	\limsup_\alpha |\,\omega(T-T_\alpha)| \ \leq\ 
	\bigl(\,\|T\|+\sup_\alpha\|T_\alpha\|\,\bigr)\,\sum_{n=N}^\infty \|x_n\|^2.
\end{equation*}
Since the tail $\sum_{n=N}^\infty\|x_n\|^2$
converges to~$0$ as~$N\to \infty$,
$\limsup_\alpha \left|\omega(T-T_\alpha)\right|=0$.
Hence $\omega(T)=\lim_\alpha \omega(T_\alpha)$,
and so~$(T_\alpha)_\alpha$ converges ultraweakly to~$T$.\qed
\end{point}
\end{point}
\end{parsec}
\begin{parsec}{770}%
\begin{point}{10}[vn-complete]{Theorem}%
	\index{von Neumann algebra!is ultrastrongly complete}%
	\index{von Neumann algebra!is bounded ultraweakly complete}%
	\index{ultraweak and ultrastrong!completeness}
A von Neumann algebra~$\scrA$ is ultrastrongly complete
and bounded ultraweakly complete.
\begin{point}{20}{Proof}%
Let~$\Omega$ be the set of all np-functionals
on~$\scrA$.
Recall from~\sref{ngns-proof}
that $\varrho_\Omega$
gives an nmiu-isomorphism
onto the 
von Neumann algebra~$\scrR:=\varrho_\Omega(\scrA)$ of operators
on the Hilbert space~$\scrH_\Omega$.
Since~$\scrB(\scrH_\Omega)$ is ultrastrongly complete
(\sref{bh-us-complete}),
and~$\scrR$ is  ultrastrongly closed in~$\scrB(\scrH_\Omega)$
(see~\sref{vnsac}),
we see that $\scrR$
is complete with respect to the ultrastrong
topology of~$\scrB(\scrH_\Omega)$,
but since any np-functional~$\omega\colon \scrR\to\C$
is of the form~$\omega\equiv \left<x,(\,\cdot\,)x\right>$
for some~$x\in\scrH_\Omega$,
and therefore  the ultrastrong topology on~$\scrB(\scrH_\Omega)$
coincides on~$\scrR$ with the ultrastrong topology of~$\scrR$,
we see that~$\scrR$ (and therefore~$\scrA$)
is complete with respect to its own ultrastrong topology.
Since similarly~$\scrB(\scrH_\Omega)$
is bounded ultraweakly complete (\sref{bh-bounded-uw-complete}),
the ultraweak topology on~$\scrB(\scrH_\Omega)$
coincides on~$\scrR$ with the ultraweak topology on~$\scrR$,
and~$\scrR$ is ultraweakly closed
in~$\scrB(\scrH_\Omega)$
(by~\sref{vnsac}),
we see that~$\scrR$ is bounded ultraweakly complete.\qed
\end{point}
\end{point}
\begin{point}{30}[vn-ball-compact]{Theorem}%
The ball $(\scrA)_1$
of a von Neumann algebra~$\scrA$ is ultraweakly compact.
\begin{point}{40}{Proof}%
Writing~$\Omega$ for
the set of npu-maps $\omega\colon \scrA\to\C$,
the map~$\kappa\colon \scrA\to \C^\Omega$
given by~$\kappa(a)=(\omega(a))_\omega$ for all~$a\in\scrA$
is clearly a linear homeomorphism from~$\scrA$ with the  ultraweak topology 
onto~$\kappa(\scrA)\,\subseteq \C^\Omega$ endowed
with the product topology.
Since~$\kappa$ restricts
to an isomorphism of uniform spaces
$(\scrA)_1\to \kappa(\,(\scrA)_1\,)$,
and $(\scrA)_1$ is ultraweakly complete 
(being a norm-bounded ultraweakly closed
subset of the bounded ultraweakly complete space~$\scrA$,
see~\sref{vn-complete}),
we see that $\kappa(\,(\scrA)_1\,)$
is complete,
and thus closed in~$\C^\Omega$.
Now note that~$\kappa(\,(\scrA)_1\,)$ is a closed subset 
of the (by Tychonoff's theorem) compact
space~$((\C)_1)^\Omega$, 
because $\left|\omega(a)\right|\leq 1$ for all~$a\in(\scrA)_1$
and $\omega\in\Omega$.
But then~$\kappa(\,(\scrA)_1\,)$,
being a closed subset of a compact Hausdorff space,
is compact,
and so $(\scrA)_1$ (being homeomorphic to it) is compact too.\qed
\end{point}
\end{point}
\begin{point}{50}[vn-extension]{Proposition}%
Given an ultraweakly dense $*$-subalgebra~$\scrS$
of a von Neumann algebra~$\scrA$,
any ultraweakly continuous and bounded linear map~$f\colon \scrS\to\scrB$
can be extended uniquely
to an ultraweakly continuous map~$g\colon \scrA\to\scrB$.

Moreover, $g$ is bounded,
and in fact,  $\|g\|=\|f\|$.
\begin{point}{60}{Proof}%
As the uniqueness of~$g$ is rather obvious
we concern ourselves only with its existence.
Let~$a\in\scrA$ be given
in order to define~$g(a)$.
Let also~$\varepsilon>0$ be given.
Note that by~\sref{dense-subalgebra} 
there is a net~$(s_\alpha)_\alpha$
in~$\scrS$
that converges ultrastrongly (and so ultraweakly too)
to~$a$
with~$\|s_\alpha \|\leq(1+\varepsilon)\|a\|$
for all~$\alpha$.
Now,
since the net~$(s_\alpha)_\alpha$
is bounded an ultraweakly Cauchy,
and~$f$ is bounded and (uniformly) ultraweakly continuous,
the net
$(f(s_\alpha))_\alpha$
is bounded and ultraweakly Cauchy too,
and thus converges (by~\sref{vn-complete})
to some element
$\uwlim_\alpha f(s_\alpha)$
of~$\scrB$.
\begin{point}{70}%
Of course we'd like to define $g(a):=\uwlim_\alpha
f(s_\alpha)$,
but must first check
that $\uwlim_\alpha f(s_\alpha')
=\uwlim_\alpha f(s_\alpha)$
when~$(s_\alpha')_\alpha$ is a second net with the same properties
as~$(s_\alpha)_\alpha$.
Let us for simplicity's sake
assume that $(s_\alpha')_\alpha$ and~$(s_\alpha)_\alpha$
have the same index set
--- matters can always be arranged this way.
Then as the difference $s_\alpha-s_\alpha'$
converges ultraweakly to~$0$ in~$\scrA$ as~$\alpha\to\infty$,
$\uwlim_\alpha f(s_\alpha-s_\alpha')=0$,
implying that $\uwlim_\alpha f(s_\alpha)
= \uwlim_\alpha f(s_\alpha')$.
\end{point}
\begin{point}{80}%
In this way
we obtain a map $g \colon \scrA\to\scrB$
--- which is clearly linear.
The map~$g$ is also bounded,
because since~$\|s_\alpha\|\leq (1+\varepsilon)\|a\|$
for all~$\alpha$,
where $(s_\alpha)_\alpha$ and~$t$ are as before,
we have $\|f (s_\alpha)\|\leq (1+\varepsilon) \|f\|\|a\|$
for all~$\alpha$,
and so~$\|g(a)\|=\|\uwlim_\alpha f(s_\alpha)\|
\leq (1+\varepsilon) \|f \|\|a\|$.
More precisely, $\|g\|\leq (1+\varepsilon)\|f\|$,
and---as~$\varepsilon>0$ was arbitrary---in 
fact~$\|g\|\leq \|f\|$, and so~$\|g\|=\|f\|$.

That, finally, $g$ is ultraweakly continuous
follows by a standard but abstract argument from the fact
that $f$ is \emph{uniformly} ultraweakly continuous.
We'll give a concrete version of this argument here.
To begin, note that it suffices to show that~$\omega\circ g$
is ultraweakly continuous
at~$0$
where~$\omega\colon \scrB\to\C$
is an np-functional.
Let~$\varepsilon>0$ be given.
Since~$f$ 
is ultraweakly continuous,
and 
thus $\omega\circ f$ is too,
there is~$\delta>0$ and an np-functional~$\nu\colon \scrA\to\C$
such that $\left|\nu(s)\right|\leq \delta \implies 
\left|\omega(f(s))\right|\leq \varepsilon$
for all~$s\in\scrS$.
We claim that $\left|\nu(a)\right|\leq \delta/2 \implies
\left|\omega(g(a))\right|\leq 2\varepsilon$
for all~$a\in\scrA$,
which implies, of course,
that $\omega\circ g$ is ultraweakly continuous on~$0$.
So let~$a\in\scrA$ with $\left|\nu(a)\right|\leq \delta/2$
be given.
Pick (as before) a bounded net~$(s_\alpha)_\alpha$
in~$\scrS$
such that $f(s_\alpha)$
converges to~$a$ as~$\alpha\to\infty$, and observe that,
for all~$\alpha$,
\begin{equation*}
\left|\omega(g (a))\right|
\ \leq\ 
\left|\omega(g(a)-f(s_\alpha))\right|
\,+\,\left|\omega(f(s_\alpha))\right|.
\end{equation*}
The first term on the right-hand side above will 
vanish as~$\alpha\to\infty$ (since 
$g (a)=\uwlim_\alpha f (s_\alpha)$),
and will thus be smaller than~$\varepsilon$ 
for sufficiently large~$\alpha$.
Since~$\lim_\alpha \left|\nu(s_\alpha)\right|=\left|\nu(a)\right|
\leq \delta/2<\delta$
we see 
that for sufficiently large~$\alpha$
we'll have~$\left|\nu(s_\alpha)\right|\leq \delta$ 
and with it~$\left|\omega(f (s))\right|\leq \varepsilon$.
Combined,
we get $\left|\omega(g (a))\right|\leq 2\varepsilon$,
and so~$g$ is ultraweakly continuous.\qed
\end{point}
\end{point}
\end{point}
\end{parsec}
\section{Division}
\label{S:division}
\begin{parsec}{780}%
\begin{point}{10}%
Using the ultrastrong completeness of von Neumann algebras
(see~\sref{vn-complete})
we'll address
the question of division:
given elements~$a$ and~$b$ of a von Neumann algebra~$\scrA$,
when is there an element~$c\in\scrA$
with $a=cb$?
Surely, 
such~$c$ can not always exist,
because
its presence
implies
\begin{equation}
	\label{douglas-ineq}
	a^*a \ \leq \ B\, b^*b,
\end{equation}
where~$B=\|c\|^2$;
but this turns out to be the only restriction:
we'll see in~\sref{douglas} that if~\eqref{douglas-ineq}
holds for some~$B\in [0,\infty)$,
then~$a=cb$ for some
unique~$c\in\scrA$
with~$\ceill{c}\leq \ceilr{b}$,
which we'll denote by~$a/b$.

The main application of this division in our work
is a universal property
for the map $b\mapsto \sqrt{a}b\sqrt{a}\colon
\scrA\to\ceil{a}\!\scrA\!\ceil{a}$
where~$a$ is a positive element of a von Neumann algebra~$\scrA$.
Indeed,
we'll show that
for every np-map $f\colon \scrB\to \scrA$
with $f(1)\leq a$
there is a (unique) np-map $g\colon\scrB\to\ceil{a}\!\scrA\!\ceil{a}$
with~$f(b)=\sqrt{a}g(b)\sqrt{a}$ for all~$b\in\scrB$
---
by taking $g(b)=\sqrt{a}\backslash  (f(b) / \sqrt{a})$,
see~\sref{canonical-filter}.
This does not give a complete description
of the map~$b\mapsto \sqrt{a}b\sqrt{a}$,
though,
since it shares its universal property
with all the maps
$b\mapsto c^*bc,\,\scrA\to\ceil{a}\!\scrA\!\ceil{a}$
where~$c\in\scrA$ with $c^*c = a$,
but that
is a  challenge for the next chapter.

Returning to division again,
another
application
is the polar decomposition
of an element~$a$ of a von Neumann algebra~$\scrA$,
see~\sref{polar-decomposition},
which is simply
\begin{equation*}
	a \,=\, (a/ \sqrt{a^*a})\, \sqrt{a^*a}.
\end{equation*}

Before we get down to business,
let us indicate the difficulty
in defining~$a/b$
for~$a$ and~$b$ that obey~\eqref{douglas-ineq}.
Surely, if~$b$ is invertible,
then we could simply put $a/b:=ab^{-1}$;
and also if~$b$ is just \emph{pseudoinvertible}
in the sense that $b^{\sim 1 }b=\ceill{b}$
and $bb^{\sim 1}=\ceilr{b}$
for some~$b^{\sim 1}$
the formula $a/b:=ab^{\sim 1}$ would work.
But,
of course,
$b$ need not be pseudoinvertible.
The  ideal of~$b^{\sim 1}$
can however be approximated
in an appropriate sense by a formal series $\sum_n t_n$ 
(which we call an \emph{approximate pseudoinverse})
so that we can take
$a/b:= \sum_n a t_n$
(using ultrastrong completeness to see that
the series converges.)
\end{point}
\end{parsec}
\subsection{(Approximate) Pseudoinverses}
\begin{parsec}{790}%
\begin{point}{10}[dfn-pseudoinverse]{Definition}%
Let~$a$ be an element of a von Neumann algebra~$\scrA$.
We'll say that~$a$ is \Define{pseudoinvertible}%
\index{pseudoinverse}
if it has a \Define{pseudoinverse},
that is,
an element~$t$ of~$\scrA$
with 
$ta=\ceill{a}=\ceilr{t}$
and $at=\ceill{t}=\ceilr{a}$.
When such~$t$ exists,
it is unique (by~\sref{mult-cancellation}),
and we'll denote it by~\Define{$a^{\sim 1}$}.%
\index{*sim1@$a^{\sim 1}$, pseudoinverse of~$a$}
If~$a^{\sim1}=a^*$,
we say that~$a$ is a \Define{partial isometry}%
\index{partial isometry!in a von Neumann algebra}
(see~\sref{partial-isometry-equivalents}).
\end{point}
\begin{point}{20}[pseudoinverse-equivalents]{Lemma}%
For elements $a,t$ of a von Neumann algebra
the following are equivalent.
\begin{enumerate}
\item
\label{pseudoinverse-1}
$ta$ is a projection, and~$\ceill{t}=\ceilr{a}$.
\item
\label{pseudoinverse-2}
$ata=a$, and $\ceill{t}\leq \ceilr{a}$ and~$\ceilr{t}\leq \ceill{a}$.
\item
\label{pseudoinverse-3}
$at$ is a projection, and $\ceill{a}=\ceilr{t}$.
\item
\label{pseudoinverse-4}
$tat=t$, and~$\ceill{a}\leq \ceilr{t}$ and~$\ceilr{a}\leq \ceill{t}$.
\item
\label{pseudoinverse-5}
$t$ is a pseudoinverse of~$a$.
\item
\label{pseudoinverse-6}
$a$ is a pseudoinverse of~$t$.
\end{enumerate}
\spacingfix%
\begin{point}{30}{Proof}%
\grayed{(\ref{pseudoinverse-5}$\iff$%
\ref{pseudoinverse-6})}\ 
is clear.
For the remainder we  make two loops.
\grayed{(\ref{pseudoinverse-1}$\Longrightarrow$\ref{pseudoinverse-2})}\ 
We have $\ceill{t}\leq \ceilr{a}$ by assumption,
and $\ceilr{t}=\ceilr{t\ceill{t}}
=\ceilr{t\ceilr{a}}
=\ceilr{ta}=ta=\ceill{ta}\leq\ceill{a}$.
Further, $ata=a$
by~\sref{mult-cancellation},
because $tata=ta$ (since~$ta$ is a projection)
and~$\ceilr{ata}\leq \ceilr{a}\leq\ceill{t}$.
\grayed{(\ref{pseudoinverse-3}$\Longrightarrow$\ref{pseudoinverse-4})} 
follows along the same lines.
\grayed{(\ref{pseudoinverse-2}$
	\Longrightarrow$%
\ref{pseudoinverse-5})}\ 
We have $ta=\ceill{a}$ by~\sref{mult-cancellation},
because $ata=a=a\ceill{a}$,
	and $\ceilr{ta}\leq\ceilr{t}\leq\ceill{a}$.
Also, $at=\ceilr{a}$,
(because~$ata=a=\ceilr{a}a$,
and~$\ceill{at}\leq\ceill{t}\leq \ceilr{a}$).
Further, $\ceill{t}=\ceilr{a}$,
because 
	$\ceilr{a}=at=\ceill{at}\leq \ceill{t}\leq\ceilr{a}$;
and, similarly, $\ceill{a}=\ceilr{t}$.
\grayed{(\ref{pseudoinverse-4}$
	\Longrightarrow$%
\ref{pseudoinverse-5})}
is proven by the same principles, and
\grayed{(\ref{pseudoinverse-5}$
	\Longrightarrow$%
\ref{pseudoinverse-1},\ref{pseudoinverse-3})}
is rather obvious.\qed
\end{point}
\end{point}
\begin{point}{40}[partial-isometry-equivalents]{Exercise}%
Show that an element~$u$ of a von Neumann algebra
is a partial isometry iff
$u^*u$ is a projection
iff $uu^*u=u$
iff $uu^*$ is a projection
iff $u^*uu^*=u^*$
iff~$u^*$ is the pseudoinverse of~$u$.
(Hint: use~\sref{pseudoinverse-equivalents},
or give a direct proof.)
\end{point}
\begin{point}{50}[pseudoinverse-basic]{Exercise}%
Let~$a$ and~$b$ be a elements of a von Neumann algebra~$\scrA$.
\begin{enumerate}
\item
Show that~$a$ is pseudoinvertible
iff~$a^*$ is pseudoinvertible,
and, in that case, $(a^*)^{\sim1}=(a^{\sim1})^*$.
\item
Assuming that~$a$ and~$b$ are pseudoinvertible,
and~$\ceilr{b}=\ceill{a}$,
show that $ab$ is pseudoinvertible,
and~$(ab)^{\sim1}=b^{\sim1}a^{\sim1}$.
\item
Show that~$a$ is pseudoinvertible
iff~$a^*a$ is pseudoinvertible,
and, in that case, $a^{\sim 1} = (a^*a)^{\sim1}a^*$
and~$(a^*a)^{\sim1}=a^{\sim1}(a^{\sim1})^*$.
\end{enumerate}
\spacingfix
\end{point}%
\begin{point}{60}[pseudoinverse-basic-2]{Exercise}%
Let~$a$ be a positive element of a von Neumann algebra~$\scrA$.
\begin{enumerate}
\item
Show that~$a$ is pseudoinvertible iff~$a$
is invertible in~$\ceil{a}\!\scrA\!\ceil{a}$
iff $at=\ceil{a}$ for some~$t\in\scrA_+$.
Show, moreover, that $at=ta$ for such~$t$.
\item
Show that~$a$ is pseudoinvertible iff
there is~$\lambda>0$ with $\lambda \ceil{a}\leq a$.
\item
Assume that~$a$ is pseudoinvertible.

Show that~$\ceil{a^{\sim1}}=\ceil{a}$.

Show that if~$b\in \scrA$ commutes with~$a$,
then $b$ commutes with~$a^{\sim 1}$.

(In other words, $a^{\sim1}\in\{a\}^{\square\square}$.)
\item
Show that~$c^{\sim1}\leq b^{\sim1}$
when~$b\leq c$ are pseudoinvertible positive \emph{commuting}
elements of~$\scrA$.
(The statement is still true without
the requirement that~$b$ and~$c$ commute,
but also much harder to prove.)

\item
Show that $(0,0,1,\frac{1}{2}, \frac{1}{3},\dotsc)$
        is not pseudoinvertible in~$\ell^\infty(\N)$.
\end{enumerate}
\end{point}
\end{parsec}
\begin{parsec}{800}%
\begin{point}{10}{Remark}%
Note that the obvious candidate
for the pseudoinverse of~$(0,0,1,\frac{1}{2},\frac{1}{3},\dotsc)$
    from~$\ell^\infty(\N)$
being~$(0,0,1,2,3,\dotsc)$
is not bounded,
    and therefore not an element of~$\ell^\infty(\N)$.
We can nevertheless approximate~$(0,0,1,2,3,\dotsc)$
by the elements
\begin{equation*}
(0,0,1,0,0,\dotsc),\ 
(0,0,1,2,0,\dotsc),\  \dots 
\end{equation*}
    of~$\ell^\infty(\N)$
forming what we will call ``approximate pseudoinverse'' for 
$(0,0,1,\frac{1}{2},\frac{1}{3},\dotsc)$.
That this can also be done for an arbitrary element
of a von Neumann algebra
is what we'll see next.
\end{point}
\begin{point}{20}[approximate-pseudoinverse-def]{Definition}%
An \Define{approximate pseudoinverse}%
\index{approximate pseudoinverse}%
\index{pseudoinverse!approximate}
of an element~$a$ of a von Neumann algebra~$\scrA$
is a sequence~$t_1,t_2,\dotsc$
of elements of~$\scrA$
such that~$t_1a,\ t_2a,\ \dotsc, at_1,\ at_2,\ \dotsc$
are projections with $\sum_n t_na = \ceill{a}=\sum_n \ceilr{t_n}$
and~$\sum_n at_n =\ceilr{a}=\sum_n \ceill{t_n}$.
\end{point}
\begin{point}{30}[approximate-pseudoinverse-reduction]{Exercise}%
Let~$b$ be an element of a von Neumann algebra~$\scrA$,
and let~$t_1,t_2,\dotsc$
be an approximate pseudoinverse
of~$b^*b$.
Show that $t_1b^*,\ t_2b^*,\ \dotsc$
is an approximate pseudoinverse of~$b$.
\end{point}
\begin{point}{40}[approximate-pseudoinverse]{Theorem}%
Every element~$a$ of a von Neumann algebra~$\scrA$
has an approximate pseudoinverse.
\begin{point}{50}{Proof}%
By~\sref{approximate-pseudoinverse-reduction},
it suffices to consider the case that~$a$ is positive.
When~$a=0$ the sequence~$0,0,0,\dotsc$
clearly yields an approximate pseudoinverse for~$a$,
so let us disregard this case,
and assume that~$a$ is positive and non-zero.

Note that $a-1 \,\leq\, a-\frac{1}{2}\,\leq\, a - \frac{1}{3}\,\leq\, \dotsb$
converges in the norm to~$a\equiv a_+$,
and so does $(a-1)_+\,\leq\,(a-\frac{1}{2})_+\,\leq\,\dotsc$,
which converges also ultraweakly to~$\bigvee_n(a-\frac{1}{n})$,
so that~$a=\bigvee_n (a-\frac{1}{n})_+$,
and thus $\ceil{a}=\bigcup_n \ceil{(a-\frac{1}{n})_+}$
by~\sref{ceil-supremum}.

Writing $q_n=
\ceil{\smash{(a-\frac{1}{n})_+}}$
--- and picturing it as the places where~$a\geq \frac{1}{n}$ ---
we have $(a-\frac{1}{n})q_n = (a-\frac{1}{n})_+\geq 0$
(because~$b\ceil{b_+}=b_+$
for a positive element~$b$ of a von Neumann algebra,
by~\sref{ceil-pos-part}),
and so $\frac{1}{n}q_n \leq  aq_n$ for all~$n>0$.

Writing $e_{n}=q_{n+1}-q_n$
for all~$n$ (taking~$q_0:=0$)
--- and thinking of it as the places
where $\frac{1}{n+1}\leq a < \frac{1}{n}$ --- 
we get a sequence of (pairwise orthogonal) projections
$e_1,e_2,\dotsc$ in~$\{a\}^{\square\square}$
with $\sum_n e_n = \ceil{a}$.
By an easy computation
involving the facts that $\frac{1}{n+1}\leq \frac{1}{n}$
and~$aq_n\leq aq_{n+1}$,
we get $\frac{1}{n+1}e_n \leq ae_n \leq  \frac{1}{n}e_n$.

We claim that~$\ceil{ae_n}=\ceil{e_n}$ for any~$n$.
Indeed, on the one hand~$ae_n=e_nae_n\leq \|a\|e_n$
(as~$e_n\in\{a\}^{\square\square}$)
and so~$\ceil{ae_n}\leq \ceil{\|a\|e_n}=e_n$
(using here that $\|a\|\neq 0$),
while on the other hand, $\frac{1}{n+1}e_n\leq ae_n$
gives $e_n\equiv \ceil{\smash{\frac{1}{n+1}e_n}}\leq \ceil{ae_n}$.
In particular, $\frac{1}{n+1}\ceil{ae_n}=\frac{1}{n+1}e_n \leq ae_n$,
	so that~$ae_n$ is pseudoinvertible (by~\sref{pseudoinverse-basic-2}).

Writing~$t_n := (ae_n)^{\sim1}$,
we have~$\ceil{t_n}=e_n$
(since~$\ceil{ae_n}=e_n)$.
Then $t_na=t_n\ceil{t_n}a=t_ne_na=\ceil{ae_n}=e_n$,
and similarly, $at_n = e_n$,
so that $\sum_n at_n = \sum_n t_n a 
=\sum_n e_n =\ceil{a}
= \sum_n \ceil{t_n}$,
making $t_1,t_2,\dotsc$ an approximate pseudoinverse of~$a$.\qed
\end{point}
\end{point}
\end{parsec}
\subsection{Division}
\begin{parsec}{810}%
\begin{point}{10}[division]{Definition}%
Let~$b$ be an element of a von Neumann algebra~$\scrA$,
and let~$a$ be an element of~$\scrA b$
--- so~$a\equiv cb$ for some~$c\in\scrA$.
We denote by~\Define{$a/b$}%
\index{*slash@$a/b$!in a von Neumann algebra}
the (by~\sref{mult-cancellation}) unique
element~$c$ of~$\scrA\!\ceilr{b}$
with~$a=cb$,
and, dually,
given an element~$a$ of~$b\scrA$ 
we denote by~\Define{$b\backslash a$} 
the unique
element~$c$ of~$\ceill{b}\!\scrA$
with~$a=bc$.
\end{point}
\begin{point}{20}{Exercise}%
Let~$a$ and~$b$
be elements of a von Neumann algebra~$\scrA$.
\begin{enumerate}
\item
Show that $c/b$ is an element of~$\ceilr{c}\!\scrA\!\ceilr{b}$
for every element~$c$ of~$b\scrA$.
\item
Show that $(ab)/b = a\ceilr{b}$
and~$b\backslash (ba)=\ceill{b}a$.
\item
Let~$c$ be an element of $a\scrA b$.
Show that
$a \backslash c \,\in\, \scrA b$,
and
$c/b\,\in\,a\scrA$, and
\begin{equation*}
(a\backslash c)/b\ =\  a\backslash(c/b)
\quad =:\  \Define{a\backslash c/b}.
\end{equation*}
\index{*slashslash@$a\backslash c/b$!in a von Neumann algebra}
Show that~$a\backslash c/b$
is the unique element~$d$ of~$\ceill{a}\scrA\ceilr{b}$
with~$c=adb$. 
\item
Let~$c$ be an element of~$\scrA b$
and let~$d$ be an element of~$a\scrA$.

Show that~$dc\in a\scrA b$,
and~$a\backslash (dc) / b  = (a\backslash d)\,(c/b)$.
\item
Let~$c$ be an element of~$\scrA b$.
Show that~$c^* \,\in\, b^* \scrA$
and $b^*\backslash c^* = (c/b)^*$.
\end{enumerate}%
\spacingfix%
\end{point}%
\begin{point}{30}[proto-douglas]{Lemma}%
Given elements $a$ and~$b$ of a von Neumann algebra~$\scrA$
with~$a^*a \leq b^*b$
we have $a\,\in\, \scrA b$.
Moreover,
given an approximate pseudoinverse $t_1,t_2,\dotsc$
of~$b$,
the series $\sum_nat_n$ converges ultrastrongly to~$a/b$,%
\index{*slash@$a/b$!in a von Neumann algebra}
and uniformly so in~$a$.
\begin{point}{40}{Proof}%
To show that~$\sum_{n=0}^N at_n$
converges ultrastrongly as~$N\to \infty$
it suffices
to show that~$(\,\sum_{n=0}^N at_n\,)_N$
is ultrastrongly Cauchy
(because~$\scrA$ is ultrastrongly complete, by~\sref{vn-complete}).
To this end, note that
\begin{alignat*}{3}
\textstyle 
(\, \sum_{n=M}^N at_n\,)^*\ \,\sum_{n=M}^N\,at_n
\ &= \textstyle\ (\sum_{n=M}^Nt_n^*) \,a^*a\, (\sum_{n=M}^N t_n)
\\
\ &\leq\ \textstyle 
 (\sum_{n=M}^N t_n^*) \,b^*b\, (\sum_{n=M}^Nt_n)
\\
\ &= \ \textstyle 
\sum_{n,m=M}^N
t_n^*b^*bt_m 
\\
\ &=\ \textstyle 
\sum_{m=M}^N
bt_m,
\end{alignat*}
where we've used that
$bt_1,\ bt_2,\ \dotsc$
are pairwise orthogonal projections
--- but then the series~$\sum_{n=0}^\infty bt_m$
converges ultraweakly by~\sref{sum-of-orthogonal-projections}.
This, coupled with the
inequality above,
gives us that~$\sum_{n=0}^N at_n$
is ultrastrongly Cauchy,
and therefore converges ultrastrongly --- 
and even uniformly so in~$a$,
because~``$a$'' does not appear in
the expression
``$\sum_{m=M}^N bt_m$'' that gave the bound.

Define~$c:=\sum_{n=0}^\infty at_n$.
Since~$a^*a \leq b^*b$,
we have $\ceill{a}\leq \ceill{b}$,
and so~$a=a\ceill{b} = a\sum_n t_n b = \sum_n at_n b = cb$. 
So to get~$c=a/b$
we only need to prove that~$\ceill{c}\leq \ceilr{b}$,
that is, $c\ceilr{b}=c$.
To this end,
recall that~$\sum_n \ceill{t_n}=\ceilr{b}$,
so that~$\ceill{t_n}\leq \ceilr{b}$,
and~$t_n\ceilr{b}=t_n$,
which implies
that~$at_n\ceilr{b}=at_n$,
and so~$c\ceilr{b}=\sum_n at_n\ceilr{b}
= \sum_n at_n=c$.\qed
\end{point}
\end{point}
\begin{point}{50}[douglas]{Exercise}%
\index{Douglas' Lemma}
Let~$a$ and~$b$ be elements of a von Neumann algebra~$\scrA$.
\begin{enumerate}
\item
Let~$\lambda\geq 0$ be given,
and recall that $(\scrA)_\lambda 
= \{c\in\scrA\colon \|c\|\leq\lambda\}$.

Show that $a$ is in $(\scrA)_\lambda b$ iff $a^*a\leq \lambda^2 b^*b$,
and then~$\|a/b\|\leq \lambda$.

(Compare this with ``Douglas'~Lemma'' from~\cite{douglas}.)
\item
Show that $a \in \scrA\!\ceilr{b}$
need not entail that~$a\in \scrA b$.
\end{enumerate}
\spacingfix%
\end{point}%
\begin{point}{60}[sequential-douglas]{Exercise}%
Let~$b$ be an element of a von Neumann algebra~$\scrA$.
\begin{enumerate}
\item
Let~$a$ be a positive element of~$\scrA$,
and let~$\lambda\geq 0$.

Show that
 $a\in b^*(\scrA)_\lambda b$
iff~$a\leq \lambda b^*b$,
and then~$\|b^*\backslash a / b \| \leq\lambda$.
\item
Show that~$b^*\backslash a / b$ is positive
for every positive element~$a$ of $ b^* \scrA b$.

(Hint: prove that $(b^*\backslash \sqrt{a})\,(\sqrt{a}/b)
= b^*\backslash a/b$.)
\end{enumerate}%
\spacingfix%
\end{point}%
\begin{point}{70}[div-approx]{Exercise}%
Given elements~$b$ and~$c$ of a von Neumann algebra~$\scrA$,
an approximate pseudoinverse $t_1,t_2,\dotsc$
of~$b$,
and an approximate pseudoinverse of~$s_1,s_2,\dotsc$
of~$c$,
show that
$(\sum_{n=1}^N s_n) \,a\, (\sum_{m=1}^Nt_m)$,
converges ultrastrongly to~$c\backslash a /b$%
\index{*slashslash@$a\backslash c/b$!in a von Neumann algebra}
as~$N\to\infty$
(and uniformly so) for~$a\in c(\scrA)_1 b$.
\end{point}
\begin{point}{80}[sequential-quotient]{Exercise}%
Show that for positive elements~$a$ and~$b$ of a von Neumann
algebra~$\scrA$,
the following are equivalent.
\begin{enumerate}
\item
$a\leq \lambda b$ for some~$\lambda\geq 0$;
\item
$a=\sqrt{b}c\sqrt{b}$
for some positive~$c\in\scrA$.
\end{enumerate}
In that case, there is a unique~$c\in\scrA_+$
with $a=\sqrt{b}c\sqrt{b}$
and~$\ceil{c}\leq \ceil{b}$.
Moreover,
if~$t_1,t_2,\dotsc$
is an approximate pseudoinverse of~$\sqrt{b}$,
then~$\sum_{m,n} t_m a t_n$
converges ultraweakly to such~$c$.
\end{point}
\begin{point}{90}[div-usc]{Lemma}%
Given elements~$b$ and~$c$  of a von Neumann algebra~$\scrA$
the maps 
\begin{equation*}
	a\mapsto a/b\colon \  (\scrA)_1b \to \scrA
\qquad\text{and}\qquad
a\mapsto c\backslash a/b\colon\  c(\scrA)_1b\to \scrA
\end{equation*}
are ultrastrongly continuous
(where $(\scrA)_1$ is the unit ball).
\begin{point}{100}{Proof}%
By~\sref{proto-douglas}
the series~$\sum_n at_n$ converges ultraweakly to~$a/b$,
where~$t_1,t_2,\dotsc$
is an approximate pseudoinverse of~$b$,
and in fact uniformly so for~$a\in(\scrA)_1 b$
(because $a^*a\leq b^*b$ for such~$a$).
Since~$a\mapsto \sum_{n=1}^N at_n,\ (\scrA)_1b\to\scrA$
is ultrastrongly continuous (by~\sref{mult-uws-cont})
--- and the uniform limit of continuous functions is continuous ---
we see that~$a\mapsto a/b,
\ (\scrA)_1b\to\scrA$ is ultrastrongly continuous.
It follows that
$(\,\cdot\,)/b\colon\,c(\scrA)_1b\to c(\scrA)_1$
and~$c\backslash(\,\cdot\,)\colon\, c(\scrA)_1\to\scrA$
are ultrastrongly continuous;
as must be their composition $c\backslash\,\cdot\,/b\colon\,
c(\scrA)_1b\to\scrA$.\qed
\end{point}
\begin{point}{110}{Remark}%
The map $a\mapsto a/b$ might not give
an ultrastrongly
continuous map on the larger domain~$\scrA b$,
because, for example, 
upon applying $(\,\cdot\,)/(1,\frac{1}{2},\frac{1}{3},\dotsc)$
to the ultrastrongly Cauchy 
sequence  $(1,0,0,\dotsc),\ (1,1,0,\dotsc),\ \dotsc$
    in~$\ell^\infty(\N)$
we get the 
sequence $(1,0,0,\dotsc),\ (1,2,0,\dotsc),\ \dotsc$,
which is not ultrastrongly Cauchy.
\end{point}
\end{point}
\end{parsec}
\subsection{Polar Decomposition}
\begin{parsec}{820}%
\begin{point}{10}[polar-decomposition]{Proposition (Polar Decomposition)}%
\index{polar decomposition!of an element of a von Neumann algebra}%
Any element~$a$ of a von Neumann
algebra~$\scrA$
can be uniquely written as~$a=\Define{[a]}\sqrt{a^*a}$,%
    \index{*brackets@$[\,\cdot\,]$!$[a]$, partial isometry from the polar decomposition of~$a$}
where~$[a]$
is an element of~$\scrA\!\ceill{a}$.
Moreover,
\begin{enumerate}
\item
$[a]$ is partial isometry
with $[a]^*[a] = \ceil{a^*a} \equiv \ceill{a}$
and $[a][a]^*=\ceil{aa^*} \equiv \ceilr{a}$,
\item
and $[a^*]=[a]^*$, so that~$\sqrt{aa^*}[a]=a=[a]\sqrt{a^*a}$.
\end{enumerate}
\spacingfix%
\begin{point}{20}{Proof}%
Since $a^*a\leq \sqrt{a^*a}\sqrt{a^*a}$,
the existence and uniqueness
of an element~$[a]$ of~$\scrA$
with $a=[a]\sqrt{a^*a}$
and~$\ceill{\,[a]\,}\leq
\ceill{a}\equiv \ceilr{\,\sqrt{a^*a}\,}$
is provided by~\sref{douglas},
and we get~$\ceilr{\,[a]\,}\leq \ceilr{a}$
to boot!
Note that~$[a]^*[a]=\ceil{a^*a}$,
by~\sref{mult-cancellation},
because 
\begin{equation*}
	\sqrt{a^*a}\,[a]^*[a]\,\sqrt{a^*a}
\ =\ 
a^*a
\ =\ 
\sqrt{a^*a}\,\ceil{a^*a}\,\sqrt{a^*a},
\end{equation*}
and $\ceil{\,[a]^*[a]\,}\leq \ceill{a}=\ceil{\,\sqrt{a^*a}\,}$.
In particular, $[a]$ is a partial isometry
(by~\sref{partial-isometry-equivalents}).

Let us prove that $[a][a]^*=\ceilr{a}$.
Note that~$[a][a]^*$
is a projection (because $[a]$ is a partial isometry, 
by~\sref{partial-isometry-equivalents}).
We already know that $[a][a]^* = \ceilr{\,[a]\,}\leq \ceilr{a}$.
Concerning the other direction,
$aa^*=[a]\sqrt{a^*a}\,\sqrt{a^*a}[a]^*=[a]\,a^*a\,[a]^*$,
so that
$\ceilr{a}=\ceil{aa^*}
= \ceil{\,[a]a^*a[a]^*\,}
\leq \ceil{\,\|a\|^2[a][a]^*\,}
	=\ceil{[a][a]^*}
\leq [a][a]^*$.

To prove that~$a=\sqrt{aa^*}[a]$,
we'll first show that~$\sqrt{aa^*} = [a] \sqrt{a^*a} [a]^*$.
Indeed, since~$[a]^*[a]=\ceil{\,\sqrt{a^*a}\,}$,
we have
$[a]\sqrt{a^*a}[a]^*[a]\sqrt{a^*a}[a]^*
= [a] \sqrt{a^*a} \,\sqrt{a^*a} [a]^* 
= aa^*$ --- now take the square root.
It follows that~$\sqrt{aa^*}[a]
= [a]\sqrt{a^*a}[a]^*[a]
= [a]\sqrt{a^*a}=a$.
Finally, upon applying~$(\,\cdot\,)^*$,
we see that
$a^*=[a]^*\sqrt{aa^*}$,
and thus~$[a^*]=[a]^*$,
by uniqueness of~$[a^*]$,
because $\ceill{\,[a]^*\,}=\ceilr{\,[a]\,}=  \ceilr{a}
=\ceill{a^*}$.\qed
\end{point}
\end{point}
\end{parsec}
\begin{parsec}{830}%
\begin{point}{10}%
Recall from~\sref{cceil-fundamental}
that the least central projection~$\cceil{e}$
above a projection~$e$ of a von Neumann algebra~$\scrA$
is given by~$\cceil{e}=\bigcup_{a\in \scrA}\ceil{a^* e a}$.
Using the polar decomposition
we can give a more economical
description of~$\cceil{e}$,
see~\sref{cceil-sum}.
\end{point}
\begin{point}{20}[vmleq]{Proposition}%
Given projections $e'$ and~$e$
of a von Neumann algebra~$\scrA$,
the following are equivalent.
\begin{enumerate}
\item
\label{vmleq-1}
$e'= \ceil{a^*e a}$ for some~$a\in\scrA$;
\item
\label{vmleq-2}
$e' = \ceill{a}$ and~$\ceilr{a}\leq e$
for some~$a\in\scrA$;
\item
\label{vmleq-3}
$e' = u^*u$ and~$uu^*\leq e$
for some partial isometry~$u$.
\end{enumerate}
In that case
we write~$\Define{e'\vmleq e}$%
\index{((vmleq@$\vmleq$, Murray--von Neumann preorder}
(and say~$e'$ is \Define{Murray--von Neumann below}~$e$).%
\index{Murray--von Neumann preorder}
\begin{point}{30}{Proof}%
That~\ref{vmleq-3}
implies~\ref{vmleq-2} is clear.
(\ref{vmleq-2}$\Rightarrow$\ref{vmleq-1})\ 
Since~$\ceilr{a}\leq e$,
we have $ea=a$,
and so~$\ceil{a^*ea}=\ceil{a^*a}=\ceill{a}=e'$.
(\ref{vmleq-1}$\Rightarrow$\ref{vmleq-3})\ 
By the polar decomposition (see~\sref{polar-decomposition})
we get a partial isometry~$u:=[ea]$
for which $u^*u = [ea]^*[ea]=\ceil{(ea)^*ea}=e'$
and~$uu^*=\ceil{eaa^*e}\leq e$.\qed
\end{point}
\end{point}
\begin{point}{40}[mvn-preorders]{Exercise}%
Show that~$\vmleq$  preorders
the projections of a von Neumann algebra.
\end{point}
\begin{point}{50}[cceil-sum]{Lemma}%
\index{*cceil@$\cceil{\,\cdot\,}$!$\cceil{a}$, central support}%
Given a projection~$e$
of a von Neumann algebra~$\scrA$
there is a family
$(e_i)_i$ of non-zero projections with $\cceil{e}=\sum_i e_i$,
and $e_i\vmleq e$ for all~$i$.
\begin{point}{60}{Proof}%
Let~$(e_i)_i$
be a maximal set of non-zero pairwise orthogonal projections
in~$\scrA$ with
$e_i\vmleq e$ for all~$i$.
Our goal is to show that $\sum_i e_i \equiv \bigcup_i e_i=\cceil{e}$.

Let~$u_i$ be a partial isometry
with $u_i^*u_i=e_i$ and $u_iu_i^* \leq e$.
Since $e_i=u_i^* u_i = u_i^* u_i u_i^* u_i
\leq u_i^* e u_i 
\leq \bigcup_{a\in \scrA} \ceil{a^* e a}=\cceil{e}$,
we have
$\bigcup_i e_i \leq \cceil{e}$.

Suppose that~$\bigcup_i e_i < \cceil{e}$
(towards a contradiction).
Then since $p:=\cceil{e}-\bigcup_i e_i $ 
is a non-zero projection,
and~$p=p\cceil{e} p
=\bigcup_{a\in\scrA} \ceil{p\ceil{a^*ea}p}
= \bigcup_{a\in\scrA} \ceil{(eap)^*eap}$,
there must be~$a\in\scrA$ with $(eap)^*eap\neq 0$.
The polar decomposition (see~\sref{polar-decomposition})
of~$eap$
gives us a partial isometry $u:=[eap]$
with $uu^* = \ceil{eap(eap)^*}=\ceil{eapa^*e}\leq e$
and $u^*u= \ceil{(eap)^*eap}\leq p$,
so that $u^*u$ is a non-zero projection,
orthogonal to all~$e_i$ with $u^*u \vmleq e$.
In other words,
$e$ could have been added to~$(e_i)_i$,
contradicting its maximality.
Hence $\bigcup_i e_i = \cceil{e}$.
\qed
\end{point}
\end{point}
\end{parsec}
\begin{parsec}{840}%
\begin{point}{10}%
Using~\sref{vmleq} we can classify
all finite-dimensional $C^*$-algebras.
\end{point}
\begin{point}{20}[fdcstar]{Theorem}%
\index{Cstar-algebra@$C^*$-algebra!finite dimensional}%
\index{von Neumann algebra!finite dimensional}%
Any finite-dimensional $C^*$-algebra~$\scrA$
is
a direct sum
of full matrix algebras,
that is,
$\scrA\cong \bigoplus_m M_{N_m}$
for some  $N_1,\dotsc,N_M\in\N$.
\begin{point}{30}{Proof}%
Let~$e_1,\dotsc,e_N$ be a basis for~$\scrA$.
We'll first show that~$\scrA$ is a von Neumann algebra,
and for this we'll 
need the fact that the unit ball $(\scrA)_1$
is compact with respect to the \emph{norm} on~$\scrA$.
For this it suffices
to show that~$\|\,\cdot\,\|$
is equivalent to the norm~$\|\,\cdot\,\|'$ on~$\scrA$
given by~$\|a\|'=\sum_n\left|z_n\right|$
for all $a\equiv \sum_n z_n e_n$
where~$z_1,\dotsc,z_N\in\C$,
(because the unit $\|\,\cdot\,\|'$-ball is clearly compact
being homeomorphic to the unit ball of~$\C^N$.)
Since for such~$a\equiv \sum_n z_n e_n$
we have
\begin{equation*}
	\textstyle
\|a\|\,\leq\, \sum_n\left|z_n\right| \left\|e_n\right\|
\,\leq\,\sum_n\left|z_n\right| \,\sup_n\left\|e_n\right\|
	\,=\, \|a\|'\,\sup_n\left\|e_n\right\|
\end{equation*}
we see that~$a\mapsto a\colon \scrA\to\scrA$
is continuous from~$\|\,\cdot\,\|'$
to~$\|\,\cdot\,\|$.
For the converse
it suffices to show that
$f_m\colon\, a\equiv \sum_n z_n \mapsto z_m,\, 
\scrA\to\C$
is bounded with respect to~$\|\,\cdot\,\|$,
because then
\begin{equation*}
	\textstyle
	\|a\|'\,\equiv\, \|\sum_n f_n (a)e_n\|'
	\,\leq\, \sum_n\left|f_n(a)\right| 
	\,\leq\, (\sum_n \|f_n\|)\,\|a\|.
\end{equation*}
In fact,
we'll show that any linear functional on~$\scrA$
is bounded.
Since the bounded linear functionals
form a linear subspace~$\scrA^*$
of $N$-dimensional vector space of all linear functionals 
on~$\scrA$
it suffices to show that~$\scrA^*$ has dimension~$N$.
So let~$f_1,\dotsc,f_M$ be a basis for~$\scrA^*$;
we must show that~$N\leq M$.
Since the states of~$\scrA$
(see~\sref{states-order-separating})
and thus all linear functionals on~$\scrA$
form a separating collection,
the functionals $f_1,\dotsc,f_N$
form a separating set too;
since therefore
\begin{equation*}
a\mapsto (f_1(a),\,\dotsc,\,f_M(a))\colon\,
\scrA\to\C^M
\end{equation*}
is a linear injection
from the~$N$-dimensional
space~$\scrA$ to the $M$-dimensional space~$\C^M$
we get~$N\leq M$.
Whence all linear functionals on~$\scrA$
are bounded, the norms $\|\,\cdot\,\|$ and~$\|\,\cdot\,\|'$
are equivalent,
and $(\scrA)_1$
is norm compact.
\begin{point}{40}[suprema-in-fdvna]{$\scrA$ is a von Neumann algebra}%
First we need to show that every bounded directed
set~$D$ of self-adjoint elements
of~$\scrA$ has a supremum (in $\Real{\scrA}$).
We may assume without loss of generality that~$\|d\|\leq 1$
for all~$d\in D$, and so~$D\subseteq (\scrA)_1$.
Since~$(\scrA)_1$
is norm compact
there is a cofinal subset~$D'$ of~$D$
that norm converges to some~$a\in\scrA$,
and thus~$D$ norm converges to~$a$ itself.
It's easily seen that~$a$ is the supremum of~$D$.
Indeed, given $d_0\in D$
we have~$d_0\leq d$ for all~$d\geq d_0$,
and so~$d_0\leq \lim_{d\geq d_0} d=a$.
Hence~$a$ is an upper bound for~$D$;
and if~$b$ is an upper bound for~$D$,
then~$d\leq b$ for all~$d\in D$, and so $a=\lim_d d \leq  b$.

Since in this finite-dimensional setting~$\bigvee D$
is apparently the norm limit of~$(d)_{d\in D}$,
any positive functional $f$ on~$\scrA$
will map~$\bigvee D$ to the limit of
$(f(d))_{d\in D}$,
which is~$\bigvee_{d\in D} f(d)$,
and so~$f(\bigvee D)=\bigvee_{d\in D} f(d)$.
Whence every positive functional on~$\scrA$
is normal; and since the positive functionals on~$\scrA$
form a separating collection,
$\scrA$ is a von Neumann algebra.
\end{point}
\begin{point}{50}{Reduction to a factor}%
Since pairwise orthogonal non-zero projections
are easily seen to be linearly independent,
and~$\scrA$ is finite dimensional,
every orthogonal set of projections in~$\scrA$ is finite.
In particular,
any descending sequence of non-zero projections must eventually become
constant.
It follows that below every (central) projection~$p$ in~$\scrA$
there is a minimal (central) projection,
and even that~$p$ is the finite sum of minimal (central) projections.
In particular,
the unit~$1$ of~$\scrA$ can be written
as $1=\sum_n z_n$ where~$z_1,\dotsc,z_M$ are minimal central projections
of~$\scrA$.
By~\sref{central-projections-sums}
we know that~$z_m\scrA$ is a von Neumann algebra for each~$m$,
and that~$\scrA$ is nmiu-isomorphic
to the direct sum $\bigoplus_m z_m \scrA$
of these von Neumann algebras
via~$a\mapsto (z_ma)_m$.
Since~$z_m$ is a minimal central projection,
the von Neumann algebra~$z_n \scrA$
has no non-trivial central projections.
\end{point}
\begin{point}{60}{When~$\scrA$ is a factor}%
Let~$e$ be a minimal
projection of~$\scrA$
(which exists by the previous discussion).
Since~$e\neq 0$,
and~$\scrA$ has no non-trivial central projections,
we have~$\cceil{e}=1$.
By~\sref{cceil-sum}
we have $1\equiv \cceil{e}=\sum_k e_k$
for some non-zero projections $e_1,\dotsc,e_K$
in~$\scrA$
with~$e_k\vmleq e$.
So there are partial isometries
$u_1,\dotsc,u_K\in\scrA$
with $u_k^* u_k = e_k$
and~$u_ku_k^*\leq e$ for all~$k$.
In fact,
since~$e$ is minimal,
we have~$u_ku_k^* = e$.
Thinking of~$u_k$
as~$\ketbra{0}{k}$
define~$u_{k\ell } = u_k ^* u_\ell$;
we'll show that
$\varrho\colon A\mapsto \sum_{k\ell} A_{k\ell}  u_{ k \ell}
\colon M_K\to\scrA$
is an miu-isomorphism.
It's easy to see that~$\varrho$
is linear, involution preserving and unital.
To see that~$\varrho$ is multiplicative,
first note that $u_j u_k^*$ equals~$e$ when~$j=k$
and is zero otherwise.
It follows
that~$u_{ij} u_{k\ell}$
equals $u_{i\ell}$
when~$k=j$ and is zero otherwise.
Whence
\begin{equation*}
	\textstyle
\varrho(A)\varrho(B)
\ =\  
	\sum_{ijk\ell}
	A_{ij}u_{ij} B_{k\ell}u_{k\ell}
\ = \ 
	\sum_{i\ell} ( \sum_k A_{ik} B_{k\ell}) u_{i\ell}
	\ = \ \varrho( AB)
\end{equation*}
for all matrices $A,B\in M_K$,
and so~$\varrho$ is multiplicative.

It remains to be shown that~$\varrho$ is a bijection.
To see that~$\varrho$ is injective,
first note that~$\varrho$
is normal,
because using the fact that~$\varrho$ is positive
and thus bounded,
we can show that~$\varrho$ preserves suprema of bounded directed
sets in much the same way we showed
that all np-functionals on~$\scrA$ are bounded.
We can thus speak of the central carrier~$\cceil{\varrho}$
of~$\varrho$,
and thus to show that~$\varrho$ is injective
it suffices to show that~$\cceil{\varrho}=1$.
Since~$M_K$ is a factor
(see \sref{central-examples})
the only alternative
is~$\cceil{\varrho}=0$
i.e.~$\varrho=0$,
which is clearly absurd
unless~$\scrA=\{0\}$
in which case
we'd already be done.
Hence~$\varrho$ is injective.

To see that~$\varrho$ is surjective
let~$a\in\scrA$ with $a\neq 0$ be given.
Since~$a\equiv \sum_{k,\ell} e_k a e_\ell
= \sum_{k,\ell}u_{k1} u_{1k}au_{\ell 1} u_{1\ell}$,
and~$u_{k1}$ and~$u_{1\ell}$ are in the range of~$\varrho$
it suffices
to show that $u_{1k}a u_{\ell 1}$
is in the range of~$\varrho$
for all~$k$ and~$\ell$.
In other words,
we may assume without loss of generality
that~$e a e = a$,
where~$e$ is the minimal projection in~$\scrA$ we started with.
Since~$e(\Real{a})_+e= (\Real{a})_+$, and so on,
we may assume that~$a$ is positive.
By scaling,
we may also assume that~$\|a\|\leq \nicefrac{1}{3}$.
Since $\ceil{\|a\|e -a}\leq e$,
and~$e$ is minimal,
we either have~$\ceil{\|a\|e-a}=e$
or~$\ceil{\|a\|e-a}=0$.

The former case is impossible:
indeed,
if $e=\ceil{\|a\|e-a}\equiv \bigvee_n (\|a\|e-a)^{\nicefrac{1}{2^n}}$
(see~\sref{vna-ceil}),
then~$(\|a\|e-a)^{\nicefrac{1}{2^n}}$
norm converges to~$\ceil{\|a\|e-a}=e$
(cf.~\sref{suprema-in-fdvna}),
and so $\|\|a\|e-a\|^{\nicefrac{1}{2^n}}$
converges to $\|e\|=1$.
Then~$\|\|a\|e-a\|=1$,
while~$\|\|a\|e-a\|\leq \|a\|\|e\|+\|a\|\leq  \frac{2}{3}$,
which is absurd.

Hence~$\ceil{\|a\|e-a}=0$,
and so~$a=\|a\|e$.
In particular, $a$ is in the range of~$\varrho$.
Whence~$\varrho$ is surjective,
and thus an miu-isomorphism $M_N\to\scrA$.\qed
\end{point}
\end{point}
\end{point}
\end{parsec}
\begin{parsec}{841}%
\begin{point}{10}[cstar-no-pu-equalisers-example]{Example}%
Using the description of all finite-dimensional $C^*$-algebras
    from~\sref{fdcstar}
we can prove the  claim made at the start of this thesis,
in~\sref{cstar-no-pu-equalisers},
that in~$\Cstar{pu}$ there's no equaliser
for the maps $f,g\colon \C^4\to \C$
given  by
\begin{equation*}
\textstyle
f(a,b,c,d)\,=\, \frac{1}{2}(a+b),
\quad \text{and}\quad
g(a,b,c,d)\,=\, \frac{1}{2}(c+d).
\end{equation*}
Indeed, suppose towards a contradiction
that $f$ and~$g$ do have an
equaliser $e\colon \scrE\to\C^4$ in~$\Cstar{pu}$,
and let~$\scrS$
denote the set-theoretic equaliser:
\begin{alignat*}{3}
    \scrS &:=\ \{\,(a,b,c,d)\in\C^4\colon\, f(a,b,c,d)=g(a,b,c,d)\,\}\\
    &=\ \{\,(a,b,c,d)\in\C^4\colon\,a+b\,=\,c+d\,\}.
\end{alignat*}
Note that the elements~$s\in \scrS$
with $0\leq s\leq 1$
form a convex subset of~$\C^4$ that is isomorphic
to an octahedron---this will be essential later.

We claim that the range of~$e\colon \scrE\to\C^4$
is simply the set-theoretic equaliser, $e(\scrE)=\scrS$.
Indeed, surely, $e(\scrE) \subseteq \scrS$.
For the other direction,
let~$v\in \scrS$ be given; we must find~$a\in\scrE$ with~$e(a)=v$.
Since~$e$ is involution preserving,
and so~$\Real{v},\Imag{v}\in\scrS$,
we may assume without loss of generality that~$v$ is self-adjoint.
Since~$v+\|v\|\geq 0$, and~$e$ is unital, 
    we may assume that~$v$ is positive too.
By scaling~$v$ if necessary,
we may assume also that~$0\leq v\leq 1$.
Now, to use the universal property of~$e\colon \scrE\to\C^4$
consider the unique pu-map $p\colon \C^2\to\C^4$
    given by~$p(1,0)=v$.
Since~$v\in \scrS$ we have $f\circ p = g\circ p$,
and so there is a unique $q\colon \C^2\to\scrE$
with $p = e\circ q$.  Then~$v=p(1,0)=e(q(1,0))$,
    and so~$e(\scrE)=\scrS$.

The next thing to note is that~$e$
is injective,
and for this it suffices to show
that~$e$ is injective on $[0,1]_\scrE$.
So let~$a,b\in[0,1]_\scrE$ with~$e(a)=e(b)$ be given;
we must show that~$e(a)=0$.
Let~$p,q\colon \C^2\to\C^4$ be the unique
pu-maps given by~$p(1,0)=a$ and~$q(1,0)=b$,
and note that~$e\circ p = e\circ q$.
Since equalisers are mono, we get~$p=q$, and so~$a=b$.

Thus~$e\colon \scrE\to \C^4$
gives a linear isomorphism from~$\scrE$
onto the 3-dimensional linear subspace~$\scrS$ of~$\scrA$,
so~$\scrE$ is 3-dimensional too.
By the classification of finite-dimensional
$C^*$-algebras, \sref{fdcstar}, 
$\scrE$ must be miu-isomorphic to~$\C^3$.

The map $e\colon \scrE\to\C^4$
is not only injective,
but in fact bipositive (see~\sref{cstar-isometry}).
Indeed,  if $e(a)\leq 0$ for some $a\in[0,1]_\scrE$
we can, as before,
find~$q\colon \C^2 \to \scrE$
with $e(q(1,0))=e(a)$,
and so~$a=q(1,0)\geq 0$.
It follows that~$e$ gives a linear order isomorphism
between~$\scrE$ and the subspace~$\scrS$ of~$\C^4$,
and so $[0,1]_\scrE$ is as convex space
isomorphic to~$\scrS\cap[0,1]_{\C^4}$.
This is problematic,
because on the one
hand
the convex space $[0,1]_\scrE$ being a cube
(because
$\scrE$ is miu-isomorphic to~$\C^3$)
has eight extreme points,
while on the other 
hand~$\scrS\cap[0,1]_{\C^4}$ being
an octahedron
has six extreme points: a contradiction.
\end{point}
\end{parsec}
\subsection{Hereditarily Atomic Von Neumann Algebras}
\begin{parsec}{842}
\begin{point}{10}%
We've seen in~\sref{fdcstar}
that every finite-dimensional von Neumann algebra
is the product of finitely many full matrix algebras.
For our purposes this class is too small,
not admitting interpretations of infinite-dimensional
datatypes,
so we've focused on all von Neumann algebras instead.
There is,
however, a rather modest but very promising subclass
of von Neumann algebras
that does sate our desire for the infinite:
following Kornell
we'll call a von Neumann algebra that is the 
    product of a possibly infinite set
of full matrix algebras
\emph{hereditarily atomic}, see~\sref{def:hereditarily-atomic}.
In his recent paper~\cite{kornell2018quantum}
Kornell develops the position
that these hereditarily atomic von Neumann 
algebras are the ``correct'' quantum generalisation
of sets, and---which is especially relevant to our work---observes
that the category of hereditarily atomic von Neumann 
algebras and the nmiu-maps between
them endowed with the regular 
tensor is monoidal coclosed (see~\cite{kornell2018quantum}, Theorem 9.1.)
This will allow us to build a model of the quantum lambda calculus
not only using all von Neumann algebras,
but also just from the hereditarily atomic ones.

Hereditarily atomic von Neumann algebras
have garnered attention for a completely different reason too:
Selinger observed in Example~2.7 of~\cite{selinger2004towards} 
that the effects of~$M_N$
(and thus of every hereditarily atomic von Neumann algebra)
form a continuous dcpo.
Furber and Weaver have proven recently that the converse also holds:
that every von Neumann algebra~$\scrA$
for which $[0,1]_\scrA$ as dcpo is continuous,
is hereditarily atomic,
see Theorem III.15 of~\cite{furbercontinuous}.
\end{point}
    \begin{point}{20}[def:hereditarily-atomic]{Definition}%
A von Neumann algebra is called \Define{hereditarily atomic}%
\index{von Neumann algebra!hereditarily atomic}%
\index{hereditarily atomic!von Neumann algebra}
if it is nmiu-isomorphic
to a direct sum
$\bigoplus_{i\in I} M_{N_i}$
of possibly infinitely many $M_{N_i}$'s.

We denote by~$\Define{\haW{miu}}$ and  $\Define{\haW{cpsu}}$%
    \index{haWmiu@$\haW{miu},\,\haW{cpsu},\,\dotsc$}
the full subcategories of $\W{miu}$ and $\W{cpsu}$, respectively,
of hereditarily atomic von Neumann algebras.
\end{point}
\begin{point}{30}{Proposition}%
A von Neumann subalgebra~$\scrB$ of a hereditarily atomic
von Neumann algebra~$\scrA$ is itself hereditarily atomic.
\begin{point}{40}{Proof}%
Since~$\scrA$ is hereditarily atomic, 
we may assume without loss of generality
that~$\scrA \equiv \bigoplus_{i \in I} M_{N_i}$
for some family of natural numbers $(N_i)_{i\in I}$.

Note that to show that~$\scrB$ is hereditarily atomic,
it suffices to find a orthogonal family of central
projections $(c_j)_{j\in J}$ in~$\scrB$
with $\sum_j c_j=1$ 
such that  each $c_j\scrB$ is finite-dimensional.
Indeed, then each~$c_j\scrB$ is hereditarily atomic,
and so will be~$\scrB\cong \bigoplus_{j\in J} c_j\scrB$
(see~\sref{central-projections-sums}).

It's even enough to find a family of central projections
$(d_k)_{k\in K}$ in~$\scrB$,
not necessarily orthogonal, but
with $\bigcup_{k\in K} d_k =1$
and each $d_k\scrB$  finite-dimensional.
Indeed, any maximal orthogonal 
family $(c_j)_{j\in J}$
of non-zero central projections in~$\scrB$
for which each~$c_j$ is below some~$d_k$
will have the properties
that $\sum_{j\in J} c_j =1$
and~$c_j\scrB$ be finite-dimensional for every~$j$.

Define $d_j := \ceil{\pi_j\circ e}$
to be the carrier (see~\sref{carrier})
of the inclusion $e\colon \scrB\to\scrA$
followed by the $j$-th projection
$\pi_j\colon \scrA\equiv\bigoplus_{i\in I}M_{N_i}\to M_{N_j}$.
Since $\pi_j\circ e$ is an nmiu-map,
    $d_j$ is a central projection by~\sref{carrier-miu}.
Since there are fewer projection in~$\scrB$ than in~$\scrA$,
we have $\ceil{\pi_j}\leq \ceil{\pi_j\circ e}\equiv d_j$.
Now, since clearly $\sum_{i\in I}\ceil{\pi_i}=1$,
this implies that $\bigcup_{i\in I} d_i =1$.

Let~$i\in I$ be given.
It remains to be shown that~$d_i\scrB$ is finite-dimensional.
To see this, simply note that the restriction
of $\pi_i\circ e$ to a map $d_i\scrB \to M_{N_i}$
    is an injection (by~\sref{carrier-miu})
    into the finite-dimensional space~$M_{N_i}$,
    and so~$d_i\scrB$ is finite-dimensional too.\qed
\end{point}
\end{point}
\begin{point}{50}[ha-equalisers]{Corollary}%
Given nmiu-maps
$f,g\colon \scrA\to\scrB$
between hereditarily atomic
von Neumann algebras~$\scrA$ and~$\scrB$,
the von Neumann subalgebra
\begin{equation*}
    \scrE \ :=\  \{\,a\in\scrA\colon\, f(a)=g(a)\,\}
\end{equation*}
of~$\scrA$
is hereditarily atomic,
and the inclusion~$e\colon \scrE\to\scrA$
is an equaliser of~$f$ and~$g$ both in~$\haW{miu}$
and~$\haW{cpsu}$.
\end{point}
\begin{point}{60}{Remark}%
It follows that~$\haW{miu}$
is the least full subcategory of~$\W{miu}$
closed under limits
that contains all finite-dimensional von Neumann algebras.
\end{point}
\end{parsec}
\section{Normal Functionals}
\begin{parsec}{850}%
\begin{point}{10}%
For our study of the category of von Neumann algebras
we need two more technical results
concerning the normal functionals
on a von Neumann algebra.

The first one,
that a net $(b_\alpha)_\alpha$
in a von Neumann algebra~$\scrA$ is (norm) bounded
provided that $(\omega(b_\alpha))_\alpha$
is bounded for each np-functionals $\omega\colon \scrA\to\C$
(see~\sref{ultraweakly-bounded-implies-bounded}),
ultimately
follows from a type of polar decomposition for 
ultraweakly linear functionals (see~\sref{polar-decomposition-of-functional}).

The second one,
that the ultraweak topology of a von Neumann subalgebra
coincides with the ultraweak topology of the surrounding space
(see~\sref{functional-permanence}),
is proven using the double commutant theorem (\sref{double-commutant})
and requires a lot of hard work.
\end{point}
\end{parsec}
\subsection{Ultraweak Boundedness}
\begin{parsec}{860}%
\begin{point}{10}%
To get a better handle on the normal positive functionals
on a von Neumann algebra,
we first analyse  the not-necessarily-positive normal functionals
in greater detail.
\end{point}
\begin{point}{20}[positive-functional-criterion]{Lemma}%
\index{functional!positive}
A linear map $f\colon \scrA\to \C$
on a $C^*$-algebra~$\scrA$
is positive iff $\|f\|\leq f(1)$.
\begin{point}{30}{Proof}%
(Based on Theorem 4.3.2 of~\cite{kr}.)

If~$f(1)=0$, then~$f=0$ in both cases 
(viz.~$f$ is positive, and~$\|f\|\leq f(1)$),
so
we may assume that~$f(1)\neq 0$.
The problem is easily reduced farther to the case 
that~$f(1)=1$
by replacing~$f$ by~$f(1)^{-1}f$
(noting that~$f(1)\geq 0$ in both cases),
so we'll assume that~$f(1)=1$.
\begin{point}{40}{$f$ positive $\implies$ $\|f\|\leq 1$}%
This follows
immediately from~\sref{cp-russo-dye} and~\sref{cp-commutative},
but here's a concrete proof:
Let~$a\in \scrA$ be given.
Pick~$\lambda\in \C$ with $\left|\lambda\right|=1$
and~$\lambda f(a)\geq 0$.
Then $\left|f(a)\right|=f(\lambda a) 
= \Real{f(\lambda a)}
= f(\Real{(\lambda a)})
\leq f(\|a\|)=\|a\|$,
because $\Real{(\lambda a)}
\leq \|\Real{(\lambda a)}\|
\leq \|\lambda a\|=\|a\|$,
and $f$ is positive.
Hence~$\|f\|\leq 1$.
\end{point}
\begin{point}{50}{$\|f\|\leq 1$ $\implies$ $f$ is positive}%
Let~$a\in [0,1]_\scrA$ be given.
To prove that~$f$ is positive, it suffices to show that $f(a)\geq 0$.
Since $(\Real{f(a)})^\perp
=\Real{(f(a)^\perp)}
\leq \left|f(a)^\perp \right|
=\left|f(a^\perp)\right|\leq 1$,
and therefore $\Real{f(a)}\geq 0$,
we just need to show that~$\Imag{f(a)}=0$.

The trick is to consider $b_n := (a - \Real{f(a)})+ni\Imag{f(a)}$.
Indeed, since $(n+1)^2(\Imag{f(a)})^2
= \left|f(b_n)\right|^2 \leq \|b_n\|^2 = 
\|b_n^*b_n\| \leq \|a-\Real{f(a)}\|^2 + n^2(\Imag{f(a)})^2$,
one sees that $(2n+1)(\Imag{f(a)})^2\leq \|a-\Real{f(a)}\|^2$
for all~$n$,
which is impossible unless~$(\Imag{f(a)})^2=0$,
that is, $\Imag{f(a)}=0$.\qed
\end{point}
\end{point}
\end{point}
\begin{point}{60}[vn-ball-extreme-point]{Lemma}%
\index{unit ball!of a $C^*$-algebra!extreme points}%
An extreme point~$u$ of the unit ball~$(\scrA)_1$
of a $C^*$-algebra~$\scrA$
is a partial isometry with $(uu^*)^\perp\scrA(u^*u)^\perp = \{0\}$.
\begin{point}{70}{Remark}%
The converse (viz.~every such partial isometry
is extreme in~$(\scrA)_1$)
also holds, but we won't need it.
\end{point}
\begin{point}{80}{Proof}%
(Based on Theorem~7.3.1 of~\cite{kr}.)

To show~$u$ is a partial isometry
it suffices to prove that~$u^*u$ is a projection.
Suppose towards a contradiction that~$u^*u$ is not a projection.
Then $u^*u$,
represented
as continuous function (on~$\spec(u^*u)$ cf.~\sref{functional-calculus}),
takes neither the value~$0$ nor~$1$ on a neighbourhood
of some point,
and so by considering a
positive continuous function,
which is sufficiently small but non-zero on this neighbourhood
and zero elsewhere, 
we can find a non-zero element~$a$ 
of the (commutative) $C^*$-subalgebra generated by~$u^*u$
with $0\leq a\leq u^*u$
and $\|u^*u (1\pm a)^2\|\leq 1$,
so that~$\|u(1\pm a)\|\leq 1$.
Since~$u$ is extreme in~$(\scrA)_1$,
and~$u=\frac{1}{2}u(1+a)\,+\,\frac{1}{2}u(1-a)$,
we get~$ua=0$,
and so $0\leq a^2\leq \sqrt{a}u^*u\sqrt{a}=u^*ua=0$,
which contradicts $a\neq 0$.

Let~$a\in (uu^*)^\perp \scrA (u^*u)^\perp$
be given; we must show that $a=0$.
Assume (without loss of generality)
that~$\|a\|\leq 1$.
We'll show that $\|u\pm a\|\leq 1$,
because,
since~$u$ is extreme in~$(\scrA)_1$,
$u\equiv \frac{1}{2}(u+a)+\frac{1}{2}(u-a)$
implies that~$u=u+a$, and so $a=0$.
Note that $a^*a \leq (u^*u)^\perp$ (because $a(u^*u)^\perp=a$)
and $u^*a = 0$ (because $(uu^*)^\perp a=a$).
Thus $(u\pm a)^*(u\pm a)
=u^*u \pm u^* a \pm a^* u + a^*a
= u^*u + a^*a \leq u^*u + (u^*u)^\perp = 1$,
so $\|u\pm a\|\leq 1$. \qed
\end{point}
\end{point}
\begin{point}{90}[polar-decomposition-of-functional]%
	{Theorem (Polar decomposition of functionals)}%
\index{polar decomposition!of a functional}
Every functional $f\colon \scrA\to \C$ on a von Neumann algebra~$\scrA$
which is ultraweakly continuous on
the unit ball~$(\scrA)_1$
is of the form $f\equiv f(uu^*(\,\cdot\,)) = f((\,\cdot\,)u^*u)$
for some partial isometry~$u$ on~$\scrA$
such that $f(u(\,\cdot\,))$
and $f((\,\cdot\,)u)\colon \scrA\to\C$
are positive.
\begin{point}{100}{Proof}%
(Based on Theorem~7.3.2 of~\cite{kr}.) 
\begin{point}{110}%
We'll first show that~$f$ takes the value~$\|f\|$
at some extreme point~$u$ of~$(\scrA)_1$.
To begin, since~$(\scrA)_1$ is ultraweakly compact (\sref{vn-ball-compact}),
and~$f$ is ultraweakly continuous
the subset $\{\,f(a)\colon\,a\in(\scrA)_1\,\}$
of~$\R$ is compact,
and therefore has a largest element, 
which must be~$\|f\|$.
Thus the convex
set~$F:=\{\,a\in(\scrA)_1\colon\, f(a)=\|f\|\,\}$
is non-empty.
Since~$F$ is ultraweakly compact (being an ultraweakly closed
subset of the ultraweakly compact~$(\scrA)_1$),
$F$ has at least one extreme point
by the Krein--Milman Theorem
(see e.g.~Theorem~V7.4 of~\cite{conway2013}), say~$u$.
Note that~$F$ is a face of~$(\scrA)_1$:
if~$\frac{1}{2}a+\frac{1}{2}b\in F$ for some~$a,b\in(\scrA)_1$,
then $\frac{1}{2}f(a)+\frac{1}{2}f(b) = \|f\|$,
so~$f(a)=f(b)=\|f\|$
(since~$\|f\|$ is extreme in $(\C)_{\|f\|}$)
and thus~$a,b\in F$.
It follows that~$u$ is not only extreme in~$F$, but also in~$(\scrA)_1$,
so that~$u$ is an partial isometry with $(uu^*)^\perp\scrA(u^*u)^\perp=\{0\}$
by~\sref{vn-ball-extreme-point}.

Note that $f(u(\,\cdot\,))$
is positive by~\sref{positive-functional-criterion}, because
$\|f(u(\,\cdot\,))\||\leq\|f\|\|u\|\leq\|f\|=f(u)=f(u(1))$.
By a similar argument~$f((\,\cdot\,)u)$
is positive.

Let~$a\in\scrA$ be given.
It remains to be shown that $f(a)=f(uu^*a)=f(au^*u)$.
First note that  $u(u^*u)^\perp = 0$ (since~$u$ is an isometry)
and so $f(u(u^*u)^\perp)=0$,
that is,  $u^*u\geq \ceil{f(u(\,\cdot\,))}$.
This entails that $f(ubu^*u)=f(ub)$ for all~$b\in\scrA$
by~\sref{carrier-fundamental}, and in particular $f(uu^*au^*u)=f(uu^*a)$.

Now, since~$(uu^*)^\perp\scrA(u^*u)^\perp = \{0\}$,
we have $uu^* a u^*u + a = uu^*a + au^*u$,
and thus $f(a)+f(uu^*a)=f(a)+f(uu^*au^*u)=f(uu^*a)+f(au^*u)$,
which yields $f(a)=f(au^*u)$.
By a similar reasoning we get $f(uu^*a)=f(a)$.\qed
\end{point}
\end{point}
\end{point}
\begin{point}{120}[uwcont-on-ball]{Corollary}%
\index{functional!ultraweakly continuous}
A functional $f\colon \scrA\to\C$
on a von Neumann algebra~$\scrA$
is ultraweakly continuous
when it is ultraweakly continuous
on the unit ball~$(\scrA)_1$.
\begin{point}{130}{Proof}%
By~\sref{polar-decomposition-of-functional}
there is a partial isometry~$u$
such that $f(uu^*(\,\cdot\,))=f$
and~$f(u(\,\cdot\,))$ is positive.
Recall from~\sref{p-uwcont} that such a positive functional~$f(u(\,\cdot\,))$
is normal when it is ultraweakly continuous
on~$[0,1]_{\scrA}$;
which it is, 
because $a\mapsto ua$ is ultraweakly continuous (see \sref{mult-uws-cont}),
maps~$[0,1]_{\scrA}$ into~$(\scrA)_1$,
and $f$ is ultraweakly continuous on~$(\scrA)_1$.
But then~$f\equiv f(uu^*(\,\cdot\,))$
being the composition of the ultraweakly continuous maps
$f(u(\,\cdot\,))$ and $a\mapsto u^*a$ 
is ultraweakly continuous on~$\scrA$ too.\qed
\end{point}
\end{point}
\begin{point}{140}[functional-norm]{Lemma}%
Let~$f\colon \scrA\to\C$
be a normal functional
on a von Neumann algebra~$\scrA$,
and let~$u$ be a partial isometry in~$\scrA$
such that~$f(u(\,\cdot\,))$ is positive,
and $f=f(uu^*(\,\cdot\,))$.
Then~$\|f\|=f(u)$.\qed
\begin{point}{150}{Proof}%
Since~$f(u(\,\cdot\,))$ is positive,
we have $\|f(u(\,\cdot\,))\|=f(u)$
by~\sref{cp-russo-dye};
hence $\|f\|=\|f(uu^*(\,\cdot\,))\|
\leq \|f(u(\,\cdot\,))\|\|u^*\|
\equiv f(u) \leq \|f\|$,
and thus~$\|f\|=f(u)$.
\end{point}
\end{point}
\end{parsec}
\begin{parsec}{870}%
\begin{point}{10}{Definition}%
Given a von Neumann algebra~$\scrA$,
the vector space of ultraweakly continuous 
linear maps $f\colon \scrA\to \C$
endowed with the operator norm
is denoted by~$\Define{\scrA_*}$,%
\index{*understar@$\scrA_*$, predual of~$\scrA$}
and called the \Define{predual}%
\index{predual}
of~$\scrA$.
\begin{point}{20}{Remark}%
The reason
that the space~$\scrA_*$ is called the \emph{predual}
of~$\scrA$
is the non-trivial fact due to Sakai~\cite{sakai} (which we 
don't need and therefore won't prove),
that the obvious 
map $\scrA\to(\scrA_*)^*$,
where $(\scrA_*)^*$ is the \emph{dual} of~$\scrA_*$ --- the
vector space of bounded linear maps $\scrA_*\to\C$
endowed with the operator norm ---,
is a surjective isometry,
and so $\scrA$ ``is'' the dual of~$\scrA_*$,
(albeit only as normed space,
because~$(\scrA_*)^*$ doesn't come equipped with a multiplication.)

We will need this:
\end{point}
\end{point}
\begin{point}{30}[predual-complete]{Proposition}%
The predual~$\scrA_*$ of a von Neumann algebra~$\scrA$
is complete (with respect to the operator norm).
\begin{point}{40}{Proof}%
Let $f_1,f_2,\dotsc$ be a sequence in~$\scrA_*$
which is Cauchy with respect to the operator norm.
We already know (from~\sref{operator-norm-complete})
that $f_1,f_2,\dotsc$ converges to a bounded linear 
functional~$f\colon \scrA\to\C$;
so we only need to prove that~$f$ is ultraweakly continuous
to see that~$\scrA_*$ is complete,
and for this,
we only need to show 
(by~\sref{uwcont-on-ball})
that~$f$
is ultraweakly continuous
on the unit ball~$(\scrA)_1$ of~$\scrA$.
So let~$(b_\alpha)_\alpha$ be 
a net in~$(\scrA)_1$
which converges ultraweakly to~$0$;
we must show that~$\lim_\alpha f(b_\alpha)=0$.
Now, note that for every~$n$ and~$\alpha$
we have the bound
\begin{equation*}
	\left|f(b_\alpha)\right|
	\ \leq\ \left|(f-f_n)(b_\alpha)\right|
	\,+\,\left|f_n(b_\alpha)\right|
	\ \leq\ 
	\|f-f_n\|\,+\,\left|f_n(b_\alpha)\right|.
\end{equation*}
From this,
and $\lim_n\|f-f_n\|=0$,
and $\lim_\alpha f_n(b_\alpha)=0$
for all~$n$,
one easily deduces that~$\lim_\alpha f(b_\alpha)=0$.
Thus~$f$ is ultraweakly continuous,
and so~$\scrA_*$ is complete.\qed
\end{point}
\end{point}
\begin{point}{50}%
Note
that for a self-adjoint
element~$a$ of a von Neumann algebra~$\scrA$
we have $\|a\|=\sup_\omega \left|\omega(a)\right|$
where~$\omega$ ranges
over the npsu-functionals,
but that the same identity does not need to hold
for arbitrary (not necessarily self-adjoint) $a\in \scrA$.
The following lemma shows
that this restriction to self-adjoint elements
can be lifted 
by letting~$\omega$
range over all of~$\scrA_*$.
\end{point}
\begin{point}{60}[norm-predual]{Lemma}%
We have $\|a\|=\sup_{f\in(\scrA_*)_1}\left|f(a)\right|$
for every element~$a$ of a von Neumann algebra~$\scrA$.
\begin{point}{70}{Proof}%
It's clear that $\sup_{f\in(\scrA_*)_1}\left|f(a)\right|\leq \|a\|$.

For the other direction,
write $a\equiv [a] \sqrt{a^*a}$ (see~\sref{polar-decomposition})
and note that $\|a\|=\|\sqrt{a^*a}\|=\sup_{\omega\in \Omega} 
\left|\omega (\,\sqrt{a^*a}\,)\right|$,
where $\Omega$ is the set of  npu-maps $\scrA\to \C$
(which is order separating).
Let~$\omega\in \Omega$ be given.
Since~$[a]^*a=\sqrt{a^*a}$ 
we have $\omega(\,\sqrt{a^*a}\,)=\omega([a]^*a)=f(a)$,
where $f:=\omega([a]^*(\,\cdot\,))\in (\scrA_*)_1$,
and so $\|a\|=\sup_{\omega\in\Omega} 
\omega(\,\sqrt{a^*a}\,) \leq \sup_{f\in (\scrA_*)_1}
\left|f(a)\right|$.\qed
\end{point}
\end{point}
\begin{point}{80}[ultraweakly-bounded-implies-bounded]{Theorem}%
A net  $(b_\alpha)_\alpha$
 in a von Neumann algebra~$\scrA$
is norm bounded
(that is,~$\sup_\alpha \|b_\alpha\|<\infty$)
provided it is \Define{ultraweakly bounded}, i.e.,%
\index{ultraweakly bounded net}
$\sup_\alpha \left|\omega(b_\alpha)\right|<\infty$
	for every np(u)-map $\omega\colon \scrA\to \C$.
\begin{point}{90}{Proof}%
Note that $f\mapsto f(b_\alpha)$
gives a linear map $(\,\cdot\,)(b_\alpha)\colon \scrA_*\to\C$
with $\|(\,\cdot\,)(b_\alpha)\|=\|b_\alpha\|$ 
by~\sref{norm-predual}
for each~$\alpha$.
So to prove that~$(b_\alpha)_\alpha$
is norm bounded, viz.~%
$\sup_\alpha \|b_\alpha\|\equiv \sup_\alpha \|(\,\cdot\,)(b_\alpha)\|<\infty$,
it suffices to show
(by the principle of uniform boundedness, \sref{pub},
using that~$\scrA_*$ is complete, \sref{predual-complete}),
that $\sup_\alpha \left|f(b_\alpha)\right|<\infty$ for all~$f\in\scrA_*$.

Since such~$f\in\scrA_*$
can be written as $f\equiv \sum_{k=0}^3 i^k \omega_k$
where $\omega_k\colon \scrA\to\C$ are np-maps
(by~\sref{normal-functionals-lemma}),
we see that $\sup_\alpha\left|f(b_\alpha)\right|
\leq \sum_{k=0}^3\sup_\alpha \left|\omega_k(b_\alpha)\right|
<\infty$, because~$(b_\alpha)_\alpha$
is ultraweakly bounded.
Thus~$(b_\alpha)_\alpha$
is norm bounded.\qed
\end{point}
\end{point}
\end{parsec}
\subsection{Ultraweak Permanence}
\begin{parsec}{880}%
\begin{point}{10}%
We turn to 
a subtle, and surprisingly difficult matter:
it is not immediately clear
that the ultraweak topology on a von Neumann
subalgebra~$\scrA$ of a von Neumann algebra~$\scrB$,
coincides (on~$\scrA$) with the ultraweak topology on~$\scrB$.
While it is easily seen that the former is finer
(that is, a net in~$\scrA$ which converges ultraweakly in~$\scrA$,
converges ultraweakly in~$\scrB$ too, because any np-map 
$\omega\colon \scrB\to\C$ is also an np-map restricted to~$\scrA$),
it is not obvious that an np-map $\omega\colon \scrA\to\C$
can be extended to an np-map on~$\scrB$
--- but it can, as we'll see~\sref{functional-permanence}.
We'll call this independence of the ultraweak topology
from the surrounding space
\emph{ultraweak permanence}
being not unlike the independence
of the spectrum of an operator from the surrounding space
known as spectral permanence (\sref{spectral-permanence}).

It is tempting to think that the extension of an np-map~$\omega\colon \scrA
\to\C$
on a von Neumann
subalgebra~$\scrA$ of a von Neumann algebra~$\scrB$
to~$\scrB$ is simply a matter of applying Hahn--Banach to~$\omega$,
but this approach presents two problems:
it yields a normal but not necessarily positive extension
of~$\omega$;
and it not clear that~$\omega$ is ultraweakly continuous
on~$\scrA$ (that is, whether Hahn--Banach applies).

Instead of applying general techniques we feel forced
to delve deeper
into the particular structure 
provided to us by von Neumann algebras
(namely the commutant, \sref{commutant})
to show that any np-map
$\omega\colon \scrA\to\C$
on a von Neumann algebra~$\scrA$
of bounded operators on a Hilbert space~$\scrH$
can be extended to an np-map on~$\scrB(\scrH)$,
and in fact, is of the form
$\omega\equiv \sum_n \left<x_n,(\,\cdot\,)x_n\right>$
for some $x_1,x_2,\dotsc \in\scrH$, see~\sref{normal-functional}.
\end{point}
\begin{point}{20}[commutant-ceil]{Proposition}%
Let~$S$ be a subset of a von Neumann algebra~$\scrA$
that is closed under multiplication, involution, and contains~$1$.
Let~$e$ be a projection in~$\scrA$.
Then $\Define{\ceil{e}_{S^\square}}= 
\bigcup_{a\in S} \ceil{a^* e a}$%
\index{*ceils@$\ceil{e}_S$}
is the least projection in~$S^\square$
above~$e$.

(Compare this with the paragraph ``Subspaces'' of \S2.6 of~\cite{kr}.)
\begin{point}{30}{Proof}%
Let us first show that~$p:= \ceil{e}_{S^\square}$
is in~$S^\square$.
Let $b\in S$ be given;
we must show that $pb=bp$.
We may may assume without loss of generality that~$\|b\|\leq 1$.
Since~$b^*(\,\cdot\,)b\colon \scrA\to\scrA$
is normal and completely positive,
and $p=\bigcup_{a\in S} \ceil{a^*ea}$,
we have $b^*pb\leq \ceil{b^*pb} = 
\bigcup_{a\in S} \ceil{ b^* \ceil{a^* e a} b }
= \bigcup_{a\in S} \ceil{(ab)^* \,e\, ab} \leq p$
by~\sref{ncp-union} and~\sref{ncp-ceil}.
Applying $p^\perp(\,\cdot\,)p^\perp$,
we get $p^\perp b^*pb p^\perp 
\leq p^\perp p p^\perp = 0$,
so that $pbp^\perp=0$,
and thus $pbp=pb$.
Since similarly $pb^* =pb^*p$,
we get  $bp=pbp=pb$
(upon applying~$(\,\cdot\,)^*$) and so~$p\in S^\square$.

Note that~$e\leq \ceil{1^* e 1}\leq p$, because $1\in S$.
It remains to be shown that~$p$ is the least projection in~$S^\square$
above~$e$, so let~$q$ be a projection in~$S^\square$ above~$e$.
Since for~$a\in S$,
we have~$aq^\perp a^*=  q^\perp aa^* q^\perp 
\leq \|a\|^2q^\perp \leq \|a\|^2e^\perp$,
and so $a^*ea\leq \|a\|^2q$ 
we get~$\ceil{a^*ea}\leq q$
for all~$a\in S$,
and thus~$p=\bigcup_{a\in S}\ceil{a^*ea} \leq q$.\qed
\end{point}
\end{point}
\begin{point}{40}[carrier-vector-state]{Exercise}%
Show that given a vector~$x$ of Hilbert space~$\scrH$,
and a collection~$S$ of bounded operators on~$\scrH$
that is closed under addition, (scalar) multiplication,
involution, and contains the identity operator,
the following coincide.
\begin{enumerate}
\item
$\ceil{\,\ketbra{x}{x}\,}_{S^\square}$,
the least projection in~$S^\square$
above $\ceil{\,\ketbra{x}{x}\,}$;
\item
$\ceil{\,\left<x,(\,\cdot\,)x\right>|S^\square\,}$,
the carrier of the vector functional on~$S^\square$
given by~$x$;
\item
$\bigcup_{a\in S} \ceil{\,\ketbra{ax}{ax}\,}$; and 
\item
the projection on~$\overline{S x}$.
\end{enumerate}
Conclude that $\overline{S^{\square\square}x}
=\overline{S x}$.
(Hint: $S^{\square\square\square}=S^\square$.)
\begin{point}{50}[proto-double-commutant]%
Now consider
(instead of~$x$)
an np-map $\omega\colon \scrB(\scrH)\to \C$,
which we know must be of the form
$\omega\equiv \sum_n \left<x_n,(\,\cdot\,)x_n\right>$
(by~\sref{bh-np})
and is therefore given by 
an element $x'\equiv (x_1,x_2,\dotsc)$ of the $\N$-fold
product~$\scrH':=\bigoplus_n\scrH$ of~$\scrH$.
\begin{enumerate}
\item
Show that  $\omega(t)=\left<x',\varrho'(t)x'\right>$,
where $\varrho'\colon \scrB(\scrH)\to\scrB(\scrH')$
is the nmiu-map given by 
$\varrho'(t)y=(ty_n)_n$
for all~$t\in\scrB(\scrH)$ and~$y\in\scrH'$.

Prove that $\varrho'(t)=\sum_n P_n^* t P_n$,
where $P_n:=\pi_n\colon \scrH'\equiv \bigoplus_n\scrH\to\scrH$
is the $n$-th projection.
\item

Let~$t\in S^{\square\square}$ be given
(with~$S$ as above).
Show that $\varrho'(t)\in \varrho'(S)^{\square\square}$.

(Hint: first show 
	$P_n aP_m^*\in S^\square$ for all~$m$, $n$, and
		$a\in \varrho'(S)^\square$.)

Conclude that  $\varrho'(t)x'\,\in\,\overline{\varrho'(S)^{\square\square}x'}
\equiv \overline{\varrho'(S)x'}$.

Whence for every~$\varepsilon>0$
one can find~$a\in S$ with $\|t-a\|_\omega \leq \varepsilon$.

\item
Deduce that $S^{\square\square}$
is contained in the ultrastrong closure of~$S$.
\end{enumerate}
\spacingfix
\end{point}%
\end{point}%
\begin{point}{60}[double-commutant]{Double Commutant Theorem}%
\index{Double Commutant Theorem}%
\index{Bicommutant Theorem}
For a collection~$S$ of bounded operators
on a Hilbert space~$\scrH$
that is closed under addition, (scalar) multiplication,
involution, and contains the identity operator
the following are the same.
\begin{enumerate}
\item
$S^{\square\square}$, the ``double commutant'' of~$S$
in~$\scrB(\scrH)$;
\item
$\mathrm{us}\text{-}\mathrm{cl}(S)$,
the ultrastrong closure of~$S$ in~$\scrB(\scrH)$;
\item
$\mathrm{uw}\text{-}\mathrm{cl}(S)$,
the ultraweak closure of~$S$ in~$\scrB(\scrH)$;
\item
$W^*(S)$,
the least von Neumann subalgebra of~$\scrB(\scrH)$
that contains~$S$.
\end{enumerate}
\begin{point}{70}{Proof}%
(Based on Theorem~5.3.1 of~\cite{kr}.) 

Note that: $\mathrm{us}\text{-}\mathrm{cl}(S)
\subseteq  \mathrm{uw}\text{-}\mathrm{cl}(S)$,
because ultrastrong convergence implies ultraweak convergence;
and
$\mathrm{uw}\text{-}\mathrm{cl}(S)
\subseteq W^*(S)$,
because~$W^*(S)$ 
is ultraweakly closed in~$\scrB(\scrH)$ by~\sref{vnsac};
and
$W^*(S)\subseteq S^{\square\square}$,
because 
$S^{\square\square}$ is a von Neumann subalgebra
of~$\scrB(\scrH)$ by~\sref{commutant-basic};
and, finally, $S^{\square\square}\subseteq \mathrm{us}\text{-}\mathrm{cl}(S)$
by~\sref{proto-double-commutant}.\qed
\end{point}
\end{point}
\begin{point}{80}[centre-commutant]{Exercise}
Show that central elements of
a von Neumann algebra~$\scrA$
of bounded operators on a Hilbert space~$\scrH$
coincide with the central elements of the commutant~$\scrA^\square$,
that is, $Z(\scrA)=Z(\scrA^\square)$.
(Hint: $\scrA^{\square\square}=\scrA$ by~\sref{double-commutant}.)
\begin{point}{90}[commutant-cceil]%
Deduce that~$\cceil{f|\scrA}=\cceil{f|\scrA^\square}$
for every np-map~$f\colon \scrB(\scrH)\to\scrB$ 
into a von Neumann algebra~$\scrB$.
\end{point}
\end{point}
\end{parsec}%
\begin{parsec}{890}%
\begin{point}{10}[gns-mapping-property]{Lemma}%
Let~$\omega\colon \scrA\to\C$
be an np-map on a von Neumann algebra~$\scrA$,
which is represented by nmiu-maps
$\varrho\colon \scrA\to\scrB(\scrH)$
and $\pi\colon \scrA\to\scrB(\scrK)$
on Hilbert spaces~$\scrH$ and~$\scrK$.
If $\left<x,\varrho(\,\cdot\,)x\right>
=\omega=\left<y,\pi(\,\cdot\,)y\right>$
for some  $x\in\scrH$ and~$y\in\scrK$,
then there is a bounded operator $U\colon \scrK\to\scrH$
for which~$UU^*$ is the projection
on~$\overline{\varrho(\scrA)x}$,
$U^*U$ is the projection
on~$\overline{\pi(\scrA)y}$,
and~$U\pi(a)=\varrho(a)U$
for all~$a\in\scrA$.
\begin{point}{20}{Proof}%
(Compare this with Proposition 4.5.3 of~\cite{kr}.)

Since $\|\varrho(a)x\|^2
= \left<x,\varrho(a^*a)x\right>
=\omega(a^*a)=\left<y,\pi(a^*a)y\right>
= \|\pi(a)y\|^2$ for all~$a\in\scrA$,
there is a unique bounded operator $V\colon \overline{\pi(\scrA)y}
\to \overline{\varrho(\scrA)x}$
with $V\pi(a)y = \varrho(a)x$ for all~$a\in\scrA$.
A moment's thought reveals
	that $V$ is a unitary (and so~$V^*V=1$ and $VV^*=1$.)
Now, define $U:=EVF^*$
where $E\colon \overline{\varrho(\scrA)x}\to\scrH$
and $F\colon \overline{\pi(\scrA)y}\to\scrK$
are the inclusions
(and so~$E^*E=1$ and~$F^*F=1$).
Then $UU^*= EVF^*FV^*E^*=EVV^*E^*=EE^*$
is the projection onto~$\overline{\varrho(\scrA)x}$,
and $UU^*=FF^*$
is the projection onto~$\overline{\pi(\scrA)y}$.

Let~$a\in\scrA$ be given.
It remains to be shown that
$U\pi(a)=\varrho(a)U$.
To this end,
observe that
$V F^* \pi(a) F = E^* \varrho(a) E V$
(because 
these two bounded linear maps
are easily seen to
agree on the dense subset $\pi(\scrA)y$
of $\overline{\pi(\scrA)y}$);
and $\varrho(a)E = EE^*\varrho(a)E$
(because $\varrho(a)$ maps $\varrho(\scrA)x$ into~$\varrho(\scrA)x$);
and similarly $\varrho(a^*)F=FF^*\varrho(a^*) F$,
so that $F^*\varrho(a) = F^* \varrho(a) FF^*$
(upon application of the~$(\,\cdot\,)^*$).
By these observations,  $U\pi(a)=
EVF^*\pi(a)=EVF^*\pi(a) FF^*
= EE^*\varrho(a)EVF^*
= \varrho(a)EVF^*
= \varrho(a)U$.\qed
\end{point}
\end{point}
\begin{point}{30}[summing-partial-isometries]{Exercise}%
It is not too difficult 
to see that the (ultraweak) sum~$\sum_i u_i$
of a collection $(u_i)_i$ 
of partial isometries from some von Neumann algebra
is again a partial isometry, 
provided that the initial projections $u_i^*u_i$
are pairwise orthogonal,
and the final projections~$u_iu_i^*$ are pairwise orthogonal.
In this exercise, you'll establish a similar result,
but for partial isometries between two different Hilbert spaces,
and avoiding the use of an analogue of 
the ultraweak topology for such operators.
\begin{point}{40}%
Let~$\scrH$ and~$\scrK$ be Hilbert spaces,
and 
let $U_i\colon \scrH\to\scrK$
be a bounded operator
for every element~$i$ from some set~$I$.
Assume that the operators~$U_i^*U_i$
are pairwise orthogonal projections in~$\scrB(\scrK)$,
and that~$U_iU_i^*$ are pairwise orthogonal projections in~$\scrB(\scrH)$.
\begin{enumerate}
\item
Let~$x\in\scrH$ and~$y\in\scrK$ be given.

Show that
$\left|\left<x,U_iy\right>\right|
\leq \|U_i^*x\| \|U_i y \|$
for each~$i$
(perhaps by first proving that  $U_i = U_i U_i^* U_i$).

Show that~$\sum_i \|U_i y\|^2 \leq \|y\|^2$
and $\sum_i \|U_i^* x\|^2\leq \|x\|^2$,
and deduce from this
that $\sum_i \left|\left<x,U_i y\right>\right| \leq \|x\|\|y\|$.

Now use~\sref{chilb-form-representation}
to show that there is a bounded operator~$U\colon \scrK\to\scrH$
with $\left<x,Uy\right>
= \sum_i \left<x,U_iy\right>$
for all~$x\in\scrH$ and~$y\in \scrK$.
\item
Show that~$U_i^*U_j = 0$ when~$i\neq j$.
Deduce from this that~$U^*U = \sum_i U_i^* U_i$.

Prove that $UU^* = \sum_i U_iU_i^*$.
\end{enumerate}
\spacingfix%
\end{point}%
\end{point}%
\begin{point}{50}[sigma-weak-lemma-2]{Lemma}%
Let~$\Omega$ be a collection of np-maps $\omega\colon \scrA\to\C$
on a von Neumann algebra~$\scrA$ 
whose central carriers,  $\cceil{\omega}$, are pairwise orthogonal to one another,
and let~$\scrH$ and~$\scrK$ be Hilbert spaces
on which~$\scrA$ is represented such
that each~$\omega\in\Omega$ 
is given by vectors $x_\omega\in\scrH$ and $y_\omega\in\scrK$,
that is,
$\left<x_\omega,\varrho(\,\cdot\,)x_\omega\right>
=\omega = \left<y_\omega,\pi(\,\cdot\,)y_\omega\right>$,
where $\varrho\colon \scrA\to\scrB(\scrH)$
and $\pi\colon \scrA\to\scrB(\scrK)$
are nmiu-maps.

Then there is a bounded operator $U\colon \scrK\to\scrH$
which intertwines~$\pi$ and~$\varrho$
in the sense that $U \pi(a)=\varrho(a) U$
for all~$a\in \scrA$
such that~$U^*U$
is a projection in
$\pi(\scrA)^\square$
with
$\cceil{U^*U}_{\pi(\scrA)^\square}=\pi(\sum_\omega \cceil{\omega})$,
and 
$UU^*$
is projection in
$\varrho(\scrA)^\square$
with
$\cceil{UU^*}_{\varrho(\scrA)^\square}=\varrho(\sum_\omega \cceil{\omega})$.
\begin{point}{60}{Proof}%
Given~$\omega\in\Omega$,
let $\sigma_\omega\colon \varrho(\scrA)\to \C$
and $\sigma_\omega'\colon \varrho(\scrA)^\square\to\C$
denote the restrictions
of the vector functional $\left<x_\omega,(\,\cdot\,)x_\omega\right>
\colon \scrB(\scrH)\to\C$,
and let  $\tau_\omega\colon \pi(\scrA)\to \C$
and $\tau_\omega'\colon \pi(\scrA)^\square\to\C$
be similar restrictions
of $\left<y_\omega,(\,\cdot\,)y_\omega\right>$.
We already know (by~\sref{gns-mapping-property}
and~\sref{carrier-vector-state})
that
there is
a bounded operator $U_\omega\colon \scrK\to\scrH$
with $U_\omega^*U_\omega = \ceil{\tau_\omega'}$,
$U_\omega U_\omega^* = \ceil{\sigma_\omega'}$,
and $U_\omega \pi(a) = \varrho(a) U_\omega$
for all~$a\in\scrA$.

We'll combine these $U_\omega$s into one operator~$U$
using~\sref{summing-partial-isometries}, but for this we must
verify that the projections $U_\omega U_\omega^*=\ceil{\sigma_\omega'}$
are pairwise orthogonal,
and that the projections $U_\omega^* U_\omega$
are pairwise orthogonal too.
To this end note that
$\cceil{\sigma_\omega}=\cceil{\sigma_\omega'}$
by~\sref{commutant-cceil}.
Thus, since the projections~$\cceil{\omega}$
are orthogonal to one another,
and $\ceil{\sigma_\omega'}\leq \cceil{\sigma_\omega'}
= \cceil{\sigma_\omega} =\varrho(\cceil{\omega})$,
we see that the projections~$U_\omega U^*_\omega \equiv 
\ceil{\sigma_\omega'}$
are indeed pairwise orthogonal.
Since for a similar reason
the projections $U^*_\omega U_\omega
\equiv \ceil{\tau_\omega'}$ are pairwise orthogonal too,
there is by~\sref{summing-partial-isometries}
a bounded operator $U\colon \scrK\to\scrH$
with 
$U^*U = \sum_\omega U_\omega^* U_\omega$,
$UU^* = \sum_\omega U_\omega U_\omega^*$,
and $\left<x,Uy\right>=\sum_\omega \left<x,U_\omega y\right>$
for all~$x\in\scrH$ and~$y\in \scrK$.

Let us check that~$U$ has the desired properties.
To begin, since the projections
$\cceil{U_\omega U^*_\omega}=\cceil{\sigma_\omega'}=\varrho(\cceil{\omega})$
are pairwise orthogonal,
we have $\cceil{U U^*}=\sum_\omega \cceil{U_\omega U^*_\omega}
=  \varrho(\sum_\omega \cceil{\omega})$
by~\sref{cceil-basic} and
\sref{sum-of-orthogonal-projections}.
Similarly, $\cceil{U^*U}= \pi(\sum_\omega \cceil{\omega})$.

Finally,
given~$a\in\scrA$
we have $U\pi(a)=\varrho(a)U$,
because $\left<x,U\pi(a)y\right>
= \sum_\omega \left<x,U_\omega \pi(a)y\right>
= \sum_\omega \left<x,\varrho(a) U_\omega y\right>
= \sum_\omega \left<\varrho(a)^* x, U_\omega y\right>
= \left<\varrho(a)^* x, U y\right>
= \left<x, \varrho(a) U y\right>$
for all~$x\in\scrH$ and~$y\in\scrK$.\qed
\end{point}
\end{point}
\begin{point}{70}[sigma-weak-lemma]{Corollary}%
Let~$\scrA$ be a von Neumann 
of bounded operators on some Hilbert space~$\scrH$,
and let~$\varrho\colon \scrA\to\scrB(\scrH)$
denote the inclusion.
Let~$\Omega$ be the collection of all np-maps $\scrA\to\C$,
and let $\varrho_\Omega\colon \scrA\to\scrB(\scrH_\Omega)$
be as in~\sref{gelfand-naimark-representation}.

There is a bounded operator $U\colon \scrH_\Omega\to\scrH$
such that $U^*U$ is a projection
in~$\varrho_\Omega(\scrA)^\square$
with 
$\cceil{U^*U}_{\varrho_\Omega(\scrA)^\square}=1$
and
$U\varrho_\Omega(a)= \varrho(a) U$
for all~$a\in\scrA$.
\begin{point}{80}{Proof}%
Let $\{x_i\}_i$ be a maximal set of vectors
in~$\scrH$ 
such that the 
central carriers~$\cceil{\omega_i}$
of the corresponding vector functionals
$\omega_i :=\left<x_i,\varrho(\,\cdot\,)x_i\right>$
on~$\scrA$
are pairwise orthogonal;
so that we'll 
have $\sum_i \cceil{\omega_i}=1$.
Now, the point of~$\scrH_\Omega$ 
is that there are vectors $y_i\in\scrH_\Omega$
with $\omega_i=\left<y_i,\varrho_\Omega(\,\cdot\,)y_i\right>$
for each~$i$.
Now apply~\sref{sigma-weak-lemma-2}
to get a map~$U\colon \scrH_\Omega\to \scrH$
with the desired properties.\qed
\end{point}
\end{point}
\begin{point}{90}[normal-functional]{Theorem}%
\index{functional!positive!normal}%
\index{normal!positive functional}%
Every np-map $\omega\colon \scrA\to\C$
on a von Neumann subalgebra~$\scrA$
of~$\scrB(\scrH)$,
where~$\scrH$ is some Hilbert space,
is of the 
form~$\omega \equiv \sum_n \left<x_n,(\,\cdot\,)x_n\right>$
for some $x_1,x_2,\dotsc \in\scrH$
(with $\sum_n\|x_n\|^2<\infty$).
\begin{point}{100}{Proof}%
(Based on Theorem~7.1.8 of~\cite{kr}.)

Let~$\varrho\colon \scrA\to\scrB(\scrH)$
denote the inclusion,
and let~$U\colon \scrH_\Omega\to\scrH$
be as in~\sref{sigma-weak-lemma}.
Since~$\omega\in\Omega$,
there is~$y\in\scrH_\Omega$
with $\omega = \left<y,\varrho_\Omega(\,\cdot\,)y\right>$.
We're going to `transfer'~$y$ from~$\scrH_\Omega$ to~$\scrH$
using the following device.
Since~$1=\cceil{U^*U}_{\varrho_\Omega(\scrA)^\square}$,
we can 
(by~\sref{cceil-sum})
find partial isometries $(v_i)_i$ in
$\varrho_\Omega(\scrA)^\square$
with $1=\sum_i v_i^*v_i$
and $v_iv_i^*\leq U^*U$
for all~$i$.
Then for every~$a\in\scrA$,
\begin{alignat*}{3}
	\textstyle
\omega(a)\  &=\
	\textstyle
	\left<\,y, \,\varrho_\Omega(a)y\,\right> 
	\\
&= 
	\textstyle
	\ \sum_i \left<\,y,\,  v_i^*v_i\, \varrho_\Omega(a) y \,\right>
	&&\qquad 
	\textstyle
	\text{since $1=\sum_i v_i^*v_i$}
	\\
&= 
	\textstyle
	\ \sum_i \left<\,y,\,  v_i^* U^*Uv_i \,\varrho_\Omega(a)y\, \right>
	&&\qquad\text{since $v_iv_i^*\leq U^*U$}
	\\
	&= 
	\textstyle
	\ \sum_i \left<\,Uv_i y,\, U \varrho_\Omega(a)v_i y\, \right>
	&&\qquad\text{since $v_i \in \varrho(\scrA)^\square$}
	\\
	&= 
	\textstyle
	\ \sum_i \left<\,Uv_i y,\, \varrho(a) \, Uv_i y\, \right>
	&&\qquad\text{since $U\varrho_\Omega(a)=\varrho(a) U$}.
\end{alignat*}
In particular, $\omega(1)=\sum_i \|Uv_i y\|^2$,
so at most countably many $Uv_iy$'s are non-zero;
and denoting those by~$x_1,x_2,\dotsc$,
we get $\omega = \sum_n \left<x_n,(\,\cdot\,)x_n\right>$.\qed
\end{point}
\end{point}
\begin{point}{110}[functional-permanence]{Corollary}%
Let~$\scrA$ be a von Neumann subalgebra
of a von Neumann algebra~$\scrB$.
\begin{enumerate}
\item
For every np-map $\omega\colon \scrA\to\C$
there is an np-map $\xi\colon \scrB\to\C$
with $\xi|\scrA=\omega$.
\item%
\index{ultraweak and ultrastrong!permanence}
\Define{Ultraweak permanence:}\ 
the restriction of the ultraweak topology on~$\scrB$
to~$\scrA$ coincides with the ultraweak topology on~$\scrA$.
\item
\Define{Ultrastrong permanence:}\ 
the restriction of the ultrastrong topology on~$\scrB$
to~$\scrA$ coincides with the ultrastrong topology on~$\scrA$.
\end{enumerate}%
\spacingfix%
\end{point}%
\begin{point}{120}[functional-extension]{Exercise}%
Let $\varrho\colon \scrA\to\scrB$
be an injective nmiu-map.

Show 
using~\sref{injective-nmiu-iso-on-image}
that any np-functional
$\omega\colon \scrA\to\C$
can be extended along~$\varrho$,
that is,
there is an np-functional
$\omega'\colon \scrB\to\C$
with $\varrho\circ \omega' = \omega$.
\end{point}
\end{parsec}
\begin{parsec}{900}%
\begin{point}{10}%
We end the chapter 
with another corollary to~\sref{normal-functional}:
 that  the np-functionals
on a von Neumann algebra
are generated (in a certain sense)
by any centre separating collection
of functionals. This fact plays
an important role
in the next chapter
for our definition of the tensor product
of von Neumann algebras
(on which the product functionals are to be centre separating,
\sref{tensor}).
\end{point}
\begin{point}{20}[vn-center-separating-fundamental]{Proposition}%
\index{centre separating collection!of np-functionals}
Given a centre separating collection~$\Omega$ of np-functionals
on a von Neumann algebra~$\scrA$,
and an ultrastrongly dense subset~$S$ of~$\scrA$
\begin{enumerate}
\item
	\label{vn-center-separating-fundamental-1}
$\Omega':= \{\,\omega(s^*(\,\cdot\,)s)\colon\,
\omega\in\Omega,\,s\in S\,\}$
is order separating, and
\item
	\label{vn-center-separating-fundamental-2}
$\Omega'':=\{\,\sum_n\omega_n\colon\, \omega_1,\dotsc,\omega_N\in\Omega'\,\}$ is operator norm dense in $(\scrA_*)_+$.
\end{enumerate}
\spacingfix%
\begin{point}{30}[vn-center-separating-fundamental-i]{Proof}%
We tackle~\ref{vn-center-separating-fundamental-1}
first. We already know from~\sref{proto-gelfand-naimark}
that the collection $\Xi:=\{\, \omega(a^*(\,\cdot\,)a)\colon\, 
\omega\in\Omega,\, a\in\scrA\,\}$,
which contains~$\Omega'$, is order separating;
so to prove that~$\Omega'$ is itself order separating
it suffices by~\sref{order-separating-dense-subset} to show that~$\Omega'$
is norm dense in~$\Xi$.
This is indeed the case
since given~$a\in\scrA$
and~$\omega\in\Omega$,
and a net~$(s_\alpha)_\alpha$ in~$S$
that converges ultrastrongly to~$a$,
the functionals~$s_\alpha \ast\omega
\equiv \omega(s_\alpha^*(\,\cdot\,)s_\alpha)$
converge in norm to~$a\ast \omega$
as $\alpha\to\infty$
by~\sref{bstaromega-basic}.%
\begin{point}{40}{Concerning~\ref{vn-center-separating-fundamental-2}}%
Let~$f\colon \scrA\to\C$ be an np-map;
we must show that~$f$ is in 
the norm closure~$\overline{\Omega''}$ of~$\Omega''$.
Note that since~$\Omega$ is centre separating,
the map $\varrho_\Omega\colon \scrA\to\scrB(\scrH_\Omega)$
from~\sref{proto-gelfand-naimark}
is injective,
and in fact restricts
to an nmiu-isomorphism
from~$\scrA$ onto $\varrho_\Omega(\scrA)$
(cf.~\sref{ngns}).
So by~\sref{normal-functional}
$f$ is of the form $f\equiv 
\sum_n \left<x_n,\varrho_\Omega(\,\cdot\,)x_n\right>$
for some $x_1,x_2,\dotsc \in\scrH_\Omega$
with $\sum_n \|x_n\|^2<\infty$,
so that the partial sums
$\sum_{n=1}^N \left<x_n,\varrho_\Omega(\,\cdot\,)x_n\right>$
converge with respect to the operator norm to~$f$
(by~\sref{vector-functional-convergence}).
Thus to show that~$f$ is in~$\overline{\Omega''}$
it suffices to show that each~$\left<x_n,\varrho_\Omega(\,\cdot\,)x_n\right>$
is in~$\overline{\Omega''}$
(since~$\overline{\Omega''}$ is clearly closed
under finite sums and norm limits).
In effect
we may assume 
without loss of generality
that~$f\equiv \left<x,\varrho_\Omega(\,\cdot\,)x\right>$
for some~$x\in \scrH_\Omega$.
We reduce the problem some more.
By definition 
of $\scrH_\omega\equiv \bigoplus_{\omega\in\Omega}\scrH_\omega$
and~$\varrho_\Omega$,
we have $f=\left<x,\varrho_\Omega(\,\cdot\,)x\right>
= \sum_{\omega\in \Omega} 
\left<x_\omega,\varrho_\omega(\,\cdot\,)x_\omega\right>$;
and so we may, by the same token, assume without loss of generality
that $f=\left<x,\varrho_\omega(\,\cdot\,)x\right>$
for some~$\omega\in \Omega$
and~$x\in \scrH_\omega$.
Since such~$x$ (by definition of~$\scrH_\omega$, 
\sref{gns})
is the norm limit of
a sequence $\eta_\omega(a_1),\,\eta_\omega(a_2),\,
\dotsb$,
where~$a_1,a_2,\dotsc\in\scrA$,
the np-maps $a_n\ast \omega\equiv 
\left<\eta_\omega(a_n),\varrho_\omega(\,\cdot\,)\eta_\omega(a_n)\right>$
converge to~$\left<x,\varrho_\omega(\,\cdot\,) x\right>=f$
in the operator norm as $n\to\infty$ by~\sref{vector-functional-convergence};
and so we may assume without loss of generality
that~$f=a\ast \omega$ for some~$a\in\scrA$ and~$\omega\in \Omega$.
Since~$S$ is ultrastrongly dense in~$\scrA$
we can find a net~$(s_\alpha)_\alpha$
in~$S$ that converges ultrastrongly to~$a$.
As the np-functionals~$s_\alpha \ast \omega$
in $\Omega'\subseteq \Omega''$
will then
operator-norm converge
to~$f= a\ast \omega$
as $\alpha\to\infty$
by~\sref{bstaromega-basic},
we conclude that~$f\in \overline{\Omega''}$.\qed
\end{point}
\end{point}
\end{point}
\end{parsec}
\begin{parsec}{910}
\begin{point}{10}
With this chapter
ends perhaps the most hairy part of this thesis:
we've developed the theory of von Neumann algebras
starting from Kadison's characterisation (see~\sref{vna})
to the point that we have a sufficiently firm hold on the normal functionals
(see e.g.~\sref{polar-decomposition-of-functional}, \sref{normal-functional}),
the ultraweak and ultrastrong topologies 
(e.g.~\sref{kaplansky}, \sref{functional-permanence}, 
	\sref{vn-center-separating-fundamental}),
the projections (\sref{vna-ceil}, \sref{ceill}, 
\sref{projections-norm-dense}), 
and the division structure 
(\sref{douglas}, \sref{polar-decomposition})
on a von Neumann algebra.
In the next chapter we reap the benefits of our labour
when we study an assortment of structures
in the category~$\W{cpsu}$ of von Neumann algebras
and ncpsu-maps.
\end{point}
\end{parsec}

\chapter{Assorted Structure in~$\W{cpsu}$}
\begin{parsec}{920}%
\begin{point}{10}%
In the previous two chapters
we have travelled through
charted territory
when developing the theory of $C^*$-algebras
and von Neumann algebras
adding some new landmarks and shortcuts
of our own along the way.
In this chapter
we properly break new ground
by revealing
two entirely new features
of the category~$\W{cpsu}$
of von Neumann algebras
and the normal completely positive sub-unital
linear maps between them,
namely, 
\begin{enumerate}
\item 
that the binary operation~$\ast$
on the effects of a von Neumann algebra~$\scrA$
given by~$p\ast q = \sqrt{p}q\sqrt{p}$
(representing measurement of~$p$)
can be axiomatised,
and 
\item 
	that the category~$\W{cpsu}$
has all the bits and pieces
needed to be a model of Selinger and Valiron's quantum lambda calculus.
\end{enumerate}
We'll deal with the first matter directly after this introduction
in Section~\ref{S:measurement}.
The second matter is treated in Section~\ref{S:model},
but only after
we have given the tensor product of von Neumann 
algebras a complete overhaul
in Section~\ref{S:tensor}.
Finally, as an offshoot of our model
of the quantum lambda calculus
we'll study all
 von Neumann algebras
that admit a `duplicator'
in Section~\ref{S:duplicable}
---
surprisingly, they're
all of the form $\ell^\infty(X)$.
\end{point}
\end{parsec}

\section{Measurement}
\label{S:measurement}
\begin{parsec}{930}%
\begin{point}{10}%
The maps on a von Neumann algebra~$\scrA$
of the form
$a\mapsto \sqrt{p}a\sqrt{p}\colon\,\scrA\to\scrA$,
where~$p$ is an effect of~$\scrA$,
represent measurement of~$p$,
and are called \emph{assert maps} in~\cite{newdirections}.
The importance of these maps 
to any logical description of
quantum computation is not easily overstated.
On the effects of~$\scrA$
these maps are also  studied 
in the guise
of the binary operation
$p\ast q=\sqrt{p} q \sqrt{p}$
called the \emph{sequential product}
(see e.g.~\cite{gudder2002sequential}).
We'll axiomatise this operation
in this section
in terms of the properties 
of the underlying assert maps.

Our first observation
to this end
is that any assert map factors as
\begin{equation*}
\xymatrix@C=10em{
\scrA
\ar[r]^-{\pi\colon a\mapsto \ceil{p}a\ceil{p}}
&
\ceil{p}\!\scrA\!\ceil{p}
\ar[r]^-{c\colon a\mapsto \sqrt{p}a\sqrt{p}}
&
\scrA
},
\end{equation*}
where both~$\pi$ and~$c$ obey a universal property:
$c$ is a \emph{filter} of~$p$, see~\sref{filter},
and~$\pi$ is a \emph{corner} of~$\ceil{p}$, see~\sref{corner}.
Such maps
that are the composition of a filter and a corner
will be called \emph{pure}, see~\sref{pure},
Since not only assert maps turn out to be pure, but also maps of the form
$b^*(\,\cdot\,)b\colon \scrA\to\scrA$ for an arbitrary element~$b$
of~$\scrA$,
we need another property of assert maps, namely
that
\begin{equation*}
	\sqrt{p}\,e_1\,\sqrt{p}\ \leq\  e_2^\perp 
	\qquad\iff\qquad
	\sqrt{p}\,e_2 \,\sqrt{p}\ \leq\  e_1^\perp
\end{equation*}
for all projections~$e_1$ and~$e_2$ of~$\scrA$---which we'll
describe by saying that 
	\begin{equation*} \sqrt{p}(\,\cdot\,)\sqrt{p}\colon \scrA\to\scrA
	\end{equation*}
is \emph{$\diamond$-self-adjoint}.
Judging only by the name
it may not surprise you that the map $b(\,\cdot\,)b\colon \scrA\to\scrA$
where~$b\in \scrA$ is self-adjoint (but not necessarily positive)
turns out to be $\diamond$-self-adjoint too,
so that as a final touch we introduce the notion
of \emph{$\diamond$-positive} maps $f\colon \scrA\to\scrA$
that are simply maps of the form~$f\equiv gg$ for some $\diamond$-self-adjoint~$g$.

The main technical result, then, of this section
is that any $\diamond$-positive map $f\colon\scrA\to\scrA$
is of the form~$f=\sqrt{p}(\,\cdot\,)\sqrt{p}$
where~$p=f(1)$;
and, accordingly, our axioms 
(in~\sref{uniqueness-sequential-product})
that uniquely
determine the sequential product~$\ast$
on the effects of a von Neumann algebra~$\scrA$ are:
for every effect~$p$ of~$\scrA$,
\begin{enumerate}
\item
$p\ast 1=p$,
\item
$p\ast q = f(q)$
for all~$q\in [0,1]_\scrA$
for some pure map~$f\colon \scrA\to\scrA$,
\item
$p=q\ast q$ for some $q$ from~$[0,1]_\scrA$,
\item
$p \ast (p \ast q) = (p\ast p)\ast q$
for all~$q\in[0,1]_\scrA$,
\item
$p \ast e_1 \leq e_2^\perp\iff
p \ast e_2 \leq e_1^\perp$
for all projections $e_1,e_2$ of~$\scrA$.
\end{enumerate}%
While I would certainly not like
to undersell the results mentioned above,
I suspect that the notion of purity exposed along the way
might turn out to be of far greater significance
for the following reason.
Our notion of purity can be described in
wildly different terms:
a map~$f\colon \scrA\to\scrB$ is pure when given its
\emph{Paschke dilation}
$\xymatrix{\scrA
	\ar[r]|-\varrho
&
	\scrP\ar[r]|-c
&
\scrB}$
the map $\varrho$ is surjective
(see~\sref{paschke-pure} and~\cite{wwpaschke}).
Because of my faith in our notion of purity I've allowed myself
to address some theoretical questions concerning it
here that are not required for the main results of this thesis,
but suppose a general interest in purity:
I'll show that every pure map~$f\colon\scrA\to\scrB$
is extreme among the ncp-maps~$g\colon \scrA\to\scrB$ with~$f(1)=g(1)$,
and, in fact, enjoys the possibly stronger property 
of being~\emph{rigid} (see~\sref{rigid} and~\sref{pure-is-rigid}).
\end{point}
\end{parsec}
\subsection{Corner and Filter}
\begin{parsec}{940}%
\begin{point}{10}{Definition}%
Given a projection~$e$ of a von Neumann algebra~$\scrA$,
the \Define{corner}%
	\index{corner (von Neumann algebra)}
of~$e$
is the subset~$e\scrA e$%
\index{*eAe@$e\scrA e$, corner}
of~$\scrA$ 
(consisting of the elements of~$\scrA$
of the form~$eae$ with~$a\in\scrA$).
In this context,
the obvious map~$e\scrA e\to\scrA$
is called the \Define{inclusion}%
\index{inclusion!of a corner}
and the map $a\mapsto eae,\ \scrA\to e\scrA e$
is called the \Define{projection}.%
\index{projection!onto a corner}
\end{point}
\begin{point}{20}[corner-vna-basic]{Exercise}%
Let~$e$ be a projection from a von Neumann algebra~$\scrA$.
\begin{enumerate}
\item
Show that~$a\in\scrA$ 
is an element of~$e\scrA e$ iff~$eae=a$
iff $\ceilr{a}\cup\ceill{a} \leq e$.
\item
Show that the corner~$e\scrA e$
is closed under addition, (scalar) multiplication,
and involution.
\item
Show that~$e$ is a unit for~$e\scrA e$,
that is, $ea=ae=a$ for all~$a\in e\scrA e$.
\item
Show that~$e\scrA e$ is norm and ultraweakly closed.\\
(Hint: use the fact that $e(\,\cdot\,)e\colon \scrA\to\scrA$
is normal and bounded.)
\item
Show that~$e\scrA e$ --- 
endowed with the addition, (scalar) multiplication,
involution and norm from~$\scrA$,
and with~$e$ as its unit ---  is a $C^*$-algebra.
\item
Show that the supremum of a bounded directed
set~$D$ of self-adjoint elements of~$e\scrA e$
computed in~$\scrA$
is itself in~$e\scrA e$,
and, in fact, the supremum of~$D$ in~$e\scrA e$.
\item
Show that the inclusion $e\scrA e\to\scrA$
is an ncpsu-map.
\item
Deduce from this that the restriction of an np-map
$\omega\colon \scrA\to\C$ to
a map $e\scrA e\to\C$
is an np-map.

Conclude that~$e\scrA e$ is a von Neumann algebra.
\item
Show that the projection $a\mapsto eae,\ \scrA\to e\scrA e$
is an ncpu-map.
\item
Show that every np-map $\omega\colon e\scrA e\to\C$
is the restriction
of the np-map $\omega(e(\,\cdot\,)e)\colon \scrA\to\C$.
Deduce from this that the ultraweak topology of~$e\scrA e$
coincides (on $e\scrA e$) with the ultraweak topology on~$\scrA$.
Show that the ultrastrong topologies on~$e\scrA e$ and~$\scrA$
coincide in a similar fashion.
\end{enumerate}
\spacingfix
\end{point}%
\begin{point}{30}[ad-ncp]{Exercise}%
Let~$a$ be an element of a von Neumann algebra~$\scrA$,
and let~$p$ and~$q$ be projections
of~$\scrA$ with $a^*pa\leq q$.
\begin{enumerate}
\item
Show that $a^*ba\in q\scrA q$
for every~$b\in p\scrA p$.
\item
Show that~$a^*(\,\cdot\,)a$
gives an ncp-map $p\scrA p\to q\scrA q$.
\end{enumerate}
\spacingfix%
\end{point}%
\end{parsec}%
\begin{parsec}{950}%
\begin{point}{10}[corner]{Definition}%
Let~$p$ be an effect of a von Neumann algebra~$\scrA$.
A \Define{corner}%
	\index{corner (map)}
of~$p$ is an
ncp-map $\pi\colon \scrA\to\scrC$
to a von Neumann algebra~$\scrC$
with~$\pi(p^\perp)=0$,
which is initial among such maps 
in the sense
that every ncp-map $f\colon \scrA\to\scrB$
with~$f(p^\perp)=0$
factors as $f=g\circ\pi$
for some unique ncp-map $g\colon \scrC\to\scrB$.

While most corners
that we'll deal with are unital,
there are also corners which are not unital
(because there are non-unital
ncp-isomorphisms).
When we write ``corner'' we shall
always mean a ``unital corner''%
\index{corner (map)!unital}
unless explicitly stated otherwise.
\end{point}
\begin{point}{20}[prop-corner]{Proposition}%
Given an effect~$p$ of a von Neumann algebra~$\scrA$,
and a partial isometry~$u$ of~$\scrA$
with $\floor{p}=uu^*$,
the map $\pi\colon \scrA\to u^*u \scrA u^*u$
given by~$\pi(a)=u^*au$ is a corner of~$p$.
\begin{point}{30}{Proof}%
By~\sref{ad-ncp}, $\pi$ is an ncp-map.
To see that~$\pi(p^\perp)\equiv u^*p^\perp u =0$,
note that since~$u^*u=u^*\,u u^*\,u$,
we have $0=u^*(uu^*)^\perp u =u^*\smash{\floor{p}}^\perp u
= u^*\ceil{\smash{p^\perp}} u$,
and so
$0 = \ceil{u^* \ceil{\smash{p^\perp}} u }
=\ceil{u^* p^\perp u}$
by~\sref{ceil-fundamental},
giving~$u^*p^\perp u=0$
by~\sref{ceil-basic}.

Let~$\scrB$ be a von Neumann algebra,
and let~$f\colon \scrA\to\scrB$ be an ncp-map
with $f(p^\perp)=0$.
To show that~$\pi$ is a corner,
we must show that there is a unique ncp-map
$g\colon u^*u \scrA u^*u\to\scrB$
with $f=g\circ \pi$.
Uniqueness follows
from surjectivity of~$\pi$.
Concerning existence,
define~$g:= f\circ \zeta$,
where $\zeta\colon  u^*u\scrA u^*u\to \scrA$
is the ncp-map given by~$\zeta(a)=uau^*$
for~$a\in\scrA$ (see~\sref{ad-ncp}),
so that it is immediately clear that~$g$ is an ncp-map.
It remains to be shown~$f=g\circ \pi$,
that is,
$f(a)=f(uu^*\,a\,uu^*)$ for all~$a\in\scrA$.
This follows from~\sref{cp-comprehension}
because~$f(\smash{(uu^*)^\perp})=0$,
since~$\ceil{\smash{f(\,\smash{(uu^*)^\perp}\,)}}
=\ceil{\smash{f(\smash{\floor{p}}^\perp)}}
=\ceil{\smash{f(\ceil{\smash{p^\perp}})}}
= \ceil{\smash{f(p^\perp)}}=\ceil{0}=0$.\qed
\end{point}
\end{point}
\end{parsec}
\begin{parsec}{960}%
\begin{point}{10}[filter]{Definition}%
A \Define{filter}%
\index{filter}
is an ncp-map $c\colon \scrC\to\scrA$
between von Neumann algebras
such that every ncp-map $f\colon \scrB\to\scrA$
with~$f(1)\leq c(1)$
factors as $f=c\circ g$
for some unique ncp-map $g \colon \scrB\to\scrC$.
We'll say that~$c$ is a \Define{filter for}~$c(1)$.
\index{filter!for~$p$}
\begin{point}{11}{Remark}%
In the abstract setting of effectus theory,
it makes sense to call these filters
``quotients'', as we do in~\cite{cho2015quotient};
but since in the concrete setting of von Neumann algebras
``quotient'' has a pre-existing and unrelated meaning,
we chose to use the word ``filter'' instead (as in ``polarising filter''),
an idea borrowed from~\cite{wilce2016royal}.
\end{point}
\end{point}
\begin{point}{20}%
To show that there is a filter
for every positive element of a von Neumann algebra
we need the following result
concerning ultraweak limits of ncp-maps.
\end{point}
\begin{point}{30}[ncp-uwlim]{Lemma}%
Given von Neumann algebras~$\scrA$
and~$\scrB$
the pointwise ultraweak limit
$f\colon \scrA\to\scrB$
of a net of  positive linear maps $f_\alpha\colon \scrA\to\scrB$
is positive, and, 
\begin{enumerate}
\item
$f$ is completely positive provided
that the $f_\alpha$ are completely positive, and
\item
$f$ is normal provided that the $f_\alpha$ are normal
and the ultraweak convergence of the~$f_\alpha$ to~$f$
is uniform on~$[0,1]_\scrA$.
\end{enumerate}
\spacingfix%
\begin{point}{40}{Proof}%
Since given~$a\in \scrA$ the element~$f(a)$
is the ultraweak limit of the positive elements~$f_\alpha(a)$,
and therefore positive (by~\sref{vn-positive-basic}),
we see that~$f$ is positive.

Suppose that each~$f_\alpha$ is completely positive.
To show that~$f$ is completely positive,
we must prove, given~$a_1,\dotsc,a_n\in\scrA$
and~$b_1,\dotsc,b_n\in\scrB$,
that 
the element $\sum_{i,j} b_i^* f(a_i^*a_j)b_j$
of~$\scrB$
is positive.
And indeed it is,
being the ultraweak limit of
the positive elements $\sum_{i,j} b_i^* f_\alpha (a_i^* a_j)b_j$,
because  $f_\alpha(a_i^* a_j)$
converges ultraweakly to~$f(a_i^* a_j)$,
and~$b_i^*(\,\cdot\,)b_j\colon \scrB\to\scrB$
is ultraweakly continuous
(\sref{mult-uws-cont})
for any~$i$ and~$j$.

If the~$f_\alpha$ 
are normal,
and converge uniformly on~$[0,1]_\scrA$ ultraweakly
to~$f$,
then~$f$ is ultraweakly continuous
on~$[0,1]_\scrA$
(because the uniform limit of continuous functions is continuous),
and thus normal (by~\sref{p-uwcont}).\qed
\end{point}
\end{point}
\begin{point}{50}[canonical-filter]{Proposition}%
Given an element~$d$ of a von Neumann algebra~$\scrA$,
the map $c\colon \ceilr{d}\!\scrA\!\ceilr{d}\to\scrA$
given by~$c(a)=d^*ad$
is a filter.
\begin{point}{60}{Proof}%
Note that~$c$ is an ncp-map by~\sref{ad-ncp}.
Let~$\scrB$ be a von Neumann algebra,
and let~$f\colon \scrB\to\scrA$ be an ncp-map
with $f(1)\leq c(1)$.
To show that~$c$ is a filter,
we must show that there is a unique ncp-map
$g\colon \scrB\to
\ceilr{d}\!\scrA\!\ceilr{d}$
with~$f=c\circ g$.
Uniqueness of~$g$ follows from the observation
that~$c$ is injective by~\sref{mult-cancellation}.

To establish the existence of such~$g$,
note that~$f(b)$ is an element of~$d^*\scrA d$,
when~$b$ is positive
by~\sref{sequential-douglas}
because~$0\leq f(b)\leq \|b\|f(1)\leq \|b\| c(1)=\|b\|d^*d$,
and thus for arbitrary~$b\in\scrB$ too
(being a linear combination
of positive elements).
We can thus define $g\colon \scrB\to \ceilr{d}\!\scrA\!\ceilr{d}$
by~$g(b)=d^*\backslash f(b)/d$
for all~$b\in\scrB$.
It is clear that~$g$ is linear and positive,
and~$c\circ g=f$.

To see that~$g$ is normal,
note that
$d^*\backslash\,\cdot\,/d\colon
d^*(\scrA)_1 d\to\scrA$
is ultrastrongly continuous by~\sref{div-usc},
as is~$f$ by~\sref{cp-uscont}
(also) as map from~$(\scrB)_1$ to~$d^*(\scrA)_1 d$,
so that~$g$ is ultrastrongly continuous on~$(\scrB)_1$,
and therefore normal by~\sref{p-uwcont}.

Finally, $g$ is completely positive
by~\sref{ncp-uwlim},
because it is by~\sref{div-approx}
the uniform ultrastrong limit
of the by~\sref{ad-ncp} completely positive maps
$(\sum_{n=1}^Nt_n)^* \,f(\,\cdot\,)\,(\sum_{n=1}^N t_n)$,
where~$t_1,t_2,\dotsc$
is an approximate pseudoinverse of~$d$.\qed
\end{point}
\end{point}
\end{parsec}
\begin{parsec}{970}
\begin{point}{10}
Before exploring their more technical aspects,
we'll explain how
corners and filters can be made to appear at opposite ends
of a chain of adjunctions:
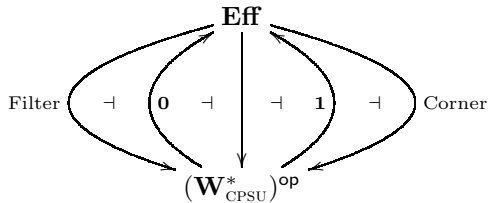
\begin{equation*} 
\xymatrix@R=5em{ 
	\Cat{Eff}\ar[d]
	\ar@/_3ex/@{{}{ }{}}[d]|\dashv
	\ar@/^3ex/@{{}{ }{}}[d]|\dashv
	\ar@/_11.5ex/@{{}{ }{}}[d]|\dashv
	\ar@/^11.5ex/@{{}{ }{}}[d]|\dashv
	\ar@/_15ex/[d]_{\text{Filter}}
	\ar@/^15ex/[d]^{\text{Corner}}
	\\ 
	\op{(\W{cpsu})}\ar@/^8ex/[u]_{\mathbf{0}}\ar@/_8ex/[u]^{\mathbf{1}} 
}  
\end{equation*} 
The category $\Define{\Cat{Eff}}$%
\index{Eff@$\Cat{Eff}$}%
\index{quotient--comprehension chain}
has as objects pairs $(\scrA, p)$,
where~$\scrA$ is a von Neumann algebra, and~$p\in[0,1]_\scrA$
is an effect from~$\scrA$.
A morphism  $(\scrA,p)\longrightarrow (\scrB,q)$
in~$\Cat{Eff}$
is an ncpsu-map $f\colon \scrB\to \scrA$ with $p\leq f(q)+f(1)^\perp$
--- that is,
\begin{equation*}
	\omega(p)\ \leq\ \omega(f(q))\,+\,\omega(f(1))^\perp
	\qquad\text{for every normal state } \omega\colon \scrA\to \C.
\end{equation*}
The functor $\Cat{Eff}\longrightarrow \op{(\W{cpsu})}$
in the middle of the diagram above 
	maps a morphism $f\colon (\scrA,p)\to(\scrB,q)$
to the underlying map $f\colon \scrB\to\scrA$.
The functors~$\mathbf{0}$ and~$\mathbf{1}$
on its sides map a von Neumann algebra~$\scrA$ to
$(\scrA,0)$ and~$(\scrA,1)$, respectively,
and send an ncpsu-map $f\colon \scrA\to\scrB$ to itself; 
this is possible since
\begin{equation*}
0\,\leq\,  f(0) + f(1)^\perp
	\qquad\text{and}\qquad
1\,\leq\, f(1)+f(1)^\perp.
\end{equation*}
That~$\mathbf{1}$
is right adjoint to the functor $\Cat{Eff}\longrightarrow \op{(\W{cpsu})}$
follows from the observation that
an ncpsu-map $f\colon \scrB\to\scrA$
is always a morphism $(\scrA,p)\to(\scrB,1)$,
whatever $p\in[0,1]_\scrA$ may be,
because  $p \leq f(1)+f(1)^\perp$.
For a similar reason  $\mathbf{0}$ is left adjoint to
$\Cat{Eff}\longrightarrow \op{(\W{cpsu})}$.

On the other hand,
a morphism $(\scrA,1)\to(\scrB,q)$ where $q\in[0,1]_\scrB$
is not just any ncpsu-map $f\colon \scrB\to\scrA$,
but one
with $1\leq f(q)+f(1)^\perp$,
that is, $f(q^\perp)=0$.
It's no surprise then that
a corner $\pi\colon \scrB\to \scrC$ for~$q\in[0,1]_\scrB$
considered as morphism $(\scrC,1)\to(\scrB,q)$
is a universal arrow from~$\mathbf{1}$ to~$(\scrB,q)$.

On the other side there's a twist:
a morphism $(\scrA,p)\to(\scrB,0)$ where~$p\in[0,1]_\scrA$
is an ncpsu-map $f\colon \scrA\to\scrB$
with $p\leq f(0)+f(1)^\perp$,
that is, $f(1)\leq p^\perp$.
It follows that any filter $c\colon \scrC\to \scrA$ for~$p^\perp$,
when considered as morphism $(\scrA,p)\to (\scrC,0)$,
is a universal arrow from~$(\scrA,p)$ to~$\mathbf{0}$.

This chain of adjunctions not only exposes
a hidden symmetry between filters and corners,
but such chains appear
in many other categories as well, see~\cite{cho2015quotient}.
\end{point}
\end{parsec}
\begin{parsec}{980}%
\begin{point}{10}[dfn-standard-corner-and-filter]{Definition}%
Let~$\scrA$ be a von Neumann algebra.
\begin{enumerate}
\item
Given a positive element~$p$
of~$\scrA$
we denote
by $\Define{c_p}\colon \ceil{p}\!\scrA\!\ceil{p}\to\scrA$%
\index{cp@$c_p$, standard filter for~$p$}
the \Define{standard filter}%
\index{filter!standard}
for~$p$
given by~$c_p(a)=\sqrt{p}a\sqrt{p}$
for all~$a\in\ceil{p}\!\scrA\!\ceil{p}$.
\item
Given an effect~$p$ of~$\scrA$
we denote
by $\Define{\pi_p}\colon \scrA\to\floor{p}\!\scrA\!\floor{p}$%
		\index{pip@$\pi_p$, standard corner of~$p$}
the \Define{standard corner}%
\index{corner (map)!standard}
		of~$p$
given by~$\pi_p(a)=\floor{p}\!a\!\floor{p}$.
\end{enumerate}
\spacingfix%
\end{point}%
\begin{point}{20}[filter-basic]{Exercise}%
Let~$c\colon \scrC\to\scrA$ be a filter,
where~$\scrC$ and~$\scrA$ are von Neumann algebras.
\begin{enumerate}
\item
Show that, writing~$p:=c(1)$,
there is a unique
ncp-map $\alpha \colon \scrC\to \ceil{p}\!\scrA\!\ceil{p}$
with $c = c_p \circ \alpha$;
and that this~$\alpha$ is a unital ncp-isomorphism.
\item
Show that~$c$ is injective
(by proving first that~$c_p$ is injective
using~\sref{mult-cancellation}).

Conclude that~$c$
is faithful (so $\ceil{f}=1$), and that~$c$ is mono in~$\W{CP}$.
\item
Show that~$c$ is bipositive
(by proving first that~$c_p$
is bipositive using~\sref{sequential-douglas}).
\end{enumerate}
\spacingfix%
\end{point}%
\begin{point}{30}[filters-composition]{Exercise}%
Show that the composition~$d\circ c$
of filters~$c\colon\scrC\to\scrD$
and~$d\colon \scrD\to\scrA$ 
between von Neumann algebras
is again a filter.
\end{point}
\begin{point}{40}[corner-basic]{Exercise}%
Let~$p$ be an effect of a von Neumann algebra~$\scrA$,
and let~$\pi\colon \scrA\to\scrC$ be a corner of~$p$.
\begin{enumerate}
\item
Show that there is a unique ncp-map
$\beta \colon \floor{p}\!\scrA\!\floor{p}\to\scrC$
with~$\pi = \beta\circ \pi_p$;
and that this~$\beta$ is unital and an ncp-isomorphism.
\item
Show that~$\pi$ is surjective, and that~$\pi$ is epi in~$\W{cp}$.
\end{enumerate}
\spacingfix%
\end{point}%
\begin{point}{50}[corners-floor]{Exercise}%
Show that an ncpu-map $\pi\colon \scrA\to\scrB$
between von Neumann algebras
is a corner for an effect~$p$ of~$\scrA$
iff~$\pi$ is a corner for~$\floor{p}$;
in which case~$\ceil{\pi}=\floor{p}$.

Thus a corner~$\pi$ is a corner for~$\ceil{\pi}$.
\end{point}
\begin{point}{60}[corners-composition]{Exercise}%
Show that the composition~$\tau\circ \pi$
of corners~$\pi\colon \scrA\to\scrB$
and~$\tau\colon \scrB\to\scrC$
between von Neumann algebras
is again a corner.\\
(Hint:
prove
and use the inequality
$\ceil{\tau}\leq \ceil{\smash{\pi(\ceil{\tau\circ \pi}^\perp)}}^\perp$.)
\end{point}
\begin{point}{70}[filter-corner]{Theorem}%
Given an ncp-map $f\colon\scrA\to\scrB$
between von Neumann algebras,
a projection~$e$ of~$\scrA$
with~$\ceil{f}\leq e$,
and a positive element~$p$
of~$\scrB$ with~$f(1) \leq p$,
there is a unique ncp-map
$g \colon e\scrA e
\to \ceil{p}\!\scrB\!\ceil{p}$
such that
\begin{equation*}
\xymatrix@C=4em{
\scrA
\ar[r]^-f
\ar[d]_{\pi_e}
&
\scrB
\\
e\scrA e
\ar[r]_-g
& 
\ceil{p}\!\scrB\!\ceil{p}
\ar[u]_{c_p}
}
\end{equation*}
commutes,
and it is given by
$g(a)=\sqrt{p}\backslash f(a)/\!\sqrt{p}$
for all~$a\in e\scrA e$.
\begin{point}{80}{Proof}%
Uniqueness of~$g$ follows from the facts
that~$\pi_e$ is epi and~$c_p$ is mono
in~$\W{cp}$,
see~\sref{corner-basic} and~\sref{filter-basic}.

Concerning existence, 
since~$\pi_e$ is a corner of~$e$,~\sref{corner},
and~$\ceil{f}\leq e$,
or in other words, $f(e^\perp)=0$,
there is a unique ncp-map $h\colon e\scrA e\to \scrB$
with $h \circ \pi_e = f$.
Note that~$h(a)=f(a)$ for all~$a$ from~$e\scrA e$.

As~$h(1)=h(\pi_e(1))=f(1)\leq p=c_p(1)$,
and~$c_p$ is a filter,~\sref{filter},
there is a unique ncp-map
$g\colon e\scrA e \to p \scrB p$
with $c_p\circ g = h$,
which is (by the proof of \sref{canonical-filter}) given by
$g(a)=\sqrt{p}\backslash h(a)/\sqrt{p}
\equiv \sqrt{p}\backslash f(a)/\sqrt{p}$
for all~$a$ from~$e\scrA e$.
Then~$c_p\circ g\circ \pi_e = h\circ \pi_e = f$.\qed
\end{point}
\end{point}
\begin{point}{90}[square-f]{Corollary}%
Given an ncp-map $f\colon \scrA\to\scrB$
between von Neumann algebras,
there is a unique ncp-map $\Define{[f]}\colon 
\ceil{f}\!\scrA\!\ceil{f}
\to
\ceil{f(1)}\!\scrB\!\ceil{f(1)}$%
\index{*brackets@$[\,\cdot\,]$!$[f]$, for an ncp-map}
such that 
\begin{equation*}
\xymatrix@C=4em{
\scrA
\ar[r]^-f
\ar[d]_{\pi_{\ceil{f}}}
&
\scrB
\\
\ceil{f}\!\scrA\!\ceil{f}
\ar[r]_-{[f]}
& 
\ceil{f(1)}\!\scrB\!\ceil{f(1)}
\ar[u]_{c_{f(1)}}
}
\end{equation*}
commutes;
and it is given by~$[f](a)=\sqrt{f(1)}\backslash f(a)/\!\sqrt{f(1)}$
for all~$a$ from $\ceil{f}\!\scrA\!\ceil{f}$.

Moreover, 
$[f]$ is unital and faithful.
\end{point}
\begin{point}{100}{Example}%
For any faithful unital ncp-map $f\colon \scrA\to \scrB$
we have~$[f]=f$.
Such a map need not be an isomorphism;
	as one may take $f\colon (\lambda,\mu)\mapsto \frac{1}{2}(\lambda+\mu),
\C^2\to\C$.
\end{point}
\begin{point}{110}[ad-pure]{Example}%
In the concrete case
that $f\equiv a^*(\,\cdot\,)a \colon
s\scrA s\to t\scrA t$,
where~$a$ is an element
of a von Neumann algebra,
and $s$ and~$t$ are projections of~$\scrA$
with
$\ceilr{a}\leq s$
and~$\ceill{a}\leq t$,
the map~$[f]$ 
is closely related to the
polar decomposition $a\equiv [a]\sqrt{a^*a}
= \sqrt{aa^*}[a]$ of~$a$,
where $[a]=a/\sqrt{a^*a}$
(see~\sref{polar-decomposition}).

Indeed,
since  $\ceil{f}=\ceilr{a}$,
$f(1)=a^*a$,
and~$[f]\equiv \sqrt{a^*a}\backslash a^*(\,\cdot\,)a/\sqrt{a^*a}
\equiv [a](\,\cdot\,)[a]^*$,
the picture becomes:
\begin{equation*}
\xymatrix@C=10em{
s\scrA s
\ar[r]^-{f\,=\,a^*\,(\,\cdot\,)\,a}
\ar[d]_{\pi_{\ceilr{a}}}
&
t\scrA t
\\
\ceilr{a}\!\scrA\!\ceilr{a}
\ar[r]_-{[f] \,=\,  [a]\,(\,\cdot\,)\,[a]^*}
& 
\ceill{a}\!\scrA\!\ceill{a}
\ar[u]_{c_{a^*a}}
}
\end{equation*}
Note that in this example
$[f]$ is an ncpu-isomorphism,
because~$[a]$ is a partial isometry
with initial projection~$\ceill{a}$
and final projection~$\ceilr{a}$.
Thus one can think of the diagram above
as an isomorphism theorem of sorts,
which applies only to certain  ncp-maps
that'll be called \emph{pure} in a moment (see~\sref{pure-fundamental}).
\end{point}
\end{parsec}
\subsection{Isomorphism}
\begin{parsec}{990}%
\begin{point}{10}%
In case you were wondering,
the ncpu-isomorphism
we encountered in~\sref{ad-pure}
is simply an nmiu-isomorphism 
(see~\sref{iso}), which follows
from the following characterisation of multiplicativity.
\end{point}
\begin{point}{20}[gardner]{Proposition}%
\index{multiplicative!ncpsu-map}
For an ncpu-map $f\colon \scrA\to\scrB$
between von Neumann algebras
the following are equivalent.
\begin{enumerate}
\item
\label{gardner-1}
$f$ is multiplicative.
\item
\label{gardner-2}
$f(a)f(b)=0$
for all $a,b\in\scrA$ with $ab=0$.
\item
\label{gardner-3}
$\ceil{f(p)}\ceil{f(q)}=0$
for all projections $p$ and~$q$ of~$\scrA$ with $pq=0$.
\item
\label{gardner-4}
$f$ maps projections to projections.
\item
\label{gardner-5}
$\ceil{f(a)}=f(\ceil{a})$
for all~$a\in\scrA_+$.
\end{enumerate}
\spacingfix%
\begin{point}{30}{Proof}%
(Based in part on the work of Gardner in~\cite{gardner}).
\begin{point}{40}{\sref{gardner-1}$\Longrightarrow$\sref{gardner-4}
 and \sref{gardner-5}$\Longrightarrow$\sref{gardner-4}}
	are rather obvious.
\end{point}
\begin{point}{50}{\sref{gardner-4}$\Longrightarrow$\sref{gardner-5}}%
 $\ceil{f(a)}
\smash{\overset{\sref{ncp-ceil}}{=\joinrel=\joinrel=}}
\ceil{f(\ceil{a})}
=f(\ceil{a})$
since~$f(\ceil{a})$ is a projection.
\end{point}
\begin{point}{60}{\sref{gardner-4}$\Longrightarrow$\sref{gardner-3}}%
Let~$p$ and~$q$ be projections of~$\scrA$ with~$pq=0$.
Then~$p\leq q^\perp$, and so~$f(p)\leq f(q^\perp)=f(q)^\perp$,
which implies that $\ceil{f(p)}\ceil{f(q)}
=f(p)f(q)=0$ since~$f(p)$ and~$f(q)$ are projections.
\end{point}
\begin{point}{70}{\sref{gardner-3}$\Longrightarrow$\sref{gardner-2}}%
Let~$a,b\in\scrA$ with~$ab=0$ be given.
We must show that~$f(a)f(b)=0$,
and for this it suffices to show that
$\ceill{f(a)}\ceilr{f(b)}=0$,
because $f(a)f(b)=f(a)\ceill{f(a)}\ceilr{f(b)}f(b)$.
Since~$ab=0$,
we have~$\ceill{a}\ceilr{b}=0$ by~\sref{mult-cancellation},
and so~$\ceil{f(\ceill{a})}\ceil{f(\ceilr{a})}=0$.
Now,
since $\ceil{f(\ceill{a})}\leq \ceill{f(a)}$
	and $\ceil{f(\ceilr{a}}\leq \ceilr{f(a)}$
	by~\sref{ncp-ceill},
we get $\ceill{f(a)}\ceilr{f(b)}
	= \ceill{f(a)} \ceil{f(\ceill{a})}
	\ceil{f(\ceilr{a})}
	\ceilr{f(a)}
	=0$.
\end{point}
\begin{point}{80}{\sref{gardner-2}$\Longrightarrow$\sref{gardner-1}}%
We must show that~$f(a)f(b)=f(ab)$
for all~$a,b\in \scrA$.
Since the linear span of projections is norm-dense in~$\scrA$,
it suffices to show that $f(a)f(e)=f(ae)$
for any $a\in\scrA$ and a projection~$e$ of~$\scrA$.
Given such~$a$ and~$e$,
we on the one hand have $ae^\perp\, e=0$,
so that~$f(ae^\perp)f(e)=0$,
that is, $f(a)f(e)=f(ae)f(e)$;
and on the other hand
we have $ae\,e^\perp=0$,
so that~$f(ae)f(e^\perp)=0$,
that is, $f(ae)=f(ae)f(e)$;
so that we reach~$f(ae)=f(a)f(e)$ as sum total,
and the result that~$f$ is multiplicative.\qed
\end{point}
\end{point}
\end{point}
\begin{point}{90}[iso]{Theorem}%
An ncpsu-isomorphism $f\colon \scrA\to\scrB$
between von Neumann algebras 
(so both~$f$ and~$f^{-1}$ are ncpsu-maps)
is an nmiu-isomorphism.
\begin{point}{100}{Proof}%
Since~$f^{-1}(1)\leq 1$
and so~$1=f(f^{-1}(1))\leq f(1)\leq 1$,
we see that~$f(1)=1$, so both $f$ and $f^{-1}$ are unital.
It remains to be shown that~$f$ and~$f^{-1}$ are multiplicative.
Since by~\sref{projection-order-sharp} an effect~$a$ of~$\scrA$
is a projection iff~$0$ is the infimum of~$a$ and~$a^\perp$,
	and~$f$ (as ncpu-isomorphism) preserves $(\,\cdot\,)^\perp$
	and order,
we see that~$f$ maps projections to projections,
and is thus multiplicative, by~\sref{gardner}.
It follows automatically that~$f^{-1}$ is multiplicative too.\qed
\end{point}
\end{point}
\begin{point}{110}{Exercise}%
Show that any filter of a projection is multiplicative.\\
(Hint: the filter is
a standard filter
up to an
ncpu-isomorphism, \sref{filter-basic},
which is an nmiu-isomorphism by~\sref{iso}.)
\end{point}
\begin{point}{120}[sharp-multiplicative]{Exercise}%
\index{multiplicative!ncp-map}
Show that for an ncp-map $f\colon \scrA\to\scrB$
between von Neumann algebras
the following are equivalent.
\begin{enumerate}
\item
$f$ is multiplicative.
\item
$f$ sends projections to projections.
\item
$\ceil{f(a)}=f(\ceil{a})$
for all~$a \in\scrA_+$.
\end{enumerate}
(Hint: factor~$f=\zeta \circ h$
where~$\zeta$ is a filter for~$f(1)$
and~$h$ is an ncp-map.)
\end{point}
\end{parsec}
\subsection{Purity}
\begin{parsec}{1000}%
\begin{point}{10}[pure]{Definition}%
Filters, corners,
and their compositions we'll call \Define{pure}.%
\index{pure map}
\end{point}
\begin{point}{20}{Exercise}%
Show that the following maps are pure.
\begin{enumerate}%
\item
An ncp-isomorphism between von Neumann algebras.
\item
The identity map~$\id\colon \scrA\to\scrA$
on a von Neumann algebra~$\scrA$.
\item
The map $a^*\,(\,\cdot\,)\,a\colon \scrA\to\scrA$
for any element~$a$ of a von Neumann algebra~$\scrA$.
\end{enumerate}
\spacingfix%
\end{point}%
\begin{point}{30}[pure-fundamental]{Proposition}%
For an ncp-map $f\colon \scrA\to\scrB$ between von Neumann algebras
the following are equivalent.
\begin{enumerate}
\item 
\label{pure-fundamental-1}
	$f$ is pure, i.e., $f$ is the composition
	of (perhaps many) filters and corners.
\item
\label{pure-fundamental-2}
	$f = c\circ \pi$ for a filter $c\colon \scrC\to\scrB$
	and a corner $\pi\colon \scrA\to\scrC$.
\item
\label{pure-fundamental-3}
	$[f]$ from~\sref{square-f} is an ncpu-isomorphism.
\end{enumerate}
\spacingfix%
\begin{point}{40}{Proof}%
\ref{pure-fundamental-3}$\Longrightarrow$\ref{pure-fundamental-2}
and \ref{pure-fundamental-2}$\Longrightarrow$\ref{pure-fundamental-1}
are rather obvious.
\begin{point}{50}{\ref{pure-fundamental-1}$\Longrightarrow$%
\ref{pure-fundamental-2}}%
Calling $f$ \emph{properly pure}
when~$f\equiv c\circ \pi$
for some filter~$c$ and corner~$\pi$,
we must show that every pure map is properly pure.
For this it suffices to show that the composition of properly
pure maps is again properly pure;
which,
since filters are closed under composition
(by~\sref{filters-composition}),
and corners are closed under composition
(by~\sref{corners-composition}),
boils down to proving that the composition
$\pi\circ c$ of a corner~$\pi$ and a filter~$c$
is properly pure.
Since~$\pi\equiv \alpha\circ \pi_{\ceil{\pi}}$
and~$c\equiv c_{c(1)}\circ \beta$
for ncpu-isomorphisms~$\alpha$ and~$\beta$
(see~\sref{filter-basic}
and~\sref{corner-basic})
it suffices to show that
$f:=\pi_{s} c_{p}$ is properly pure
for a positive element~$p$ and a projection~$s$
of a von Neumann algebra~$\scrA$.
Since such~$f$ is of the form $f=s\sqrt{p}(\,\cdot\,)\sqrt{p}s
\colon \ceil{p}\!\scrA\!\ceil{p}\to s\scrA s$,
we know by~\sref{ad-pure}
that~$[f]$ is an ncpu-isomorphism,
and thus that~$f\equiv c_{f(1)}\circ [f]\circ \pi_{\ceil{f}}$ is properly pure.
\end{point}
\begin{point}{60}{\ref{pure-fundamental-2}$\Longrightarrow$%
\ref{pure-fundamental-3}}%
Recall that $[f]$
is by definition the unique ncp-map
with~$f = c_{f(1)} [f] \pi_{\ceil{f}}$,
see~\sref{square-f}.
Note that since~$f=c\circ \pi$,
	we have~$\ceil{f}=\ceil{\pi}$ (because~$\ceil{c}=1$ 
	by~\sref{filter-basic}),
and~$f(1)=c(1)$ (because~$\pi(1)=1$).
Since there are ncpu-isomorphisms~$\alpha$ and~$\beta$
with $\pi= \alpha \pi_{\ceil{\pi}}$ and  $c=c_{c(1)} \beta$,
we see that~$f=c_{c(1)} \beta\alpha \pi_{\ceil{\pi}}$,
and so~$[f]=\beta\alpha$
by definition of~$[f]$,
so~$[f]$ is an ncpu-isomorphism.\qed
\end{point}
\end{point}
\end{point}
\begin{point}{70}[special-pure-maps]{Exercise}%
Use~\sref{pure-fundamental} to show that 
\begin{enumerate}
\item
a faithful pure map is a filter,
\item
a unital pure map is a corner, and
\item
a unital and faithful pure map is an ncpu-isomorphism.
\end{enumerate}
\spacingfix
\end{point}%
\end{parsec}%
\subsection{Contraposition}
\begin{parsec}{1010}%
\begin{point}{10}{Definition}%
Given an ncp-map $f\colon \scrA\to\scrB$
between von Neumann algebras
we define
$\Define{f^\diamond}\colon \Proj(\scrA)\to \Proj(\scrB)$%
\index{*diamond@$f^\diamond$}
by~$f^\diamond(e)=\ceil{f(e)}$
for all~$e\in\Proj(\scrA)$.
\begin{point}{11}{Remark}%
The significance of the symbol~``$\diamond$''
in~$f^\diamond$
is in 
accommodating the notation~$f^\BOX(e):=f^\diamond(e^\perp)^\perp$
used
    in the next thesis, in~\sref{diamond-basics}.
\end{point}
\end{point}
\begin{point}{20}{Proposition}%
Given an ncp-map $f\colon \scrA\to\scrB$
between von Neumann algebras
and a projection~$e$ from~$\scrB$
there is a least projection~$\Define{f_\diamond(e)}$ from~$\scrA$%
\index{*diamondsub@$f_\diamond$}
with~$\ceil{f(\,f_\diamond(e)^\perp\,)}\leq e^\perp$,
namely~$f_\diamond(e) =\ceil{\,ef(\,\cdot\,)e\,}$
(being the carrier 
	of the ncp-map $ef(\,\cdot\,)e$ from~\sref{carrier});
giving a map $\Define{f_\diamond}\colon \Proj(\scrB)\to\Proj(\scrA)$.
\begin{point}{30}{Proof}%
Since by definition $\ceil{\,ef(\,\cdot\,)e\,}$
is the greatest projection~$s$ of~$\scrA$
with $ef(s^\perp)e=0$ (see~\sref{carrier});
and~$ef(s^\perp )e=0$ iff~$\ceil{f(s^\perp)}
\leq\ceil{e(\,\cdot\,)e}^\perp\equiv
e^\perp$;
the projection
$\ceil{\,ef(\,\cdot\,)e\,}$
satisfies the description of~$f_\diamond(e)$.\qed
\end{point}
\end{point}
\begin{point}{40}[diamond-suprema]{Exercise}%
Let~$f\colon \scrA\to\scrB$ be an ncp-map between von Neumann algebras.
\begin{enumerate}
\item
Show that $f^\diamond(s)\leq t^\perp$
iff $f_\diamond(t)\leq s^\perp$,
for all~$s\in\Proj(\scrA)$ and~$t\in\Proj(\scrB)$.
\item
Show that $f^\diamond(\,\bigcup E\,)
= \bigcup_{e\in E} f^\diamond(e)$
for every set of projections~$E$ from~$\scrA$.
\end{enumerate}%
\spacingfix%
\end{point}%
\begin{point}{50}{Exercise}%
Show that for ncp-maps $f,g\colon\scrA\to\scrB$
between von Neumann algebras $f^\diamond = g^\diamond$
iff $f_\diamond = g_\diamond$.
In that case we say that $f$ and~$g$ are \Define{equivalent}.%
	\index{equivalent ncp-maps}
    \begin{point}{60}[contraposed]%
Show that for ncp-maps $f\colon \scrA\to\scrB$
and~$g\colon \scrB\to\scrA$ we have
$f^\diamond=g_\diamond$ iff $f_\diamond = g^\diamond$
iff $\ceil{f(s)}\leq t^\perp\iff \ceil{g(t)}\leq s^\perp$
for all projections $s$ from~$\scrA$ and~$t$ from~$\scrB$.

In that case we say that~$f$ and~$g$ are \Define{contraposed}.%
	\index{contraposed}
\end{point}
\end{point}
\begin{point}{70}[equivalent-examples]{Examples}%
\begin{enumerate}
\item
Given an element~$a$ of a von Neumann algebra~$\scrA$,
the maps $a^*(\,\cdot\,)a$ and~$a(\,\cdot\,)a^*$
on~$\scrA$ are contraposed.

If~$p$ and~$q$ are projections of~$\scrA$
with $a^*pa\leq q$
(as in~\sref{ad-ncp}),
then the maps
$a^*(\,\cdot\,)a \colon p\scrA p\to q\scrA q$
and~$a(\,\cdot\,)a^*\colon q\scrA q \to p\scrA p$
are contraposed.

In particular,
the standard corner $\pi_s\colon \scrA\to s \scrA s$
and the standard filter $c_s\colon s\scrA s\to \scrA$
for a projection~$s$ from~$\scrA$
are contraposed.
\item
An ncp-isomorphism $f\colon \scrA\to\scrB$
between von Neumann algebras
is contraposed to its inverse~$f^{-1}\colon \scrB\to\scrA$.
\item
There may be many maps equivalent to a given ncp-map $f\colon \scrA\to\scrB$
between von Neumann algebras:
show that~$(zf)^\diamond = f^\diamond$
for every positive central element~$z$ of~$\scrB$ with~$\ceil{z}=1$.
\end{enumerate}
\spacingfix%
\end{point}%
\begin{point}{80}[diamond-composition]{Exercise}%
Let $\xymatrix{
	\scrA\ar[r]|-{f}&
	\scrB\ar[r]|-{g}&
\scrC}$
be ncp-maps between von Neumann algebras~$\scrA$,
$\scrB$ and~$\scrC$.
\begin{enumerate}
\item
Show that $(g\circ f)^\diamond = g^\diamond\circ f^\diamond$
(using~\sref{ncp-ceil}),
and $(g\circ f)_\diamond = f_\diamond\circ g_\diamond$.

\item
Assuming that $f$ is equivalent 
to an ncp-map $f'\colon \scrA\to\scrB$
and~$g$ is equivalent to
an ncp-map~$g'\colon \scrB\to\scrC$,
show that~$g\circ f$ is equivalent to~$g'\circ f'$.
\item
Assuming that $f$ is contraposed to
an ncp-map $f'\colon \scrB\to\scrA$
and~$g$ is contraposed to
an ncp-map $g'\colon \scrC\to\scrB$,
show that~$g\circ f$ is contraposed to~$f'\circ g'$.
\end{enumerate}
\spacingfix
\end{point}%
\begin{point}{90}[diamond-sum]{Proposition}%
Given ncp-maps~$f,g\colon \scrA\to\scrB$
between von Neumann algebras
\begin{equation*}
(f+g)^\diamond(s) \,=\, f^\diamond(s)\, \cup\, g^\diamond(s)
\qquad
\text{and}
\qquad (f+g)_\diamond(t)\, =\, f_\diamond(t) \,\cup\, g_\diamond(t)
\end{equation*}
for all~$s\in \Proj(\scrA)$ and~$t\in \Proj(\scrB)$.
\begin{point}{100}[diamond-sum-proof]{Proof}%
Note that $(f+g)^\diamond(s)
=  \ceil{(f+g)(s)}
= \ceil{f(s)+g(s)}
= \ceil{f(s)}\cup \ceil{g(s)}
= f^\diamond(s) \cup g^\diamond(s)$
by~\sref{ceil-basic}.
Since~$(f+g)_\diamond(t)\leq s^\perp$
iff~$f^\diamond(s)\cup g^\diamond(s)\equiv (f+g)^\diamond(s)\leq t^\perp$
iff both $f^\diamond(s)\leq t^\perp$ and~$g^\diamond(s)\leq t^\perp$
iff both $f_\diamond(t)\leq s^\perp$ and~$g_\diamond(t)\leq s^\perp$
iff~$f_\diamond(t)\cup g_\diamond(t)\leq s^\perp$,
we see that $(f+g)_\diamond(t)=f_\diamond(t)\cup g_\diamond(t)$.\qed
\end{point}
\end{point}
\begin{point}{110}[carrier-f-dagger-f]{Lemma}%
Given contraposed
maps~$f\colon \scrA\to\scrB$
and~$g\colon \scrB\to\scrA$ between von Neumann algebras,
we have $\ceil{f}=\ceil{gf}$.
\begin{point}{120}{Proof}%
$\ceil{gf}=(gf)_\diamond(1)
= f_\diamond(g_\diamond(1))
= g^\diamond(\ceil{g})
= g^\diamond(1)=f_\diamond(1)=\ceil{f}$.\qed
\end{point}
\end{point}
\end{parsec}
\subsection{Rigidity}
\begin{parsec}{1020}%
\begin{point}{10}%
We now turn to  a remarkable property shared
by pure and nmiu-maps.
\end{point}
\begin{point}{20}[rigid]{Definition}%
We say that an ncp-map $f\colon \scrA\to\scrB$
between von Neumann algebras is \Define{rigid}%
\index{rigid ncp-map}
when the only ncp-map $g\colon \scrA\to\scrB$
with $g(1)=f(1)$ and $\ceil{f(p)}=\ceil{g(p)}$ for all
projections~$p$ from~$\scrA$ is~$f$ itself.
\end{point}
\begin{point}{30}[rigid-ncp-extreme]{Proposition}%
A rigid map $f\colon \scrA\to\scrB$
between von Neumann algebras
is extreme among the ncp-maps $g\colon \scrA\to\scrB$
with $g(1)=f(1)$.
\begin{point}{40}{Proof}%
Given $f\equiv \lambda g_1 + \lambda^\perp g_2$
where~$\lambda\in(0,1)$ 
and~$g_1,g_2\colon \scrA\to\scrB$
are ncp-maps with $g_i(1)=f(1)$,
we must show that~$f=g_1=g_2$.
Note that for every projection~$s$
of~$\scrA$
we have~$f^\diamond(s) = (\lambda g_1+\lambda^\perp g_2)^\diamond(s)
= g_1^\diamond(s)\cup g_2^\diamond(s)$
by~\sref{diamond-sum} and~\sref{equivalent-examples};
and in particular~$g_1^\diamond(s)\leq f^\diamond(s)$.
Then for $h:=\lambda  g_1 + \lambda^\perp f$
we have $h(1)=f(1)$
and~$h^\diamond(s) = g_1^\diamond(s)\cup f^\diamond(s)
= f^\diamond(s)$,
so that~$\lambda g_1 + \lambda^\perp f \equiv 
h=f = \lambda g_1 +\lambda^\perp g_2$ by rigidity of~$f$;
and thus~$f=g_2$.
Similarly, $f=g_1$.\qed%
\end{point}
\end{point}
\begin{point}{50}[nmiu-rigid]{Proposition}%
A nmiu-map $\varrho\colon \scrA\to\scrB$
between von Neumann algebras is rigid.
\begin{point}{60}{Proof}%
Let~$g\colon \scrA\to\scrB$
be an ncpu-map
with~$\ceil{\varrho(p)}=\ceil{g(p)}$
for every projection~$p$ of~$\scrA$.
To show that~$\varrho$ is rigid,
we must show that~$g=\varrho$,
and for this, it suffices to prove that $g(p)=\varrho(p)$
for every projection~$p$ of~$\scrA$.
To this end, we'll show that~$g$ is multiplicative,
because then~$g$ maps projections to projections,
so that $g(p)=\ceil{g(p)}=\ceil{\varrho(p)}=\varrho(p)$
for every projection~$p$ of~$\scrA$.
We'll show that $g$ is multiplicative
using~\sref{gardner}
by proving that
$\ceil{g(p)}\ceil{g(q)}=0$
for projections $p$ and~$q$ of~$\scrA$
with~$pq=0$.
Indeed,
$\ceil{g(p)}\ceil{g(q)}=\ceil{\varrho(p)}\ceil{\varrho(q)}
=\varrho(p)\varrho(q)=\varrho(pq)=\varrho(0)=0$.\qed
\end{point}
\end{point}
\begin{point}{70}[canonical-quotient-rigid]{Lemma}%
Given an element~$b$ of a von Neumann algebra~$\scrA$
the ncp-map $a\mapsto b^* a b,\ \ceilr{b}\!\scrA\!\ceilr{b}\to\scrA$
is rigid.
\begin{point}{80}{Proof}%
Let~$g\colon \ceilr{b}\!\scrA\!\ceilr{b}\to\scrA$
be an ncp-map with~$g(1)=b^*b$
and $\ceil{b^*pb}=\ceil{g(p)}$ for every projection~$p$ 
of $\ceilr{b}\!\scrA\!\ceilr{b}$.
To prove that~$c:=b^*(\,\cdot\,) b
\colon \ceilr{b}\!\scrA\!\ceilr{b}\to \scrA$ is rigid,
we must show that~$g=c$.
Since~$c$ is a filter
(by~\sref{canonical-filter})
and~$g(1)=b^*b$
there is a unique ncp-map~$h\colon \ceilr{b}\!\scrA\!\ceilr{b}
\to\ceilr{b}\!\scrA\!\ceilr{b}$
with~$g=c\circ h$.
Our task then is to show that~$h=\id$,
and for this it suffices to show that,
for all~$a\in\ceilr{b}\!\scrA\!\ceilr{b}$,
\begin{equation}
\label{filter-rigid-1}
e_n\, h(\,e_n\, a\, e_n\,)\, e_n
\ = \ e_n \, a\,  e_n
\end{equation}
for some sequence of projections $e_1,e_2,\dotsc$
of~$\ceilr{b}\!\scrA\!\ceilr{b}$
that converges ultrastrongly to~$\ceilr{b}$,
because by~\sref{mult-jus-cont} the left-hand side of the equation above 
converges ultrastrongly to~$g(a)$,
while the right-hand side converges ultrastrongly to~$a$.
We'll take $e_N := \sum_{n=1}^N \ceill{t_n}$,
where~$t_1,t_2,\dotsc$
is an approximate pseudoinverse for~$b$,
because $\ceilr{b} = \sum_n\ceill{t_n}$.

Since the identity on~$e_n \scrA e_n$ is rigid
by~\sref{nmiu-rigid},
it suffices (for~\eqref{filter-rigid-1})
to show that 
$e_n h(e_n) e_n = e_n$
and 
$\ceil{e_nh(p)e_n} = p$
for every projection $p$ from $e_n\scrA e_n$.
Writing~$s_N:=\sum_{n=1}^N t_n$,
we have $bs_n = e_n$,
and so
$
\ceil{e_n h(p) e_n}
=
\ceil{s_n^* b^* h(p) b s_n}
=
\ceil{s_n^*g(p) s_n}
=
\ceil{s_n^* \ceil{g(p)} s_n}
=
\ceil{s_n^* \ceil{b^* p b} s_n}
=
\ceil{s_n^* b^* p b s_n}
=
\ceil{e_n p e_n }
$
for every 
projection~$p$ from~$\ceilr{b}\!\scrA\!\ceilr{b}$.
In particular, $\ceil{e_n h(p)e_n} = p$
when~$p$ is from~$e_n \scrA e_n$;
and we see $\ceil{e_n h(e_n^\perp) e_n}=\ceil{e_n e_n^\perp e_n}=0$
when we take~$p=e_n^\perp$,
so that~$e_n h(e_n^\perp) e_n =0$,
which yields $e_nh(e_n)e_n = e_n$.\qed
\end{point}
\end{point}
\begin{point}{90}[pure-is-rigid]{Theorem}%
\index{pure map!is rigid}
Every pure map between von Neumann algebras is rigid.
\begin{point}{100}{Proof}%
Let~$f\colon \scrA\to\scrB$ be a pure map between von Neumann algebras,
and let~$g\colon \scrA\to\scrB$ be an ncp-map
with $f(1)=g(1)$
and $f^\diamond = g^\diamond$.
To show that~$f$
is rigid,
we must prove that~$f=g$.
We know by~\sref{square-f}
that $f$ can be written as $f\equiv c_{f(1)} \circ [f]\circ \pi_{\ceil{f}}$,
and that~$c_{f(1)}$ is rigid,
by~\sref{canonical-quotient-rigid},
which we'll use shortly.
To this end,
note that since~$f^\diamond = g^\diamond$,
we have $f_\diamond = g_\diamond$,
and so $\ceil{f}=f_\diamond(1)=g_\diamond(1)=\ceil{g}$.
As~$\pi_{\ceil{f}}$ is a corner of~$\ceil{f}=\ceil{g}$,
there is a unique ncp-map $h\colon \ceil{f}\!\scrA\!\ceil{f}\to\scrB$
with $h\circ \pi_{\ceil{f}} =g$. 
Since then
$h^\diamond \circ \pi_{\ceil{f}}^\diamond
= g^\diamond = f^\diamond 
= c_{f(1)}^\diamond
\circ [f]^\diamond \circ \pi_{\ceil{f}}^\diamond$,
 and $\pi_{\ceil{f}}^\diamond$ is clearly surjective,
we get~$h^\diamond = c_{f(1)}^\diamond\circ [f]^\diamond$,
and thus  $(h\circ [f]^{-1})^\diamond = c_{f(1)}^\diamond$,
using here that~$[f]$ is invertible,
because~$f$ is pure.
Now,
using that~$c_{f(1)}$
is rigid,
and $h([f]^{-1}(1))=h(1)=h(\pi_{\ceil{f}}(1))=g(1)=f(1)=c_{f(1)}(1)$,
we get~$h\circ [f]^{-1}=c_{f(1)}$,
which yields
$g=h\circ \pi_{\ceil{f}} 
=h\circ [f]^{-1}\circ [f]\circ \pi_{\ceil{f}}
= c_{f(1)} \circ [f] \circ \pi_{\ceil{f}} = f$,
and thus~$f$ is rigid.\qed
\end{point}
\end{point}
\end{parsec}
\subsection{\texorpdfstring{$\diamond$-P}{Diamond p}ositivity}
\begin{parsec}{1030}%
\begin{point}{10}{Definition}%
We'll call an ncp-map $f\colon \scrA\to\scrA$
between von Neumann algebras
\begin{enumerate}
\item
\Define{$\diamond$-self-adjoint}%
\index{$\diamond$-self-adjoint}
if~$f$ is pure and contraposed to itself
($f^\diamond = f_\diamond$), and
\item
\Define{$\diamond$-positive}%
\index{$\diamond$-positive}
if~$f\equiv gg$
for some $\diamond$-self-adjoint map
$g\colon \scrA\to\scrA$.
\end{enumerate}
We added ``$\diamond$-'' to ``positive''
not only to distinguish it from 
the pre-existing notion of positivity for maps between $C^*$-algebras,
but also to contrast it with
	the notion of ``$\dagger$-positivity''
	that appears in the following thesis (see~\sref{dagger-effectus}).
\end{point}
\begin{point}{20}[purely-positive-examples]{Examples}%
Let~$\scrA$ be a von Neumann algebra.
\begin{enumerate}
\item
Given~$a\in\Real{\scrA}$ the map~$a(\,\cdot\,)a\colon \scrA\to\scrA$ is 
$\diamond$-self-adjoint.
\item
Given~$a\in\scrA_+$
the map $a(\,\cdot\,)a\colon \scrA\to\scrA$
is $\diamond$-positive.
\end{enumerate}
\spacingfix
\end{point}%
\begin{point}{30}[purely-positive-basic]{Exercise}%
Let~$f\colon \scrA\to\scrA$
be an ncp-map,
where~$\scrA$ is a von Neumann algebra.
\begin{enumerate}
\item
Show that~$\ceil{f}=\ceil{f(1)}$ when~$f$
is  $\diamond$-self-adjoint.
\item
Assuming~$f$ is $\diamond$-self-adjoint,
show that~$ff$ is $\diamond$-self-adjoint,
and show that~$\ceil{ff}=\ceil{f}$ (cf.~\sref{carrier-f-dagger-f}).
\item
Show that~$f$ is $\diamond$-self-adjoint
when~$f$ is $\diamond$-positive.
\end{enumerate}
\spacingfix%
\end{point}
\end{parsec}
\begin{parsec}{1040}%
\begin{point}{10}%
We now turn 
	to the question
	roughly speaking 
to what extent a filter~$c$ is determined by
its action~$c^\diamond\colon e\mapsto \ceil{c(e)}$ on projections;
we will see (essentially in~\sref{positive-quotients-centrally-similar})
that two filters $c_1$ and~$c_2$
are equivalent, $c_1^\diamond = c_2^\diamond$,
if and only if~$c_1(1)$ and~$c_2(1)$
are equal up to some central elements,
that is, \emph{centrally similar}.
\end{point}
\begin{point}{20}{Definition}%
We say that positive elements $p$ and~$q$ of a von Neumann algebra~$\scrA$
are \Define{centrally similar}
if~$cp=dq$ for some positive central elements~$c$ and~$d$ of~$\scrA$
with~$\ceil{p}\leq \ceil{c}$
and~$\ceil{q}\leq \ceil{d}$.
\end{point}
\begin{point}{30}[centrally-similar-basic]{Exercise}%
Let~$p$ and~$q$ be positive elements
of a von Neumann algebra~$\scrA$.
\begin{enumerate}
\item
Show that when~$p$ and~$q$ are centrally similar,
every element~$a$ of~$\scrA$ that commutes
with~$p$ commutes with~$q$ too;
and in particular, $pq=qp$.
\item
Show that when~$p$ and~$q$ are centrally similar,
$\ceil{p}=\ceil{q}$.
\item[2a.]
Assuming that $p\leq Bq$ for some $B\in [0,\infty)$,
show that~$p$ and~$q$ are centrally similar iff $p/q$ is central.

Show that~$p$ is centrally similar to~$1$ iff~$p$ is central.

Show that~$p$ is centrally similar to~$p^2$ iff~$p$ is central.
\item
Show that when~$p$ and~$q$ commute,
and both $\frac{p\wedge q}{p}$ 
and~$\frac{p\wedge q}{q}$
are central,
$p$ and~$q$ are centrally similar.
\item
Show that when~$p$ and~$q$ are pseudoinvertible,
then:
$p$ and~$q$ are centrally similar iff
$pq^{\sim 1}$ is central
iff $qp^{\sim 1}$ is central
iff both $(p\wedge q)p^{\sim 1}$
and~$(p\wedge q)q^{\sim 1}$ are central.
\item
Assuming that $p$ and~$q$ commute
and $e_1 \leq e_2 \leq \dotsb$
are projections commuting with~$p$ and~$q$
with~$\bigcup_n e_n=\ceil{p}$
such that the~$e_np$ and~$e_nq$
are pseudoinvertible,
and centrally similar,
show that $p$ and~$q$ are centrally similar
on the 
grounds that both  $\frac{p\wedge q}{p}$
and~$\frac{p\wedge q}{q}$ are central.

(Hint:
$\smash{e_n \frac{p\wedge q}{p} = \frac{(e_np)\wedge(e_nq)}{e_np}}$
are central,
and
converge ultraweakly to $\frac{p\wedge q}{p}$.)
\end{enumerate}
\spacingfix%
\end{point}%
\begin{point}{40}[centrally-similar-fundamental]{Lemma}%
Suppose that $\ceil{q \, \vartheta(e)\, q}\leq e$
and~$\ceil{q \,\vartheta(e^\perp)\, q} \leq e^\perp$,
where~$e$ is a projection of a von Neumann algebra~$\scrA$,
$q$ is a positive element of~$\scrA$,
and~$\vartheta\colon \scrA\to\scrA$ is an miu-map.
Then~$eq=qe$ and~$\vartheta(e)=e$.
\begin{point}{50}{Proof}%
We have $\vartheta(e)qe=\vartheta(e)q$,
because~$e\geq  \ceil{q\,\vartheta(e)\,q}
\equiv \ceill{\vartheta(e)q}$ (see~\sref{ceill-basic}).
Similarly, $\vartheta(e^\perp)qe^\perp = \vartheta(e^\perp)q$,
because $e^\perp \geq \ceil{q\,\vartheta(e^\perp)\,q}
\equiv \ceill{\vartheta(e^\perp)q}$,
and so~$\vartheta(e^\perp)qe=0$,
which implies $\vartheta(e)qe=qe$.
Thus~$qe=\vartheta(e)qe=\vartheta(e)q$,
and so $q^2e=q\vartheta(e)q$ is self-adjoint,
which gives us that $q^2e=(q^2e)^*=eq^2$.
Since~$q^2$ commutes with~$e$,
$q=\smash{\sqrt{q^2}}$ commutes
with~$e$ too (see~\sref{sqrt}).
Finally, $\vartheta(e)q=qe=eq$
and~$\ceil{q}=1$
imply that~$\vartheta(e)=e$ by~\sref{mult-cancellation}.\qed
\end{point}
\end{point}
\begin{point}{60}[centrally-similar-corollary]{Corollary}%
A positive element~$q$
of a von Neumann algebra~$\scrA$
with~$\ceil{q}=1$
is central provided
there is an miu-map~$\vartheta\colon \scrA\to\scrA$
with $\ceil{q\,\vartheta(e)\,q}\leq e$
for every projection~$e$ from~$\scrA$;
and in that case~$\vartheta=\id$.
\end{point}
\begin{point}{70}[positive-quotients-centrally-similar]{Proposition}%
Positive elements~$p$ and~$q$
of a von Neumann algebra~$\scrA$
with~$\ceil{p}=\ceil{q}=1$
are centrally similar 
when there is an miu-isomorphism
$\vartheta\colon \scrA\to\scrA$
with~$\ceil{pep}=\ceil{q\,\vartheta(e)\,q}$
for all projections~$e$ of~$\scrA$;
and in that case  $\vartheta=\id$.
\begin{point}{80}{Proof}%
Let~$e$ be a projection from~$\scrA$ with~$ep=pe$.
Since~$1=\ceil{p}=\ceil{\smash{p^2}}$
we have $e=\ceil{e\ceil{\smash{p^2}}e}
=\ceil{e\smash{p^2} e}=\ceil{pep}=\ceil{q\,\vartheta(e)\,q}$.
Since~$e^\perp$ commutes with~$p$ too,
we get~$e^\perp = \ceil{\smash{q\,\vartheta(e^\perp)\,q}}$
by the same token;
and thus~$eq=qe$ and~$\vartheta(e)=e$ 
by~\sref{centrally-similar-fundamental}.
Since~$p$ is the norm limit
of linear combinations of such projections~$e$,
we get $pq=qp$ and~$\vartheta(p)=p$.

Since~$p$ and~$q$ commute,
we can find a sequence
of projections~$e_1\leq e_2 \leq \dotsb$
that commute with~$p$ and~$q$
with~$\bigcup_n e_n =\ceil{p}$
and such that $pe_n$ and~$qe_n$
are pseudoinvertible --- one may,
for example,
take $e_N:=\sum_{n=1}^N \ceil{t_n}$
where~$t_1,t_2,\dotsc$
is an approximate pseudoinverse
of~$p\wedge q$ (see~\sref{approximate-pseudoinverse}).
Note that to prove that~$p$ and~$q$ are centrally similar,
it suffices to show that $pe_n$ and~$qe_n$ are centrally similar,
by~\sref{centrally-similar-basic}.
Further, to prove that~$\vartheta(a)=a$
for some~$a\in\scrA$,
it suffices to show that~$\vartheta( e_n a e_n  ) = e_n a e_n$,
because $e_n a e_n$ converges ultrastrongly to~$a$
	by~\sref{mult-jus-cont}.
Note that~$\vartheta(e_n)=e_n$,
because~$e_np=pe_n$,
and so~$\vartheta$ maps~$e_n\scrA e_n$ into~$e_n \scrA e_n$.

Thus, by considering~$e_n \scrA e_n$ 
instead of~$\scrA$,
and the restriction of~$\vartheta$ to~$e_n\scrA e_n$
instead of~$\vartheta$,
and~$pe_n$ and~$qe_n$
instead of~$p$ and~$q$,
we reduce the problem to the case that~$p$ and~$q$ are invertible;
and so we may assume without loss of generality that~$p$ and~$q$
are invertible to start with.
Given a projection~$e$ from~$\scrA$
we have $\ceil{p^{-1}q \,\vartheta(e)\, q p^{-1}}
= \ceil{p^{-1}\ceil{q\,\vartheta(e)\,q}p^{-1}}
= \ceil{p^{-1}\ceil{pep}p^{-1}}=e$;
so 
by~\sref{centrally-similar-corollary},
we get that
$\vartheta=\id$
and
$p^{-1}q$ is central;
and so
$p$ and~$q$ are centrally similar (by~\sref{centrally-similar-basic}).
\qed
\end{point}
\end{point}
\begin{point}{90}[faithful-positive-map-uniqueness]{Proposition}%
A faithful $\diamond$-positive map $f\colon \scrA\to\scrA$
on a von Neumann algebra~$\scrA$
is of the form~$f=\sqrt{p}(\,\cdot\,)\sqrt{p}$
where $p:=f(1)$.
\begin{point}{100}{Proof}%
Note that~$f$,
being faithful and pure,
is a filter
(by~\sref{special-pure-maps}),
and thus of the form $f\equiv \sqrt{p}\,\vartheta(\,\cdot\,)\,\sqrt{p}$
for some isomorphism~$\vartheta\colon \scrA\to\scrA$.
Our task then is to show that~$\vartheta=\id$,
and for this
it suffices, by~\sref{positive-quotients-centrally-similar},
to find some positive~$q$ in~$\scrA$ with~$\ceil{q}=1$
and~$f^\diamond(e)\equiv\ceil{\sqrt{p}\,\vartheta(e)\,\sqrt{p}}
= \ceil{qeq}$ for all projections~$e$ in~$\scrA$.

Since~$f$ is $\diamond$-positive,
we have~$f\equiv \xi \xi$ for some $\diamond$-self-adjoint
map~$\xi\colon \scrA\to\scrA$.
Since~$1=\ceil{f}=f_\diamond(1)=
\xi_\diamond(\xi_\diamond(1))\leq \xi_\diamond(1)=\ceil{\xi}$
we have~$\ceil{\xi}=1$,
and so, $\xi$, being pure and faithful,
is a filter (by~\sref{special-pure-maps}).
Furthermore,
as~$\tilde\xi:=\sqrt{\xi(1)}(\,\cdot\,)\sqrt{\xi(1)}\colon \scrA\to\scrA$
is a filter of~$\xi(1)$ too,
there is an isomorphism~$\alpha\colon \scrA\to\scrA$
with $\xi=\tilde\xi\alpha$.
Now, $ {\tilde\xi}^\diamond \alpha^\diamond
={\xi}^\diamond={\xi}_\diamond
=\alpha_\diamond\tilde \xi_\diamond
= (\alpha^\diamond)^{-1}{\tilde\xi}^\diamond$
implies~${\tilde\xi}^\diamond = \alpha^\diamond 
{\tilde \xi}^\diamond \alpha^\diamond$,
and 
$f^\diamond= (\xi\xi)^\diamond
= {\tilde \xi}^\diamond\alpha^\diamond{\tilde \xi}^\diamond\alpha^\diamond
={\tilde \xi}^\diamond{\tilde \xi}^\diamond=(\tilde \xi\tilde \xi)^\diamond$.
In other words,
$\ceil{\sqrt{p}\,\vartheta(e)\,\sqrt{p}}
=f^\diamond(e)=(\tilde\xi\tilde\xi)^\diamond(e)
= \ceil{\xi(1)\,e\,\xi(1)}$
for all projections~$e$ of~$\scrA$,
which implies that~$\vartheta=\id$
by~\sref{positive-quotients-centrally-similar},
and hence that~$f=\sqrt{p}\,(\,\cdot\,)\,\sqrt{p}$.\qed
\end{point}
\end{point}
\end{parsec}
\begin{parsec}{1050}%
\begin{point}{10}%
To strip 
from~\sref{faithful-positive-map-uniqueness}
the assumption 
that~$f$
be faithful 
we employ this device:
\end{point}
\begin{point}{20}[chevron-f]{Definition}%
Given an ncp-map $f\colon \scrA\to\scrB$
between von Neumann algebras
we denote by
$\Define{\left<f\right>}\colon \ceil{f}\!\scrA\!\ceil{f}
\to \ceil{f(1)}\!\scrB\!\ceil{f(1)}$%
\index{*diamondf@$\left<f\right>$}
the unique ncp-map
such that 
\begin{equation*}
\xymatrix@C=6em{
\scrA
\ar[r]^f
\ar[d]_{\pi_{\ceil{f}}}
&
\scrB
\\
\ceil{f}\!\scrA\!\ceil{f}
\ar[r]^{\left<f\right>}
&
\ceil{f(1)}\!\scrB\!\ceil{f(1)}
\ar[u]_{c_{\ceil{f(1)}}}
}
\end{equation*}
commutes.
(Compare this with the definition of~$[f]$ in~\sref{square-f}.)
\end{point}
\begin{point}{30}[chevron-f-basic]{Exercise}%
Let~$f\colon \scrA\to\scrB$ be an ncp-map.
\begin{enumerate}
\item
Show that~$\left<f\right>
= \pi_{\ceil{f(1)}}\circ f \circ c_{\ceil{f}}$
(using, perhaps, that $\pi_{\ceil{f}}\circ c_{\ceil{f}}=\id$).
\item
Show that
$\left<f\right> = \pi_{\ceil{f(1)}} \circ c_{f(1)}\circ [f]$.

(Thus $\left<f\right>\!(a) = 
\sqrt{f(1)}\ [f]\!(a)\ \sqrt{f(1)}$
for all $a$ from~$\ceil{f}\!\scrA\!\ceil{f}$.)
\item
Show that~$\left<f\right>$
is faithful,
and~$\left<f\right>\!(1)=f(1)$.
\item
Assuming that~$f$ is pure,
show that~$\left<f\right>$ is pure,
and hence a filter (by~\sref{special-pure-maps}).
\end{enumerate}
\spacingfix%
\end{point}%
\begin{point}{40}[chevron-f-purely-positive]{Exercise}%
Let~$f\colon \scrA\to\scrA$
be an ncp-map, where~$\scrA$ is a von Neumann algebra.
\begin{enumerate}
\item
Suppose that~$f$ is $\diamond$-self-adjoint.

Recall that~$\ceil{f}=\ceil{f(1)}$,
and so $\left<f\right>\colon \ceil{f}\!\scrA\!\ceil{f}
\to \ceil{f}\!\scrA\!\ceil{f}$.

Prove that~$\left<f\right>$
is $\diamond$-self-adjoint.
\item
Suppose again that~$f$ is $\diamond$-self-adjoint,
and recall from~\sref{purely-positive-basic} that $f^2$
is $\diamond$-self-adjoint, and~$\ceil{f^2} = \ceil{f}$.
Show that $\left<f^2\right> = \left<f\right>^2$.
\item
Assuming that~$f$ is $\diamond$-positive,
show that~$\left<f\right>$ is $\diamond$-positive.
\end{enumerate}
\spacingfix
\end{point}%
\begin{point}{50}[positive-map-uniqueness]{Theorem}%
Given a positive element~$p$ of a von Neumann algebra~$\scrA$
there is a unique $\diamond$-positive map $f\colon \scrA\to\scrA$
with~$f(1)=p$,
namely~$f=\sqrt{p}(\,\cdot\,)\sqrt{p}$.
\begin{point}{60}{Proof}%
We've already seen in~\sref{purely-positive-examples}
that $f=\sqrt{p}(\,\cdot\,)\sqrt{p}\colon \scrA\to\scrA$
is a $\diamond$-positive map with~$f(1)=p$.
Concerning uniqueness,
(given arbitrary~$f$)
the map~$\left<f\right>\colon \ceil{p}\!\scrA\!\ceil{p}
\to \ceil{p}\!\scrA\!\ceil{p}$
from~\sref{chevron-f}
is $\diamond$-positive by~\sref{chevron-f-purely-positive},
and faithful by~\sref{chevron-f-basic},
and so of the form
$\left<f\right>=\sqrt{p}(\,\cdot\,)\sqrt{p}$
by~\sref{faithful-positive-map-uniqueness}
(since~$\left<f\right>\!(1)=f(1)=p$);
implying that
$f= c_{\ceil{p}}\circ\left<f\right>\circ \pi_{\ceil{p}}
= \sqrt{p}\ceil{p}(\,\cdot\,)\ceil{p}\sqrt{p}
= \sqrt{p}(\,\cdot\,)\sqrt{p}$.\qed
\end{point}
\end{point}
\begin{point}{70}[sqrt-axiom]{Corollary (``Square Root Axiom'')}%
\index{square root axiom}
Given a positive element~$p$ of a von Neumann algebra~$\scrA$
there is a unique $\diamond$-positive map~$g\colon \scrA\to\scrA$
with~$g(g(1))=p$, namely
$g=\sqrt[4]{p}\,(\,\cdot\,)\,\sqrt[4]{p}$.
\begin{point}{80}{Proof}%
Any $\diamond$-positive map~$g\colon \scrA\to\scrA$ with~$g(g(1))=p$
will be of the form
$g=\smash{\sqrt{g(1)}\,(\,\cdot\,)\,\sqrt{g(1)}}$
by~\sref{positive-map-uniqueness};
so that~$p=g(g(1))=g(1)^2$
implies that~$g(1)=\sqrt{p}$
by~\sref{sqrt},
thereby giving~$g=\sqrt[4]{p}\,(\,\cdot\,)\,\sqrt[4]{p}$.\qed
\end{point}
\end{point}
\end{parsec}
\begin{parsec}{1060}%
\begin{point}{10}[uniqueness-sequential-product]{Theorem}%
\index{sequential product}
On the effects of every von Neumann algebra~$\scrA$
there is a unique binary operation~$\ast$
such that for all~$p$ from~$[0,1]_\scrA$,
\begin{enumerate}[A.]
\item \label{ax1}
$p\ast 1 = p$,
\item
	\label{ax2}
$p\ast q = f(q)$
for all~$q$ from~$[0,1]_\scrA$
for some pure map~$f\colon \scrA\to\scrA$,
\item\label{ax3}
$p\ast (p\ast q)=(p\ast p)\ast q$
for all~$q$ from $[0,1]_\scrA$,
\item\label{ax4}
$p=q\ast q$ for some~$q$ from~$[0,1]_\scrA$,
\item\label{ax5}
$p \ast e_1 \leq e_2^\perp
\iff p\ast e_2 \leq e_1^\perp$
for all projections~$e_1,e_2$ from~$\scrA$;
\end{enumerate}
namely, the sequential product,
given by
$p\ast q = \sqrt{p}q\sqrt{p}$
for all~$p, q$ from~$[0,1]_\scrA$.
\begin{point}{20}{Proof}%
Let~$p$  from $[0,1]_\scrA$ be given.
Pick~$p'$ from~$[0,1]_\scrA$
with $p = p'\ast p'$
using~\ref{ax4},
and find a pure map~$f\colon \scrA\to\scrA$
with~$f(q)=p'\ast q$ for all~$q$ from $[0,1]_\scrA$
using~\ref{ax2}.
Then~$f$ is $\diamond$-self-adjoint by~\ref{ax5},
and so~$ff$ is $\diamond$-positive.
Since~$f(f(1))=p'\ast (p'\ast 1)
= p'\ast p'=p$ by~\ref{ax1},
$ff=\sqrt{p}(\,\cdot\,)\sqrt{p}$
by~\sref{positive-map-uniqueness},
so $p\ast q
= (p'\ast p')\ast q
= p'\ast (p' \ast q) = f(f(q))=\sqrt{p}q\sqrt{p}$
for all~$q\in [0,1]_\scrA$ by~\ref{ax3}.\qed
\end{point}
\end{point}
\begin{point}{30}{Exercise}%
None of the axioms
from~\sref{uniqueness-sequential-product}
may be omitted (except perhaps~\ref{ax4},
see~\sref{fourth-axiom}):
\begin{enumerate}
\item
Show that
$p\ast q := \ceil{p}q\ceil{p}$
satisfies all axioms of~\sref{uniqueness-sequential-product}
except~\ref{ax1}.
\item
Show that  $p\ast q := \floor{p}q\floor{p}\ +\ \smash{\sqrt{p-\floor{p}}\,q\,
\sqrt{p-\floor{p}}}$
satisfies all axioms except~\ref{ax2}.
\item
Show that if for every effect~$p$ of~$\scrA$
we pick a unitary~$u_p$ from~$\ceil{p}\!\scrA\!\ceil{p}$
then~$\ast$ given by
$p\ast q= \sqrt{p}u_p^* \,q\, u_p\sqrt{p}$
satisfies~\ref{ax1} and~\ref{ax2}.

Show that this~$\ast$ obeys~\ref{ax3} when~$u_p^2=u_{p^2}$,
and~\ref{ax4} when $pu_p=u_p p$,
and~\ref{ax5} when~$u_p^*=u_p$.

Conclude that when $u_p$ is defined by $u_p:=g(p)$,
where~$g\colon [0,1]\to\{-1,1\}$
is any Borel function with $g(\nicefrac{2}{3})=1$
and~$g(\nicefrac{4}{9})=-1$
the operation~$\ast$ (defined by~$u_p$ as above) satisfies
all conditions of~\sref{uniqueness-sequential-product} except~\ref{ax3}.
\item
Show that there is a Borel
function~$g\colon[0,1]\to S^1$
with $g(\nicefrac{1}{2})\neq 1$
and~$g(\lambda^2)=g(\lambda)^2$ for all~$\lambda\in [0,1]$,
and that~$\ast$ given by~$p\ast q = \sqrt{p} g(p)^* \,q \,g(p)\sqrt{p}$
satisfies all conditions of~\sref{uniqueness-sequential-product}
except~\ref{ax5}.
\end{enumerate}
\spacingfix%
\end{point}%
\begin{point}{40}[fourth-axiom]{Problem}%
Is there a binary operation~$\ast$
on the effects~$[0,1]_\scrA$
of a von Neumann algebra~$\scrA$
that satisfies all axioms of~\sref{uniqueness-sequential-product}
except~\sref{ax4}?
\end{point}
\begin{point}{50}{Remark}%
The axioms for the sequential product
(on a single von Neumann algebra)
presented here (in~\sref{uniqueness-sequential-product})
evolved
from the following axioms for all
sequential products
on von Neumann algebras
$(\ast_\scrA)_\scrA$
we previously published in~\cite{westerbaan2016universal}.
\begin{enumerate}
\item[Ax.1]%
For every effect~$p$ of a von Neumann algebra~$\scrA$
there is a filter $c\colon \scrC\to\scrA$ of~$p$
and a corner~$\pi\colon \scrA\to\scrC$
		of~$\floor{p}$ with~$p\ast_\scrA q = c(\pi(q))$
		for all $q\in[0,1]_\scrA$.
\item[Ax.2]%
$p\ast_\scrA(p\ast_\scrA q)= (p \ast_\scrA p )\ast_\scrA q$
for all effects~$p$ and~$q$ from a von Neumann algebra~$\scrA$.
\item[Ax.3]%
$f(p\ast_\scrA q)= f(p)\ast_\scrB f(q)$
for every nmisu-map $f\colon \scrA\to\scrB$
between von Neumann algebras
and all effects~$p$ and~$q$ from~$\scrA$.
\item[Ax.4]%
$p\ast_\scrA e_1 \leq e_2^\perp
\iff p\ast_\scrA e_2 \leq e_1^\perp$
for every effect~$p$ from a von Neumann algebra~$\scrA$
and projections $e_1$ and~$e_2$ from~$\scrA$.
\end{enumerate}
Note that~Ax.2 and Ax.4 are mutatis mutandis 
the same as axioms~C and~E, respectively,
and Ax.1 is essentially the same as the combination 
of axioms~A and~B.
In other words,
we managed to get rid of~Ax.3---and with it
the need to axiomatise all sequential products simultaneously---at
the slight cost of adding axiom~D,
though that one might be superfluous as well
	(see~\sref{fourth-axiom}).

We refer to \S{}VI of~\cite{westerbaan2016universal}
for comments on the relation
of our axioms
	with those of Gudder 
	and Lat\'emoli\`ere\cite{gudder2008characterization}
and for some more pointers to the literature.
\end{point}
\end{parsec}

\section{Tensor product}
\label{S:tensor}
\begin{parsec}{1070}%
\begin{point}{10}%
The tensor
product of von Neumann algebras~$\scrA$
and~$\scrB$
represented on Hilbert spaces~$\scrH$ and~$\scrK$,
respectively,
is usually defined as the von Neumann subalgebra
of~$\scrB(\scrH\otimes\scrK)$
generated by the operators
on~$\scrH\otimes\scrK$
of the form~$A\otimes B$ where~$A\in\scrA$
and~$B\in\scrB$.
In line with the representation-avoiding
treatment of von Neumann algebras
from the previous chapter
we'll take an entirely different
approach
by defining the tensor product of von Neumann
algebras $\scrA$ and~$\scrB$ abstractly 
as an miu-bilinear map $\otimes\colon \scrA\times\scrB\to\scrA\otimes\scrB$
whose range generates~$\scrA\otimes\scrB$
and admits sufficiently many product functionals (see~\sref{tensor});
we'll only resort to the concrete representation
of the tensor product  mentioned above to show that 
such an abstract tensor product actually exists (see~\sref{special-tensor}).

Moreover, we'll show that the tensor product has a universal
property~\sref{tensor-universal-property}
yielding bifunctors on~$\W{cpsu}$
and~$\W{miu}$ (see~\sref{tensor-functor})
turning them into a monoidal categories
(see~\sref{vn-smc}).
In the next chapter,
    we'll see that~$\op{(\W{miu})}$ is even
monoidal closed (see~\sref{tensor-closed}).
This fact is one ingredient of our model for the quantum lambda calculus
from~\cite{model}
built of von Neumann algebras,
but more of that later.
\end{point}
\end{parsec}
\subsection{Definition}
\begin{parsec}{1080}%
\begin{point}{10}[bilinear-basic]{Definition}%
A bilinear map
	$\beta\colon \scrA\times \scrB\to\scrC$
between von Neumann algebras is
\begin{enumerate}
\item
\Define{\textbf{u}nital}%
\index{unital!bilinear map}
when~$\beta(1,1)=1$,
\item
\Define{\textbf{m}ultiplicative}%
\index{multiplicative!bilinear map}
if~$\beta(ab,cd)=\beta(a,c)\beta(b,d)$
for all~$a,b\in\scrA$, $c,d\in\scrB$,
\item
\Define{\textbf{i}nvolution preserving}
\index{involution preserving!bilinear map}
if~$\beta(a,b)^*=\beta(a^*,b^*)$
for all~$a\in\scrA$, $b\in\scrB$.
\item
(This list is extended in~\sref{tensor-extra}.)
\end{enumerate}
We abbreviate these properties as in~\sref{maps},
and say, for instance, that $\beta$ is \Define{miu-bilinear}%
\index{miu-bilinear}
when it is unital, multiplicative and involution preserving.

\end{point}
\begin{point}{20}[tensor]{Definition}%
A miu-bilinear map $\gamma\colon \scrA\times\scrB\to\scrT$
between von Neumann algebras
is a \Define{tensor product}%
\index{tensor product!of von Neumann algebras}
of~$\scrA$ and~$\scrB$
when it obeys the following three conditions.
\begin{enumerate}
\item
\label{tensor-1}
The range of~$\gamma$ generates~$\scrT$
(which  means in this case that the linear span
of the range of~$\gamma$ is ultraweakly dense in~$\scrT$.)

This implies
that for all $f \in \scrA_*$
and $g\in \scrB_*$
there is at most one $h\in \scrT_*$
with, for all~$a\in\scrA$ and $b\in\scrB$,
\begin{equation*}
	h(\gamma(a,b))\ =\  f(a)\,g(b),
\end{equation*}
which we'll call the \Define{product functional}%
\index{product functional}
for~$f$ and~$g$,
and denote by~$\Define{\gamma(f,g)}$%
\index{*gammafg@$\gamma(f,g)$, product functional}
(when it exists).
\item
\label{tensor-2}
For all np-functionals
$\sigma\colon \scrA\to\C$
and $\tau\colon \scrB\to\C$
the product functional $\gamma(\sigma,\tau)\colon \scrT\to\C$
exists and is positive.
\item
\label{tensor-3}
The product functionals 
$\gamma(\sigma,\tau)$
of np-functionals~$\sigma$ and~$\tau$ form a faithful collection
of np-functionals on~$\scrT$.
\end{enumerate}
(We'll see a slightly different characterisation
of the tensor
in which not all product functionals of np-functionals
are required to exist upfront 
in~\sref{tensor-characterization}.)
\end{point}
\begin{point}{30}{Remark}%
This compact definition of the tensor product
leaves four questions unanswered:
whether such a tensor product of two von Neumann algebras
always exists,
whether it has some  universal property,
whether it is unique in some way,
and whether it coincides with the usual definition.
We'll shortly address all four questions.
\end{point}
\end{parsec}
\subsection{Existence}
\begin{parsec}{1090}%
\begin{point}{10}%
We'll start with the existence
of a tensor product of von Neumann algebras
for which we'll first need the tensor product
of Hilbert spaces.
\end{point}
\begin{point}{20}{Definition}%
We'll call a bilinear map $\gamma\colon \scrH\times\scrK\to\scrT$
between Hilbert spaces
a \Define{tensor product}%
\index{tensor product!of Hilbert spaces}
when it obeys the following two conditions.
\begin{enumerate}
\item
The linear span of the range of~$\gamma$ is dense in~$\scrT$.
\item
$\left<\gamma(x,y),\gamma(x',y')\right>
= \left<x,x'\right> \left<y,y'\right>$
for all~$x,x'\in\scrH$ and~$y,y'\in\scrK$.
\end{enumerate}
\spacingfix
\end{point}%
\begin{point}{30}{Exercise}%
\index{tensor product!of Hilbert spaces!exists}
We're going to prove that every pair of Hilbert
spaces~$\scrH$ and~$\scrK$ admits
a tensor product.
\begin{enumerate}
\item
Given sets~$X$ and~$Y$
show that  $\gamma\colon\ell^2(X)\times \ell^2(Y)\to \ell^2(X\times Y)$
given by
\begin{equation*}
	\gamma(f,g)\ = \ (\,f(x)\,g(y)\,)_{x\in X,y\in Y}
\end{equation*}
is a tensor product of~$\ell^2(X)$ and~$\ell^2(Y)$.
\item
Show that a subset~$\scrE$ of a Hilbert space~$\scrH$
is an orthonormal basis (see~\sref{orthonormal}) iff
the map $T\colon \ell^2(\scrE)\to\scrH$
given by~$T(x)=\sum_{e\in\scrE} x_e e$
is an isometric isomorphism.
\item
Show that any pair~$\scrH$ and~$\scrK$
of Hilbert spaces has a tensor product
(using the fact that every Hilbert space
has an orthonormal basis).
\end{enumerate}
\spacingfix%
\end{point}%
\begin{point}{40}[hilb-tensor-basic]{Proposition}%
Let~$\gamma\colon \scrH\times\scrK\to\scrT$
be a tensor product
of Hilbert spaces.
\begin{enumerate}
\item
\label{hilb-tensor-basic-1}
We have $\|\gamma(x,y)\|=\|x\|\|y\|$
for all~$x\in\scrH$ and~$y\in\scrK$.
\item
\label{hilb-tensor-basic-2}
Given orthonormal bases~$\scrE$ and~$\scrF$
of~$\scrH$ and~$\scrK$, respectively,
the set
\begin{equation*}
	\scrG\,:=\,\{\,\gamma(e,f)\colon \,e\in \scrE,\,f\in\scrF\,\}
\end{equation*}
is an orthonormal basis for~$\scrT$.%
\end{enumerate}
\spacingfix%
\begin{point}{50}{Proof}%
\ref{hilb-tensor-basic-1}\ 
We have $\|\gamma(x,y)\|^2
= \left<\gamma(x,y),\gamma(x,y)\right>
=\left<x,x\right>\left<y,y\right>
= \|x\|^2\|y\|^2$.

\ref{hilb-tensor-basic-2}\ 
Since $\left<\gamma(e,e'),\gamma(f,f')\right>
=\left<e,e'\right>\left<f,f'\right>$
where~$e,e'\in \scrE$
and~$f,f'\in\scrF$,
the set~$\scrG$ is clearly orthonormal.
To see that~$\scrG$ is maximal (and thus a basis)
it suffices
to show that the span of~$\scrG$ 
is dense in~$\scrT$,
and for this it suffices
to show that each~$\gamma(x,y)$
where~$x\in\scrH$ and~$y\in\scrK$
is in the closure of the span of~$\scrG$.
Now,
since~$y=\sum_{f\in\scrF}\left<f,y\right> f$,
by~\sref{orthonormal}
and~$\left<x,(\,\cdot\,)\right>$
is bounded by~\ref{hilb-tensor-basic-1}
we have~$\gamma(x,y)= \sum_{f\in\scrF}
\left<y,f\right>\gamma(x,f)$.
Since similarly $\gamma(x,f)=\sum_{e\in\scrE} \left<e,x\right>\gamma(e,f)$ 
for all~$f\in\scrF$,
we see that~$\gamma(x,y)$
is indeed in the closure of the span of~$\scrG$.\qed
\end{point}
\end{point}
\end{parsec}
\begin{parsec}{1100}%
\begin{point}{10}{Definition}%
We'll say that a bilinear map $\beta\colon\scrH\times\scrK\to\scrL$
between Hilbert spaces
is \Define{$\ell^2$-bounded}%
\index{l2-bounded@$\ell^2$-bounded bilinear map}
by~$B\in [0,\infty)$
when
\begin{equation*}
	\|\sum_{i}\beta(x_i,y_i)\|^2
\leq B^2 \sum_{i,j}\left<x_i,x_j\right>\left<y_i,y_j\right>
\end{equation*}
for all~$x_1,\dotsc,x_n\in \scrH$
and~$y_1,\dotsc,y_n\in\scrK$.
\begin{point}{20}{Remark}%
We added the prefix ``$\ell^2$-''
to clearly distinguish
it from the boundedness
of (sesquilinear) forms from~\sref{chilb-form},
which one might call ``$\ell^\infty$-boundedness.''%
\index{linfty-bounded@$\ell^\infty$-bounded bilinear map}

This distinction is needed
since for example given a Hilbert space~$\scrH$ the
bilinear map $(f,x)\mapsto f(x)\colon \scrH^*\times \scrH\to\C$
is always $\ell^\infty$-bounded in the 
sense that~$\left|f(x)\right|\leq \|f\|\|x\|$
for all~$f\in\scrH^*$ and~$x\in\scrH$,
but it is not~$\ell^2$-bounded
when~$\scrH$ is infinite dimensional
\end{point}
\end{point}
\begin{point}{30}[hilb-tensor-universal-property]{Theorem}%
\index{tensor product!of Hilbert spaces!universal property}
\index{tensor product!of Hilbert spaces!is $\ell^2$-bounded}
A tensor product $\gamma\colon \scrH\times\scrK\to\scrT$
of Hilbert spaces
is $\ell^2$-bounded,
and initial as such 
in the sense
that for any by $B\in [0,\infty)$ 
$\ell^2$-bounded bilinear map $\beta\colon \scrH\times\scrK \to\scrL$
into a Hilbert space~$\scrL$
there is a unique bounded linear map $\beta_\gamma\colon\scrT
\to\scrL$
with $\beta_\gamma(\gamma(x,y))=\beta(x,y)$
for all~$x\in\scrH$ and~$y\in\scrK$.
Moreover, $\|\beta_\gamma\|\leq B$
for such~$\beta$.
\begin{point}{40}{Proof}%
Note that $\gamma$ is $\ell^2$-bounded,
since
for all~$x_1,\dotsc,x_n\in\scrH$,
$y_1,\dotsc,y_n\in\scrK$,
we have
$\|\sum_i\gamma(x_i,y_i)\|^2
= \sum_{i,j}\left<\gamma(x_i,y_i),
\gamma(x_j,y_j)\right>
= \sum_{i,j} \left<x_i,x_j\right>\left<y_i,y_j\right>$.

Let~$\scrE$ and~$\scrF$
be orthonormal bases
for~$\scrH$ and~$\scrK$, respectively.
Then since~$\{\,\gamma(e,f)\colon\,e\in\scrE,\,f\in\scrF\,\}$
is an orthonormal basis for~$\scrT$
by~\sref{hilb-tensor-basic},
and~$\beta_\gamma$
is fixed on it by~$\beta_\gamma(\gamma(e,f))=\beta(e,f)$,
uniqueness of~$\beta_\gamma$ is clear.

Concerning existence of~$\beta_\gamma$,
note that
since $t=\sum_{e\in\scrE,f\in\scrF} \left<\gamma(e,f),t\right>\gamma(e,f)$
for all~$t\in\scrT$
by~\sref{orthonormal},
 we'd like to define~$\beta_\gamma$ by
\begin{equation}
	\label{tensor-universal-property-1}
	\beta_\gamma(t) \ =\  
	\hspace{-1em}
	\sum_{e\in\scrE,\,f\in\scrF}
	\hspace{-.5em}
	\left<\gamma(e,f),t\right>\,\gamma(e,f);
\end{equation}
but before we can do this
we must first check that this series converges.
To this end,
note that since~$\beta$ is $\ell^2$-bounded by~$B$
we have, given~$t\in \scrT$,
\begin{alignat*}{3}
	\bigl\|\hspace{-.75em}\sum_{e\in E,\,f\in F}
	\hspace{-.5em}
\left< \gamma(e,f),t\right>\beta(e,f)
\,\bigr\|^2
\ &=\  
\bigl\|\hspace{-.75em}
\sum_{e\in E,\,f\in F}
\hspace{-.5em}
\beta(e,\,\left<\gamma(e,f),t\right>f\,)
\,\bigr\|^2
\\
\ &\leq \  
B^2\hspace{-1.5em}\sum_{e'\!\!,e\in E,\ f'\!\!,f\in F}
\hspace{-1em}
\left<e',e\right>\,
\,\left<t,\gamma(e',f')\right>\,
\left<f',f\right>\,\left<\gamma(e,f),t\right>
\\
\ &=\ 
B^2\hspace{-.75em}\sum_{e\in E,\,f\in F}
\hspace{-.5em}
\left|\left<\gamma(e,f),t\right>\right|^2
\end{alignat*}
for all finite subsets $E\subseteq\scrE$
and~$F\subseteq \scrF$.
Since $\|t\|^2=\sum_{e\in \scrE,\,f\in\scrF}
\left|\left<\gamma(e,f),t\right>\right|^2$
by Parseval's identity (\sref{orthonormal}),
we see
that the series
from~\eqref{tensor-universal-property-1}
converges defining~$\beta_\gamma(t)$,
and, moreover,
that~$\|\beta_\gamma(t)\|^2\leq B^2\|t\|^2$.

The resulting map $\beta_\gamma\colon \scrT\to\scrL$
is clearly linear,
and
bounded by~$B$.
Further,
$\beta_\gamma(\gamma(e,f))=\beta(e,f)$
for all~$e\in\scrE$ and~$f\in\scrF$
implies
that $\beta_\gamma(\gamma(x,y))=\beta(x,y)$
for all~$x\in \scrH$ and~$y\in\scrK$,
and so we're done.\qed
\end{point}
\end{point}
\begin{point}{50}{Exercise}%
Show that the tensor product of Hilbert spaces~$\scrH$
and~$\scrK$ is unique in the sense
that given tensor products $\gamma\colon \scrH\times\scrK\to\scrT$
and $\gamma'\colon\scrH\times\scrK\to\scrT'$
there is a unique isometric linear isomorphism
$\varphi\colon \scrT\to\scrT'$
with $\gamma'(x,y) = \varphi(\gamma(x,y))$
for all~$x\in\scrH$
and~$y\in\scrK$.
\end{point}
\begin{point}{60}{Notation}%
Now that we've established
that the tensor product
of Hilbert spaces~$\scrH$ and~$\scrK$
exists and is unique (up to unique isomorphism)
we just pick one and denote it by $\Define{\otimes}\colon
\scrH\times\scrK\to\Define{\scrH\otimes\scrK}$.%
\index{*tensor@$\otimes$, tensor product!$\scrH\otimes \scrK$, of Hilbert spaces}
\index{*tensor@$\otimes$, tensor product!$x \otimes y$, of elements of Hilbert spaces}
\end{point}
\end{parsec}
\begin{parsec}{1110}%
\begin{point}{10}%
Essentially to turn~$\otimes$
into a functor on
the category of Hilbert spaces
in~\sref{hilb-tensor-functor},
we'll need
the following result
(known as part of \emph{Schur's product theorem}),%
\index{Schur's Product Theorem}
which will be useful several times later on.
\end{point}
\begin{point}{20}[schur]{Lemma}%
For any natural number~$N$
the entrywise product 
$(a_{nm}b_{nm})$
of positive  $N\times N$-matrices
$(a_{nm})$ and~$(b_{nm})$
over~$\C$
is positive.
\begin{point}{30}{Proof}%
Let~$z_1,\dotsc,z_N\in\C$ be given.
To show that $(a_{nm}b_{nm})$
is positive,
it suffices by~\sref{when-a-matrix-over-a-cstar-algebra-is-positive} 
to prove that
$\sum_{n,m} \overline{z}_n a_{nm} b_{nm} z_m \geq 0$
for all~$n,m$.
Since~$(a_{nm})$
is a positive element
of the $C^*$-algebra~$M_N$
it's of the form~$(a_{nm})=C^*C$
for some $N\times N$-matrix $C\equiv(c_{nm})$
over~$\C$,
so  $a_{nm} = \sum_k \overline{c}_{kn} c_{km}$
for all~$n,m$.
Similarly, there a $N\times N$-matrix
$(d_{nm})$ over~$\C$
with $b_{nm} = \sum_\ell \overline{d}_{\ell n}
d_{\ell m}$
for all~$n,m$.
Then
\begin{alignat*}{3}
	\sum_{n,m} \overline{z}_n
a_{nm} b_{nm} z_m
\ &=\  \hspace{-.5em} \sum_{n,m,k,\ell}
\overline{z}_n \,\overline{c}_{kn}
c_{km} \,\overline{d}_{\ell n} d_{\ell m} z_m\,
\\
\ &=\  \sum_{k,\ell}
\Bigl(\sum_n
\overline{z}_n\overline{c}_{kn}
\overline{d}_{\ell n} \Bigr)
\,
\Bigl(
\sum_m
z_m
c_{km} 
d_{\ell m} 
\Bigr)
\\
\ &=\  \sum_{k,\ell}\ 
\Bigl|\sum_n z_n c_{kn} d_{\ell n} \Bigr|^2 \geq 0,
\end{alignat*}
and so~$(a_{nm}b_{nm})$ is positive.
\qed
\end{point}
\end{point}
\begin{point}{40}[mult-completely-monotone]{Exercise}%
Given square matrices 
$(a_{nm})\leq (\tilde{a}_{nm})$
and $(b_{nm})\leq (\smash{\tilde{b}_{nm}})$
over~$\C$
of the same dimensions,
show that $(\,a_{nm}b_{nm}\,)\leq (\,\tilde{a}_{nm}\smash{\tilde{b}_{nm}}\,)$.
\end{point}
\begin{point}{50}[hilb-tensor-functor]{Proposition}%
Given bounded linear maps
$A\colon \scrH\to\scrH'$
and~$B\colon \scrK\to\scrK'$
between Hilbert spaces
there is a unique bounded linear
map 
\begin{equation*}
	\Define{A\otimes B}\colon \scrH\otimes\scrK\to\scrH'\otimes\scrK'
\end{equation*}%
\index{*tensor@$\otimes$, tensor product!$A\otimes B$, 
of operators between Hilbert spaces}
with $(A\otimes B)(x\otimes y)=(Ax)\otimes (By)$
for all~$x\in\scrH$
and~$y\in\scrK$.
\begin{point}{60}{Proof}%
In view of~\sref{hilb-tensor-universal-property}
the only thing we need
to prove is that the bilinear map
$\otimes\circ (A\times B)\colon \scrH\times \scrK
\to\scrH\otimes\scrK$
is $\ell^2$-bounded
(for then~$A\otimes B=(\,\otimes\circ (A\times B)\,)_\otimes$.)
So let~$x_1,\dotsc,x_n\in\scrH$
and~$y_1,\dotsc,y_n\in\scrK$
be given,
and note that
\begin{alignat*}{3}
	\|\sum_{i}(\otimes\circ(A\times B))(x_i,y_i)\,\|^2
	\ &=\ 
	\|\sum_{i} (Ax_i)\otimes (By_i)\|^2
	\\
	\ &=\ 
	\sum_{i,j} 
	\left<Ax_i,Ax_j\right>
	\,\left<By_i,By_j\right>
	\\
	\ &\leq\ 
	\|A\|^2\|B\|^2
	\sum_{i,j} 
	\left<x_i,x_j\right>
	\left<y_i,y_j\right>,
\end{alignat*}
so~$\otimes\circ(A\times B)$
is bounded by~$\|A\|\|B\|$.
The last step in the display above
is justified by~\sref{mult-completely-monotone},
and the inequalities
$(\,\left<Ax_i,A x_j\right>\,)\leq
(\, \|A\|^2 \left<x_i,x_j\right>\,)$
and
$(\,\left<By_i,B y_j\right>\,)\leq
(\, \|B\|^2 \left<y_i,y_j\right>\,)$.\qed
\end{point}
\end{point}
\begin{point}{70}[special-tensor]{Theorem}%
\index{spacial tensor product}
Let~$\scrA$
and~$\scrB$
be von Neumann algebras
of bounded operators
on Hilbert spaces~$\scrH$ and~$\scrK$,
respectively.
Sending operators $A\in\scrA$ and~$B\in\scrB$
to $A\otimes B\colon \scrH\otimes \scrK\to\scrH\otimes \scrK$
from~\sref{hilb-tensor-functor} gives
an miu-bilinear 
map 
\begin{equation*}
	\otimes\colon \scrA\times \scrB\longrightarrow
	\scrB(\scrH\otimes\scrK).
\end{equation*}
Letting~$\scrT$ 
be the von Neumann subalgebra
of $\scrB(\scrH\otimes\scrK)$
generated by the range of~$\otimes$,
the restriction~$\gamma\colon \scrA\times \scrB\to\scrT$
of~$\otimes$
is a tensor product of~$\scrA$ and~$\scrB$.
\begin{point}{80}{Proof}%
We'll check that the three conditions of~\sref{tensor} hold;
we leave it to the reader to verify that~$\otimes$
is miu-bilinear.
\begin{point}{90}{Condition~\ref{tensor-1}}%
The range of~$\gamma$ being the same as the range of~$\otimes$
generates~$\scrT$ 
simply by the way~$\scrT$ was defined.
\end{point}
\begin{point}{100}{Condition~\ref{tensor-2}}%
Let~$\sigma\colon \scrA\to\C$
and~$\tau\colon \scrB\to\C$ be np-maps.
We must find an np-functional~$\omega$ on~$\scrT$
with $\omega(A\otimes B)=\sigma(A)\tau(B)$
for all~$A\in \scrA$, $B\in \scrB$.
Note that by~\sref{normal-functional}
$\sigma$ and~$\tau$ are of the form
$\sigma\equiv\sum_n \left<x_n,(\,\cdot\,)x_n\right>$
and~$\tau\equiv\sum_n \left<y_n,(\,\cdot\,)y_n\right>$
for some~$x_1,x_2,\dotsc\in\scrH$
and~$y_1,y_2,\dotsc\in\scrK$
with~$\sum_n\|x_n\|^2<\infty$
and~$\sum_m\|y_m\|^2 <\infty$.
So as~$\sum_{n,m}\|x_n\otimes y_m\|^2\equiv\sum_n \|x_n\|^2\,
\sum_m\|y_m\|^2 <\infty$,
we can define
an np-functional~$\omega$ on~$\scrT$
by $\omega(T):= \sum_{n,m} \left<x_n\otimes y_m,
T\,x_n\otimes y_m\right>$;
which does the job:  $\omega(A\otimes B)
= \sum_{n,m} \left<x_n, Ax_n\right> \left<y_m,By_m\right>
= \sigma(A)\tau(B)$
for all~$A\in\scrA$ and~$B\in\scrB$.
\end{point}
\begin{point}{110}{Condition~\ref{tensor-3}}%
It remains to be shown that the
product functionals
on~$\scrT$ form a faithful collection.
These functionals are---as we've just seen---all
of the form
$\sum_{m,n}\left<x_n\otimes y_n,(\,\cdot\,)\,x_n\otimes y_m\right>$
for some $x_1,x_2,\dotsc\in\scrH$
and~$y_1,y_2,\dotsc\in\scrK$
(and, conversely, it's easily seen
that a functional
of that form
is a product functional).
It suffices,
then,
to show that the subset of
product functionals
of the  form $\left<x\otimes y,(\,\cdot\,)x\otimes y\right>$
where~$x\in\scrH$ and~$y\in\scrK$ is faithful.
To this end,
let~$T\in\scrT_+$
with $\left< x\otimes y, Tx\otimes y\right>=0$
for all~$x\in\scrH$ and~$y\in\scrK$ be given
in order to show that~$T=0$.
Note that since
$\|\sqrt{T}\,x\otimes y\|^2=\left<x\otimes y,T\,x\otimes y\right>
=0$, and so~$\sqrt{T}\,x\otimes y=0$
for all~$x\in\scrH$, $y\in \scrK$,
we have~$\sqrt{T}=0$ (since the linear
span of the $x\otimes y$ is dense in~$\scrH\otimes\scrK$),
and thus~$T=0$.\qed
\end{point}
\end{point}
\end{point}
\begin{point}{120}{Exercise}%
\index{tensor product!of von Neumann algebras!exists}%
Given von Neumann algebras~$\scrA$ and~$\scrB$
(which are not a priori represented on Hilbert spaces)
construct a tensor product~$\gamma\colon \scrA\times \scrB\to\scrT$
of~$\scrA$ and~$\scrB$
using~\sref{ngns} and~\sref{special-tensor}.
\end{point}
\end{parsec}
\subsection{Universal Property}
\begin{parsec}{1120}%
\begin{point}{10}%
Before we bring our categorical faculties
to bear upon the tensor product for von Neumann algebras
we quickly review
the (algebraic)
tensor product of plain vector spaces~$V$ and~$W$ first ---
    it is a vector space~$\Define{V\odot W}$%
\index{*tensora@$\odot$, algebraic tensor product}%
\index{tensor product!algebraic = of vector spaces}
equipped
    with a bilinear mapping~$\Define{\odot}\colon V\times W\to V\odot W$
which is  universal  
in the sense that for every bilinear mapping~$\beta\colon V\times W\to Z$
into some vector space~$Z$
    there is a unique linear map~$\Define{\beta_\odot}\colon V\odot W\to Z$%
\index{*tensorasup@$\beta_\odot$}
with~$\beta_\odot(v\odot w)=\beta(v,w)$
for all~$v\in V$ and~$w\in W$.
This property
uniquely determines the algebraic tensor product in the sense
that for any bilinear map~$\mathbin{\tilde\odot}\colon
 V\times W\to V\mathbin{\tilde\odot} W$
 into a vector space~$V\mathbin{\tilde \odot} W$
which shares this property
there is a unique linear isomorphism $\varphi\colon V\odot W\to V
\mathbin{\tilde \odot} W$
with $\varphi(v\odot w) = v\mathbin{\tilde\odot} w$
for all~$v\in V$ and~$w\in W$.

In fact, one may take this property as a neat abstract 
definition of the algebraic
tensor product.
However, to  see that the darn thing actually exists,
one still needs a concrete description
such as this one:
take given a basis~$B$ of~$V$ and a basis~$C$ of~$W$
the bilinear map $\odot$ on~$V\times W$
to the vector space $(B\times C)\cdot \C$ with basis~$B\times C$
determined by~$b\odot c = (b,c)$
for~$b\in B$ and~$c\in C$.
This shows us not only that the algebraic tensor product
exists,
but also 
that~$\odot$ is injective (among other things).

This is all, of course, well known,
and we already saw in~\sref{hilb-tensor-universal-property}
that the tensor product for Hilbert spaces
has a similar universal property;
the interesting thing here is that
with some work one can see that
a tensor product
$\gamma\colon \scrA\times \scrB\to \scrT$
of von Neumann algebras
$\scrA$ and~$\scrB$
has a similar universal property too!
We'll see that any bilinear map $\beta\colon \scrA\times \scrB\to\scrC$
into a von Neumann algebra~$\scrC$
which is sufficiently regular
extends uniquely along~$\gamma$ to a ultraweakly continuous
map $\beta_\gamma \colon  \scrT\to\C$,
where regular will mean 
that the extension $\beta_\odot\colon \scrA\odot\scrB\to\scrC$
from the \emph{algebraic} tensor product
is ultraweakly continuous
and bounded
with respect to the norm and ultraweak topology
induced on~$\scrA\odot\scrB$ by~$\scrT$
via~$\gamma$.

To prevent a circular description
here,
we'll first describe the norm and ultraweak topology
that the tensor product induces on~$\scrA\odot \scrB$
directly,
which turns out to be independent (as it should)
from the choice of~$\gamma$.
This description
is essentially based
on the fact that the product functionals on~$\scrT$
are centre separating;
and that this determines both norm and 
ultraweak topology
is just a general 
observation concerning
centre separating sets, as we saw in~\sref{vn-center-separating-fundamental}.
\end{point}
\begin{point}{20}[tensor-extra]{Definitions}%
Let~$\scrA$ and~$\scrB$ be von Neumann algebras.
\begin{enumerate}
\item
A \Define{basic functional}%
\index{functional!basic, on~$\scrA\odot\scrB$}
is 
a map $\omega \colon \scrA\odot\scrB\to\C$
with
$\omega\equiv (\sigma\odot \tau)(t^*(\,\cdot\,)t)$
for some np-maps
$\sigma\colon \scrA\to\C$, $\tau\colon \scrB\to\C$,
and
$t\in \scrA\odot\scrB$.

A \Define{simple functional}%
\index{functional!simple, on~$\scrA\odot\scrB$}
is a finite sum of basic functionals.
\item
Each basic functional $\omega \colon \scrA\odot\scrB\to\C$
gives us an operation~$\Define{[\,\cdot\,,\,\cdot\,]_\omega}$,
that will turn out to be an inner product in~\sref{basic-state-inner-product}
 by
$\Define{[s,t]_\omega}:=\omega(s^*t)$
(cf.~\sref{state-inner-product}),
and an associated semi-norm
denoted by~$\Define{\|t\|_\omega}:=[t,t]_\omega^{\smash{\nicefrac{1}{2}}} 
= \omega(t^*t)^{\nicefrac{1}{2}}$.

The \Define{tensor product norm}
on~$\scrA\odot \scrB$
is the norm (see~\sref{tensor-product-norm})
given by
\begin{equation*}
	\textstyle
	\Define{\|t\|}\ =\ \sup_\omega \|t\|_\omega,
\end{equation*}
where~$\omega$ ranges over all basic functionals
on~$\scrA\odot\scrB$
with~$\omega(1)\leq 1$.
\item
Note that having endowed~$\scrA\odot \scrB$
with the tensor product norm
we can speak of bounded functionals on~$\scrA\odot \scrB$,
and the operator norm on them;
and note that
the basic and simple  functionals are bounded.

The \Define{ultraweak tensor product topology}%
\index{ultraweak tensor product topology}
is the least topology on~$\scrA\odot\scrB$
that makes all operator norm limits
of simple functionals continuous.
\item
A bilinear map $\beta\colon \scrA\times \scrB\to\scrC$
to a von Neumann algebra~$\scrC$
is called
\begin{enumerate}
\item
(continues the list from~\sref{bilinear-basic})
\item
\Define{bounded}%
\index{bilinear map!bounded}
when the unique extension $\beta_\odot \colon \scrA\odot \scrB\to \scrC$
is bounded,
\item 
\Define{\textbf{n}ormal}%
\index{bilinear map!normal}%
\index{normal!bilinear map}
when~$\beta_\odot$
is continuous with respect to the ultraweak tensor product topology 
on~$\scrA\odot\scrB$
and the ultraweak topology on~$\scrC$,
\item
\Define{\textbf{c}ompletely \textbf{p}ositive}%
\index{bilinear map!completely positive}%
\index{completely positive!bilinear map}%
\index{positive!completely $\sim$ bilinear map}
when~$\sum_{i,j}c_i^*\,\beta(a_i^*a_j,b_i^*b_j)\,c_j\geq 0$
for all tuples $a_1,\dotsc,a_N\in\scrA$,
$b_1,\dotsc,b_N\in\scrB$,
and $c_1,\dotsc,c_N\in\scrC$.
\end{enumerate}
\end{enumerate}
\spacingfix%
\begin{point}{21}{Warning}%
While we'll be able to see shortly that any bilinear map
$\beta\colon \scrA\times\scrB\to\scrC$
between von Neumann algebras
that is normal
is jointly ultraweakly continuous as well,
(as a consequence of~\sref{tensor-basic},)
we do not know---but doubt---that the converse holds.
So to clearly differentiate between these two possibly different properties,
we decided to call the former ``normality'' instead of 
    the more likely ``ultraweak continuity'',
stretching the use of the word ``normal'' beyond its usual 
    domain of positive (bilinear) maps.
\end{point}
\end{point}
\begin{point}{30}[product-state-positive]{Lemma}%
Given $C^*$-algebras~$\scrA$ and~$\scrB$
we have~$(\sigma\odot \tau) (t^*t)\geq 0$
for all  $t\in\scrA\odot\scrB$
and p-maps  $\sigma\colon \scrA\to\C$
and~$\tau\colon \scrB\to\C$.
\begin{point}{40}{Proof}%
Note that writing~$t\equiv \sum_n a_n \odot b_n$,
where~$a_1,\dotsc,a_N\in \scrA$, $b_1,\dotsc,b_N\in \scrB$,
we have
$(\sigma\odot\tau)(t^*t)
= \sum_{n,m} \sigma(a_n^*a_m)\,\tau(b_n^*b_m)$.
Since~$(a_n^*a_m)$
is a positive matrix over~$\scrA$,
and~$\sigma\colon \scrA\to\C$
is completely positive (by~\sref{cp-commutative}),
the matrix~$(\sigma(a_n^*a_m))$ is positive.
Since~$(\tau(b_n^*b_m))$
is positive by the same token,
the entrywise product
$(\,\sigma(a_n^*a_m)\,\tau(b_n^*a_m)\,)$
is positive too (by~\sref{schur}).
Whence
$(\sigma\odot\tau)(t^*t)
=\ \sum_{n,m} \sigma(a_n^*a_m)\,\tau(b_n^*b_m) \geq 0$.\qed
\end{point}
\end{point}
\begin{point}{50}[basic-state-inner-product]{Exercise}%
Use~\sref{product-state-positive} 
to show that
 $[\,\cdot\,,\,\cdot\,]_\omega$
from~\sref{tensor-extra}
is an inner product.
\end{point}
\begin{point}{60}{Lemma}%
Product functionals on~$\scrA\odot\scrB$
formed from 
separating
collections~$\Omega$ and~$\Xi$ 
of linear functionals
on $C^*$-algebras~$\scrA$ and~$\scrB$,
respectively,
are separating
in the sense that given~$t\in\scrA\odot\scrB$
the condition that $(\sigma\odot \tau)(t)=0$
for all~$\sigma\in \Omega$ and~$\tau\in\Xi$
entails that~$t=0$.
\begin{point}{70}{Proof}%
Write~$t\equiv \sum_n a_n\odot b_n$
for some  $a_1,\dotsc,a_N\in\scrA$
and~$b_1,\dotsc,b_N\in \scrB$.
	Note that (by replacing them if necessary)
we may assume that the~$a_1,\dotsc,a_N$
are linearly independent.
Let~$\tau\in\Xi$ be given.
Since~$0=(\sigma\odot \tau)(t)
= \sum_n\sigma(a_n)\tau(b_n)
= \sigma(\,\sum_n a_n\tau(b_n)\,)$
for all $\sigma$ from the separating collection~$\Omega$,
we have~$0=\sum_n a_n\tau(b_n)$,
and so---$a_1,\dotsc,a_N$ being linearly independent---we get
 $0=\tau(b_1)=\dotsb = \tau(b_N)$.
Since this holds for any~$\tau$
in the separating collection~$\Xi$
we get~$0=b_1=\dotsb=b_N$,
and thus~$t=\sum_n a_n\odot b_n=0$.\qed
\end{point}
\end{point}
\begin{point}{80}[tensor-product-norm]{Exercise}%
Show that the tensor product norm
from~\sref{tensor-extra}
is, indeed, a norm.
\end{point}
\begin{point}{90}[product-functional]{Exercise}%
Note that given np-functionals 
$\sigma \colon \scrA\to\C$
and~$\tau\colon \scrB\to\C$
on von Neumann algebras,
the functional $\sigma\odot\tau\colon \scrA\odot\scrB\to\C$
is ultraweakly continuous and bounded,
almost by definition.

Show that $f\odot g$ is
bounded and ultraweakly continuous too
for all~$f\in\scrA_*$ and~$g\in\scrB_*$
(perhaps using~\sref{luws}).
\end{point}
\begin{point}{100}[tensor-basic]{Exercise}%
We're going to show that 
the ultraweak tensor product topology
and tensor product norm from~\sref{tensor-extra}
actually describe the norm and ultraweak topology
on~$\scrA\odot\scrB$ induced by a tensor product
$\scrA\times\scrB\to\scrT$ (via~$\gamma_\odot$)
by establishing
the two closely related facts that
$\gamma_\odot\colon\scrA\odot\scrB\to\scrT$
is an isometry
and an ultraweak embedding,
and that certain
functionals 
$\omega\colon \scrA\odot\scrB\to\C$
can be extended uniquely to~$\scrT$ along~$\gamma_\odot$.
\begin{enumerate}
\item
Show using~\sref{vn-center-separating-fundamental}
that the collection~$\Omega$
of np-functionals
on~$\scrT$ of the form 
$\gamma(\sigma,\tau)(\gamma_\odot(s)^*(\,\cdot\,)\gamma_\odot(s))$,
where~$\sigma\colon \scrA\to\C$,
$\tau\colon \scrB\to\C$
are np-functionals
and~$s\in\scrA\odot\scrB$,
is order separating,
and that
every np-functional on~$\scrT$
is the operator norm limit of finite sums
of functionals from~$\Omega$.

Show that~$\omega\circ \gamma_\odot$ is a basic functional
(see \sref{tensor-extra})
for every~$\omega\in\Omega$,
and that every basic functional is of this form
for some unique~$\omega\in\Omega$.

\item
Show that the subset~$\Omega_1$ of~$\Omega$
of unital maps is order separating,
and so determines the norm on~$\scrT$
via~$\|a\|^2=\|a^*a\| = \sup_{\omega\in\Omega_1} \omega(a^*a)$
for all~$a\in\scrT$
(see~\sref{order-separating-norm}).

Prove that $\|\gamma_\odot(s)\|
=\sup_{\omega\in\Omega_1} \omega(s^*s)^{\nicefrac{1}{2}}
=\sup_{\omega\in\Omega_1}
\|s\|_{\omega\circ \tau_\odot}=\|s\|$
for all~$s\in\scrA\odot\scrB$,
and 
conclude
that~$\gamma_\odot$ is an isometry.

\item
Show that~$\|f\circ\gamma_\odot\|\leq \|f\|$
for every $f\in\scrT_*$,
and deduce from this
that when~$\omega\colon \scrT\to\C$
is an np-functional
its restriction $\omega\circ\gamma_\odot$
is the operator norm limit
of simple functionals on~$\scrA\odot\scrB$
implying that~$\omega\circ\gamma_\odot$---and
thus $\gamma_\odot$ itself---is ultraweakly continuous.

\item
In order to show that~$\gamma_\odot$ is an ultraweak embedding,
we'll need the equality $\|f\circ \gamma_\odot\|=\|f\|$
for all~$f\in\scrT_*$.

In order to show this in turn,
recall 
(from \sref{polar-decomposition-of-functional})
that there is a partial isometry~$u$ in~$\scrT$
with~$f(u)=\|f\|$ (see~\sref{functional-norm}). 

Show that given~$\varepsilon>0$ 
there is a net~$(s_\alpha)_\alpha$
in~$\scrA\odot\scrB$
with $\|s_\alpha\|\leq 1+\varepsilon$ for all~$\alpha$
such that $\gamma_\odot(s_\alpha)$
converges ultrastrongly to~$t$ as~$\alpha\to\infty$
(cf.~\sref{dense-subalgebra}).

Deduce that $\|f\|=f(u)=\left|f(u)\right|
=\lim_\alpha \left|f(\gamma_\odot(s_\alpha))\right|
\leq \|f\circ \gamma_\odot\| (1+\varepsilon)$,
and
conclude that~$\|f\|=\|f\circ\gamma_\odot\|$.

\item
Show that any functional $\omega'\colon \scrA\odot\scrB\to\C$
that is the operator norm limit of simple functionals
on~$\scrA\odot\scrB$
can be extended uniquely along~$\gamma_\odot$
to an np-functional on~$\scrT$
(using the fact that the operator
norm limit of np-functionals is an np-functional again,
see~\sref{predual-complete}).

Deduce from this that~$\gamma_\odot$ is a ultraweak topological embedding.

(Note that by~\sref{vn-extension} any 
bounded ultraweakly continuous functional
on~$\scrA\odot\scrB$ can be extended uniquely
to a normal functional on~$\scrT$.)
\end{enumerate}
\spacingfix%
\end{point}
\begin{point}{110}[tensor-universal-property]{Theorem}%
\index{tensor product!of von Neumann algebras!universal property}
A tensor product~$\gamma\colon \scrA\times\scrB\to\scrT$
of von Neumann algebras~$\scrA$ and~$\scrB$
has this universal property:
for every normal bounded bilinear map $\beta\colon \scrA\times \scrB\to\scrC$
to a von Neumann algebra~$\scrC$
there is a unique ultraweakly continuous
map  $\Define{\beta_\gamma}\colon \scrT
\to \scrC$%
\index{*tensorsup@$\beta_\gamma$}
with $\beta_\gamma\circ \gamma  = \beta$.
Moreover, $\|\beta_\gamma\|=\|\beta_\odot\|$.
\begin{point}{120}{Proof}%
Since $\beta_\odot\colon \scrA\odot\scrB\to\scrC$
is ultraweakly continuous and bounded,
and
$\scrA\odot\scrB$
can by~\sref{tensor-basic} be
considered an ultraweakly dense $*$-subalgebra
of~$\scrT$ via~$\gamma_\odot$,
the theorem follows from~\sref{vn-extension}
except for some trivial details.\qed
\end{point}
\end{point}
\end{parsec}%
\begin{parsec}{1130}%
\begin{point}{10}%
\index{completely positive!bilinear map}%
We'll need some observations
concerning completely positive bilinear maps.
\end{point}
\begin{point}{20}{Exercise}%
\index{Schur's Product Theorem}
Show that a mi-bilinear
map $\beta\colon \scrA\times\scrB\to\scrC$ between von Neumann 
algebras is completely positive.
\end{point}
\begin{point}{30}{Notation}%
Given a bilinear map
$\beta\colon \scrA\times\scrB\to\scrC$
between von Neumann algebras,
we define
$\Define{M_N\beta}\colon M_N\scrA\times M_N\scrB\to M_N\scrC$
by  $(M_N\beta)(A,B) = 
(\beta(A_{ij},B_{ij}))_{ij}$%
\index{Mnbeta@$M_n\beta$, for bilinear~$\beta$}
for each~$N$.
\end{point}
\begin{point}{40}[cp-bilinear]{Exercise}%
Show that for a bilinear map $\beta\colon \scrA\times\scrB\to\scrC$
between von Neumann algebras
the following are equivalent.
\begin{enumerate}
\item
$\beta$ is completely positive.
\item
$M_N\beta$ is completely positive
for each~$N$.
\item
$(M_N\beta)(A,B)\geq 0$
for all~$A\in M_N(\scrA)_+$, $B\in M_N(\scrB)_+$ and~$N$.
\end{enumerate}
Deduce 
as a corollary that $h\circ \beta \circ (f\times g)$
is completely positive
when~$f\colon \scrA'\to\scrA$,
$g\colon \scrB'\to\scrB$
and~$h\colon \scrC\to\scrC'$
are cp-maps between von Neumann algebras.
\end{point}
\end{parsec}
\begin{parsec}{1140}%
\begin{point}{10}[tensor-universal-property-extra]{Exercise}%
\index{tensor product!of von Neumann algebras!universal property}
Let $\gamma\colon \scrA\times \scrB\to\scrT$
be a tensor product
of von Neumann algebras,
  $\beta\colon \scrA\times \scrB\to\scrC$
 a normal bounded bilinear map,
 and~$\beta_\gamma\colon \scrT\to\scrC$
 its extension along~$\gamma_\odot$
from~\sref{tensor-universal-property}.
Show that
\begin{enumerate}
\item
$\beta_\gamma$
is multiplicative iff~$\beta$ is multiplicative
(see~\sref{tensor-extra});
\item 
$\beta_\gamma$ 
is involution preserving iff~$\beta$ is involution preserving;
\item
$\beta_\gamma$ is unital
iff $\beta$ is unital;
\item
$\beta_\gamma$ is positive
iff
$\sum_{i,j} \beta(a_i^*a_j,b_i^*b_j) \geq 0$
for all tuples $a_1,\dotsc,a_N$
from~$\scrA$ and $b_1,\dotsc,b_N$ from $\scrB$;
\item
$\beta_\gamma$ is completely positive iff
$\beta$ is completely positive.
\end{enumerate}
\spacingfix%
\end{point}%
\begin{point}{20}[tensor-uniqueness]{Exercise}%
\index{tensor product!of von Neumann algebras!uniqueness}
Show that the tensor product of von Neumann algebras~$\scrA$
and~$\scrB$ is unique in
the sense
that when~$\gamma\colon \scrA\times \scrB\to\scrT$
and~$\gamma'\colon \scrA\times \scrB\to\scrT'$
are tensor products of~$\scrA$ and~$\scrB$,
then there is a unique
nmiu-isomorphism $\varphi\colon \scrT\to\scrT'$
with $\varphi(\gamma(a,b))=\gamma'(a,b)$
for all~$a\in\scrA$ and~$b\in\scrB$.
\end{point}
\end{parsec}
\subsection{Functoriality}
\begin{parsec}{1150}%
\begin{point}{10}{Notation}%
Now that we've established
that that the tensor product
of von Neumann algebras~$\scrA$ and~$\scrB$
exists and is unique (up to unique nmiu-isomorphism)
we just pick one and denote it by $\Define{\otimes}\colon
\scrA\times\scrB\to\Define{\scrA\otimes\scrB}$.%
\index{*tensor@$\otimes$, tensor product!$\scrA\otimes\scrB$, of von Neumann algebras}%
\index{*tensor@$\otimes$, tensor product!$a\otimes b$,
of elements of von Neumann algebras}
\end{point}
\begin{point}{20}[tensor-functorial]{Proposition}%
Given ncp-maps $f\colon \scrA\to\scrC$
and $g\colon \scrB\to\scrD$
between von Neumann algebras
there is a unique ncp-map
$\Define{f\otimes g}\colon \scrA\otimes\scrB\to\scrC\otimes\scrD$
with 
\begin{equation*}
	(f\otimes g)(a\otimes b) \ =\ f(a)\otimes f(b)
\end{equation*}%
\index{*tensor@$\otimes$, tensor product!$f\otimes g$, of np-maps}
for all~$a\in\scrA$ and~$b\in\scrB$.
Moreover,
\begin{enumerate}
\item
$f\otimes g$ is multiplicative 
when~$f$ and~$g$ are multiplicative;
\item
$f\otimes g$ is involution preserving
when~$f$ and~$g$ are involution preserving; and
\item
$f\otimes g$ is (sub)unital 
when~$f$ and~$g$ are (sub)unital.
\end{enumerate}
\spacingfix%
\begin{point}{30}{Proof}%
As uniqueness of~$f\otimes g$ is rather obvious,
we leave it at that.
To establish
existence of~$f\otimes g$,
it suffices to show
that the bilinear map~$\beta\colon \scrA\times\scrB\to\scrC\otimes\scrD$
given by~$\beta(a,b)=f(a)\otimes g(b)$,
which is completely positive by~\sref{cp-bilinear},
is bounded and normal;
because then we may take $f\otimes g:=\beta_\otimes$%
\index{*tensorsupexact@$\beta_\otimes$}
as in~\sref{tensor-universal-property}
and all
the properties claimed for~$f\otimes g$ will then  follow 
with the very least of effort
from~\sref{tensor-universal-property-extra}.

To see that~$\beta$ is bounded,
we'll prove that~$\|\beta_\odot(s)\| \leq \|f\|\|g\| \|s\|$
given an element~$s$ of~$\scrA\otimes \scrB$,
and for this
it suffices (by the definition
of the tensor product norm, \sref{tensor-extra}) 
to show that
$\omega(\beta_\odot(s)^*\beta_\odot(s))
\leq \|f\|^2\|g\|^2\|s\|^2$
given a basic functional~$\omega$
on~$\scrA\odot\scrB$ with~$\omega(1)\leq 1$.
We'll prove in a moment that
$\|\omega\circ\beta_\odot\|\leq \|f\|\|g\|$
and~$\beta_\odot(s)^*\beta_\odot(s)\leq
\|f\|\|g\|\beta_\odot(s^*s)$,
because with these two claims
we get 
$\omega(\beta_\odot(s)^*\beta_\odot(s))
\leq \|f\|\|g\|\omega(\beta_\odot(s^*s))
\leq \|f\|\|g\| \|\omega\circ \beta_\odot\| \|s\|^2
\leq \|f\|^2\|g\|^2\|s\|^2$
 --- which is the result desired.

Concerning the first
promise, that $\|\omega\circ \beta_\odot\|\leq \|f\|\|g\|$,
note that writing
$\omega\equiv (\sigma\odot\tau)(t^*(\,\cdot\,)t)$,
where~$\sigma$ and~$\tau$
are np-maps on~$\scrC$ and~$\scrD$, respectively,
and~$t\equiv\sum_{ij} c_i \odot d_i$
is from~$\scrC\odot\scrD$,
we have 
\begin{equation*}
	\omega\circ\beta_\odot
= \textstyle \sum_{ij} \sigma(c_i^* f(\,\cdot\,)c_j)\,\odot\,
\tau(d_i^* g(\,\cdot\,)d_j),
\end{equation*}
and so~$\omega\circ \beta_\odot$
is ultraweakly continuous and bounded
by~\sref{product-functional},
because the $\sigma(c_i^*f(\,\cdot\,)c_j)$
and~$\tau(d_i^*g(\,\cdot\,)d_j)$
are bounded ultraweakly continuous functionals.
Although the bound for~$\omega\circ\beta_\odot$ 
thus obtained is in all probability nowhere near~$\|f\|\|g\|$,
it does allow
us by~\sref{tensor-universal-property}
to extend~$\omega\circ\beta_\odot$
to an ultraweakly continuous 
functional $\omega':=(\omega\circ\beta)_\otimes$ on~$\scrC\otimes \scrD$
with the same norm, $\|\omega'\|=\|\omega\circ\beta_\odot\|$.
Since this extension~$\omega'$
is completely positive
(because~$\beta$ and thus~$\omega\circ \beta$ are completely positive,
see~\sref{cp-bilinear})
its norm is by~\sref{cp-russo-dye} given by~$\|\omega'\|=\omega'(1)\equiv
\omega(f(1)\otimes g(1))\leq \|f\|\|g\|$,
where we used that~$\omega(1)\leq 1$.
Thus $\|\omega\circ \beta_\odot\|=\|\omega'\| \leq \|f\|\|g\|$,
as was claimed.

Incidentally,
since each~$\omega\circ \beta_\odot$
is ultraweakly continuous,
so is~$\beta_\odot$,
and thus~$\beta$ is normal.
The only thing that remains
is to make good on our last promise,
that~$\beta_\odot(s)^*\beta_\odot(s)\leq 
\|f\|\|g\|\beta_\odot(s^*s)$.
To this end,
write $s\equiv \sum_i a_i\odot b_i$,
and consider
the matrices~$A$ and~$B$ given by
\begin{equation*}
A\ :=\  \begin{pmatrix}
a_1 & a_2 & \dotsb & a_n \\
0 & 0 & \dotsb & 0 \\
\vdots & \vdots & \ddots&\vdots \\
0 & 0 & \dotsb & 0
\end{pmatrix}
\qquad\qquad
B\ :=\  \begin{pmatrix}
b_1 & b_2 & \dotsb & b_n \\
0 & 0 & \dotsb & 0 \\
\vdots & \vdots & \ddots&\vdots \\
0 & 0 & \dotsb & 0
\end{pmatrix},
\end{equation*}
and the cp-map $h\colon M_n(\scrC\otimes \scrD)\to \scrC\otimes \scrD$
given by~$h(C)= \left<(1,\dotsc,1), C(1,\dotsc,1)\right> =\sum_{ij}C_{ij}$.
We make these arrangements
so that we
may apply the inequality
$(M_nf)(A)^*  (M_nf)(A)
\leq \|(M_nf)(1)\| (M_nf)(A^*A)$
easily derived from~\sref{cp-cs}.
Indeed,
noting also~$\|(M_nf)(1)\|=\|f(1)\|=\|f\|$,
we have
\begin{alignat*}{3}
	\beta_\odot(s)^*\beta_\odot(s) 
\ &=\ \textstyle \sum_{ij}f(a_i)^*f(a_j)\otimes g(b_i)^*g(b_j)\\
\ &=\ h(\ 
(M_nf)(A)^* (M_nf)(A) \ \ (M_n \otimes)\ \ 
(M_ng)(B)^*(M_ng)(B) \ ) \\
\ &\leq\  
\|f\|\|g\| 
h(\ (M_nf)(A^*A) \ \ (M_n \otimes)\ \ 
(M_ng)(B^*B) \ ) 
\\
\ &=\  
\|f\|\|g\| \,\textstyle \sum_{ij} f(a_i^*a_j)\otimes g(b_i^*b_j) \\
\ &= \ \|f\|\|g\|\,\beta_\odot(s^*s),
\end{alignat*}
which concludes this proof.\qed
\end{point}
\end{point}
\begin{point}{40}[tensor-functor]{Exercise}%
Show that the assignments
$(\scrA,\scrB)\mapsto \scrA\otimes \scrB$,
and~$(f,g)\mapsto f\otimes g$
    give a bifunctor~$\Define{\otimes}\colon \Cat{C}\times\Cat{C}
\to\Cat{C}$
\index{*tensor@$\otimes$, tensor product!bifunctor 
on~$\W{miu},\dotsc$}%
\index{tensor product!of von Neumann algebras!functorial}
where~$\Cat{C}$
can be~$\W{miu}$, $\W{cp}$, $\W{cpu}$
or $\W{cpsu}$.
\end{point}
\begin{point}{50}[tensor-injective]{Proposition}%
Given injective nmiu-maps
$f\colon \scrA\to\scrC$ and~$g\colon \scrB\to\scrD$,
the nmiu-map $f\otimes g\colon \scrA\otimes \scrB\to
\scrC\otimes\scrD$ is injective.
\begin{point}{60}{Proof}%
The trick
is to consider the von Neumann subalgebra~$\scrT$
generated by
the elements of~$\scrC\otimes \scrD$
of the form~$f(a)\otimes g(b)$
where~$a\in\scrA$ and~$b\in\scrB$,
and to show that  the miu-bilinear map
$\gamma\colon \scrA\times\scrB\to\scrT$
given by~$\gamma(a,b)=f(a)\otimes g(b)$
is a tensor product of~$\scrA$ and~$\scrB$.
Indeed,
if this is achieved,
then
there is, by~\sref{tensor-uniqueness},
a unique nmiu-map $\varphi\colon \scrA\otimes\scrB\to\scrT$
with $\varphi(a\otimes b)=\gamma(a,b)
=f(a)\otimes g(b)$,
so that the following diagram commutes.
\begin{equation*}
\xymatrix@C=4em{
\scrA\times \scrB
\ar[rr]^-{f\times g}
\ar[rd]^-\gamma
\ar[d]_-\otimes
&
&
\scrC\times\scrD 
\ar[d]^-\otimes
\\
\scrA\otimes\scrB
\ar[r]|-\varphi
&
\scrT
\ar[r]|-\subseteq
&
\scrC\otimes\scrD
}
\end{equation*}
The map on the bottom side of this rectangle above is
none other than~$f\otimes g$,
and is thus,
being
the composition of the isomorphism~$\varphi$
with the inclusion $\scrT\subseteq \scrC\otimes\scrD$,
injective.

It remains to be shown that~$\gamma$
is a tensor product,
that is, obeys the conditions from~\sref{tensor}.
Condition~\ref{tensor-1}
holds simply by definition of~$\scrT$.
To see that~$\gamma$
obeys condition~\ref{tensor-2},
let np-functionals
$\tilde\sigma\colon \scrA\to\C$
and $\tilde\tau\colon \scrB\to\C$ be given;
we must find an np-functional $\gamma(\tilde\sigma,
\tilde\tau)$ on~$\scrT$
with $\gamma(\tilde\sigma,\tilde\tau)(a\otimes b)
= \gamma(a,b)$.

By ultraweak permanence
$\tilde\sigma$ and~$\tilde\tau$
can be extended along~$f$ and~$g$, respectively,
see~\sref{functional-extension},
giving us np-functionals $\sigma\colon \scrC\to\C$
and $\tau\colon\scrD\to\C$
with~$\tilde\sigma = \sigma\circ f$
and $\tilde\tau = \tau\circ g$.
Now simply take $\gamma(\tilde\sigma,\tilde\tau)$
to be the restriction
of $\sigma\otimes \tau$
to~$\scrT$,
which does the job.

Finally,
concerning
condition~\ref{tensor-3},
let~$z$ be a central projection of~$\scrT$
with $\gamma(\tilde\sigma,\tilde\tau)(z)=0$
for all~$\tilde\sigma$
and~$\tilde\tau$ of aforementioned type.
We must show that~$z=0$,
and for this
it suffices to show that
$(\sigma\otimes \tau)(z)=0$
for all  np-functionals
$\sigma$ and~$\tau$ on~$\scrC$ and~$\scrD$,
respectively.
Since for such~$\sigma$ and~$\tau$
we have
$\gamma(\tilde\sigma,\tilde\tau)
(\gamma(a,b))
 = 
 \sigma(f(a))\,\tau(g(b))
 = (\sigma\otimes\tau)(\gamma(a,b))$
 for all~$a\in\scrA$ and~$b\in\scrB$,
 we have $\gamma(\tilde\sigma,\tilde\tau)
 (t) = (\sigma\otimes \tau)(t)$
 for all~$t\in\scrT$,
 and, in particular,
 $0=\gamma(\tilde\sigma,\tilde\tau)(z)
 = (\sigma\otimes\tau)(z)$.
 Hence~$z=0$.\qed
\end{point}
\end{point}
\end{parsec}%
\subsection{Miscellaneous Properties}
\begin{parsec}{1160}%
\begin{point}{10}[product-functional-norm]{Lemma}%
Given von Neumann algebras~$\scrA$ and~$\scrB$,
we have $\|f\otimes g\|=\|f\|\|g\|$
for all~$f\in\scrA_*$ and~$g\in\scrB_*$.%
\index{*tensor@$\otimes$!$f\otimes g$, of normal functionals}
\begin{point}{20}{Proof}%
The trick is to use the polar decomposition
for normal functionals, \sref{polar-decomposition-of-functional}.
On its account we can find partial isometries $u\in\scrA$
and $v\in\scrB$
such that
$f(u(\,\cdot\,))$
and~$g(v(\,\cdot\,))$ are positive,
and
$f\equiv f(uu^*(\,\cdot\,))$,  $g\equiv g(vv^*(\,\cdot\,))$.
Then $u\otimes v$ is a partial isometry
such that $(f\otimes g)((u\otimes v)(\,\cdot\,))$
is positive,
and $f\otimes g
= (f\otimes g)(\,(u\otimes v)\, (u\otimes v)^*\,(\,\cdot\,)\,)$
so that $\|f\otimes g\|=(f\otimes g)(u\otimes v)
= f(u)g(v)=\|f\|\|g\|$ by~\sref{functional-norm}.\qed
\end{point}
\end{point}
\begin{point}{30}[tensor-simple-facts]{Exercise}%
There are some easily obtained facts
concerning the tensor product~$\scrA\otimes\scrB$
of von Neumann algebras
that nevertheless deserve explicit mention.
\begin{enumerate}
\item
Show that~$a\otimes b\geq 0$
for all~$a\in\scrA_+$ and~$b\in\scrB_+$;
and conclude that $a_1\otimes b_1 \leq a_2\otimes b_2$
for all
$a_1\leq a_2$ from~$\scrA$ and $b_1\leq b_2$ from~$\scrB$.
\item
Show that $\|a\otimes b\| = \|a\|\|b\|$
for all~$a\in\scrA$ and~$b\in\scrB$.

Conclude that $\otimes\colon \scrA\times\scrB\to\scrA\otimes \scrB$
is norm continuous.

(Warning: as~$\otimes$ is not linear
this is not entirely trivial.)
\item
Show that $\otimes\colon \scrA_*\times\scrB_*\to(\scrA\otimes\scrB)_*$
is norm continuous
(using~\sref{product-functional-norm}).
\item
Show that~$\otimes\colon \scrA\times \scrB\to \scrA\otimes \scrB$
is ultraweakly continuous.

(Hint: since we
already know that~$\otimes_\odot\colon \scrA\odot\scrB\to\scrA\otimes\scrB$
is ultraweakly continuous, by~\sref{tensor-basic},
an equivalent question 
is whether~$\odot\colon \scrA\times\scrB\to\scrA\odot\scrB$
is ultraweakly continuous,
which may be boiled down
to the fact
that $(a,b)\mapsto \sum_{ij} \sigma(a_i^* a a_j)\,\tau(b_i^* b b_j)\colon
\ \scrA\times \scrB\to\C$ is ultraweakly continuous,
where~$\sigma$ and~$\tau$ are np-functionals
on~$\scrA$ and~$\scrB$, respectively,
and
$a_1,\dotsc,a_n\in\scrA$, and~$b_1,\dotsc,b_n\in\scrB$.)
\item
Show that $a\otimes(\,\cdot\,)\colon \scrB\to \scrA \otimes \scrB$
is a ncp-map for every~$a\in\scrA$,
and that $1\otimes(\,\cdot\,)\colon \scrB\to \scrA\otimes \scrB$
is an nmiu-map.
\end{enumerate}
\spacingfix%
\end{point}%
\begin{point}{31}%
The following observation will come in very handy later on when we prove that
\begin{equation*}
\scrA\otimes (\scrB\otimes \scrC)
    \cong (\scrA\otimes\scrB)\otimes \scrC,
    \qquad\text{and}\qquad
\scrA\otimes (\scrB\oplus\scrC)\cong\scrA\otimes \scrB
    \,\oplus\,\scrA\otimes \scrC,
\end{equation*}
see \sref{associator}, and \sref{tensor-distributes-over-sums}.
\end{point}
\begin{point}{40}[tensor-generation]{Proposition}
Let~$\scrA$ and~$\scrB$ 
be von Neumann algebras.
\begin{enumerate}
\item
\label{tensor-generation-1}
If~$S$ and~$T$ are subsets of~$\scrA$
and~$\scrB$, respectively,
whose linear span is ultraweakly dense,
then the linear span
of $\{\,s\otimes t\colon\,s\in S,\,t\in T\,\}$
is ultraweakly dense in~$\scrA\otimes\scrB$.
\item
\label{tensor-generation-2}
If~$\Omega$ and~$\Theta$ are centre separating collections
of np-functionals on~$\scrA$ and~$\scrB$, respectively,
then~$\{\,\omega\otimes\vartheta\colon\, 
\omega\in\Omega,\,\vartheta\in\Theta\,\}$
is centre separating for~$\scrA\otimes\scrB$.
\end{enumerate}%
\spacingfix%
\begin{point}{50}{Proof}%
Concerning~\ref{tensor-generation-1}:
Let~$S'$ and~$T'$ denote the linear spans of~$S$
and~$T$, respectively.
Since linear combinations of elements of~$\scrA\otimes\scrB$
of the form $a\otimes b$
lie ultraweakly dense in~$\scrA\otimes\scrB$ 
where~$a\in\scrA$ and~$b\in\scrB$,
it suffices to show that such element~$a\otimes b$
is the ultraweak limit of elements of the form $s'\otimes t'$
where~$s'\in S'$ and~$t'\in T'$
(because such~$s'\otimes t'$ are, of course, a linear combinations
of elements of the form $s\otimes t$ where $s\in S$
    and~$t\in T$.)

This is indeed the case
as there are nets $(s_\alpha')_\alpha$
and~$(t_\beta')_\beta$ in~$S'$ and~$T'$
that converge ultraweakly to~$a$ and~$b$, respectively,
and so, 
because~$\otimes$ is ultraweakly continuous by~\sref{tensor-simple-facts},
we see that~$s_\alpha'\otimes t_\beta'$
converges ultraweakly to~$a\otimes b$
as~$\alpha,\beta\to\infty$.

Concerning~\ref{tensor-generation-2},
let~$t$ be a positive element of~$\scrA\otimes \scrB$
with $(\omega\otimes \vartheta)(s^*ts)=0$
for all~$\omega\in\Omega$, $\vartheta\in\Theta$,
and~$s\in\scrA\otimes \scrB$;
we must show that~$t=0$.
For this it suffices
to show that~$(\sigma\otimes\tau)(t)=0$
for all np-functionals  $\sigma\colon \scrA\to\C$
and 
$\tau\colon\scrB\to\C$
(since the product functionals $\sigma\otimes\tau$
	form a faithful collection.)
Now, since~$\Omega$ is centre separating
such $\sigma$ may by~\sref{vn-center-separating-fundamental} be obtained
as operator norm limit of finite sums
of functionals of the
form~$\omega(a^*(\,\cdot\,)a)$
where~$\omega\in\Omega$ and~$a \in\scrA$.
Since an np-functional $\tau\colon \scrB\to\C$ 
can be obtained in a similar fashion
from~$\Theta$,
and~$\otimes \colon \scrA_*\otimes \scrB_*\to(\scrA\otimes\scrB)_*$
is operator norm continuous (by~\sref{tensor-simple-facts}),
we see that a product functional~$\sigma\otimes \tau$
can be obtained as the operator norm limit
of finite sums of functionals
of the form $\omega(a^*(\,\cdot\,)a)\,\otimes\,
\vartheta(b^*(\,\cdot\,)b)
\,\equiv\, (\omega\otimes\vartheta)(\,(a\otimes b)^*\,(\,\cdot\,)
\,(a\otimes b)\,)$;
and since those functionals 
map~$t$ to~$0$,
by assumption,
we conclude that~$(\sigma\otimes \tau)(t)=0$ too.\qed
\end{point}
\end{point}
\begin{point}{60}%
To obtain certain examples
the following characterisation of
the tensor product of von Neumann algebras
proves useful.
\end{point}
\begin{point}{70}[tensor-characterization]{Theorem}%
Given centre separating collections~$\Sigma$ and~$\Gamma$
of np-functionals
on von Neumann algebras~$\scrA$ and~$\scrB$,
respectively,
an miu-bilinear map
$\gamma\colon \scrA\times\scrB\to\scrT$
is a tensor product iff all of the following conditions hold.
\begin{enumerate}
\item
The range of~$\gamma$ generates~$\scrT$.
\item
For all~$\sigma\in \Sigma$ and $\tau\in\Gamma$
the product functional $\gamma(\sigma,\tau)\colon \scrT\to\C$
exists
(see~\sref{tensor})
 and is positive.
\item
The set
$\{\,\gamma(\sigma,\tau)\colon\,\sigma\in\Sigma,\,\tau\in \Gamma\,\}$
is centre separating for~$\scrT$.
\end{enumerate}
\spacingfix%
\begin{point}{80}{Proof}%
A tensor product~$\gamma$
obeys these conditions
by definition and by~\sref{tensor-generation},
so we only need to show that
a~$\gamma$ that obeys these conditions is a tensor product,
and for this it
suffices to show that~$\gamma$ can be extended
to an nmiu-isomorphism $\gamma_\otimes \colon \scrA\otimes\scrB\to \scrT$.
To extend~$\gamma$ to just an miu-map~$\gamma_\otimes$
(to begin with)
it suffices by~\sref{tensor-universal-property}
and~\sref{tensor-universal-property-extra} 
to show that~$\gamma_\odot\colon \scrA\odot\scrB\to\scrT$
is bounded with respect to the tensor product norm
on~$\scrA\odot\scrB$
and continuous with respect to the tensor product
topology on~$\scrA\odot \scrB$ and the ultraweak topology on~$\scrT$.

To see that~$\gamma_\odot$
is bounded, 
let~$t\in \scrA\odot\scrB$ be given;
we'll show
that~$\|\gamma_\odot(t)\|^2
\equiv \|\gamma_\odot(t^*t)\|\leq \|t\|^2$
where~$\|t\|$ is the tensor product norm of~$t$.
Since by~\sref{vn-center-separating-fundamental}
the np-functionals on~$\scrT$
of the form
\begin{equation}
	\label{tensor-characterization-1}	
	\gamma(\sigma,\tau)(\,\gamma_\odot(s)^*\,
(\,\cdot\,)\,\gamma_\odot(s)\,)
\end{equation}
where~$\sigma\in\Sigma$, $\tau\in\Gamma$
and~$s\in\scrA\odot\scrB$,
are order separating,
also with  the restriction
that~$1=\gamma(\sigma,\tau)(\gamma_\odot(s^*s))\equiv
(\sigma\odot\tau)(s^*s)$,
and therefore determine the norm of~$t^*t$ as in~\sref{order-separating-norm},
it suffices to show that
$\gamma(\sigma,\tau)(\gamma_\odot(s)^* \gamma_\odot(t^*t)
\gamma_\odot(s)) \leq \|t\|^2$
given such~$\sigma$, $\tau$, and~$s$
(with $(\sigma\odot \tau)(s^*s)=1$).
But since $
\gamma(\sigma,\tau)(\gamma_\odot(s)^* \gamma_\odot(t^*t)\gamma_\odot(s)) =
(\sigma\odot\tau)(s^* t^*t s)
= \|t\|_{(\sigma\odot\tau)(s^*(\,\cdot\,)s)}^2
\leq \|t\|^2$
by the definition of the tensor product norm
(see~\sref{tensor-extra}),
this is indeed the case.

To see that $\gamma_\odot\colon \scrA\odot\scrB\to\scrT$
is ultraweakly continuous
it suffices to show that  $\omega\circ \gamma_\odot$
is the operator norm limit of finite sums of basic functionals
on~$\scrA\odot\scrB$
(see~\sref{tensor-extra})
given any np-functional $\omega\colon \scrT\to\C$.
Since by~\sref{vn-center-separating-fundamental} such~$\omega$
is the norm limit
of finite sums of
functionals on~$\scrT$
of the form displayed in~\eqref{tensor-characterization-1},
and~$\gamma_\odot$
is bounded,
we may assume without loss of generality
that~$\omega$
itself is as shown in~\eqref{tensor-characterization-1}.
Since $\omega\circ \gamma_\odot 
\equiv (\sigma\odot\tau)(s^*(\,\cdot\,)s)$
is then a basic functional $\gamma_\odot$
is ultraweakly continuous.

Having established boundedness and continuity
of~$\gamma_\odot$
we obtain our  nmiu-map $\gamma_\otimes\colon \scrA\otimes\scrB
\to\scrT$
with~$\gamma_\otimes(a\otimes b)=\gamma(a,b)$
for all~$a\in\scrA$ and~$b\in \scrB$.
To show that~$\gamma$ is a tensor product,
it suffices
to show that~$\gamma_\otimes$
is an nmiu-isomorphism,
and for this,
it suffices
to show that~$\gamma_\otimes$ is a bijection.
In fact, 
we only need to show that~$\gamma_\otimes$ is injective,
because since the elements
of~$\scrT$ of the form~$\gamma(a,b)\equiv\gamma_\otimes(a\otimes b)$
generate~$\scrT$ (by assumption),
and are in the range of~$\gamma_\otimes$
(which is a von Neumann subalgebra of~$\scrT$
by~\sref{injective-nmiu-iso-on-image}),
$\gamma_\otimes$ will be surjective.

To show that~$\gamma_\otimes$ is injective,
it suffices to show that~$\ceil{\gamma_\otimes}\equiv\cceil{\gamma_\otimes}
=1$
(see~\sref{carrier-miu}).
Since the product functionals on $\scrA\otimes\scrB$
of the form $\sigma\otimes \tau$
where~$\sigma\in\Sigma$ and~$\tau\in\Gamma$
are centre separating (by~\sref{tensor-generation}),
and~$\cceil{\gamma_\otimes}$ is central,
it suffices to show that
$(\sigma\otimes\tau)(\,\cceil{\gamma_\otimes}^\perp\,)=0$
given~$\sigma\in\Sigma$ and~$\tau\in\Gamma$.
But this is easy ---
$(\sigma\otimes \tau)(\,\cceil{\gamma_\otimes}^\perp \,)
=
\gamma(\sigma, \tau)(\gamma_\otimes(\,\cceil{\gamma_\otimes}^\perp \,))
= 0$.
Whence~$\gamma$ is a tensor product.\qed
\end{point}
\end{point}
\end{parsec}
\begin{parsec}{1170}%
\begin{point}{10}%
Using the characterization from~\sref{tensor-characterization}
it is pretty 
easy to see that the tensor product distributes
over (infinite) direct sums
(see~\sref{tensor-distributes-over-sums})
after some unsurprising
observations regarding direct sums
(in~\sref{sum-generation}).
\end{point}
\begin{point}{20}[sum-generation]{Exercise}%
Let~$(\scrA_i)_{i\in I}$
be a collection of von Neumann algebras.
\begin{enumerate}
\item
Show that
given a generating subset~$A_i$ 
for each von Neumann algebra~$\scrA_i$
the set $\bigcup_{i\in I} \kappa_i(A_i)$ 
generates~$\bigoplus_{i\in I} \scrA_i$,
where~$\kappa_i\colon \scrA_i\to\bigoplus_{i\in I} \scrA_i$
denotes
the np-map given by~$(\kappa_i(a))_i=a$ and~$(\kappa_i(a))_j=0$
when~$j\neq i$.
\item
Show that
given a centre separating collection~$\Omega_i$
of np-functionals on~$\scrA_i$
for each~$i\in I$
the collection $\{\,\omega\circ \pi_i\colon \,
\omega\in\Omega_i,\,i\in I\,\}$
is centre separating for~$\bigoplus_{i\in I}\scrA_i$.
\end{enumerate}
\spacingfix%
\end{point}%
\begin{point}{30}[tensor-distributes-over-sums]{Proposition}%
Given von Neumann algebras
$\scrA$ and~$(\scrB_i)_{i\in I}$
the bilinear map
\begin{equation*}
	\textstyle
	\gamma\colon 
	\ 
	\scrA \times \bigoplus_i \scrB_i
	\longrightarrow \bigoplus_i \scrA\otimes \scrB_i,\ 
	(a,b)\mapsto (a_i\otimes b)_i
\end{equation*}
is a tensor product.
(Whence~$\scrA\otimes \bigoplus_i \scrB_i
\cong \bigoplus_i\scrA\otimes \scrB_i$.)
\begin{point}{40}{Proof}%
We use~\sref{tensor-characterization}
to show that~$\gamma$ is a tensor product.
Note that
$\gamma$ is clearly miu-bilinear,
and that
the elements of the form $\gamma(a,\kappa_i(b))
= \kappa(a\otimes b)$
from the range of~$\gamma$
where~$a\in \scrA$, $i\in I$, and~$b\in\scrB_i$ 
generate~$\bigoplus_i\scrA \otimes\scrB_i$ by~\sref{sum-generation}.
Further, since
given $i\in I$ and np-functionals
$\sigma\colon \scrA\to\C$ and
$\tau\colon \scrB_i\to\C$
the product functional $\gamma(\sigma,\tau\circ \pi_i)$
exists being simply~$(\sigma\otimes \tau)\circ\pi_i
\colon \bigoplus_i \scrA\otimes \scrB_i\to\C$,
and such product functionals form  a centre separating collection 
by~\sref{sum-generation},
we see that~$\gamma$ is indeed a tensor product.\qed
\end{point}
\end{point}
\end{parsec}
\begin{parsec}{1180}%
\begin{point}{10}%
The tensor interacts
with projections as expected.
\end{point}
\begin{point}{20}{Lemma}%
Let~$\scrA$ and~$\scrB$ be von Neumann algebras.
\begin{enumerate}
\item
We have $\ceil{a\otimes b} = \ceil{a}\otimes \ceil{b}$
for all~$a\in\scrA_+$ and~$b\in\scrB_+$.
\item
We have~$\cceil{a\otimes b} = \cceil{a}\otimes \cceil{b}$
for all~$a\in\scrA$ and~$b\in\scrB$.
\end{enumerate}
\spacingfix%
\begin{point}{30}{Proof}%
Let~$a\in\scrA_+$ and~$b\in\scrB_+$ be given.
Since the map~$(\,\cdot\,)\otimes b\colon \scrA\to \scrA\otimes \scrB$
is np,
$\ceil{a\otimes b}
\smash{\overset{\sref{ncp-ceil}}{=\joinrel=\joinrel=}}
\ceil{\ceil{a}\otimes b}$.
Since similarly~$\ceil{\ceil{a}\otimes b}
=\ceil{\,\ceil{a}\otimes \ceil{b}\,}
\equiv \ceil{a}\otimes\ceil{b}$
using here that~$\ceil{a}\otimes\ceil{b}$
is already a projection,
we get~$\ceil{a}\otimes \ceil{b}=\ceil{a\otimes b}$.

Let~$a\in\scrA$ and~$b\in\scrB$ be given
in order to prove that~$\cceil{a\otimes b}
=\cceil{a}\otimes\cceil{b}$.
Since~$\cceil{a}\otimes 1$ 
commutes with all elements of~$\scrA\otimes \scrB$
of the form~$a'\otimes b'$,
and thus with all elements of~$\scrA\otimes \scrB$,
we see that~$\cceil{a}\otimes 1$ is central.
Since similarly $1\otimes\cceil{b}$ is central,
we see that~$\cceil{a}\otimes \cceil{b}
= (\cceil{a}\otimes 1)\otimes(1\otimes \cceil{b})$
is central too.
Since
in addition 
$\cceil{a}\otimes\cceil{b}$
is a projection,
and $(\cceil{a}\otimes\cceil{b})\,(a\otimes b)
= (\cceil{a}a)\otimes (\cceil{b}b)
= a\otimes b$
we see that~$\cceil{a\otimes b}\leq \cceil{a}\otimes \cceil{b}$
(by definition, see~\sref{central support}).

So all that remains is to show that~$\cceil{a}\otimes\cceil{b}
\leq \cceil{a\otimes b}$.
Recall that~$\cceil{a}= \bigcup_{\tilde{a}\in\scrA}
\ceil{\tilde{a}^* a^*a\tilde{a}}$
by~\sref{cceil-fundamental}.
Using this,  a similar
expression for~$\cceil{b}$,
and~\sref{ncp-union},
we see that~$\cceil{a}\otimes \cceil{b}
=\bigcup_{\tilde{a}\in\scrA}
\bigcup_{\tilde{b}\in\scrB}
\ceil{ \smash{(\tilde{a}^*a^*a\tilde{a})
\otimes (\tilde{b}^*b^*b\tilde{b})} }$,
and so it suffices to show that 
$\ceil{\smash{(\tilde{a}^*a^*a\tilde{a})
\otimes (\tilde{b}^*b^*b \tilde{b})}}
\leq \cceil{a\otimes b}$
given~$\tilde{a}\in\scrA$
and~$\tilde{b}\in\scrB$.
This is indeed the case
since
$\ceil{\smash{(\tilde{a}^*a^*a\tilde{a})
\otimes (\tilde{b}^*b^*b \tilde{b})}} 
=
\ceil{\smash{(\tilde{a}\otimes\tilde{b})^*\,
	(a \otimes b)^*(a\otimes b) \,
(\tilde{a}\otimes \tilde{b})}}
\leq \cceil{a\otimes b}$
(by~\sref{cceil-fundamental}, again.)\qed
\end{point}
\end{point}
\begin{point}{40}[carrier-tensor]{Exercise}%
Let~$f\colon \scrA\to\scrB$
and~$g\colon \scrC\to\scrD$ be np-maps
between von Neumann algebras.
We're going to prove that~$\ceil{f\otimes g}
=\ceil{f}\otimes \ceil{g}$.
\begin{enumerate}
\item
Show that~$(f\otimes g)(\ceil{f}\otimes\ceil{g})=1\otimes 1$,
and conclude that~$\ceil{f\otimes g} \leq\ceil{f}\otimes\ceil{g}$.
\item
Assume for the moment that~$\scrA$
and~$\scrC$ are von Neumann algebras
of bounded operators on Hilbert spaces~$\scrH$
and~$\scrK$, respectively,
and that~$f$ and~$g$ are vector functionals,
that is, $\scrB=\scrD=\C$,
and~$f=\left<x,(\,\cdot\,)x\right>$
for some~$x\in\scrH$,
and
$g=\left<y,(\,\cdot\,)y\right>$
for some
$y\in \scrK$.

Show that~$\ceil{f}=\bigcup_{a\in\scrA^\square}
\ceil{\,a^* \ketbra{x}{x} a\,}$
using~\sref{carrier-vector-state}
and~\sref{double-commutant}.

\item
With the same assumptions as in the previous
point,
suppose, furthermore, without loss of generality
that~$\scrA\otimes \scrB$
is
given as the von Neumann subalgebra of~$\scrB(\scrH\otimes \scrK)$
generated by the operators~$A\otimes B$
where~$A\in\scrA$ and~$B\in\scrB$
(cf.~\sref{special-tensor}).

Show that~$f\otimes g = \left<x\otimes y,(\,\cdot\,)x\otimes y\right>$.

Given~$a\in\scrA^\square$ and~$b\in\scrB^\square$
show that $a\otimes b \in (\scrA\otimes \scrB)^\square$,
and thus
\begin{equation*}
	\ceil{a^* \ketbra{x}{x}a}\otimes
\ceil{b^* \ketbra{y}{y}b} \,\leq\, \ceil{f\otimes g}.
\end{equation*}

Deduce from this that~$\ceil{f}\otimes \ceil{g}\leq \ceil{f\otimes g}$,
so~$\ceil{f}\otimes \ceil{g}= \ceil{f\otimes g}$.
\item
Let~$f$ and~$g$ be arbitrary again,
and assume now that~$f$ and~$g$ are functionals,
that is, $\scrB=\scrD=\C$.
Show that~$\ceil{f\otimes g}=\ceil{f}\otimes \ceil{g}$.
\item
Let~$f$ and~$g$ be arbitrary again,
and recall
from~\sref{ultracyclic-basic}
that~$1=\bigcup_{\sigma}\ceil{\sigma}$
when~$\sigma$ ranges over the np-functionals
$\sigma$ on~$\scrB$.

Show that~$1\otimes 1 
= \bigcup_{\sigma,\tau}\ceil{\sigma \otimes \tau}$
where~$\sigma$ and~$\tau$
range over the np-functionals
on~$\scrB$ and~$\scrD$, respectively.

Show using~\sref{diamond-suprema}
and~\sref{diamond-composition}
that
$\ceil{f\otimes g} \equiv
(f\otimes g)_\diamond(1\otimes 1)
= \ceil{f}\otimes \ceil{g}$.
\item
Show that $(f\otimes g)_\diamond (s\otimes t)
= f_\diamond(s)\otimes g_\diamond(t)$
for projections
$s\in \scrB$
and~$t\in \scrD$.
\end{enumerate}
\spacingfix%
\end{point}%
\end{parsec}%
\subsection{Monoidal Structure}
\begin{parsec}{1190}%
\begin{point}{10}%
Up to this point
we have only written about the tensor product~$\scrA\otimes \scrB$
of \emph{two} von Neumann algebras
(to save ink),
but all of it,
as you will no doubt have observed already,
can be easily adapted 
to deal with
a tensor product
$\otimes\colon \scrA_1\times\dotsc\times \scrA_n
\to\scrA_1\otimes \dotsb\otimes \scrA_n$
of a tuple $\scrA_1,\dotsc,\scrA_n$ of von Neumann algebras,
which will then, of course, be a multilinear map
instead of a bilinear map, etc..

What is less obvious
is that there should be any relation
between 
$(\scrA\otimes\scrB)\otimes \scrC$,
and 
$\scrA\otimes (\scrB\otimes \scrC)$
and~$\scrA\otimes\scrB\otimes\scrC$;
but there is.
\end{point}
\begin{point}{20}{Proposition}%
Given von Neumann algebras~$\scrA$, $\scrB$ and~$\scrC$,
the trilinear map $\gamma\colon (a,b,c)\mapsto (a\otimes b)\otimes c,\ 
\scrA\times\scrB\times\scrC \to (\scrA\otimes \scrB)\otimes\scrC$
is a tensor product.
\begin{point}{30}{Proof}%
We need to verify the three conditions
from~\sref{tensor} (adapted
to trilinear maps).
The first condition,
that the elements of the form~$(a\otimes b)\otimes c$
generate $(\scrA\otimes \scrB)\otimes \scrC$
follows
by~\sref{tensor-generation}
since
the elements of the form~$a\otimes b$
generate~$\scrA\otimes \scrB$
(and~$\scrC$ generates~$\scrC$).
The second condition
is met by defining
$\gamma(\sigma,\tau,\upsilon):= (\sigma\otimes\tau)\otimes\upsilon$
for all np-functionals
$\sigma\colon \scrA\to\C$,
$\tau\colon \scrB \to\C$
and~$\upsilon\colon \scrC\to\C$.
Finally,
these product functionals
$\gamma(\sigma,\tau,\upsilon)$
are centre separating by~\sref{tensor-generation}
because
the functionals on~$\scrA\otimes \scrB$
of the form $\sigma\otimes\tau$
are centre separating (and so is
the set of all np-functionals on~$\scrC$),
which was the third condition.\qed
\end{point}
\end{point}
\begin{point}{40}[associator]{Corollary}%
There is a unique nmiu-isomorphism
\begin{equation*}
    \Define{\alpha_{\scrA,\scrB,\scrC}}\colon 
    \, \scrA\otimes(\scrB\otimes \scrC)
\longrightarrow  (\scrA \otimes \scrB)\otimes\scrC,
\end{equation*}%
    called an \Define{associator}\index{associator}%
\index{alpha@$\alpha_{\scrA,\scrB,\scrC}$, associator},
with $\alpha_{\scrA,\scrB,\scrC}
    (\,a\otimes(b\otimes c)\,)=(a\otimes b)\otimes c$
for all~$a\in\scrA$, $b\in\scrB$, $c\in\scrC$,
for any von Neumann algebras
$\scrA$, $\scrB$, $\scrC$.
\end{point}
\begin{point}{41}%
If the point above means that~$\otimes$
is associative,
then the following two points
mean that~$\otimes$ has~$\C$ as its unit,
and~$\otimes$ is commutative, respectively.
\end{point}
\begin{point}{42}{Exercise}%
Show that given a von Neumann algebra~$\scrA$
the bilinear maps $(z,a)\mapsto za\colon \C\times \scrA\to\scrA$
    and $(a,z)\mapsto za\colon \scrA\times \C\to\scrA$
are tensor products,
and deduce from this that there are unique nmiu-isomorphisms
\begin{equation*}
\Define{\lambda_\scrA}\colon \C\otimes\scrA\longrightarrow\scrA,
\qquad
    \text{and,}\qquad
\Define{\varrho_\scrA}\colon \scrA\otimes\C\longrightarrow\scrA,
\end{equation*}%
    called a \Define{left} and \Define{right unitor}\index{unitor}%
\index{rho@$\varrho_\scrA$, right unitor}%
\index{lambda@$\lambda_\scrA$, left unitor},
respectively,
    with~$\lambda_\scrA(z\otimes a)=za=\varrho_\scrA(a\otimes z)$ for all~$a\in\scrA$
and~$z\in\C$.
\end{point}
\begin{point}{43}{Exercise}
Show that given von Neumann algebras~$\scrA$ and~$\scrB$
the bilinear map $(a,b)\mapsto b\otimes a\colon \scrA \otimes \scrB
\longrightarrow \scrB\otimes \scrA$
is a tensor product,
and deduce from this that there is a unique nmiu-isomorphism
\begin{equation*}
    \Define{\gamma_{\scrA,\scrB}}\colon \scrA\otimes\scrB\longrightarrow
    \scrB\otimes\scrA,
\end{equation*}
called a \Define{braiding}\index{braiding}%
\index{gamma@$\gamma_{\scrA,\scrB}$, braiding},
    with $\gamma_{\scrA,\scrB}(a\otimes b)=b\otimes a$
    for all~$a\in\scrA$ and~$b\in\scrB$.
\end{point}
\begin{point}{50}[vn-smc]{Theorem}%
Endowed with the tensor product,
$\W{miu}$, $\W{cp}$, $\W{cpu}$,
and~$\W{cpsu}$ are symmetric monoidal categories\cite{maclane}
with~$\C$ as unit.
\begin{point}{51}{Proof}%
The first order of business is showing that
the associators $\alpha_{\scrA,\scrB,\scrC}$
(from~\sref{associator})
form a natural transformation
    in~$\W{cp}$
    (and thus in~$\W{miu}$,
    $\W{cpu}$, and $\W{cpsu}$ too,
    as the $\alpha_{\scrA,\scrB,\scrC}$'s are nmiu),
that is, 
that the following diagram commutes
\begin{equation}
    \label{eq:tensor-alpha-natural}
    \xymatrix@C=4em{
\scrA\otimes(\scrB\otimes\scrC)
    \ar[r]^-{f\otimes(g\otimes h)}
    \ar[d]_-{\alpha_{\scrA,\scrB,\scrC}}
&
\scrA'\otimes(\scrB'\otimes\scrC')
    \ar[d]^-{\alpha_{\scrA',\scrB',\scrC'}}
\\
(\scrA\otimes\scrB)\otimes\scrC
    \ar[r]^-{(f\otimes g)\otimes h}
&
(\scrA'\otimes \scrB')\otimes \scrC'
}
\end{equation}
for all ncp-maps $f\colon \scrA\to \scrA'$,
    $g\colon \scrB\to\scrB'$,
    and $h\colon \scrC\to\scrC'$.
Note that by both routes through this diagram
$a\otimes (b\otimes c)$
gets mapped to $(f(a)\otimes g(b))\otimes h(c)$
for all~$a\in\scrA$, $b\in\scrB$, and~$c\in \scrC$.
Since the linear span of such
$a\otimes (b\otimes c)$'s
is ultraweakly dense in~$\scrA\otimes(\scrB\otimes \scrC)$
    (by~\sref{tensor-generation},)
this entails that~\eqref{eq:tensor-alpha-natural} commutes.

By a similar but simpler argument one sees
that the braidings ($\gamma_{\scrA,\scrB}$) and unitors
($\lambda_\scrA$ and $\varrho_\scrA$)
give natural transformations.
\begin{point}{52}%
It remains to be shown
that the appropriate coherence relations hold.
Given von Neumann algebras $\scrA$, $\scrB$,
$\scrC$, and~$\scrD$,
the pentagon
    \newcommand{\wrapcell}[1]{\hbox to 2em{\hss $#1$\hss}}
\begin{equation*}
\xymatrix@C=.1em@R=3em{
    && \wrapcell{(\scrA\otimes \scrB)\otimes (\scrC\otimes \scrD) } 
    \ar[rrd]^{\alpha_{\scrA\otimes\scrB,\scrC,\scrD}}
    && \\
    \wrapcell{\scrA\otimes(\scrB\otimes (\scrC \otimes \scrD))}
    \ar[rru]^{\alpha_{\scrA,\scrB,\scrC\otimes\scrD}}
    \ar[ddr]|{\id_\scrA\otimes\alpha_{\scrB,\scrC,\scrD}}
    &&&&
    \wrapcell{((\scrA\otimes \scrB)\otimes \scrC)\otimes \scrD}
    \\
    \\
    & 
    {\scrA\otimes((\scrB\otimes \scrC)\otimes\scrD)}
    \ar[rr]_{\alpha_{\scrA,\scrB\otimes\scrC,\scrD}}
    &&
    {(\scrA\otimes (\scrB\otimes \scrC))\otimes\scrD}
    \ar[ruu]|{\alpha_{\scrA,\scrB,\scrC}\otimes \id_\scrD}
}
\end{equation*}
commutes,
since by both routes
the elements in~$\scrA\otimes(\scrB\otimes (\scrC\otimes \scrD))$
of the form $a\otimes (b\otimes (c\otimes d))$
(whose linear span is ultraweakly dense)
get sent to $((a\otimes b)\otimes c)\otimes d$.

By similar arguments the diagrams
\begin{equation*}
\xymatrix{
    \scrA\otimes(\C\otimes \scrC)
    \ar[rr]^{\alpha_{\scrA,\C,\scrC}}
    \ar[rd]_{\id_\scrA\otimes \lambda_\scrC}
    &
    &
    (\scrA\otimes \C)\otimes \scrC
    &
    \C\otimes \C 
    \ar@/_1em/[d]_{\varrho_\C}
    \ar@/^1em/[d]^{\lambda_\C}
    \\
    &
    \scrA\otimes \scrC
    \ar[ru]_{\varrho_\scrA\otimes \id_\scrC}
    &
    &
    \C
}
\end{equation*}
commute, as does the diagram
\begin{equation*}
\xymatrix@C=4em{
    \scrA\otimes(\scrB\otimes \scrC)
    \ar[r]^{\alpha_{\scrA,\scrB,\scrC}}
    \ar[d]_{\id_\scrA\otimes\gamma_{\scrB,\scrC}}
    &
    (\scrA\otimes\scrB)\otimes\scrC
    \ar[r]^{\gamma_{\scrA\otimes\scrB,\scrC}}
    &
    \scrC\otimes(\scrA\otimes\scrB)
    \ar[d]^{\alpha_{\scrC,\scrA,\scrB}}
    \\
    \scrA\otimes(\scrC\otimes \scrB)
    \ar[r]_{\alpha_{\scrA,\scrC,\scrB}}
    &
    (\scrA\otimes \scrC)\otimes \scrB
    \ar[r]_{\gamma_{\scrA,\scrC}\otimes \id_\scrB}
    &
    (\scrC\otimes \scrA)\otimes \scrB
},
\end{equation*}
and do the diagrams
\begin{equation*}
\xymatrix{
\scrA\otimes\scrB
    \ar[r]^{\gamma_{\scrA,\scrB}}
    \ar[rd]_{\id_{\scrA\otimes \scrB}}
&
\scrB\otimes \scrA
    \ar[d]^{\gamma_{\scrB,\scrA}}
&
    &
\scrB\otimes\C
    \ar[r]^{\gamma_{\scrB,\C}}
    \ar[rd]_{\varrho_{\scrB}}
&
\C\otimes \scrB
    \ar[d]^{\lambda_\scrB}
\\
&
\scrA\otimes \scrB
&
    &
&
\scrB
}.
\end{equation*}
Thus $\W{miu}$, $\W{cp}$, $\W{cpu}$,
and~$\W{cpsu}$ are symmetric monoidal categories.\qed
\end{point}
\end{point}
\end{point}
\end{parsec}
\section{Quantum Lambda Calculus}
\label{S:model}
\begin{parsec}{1200}%
\begin{point}{10}[qlcm-intro]%
In this section we provide
the parts needed
to build a model
	of the quantum lambda calculus\index{quantum lambda calculus}
using von Neumann algebras.
We will not venture
to describe the quantum lambda calculus
in all its details here,
nor will we describe how to build the model
from these parts
(as we did in~\cite{model});
we'll just touch upon
the two key ingredients:
the interpretation of ``$\bang$'' and ``$\limp$''%
 \index{*bang@$\bang$}\index{*lollypop@$\limp$}
---
with them the expert
can easily produce the model.

Let us, nevertheless, try to give some impression to those who are not familiar 
with the quantum lambda calculus.
The quantum lambda calculus is a type theory
proposed by Selinger and Valiron in~\cite{selinger2005,selinger2006}
to describe programs for quantum computers
especially designed
to include 
not only
function types ($\limp$)
and
classical data types (such as~$\bit$),
but also quantum data types
(such as $\qbit$),
so that there can be a term such as
$\mathsf{new}\colon \bit\limp \qbit$
that represents the program
that initialises a qubit in the given state.
There are of course also terms
such as
$\mathsf{0}\colon \bit$
and~$\mathsf{1}\colon \bit$,
so that~$\mathsf{new}\,\mathsf{0}\colon \qbit$
represents a qubit in state~$\ket{0}$.
The addition of quantum data 
to a type theory
is a very delicate matter
for if one were to allow
for example
in this system
a variable to be used twice
(a thing usually beyond dispute)
it would not take much more
to construct
a program
that duplicates the contents
of a qubit,
which is nonphysical.

Still, classical data
such as a bit
can be duplicated freely,
so to accommodate this 
the type~$\bang \bit$ is used.
More precisely,
the type
$\bang A$
represents that part of the type of~$A$ that is duplicable,
so that~$\bang \bit$
is the proper type for a bit,
and $\bang\qbit$ is empty.
For example,
the term that represents the measurement
of a qubit is~$\mathsf{meas}\colon \qbit\limp \bang\bit$,
where the~$\bang$ indicates that the bit resulting from the measurement
may be duplicated freely.

The model
we alluded to
assigns to each type~$A$ a von Neumann algebra~$\sem{A}$,%
\index{*sem@$\sem{\,\cdot\,}$}
e.g.~$\sem{\qbit}=M_2$
and~$\sem{\bit}=\C^2$.
A (closed) term~$t:A$
is interpreted as an npsu-functional $\sem{t:A}\colon \sem{A}\to\C$,
so for example $\sem{\mathsf{0}\colon \bit}\colon
(x,y)\mapsto x\colon \C^2\to\C$.
When~$t:A$ has free variables $x_1:B_1,\dotsc,x_N:B_N$
the interpretation 
becomes an ncpsu-map $\sem{t}\colon \sem{A}\to\sem{B_1}\otimes 
\dotsb \otimes \sem{B_N}$,
so for example,
\begin{equation*}
	\sem{\, x\colon \qbit \vdash \mathsf{meas} x\,}
\colon (x,y)\mapsto
\smash{\bigl(
	\begin{smallmatrix}x& 0 \\ 0 & y
	\end{smallmatrix}\bigr)}
\colon \C^2\to M_2.
\end{equation*}
In short, there are no surprises here.
As said, the difficulty lies
in the definition
of $\sem{!A}$ and~$\sem{A\limp B}$,
for which we will provide the following
three ingredients.
\begin{itemize}
\item
The observation (by Kornell, \cite{kornell2012})
that
the category 
$\op{(\W{miu})}$
is monoidal closed,
that is,
that
for every von Neumann algebra~$\scrB$,
the functor
$\scrB\otimes (\,\cdot\,)\colon 
\W{miu}\to \W{miu}$
has a left adjoint
$(\,\cdot\,)^{*\scrB}$.

\item
The following two adjunctions.
\begin{equation*}
\xymatrix@C=6em{
	\Cat{Set}
	\ar@/^1.0em/[r]^-{\ell^\infty}
	\ar@{}[r]|-{\perp}
	&
	\op{(\W{miu})}
	\ar@/^1.0em/[r]^-{\subseteq}
	\ar@/^1.0em/[l]^-{\nsp:=\W{miu}(-,\C)}
	\ar@{}[r]|-{\perp}
	&
	\op{(\W{cpsu})}
	\ar@/^1.0em/[l]^-{\mathcal{F}}
}
\end{equation*}
\end{itemize}
The interpretation of~$\sem{!A}$
and~$\sem{A\limp B}$ will then be
\begin{equation*}
\sem{!A}\ =\ 
\linf(\nsp(\sem{A}))
\qquad\text{and}\qquad
\sem{A\limp B}
\ = \ 
\mathcal{F}(\sem{B})^{*\sem{A}}.
\end{equation*}%
By the universal properties of~$\mathcal{F}$
and~$(\,\cdot\,)^{*\scrB}$,
an ncpsu-map $f\colon \scrA\rightarrow \scrC\otimes \scrB$
corresponds to unique nmiu-map 
$\Lambda(f)\colon \mathcal{F}(\scrA)^{*\scrB}\longrightarrow \scrC$.
This is used to interpret 
the ``$\lambda$'':
\begin{equation*}
    \sem{\lambda x^B. t}
    \,=\, \Lambda(\sem{t})\colon \mathcal{F}(\sem{A})^{\sem{B}}
    \longrightarrow\sem{B_1}\otimes\dotsb\otimes\sem{B_N}
\end{equation*}
for any term  $t:A$  with free variables
$x_1:B_1$, \ldots, $x_N:B_N$, $x:B$,
so  
\begin{equation*}
\sem{t} \colon \sem{A}\longrightarrow 
\sem{B_1}\otimes\dotsb\otimes\sem{B_N}\otimes \sem{B}.
\end{equation*}

Note that~$\sem{!A}$
will always be a `discrete' commutative von Neumann algebra
no matter how complicated~$\sem{A}$ may be,
so that although this does the job
perhaps a more interesting interpretation of~$\bang$
may be chosen as well.
This is not the case:
in the next
section we'll show that
any von Neumann algebra
that carries 
a $\otimes$-monoid structure
(such as $\sem{!A}$)
is commutative and discrete,
and that~$\linf(\nsp(\scrA))$
is moreover the free $\otimes$-monoid
on~$\scrA$.
\end{point}
\begin{point}{20}%
In~\cite{qlc} the quantum lambda calculus
is extended with recursion via the ``$let\ rec$''
operator;
we don't know whether it's possible to interpret
    $let\ rec$ in our model.
\end{point}
\end{parsec}
\begin{parsec}{1210}%
\begin{point}{10}%
In this section,
we'll need the following result
from the literature on von Neumann algebras.
\end{point}
\begin{point}{20}[intersection-tensor]{Proposition}%
Given Hilbert spaces~$\scrH$ and~$\scrK$,
and von Neumann subalgebras~$\scrA_1$ and~$\scrA_2$
of~$\scrB(\scrH)$
and von Neumann subalgebras~$\scrB_1$ and~$\scrB_2$
of~$\scrB(\scrK)$,
we have 
\begin{equation*}
(\scrA_1 \otimes  \scrB_1)\,\cap\,
(\scrA_2 \otimes \scrB_2)
= (\scrA_1\cap\scrA_2)\,\otimes\,
(\scrB_1\cap \scrB_2).
\end{equation*}
Here
$\scrA_1 \otimes \scrB_1$
denotes not just any tensor product of~$\scrA_1$ and~$\scrB_1$,
but instead
the ``concrete'' tensor product 
of~$\scrA_1$ and~$\scrB_1$:
the least von Neumann subalgebra of~$\scrB(\scrH\otimes\scrK)$
that contains all operators
of the form $A\otimes B$ where~$A\in\scrA_1$
and~$B\in\scrB_1$.
\begin{point}{30}[intersection-tensor-proof]{Proof}%
See Corollary IV.5.10 of~\cite{Takesaki1}.\qed
\end{point}
\end{point}
\end{parsec}
\subsection{First Adjunction}
\begin{parsec}{1220}%
\begin{point}{10}{Definition}%
We write $\Define{\nsp}:=\W{miu}(\,\cdot\,,\C)$
	\index{nsp@$\nsp\colon \op{(\W{miu})}\to \Cat{Set}$}
for the functor $\op{(\W{miu})}\to \Cat{Set}$
which maps a von Neumann algebra~$\scrA$
to its set
of nmiu-functionals,
$\nsp(\scrA)$,
and sends an nmiu-map $f\colon \scrA\to\scrB$
to the map~$\nsp(f)\colon \nsp(\scrB)\to\nsp(\scrA)$
given by~$\nsp(f)(\varphi)=\varphi \circ f$
for~$\varphi \in \nsp(\scrB)$.
\end{point}
\begin{point}{20}[first-adjunction]{Proposition}%
Given a set~$X$
the map
\begin{equation*}
\eta\colon X\to\nsp(\linf(X))
\qquad\text{given by}\qquad
	\eta(x)(h)\ =\ h(x)
\end{equation*}
is universal
in the sense that for every 
map $f\colon X\to\nsp(\scrA)$,
where~$\scrA$ is a von Neumann algebra,
there is a unique
nmiu-map $g\colon \scrA\to\linf(X)$
such that
\begin{equation*}%
\xymatrix{%
X
\ar[r]^-\eta
\ar[rd]_-f
&
\nsp(\linf(X))
\ar[d]^-{\nsp(g)}
&
\linf(X)
\\
&
\nsp(\scrA)
&
\scrA
\ar@{.>}[u]_-g
}
\end{equation*}
commutes.
Moreover, and as a result, the assignment
$X\mapsto \linf(X)$
extends to a functor
$\linf\colon \Cat{Set}\to\op{(\W{miu})}$%
	\index{linfunctor@$\linf\colon \Cat{Set}\to\op{(\W{miu})}$}
that is left adjoint
to~$\nsp$,
and is given by $\linf(f)(h) = h\circ f$
for any map $f\colon X\to Y$
and~$h\in \linf(Y)$.
\begin{point}{30}{Proof}%
Note that if we identify~$\linf(X)$
with the $X$-fold product of~$\C$,
we see that~$\eta(x)\colon \linf(X)\equiv \bigoplus_{x\in X}\C\to \C$
is simply the $x$-th projection,
and thus an nmiu-map (see~\sref{vn-products}).
Hence we do indeed get a map $\eta\colon X\to\nsp(\linf(X))$.

To see that~$\eta$ has the desired universal property,
let~$f\colon X\to\nsp(\scrA)$ be given,
and define~$g\colon \scrA\to \linf(X)$
by~$g(a)(x)=f(x)(a)$.
One can now either prove directly that~$g$ is nmiu,
or reduce this in a slightly roundabout way
from the known fact that~$\linf(X)$ is the 
$X$-fold product of~$\C$ with the $\eta(x)$ as projections;
indeed~$g$ is simply the  unique nmiu-map with $\eta(x)\circ g = f(x)$
for all~$x\in X$, that is,  $g=\left<f(x)\right>_{x\in X}$.
In any case,
we see that $\nsp(g)(\eta(x)) \equiv \eta(x)\circ g = f(x)$
for all~$x\in X$, and so~$\nsp(g) \circ \eta = f$.
Concerning uniqueness of such~$g$,
note that given an nmiu-map $g'\colon \scrA\to\linf(X)$
with~$\nsp(g')\circ \eta = f$
we have~$\eta(x)\circ  g'= \nsp(g')(\eta(x))
= f(x)$ for all~$x\in X$,
and so~$g'=\left<f(x)\right>_{x\in X}=g$.

Hence~$\eta$ is a universal arrow from~$X$
to~$\nsp$.
That as a result the assignment $X\mapsto \linf(X)$
extends to a functor
$\Cat{Set}\to\op{(\W{miu})}$
by sending~$f\colon X\to Y$
to the unique nmiu-map $\linf(f)\colon \linf(Y)\to\linf(X)$
with $\nsp(\linf(f))\circ \eta_X = \eta_Y\circ f$
is a known and easily checked fact
(where~$\eta_X:= \eta$ and~$\eta_Y\colon Y\to\nsp(\linf(Y))$ 
is what you'd expect).
Finally, applying~$x\in X$ and~$h\in \linf(Y)$
we get
$ \linf(f)(h)(x)
= \eta_X(x)(\linf(f)(h)))
=\nsp(\linf(x))(\eta_X(x))(h)
=\eta_Y(f(x))(h)
=h(f(x))$.\qed
\end{point}
\end{point}
\begin{point}{40}[nmiu-functional-product]{Lemma}%
A nmiu-functional $\varphi$
on a direct sum $\bigoplus_i \scrA_i$
of von Neumann algebras
is of the form~$\varphi\equiv \varphi'\circ \pi_i$
for some~$i$ and nmiu-functional~$\varphi'$ on~$\scrA_i$.
\begin{point}{50}{Proof}%
Let~$e_j$ denote the element of~$\bigoplus_i \scrA_i$
given by $e_j(j)=1$ and~$e_j(i)=0$ for all~$i\neq j$.
Note that given~$i$ and~$j$ with $i\neq j$
we have~$e_ie_j=0$
and so~$0=\varphi(e_i e_j)=\varphi(e_i)\varphi(e_j)$;
from this we see that
there is at most one~$i$ with~$\varphi(e_i)\neq 0$.
Since for this~$i$
we have $e_i^\perp = \sum_{j\neq i} e_j$
and so~$\varphi(e_i^\perp)=\sum_{j\neq i} \varphi(e_j)=0$,
we see that $\varphi(a)=\varphi(e_i a )$
for all~$a\in \bigoplus_i \scrA_i$.
Letting  $\kappa_i\colon \scrA_i\to\bigoplus_j \scrA_j$
be the nmisu-map given by~$\kappa_i(a)(i)=a$
and~$\kappa_i(a)(j)=0$ for~$j\neq i$
we have~$\varphi = \varphi \circ \kappa_i \circ \pi_i$.
Hence taking~$\varphi':=\varphi\circ \kappa_i$
does the job.\qed
\end{point}
\end{point}
\begin{point}{60}[cor:linf-ff]{Exercise}%
Deduce from~\sref{nmiu-functional-product}
that the functor $\nsp\colon \op{(\W{miu})}\to \Cat{Set}$
preserves coproducts,
and that the 
map $\eta\colon X\to \nsp(\linf(X))$
from~\sref{first-adjunction} is a bijection.

Show that~$\linf\colon\Cat{Set}\to \op{(\W{miu})}$
is full and faithful.
	Whence~$\Cat{Set}$ is (isomorphic to)  
a coreflective subcategory of~$\op{(\W{miu})}$ 
	via~$\linf\colon \Cat{Set} \to \op{(\W{miu})}$.
\end{point}
\end{parsec}
\begin{parsec}{1230}
\begin{point}{10}{Exercise}%
	We're going to prove that
	$\ell^\infty(X\times Y)\cong \ell^\infty(X)\otimes \ell^\infty(Y)$.
\begin{enumerate}%
\item
Given an element~$x$ of a set~$X$
let~$\hat{x}$ denote
the element of~$\ell^\infty(X)$
that equals~$1$ on~$x$ and is zero elsewhere.

Show that~$\{\,\hat{x}\colon\, x\in X\,\}$
generates~$\ell^\infty(X)$.
\item
Show that the  projections $\pi_x\colon \ell^\infty(X)\equiv
\bigoplus_{y\in X}\C\to\C$
form an order separating collection of
		nmiu-functionals on~$\ell^\infty(X)$.
\item
Using this, and~\sref{tensor-characterization},
prove that given sets~$X$ and~$Y$
the map
\begin{equation*}
	\otimes\colon \ell^\infty(X)\times \ell^\infty(Y)
\to\ell^\infty(X\times Y)
\end{equation*}
given by~$(f\otimes g)(x,y)=f(x)g(y)$
is a tensor product.

Conclude that 
	$\ell^\infty(X\times Y)\cong \ell^\infty(X)\otimes \ell^\infty(Y)$.

		(In fact, it follows that~$\ell^\infty$ is strong monoidal.)
\end{enumerate}
\spacingfix%
\end{point}%
\begin{point}{20}{Exercise}%
Let~$\scrA$ and~$\scrB$
be von Neumann algebras.
We're going to show that~$\nsp(\scrA\otimes\scrB)\cong \nsp(\scrA)\times
	\nsp(\scrB)$.
\begin{enumerate}
\item
Given an nmiu-functional $\varphi\colon \scrA\otimes \scrB\to\C$
show that~$\sigma:= \varphi((\,\cdot\,)\otimes 1)$
and~$\tau:=\varphi(1\otimes (\,\cdot\,))$
are nmiu-functionals on~$\scrA$ and~$\scrB$, respectively;
		and show that~$\varphi=\sigma\otimes \tau$
		(by proving that $\varphi(a\otimes b)=\sigma(a)\tau(b)$.)
\item
Show that~$\sigma,\tau\mapsto \sigma\otimes \tau$
		gives a bijection $\nsp(\scrA)\times\nsp(\scrB)
		\to\nsp(\scrA\otimes\scrB)$.

		(This  makes~$\nsp$  strong monoidal.)
\end{enumerate}
\spacingfix%
\end{point}%
\end{parsec}%
\subsection{Second Adjunction}
\begin{parsec}{1240}%
\begin{point}{10}[vn-generation-bound]{Lemma}%
If a von Neumann algebra~$\scrA$
is generated by~$S\subseteq \scrA$,
then
\begin{equation*}
\#\scrA \ \leq\  2^{2^{\#\mathbb{C}+\#S}},
\end{equation*}
where~$\#S$ denotes the cardinality of~$S$, and so on.
\begin{point}{20}{Proof}%
Note that the
$*$-subalgebra~$S'$ of~$\scrA$
generated by~$S$
is ultraweakly dense in~$\scrA$.
Since every element of~$S'$
can be formed
from the infinite set $S\cup \mathbb{C}$ using 
the finitary operations
of
addition, multiplication,
and
involution,
$\#S'\leq \#\mathbb{C}+ \#S$.
Since every element of~$\scrA$
is the ultraweak limit of a filter
(see \cite[\S12]{willard})
on~$S'$
of which there no more than~$2^{2^{\#S'}}$,
we conclude~$\#\scrA \leq 2^{2^{\#\mathbb{C}+\#S}}$.\qed
\end{point}
\end{point}
\begin{point}{30}[second-adjunction]{Theorem}%
The inclusion $\W{miu}\to\W{cpsu}$
has a left adjoint $\Define{\mathcal{F}}\colon \W{cpsu}\to\W{miu}$.%
	\index{F@$\mathcal{F}\colon \W{cpsu}\to\W{miu}$}
\begin{point}{40}[second-adjunction-proof]{Proof}%
Note that since the category~$\W{miu}$
has all products (\sref{vn-products}),
and equalisers (\sref{vn-equalisers}),
$\W{miu}$ has all limits 
(by Theorem~V2.1 and Exercise V4.2 of~\cite{maclane}).
Moreover, the inclusion~$U\colon \W{miu}\to\W{cpsu}$
preserves these limits (see~\sref{vn-products} and~\sref{vn-equalisers}).
So by Freyd's adjoint functor theorem (Theorem~V6.1 of~\cite{maclane})
it suffices to check the \emph{solution set condition},
that is,
that 
\begin{quote}
for every von Neumann algebra~$\scrA$
there 
be a set~$I$,
and for each~$i\in I$ an ncpsu-map $f_i\colon \scrA\to\scrA_i$
into a von Neumann algebra~$\scrA_i$
such that every ncpsu-map $f\colon \scrA\to \scrB$
into some von Neumann algebra~$\scrB$
is of the form~$ f\equiv h\circ f_i$
for some~$i\in I$ and nmiu-map $h\colon \scrA_i\to\scrB$.
\end{quote}
To this end, given a von Neumann algebra~$\scrA$,
let~$\kappa:=2^{2^{\#\C+\#\scrA}}$, define
\begin{alignat*}{3}
	I\ = \ \{\ (\scrC,\gamma)\colon \scrC\text{ is }&
	\text{a von Neumann algebra
	on a subset of~$\kappa$,}\\
&\text{and $\gamma \colon \scrA\to\scrC$ is an ncpsu-map} \ \},
\end{alignat*}
and set~$f_i := \gamma$ for every $i\equiv(\scrC,\gamma)\in I$.

Let $f\colon \scrA\to\scrB$
be an ncpsu-map
into a von Neumann algebra~$\scrB$.
The von Neumann algebra~$\scrB'$ generated by~$f(\scrA)$
has cardinality below~$\kappa$
by~\sref{vn-generation-bound},
and so by relabelling the elements of~$\scrB'$
we may find a von Neumann algebra~$\scrC$
on a subset of~$\kappa$
isomorphic to~$\scrB'$
via some nmiu-isomorphism $\Phi\colon \scrB'\to\scrC$.
Then the map $\gamma\colon \scrA\to \scrC$
given by~$\gamma(a)=\Phi(f(a))$ for all~$a\in\scrA$
is ncpsu,
so that~$i:=(\scrC,\gamma)\in I$,
and, moreover,
the assignment~$c\mapsto \Phi^{-1}(c)$
gives an nmiu-map $h\colon \scrC\to \scrB$
with~$h\circ f_i \equiv h\circ \gamma=f$.
Hence~$U\colon \W{miu}\to\W{cpsu}$
obeys the solution set condition,
and therefore has a left adjoint.\qed
\end{point}
\begin{point}{50}{Remark}%
A bit more can be said
about the adjunction between the inclusion~$U\colon \W{miu}\to\W{cpsu}$
and~$\mathcal{F}$:
since~$\W{miu}$ has the same objects as~$\W{cpsu}$,
the category~$\op{(\W{cpsu})}$
is, for very general reasons, equivalent
to the Kleisli category
	of the (by the adjunction induced) monad~$\mathcal{F}U$ 
on~$\op{(\W{miu})}$
in a certain natural way
(see e.g.~Theorem~9 of~\cite{qpakm}).
\end{point}
\end{point}
\end{parsec}
\subsection{Free Exponential}
\begin{parsec}{1250}%
\begin{point}{10}%
	We'll prove Kornell's result (from~\cite{kornell2012})
that the functor
$\scrB\otimes(\,\cdot\,)\colon
\W{miu}\to\W{miu}$
has a left adjoint~$(\,\cdot\,)^{*\scrB}$
for every von Neumann algebra~$\scrB$.
Kornell original proof is rather complex,
and so is ours, unfortunately,
but we've managed
to peel off one layer of complexity from the original 
proof
by way of Freyd's Adjoint Functor Theorem,
reducing the problem
to the facts that
$\scrB\otimes(\,\cdot\,)\colon \W{miu}\to\W{miu}$
preserves products,  equalisers,
and satisfies the solution set condition.
\end{point}
\begin{point}{20}[vn-gns-bound]{Lemma}%
A von Neumann algebra~$\scrA$
can be faithfully represented
on a Hilbert space which contains no more
than~$2^{\#\scrA}$ vectors.
\begin{point}{30}{Proof}%
If~$\scrA=\{0\}$,
then the result is obvious,
so let us assume that~$\scrA\neq \{0\}$.
Then~$\scrA$ is infinite,
and so~$\aleph_0 \cdot \#\scrA  = \#\scrA$. 

Let~$\Omega$ be the set of np-functionals on~$\scrA$.
Recall that 
by the GNS-construction (see~\sref{ngns})
$\mathscr{A}$
can be faithfully represented on
the Hilbert space
$\scrH_\Omega\equiv \bigoplus_{\omega\in\Omega} \scrH_\omega$.
Since every element of~$\scrH_\omega$
is the limit of a sequence of elements from~$\scrA$,
we have $\#\scrH_\omega \leq \aleph_0^{\#\scrA} \leq (2^{\aleph_0})^{\#\scrA} 
= 2^{\#\scrA}$,
because $\aleph_0\cdot \#\scrA=\#\scrA$.
Since every normal state is a map $\omega\colon \scrA\to\C$,
we have $\#\Omega\leq \#\C^{\#\scrA}=(2^{\aleph_0})^{\#\scrA}
= 2^{\#\scrA}$, because $\aleph_0 \cdot \#\scrA = \#\scrA$.
Hence $\#\scrH = \sum_{\omega\in\Omega} \#\scrH_\omega
\leq 2^{\#\scrA}\cdot2^{\#\scrA}
=2^{\#\scrA}$.\qed
\end{point}
\end{point}
\begin{point}{40}[equaliser-lemma]{Lemma (Kornell)}%
Every nmiu-map $h\colon \scrD\to\scrA\otimes \scrC$,
where~$\scrA$, $\scrC$ and~$\scrD$ are von Neumann algebras,
factors as 
$\smash{\xymatrix@C=3em{\scrD
\ar[r]|-{\tilde{h}}
& 
\tilde\scrA\otimes \scrC
\ar[r]|-{\iota\otimes\id}
&
\scrA \otimes \scrC
}}$,
where~$\smash{\tilde\scrA}$
is a von Neumann algebra,
and~$\iota$ and~$\tilde{h}$
are nmiu-maps,
such that
for all nmiu-maps $f,g\colon \scrA\to\scrB$
into some von Neumann algebra~$\scrB$
with $(f\otimes\id)\circ h = (g\otimes \id)\circ h$
we have $f\circ \iota = g\circ \iota$.

Moreover,
$\tilde\scrA$
can be generated by
less than~$\#\scrD\,\cdot\,2^{\#\scrC}$ elements.
\begin{point}{50}{Proof}%
Assume (without loss of generality)
that
$\scrC$
is a von Neumann algebra of operators on
a Hilbert space~$\scrH$
with no more than $2^{\#\scrC}$ vectors, see~\sref{vn-gns-bound}.

For every vector~$\xi\in \scrH$
let
$r_\xi\colon \scrA\otimes \scrC\to\scrA$
be the unique np-map
given by~$r_\xi(a\otimes c) = \left<\xi,c\xi\right> a$
for all~$a\in\scrA$ and~$c\in \scrC$
(see~\sref{tensor-universal-property}
and~\sref{tensor-universal-property-extra}),
and
let~$\tilde\scrA$
be the least von Neumann subalgebra
of~$\scrA$
that contains $S:=\smash{ \bigcup_{\xi \in\scrH} r_\xi(h(\scrD))}$,
and let~$\iota\colon \tilde\scrA\to\scrA$
be the inclusion
(so~$\iota$ is nmiu).
Note that~$S$ (which generates~$\tilde\scrA$)
has no more than $\#\scrD\cdot\#\scrH\leq \#\scrD \cdot 2^{\#\scrC}$
elements.

Let~$f,g\colon \scrA\to\scrB$
be nmiu-maps
into a von Neumann algebra~$\scrB$
such that $(f\otimes\id)\circ h = (g\otimes \id)\circ h$.
We must show that~$f\circ \iota = g\circ \iota$.
By definition of~$\tilde\scrA$
(and the fact that~$f$ and~$g$ are nmiu),
it suffices to show that $f\circ r_\xi\circ h=g\circ r_\xi\circ h$
for all~$\xi\in\scrH$.
Note that given such~$\xi$,
we have $f\circ r_\xi = r_\xi' \circ (f\otimes\id)$,
where~$r_\xi'\colon \scrB\otimes\scrC\to\scrB$
is the np-map
given by~$r_\xi'(b\otimes c)=\left<\xi,c\xi\right>b$.
Since similarly,
$g\circ r_\xi = r_\xi'\circ (g\otimes \id)$,
we get~$f\circ r_\xi \circ h
= r_\xi' \circ (f\otimes \id)\circ h
= r_\xi' \circ (g\otimes \id)\circ h
= g\circ r_\xi \circ h$.

It remains only to 
be shown that $h(\scrD) \subseteq \tilde\scrA\otimes \scrC$,
because we may then simply let~$\tilde{h}$
be the restriction of~$h$ to~$\tilde\scrA\otimes \scrC$.
It is enough to prove that
$h(\scrD) \subseteq \tilde\scrA 
\otimes  \bsp(\scrH)$,
because 
$\tilde\scrA \otimes  \scrC
\,=\, (\tilde\scrA\otimes \bsp(\scrH))\,\cap\,
(\scrA\otimes \scrC)$
(see~\sref{intersection-tensor})
and we already know that $h(\scrD)\subseteq \scrA\otimes \scrC$.
Let $(e_k)_k$ be orthonormal basis of $\scrH$.
Since $1= \sum_k \ket{e_k}\!\bra{e_k}$
in $\bsp(\scrH)$,
we have, for all~$d\in\scrD$,
\begin{align*}
h(d) \ &=\ \textstyle\bigl(\sum_k 1\otimes \ket{e_k}\!\bra{e_k}\bigr)
\ h(d)\ \bigl(\sum_\ell 1\otimes \ket{e_\ell}\!\bra{e_\ell}\bigr)\\
\ &=\ \textstyle \sum_k\sum_\ell\ 
 (\,1\otimes \ket{e_k}\!\bra{e_k}\,) \ h(d)\  
 (\,1\otimes \ket{e_\ell}\!\bra{e_\ell}\,).
\end{align*}
We are done
if we can prove that,
for all~$\xi,\zeta\in\scrH$,
\begin{equation}
\label{eq:main-equaliser-todo}
(\,1\otimes \ket{\xi}\!\bra{\xi}\,) \ 
h(d)\  (\,1\otimes \ket{\zeta}\!\bra{\zeta}\,)
\ \in\ \tilde\scrA\otimes \bsp(\scrH).
\end{equation}
By an easy computation, we see that,
for all $e \,\in\,\scrA\otimes \scrC$
	of the form $e\equiv a\otimes c$,
\begin{equation*}
\label{eq:polarisation-equaliser}
(\,1\otimes \ket{\xi}\!\bra{\xi}\,) \ e
\  (\,1\otimes \ket{\zeta}\!\bra{\zeta}\,)
\ =\ 
\frac{1}{4}\sum_{k=0}^3i^k \,r_{i^k\xi+\zeta}(e)\otimes \ket\xi\!\bra\zeta.
\end{equation*}
It follows that the equation above holds for all~$e\in \scrA\otimes \scrC$.
Choosing $e=h(d)$
we see that~\eqref{eq:main-equaliser-todo}
holds,
because
$r_{i^k\xi+\zeta}(h(d)) \in\tilde\scrA$.\qed
\end{point}
\end{point}
\begin{point}{60}[tensor-equalisers]{Proposition}%
Let $e\colon\scrE\to\scrA$ be an equaliser of
nmiu-maps $f,g\colon\scrA\to\scrB$ 
between von Neumann algebras.
Then $e\otimes\id\colon\scrE\otimes\scrC\to\scrA\otimes\scrC$
is an equaliser of 
$f\otimes\id$ and $g\otimes\id$
for every von Neumann algebra~$\scrC$.
\begin{point}{70}{Proof}%
Let~$h\colon\scrD\to\scrA\otimes \scrC$
be an nmiu-map 
with $(f\otimes\id)\circ h=(g\otimes\id)\circ h$.
We must show that there is a unique
 nmiu-map $k\colon\scrD\to\scrE\otimes\scrC$
such that  $h=(e\otimes\id)\circ k$.
Note that since the equaliser map~$e$ 
is injective,
$e\otimes\id\colon\scrE\otimes\scrC\to\scrA\otimes\scrC$
is injective (by~\sref{tensor-injective})
and thus uniqueness of~$k$ is clear.
Concerning existence,
by~\sref{equaliser-lemma},
$h$ factors as
$\smash{\xymatrix@C=3em{\scrD
\ar[r]|-{\tilde{h}}
& 
\tilde\scrA\otimes \scrC
\ar[r]|-{\iota\otimes\id}
&
\scrA \otimes \scrC
}}$
where~$\tilde{h}$ and~$\iota$ are nmiu-maps,
and moreover,
we have $f\circ \iota = g\circ \iota$.
Since~$e$ is an equaliser of~$f$ and~$g$,
there is a unique nmiu-map $\tilde\iota\colon \tilde\scrA\to \scrE$
with~$e\circ \tilde\iota = \iota$.
Now, define $k:=(\tilde \iota\otimes \id)\circ \tilde{h}\colon 
\scrD\to \scrE\otimes \scrC$.
Then~$(e\otimes \id)\circ k = 
((e\circ\tilde \iota)\otimes \id)\circ \tilde{h}
= (\iota\otimes \id)\circ \tilde{h}
= h$.\qed
\end{point}
\end{point}
\begin{point}{71}%
So given a von Neumann algebra~$\scrA$
the functor
$(\,\cdot\,)\otimes \scrA\colon \W{miu}\to\W{miu}$
preserves all equalisers and products,
thus all limits, and in particular, all pullbacks.
This has the following pleasant consequence  used later on.
\end{point}
\begin{point}{72}[tensor-preimage]{Exercise}%
Given a nmiu-map~$\varrho\colon \scrB\to\scrC$ 
between von Neumann algebras~$\scrB$ and~$\scrC$,
and a von Neumann subalgebra $\scrS$ of~$\scrC$,
show that
\begin{equation*}
(\varrho\otimes\scrA)^{-1}(\scrS\otimes \scrA)\ 
    =\ \varrho^{-1}(\scrS)\otimes \scrA
\end{equation*}
for every von Neumann algebra~$\scrA$,
where for the sake of simplicity we
take $\varrho^{-1}(\scrS)\otimes \scrA$
to be the von Neumann subalgebra of~$\scrB\otimes \scrA$
    generated by 
\begin{equation*}
    \{ \,b\otimes a\colon \,b\in \varrho^{-1}(\scrS),\,
a\in\scrA\,\}.
\end{equation*}
(Hint: express $\varrho^{-1}(\scrS)$ 
as pullback in~$\W{miu}$
of $\varrho\circ \pi_1, e\circ \pi_2 \colon \scrB\oplus\scrS\to\scrC$,
    where $e\colon \scrS\to\scrC$ is the inclusion.)
\end{point}
\begin{point}{80}[tensor-closed]{Theorem (Kornell)}%
The functor $(\,\cdot\,)\otimes\scrA\colon\W{miu}\to\W{miu}$
has a left adjoint
    $\Define{(\,\cdot\,)^{*\scrA}}$
for every von Neumann algebra~$\scrA$.
	\index{free exponential}%
	\index{*AstarB@$\scrA^{*\scrB}$, free exponential}
    \begin{point}{90}[tensor-closed-proof]{Proof}%
	The category $\W{miu}$ is (small-)complete,
and
$(-)\otimes\scrA\colon\W{miu}\to\W{miu}$
preserves (small-)products and equalisers.
Thus,
by Freyd's (General)
Adjoint Functor Theorem~\cite[Thm.~V.6.2]{maclane},
it suffices to check the following Solution Set Condition
(where we've used that $\W{miu}$
is locally small).
\begin{itemize}
\item
	For each $\scrB\in\W{miu}$, there is a small subset $\mathcal{S}$ of objects in $\W{miu}$
such that every arrow $h\colon \scrB\to\scrC\otimes\scrA$
can be written as a composite $h=(t\otimes\id_{\scrA})\circ f$ for some $\scrD\in \mathcal{S}$,
$f\colon\scrB\to\scrD\otimes \scrA$, and $t\colon \scrD\to\scrC$.
\end{itemize}
Let $\scrB$ be an arbitrary von Neumann algebra.
We claim that the following set $\mathcal{S}$ satisfies the required condition:
\[
	\mathcal{S}=
\{\,\scrD\colon\,
    \text{$\scrD$ is a von Neumann algebra on a subset of~$\kappa$}\,\},
\quad\text{where}\quad\kappa=2^{2^{\#\C\cdot\#\scrB\cdot2^{\#\scrA}}}.
\]
Note that~$\kappa$
being an ordinal number
is just the set of all ordinal numbers~$\alpha<\kappa$.
To prove the claim,
suppose that $h\colon \scrB\to\scrC\otimes\scrA$ is given.
By \sref{equaliser-lemma},
$h$ factors
as
\begin{equation*}
	\xymatrix@C=4em{\scrB\ar[r] & \tilde\scrC\otimes 
\scrA\ar[r]|-{\iota\otimes \id} & \scrC\otimes\scrA },
\end{equation*}
where~$\tilde\scrC$
is a von Neumann algebra
generated by no more than~$\#\scrB\cdot2^{\#\scrA}$
elements.
It follows that~$\tilde\scrC$ has no
more than~$\kappa$ elements (by~\sref{vn-generation-bound}).
Thus we may assume without loss of generality
that~$\tilde\scrC$
is a subset of~$\kappa$,
that is,~$\tilde\scrC\in\mathcal{S}$.\qed
\end{point}
\end{point}
\begin{point}{100}[cstar-no-model]{Remark}%
It should be noted that analogues
of the first and second adjunctions
can be found in the setting of $C^*$-algebras,
which raises the question as to whether
a variation on the free exponential exist for $C^*$-algebras,
that is, is there a tensor~$\otimes$ on~$\Cstar{miu}$
such that $(-)\otimes \scrA\colon
\Cstar{miu}\to\Cstar{miu}$
has a left adjoint?

Such a tensor does not exist
if we require that
on commutative
$C^*$-algebras it is given by the product of the spectra
(as is the case for the projective and injective tensors
of $C^*$-algebras)
in the sense that there
is a natural isomorphism
$\Phi_{X,Y}\colon 
\smash{\xymatrix{
C(X)\otimes C(Y)
\ar[r]|-{\cong}
&
C(X\times  Y)}}$ 
between the obvious functors of type $\CH\times \CH\to \op{(\Cstar{miu})}$.
Indeed, if~$(-)\otimes \scrA\colon \Cstar{miu}\to\Cstar{miu}$
had a left adjoint 
and so would preserve all limits
for all $C^*$-algebras~$\scrA$,
then the functor~$(-)\times X\colon \CH\to\CH$
would preserve all colimits
for every compact Hausdorff space~$X$,
which it does not,
because if it did
the square $\beta\N \times \beta\N$ of
the Stone--\v{C}ech compactification~$\beta\N$
of the natural numbers (being the $\N$-fold
coproduct of the one-point space)
would be homeomorphic to the Stone--\v{C}ech
compactification~$\beta(\N\times \N)$
of~$\N\times \N$,
which it is not (by Theorem~1 of~\cite{glicksberg1959}).

Whence
$\Cstar{cpsu}$
does not form a model of the quantum lambda calculus
in the same way that~$\W{cpsu}$ does.
\end{point}
\end{parsec}
\subsection{Hereditarily Atomic Von Neumann Algebras}
\begin{parsec}{1251}
\begin{point}{10}
We'll argue that it's possible to modify our model of the quantum lambda 
calculus from~\cite{model} to include only 
hereditarily atomic (\sref{def:hereditarily-atomic})
von Neumann algebras
    (as suggested by Kornell on page~5 of~\cite{kornell2018quantum}.)
To this end
we must bring up that the types of the quantum lambda calculus
are generated as follows:
there's a type $\qbit$, and a type~$\top$;
and from types~$A$ and~$B$, we can form\footnote{The 
type $\bit$ discussed in~\sref{qlcm-intro} is missing from this list,
since it can defined by $\bit:=\top\oplus\top$.}
\begin{equation*}
A\oplus B,\quad A\otimes B,\quad !A,\quad \text{and}\quad 
     A\limp B.
\end{equation*}
Note that 
the interpretations,
$\sem{\qbit}=M_2$
    and $\sem{\top}=\C$,
of the ground types
are hereditarily atomic,
and that the interpretation
of
the sum,
$\sem{A\oplus B}=\sem{A}\oplus\sem{B}$, and the tensor,
$\sem{A\otimes B} = \sem{A}\otimes \sem{B}$,
are hereditarily atomic
when~$\sem{A}$ and~$\sem{B}$ are hereditarily atomic.
The interpretation
$\sem{!A} = \linf(\nsp(\sem{A}))$
is hereditarily atomic regardless
of whether $\sem{A}$ is hereditarily atomic, or not.
So whether all von Neumann algebras
in our model are hereditarily atomic
hinges only on the interpretation
of~$\limp$.
As it turns out, 
the interpretation $\sem{A\limp B} = \mathcal{F}(\sem{B})^{*\sem{A}}$
we chose is not always hereditarily atomic
when $\sem{A}$ and~$\sem{B}$
are hereditarily atomic:
    we claim (without proof) that
$\sem{\top^{\oplus 3} \limp \top}
\equiv \C^{*\C^3}$
has $\scrB(\ell_2)$  as factor,
and that $\sem{\top\limp \bit}\equiv\mathcal{F}(\C^2)$
has $L^\infty[0,1]$ as summand.
\end{point}
\begin{point}{20}
The solution is obvious:
show that the functor $(\,\cdot\,)\otimes \scrA\colon
\haW{miu}\to\haW{miu}$
has a left adjoint 
    $(\,\cdot\,)^{*_\mathrm{ha}\scrA}$
for every hereditarily atomic von Neumann 
algebra~$\scrA$,
and show that the inclusion $\haW{miu}\to\haW{cpsu}$
    has a left adjoint $\mathcal{F}_{\mathrm{ha}}$.
One may then define $\sem{\,\cdot\,}_\mathrm{ha}$
exactly the same as~$\sem{\,\cdot\,}$ except for
\begin{equation*}
    \sem{A\limp B}_{\mathrm{ha}}
    \ :=\ \mathcal{F}_\mathrm{ha}(\,
    \sem{B}_\mathrm{ha}\, )^{*_\mathrm{ha}\sem{A}_{\mathrm{ha}}}.
\end{equation*}
The benefit of using the hereditarily atomic model
is that
$\mathcal{F}_\mathrm{ha}$
and $\scrA^{*_\mathrm{ha}\scrB}$
admit a significantly more concrete
    description see~\sref{Fha-concrete}
    and~\sref{AstarhaB-concrete}
A potential drawback might be that the purely quantum mechanical is restricted
to finite dimensions, so to speak.
\end{point}%
\end{parsec}%
\begin{parsec}{1252}
\begin{point}{10}
We establish the existence of~$\mathcal{F}_\mathrm{ha}$
indirectly at first.
\end{point}
\begin{point}{20}{Proposition}%
The inclusion $\haW{miu}\to\haW{cpsu}$
has a left adjoint
\begin{equation*}
    \Define{\mathcal{F}_\mathrm{ha}}(\scrA)
    \colon \haW{cpsu}\longrightarrow \haW{miu}.
\end{equation*}%
\index{Fha@$\mathcal{F}_\mathrm{ha}\colon \haW{cpsu}\longrightarrow \haW{miu}$}%
\spacingfix{}%
\begin{point}{30}{Proof}%
Given our definition of hereditary atomicity,
    \sref{def:hereditarily-atomic},
it's  pretty clear 
that the subcategory
$\haW{miu}$
    of~$\W{miu}$
is closed under
products,
and that these products are preserved
by the inclusion functor
    $\haW{miu}\to\haW{cpsu}$.
Using~\sref{ha-equalisers}
one sees the same holds for equalisers.
Whence the proof is completed by
an application of Freyd's adjoint functor theorem,
exactly as in~\sref{second-adjunction-proof},
but with as solution set
for a hereditarily atomic von Neumann algebra~$\scrA$,
the ncpsu-maps 
$\gamma\colon \scrA\to\scrC$
for which $\scrC$
    is a \emph{hereditarily atomic} 
    von Neumann algebra on a subset of
    the cardinal $\kappa\equiv\smash{ 2^{2^{\#\C+\#\scrA}}}$.
    \qed
\end{point}
\end{point}
\end{parsec}
\begin{parsec}{1253}%
\begin{point}{10}
To give a concrete description
of the functor~$\mathcal{F}_\mathrm{ha}$
we need some notation first.
\end{point}
\begin{point}{20}
Let~$\scrA$ be a hereditarily atomic von Neumann algebra.
    We'll describe $\mathcal{F}_\mathrm{ha}(\scrA)$
    in terms of ncpsu-maps $f\colon \scrA\to M_{N_f}$
    with $W^*(f(\scrA))=M_{N_f}$.
    Let us say that two such maps $f_1\colon \scrA\to M_{N_{f_1}}$
    and $f_2 \colon \scrA\to M_{N_{f_2}}$
are \Define{miu-equivalent}\index{miu-equivalent}
when there is an nmiu-isomorphism
    $\varphi\colon M_{N_{f_1}}\to M_{N_{f_2}}$
with~$\varphi \circ f_1 = f_2$,
(which implies that  $N_{f_1}=N_{f_2}$.)
Choose a set~$\Define{R_\scrA}$\index{RA@$R_\scrA$}
of representatives
for this miu-equivalence.
\end{point}
    \begin{point}{30}[Fha-concrete]{Theorem}%
Given a hereditarily atomic von Neumann algebra~$\scrA$,
the unique nmiu-map~$\Phi$ that causes the diagram
\begin{equation*}
\xymatrix@C=6em{
    \scrA
    \ar[r]^-{\eta_\scrA}
    \ar[rd]_-{\left<r\right>_{r\in R_\scrA}}
    &
    \mathcal{F}_\mathrm{ha}(\scrA)
    \ar[d]^\Phi
    \\
    &
    \bigoplus_{r\in R_\scrA} M_{N_r}
}
\end{equation*}
to commute is an nmiu-isomorphism.
Here~$\eta$ denotes the unit of the adjunction
    between $\mathcal{F}_\mathrm{ha}$
    and the inclusion
    $\haW{miu}\to\haW{cpsu}$.
\begin{point}{40}[Fha-concrete-proof]{Proof}%
Since~$\mathcal{F}_\mathrm{ha}(\scrA)$
is hereditarily atomic,
it's nmiu-isomorphic to 
a direct sum of the form $\bigoplus_{i\in I} M_{N_i}$.
We may as well assume
that $\mathcal{F}_\mathrm{ha}(\scrA)\equiv
\bigoplus_{i\in I} M_{N_i}$.
We claim,
writing $\eta_\scrA\equiv \left<s_i\right>_{i\in I}
    \colon \scrA\to\bigoplus_{i\in I}M_{N_i}$,
that the $s_i$ form a set of representatives
for miu-equivalence as well.

The theorem follows easily from this claim.
Indeed, if given~$r\in R_\scrA$
we denote
by~$i_r$ the unique element of~$I$
    for which $s_{i_r}$ is miu-equivalent to~$r$,
    and let~$\varphi_r\colon M_{N_{i_r}}\to M_{N_{r}}$
    be a corresponding nmiu-isomorphism with~$r = \varphi_r\circ s_{i_r}$,
then one easily sees using its defining property
    that $\Phi$ is the composition of
\begin{equation*}
\xymatrix@C=5.5em{
\mathcal{F}_\mathrm{ha}(\scrA)\,\equiv\,
\bigoplus_{i \in I} M_{N_i}
    \ar[r]^-{\left<\pi_{i_r} \right>_{r\in R_\scrA}}
    & 
    \bigoplus_{r\in R_\scrA}
    M_{N_{r_i}}
    \ar[r]^-{\bigoplus_{r\in R_\scrA} \varphi_r}
    &
    \bigoplus_{r\in R_\scrA} M_{N_r}.
}
\end{equation*}
Since $r\mapsto i_r$ gives a bijection~$R_\scrA\to I$,
the first map above
is a nmiu-isomorphism.
Since the second map is clearly a nmiu-isomorphism
too, $\Phi$ is a nmiu-isomorphism.
\begin{point}{50}
Let us begin by proving that $W^*(s_i(\scrA))=M_{N_i}$
for every $i\in I$.

To this end, we'll first show 
that~$W^*(\eta_\scrA(\scrA))=\mathcal{F}_\mathrm{ha}(\scrA)$.
Let us denote by $f\colon \scrA\to W^*(\eta_\scrA(\scrA))$
the restriction of~$\eta_\scrA$.
By the universal property of~$\eta_\scrA$,
    there's a unique nmiu-map $\varrho\colon \mathcal{F}_\mathrm{ha}(\scrA)
    \to W^*(\eta_\scrA(\scrA))$
    with $f = \varrho\circ \eta_\scrA$.
Letting $e\colon W^*(\eta_\scrA(\scrA))\to \mathcal{F}_\mathrm{ha}(\scrA)$
be the inclusion,
we have $e\circ \varrho \circ \eta_\scrA = e\circ f = \eta_\scrA$.
    Since the identity $\id_{\mathcal{F}_\mathrm{ha}(\scrA)}\colon \mathcal{F}_\mathrm{ha}(\scrA)
\to \mathcal{F}_\mathrm{ha}(\scrA)$
is the unique nmiu-map $\tau\colon \mathcal{F}_\mathrm{ha}(\scrA)
\to \mathcal{F}_\mathrm{ha}(\scrA)$
with $\tau\circ \eta_\scrA = \eta_\scrA$,
    we get~$\id_{\mathcal{F}_\mathrm{ha}(\scrA)}= e\circ \varrho$.
    Since~$\id_{\mathcal{F}_\mathrm{ha}(\scrA)}$ is surjective,
    so is~$e$, and thus $W^*(\eta_\scrA(\scrA))=\mathcal{F}_\mathrm{ha}(
    \scrA)$.

Note that $W^*(\eta_\scrA(\scrA))
    \subseteq \bigoplus_{i\in I}
    W^*(s_i(\scrA))$
    because $\bigoplus_{i\in I} W^*(s_i(\scrA))$
    is a von Neumann subalgebra
    of~$\mathcal{F}_\mathrm{ha}(\scrA)$
    with $\eta_\scrA(\scrA) 
    \,\subseteq\, \bigoplus_{i\in I} W^*(s_i(\scrA))$.
    Since $W^*(\eta_\scrA(\scrA))=
    \mathcal{F}_\mathrm{ha}(
    \scrA)$,
    we get 
    $\bigoplus_{i\in I} W^*(s_i(\scrA))=
    \mathcal{F}_\mathrm{ha}(\scrA)
    \equiv \bigoplus_{ i \in I} M_{N_i}$,
    and so $W^*(s_i(\scrA))=M_{N_i}$
    for all~$i\in I$.

It remains to be shown
that given an ncpsu-map
$f\colon \scrA\to M_{N_f}$
    with $W^*(f(\scrA))= M_{N_f}$
there's a unique $i\in I$
such that~$s_i$ is miu-equivalent to~$f$.
\end{point}
\begin{point}{60}{Uniqueness}%
Let~$i,j\in I$ such that $s_i$ and~$s_j$ are miu-equivalent
be given, 
and let $\varphi\colon M_{N_i}\to M_{N_j}$
be the associated nmiu-isomorphism with~$\varphi \circ s_i = s_j$.
We must show that~$i=j$.
Recall that by the universal property of~$\eta_\scrA$,
there's a unique
nmiu-map $\varrho\colon \mathcal{F}_\mathrm{ha}(\scrA)
\to M_{N_j}$
with $s_j = \varrho\circ \eta_\scrA$.
    Surely, 
\begin{equation*} 
    \textstyle
    \pi_j\colon \mathcal{F}_\mathrm{ha}(\scrA)
    \equiv \bigoplus_{i'\in I}  M_{N_{i'}}
    \longrightarrow M_{N_j}
\end{equation*}
fits this description;
but so does $\varphi\circ \pi_i$,
since $\varphi\circ \pi_i \circ \eta_\scrA = \varphi\circ s_i = s_j$.
Hence~$\pi_j = \varphi\circ \pi_i$.
    This entails that the carriers (see~\sref{carrier-miu}) 
    of~$\pi_j$ and~$\pi_i$
are equal, so~$i=j$.
\begin{point}{70}[Fha-concrete-existence]{Existence}
Let~$f\colon \scrA\to M_{N_f}$
be an ncpsu-map with $M_{N_f}=W^*(f(\scrA))$.
We must show that there is an~$i\in I$
such that~$s_i$
is miu-equivalent with~$f$.
By the universal property of~$\eta_\scrA$
there's a unique nmiu-map $\varrho\colon \mathcal{F}_\mathrm{ha}(\scrA)
\to M_{N_f}$ with $f = \varrho\circ \eta_\scrA$.
We claim that~$\varrho$ must be of the form~$\varrho'\circ \pi_i$
for some nmiu-isomorphism $\varrho'\colon M_{N_i}\to M_{N_f}$.

First note that~$\varrho$ is surjective:
indeed, $\varrho(\mathcal{F}_\mathrm{ha}(\scrA))$
is a von Neumann subalgebra of~$M_{N_f}$ by~\sref{nmiu-image},
that contains $f(\scrA)$.
Thus $\varrho(\mathcal{F}_\mathrm{ha}(\scrA))
    \,\supseteq\, W^*(f(\scrA))\equiv M_{N_f}$,
which implies that~$\varrho$ is surjective.
Since~$\varrho$ is surjective, it maps central projections
    of $\mathcal{F}_\mathrm{ha}(\scrA)$
    to central projections
    of~$M_{N_f}$.
    For each~$i\in I$ let~$c_i$
    denote the central projection in $\mathcal{F}_\mathrm{ha}(\scrA)$
    given by $c_i(i)=1$ and $c_i(j)=0$ for all~$j\neq i$.
    Then the $c_i$ form an orthogonal family
    of  central projections with~$\sum_{i\in I}c_i =1$.
    So the $\varrho(c_i)$
    form an orthogonal family
    of central projections of~$M_{N_f}$ as well, 
    with $\sum_{i\in I} \varrho(c_i)=1$.
Since the only non-zero central projection in 
    $M_{N_f}$ (being a factor, \sref{central-examples}) is~$1$,
    it follows that there is exactly one~$i\in I$
    with $\varrho(c_i)=1$,
    and that $\varrho(c_j)=0$ for all~$j\neq i$.
    From this one easily deduces (c.f.~\sref{nmiu-factors}) 
    that~$\varrho$ must be of the form~$\varrho = \varrho'\circ \pi_i$
    for some injective nmiu-map $\varrho'\colon M_{N_i}\to M_{N_f}$.
    Since~$\varrho$ is surjective, $\varrho'$ is surjective too,
    and thus $\varrho'$ is a nmiu-isomorphism.

Now, since $\varrho = \varrho' \circ \pi_i $,
and $f=\varrho\circ \eta_\scrA$,
we get $f = \varrho' \circ \pi_i \circ \eta_\scrA\equiv
\varrho' \circ s_i$.
Since~$\varrho'$ is a nmiu-isomorphism,
we see that $f$ is miu-equivalent to~$s_i$.\qed
\end{point}
\end{point}
\end{point}
\end{point}
\begin{point}{80}{Remark}%
Given the concrete description for~$\mathcal{F}_\mathrm{ha}$
from~\sref{Fha-concrete}
it seems tempting to prove directly
that $\left<r\right>_{r\in R_\scrA}
\colon \scrA\to\bigoplus_{r\in R_\scrA}M_{N_r}$
is a universal arrow from~$\scrA$ to the inclusion
$\haW{miu}\to\haW{cpsu}$,
without presupposing the existence of~$\mathcal{F}_\mathrm{ha}$.
However, our attempts to do so have been thwarted
by our inability to prove  that
$W^*(\left<r\right>_{r\in R_\scrA})
= \bigoplus_{r\in R_\scrA} M_{N_r}$
using elementary means.

Whether this indicates an error in our proof,
or the power of the adjoint functor theorem,
remains to be seen.
\end{point}
\end{parsec}%
\begin{parsec}{1254}%
\begin{point}{10}%
To allow interpretation of~$\limp$
in~$\haW{miu}$,
we'll show that
the functor
$(\,\cdot\,)\otimes \scrB\colon \haW{miu}\to\haW{miu}$
    has a left adjoint $(\,\cdot\,)^{*_\mathrm{ha}\scrB}$ for every
hereditarily atomic von Neumann algebra~$\scrB$.
This result has already been established by Kornell
(in Theorem~9.1 of~\cite{kornell2018quantum});
we improve upon it by giving a different,
slightly more concrete, description.
    As was the case for~$\mathcal{F}_\mathrm{ha}$,
    we establish the existence
    of~$(\,\cdot\,)^{*_\mathrm{ha}\scrB}$
indirectly at first.
\end{point}
\begin{point}{20}{Proposition}%
Given a hereditarily atomic von Neumann algebra~$\scrB$,
    the functor $(\,\cdot\,)\otimes\scrB\colon
    \haW{miu}\longrightarrow \haW{miu}$
    has a left adjoint
    $\Define{(\,\cdot\,)^{*_\mathrm{ha}\scrB}}\colon
    \haW{miu}\to\haW{miu}$%
\index{*AstarBha@$\scrA^{*_\mathrm{ha}\scrB}$}.
\begin{point}{30}{Proof}%
We already know from~\sref{tensor-closed} that~$(\,\cdot\,)\otimes \scrB$
preserves limits as functor $\W{miu}\to\W{miu}$.
Since the subcategory $\haW{miu}$
of~$\W{miu}$ is closed
under products and equalisers,
    the restriction of~$(\,\cdot\,)\otimes\scrB$
    to a functor $\haW{miu}\to\haW{miu}$
    preserves limits as well.
The proof is now completed by an application
of Freyd's adjoint functor theorem,
exactly as in~\sref{tensor-closed-proof},
    but with a suitably modified solution set.\qed
\end{point}
\end{point}
\end{parsec}
\begin{parsec}{1255}
\begin{point}{10}
To describe $\scrA^{*_\mathrm{ha}\scrB}$ concretely
we need some notation.
\end{point}
\begin{point}{20}{Definition}
We say that a nmiu-map~$s\colon \scrA\to\scrC\otimes \scrB$,
where $\scrA$, $\scrB$ and~$\scrC$
are von Neumann algebras,
is \Define{$(\,\cdot\,)\otimes\scrB$-surjective}%
\index{surjective@$(\,\cdot\,)\otimes\scrB$-surjective}
when the only von Neumann subalgebra~$\scrS$
of~$\scrC$ with $s(\scrA)\subseteq \scrS\otimes\scrB$
is~$\scrS=\scrC$,
where for the sake of simplicity we regard $\scrS\otimes \scrB$
to be a von Neumann subalgebra of~$\scrC\otimes \scrB$
(c.f.~\sref{tensor-injective}).
\begin{point}{21}[tensor-map-factorisation]
    By inspecting the proof of~\sref{equaliser-lemma}
one sees that for any nmiu-map $s\colon \scrA\to\scrC\otimes \scrB$
there is a von Neumann subalgebra~$\tilde\scrC$
of~$\scrC$
such that~$s(\scrA)\subseteq \tilde\scrC\otimes \scrB$,
and the restriction of~$s$
to a a map $s\colon \scrA\to\tilde\scrC\otimes \scrB$
    is $(\,\cdot\,)\otimes\scrB$-surjective.
\end{point}
\end{point}
\begin{point}{30}[tensorBsurjectivity]{Lemma}
Given a $(\,\cdot\,)\otimes \scrB$-surjective
nmiu-map $s\colon \scrA\to\scrC\otimes \scrB$
and a nmiu-map  $\varrho\colon \scrC\to\scrD$
between von Neumann algebras,
the composition 
\begin{equation*}
\xymatrix@C=3em{
\scrA
\ar[r]^-s
    &
\scrC\otimes\scrB
    \ar[r]^-{\varrho\otimes\scrB}
    &
    \scrD\otimes\scrB
}
\end{equation*}
is $(\,\cdot\,)\otimes\scrB$-surjective
iff $\varrho$ is surjective.
\begin{point}{40}{Proof}%
Suppose that~$\varrho$ is surjective,
and let~$\scrS$ be a von Neumann subalgebra
of~$\scrD$ with 
    $(\varrho\otimes \scrB)(s(\scrA))\subseteq \scrS\otimes \scrB$.
To prove that $(\varrho\otimes \scrB)\circ s$
is $(\,\cdot\,)\otimes \scrB$-surjective,
we must show that~$\scrS=\scrD$.
Since $(\varrho\otimes \scrB)(s(\scrA))\subseteq \scrS\otimes \scrB$,
we have 
$s(\scrA)\subseteq (\varrho\otimes\scrB)^{-1}(\scrS \otimes \scrB)
\equiv \varrho^{-1}(\scrS) \otimes \scrB$,
by~\sref{tensor-preimage},
and so~$\varrho^{-1}(\scrS)=\scrC$,
because~$s$ is $(\,\cdot\,)\otimes \scrB$-surjective.
Whence~$\scrS = \varrho(\varrho^{-1}(\scrS))
\equiv \varrho(\scrC) = \scrD$,
using here that~$\varrho$ is surjective.
\begin{point}{50}
For the other direction suppose that~$(\varrho\otimes \scrB)\circ s$
is $(\,\cdot\,)\otimes\scrB$-surjective.
Since the range of~$\varrho\otimes \scrB$ 
is~$\varrho(\scrC)\otimes \scrB$,
we have 
$(\varrho\otimes\scrB)(s(\scrA))\subseteq \varrho(\scrC)\otimes \scrB$,
and so~$\varrho(\scrC)=\scrD$,
because~$(\varrho\otimes \scrB)\circ s$
is $(\,\cdot\,)\otimes\scrB$-surjective.\qed
\end{point}
\end{point}
\end{point}
\begin{point}{60}{Definition}%
Let~$\scrA$ and~$\scrB$ be hereditarily atomic von Neumann algebras.
We'll describe $\scrA^{*_\mathrm{ha}\scrB}$
in terms
of the $(\,\cdot\,)\otimes \scrB$-surjective
    nmiu-maps $f\colon \scrA\to M_{N_f}\otimes \scrB$.
Let us say that two such maps
    $f_1\colon \scrA\to M_{N_{f_1}}\otimes \scrB$
    and $f_2 \colon \scrA \to M_{N_{f_2}}\otimes \scrB$
    are 
\Define{$(\,\cdot\,)\otimes\scrB$-equivalent}%
\index{equivalent@$(\,\cdot\,)\otimes\scrB$-equivalent}
when there is a nmiu-isomorphism
    $\varphi\colon M_{N_{f_1}}\to M_{N_{f_2}}$
    with $(\varphi\otimes\scrB) \circ f_1 = f_2$
    (which implies that $N_{f_1}=N_{f_2}$.)
Pick a set of representatives
    $\Define{S_{\scrA,\scrB}}$%
    \index{SAB@$S_{\scrA,\scrB}$}
    for this 
$(\,\cdot\,)\otimes\scrB$-equivalence.
\end{point}
\begin{point}{70}[AstarhaB-concrete]{Theorem}
Let~$\scrA$ and~$\scrB$ be hereditarily atomic 
von Neumann algebras.
the unique nmiu-map
$\Phi\colon 
    \scrA^{*_\mathrm{ha}\scrB}
    \longrightarrow
    \bigoplus_{s\in S_{\scrA,\scrB}} M_{N_s}$
    that makes the diagram
\begin{equation*}
\xymatrix@C=10em{
\scrA
    \ar[r]^-{\eta_{\scrA,\scrB}}
    \ar[rdd]_-{\left<s\right>_{s\in S_{\scrA,\scrB}}}
&
\scrA^{*_\mathrm{ha}\scrB}\otimes \scrB
    \ar[d]^-{\Phi\otimes \scrB}
\\
&
\bigl(\,\bigoplus_{s\in S_{\scrA,\scrB}} M_{N_s}\,\bigr)
\otimes \scrB
    \ar[d]_{\cong}^{\left<\pi_s\otimes \scrB\right>_{s\in S_{\scrA,\scrB}}}
\\
&
\bigoplus_{s\in S_{\scrA,\scrB}} M_{N_s}\otimes \scrB
}
\end{equation*}
commute
is a nmiu-isomorphism.
Here $\eta_{(\,\cdot\,),\scrB}$
    denotes
the unit of the adjunction
between $(\,\cdot\,)^{*_\mathrm{ha}\scrB}$
and~$(\,\cdot\,)\otimes\scrB$.
\begin{point}{80}{Proof}%
We follow roughly the same lines
    as the proof of~\sref{Fha-concrete}.
Since~$\scrA^{*_\mathrm{ha}\scrB}$
is hereditarily atomic
we write~$\scrA^{*_\mathrm{ha}\scrB}
\equiv \bigoplus_{i\in I} M_{N_i}$
without loss of generality.
Note that writing $e_i = (\pi_i\otimes \scrB)\circ \eta_{\scrA,\scrB}$
the diagram
\begin{equation*}
\xymatrix@C=10em{
\scrA
    \ar[r]^-{\eta_{\scrA,\scrB}}
    \ar[rd]_-{\left<e_i\right>_{i\in I}}
&
\scrA^{*_\mathrm{ha}\scrB}\otimes \scrB
    \,\equiv\, \bigoplus_{i\in I}M_{N_i}\otimes \scrB
    \ar[d]_{\cong}^{\left<\pi_i\otimes \scrB\right>_{i \in I}}
\\
&
\bigoplus_{i \in I} M_i\otimes \scrB
}
\end{equation*}
commutes.
We claim that the $e_i$ are $(\,\cdot\,)\otimes \scrB$-surjective,
and, moreover,
form a set of representatives
for $(\,\cdot\,)\otimes \scrB$-equivalence
on the set of nmiu-maps
    $f\colon \scrA\to M_{N_f}\otimes \scrB$.
From this claim the theorem follows
with a reasoning similar to that
in~\sref{Fha-concrete-proof}, which we won't repeat here.
\begin{point}{90}
To prove that $e_i\equiv (\pi_i \otimes \scrB)\circ \eta_{\scrA,\scrB}$
is $(\,\cdot\,)\otimes\scrB$-surjective,
it suffices,
by~\sref{tensorBsurjectivity},
to show that $\eta_{\scrA,\scrB}$
is $(\,\cdot\,)\otimes \scrB$-surjective
(since~$\pi_i$ is surjective.)
So let~$\scrS$ be a von Neumann subalgebra of~$\scrA^{*_\mathrm{ha}\scrB}$
such that $\eta_{\scrA,\scrB}(\scrA)\subseteq
\scrS \otimes \scrB$.
We must show that~$\scrS = \scrA^{*_\mathrm{ha}\scrB}$.
Letting~$f\colon \scrA\to\scrS\otimes \scrB$
denote the restriction of~$\eta_{\scrA,\scrB}$,
there is, by the universal property of~$\eta_{\scrA,\scrB}$,
a unique nmiu-map $\varrho\colon \scrA^{*_\mathrm{ha}\scrB}
\to \scrS$
such that~$f= (\varrho\otimes \scrB) \circ \eta_{\scrA,\scrB}$.
Note that if we compose $\varrho$ with the inclusion
    $e\colon \scrS\to \scrA^{*_\mathrm{ha}\scrB}$,
    then
we get a nmiu-map $\sigma := e\circ \varrho\colon 
    \scrA^{*_\mathrm{ha}\scrB}
    \to\scrA^{*_\mathrm{ha}\scrB}$
    with the property
    that 
    $(\sigma\otimes\scrB)\circ \eta_{\scrA,\scrB} = \eta_{\scrA,\scrB}$.
Since the identity on~$\scrA^{*_\mathrm{ha}\scrB}$
is the only map with this property,
we get $e \circ \varrho = \id$.
This implies that~$e$ is surjective,
and thus~$\scrS=
    \scrA^{*_\mathrm{ha}\scrB}$.
    Whence~$\eta_{\scrA,\scrB}$
    is $(\,\cdot\,)\otimes \scrB$-surjective.

It remains to be show that for every nmiu-map
$f\colon \scrA\to M_{N_f}\otimes \scrB$
there is a unique $i\in I$
such that~$e_i$
is $(\,\cdot\,)\otimes \scrB$-equivalent
to~$f$.
\end{point}
\begin{point}{100}{Uniqueness}%
Suppose that $e_i$ and~$e_j$ are $(\,\cdot\,)\otimes\scrB$-equivalent
for some~$i,j\in I$;
we must show that~$i=j$.
Let~$\varphi\colon M_{N_i} \to M_{N_j}$
be an nmiu-isomorphism 
with $(\varphi\otimes \scrB ) \circ e_i = e_j$.
Note that $\pi_j\colon \scrA^{*_\mathrm{ha}\scrB}
    \equiv \bigoplus_{j'\in I} M_{N_{j'}}\longrightarrow
M_{N_j}$
    is the unique nmiu-map~$\varrho\colon \scrA^{*_\mathrm{ha}\scrB}
    \to M_{N_j}$ 
with $e_j = \varrho\otimes \scrB \circ \eta_{\scrA,\scrB}$.
    Since $(\varphi \circ \pi_i)\otimes \scrB \circ \eta_{\scrA,\scrB}
    = (\varphi\otimes \scrB) \,\circ\, (\pi_i\otimes\scrB )
    \,\circ\, \eta_{\scrA,\scrB}
    = \varphi\otimes \scrB \circ e_i = e_j$,
    we get $\varphi \circ \pi_j = \pi_i$.
    Hence~$i=j$.
\end{point}
\begin{point}{110}{Existence}%
Let $f\colon \scrA\to M_{N_f}\otimes\scrB$
be a $(\,\cdot\,)\otimes \scrB$-surjective  nmiu-map.
We must show that there is a unique $i\in I$
such that $f$ is $(\,\cdot\,)\otimes\scrB$-equivalent
to~$e_i$.
By the universal property of~$\eta_{\scrA,\scrB}$,
there's a unique nmiu-map
$\varrho\colon \scrA^{*_\mathrm{ha}\scrB}\longrightarrow
M_{N_f}$
with $(\varrho\otimes \scrB)\circ \eta_{\scrA,\scrB}
=  f$.
Note that~$\varrho$
    is surjective by~\sref{tensorBsurjectivity}.
Now, following the same reasoning as in~\sref{Fha-concrete-existence},
we see that~$\varrho\colon \scrA^{*_\mathrm{ha}\scrB}\equiv
\bigoplus_{i \in I} M_{N_i}\longrightarrow M_{N_f}$
must be of the form~$\varrho\equiv \varrho'\circ \pi_i$
for some~$i\in I$ and nmiu-isomorphism $\varrho'\colon
M_{N_i}\to M_{N_f}$.
So~$ f= (\varrho\otimes \scrB) \circ \eta_{\scrA,\scrB}
    = (\varrho'\otimes\scrB) \,\circ\, 
    (\pi_i\otimes \scrB) \,\circ\, \eta_{\scrA,\scrB}
    = (\varrho'\otimes\scrB) \circ e_i$,
    and hence $f$ is $(\,\cdot\,)\otimes\scrB$-equivalent
    to~$e_i$.\qed
\end{point}
\end{point}
\end{point}
\end{parsec}
\section{Duplicators and Monoids}
\label{S:duplicable}
\begin{parsec}{1260}%
\begin{point}{10}%
When asked for an
interpretation of the type $!A$
as a von Neumann algebra
\begin{equation}
\label{fock}
\sem{!A}
\ = \ 
\bigoplus_n \sem{A}^{\otimes n}
\end{equation}
definitely seems
like a suitable answer
given
the cue that~$!A$ should represent
as many instances of~$A$ as needed,
which makes the interpretation
we actually use in our model of the quantum lambda calculus
(namely $\sem{!A}=\ell^\infty(\nsp(\sem{A}))$)
rather
suspect.
To address such concerns
we'll show that any von Neumann algebra
that carries a $\otimes$-monoid structure
(in~$\W{miu}$ as~$\sem{!A}$ should)
must be nmiu-isomorphic
to $\ell^\infty(X)$ for some set~$X$
(see \sref{duplicable})
ruling out the interpretation~\eqref{fock}
for all but the most trivial cases.
We'll show in fact that~$\ell^\infty(\nsp(\scrA))$
is the free $\otimes$-monoid
over~$\scrA$
in~$\W{miu}$ (see \sref{thm:free-monoid-in-vNAMIU})
exonerating it in our minds from all doubts.
\end{point}
\end{parsec}
\subsection{Duplicators}
\begin{parsec}{1270}%
\begin{point}{10}[def:duplicator]{Definition}%
A von Neumann algebra~$\mathscr{A}$
is \Define{duplicable}%
	\index{duplicable von Neumann algebra}
if there is a \Define{duplicator} on~$\mathscr{A}$,%
	\index{duplicator}
that is,
an npsu-map
$\delta\colon \mathscr{A}\otimes \mathscr{A}\to\mathscr{A}$
with a \Define{unit} $u\in [0,1]_\scrA$ 
satisfying
\begin{equation*}
\delta(a\otimes u)\ =\ a\ = \ \delta(u\otimes a)
\quad\text{for all~$a\in\mathscr{A}$.}
\end{equation*}
(Note that we require of~$\delta$ neither associativity
nor commutativity.)
\end{point}
\begin{point}{20}{Remark}
The unit $u$ can be identified with
a positive subunital map $\tilde{u}\colon \C\to \scrA$ via $\tilde{u}(\lambda)=\lambda u$.
The definition is motivated by the fact that the interpretation of $\bang A$
must carry a commutative monoid structure in $\W{miu}$.
The condition is weaker, requiring the maps to be only positive subunital,
and dropping associativity and commutativity.
Nevertheless this is sufficient to prove the following.
\end{point}
\begin{point}{30}[duplicable]{Theorem}%
A von Neumann algebra~$\mathscr{A}$
is duplicable if and only if
$\mathscr{A}$ is nmiu-isomorphic to $\linf(X)$ for some set~$X$.
In that case, the duplicator $(\delta,u)$
is unique, given by
$\delta(a\otimes b) = a\cdot b$ and $u=1$.
\begin{point}{40}%
Thus, to interpret duplicable types,
we can really only use von Neumann algebras of the form $\linf(X)$.
It also follows that a von Neumann algebra
is duplicable precisely when it is a (commutative) monoid
in $\W{miu}$, or in the 
symmetric monoidal category $\W{cpsu}$ of von Neumann algebras
and normal completely positive subunital (CPsU) maps.
\end{point}
\begin{point}{50}[sec:characterisation-dup-vna]%
To prove~\sref{duplicable}
we proceed as follows.
First we prove 
in~\sref{lem:uniqueness-duplicator}
every duplicable von Neumann algebra~$\scrA$
is commutative (and that the duplicator is given by multiplication).
This reduces the problem
to a measure theoretic one,
because~$\mathscr{A} \cong \bigoplus_i L^\infty(X_i)$
for some finite complete measure spaces~$X_i$
(by~\sref{cvn}).
Since each of the~$L^\infty(X_i)$s will be duplicable
(see~\sref{cor:duplicable-product})
we may assume without loss of generality
that~$\scrA\cong L^\infty(X)$ for some finite complete measure space~$X$.
Since~$X$ splits into a discrete and a continuous
part (see~\sref{lem:measure-space-continuous-discrete}),
and the result is obviously true for discrete spaces,
we only need to show that~$L^\infty(C)=\{0\}$
for any continuous complete finite measure space~$C$
for which~$L^\infty(C)$ is duplicable.
In fact,
we'll show that~$\mu(C)=0$
for such~$C$
(see~\sref{lem:continuous-measure-space}).
\end{point}
\end{point}
\begin{point}{60}[lem:unit-duplicator]{Lemma}%
Let~$\delta$
be a duplicator 
with unit~$u$
on a von Neumann algebra~$\mathscr{A}$.
Then~$u=1$ and~$\delta(1\otimes 1)=1$.
\begin{point}{70}{Proof}%
Since~$1=\delta(u\otimes 1)\leq \delta(1\otimes 1) \leq 1$,
	we have $\delta(1\otimes 1)=1$, and so $\delta(u^\perp\otimes 1)=0$.
But,
because  $u^\perp = \delta(u^\perp \otimes u)
\leq \delta(u^\perp \otimes 1) = 0$,
we have~$u^\perp=0$, and thus~$u=1$.
Hence~$1=\delta(1\otimes u)=\delta(1\otimes 1)$.\qed
\end{point}
\end{point}
\end{parsec}
\begin{parsec}{1280}%
\begin{point}{10}%
To prove that a duplicable von Neumann 
algebra is commutative we'll need the following two classical theorems
from the theory of $C^*$-algebras.
\end{point}
\begin{point}{20}[tomiyama]{Theorem (Tomiyama)}%
\index{Tomiyama's Theorem}
Any linear surjection
 $f\colon \scrA\to\scrB$
of a von Neumann algebra~$\scrA$ onto a von Neumann subalgebra 
$\scrB\subseteq \scrA$
with~$f(f(a))=f(a)$ and~$\|f(a)\|\leq\|a\|$ 
for all~$a\in \scrA$
obeys~$bf(a)=f(ba)$
for all~$a\in\scrA$ and~$b\in\scrB$.
\begin{point}{21}{Remark}%
The usual (see e.g.~10.5.86 of~\cite{kr}) 
and original\cite{tomiyama} formulations of Tomiyama's theorem
involve $C^*$-algebras instead of von Neumann algebras,
and include the conclusion that~$f$ is positive.
Since these improvements
weren't necessary for our purposes,
we've left them out, shortening the proof.
\end{point}
\begin{point}{30}{Proof}%
(Based on II.6.10.2 of \cite{blackadar2006operator}.)

Let~$a\in\scrA$ and~$b\in\scrB$ be given.
Since $b$ is the
norm limit of linear combinations of projections
(cf.~\sref{projections-norm-dense}),
it suffices to show that $ef(a)=f(ea)$
for every projection~$e$ from~$\scrB$ and~$a\in \scrA$.
For this, in turn,
it suffices to show that $e^\perp f(ea)=0$
for every projection~$e$ from~$\scrB$,
(and thus also~$ef(e^\perp a)=0$,)
because then $f(ea)=ef(ea)=ef(a)$.

Let~$\lambda \in \R$ be given.
The trick is to obtain the following inequality.
\begin{equation}
\label{eq:tomiyama}
(1+2\lambda) \left\|e^\perp f(ea)\right\|^2 
\ \leq\  \left\|ea\right\|^2
\end{equation}
Indeed, this inequality can only hold for all~$\lambda$
when $\|e^\perp f(ea)\|=0$.
Working towards~\eqref{eq:tomiyama},
let us first note that $f(e^\perp f(ea))=e^\perp f(ea)$:
indeed, since  $e^\perp f(ea)\in \scrB$
and~$f\colon \scrA\to\scrB$ is surjective,
    there must be~$a'\in\scrA$ with~$f(a')=e^\perp f(ea)$,
    and thus  $e^\perp f(ea)=f(a')=f(f(a'))=f(e^\perp f(ea))$.
Then:
\begin{alignat*}{3}
    &(1+\lambda)^2 \left\|e^\perp f(ea)\right\|^2\\
    &= \ \left\|e^\perp f(\,ea + \lambda e^\perp f(ea)\,)\right\|^2
    \qquad && \text{since }f(e^\perp f(ea))=e^\perp f(ea) \\
    &\leq \ \left\|\,ea\,+\,\lambda e^\perp f(ea)\,\right\|^2\qquad
    &&\text{since }\|e^\perp\|\leq 1\text{ and }\|f\|\leq 1 \\
    &=\ \left\|ea\right\|^2\,+\,\lambda^2 \left\|e^\perp f(ea)\right\|^2
    \qquad&&\text{using $\|c\|^2=\|c^*c\|$ and $ee^\perp=0$}
\end{alignat*}
Subtracting $\lambda^2 \left\|e^\perp f(ea)\right\|^2$ from both
sides yields inequality~\eqref{eq:tomiyama}.\qed
\end{point}
\end{point}
\begin{point}{40}{[Moved to~\sref{russo-dye}]}%
    \begin{point}{50}{[Removed]}%
\end{point}
\end{point}
\begin{point}{60}[lem:sef-instrument]{Lemma}%
Let~$\mathscr{A}$ be a von Neumann algebra,
and let~$f\colon \mathscr{A}\oplus\mathscr{A}\to \mathscr{A}$
be a pu-map 
with $f(a,a)=a$ for all~$a\in \mathscr{A}$.
Then $p:=f(1,0)$ is central,
and
\begin{equation*}
	f(a,b) \,=\, ap\,+\, bp^\perp
\end{equation*}
for all~$a,b\in\mathscr{A}$.
\begin{point}{70}{Proof}%
(Based on Lemma~8.3 of~\cite{newdirections}.)

Note that $(c,d)\mapsto (\,f(c,d),\,f(c,d)\,)$
gives a pu-map~$f'$ from~$\scrA\oplus \scrA$
onto its von Neumann subalgebra~$\{\,(a,a)\colon\,a\in\scrA\,\}$
with~$f'(f'(c,d))=f'(c,d)$ for all~$c,d\in\scrA$.
Since~$\|f'\|\leq 1$ as a result of Russo--Dye's theorem 
(see~\sref{russo-dye-cor}),
Tomiyama's theorem (\sref{tomiyama}) 
implies that 
for all~$a,c,d\in\mathscr{A}$,
\begin{equation*}
	(a,a)\,f'(c,d)
	\ =\  f'(\,ac,\, ad\,),
	\quad\text{and so}\quad
a \,f(c,d) \ = \ f(\,ac\,,\,ad\,).
\end{equation*}
Upon taking the adjoint
we see that $f(c,d)b=f(cb,db)$
for all~$b,c,d\in\scrA$
(using here that~$f$ being positive is involution preserving,
see~\sref{cstar-p-implies-i}.)
As a result of these observations, we get
 $ap\equiv af(1,0)=f(a,0)=f(1,0)a\equiv pa$
 for all~$a\in\scrA$,
and so~$p$ is central.
Similarly, $f(0,b)=bp^\perp$
for all~$b\in\mathscr{A}$.
Then~$f(a,b)=f(a,0)+f(0,b)=ap+bp^\perp$
for all~$a,b\in\mathscr{A}$.\qed
\end{point}
\end{point}
\begin{point}{80}[lem:uniqueness-duplicator]{Lemma}%
Let $\delta\colon\mathscr{A}\otimes \mathscr{A}\to\mathscr{A}$
be a duplicator on a von Neumann algebra~$\mathscr{A}$.
Then~$\mathscr{A}$ is commutative and~$\delta(a\otimes b)=a\cdot b$
for all~$a,b\in\mathscr{A}$.%
	\index{duplicator!is multiplication}%
	\index{duplicable von Neumann algebra!is commutative}
\begin{point}{90}{Proof}%
To prove~$\scrA$ is commutative
we must show that all~$a\in\mathscr{A}$ are central,
but, of course,
it suffices
to show that all~$p\in [0,1]_\mathscr{A}$ are central
(by the usual reasoning).
Similarly, 
we only need to prove that $\delta(a\otimes p) = a\cdot p$
for all~$a\in\mathscr{A}$ and $p\in [0,1]_\mathscr{A}$.
Given such~$p\in[0,1]_\mathscr{A}$ 
define~$f\colon \mathscr{A}\oplus\mathscr{A}\to\mathscr{A}$
by $f(a,b) = \delta(a\otimes p+b\otimes p^\perp)$
for all~$a,b\in\mathscr{A}$.
Then~$f$ is positive, unital,
$f(1,0)=p$,
and 
$f(a,a)=a$
for all~$a\in \mathscr{A}$.
Thus by~\sref{lem:sef-instrument},
$p$ is central,
and  $f(a,b)=ap+bp^\perp$ for all~$a,b\in\mathscr{A}$.
Then~$a\cdot p=f(a,0)=\delta(a \otimes p)$.\qed
\end{point}
\end{point}
\begin{point}{100}{Remark}%
The special
case of~\sref{lem:uniqueness-duplicator}
in which~$\delta$ is \emph{completely} positive
can be found
in the literature,
see for example
Theorem~6 of~\cite{Maassen2010}
(where duplication is called copying,)
and \S3.3 of~\cite{clifton2003characterizing} 
    (where it's called broadcasting.)
\end{point}
\begin{point}{101}{Remark}
Note that we've not yet used that a duplicator is normal.
That a duplicator is normal will come in to play momentarily
when we show that~$L^\infty(X)$
is duplicable iff~$X$ is discrete.
\end{point}
\begin{point}{110}[cor:duplicability-multiplication]{Corollary}
Let~$\mathscr{A}$
be a von Neumann algebra.
Then~$\mathscr{A}$
is duplicable
iff there is
an np-map $\delta\colon\mathscr{A}\otimes\mathscr{A}
\to \mathscr{A}$
with $\delta(a\otimes b)=a\cdot b$ 
for all~$a,b\in\mathscr{A}$,
(and in that case~$\scrA$ is commutative.)
\begin{point}{120}{Remark}%
Thus
for a non-commutative von Neumann algebra~$\scrA$ 
multiplication $(a,b)\mapsto ab\colon \scrA\times \scrA\to\scrA$
is not a normal bilinear map
in the sense of~\sref{tensor-extra}.
\end{point}
\end{point}
\begin{point}{130}[cor:duplicable-product]{Corollary}%
When the direct sum~$\scrA\oplus\scrB$
of von Neumann algebras~$\scrA$ and~$\scrB$
is duplicable,
$\scrA$ and~$\scrB$ are duplicable
\begin{point}{140}{Proof}%
Let $\delta\colon (\mathscr{A}\oplus\mathscr{B})\otimes
(\mathscr{A}\oplus\mathscr{B})\longrightarrow
\mathscr{A}\oplus\mathscr{B}$
be a duplicator on~$\mathscr{A}\oplus\mathscr{B}$.
By~\sref{lem:uniqueness-duplicator}
$\mathscr{A}\oplus\mathscr{B}$
is commutative
and $\delta((a_1,b_1)\otimes (a_2,b_2))
= (a_1a_2,b_1b_2)$
for all $a_1,a_2\in\mathscr{A}$
and~$b_1,b_2\in\mathscr{B}$.
Let~$\kappa_1\colon \mathscr{A}\to\mathscr{A}\oplus\mathscr{B}$
be the nmiu-map
given by~$\kappa_1(a)=(a,0)$ for all~$a\in\mathscr{A}$.
Let~$\delta_\mathscr{A}$ be the composition of
$\xymatrix@C=3em{
\mathscr{A}\otimes\mathscr{A}
\ar[r]|-{\kappa_1\otimes\kappa_1}
&
(\mathscr{A}\oplus\mathscr{B})
\otimes
(\mathscr{A}\oplus\mathscr{B})
\ar[r]|-{\delta}
&
\mathscr{A}\oplus\mathscr{B}
\ar[r]|-{\pi_1}
&
\mathscr{A}
}$.
Then~$\delta_\mathscr{A}$ is normal, positive,
and
$\delta_\mathscr{A}(a_1\otimes a_2)
\,=\,  \pi_1(\delta((a_1,0)\otimes (a_2,0))) 
\,=\, \pi_1(a_1a_2,0)\,=\,a_1a_2$
for all~$a_1,a_2\in\mathscr{A}$.
Thus, by~\sref{cor:duplicability-multiplication},
$\mathscr{A}$
is duplicable.\qed
\end{point}
\end{point}
\end{parsec}
\begin{parsec}{1290}%
\begin{point}{10}%
We will now work towards
the proof that
if~$C$ is a continuous complete finite measure space,
then
$L^\infty(C)$
cannot be duplicable 
unless~$\mu(C)=0$,
see~\sref{lem:continuous-finite-measure-space-not-duplicable}.
Let us first fix some more terminology
from measure theory (see~\sref{measure-theory-recap} and~\cite{fremlin}).
\end{point}
\begin{point}{20}{Definition}
\label{def:measure-space}
Let~$X$ be a finite complete measure space.
\begin{enumerate}
\item
A measurable subset~$A$ of~$X$
is \Define{atomic}
\index{atomic subset of a measure space}
if $0<\mu(A)$
and $\mu(A')=\mu(A)$ for all~$A'\in\Sigma_X$
with~$A'\subseteq A$ and~$\mu(A')>0$.

\item
$X$ is \Define{discrete} if $X$ is covered by atomic measurable subsets.%
\index{measure space!discrete}

(This coincides with being ``purely atomic''
from 211K of~\cite{fremlin}.)
\item
$X$ is \Define{continuous}%
\index{measure space!continuous}
(or ``atomless'')
if~$X$ contains no atomic subsets.
\end{enumerate}
\spacingfix%
\end{point}%
\begin{point}{30}%
The following lemma,
which will be very useful,
is a variation on
Zorn's Lemma
(that does not require the axiom of choice).
\end{point}
\begin{point}{40}[lem:measure-zorn]{Lemma}
Let~$\mathcal{S}$
be a collection of measurable subsets of
a finite complete measure space~$X$ 
such that 
for 
every ascending countable sequence
$A_1\subseteq A_2 \subseteq \dotsb$
in~$\mathcal{S}$
there is~$A\in\mathcal{S}$
with $A_1\subseteq A_2\subseteq\dotsb \subseteq A$.

Then each element~$A\in\mathcal{S}$
is contained in some~$B\in\mathcal{S}$
that is maximal in~$\mathcal{S}$ in the sense
that $\mu(B')=\mu(B)$
for all~$B'\in\mathcal{S}$ with $B\subseteq B'$.
\begin{point}{50}{Proof}
	The trick is to consider for every $C\in\mathcal{S}$
the quantity
\begin{equation*}
	\beta_C \ = \ \sup\{\,\mu(D)\colon\, C\subseteq D
\text{ and }D\in \mathcal{S}\,\}.
\end{equation*}
Note that $\mu(C)\leq \beta_C \leq \mu(X)$
	for all~$C\in\mathcal{S}$,
	and $\beta_{C_2}\leq \beta_{C_1}$
	for all~$C_1,C_2\in\mathcal{S}$ with $C_1 \subseteq C_2$.
	To prove this lemma, it suffices to find~$B\in\mathcal{S}$
with $A\subseteq B$ and~$\mu(B)=\beta_B$.

Define~$B_1:= A$.
Pick~$B_2\in\mathcal{S}$
such that $B_1 \subseteq B_2$
and~$\beta_{B_1}-\mu(B_2)\leq \nicefrac{1}{2}$.
Pick~$B_3\in\mathcal{S}$
such that $B_2\subseteq B_3$
and~$\beta_{B_2}-\mu(B_3) \leq \nicefrac{1}{3}$.
Proceeding in this fashion,
we get a sequence $B\equiv B_1\subseteq B_2 \subseteq \dotsb$
in~$\mathcal{S}$
with $\beta_{B_{n}}-\mu(B_{n+1})\leq \nicefrac{1}{n}$
for all~$n$.
By assumption there is a~$B\in\mathcal{S}$
with $B_1\subseteq B_2\subseteq \dotsb \subseteq B$.
Note that
\begin{equation*}
\mu(B_1)\,\leq\, \mu(B_2)\,\leq\,
\dotsb \,\leq\,\mu(B)\,\leq\, \beta_B \,\leq\, \dotsb
\,\leq\, \beta_{B_2}\,\leq\, \beta_{B_1}.
\end{equation*}
Since for every~$n\in\mathbb{N}$
we have both $\mu(B_{n+1})\leq \mu(B)\leq \beta_B \leq \beta_{B_n}$
and $\beta_{B_n}- \mu(B_{n+1}) \leq \nicefrac{1}{n}$,
we get $\beta_B-\mu(B)\leq \nicefrac{1}{n}$,
and so~$\beta_B = \mu(B)$.\qed
\end{point}%
\end{point}%
\begin{point}{60}[lem:measure-space-continuous-discrete]{Lemma}%
Each finite complete  measure space~$X$ 
contains a  discrete measurable  subset~$D$
such that~$X\backslash D$ is continuous.
\begin{point}{70}{Proof}%
Since clearly the countable union
of discrete measurable subsets of~$X$
is again discrete,
there is by~\sref{lem:measure-zorn}
a discrete measurable subset~$D$ of~$X$
which is maximal in the sense that~$\mu(D')=\mu(D)$
for every discrete measurable subset~$D'$ of~$X$ with $D\subseteq D'$.
To show that~$X\backslash D$ is continuous,
we must prove that~$X\backslash D$
contains no atomic measurable subsets.
If~$A\subseteq X\backslash D$ is an atomic measurable subset
of~$X$,
then~$D\cup A$
is a discrete measurable
subset of~$X$
which contains~$D$,
and $\mu(D\cup A)=\mu(D)\cup \mu(A) > \mu(D)$.
This contradicts the  maximality of~$D$.
Thus~$X\backslash D$ is continuous.\qed
\end{point}
\end{point}
\begin{point}{80}[lem:continuous-measure-space]{Lemma}%
Given a continuous finite complete measure space~$X$,
and~$r\in [0,\mu(X)]$,
there is a measurable subset~$A$ of~$X$ with $\mu(A)=r$.
\begin{point}{90}{Proof}%
Let us quickly get rid of the case that~$\mu(X)=0$.
Indeed, then~$r=0$, and so~$A=\varnothing$ will do.
For the remainder, assume that~$\mu(X)>0$.

For starters, we show that for every~$\varepsilon >0$
and~$B\in\Sigma_X$ with~$\mu(B)>0$
there is~$A\in\Sigma_X$ with $A\subseteq B$
and  $0<\mu(A)<\varepsilon$.
Define~$A_1 := B$.
Since~$\mu(B)>0$,
and~$A_1$ is not atomic (because~$X$ is continuous)
there is~$A\in\Sigma_X$ with $A\subseteq A_1$ 
and $\mu(A)\neq \mu(A_1)$.
Since~$\mu(A)+\mu(A_1\backslash A)=\mu(A_1)$,
either $0<\mu(A)\leq \frac{1}{2}\mu(A_1)$
or $0<\mu(X\backslash A)\leq \frac{1}{2}\mu(A_1)$.
In any case,
there is~$A_2\subseteq A_1$
with $A_2\in\Sigma_X$
and $0<\mu(A_2)\leq \frac{1}{2}\mu(A_1)$.
Similarly,
since~$A_2$ is not atomic (because~$X$ is continuous),
there is~$A_3\subseteq A_2$
with~$A_3\in\Sigma_X$ and $0<\mu(A_3)\leq \frac{1}{2}\mu(A_2)$.
Proceeding in a similar fashion,
we obtain a sequence $B\equiv A_1 \supseteq A_2\supseteq \dotsb$
of measurable subsets of~$X$
with $0<\mu(A_n)\leq 2^{-n}\mu(X)$.
Then, for every $\varepsilon >0$
there is~$n\in\mathbb{N}$
such that $0<\mu(A_n)\leq \varepsilon$ and~$A_n\subseteq B$.

Now, 
let us prove that there is~$A\in\Sigma_X$ with $\mu(A)=r$.
By~\sref{lem:measure-zorn}
there is a measurable
subset~$A$ of~$X$
with $\mu(A)\leq r$
and which is maximal
in the sense that $\mu(A')=\mu(A)$
for all~$A'\in\Sigma_X$
with $\mu(A)\leq r$ and~$A\subseteq A'$.
In fact, we claim that~$\mu(A)=r$.
Indeed, suppose that~$\varepsilon := r-\mu(A)>0$
towards a contradiction.
By the previous discussion,
there is~$C\in\Sigma_X$ with $C\subseteq X\backslash A$
such that $\mu(C)\leq \varepsilon$.
Then~$A\cup C$ is measurable,
and $\mu(A\cup C)=\mu(A)+\mu(C)\leq \mu(A)+\varepsilon\leq r$,
which contradicts the maximality of~$A$.\qed
\end{point}
\end{point}
\begin{point}{100}[lem:continuous-finite-measure-space-not-duplicable]{Lemma}%
Let~$X$ be a continuous finite complete measure space
for which~$L^\infty(X)$ is duplicable.
Then~$\mu(X)=0$.
\begin{point}{110}{Proof}%
Suppose that~$\mu(X)>0$
towards a contradiction.
Let~$\delta$
be a duplicator
on~$L^\infty(X)$.
By~\sref{lem:uniqueness-duplicator}
$\delta(\mathfrak{f}\otimes \mathfrak{g})=
\mathfrak{f}\cdot \mathfrak{g}$ for 
all~$\mathfrak{f},\mathfrak{g}\in L^\infty(X)$.

Let~$\omega\colon L^\infty(X)\to \mathbb{C}$
be given by~$\omega(f^\circ)=\frac{1}{\mu(X)}\int f \,d\mu$
for all~$f\in \mathcal{L}^\infty(X)$.
Then~$\omega$ is normal, positive, unital and faithful
(cf.~\sref{Linfty-vn}).
We'll use
the product
functional 
$\omega\otimes \omega\colon L^\infty(X)\otimes L^\infty(X)\to \mathbb{C}$,
(which is also faithful,
by~\sref{carrier-tensor})
to tease out a contradiction,
but first we  need a second ingredient.

Since~$X$ is continuous,
we may partition~$X$ into two measurable
subsets of equal measure 
with the aid of~\sref{lem:continuous-measure-space},
that is,
there are measurable subsets $X_{1}$ and~$X_{2}$
of~$X$ with $X=X_{1}\cup X_{2}$, $X_{1}\cap X_{2}=\varnothing$,
and
$\mu(X_{1})=\mu(X_{2})=\frac{1}{2}\mu(X)$.
Similarly, $X_{1}$ 
can be split into two measurable subsets, $X_{11}$ and $X_{12}$,
of equal measure, and so on.
In this way,
we obtain for every word~$w$ over the alphabet~$\{1,2\}$
--- in symbols, $w\in \{1,2\}^*$ ---
a measurable subset~$X_w$ of~$X$
such that $X_w = X_{w1}\cup X_{w2}$,
$X_{w1}\cap X_{w2}=\varnothing$,
and $\mu(X_{w1})=\mu(X_{w2})=\frac{1}{2}\mu(X_w)$.
It follows that~$\mu(X_w)=\frac{1}{2^{\#w}}\mu(X)$,
where~$\#w$ is the length of the word~$w$.

Now, $p_w := \mathbf{1}_{X_w}^\circ$ is a projection in~$L^\infty(X)$,
and~$\omega(p_w)=2^{-\#w}$
for every~$w\in\{1,2\}^*$.
Moreover, $p_w = p_{w1}+p_{w2}$,
and so
\begin{alignat*}{3}
p_w\otimes p_w 
\ &=\  
p_{w1}\otimes p_{w1} \,+\,
p_{w1}\otimes p_{w2} \,+\,
p_{w2}\otimes p_{w1} \,+\,
p_{w2}\otimes p_{w2}\\
\ &\geq\ 
p_{w1}\otimes p_{w1} \,+\,
p_{w2}\otimes p_{w2}.
\end{alignat*}
Thus, if we define 
$q_N\ :=\ \sum_{w\in \{1,2\}^N}\,p_w\otimes p_w$
for every natural number~$N$,
where~$\{1,2\}^N$ is the set of words over~$\{1,2\}$ of length~$N$,
then we get a descending sequence $q_1\geq q_2\geq q_3\geq \dotsb$
of projections in~$L^\infty(X)\otimes L^\infty(X)$.
Let~$q$ be the infimum of $q_1\geq q_2 \geq \dotsb$ 
in the set of self-adjoint elements of~$L^\infty(X)\otimes
L^\infty(X)$.
Do we have~$q=0$ ?

On the one hand,
we claim that $\delta(q)=1$, and so~$q\neq 0$.
Indeed,
$\delta(p_w\otimes p_w)=p_w\cdot p_w = p_w$
for all~$w\in \{1,2\}^N$.
Thus $\delta(q_N) = \sum_{w\in \{1,2\}^N}  \delta(p_w\otimes p_w)
= \sum_{w\in\{1,2\}^N} p_w=1$ for all~$N\in \mathbb{N}$.
Hence $\delta(q)=\bigwedge_n \delta(q_N) = 1$,
because~$\delta$ is normal.
On the other hand,
we claim that $(\omega\otimes \omega)(q)=0$,
and so~$q=0$ since~$\omega\otimes \omega$ is 
faithful and $q\geq 0$.
Indeed,
$(\omega\otimes\omega)(q_N)=
\sum_{w\in\{1,2\}^N} \omega(p_w)\cdot\omega(p_w)
= \sum_{w\in\{1,2\}^N} 2^{-N}\cdot 2^{-N} = 2^{-N}$
for all~$N\in \mathbb{N}$,
	and so $(\omega\otimes\omega)(q)
=\bigwedge_N (\omega\otimes\omega)(q_N) = \bigwedge_N 2^{-N}=0$.
Thus, $q=0$ and $q\neq 0$, which is impossible.\qed
\end{point}
\end{point}
\end{parsec}
\begin{parsec}{1300}%
\begin{point}{10}%
This takes care of the continuous case.
To deal with the discrete case we first need some
simple observations.
\end{point}
\begin{point}{20}[lem:atomic-measure-space]{Lemma}%
Let~$A$ be an atomic measure space.
Then~$L^\infty(A)\cong \mathbb{C}$.
\begin{point}{30}{Proof}
Let~$f\in \mathcal{L}^\infty(A)$ be given.
It suffices to show that
there is~$z\in \mathbb{C}$
such that
$f(x)=z$ for almost all~$x\in A$.
Moreover, we only need to consider the case
that~$f$ takes its values in~$\mathbb{R}$
(because we may split~$f$ in its real and imaginary parts,
and in turn split these in positive and negative parts).

Let~$S$ be some measurable subset of~$A$.
Note that either~$\mu(S)=0$ or~$\mu(A\backslash S)=0$.
Indeed, if not~$\mu(S)=0$,
then~$\mu(S)>0$,
and so~$\mu(S)=\mu(A)$
(by atomicity of~$A$),
which entails that~$\mu(A\backslash S)=0$.

In particular,
for every real number~$t\in\R$
one of the sets
\begin{equation*}
\{\,x\in A\colon\, t\leq f(x) \,\}
\qquad
\{\,x\in A\colon\, f(x)<t \,\}
\end{equation*}
must be negligible.  Whence
either~$t\leq f^\circ$ or~$f^\circ\leq t$.
It follows that
the two closed sets $L:=\{t\in\R\colon t\leq f^\circ\}$
and $U:=\{t\in \R\colon f^\circ\leq t\}$
cover~$\R$.
Since clearly~$-\|f\| \in L$
and~$\|f\|\in  U$,
the sets~$L$ and~$U$
can't be disjoint,
because they would partition~$\R$ into two clopen non-empty sets.
For an element~$t\in L\cap U$
in the intersection
we have~$ t\leq f^\circ \leq t$,
	and so~$t=f^\circ$. Hence~$L^\infty(X)\cong \C$. \qed
\end{point}
\end{point}
\begin{point}{40}[lem:measure-space-partition]{Exercise}%
Let~$X$ be a measure space with~$\mu(X)<\infty$.
Show that~$L^\infty(X)\cong \bigoplus_{A\in\mathcal{A}} L^\infty(A)$
for every countable partition~$\mathcal{A}$
of~$X$
consisting of measurable subsets.
\end{point}
\begin{point}{50}[cor:discrete-ell-x]{Corollary}%
For every discrete 
measure space~$X$ with~$\mu(X)<\infty$
there is a  set~$Y$ with $L^\infty(X)\cong \linf(Y)$.
\end{point}
\end{parsec}
\begin{parsec}{1310}%
\begin{point}{10}%
We are now ready to prove
the main result of this section.
\end{point}
\begin{point}{20}[duplicable-proof]{Proof of~\ref{duplicable}}%
We have already seen that $\linf(X)$
can be  equipped with a commutative monoid
structure in~$\W{miu}$
for any set~$X$,
and is thus duplicable.
Conversely,
let~$\delta\colon \mathscr{A}\otimes\mathscr{A}\to\mathscr{A}$
be a duplicator with unit~$u$ on a von Neumann algebra~$\mathscr{A}$.
By~\sref{lem:unit-duplicator}, we know that~$u=1$,
and by~\sref{lem:uniqueness-duplicator},
we know that~$\mathscr{A}$
is commutative 
and~$\delta(a\otimes b)=a\cdot b$
for all~$a,b\in \mathscr{A}$.
Thus, the only thing that remains to be shown
is that~$\mathscr{A}$ is miu-isomorphic to $\linf(Y)$
for some set~$Y$.
By~\sref{cvn}
$\scrA\cong \bigoplus_i L^\infty(X_i)$
for some finite complete measure spaces~$X_i$.
So to prove that~$\mathscr{A}\cong \ell^\infty(Y)$
for some set~$Y$ it 
suffices
to find
 a set~$Y_i$
with $L^\infty(X_i)\cong \ell^\infty(Y_i)$
for each~$i$,
because then 
$\textstyle \mathscr{A}\ \cong \ 
\bigoplus_{i\in I} \ell^\infty(Y_i)\ \cong\ 
\ell^\infty\bigl(\,\bigcup_{i\in I} Y_i\,\bigr)$.

Let~$i\in I$ be given.
Since~$\mathscr{A}\cong L^\infty(X_i)\,\oplus\,\bigoplus_{j\neq i} 
L^\infty(X_j)$ is duplicable,
$L^\infty(X_i)$ is duplicable
by~\sref{cor:duplicable-product}.
By~\sref{lem:measure-space-continuous-discrete}
there is a measurable subset~$D$ of~$X_i$ such that~$D$
is discrete, and $C:=X\backslash D$ is continuous.
We have~$L^\infty(X_i)\cong L^\infty(D)\oplus L^\infty(C)$
by~\sref{lem:measure-space-partition},
and
so $L^\infty(D)$ and~$L^\infty(C)$
are duplicable
(again by~\sref{cor:duplicable-product}).
By~\sref{lem:continuous-finite-measure-space-not-duplicable},
$L^\infty(C)$
can only be duplicable if~$\mu(C)=0$,
and so~$L^\infty(C)\cong \{0\}$.
On the other hand,
since~$D$ is discrete,
we have~$L^\infty(D)\cong \ell^\infty(Y)$
for some set~$Y$
(by~\sref{cor:discrete-ell-x}).
All in all, we have $L^\infty(X_i)\cong \ell^\infty(Y)$.\qed
\end{point}
\end{parsec}
\subsection{Monoids}
\label{sec:monoids-in-vna}
\begin{parsec}{1320}%
\begin{point}{10}%
We further justify our choice,
$\sem{!A} = \linf(\nsp(\sem{A}))$,
by proving that $\linf(\nsp(\mathscr{A}))$
is the free (commutative) monoid on~$\mathscr{A}$ in~$\W{miu}$.
As a corollary, we also obtain that $\linf(\W{cpsu}(\mathscr{A},\mathbb{C}))$
is the free (commutative) monoid on~$\mathscr{A}$
in~$\W{cpsu}$.
\end{point}
\begin{point}{20}%
Let us first recall some terminology.
Given a symmetric monoidal
	category (\Define{SMC})%
\index{symmetric monoidal category (SMC)}
$\Cat{C}$, 
	a \Define{monoid} \index{monoid!in an SMC}
	in $\Cat{C}$ is an object $A$ from~$\Cat{C}$ endowed
with
a \Define{multiplication} map
$m\colon A\otimes A\to A$
and a \Define{unit} map $u\colon I\to A$
satisfying the associativity and the unit law,
i.e.~making the following diagrams commute.
\[
\xymatrix@R-.5pc{
(A\otimes A)\otimes A
\ar[d]_{\alpha}
\ar[rr]^-{m\otimes \id}
&&
A\otimes A
\ar[d]^{m}
\\
A\otimes (A\otimes A)
\ar[r]_-{\id\otimes m}
&
A\otimes A
\ar[r]_-{m}
&
A
}
\qquad
\xymatrix@R-.5pc{
I\otimes A
\ar[dr]_{\lambda}
\ar[r]^-{u\otimes \id}
&
A\otimes A
\ar[d]^{m}
&
\ar[l]_-{\id\otimes u}
A\otimes I
\ar[dl]^{\rho}
\\
&
A
&
}
\]
Here $\alpha,\lambda,\rho$ respectively
denote the associativity isomorphism, and
the left and the right unit isomorphism.
A monoid $A$ is \Define{commutative} if
$m\circ \gamma=m$,
where $\gamma\colon A\otimes A\to A\otimes A$ is the symmetry isomorphism.
A \Define{monoid morphism} between monoids $A_1$
and $A_2$ is an arrow $f\colon A_1\to A_2$
that satisfies $m_{A_2}\circ (f\otimes f)=f \circ m_{A_1}$
and $u_{A_2}=f \circ u_{A_1}$.
We denote the category of monoids
and monoid morphisms in $\Cat{C}$
by $\Define{\Mon(\Cat{C})}$.
The full subcategory of commutative
monoids is denoted by~$\Define{\CMon(\Cat{C})}$.
Recall that~$\W{miu}$ and~$\W{cpsu}$
are symmetric monoidal categories
with~$\C$ as tensor unit (see~\sref{vn-smc}),
and so we may speak about monoids
in~$\W{miu}$ and~$\W{cpsu}$.
\end{point}
\begin{point}{30}[prop:dup-vna-is-monoid]{Exercise}%
Let~$\scrA$ be a von Neumann algebra.
\begin{enumerate}
\item
Show that any monoid structure on~$\scrA$
in~$\W{cpsu}$
is a duplicator on~$\scrA$.
\item
Deduce
from this and~\sref{duplicable}
that there is a monoid structure
on~$\scrA$ in~$\W{miu}$ or in~$\W{cpsu}$
iff~$\scrA$ is duplicable
iff~$\scrA\cong \ell^\infty(X)$ for some set~$X$;
and that,
in that case
the multiplication $m\colon \scrA\otimes\scrA\to\scrA$
of the monoid
is commutative
and uniquely being fixed by~$m(a\otimes b)=a\cdot b$.
\item
Show that the monoid morphisms
in~$\W{miu}$ and in~$\W{cpsu}$
are precisely
the nmiu-maps.
\item
Conclude that
$\CMon(\W{miu})
=\Mon(\W{miu})
=\CMon(\W{cpsu})
=\Mon(\W{cpsu})$.
\item
Show  that $\Mon(\W{miu})\cong\dW{miu}\simeq\op{\Cat{Set}}$,
where $\Define{\dW{miu}}$ denotes the full subcategory
of~$\W{miu}$
consisting of duplicable von Neumann algebras.

(Hint: $\linf\colon \Cat{Set}\to \op{(\W{miu})}$ is full and faithful
by~\sref{cor:linf-ff}.)
\end{enumerate}
\spacingfix%
\end{point}%
\begin{point}{40}[thm:free-monoid-in-vNAMIU]{Theorem}%
\index{free monoid!in~$\W{miu}$}
Let~$\mathscr{A}$
be a von Neumann algebra,
and let~$\eta\colon \mathscr{A}\to\linf(\nsp(\mathscr{A}))$
be the nmiu-map
given by~$\eta(a)(\varphi)= \varphi(a)$.
Then $\linf(\nsp(\mathscr{A}))$
is the free (commutative) monoid
on~$\mathscr{A}$
in~$\W{miu}$ via~$\eta$.
\begin{point}{50}{Proof}%
Let~$\mathscr{B}$
be a monoid 
on~$\W{miu}$,
and let~$f\colon \mathscr{A}\to\mathscr{B}$
be an nmiu-map.
We must show that
there is a unique
monoid morphism
$g\colon \linf(\nsp(\mathscr{A}))
\rightarrow \mathscr{B}$
such that~$g\circ \eta = f$.
Since the monoid structure on~$\scrB$
is a duplicator on~$\scrB$
we may assume,
by~\sref{duplicable},
that~$\mathscr{B}=\linf(Y)$
for some set~$Y$.
Since~$\nsp\colon \op{(\W{miu})}\to \Cat{Set}$
is left adjoint
to~$\linf\colon \Cat{Set} \to \op{(\W{miu})}$
with unit~$\eta$ (see~\sref{first-adjunction}),
there is a unique map $h\colon Y\to \nsp(\mathscr{A})$
with $\linf(h)\circ \eta = f$.
Since~$\linf$ is full and faithful
by~\sref{cor:linf-ff},
the only thing that remains to be shown is that~$\linf(h)$
is a monoid morphism.
Indeed it is,
since the monoid multiplication
on~$\linf(\nsp(\mathscr{A}))$
and~$\linf(Y)$
is given by ordinary multiplication,
which is preserved by~$\linf(h)$ being an miu-map.\qed
\end{point}
\end{point}
\begin{point}{60}{Corollary}
Let~$\mathscr{A}$
be a von Neumann algebra.
Then $\linf(\W{cpsu}(\mathscr{A},\C))$
is the free (commutative) monoid
on~$\mathscr{A}$ in~$\W{cpsu}$.%
\index{free monoid!in~$\W{cpsu}$}
\begin{point}{70}{Proof}
By~\sref{thm:free-monoid-in-vNAMIU}
$\linf\circ\nsp$ is a left adjoint to
the forgetful functor $\Mon(\W{miu})\to\W{miu}$.
Note that
by~\sref{prop:dup-vna-is-monoid},
the forgetful functor $\Mon(\W{cpsu})\to\W{cpsu}$
factors through $\W{miu}$ as:
\[
\xymatrix{
\Mon(\W{cpsu})
\ar@{=}[r]
&
\Mon(\W{miu})
\ar[r]^-{\bot}&
\ar@/_3ex/[l]_{\linf\circ\nsp}
\W{miu}
\ar@{->}[r]^-{\bot}&
\ar@/_3ex/[l]_{\mathcal{F}}
\W{cpsu}
}
\]
where~$\mathcal{F}$ is from~\sref{second-adjunction}.
Thus the free monoid on $\scrA$ in $\W{cpsu}$ is given by:
\[
(\linf\circ\nsp\circ\mathcal{F})(\scrA)
\ =\ 
\linf(\W{miu}(\mathcal{F}\scrA,\C))
\ \cong\ 
\linf(\W{cpsu}(\scrA,\C))
\]
as was claimed.\qed
\end{point}
\end{point}
\end{parsec}

\begin{parsec}{1330}%
\begin{point}{10}[conclusion]{Conclusion}
Here ends this thesis,
but not the entire story.
There's much more to be said
about self-dual Hilbert $\scrA$-modules,
about dilations and their relation to purity,
and about the abstract theory of corners, filters,
and $\diamond$-positivity.
You'll see all this,
and more,
in the sequel,
``Dagger and dilations in the category of von Neumann algebras''\cite{bas},
brought to you by my twin brother.
\end{point}
\end{parsec}
\begin{parsec}{1340}[final-bram]
\begin{point}{10}
(Paragraphs numbered~\sref{final-bram}
and up can be found in~\cite{bas}.)
\end{point}
\end{parsec}

\backmatter

\fancyfoot[CE]{}
\fancyfoot[CO]{}
\fancypagestyle{plain}{
    \fancyfoot[CE]{}
    \fancyfoot[CO]{}
}

\printindex

\begingroup
\renewcommand\chapter[2]{\backmattertitle{Bibliography}}
\bibliography{main}{}
\endgroup

\bibliographystyle{plain}

\oldchapter{Lekensamenvatting}

Wat is wiskunde? Volgens mij is ze niets anders dan een studie van de patronen
in de wereld om ons heen in de meest \emph{strenge} en 
\emph{formele} zin.
De wiskundige
maakt zich  bewust, bijna komisch, blind 
voor zaken die een ander z\'o in het oog zouden springen.
``Ja.'' antwoorden op de vraag ``Wil je koffie of thee?'' is voor
een wiskundige volstrekt acceptabel.
Maar door deze blindheid kan de wiskundige wel
met een onge\"evenaarde nauwkeurigheid,
zekerheid en tijdloosheid
uitdrukking geven aan afzonderlijke  aspecten van de werkelijkheid.

Slechts een zeer beperkt aantal fenomenen
is vatbaar voor wiskundige \hyphenation{ana-lyse}analyse.
Soms begrijpen we een fenomeen er nog niet goed genoeg voor.
Soms is de benodigde wiskunde er nog niet voor ontwikkeld.
Maar in de meeste gevallen is het zo dat
het mes van de formele wiskunde simpelweg te scherp snijdt. 
Een wiskundige zou je niet kunnen vertellen
of een schip hetzelfde blijft als je alle planken vervangt,\footnote{Het schip
van Theseus, zie~\S23 van~\cite{theseus}.}
wat \emph{twee} is w\'el.
Het is niet zozeer dat hij of zij zich niet bezig zou willen houden
met dagelijkse begrippen, maar eerder dat deze begrippen daartegen niet
bestand zijn.
Daarom gaat de wiskunde over abstracte en ge\"idealiseerde 
objecten zoals getallen, vierkanten, cirkels, distributies,
tensoren en vari\"eteiten,
in plaats van wereldse tegenhangers zoals
respectievelijk
logistiek, akkers, hemellichamen, toekomstverwachtingen,
elektrische velden en sterrenstelsels.

Het is wonderbaarlijk hoe behulpzaam
wiskunde kan zijn ondanks en juist dankzij haar
beperkingen.
Niemand---ook geen wiskundige---kan zich het stilleven 
op de omslag van dit boekje vanuit een andere hoek nauwkeurig voorstellen.
Maar als ik vertel dat het een plaatje is van
een reflecterende bol in~$\R^{12}$
(tussen twee afgeknotte hypervlakken, met een deels reflecterend 
schaakbordpatroon) en wanneer ik 
de relevante co\"ordinaten en afmetingen geef, dan
kan ieder die \hyphenation{be-kend}bekend is met `inproducten' (vgl.~\sref{hilb-def})
een computerprogramma schrijven dat
dit tafereel vanuit een andere hoek toont.
Dankzij de wiskunde
is bovenmenselijk inzicht in de twaalfde dimensie niet nodig!

Dat sommige wiskundige begrippen zeer bruikbaar zijn,
wil niet zeggen dat \emph{elk} wiskundig begrip dit is.
Integendeel zelfs:
van alle denkbare figuren op papier
grijpt men ongetwijfeld het vaakst terug
naar de cirkel en de lijn.
Bepaalde constanten ($2\pi$, $e$, $\varphi$, $\sqrt{2}$, \ldots) 
verschijnen ook vaker in formules dan andere.
Natuurlijk kan een nuttig wiskundig begrip 
ook slecht gebruikt worden.
Men kan bijvoorbeeld 
de baan van een planeet met een groot
aantal cirkels (Ptolemaeus' `epicykels') 
omschrijven in plaats van met \'e\'en ellips.
Bij een dergelijk waardeoordeel als `goed' of `slecht' gebruik
van een begrip is er zeker sprake van persoonlijke
voorkeur, conventie en willekeur. Waarom gebruiken
we niet $2\pi=6,\!2831\dotsc$
als constante in plaats van
$\pi=3,\!1415\dotsc$? Volgens mij zou
dat veel formules eleganter maken;
zie~\cite{palais2001pi}.
Dat een wiskundig begrip (zoals epicykels)
in ongebruik geraakt is, wil bovendien niet zeggen
dat het geen \hyphenation{stijl-vol}stijlvol herintreden kan maken:
de befaamde Fourier-reeks (die in feite bestaat uit oneindig veel epicykels)
geeft de ellips weer 
het nakijken!\footnote{Of toch niet,  \cite{mayer2010drawing}?}

Dit proefschrift gaat over zo'n bruikbaar, abstract en 
ge\"idealiseerd wiskundig begrip:
de \emph{von Neumann-algebra},
bedacht door en vernoemd naar het Hongaarse genie
en de alleskunner
John von Neumann (geboren als Neumann J\'anos Lajos).
Aan hem hebben we niet alleen
de moderne computer-architectuur te danken,
maar bijvoorbeeld ook de springstoflenzen die in kernwapens gebruikt worden,
de afschrikwekkende strategie van \emph{mutual assured destruction}
en de numerieke weersvoorspelling.

Het is niet eenvoudig om uit te leggen
 wat von Neumann-algebra's precies zijn (\sref{vna}),
 maar 
ik kan wel een indruk geven waarvoor ik ze gebruik.
Uitgangspunt van dit proefschrift is dat 
de elementaire
systemen
die je tegenkomt bij het ontwerpen van een algoritme
voor een kwantumcomputer omschreven
kunnen worden door von Neumann-algebra's.
Zo wordt de klassieke bit (die~$0$ of~$1$ kan zijn)
voorgesteld door de von Neumann-algebra ``$\C^2$''
en wordt de kwantumbit (die een complexe combinatie van~$0$ of~$1$ is totdat
je haar meet)
voorgesteld door de von Neumann-algebra~``$M_2$''.
Dit geeft de hoop dat ingewikkeldere samengestelde systemen
ook een bijpassende von Neumann-algebra hebben.
Dat hangt natuurlijk af van de wijze van samenstelling:
het systeem dat bestaat uit twee kwantumbits en \'e\'en bit wordt voorgesteld
met de von Neumann-algebra $M_2\otimes M_2 \otimes \C^2$,
terwijl het systeem dat een kwantumbit of een klassieke bit bevat 
omschreven wordt met de von Neumann-algebra $M_2 \oplus \C^2$.
De bewerking ``$\otimes$'' op von Neumann-algebra's
is de interpretatie voor de  ``en''-samenstelling van systemen, 
terwijl ``$\oplus$''  invulling geeft aan ``of''-samenstelling van systemen.
Een veel complexere samenstelling van twee systemen~$A$ en~$B$
bestaat uit alle `processen' van~$A$ naar~$B$.
Meting is bijvoorbeeld zo'n proces van een kwantumbit naar een klassieke bit.
Zulke processen worden in de wereld van von Neumann-algebra's
voorgesteld door zogenaamde ``ncpsu-afbeeldingen'' (\sref{maps})
tussen von Neumann-algebra's.

E\'en van de hoofdresultaten van dit proefschrift
is een interpretatie voor deze `processen'-samenstelling,
de bewerking~$\limp$.
Door bepaalde formele kaders
was het al van tevoren duidelijk dat er hoogstens \'e\'en
interpretatie mogelijk zou zijn;
de vraag was alleen: welke?  Sterker nog: is er \"uberhaupt eentje?
Dit is vergelijkbaar met de vraag wat de kleinste\footnote{Kleinste
in de woordenboekordening.}
grammaticaal
correcte tekst is waarin
alle Nederlandse woorden voorkomen.
Om dit probleem op te lossen kun je niet zomaar alle woorden
op een rij zetten---de onderlinge samenhang moet immers kloppen.
De crux was voor mij om niet te willen proberen om de bewerking~$\limp$ 
direct te omschrijven,
maar om het bestaan ervan indirect aan te tonen,
zoals je ook kunt laten zien dat de  bovengenoemde Nederlandse tekst
uit het voorbeeld bestaat,
zonder te weten hoe hij precies is samengesteld.
(Voor een beperkte klasse van
von Neumann-algebra's, de `hereditair atomische', bleek een
directe omschrijving trouwens wel mogelijk, \sref{AstarhaB-concrete}.)

Het tweede hoofdresultaat van dit proefschrift
is de vondst van een abstracte omschrijving
van de processen
die bij een meting horen (\sref{sqrt-axiom}).
Het \hyphenation{ken-mer-ken-de}kenmerkende aspect van deze omschrijving is dat
het alleen gebruik maakt van zogenaamd `categorisch' jargon.
Dat dit mogelijk is toont aan
dat we in de categorie van von Neumann-algebra's
niet op een te abstract niveau werken:
we kunnen het nog steeds over meting hebben.
Dankzij de  categorische omschrijving  wordt het bovendien
mogelijk `meting' in andere contexten  te interpreteren.
Bij de zoektocht naar deze categorische omschrijving
was lange tijd het probleem om
een manier te vinden
om onderscheid te maken tussen de processen\footnote{Hier is~$u\neq 1$
een  unitaire
en 
$p$ een positief element van een von Neumann-algebra~$\scrA$
met $\ceil{p}=1$.}
\begin{equation*}
\sqrt{p}(\,\cdot\,)\sqrt{p}\qquad\text{en}\qquad
\sqrt{p}u^*(\,\cdot\,)u \sqrt{p}.
\end{equation*}
De linker hoort bij een meting, de rechter niet.
De oplossing was om een begrip uit de theorie 
van Hilbert-ruimten---geadjungeerdeerdheid---in een 
afgezwakte vorm over te nemen, namelijk $\diamond$-geadjungeerdheid 
(of \emph{contraposedness}, \sref{contraposed}).
Het linker proces blijkt een kwadraat te zijn van een aan zichzelf
$\diamond$-geadjungeerd proces; de rechter niet.
Beide hoofdresultaten staan in het laatste hoofdstuk.

De rest van dit proefschrift bestaat---enigszins ongebruikelijk---uit 
een grondige introductie
tot de benodigde, reeds bestaande, theorie
van $C^*$-algebra's en van von Neumann-algebra's.
Niet alleen bestond
een geschikte introductie nog niet, maar het leek me ook een goede kans 
om me verder te verdiepen in de theorie van von Neumann-algebra's.
In de eerste hoofdstukken ontwikkel ik de gehele benodigde theorie, inclusief 
bewijzen.
Het is mijn bedoeling dat eenieder die een bachelorgraad in de wiskunde heeft 
behaald deze tekst zou moeten kunnen begrijpen.
Op een groot aantal plekken wijk ik af van 
het begane pad:
soms om de tekst kort te houden, (zo ontwijk ik de
theorie van Banach algebra's volkomen),
maar meestal om te experimenteren met variaties.
Zo gebruik ik Kadison's omschrijving
van von Neumann-algebra's om
de theorie op te bouwen (wat niet eerder gedaan is.)
Zo houdt men goede wiskundige begrippen levend.

\oldchapter{About the Author}
Bram Westerbaan, born August 30, 1988,
enrolled as a physics and astronomy student at 
the Radboud University in 2006.
He obtained a bachelor's degree in mathematics in 2012 (cum laude)
with a thesis
on the simply typed $\lambda$-calculus
supervised by prof.~dr.~H.P.~Barendregt
(resulting in a publication~\cite{brambachelor}.)
The same year,
he obtained a master's degree in mathematics (summa cum laude)
with a thesis\cite{brammaster} on measure and integral
under supervision of prof.~dr.~A.C.M.~van Rooij.
His doctoral studies started in 2013
on the topic of (co)algebra guided by dr.~A.~Silva,
and shifted 
to the material presented in this thesis
in~2014
under the auspices of prof.~dr.~B.P.F.~Jacobs.
From 2018 onward Bram works as a postdoctoral researcher
at the digital security department on
an NWA project applying polymorphic pseudonymisation
to network traffic flow data.

\message{Rerun to get citations correct.}

\end{document}